\documentclass[10pt]{memoir}

\settrims{0pt}{0pt} % page and stock same size
\settypeblocksize{190mm}{110mm}{*} % {height}{width}{ratio}
\setlrmargins{*}{*}{1} % {spine}{edge}{ratio}
\setulmarginsandblock{1in}{1in}{*} % height of typeblock computed
\setheadfoot{\onelineskip}{2\onelineskip} % {headheight}{footskip}
\setheaderspaces{*}{1.5\onelineskip}{*} % {headdrop}{headsep}{ratio}
\checkandfixthelayout

\chapterstyle{bianchi}
\newcommand{\titlefont}{\normalfont\Huge\bfseries}

\usepackage{amsthm}
\usepackage{mathtools}

\usepackage{imakeidx}
\usepackage[framemethod=tikz]{mdframed}
\usepackage{footnote}
\usepackage{tablefootnote}

\usepackage[inline]{enumitem}
\setlistdepth{6}
\usepackage{ifthen}
\usepackage[utf8]{inputenc} %allows non-ascii in bib file
\usepackage{xcolor}

\usepackage{subfiles}
\usepackage[backend=biber, backref=true, maxbibnames = 10, style = alphabetic]{biblatex}
\usepackage{xr-hyper}
\usepackage[bookmarks=true,hidelinks,pdfencoding=unicode]{hyperref}
\usepackage[capitalize]{cleveref}

\usepackage{multirow}

\usepackage{xstring}

\usepackage{answers}

\usepackage{tikz}
\usepackage{varwidth}
\usepackage[prefix=tikzsym]{tikzsymbols}
\usepackage{makecell}%database table thickness
\usepackage{multirow}

\usepackage{amssymb}
\usepackage{newpxtext}
\usepackage[varg,bigdelims]{newpxmath}
\usepackage{mathrsfs}
\usepackage{dutchcal}
\usepackage{fontawesome}

\DeclareMathVersion{normal2}

\graphicspath{graphics/}

  \crefformat{enumi}{\##2#1#3}
  \crefmultiformat{enumi}{\##2#1#3}{ and~\##2#1#3}{, \##2#1#3}{, and~\##2#1#3}

  \makeindex

  \hypersetup{final}

  \setlist{nosep}

	\makesavenoteenv{tabular}

  \usetikzlibrary{
  	cd,
  	math,
  	decorations.markings,
		decorations.pathreplacing,
  	positioning,
  	arrows.meta,
  	shapes,
		shadows,
		shadings,
  	calc,
  	fit,
  	quotes,
  	intersections,
    circuits,
    circuits.ee.IEC
  }

	\tikzcdset{arrow style=tikz, diagrams={>=To}}

	\pgfdeclarelayer{edgelayer}
\pgfdeclarelayer{nodelayer}
\pgfsetlayers{edgelayer,nodelayer,main}

\tikzset{
	tick/.style={postaction={
  	decorate,
    decoration={markings, mark=at position 0.5 with
    	{\draw[-] (0,.4ex) -- (0,-.4ex);}}}
  }
}

\tikzset{
	bottom base/.style={baseline=(current bounding box.south)}
}

\tikzset{trees/.style={
	inner sep=0,
	minimum width=0,
	minimum height=0,
	level distance=.5cm,
	sibling distance=.5cm,
	edge from parent/.style={shorten <= -2pt, draw, ->},
	grow'=up,
	decoration={markings, mark=at position 0.75 with \arrow{stealth}}
	}
}

\newcommand{\idchild}{edge from parent[double, -]}

  \tikzset{
     oriented WD/.style={%everything after equals replaces "oriented WD" in key.
        every to/.style={out=0,in=180,draw},
        label/.style={
           font=\everymath\expandafter{\the\everymath\scriptstyle},
           inner sep=0pt,
           node distance=2pt and -2pt},
        semithick,
        node distance=1 and 1,
        decoration={markings, mark=at position \stringdecpos with \stringdec},
        ar/.style={postaction={decorate}},
        execute at begin picture={\tikzset{
           x=\bbx, y=\bby,
           }}
        },
     string decoration/.store in=\stringdec,
     string decoration={\arrow{stealth};},
     string decoration pos/.store in=\stringdecpos,
     string decoration pos=.7,
	 	 dot size/.store in=\dotsize,
	   dot size=3pt,
	 	 dot/.style={
			 circle, draw, thick, inner sep=0, fill=black, minimum width=\dotsize
	   },
     bbx/.store in=\bbx,
     bbx = 1.5cm,
     bby/.store in=\bby,
     bby = 1.5ex,
     bb port sep/.store in=\bbportsep,
     bb port sep=1.5,
     bb port length/.store in=\bbportlen,
     bb port length=4pt,
     bb penetrate/.store in=\bbpenetrate,
     bb penetrate=0,
     bb min width/.store in=\bbminwidth,
     bb min width=1cm,
     bb rounded corners/.store in=\bbcorners,
     bb rounded corners=2pt,
     bb small/.style={bb port sep=1, bb port length=2.5pt, bbx=.4cm, bb min width=.4cm,
bby=.7ex},
		 bb medium/.style={bb port sep=1, bb port length=2.5pt, bbx=.4cm, bb min width=.4cm,
bby=.9ex},
     bb/.code 2 args={%When you see this key, run the code below:
        \pgfmathsetlengthmacro{\bbheight}{\bbportsep * (max(#1,#2)+1) * \bby}
        \pgfkeysalso{draw,minimum height=\bbheight,minimum width=\bbminwidth,outer
sep=0pt,
           rounded corners=\bbcorners,thick,
           prefix after command={\pgfextra{\let\fixname\tikzlastnode}},
           append after command={\pgfextra{\draw
              \ifnum #1=0{} \else foreach \i in {1,...,#1} {
                 ($(\fixname.north west)!{\i/(#1+1)}!(\fixname.south west)$) +(-
\bbportlen,0)
  coordinate (\fixname_in\i) -- +(\bbpenetrate,0) coordinate (\fixname_in\i')}\fi
              \ifnum #2=0{} \else foreach \i in {1,...,#2} {
                 ($(\fixname.north east)!{\i/(#2+1)}!(\fixname.south east)$) +(-
\bbpenetrate,0)
  coordinate (\fixname_out\i') -- +(\bbportlen,0) coordinate (\fixname_out\i)}\fi;
           }}}
     },
     bb name/.style={append after command={\pgfextra{\node[anchor=north] at
(\fixname.north) {#1};}}}
  }

\tikzset{polybox/.style={
	poly/.style ={
  	rectangle split,
  	rectangle split parts=2,
		rectangle split part align={bottom},
  	draw,
  	anchor=south,
  	inner ysep=4.5,
	  prefix after command={\pgfextra{\let\fixname=\tikzlastnode}},
		append after command={\pgfextra{
			\node[inner xsep=0, inner ysep=0,
				fit=(\fixname.north west) (\fixname.two split east)]
				(\fixname_dir) {};
			\node[inner xsep=0, inner ysep=0,
				fit=(\fixname.south west) (\fixname.two split east)]
				(\fixname_pos) {};
			}}
	},
	dom/.style={
		rectangle split part fill={none, blue!10}
	},
	cod/.style={
		rectangle split part fill={blue!10, none}
	},
	constant/.style={
		rectangle split part fill={red, none}		%rectangle split part fill={red!50, none}
	},
	constant dom/.style={
		rectangle split part fill={red, blue!10}
	},
	terminal/.style={
		rectangle split part fill={red, gray}
	},
	identity/.style={
		rectangle split part fill={gray, gray}
	},
	linear/.style={
		rectangle split part fill={gray, none}
	},
	linear dom/.style={
		rectangle split part fill={gray, blue!10}
	},
	linear cod/.style={
		rectangle split part fill={gray, none}
	},
	pure/.style={
		rectangle split part fill={none, gray}
	},
	pure dom/.style={
		rectangle split part fill={none, gray}
	},
	pure cod/.style={
		rectangle split part fill={blue!10, gray}
	},
	shorten <=3pt, shorten >=3pt,
	first/.style={out=0, in=180},
	climb/.style={out=180, in=180, looseness=2},
	climb'/.style={out=0, in=0, looseness=2},
	last/.style={out=180, in=0},
	mapstos/.style={arrows={|->}},
	tos/.style={arrows={->}},
	font=\footnotesize,
	node distance=2ex and 1.5cm
}
}

\tikzset{
biml/.tip={Glyph[glyph math command=triangleleft, glyph length=.95ex]},
bimr/.tip={Glyph[glyph math command=triangleright, glyph length=.95ex]},
}

\newcommand{\bimodfrom}[1][]{
	\begin{tikzcd}[ampersand replacement=\&, cramped]\ar[r, bimr-biml, "#1"]\&~\end{tikzcd}
}

\newcommand{\treepic}{
\begin{tikzpicture}[trees, scale=1.5,
  level 1/.style={sibling distance=20mm},
  level 2/.style={sibling distance=10mm},
  level 3/.style={sibling distance=5mm},
  level 4/.style={sibling distance=2.5mm},
  level 5/.style={sibling distance=1.25mm}]
  \node[dgreen] (a) {$\bullet$}
    child {node[dgreen] {$\bullet$}
    	child {node[dgreen] {$\bullet$}
    		child {node[dgreen] {$\bullet$}
  				child {node[dgreen] {$\bullet$}
    				child {}
    				child {}
    			}
  				child {node[dyellow] {$\bullet$}
    				child {}
    				child {}
    			}
  			}
    		child {node[dyellow] {$\bullet$}
					child {node[dgreen] {$\bullet$}
      			child {}
      			child {}
     			}
    			child  {node[red] {$\bullet$}}
  			}
    	}
    	child {node[dyellow] {$\bullet$}
    		child {node[dgreen] {$\bullet$}
  				child {node[dgreen] {$\bullet$}
    				child {}
    				child {}
    			}
  				child {node[dyellow] {$\bullet$}
    				child {}
    				child {}
    			}
  			}
    		child  {node[red] {$\bullet$}}
    	}
    }
    child {node[dyellow] {$\bullet$}
    	child {node[dgreen] {$\bullet$}
    		child {node[dgreen] {$\bullet$}
  				child {node[dgreen] {$\bullet$}
    				child {}
    				child {}
    			}
  				child {node[dyellow] {$\bullet$}
    				child {}
    				child {}
    			}
  			}
    		child {node[dyellow] {$\bullet$}
					child {node[dgreen] {$\bullet$}
      			child {}
      			child {}
     			}
    			child  {node[red] {$\bullet$}}
  			}
  		}
  		child {node[red] {$\bullet$}
  		}
  	}
  ;
\end{tikzpicture}
}

\newcommand{\coverpic}{
\begin{tikzpicture}
	\node {\treepic};
	\draw[blue, densely dotted] (current bounding box.south west) rectangle
(current bounding box.north east);
\end{tikzpicture}
}%{\earpic}

\tikzset{
    mapsto/.style={->, densely dashed, shorten <=\short, shorten
        >=\short, >=stealth, thick},
    short/.store in=\short,
    short=0pt
}

\makeatletter
\AfterEndEnvironment{mdframed}{%
  \tfn@tablefootnoteprintout%
  \gdef\tfn@fnt{0}%
}
\makeatother

\newlist{longenum}{enumerate}{6}
\setlist[longenum,1]{label=\arabic*.}
\setlist[longenum,2]{label=\arabic*.}
\setlist[longenum,3]{label=\arabic*.}
\setlist[longenum,4]{label=\arabic*.}
\setlist[longenum,5]{label=\arabic*.}
\setlist[longenum,6]{label=\arabic*.}

\colorlet{theoremcolor}{white}
\colorlet{definitioncolor}{white}
\colorlet{examplecolor}{white}
\colorlet{exercisecolor}{white}

\mdfdefinestyle{theoremframe}{
    linewidth=1pt,
    roundcorner=6pt,
    nobreak=true,
    leftmargin=0,
    innerleftmargin=9pt,
    rightmargin=0,
    innerrightmargin=9pt,
    }

\mdfdefinestyle{definitionframe}{
    linewidth=1pt,
    roundcorner=6pt,
    leftmargin=0,
    innerleftmargin=9pt,
    rightmargin=0,
    innerrightmargin=9pt,
    }

\mdfdefinestyle{exampleframe}{
    linewidth=1pt,
    roundcorner=6pt,
    leftmargin=0,
    innerleftmargin=9pt,
    rightmargin=0,
    innerrightmargin=9pt,
    }

\mdfdefinestyle{exerciseframe}{
    linewidth=1pt,
    leftmargin=0,
    innerleftmargin=9pt,
    rightmargin=0,
    innerrightmargin=9pt,
    }

\newtheoremstyle{plain}
  {-\topsep}		%space above
  {}			%space below
  {\normalfont}		%body font
  {}			%indent amount
  {\bfseries}		%theorem head font
  {.}			%punctuation after theorem head
  {.5em}		%space after theorem head
  {}			%theorem head spec

  \theoremstyle{plain}
  \newmdtheoremenv[style=theoremframe]{theorem}[equation]{Theorem}
  \newmdtheoremenv[style=theoremframe]{proposition}[equation]{Proposition}
  \newmdtheoremenv[style=theoremframe]{corollary}[equation]{Corollary}
  \newmdtheoremenv[style=theoremframe]{lemma}[equation]{Lemma}

  \theoremstyle{plain}
  \newmdtheoremenv[style=definitionframe]{definition}[equation]{Definition}
  \newmdtheoremenv[style=definitionframe]{informal}[equation]{Informal Definition}
  \newmdtheoremenv[style=definitionframe]{roughDef}[equation]{Rough Definition}
  \crefname{roughDef}{Definition}{Definitions}

  \newtheorem*{axiom*}{Axiom}

  \theoremstyle{remark}
  \newtheorem{remark}[equation]{Remark}

\makeatletter
\newcommand{\nolisttopbreak}{\nobreak\@afterheading}
\makeatother

\Newassociation{solution}{Sol}{solutions}

\newmdtheoremenv[style=exampleframe]{example}[equation]{Example}

\newmdtheoremenv[style=exerciseframe]{exc-inner}[equation]{Exercise}
\newenvironment{exercise}[1][]{
  \begin{exc-inner}[#1\ifthenelse{\equal{#1}{}}{\hyperlink{sol:\theequation}{Solution here}}{; \hyperlink{sol:\theequation}{solution here}}]\hypertarget{ex:\theequation}{}~
}{
  \end{exc-inner}
}
  \crefname{exercise}{Exercise}{Exercises}

\newcommand{\adj}[5][30pt]{%[size] Cat L, Left, Right, Cat R.
\begin{tikzcd}[ampersand replacement=\&, column sep=#1]
  #2\ar[r, shift left=7pt, "#3"]
  \ar[r, phantom, "\scriptstyle\Rightarrow"]\&
  #5\ar[l, shift left=7pt, "#4"]
\end{tikzcd}
}

\newcommand{\adjr}[5][30pt]{%[size] Cat R, Right, Left, Cat L.
\begin{tikzcd}[ampersand replacement=\&, column sep=#1]
  #2\ar[r, shift left=7pt, "#3"]\&
  #5\ar[l, shift left=7pt, "#4"]
  \ar[l, phantom, "\scriptstyle\Leftarrow"]
\end{tikzcd}
}

\newcommand{\slogan}[1]{\begin{center}\textit{#1}\end{center}}

\DeclareSymbolFont{stmry}{U}{stmry}{m}{n}
\DeclareMathSymbol\fatsemi\mathop{stmry}{"23}

\DeclareFontFamily{U}{mathx}{\hyphenchar\font45}
\DeclareFontShape{U}{mathx}{m}{n}{
      <5> <6> <7> <8> <9> <10>
      <10.95> <12> <14.4> <17.28> <20.74> <24.88>
      mathx10
      }{}
\DeclareSymbolFont{mathx}{U}{mathx}{m}{n}
\DeclareFontSubstitution{U}{mathx}{m}{n}
\DeclareMathAccent{\widecheck}{0}{mathx}{"71}

\renewcommand{\ss}{\subseteq}

\DeclarePairedDelimiter{\ihom}{[}{]}

\newcommand{\lchom}[2]{\genfrac{[}{]}{0pt}{}{#1}{#2}}
\newcommand{\rchom}[1]{\overset{#1}{\frown}}

\DeclareMathOperator{\Mor}{Mor}
\DeclareMathOperator{\dom}{dom}
\DeclareMathOperator{\cod}{cod}
\DeclareMathOperator{\idy}{idy}
\DeclareMathOperator{\comp}{com}
\DeclareMathOperator*{\colim}{colim}

\DeclareMathOperator{\Ob}{Ob}

\newcommand{\const}[1]{\texttt{#1}}%a constant, or named element of a set
\newcommand{\Set}[1]{\mathsf{#1}}%a named set
\newcommand{\ord}[1]{\mathsf{#1}}%an ordinal
\newcommand{\cat}[1]{\mathcal{#1}}%a generic category
\newcommand{\Cat}[1]{\mathbf{#1}}%a named category
\newcommand{\fun}[1]{\mathrm{#1}}%a function
\newcommand{\Fun}[1]{\mathrm{#1}}%a named functor
\newcommand{\com}[1]{\mathscr{#1}}%polynomial comonoid
\newcommand{\oper}[1]{\operatorname{#1}}

\newcommand{\id}{\mathrm{id}}
\newcommand{\then}{\mathbin{\fatsemi}}

\newcommand{\iso}{\cong}
\newcommand{\too}{\longrightarrow}
\newcommand{\tto}{\rightrightarrows}
\newcommand{\To}[1]{\xrightarrow{#1}}
\renewcommand{\Mapsto}[1]{\xmapsto{#1}}
\newcommand{\Tto}[3][13pt]{\begin{tikzcd}[sep=#1, cramped, ampersand replacement=\&, text height=1ex, text depth=.3ex]\ar[r, shift left=2pt, "#2"]\ar[r, shift right=2pt, "#3"']\&{}\end{tikzcd}}
\newcommand{\Too}[1]{\xrightarrow{\;\;#1\;\;}}
\newcommand{\from}{\leftarrow}
\newcommand{\From}[1]{\xleftarrow{#1}}

\newcommand{\surj}{\twoheadrightarrow}
\newcommand{\inj}{\hookrightarrow}

\newcommand{\tickar}{\xtickar{}}
\newcommand{\imp}{\Rightarrow}

\newcommand{\fromto}{\leftrightarrows}
\newcommand{\xtickar}[1]{\begin{tikzcd}[baseline=-0.5ex,cramped,sep=small,ampersand
replacement=\&]{}\ar[r,tick, "{#1}"]\&{}\end{tikzcd}}

\newcommand{\inv}{^{-1}}
\newcommand{\op}{^\tn{op}}

\newcommand{\tn}[1]{\textnormal{#1}}
\newcommand{\ol}[1]{\overline{#1}}

\newcommand{\wh}[1]{\widehat{#1}}
\newcommand{\wc}[1]{\widecheck{#1}}

\newcommand{\LMO}[2][over]{\LMOO[#1]{\phantom{\cdot}}{#2}}
\newcommand{\LMOO}[3][over]{\ifthenelse{\equal{#1}{over}}
  {
    \begin{tikzpicture}[inner sep=0pt, solid] \node[draw, circle, label={[label distance
        = 2pt]90:{$#3$}}] {$#2$}; \end{tikzpicture}
  }{
    \begin{tikzpicture}[inner sep=0pt, solid] \node[draw, circle, label={[label distance
        = 2pt]270:{$#3$}}] {$#2$}; \end{tikzpicture}
  }
}

\newcommand{\nn}{\mathbb{N}}

\newcommand{\qq}{\mathbb{Q}}
\newcommand{\zz}{\mathbb{Z}}
\newcommand{\rr}{\mathbb{R}}

\newcommand{\smset}{\Cat{Set}}
\newcommand{\smcat}{\Cat{Cat}}
\newcommand{\catsharp}{\smcat^\sharp}
\newcommand{\comon}{\Cat{Comon}}

\newcommand{\List}{\Fun{List}}

\newcommand{\boxCD}[3][black] % [border color] {interior color} {math text}
{\fcolorbox{#1}{#2}{\begin{varwidth}{\textwidth}\centering #3\end{varwidth}}}

\newenvironment{dedication}
     {\begin{flushright}\itshape}
     {\end{flushright}}

\newcommand{\yon}{\mathcal{y}}
\newcommand{\poly}[1][]{#1\Cat{Poly}}
\newcommand{\0}{\textsf{0}}
\newcommand{\1}{\textsf{1}}
\newcommand{\2}{\textsf{2}}
\newcommand{\3}{\textsf{3}}
\newcommand{\4}{\textsf{4}}
\newcommand{\5}{\textsf{5}}
\newcommand{\6}{\textsf{6}}
\newcommand{\7}{\textsf{7}}
\newcommand{\8}{\textsf{8}}
\newcommand{\9}{\textsf{9}}

\definecolor{dgreen}{rgb}{0.0, 0.5, 0.3}
\definecolor{dyellow}{rgb}{8.0, 0.74, 0}

\newcommand{\at}{\nodepart{two}}

\newcommand{\qand}{\quad\text{and}\quad}
\newcommand{\qqand}{\qquad\text{and}\qquad}
\newcommand{\qqor}{\qquad\text{or}\qquad}

\newcommand{\lst}{\Set{List}}

\newcommand{\tr}{\Set{tree}}
\newcommand{\vtx}{\Set{vtx}}

\newcommand{\lensput}{\fun{put}}
\newcommand{\lensget}{\fun{get}}

\newcommand{\src}{\fun{src}}
\newcommand{\tgt}{\fun{tgt}}

\newcommand{\elts}{\textstyle{\int}}

\newcommand{\cofree}[1]{\com{T}_{#1}}

\newcommand{\bimod}[2]{~_{#1}\Cat{Mod}_{#2}}

\newcommand{\car}[1]{\mathfrak{#1}}

\newcommand{\cof}{\nrightarrow}

\newcommand{\tri}{\mathbin{\triangleleft}}
\newcommand{\tripow}[1]{^{\tri\,#1}}

\newcommand{\bul}[1][black]{{\color{#1}\ensuremath{\bullet}}}

 \makeatletter
\newcommand*{\bdiv}{%
  \nonscript\mskip-\medmuskip\mkern5mu%
  \mathbin{\operator@font div}\penalty900\mkern5mu%
  \nonscript\mskip-\medmuskip
}
\makeatother

\newcounter{notecounter}

\newcommand{\Red}{{\color{my-red}\const{red}}}
\newcommand{\Blue}{{\color{my-blue}\const{blue}}}
\newcommand{\Black}{\const{black}}
\newcommand{\Yellow}{{\color{my-yellow}\const{yellow}}}

\definecolor{my-blue}{HTML}{005AB5}
\definecolor{my-red}{HTML}{DC3220}
\definecolor{my-lavender}{HTML}{785EF0}
\definecolor{my-magenta}{HTML}{DC267F}
\definecolor{my-yellow}{HTML}{FFB000}

\newcommand{\rectangle}{{%
    \Large
  \ooalign{$\sqsubset\mkern3mu$\cr$\mkern3mu\sqsupset$\cr}%
}}
\newcommand{\sqiint}{{\rectangle}\kern-1.4em{\iint}}

\makeatletter
\DeclareRobustCommand{\rvdots}{%
  \vbox{
    \baselineskip4\p@\lineskiplimit\z@
    \kern-\p@
    \hbox{.}\hbox{.}\hbox{.}
  }}
\makeatother

\allowdisplaybreaks
\setsecnumdepth{subsection}
\settocdepth{subsection}
\begin{document}

\frontmatter

\title{\titlefont Polynomial Functors:\\\medskip
A Mathematical Theory of Interaction}

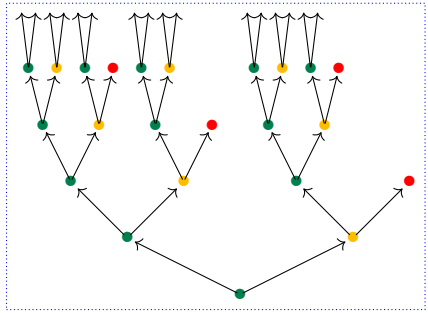
\author{
\LARGE Nelson Niu
\and
\LARGE David I. Spivak
\normalsize}

\posttitle{
\vspace{.8in}
\normalsize
\[
\coverpic
\]
  \vspace{.5in}
  \endgroup
}
\date{\vfill Last updated: \today \\ Source: \url{https://github.com/ToposInstitute/poly}}

\maketitle

\thispagestyle{empty}
\clearpage
This page intentionally left blank.
\clearpage
\clearpage

\mbox{}
\vspace{2in}
\begin{center}
\begin{tabular}{lll}
\LARGE Nelson Niu&~\hspace{.5in}~&\LARGE  David I. Spivak\\
\large University of Washington&&\large Topos Institute\\
\large Seattle, WA&&\large Berkeley, CA
\end{tabular}
\end{center}
\clearpage

\begin{dedication}
To Andr\'e Joyal
\smallskip

---D.S.
\bigskip

To my graduate cohort at UW
\smallskip

---N.N.
\end{dedication}

%------------ Chapter ------------%
\chapter*{Preface}\label{chapter.0}
\addcontentsline{toc}{chapter}{Preface}

\begin{quote}
	The proposal is also intended to [serve] equally as a foundation for the academic, intellectual, and technological, on the one hand, and for the curious, the moral, the erotic, the political, the artistic, and the sheerly obstreperous, on the other.\\
\mbox{}\hfill ---Brian Cantwell Smith\\
\mbox{}\hfill \emph{On the Origin of Objects}
\end{quote}

%\begin{quote}
%For me, though, it is difficult to resist the idea that space-time is not essentially different from matter, which we understand more deeply. If so, it will consist of vast numbers of identical units---``particles of space''---each in contact with a few neighbors, exchanging messages, joining and breaking apart, giving birth and passing away.\\
%\mbox{}\hfill ---Frank Wilczek\\
%\mbox{}\hfill \emph{Fundamentals}
%\end{quote}

\begin{quote}
And that is the way to get the greatest possible variety, but with all the order there could be; i.e. it is the way to get as much perfection as there could be.\\
\mbox{}\hfill ---Gottfried Wilhelm Leibniz\\
\mbox{}\hfill \emph{Monadology}, 58.
\end{quote}

\noindent During the Fifth International Conference on Applied Category Theory in 2022, at least twelve of the fifty-nine presentations and two of the ten posters referenced the category of polynomial functors and dependent lenses or its close cousins (categories of optics and Dialectica categories) and the way they model diverse forms of interactive behavior.
At the same time, all that is needed to grasp the construction of this category---called $\poly$ for short---is an understanding of mathematical sets and functions.
There is no need for the theory and applications of polynomial functors to remain the stuff of technical papers; $\poly$ is far too versatile, too full of potential, to be kept out of reach.

\index{applied category theory}\index{lens}\index{dialectica category}

Informally, a polynomial functor is a collection of elements we call \emph{positions} and, for each position, a collection of elements we call \emph{directions}.
There is then a natural notion of a morphism between polynomial functors that sends positions forward and directions backward, modeling two-way communication.
From these basic components, category theory allows us to construct an immense array of mathematical gadgets that model a diverse range of interactive processes.
In this book, we will establish the theory of polynomial functors and categorical constructions on them while exploring how they model interaction.

\section[Purpose and prerequisites]{Purpose and prerequisites%
  \sectionmark{}}
\sectionmark{}

A categorical theory of general interaction must be interdisciplinary by its very nature.
Already, drafts of this text have been read by everyone from algebraic geometers to neuroscientists and AI developers.
We hope to extend our reach ever further, to bring together thinkers and tinkerers from a diverse array of backgrounds under a common language by which to study interactive systems categorically.
In short---we know about $\poly$; you know about other things; but only our collective knowledge can reveal how $\poly$ could be applied to those other things.

\index{interdisciplinary}

As such, we have strived to write a friendly and accessible expository text that can serve as a stepping stone toward further investigations into polynomial functors.
We include exercises and, crucially, solutions to guide the learning process; we draw extensive analogies to provide motivation and develop intuition; we pose examples whenever necessary.
Proofs may bear far more detail than you would find in a research paper, but not so much detail that it would clutter the key ideas.
A few critical proofs are even argued through pictures, yet we contend that they are no less rigorous than the clouds of notation whose places they take.

On the other hand, there is some deep mathematical substance to the work we will discuss, drawing from the well-established theory of categories.
Although you will find, for example, a complete proof of the Yoneda lemma within these pages, we don't intend to build up everything from scratch.
There are plenty of excellent resources for learning category theory out there, catering to a variety of needs, without adding our own to the mix when our primary goal is to introduce $\poly$.
So for the sake of contributing only what is genuinely helpful, we assume a certain level of mathematical background.
You are ready to read this book if you can define the following fundamental concepts from category theory, and give examples of each:
\begin{itemize}
    \item categories,
    \item functors,
    \item natural transformations,
    \item (co)limits,
    \item adjunctions, and
    \item (symmetric) monoidal categories.
\end{itemize}
We will additionally assume a passing familiarity with the language of graph-theoretic trees (e.g.\ vertices, roots, leaves, paths).

That said, with a little investment on your part, you could very well use this book as a way to teach yourself some category theory.
If you have ever tried to learn category theory, only to become lost in abstraction or otherwise overwhelmed by seemingly endless lists of examples from foreign fields, perhaps you will benefit from a focused case study of one particularly fruitful category.
If you encounter terms or ideas that you would like to learn more about, we encourage you to look them up elsewhere, and you may find yourself spending a pleasant afternoon doing a deep dive into a new definition or theorem.
Then come back when you're ready---we'll be here.

\section[Outline]{Outline%
  \sectionmark{}}
\sectionmark{}

This book is designed to be read linearly.
\cref{part.poly} introduces the category $\poly$ and illustrates how it models interaction protocols; while \cref{part.comon} highlights a crucial operation on $\poly$, the composition product, which upgrades the theory so that it properly captures the time evolution of a dynamical system.

\bigskip

\cref{part.poly} (The category of polynomial functors) consists of:
\begin{itemize}
  \item \cref{ch.poly.rep-sets} (Representable functors from the category of sets), in which we review category-theoretic constructions on sets and the Yoneda lemma;
  \item \cref{ch.poly.obj} (Polynomial functors), in which we introduce our objects of study and present several perspectives from which to view them;
  \item \cref{ch.poly.cat} (The category of polynomial functors), in which we define $\poly$ by introducing its morphisms and demonstrate many ways to work with them, including how they model interaction protocols;
  \item \cref{ch.poly.dyn_sys} (Dynamical systems as dependent lenses), in which we use a specific class of morphisms in $\poly$ to model discrete-time dynamical systems; and
  \item \cref{ch.poly.bonus} (More categorical properties of polynomials), in which we describe a smorgasbord of additional category-theoretic structures on $\poly$.
\end{itemize}
\bigskip

\cref{part.comon} (A different category of categories) consists of:
\begin{itemize}
  \item \cref{ch.comon.comp} (The composition product), in which we examine a monoidal structure on $\poly$ given by substituting one polynomial into another;
  \item \cref{ch.comon.sharp} (Polynomial comonoids and retrofunctors), in which we show that the category of comonoids in $\poly$ with respect to the composition product is equivalent to a category of small categories we call $\smcat^\sharp$ whose morphisms are not functors;
  \item \cref{ch.comon.cofree} (Categorical properties of polynomial comonoids), in which we study the structure and utility of $\smcat^\sharp$; and
  \item \cref{ch.comon.vistas} (New horizons), in which we list open questions.
\end{itemize}

\section[Choices and conventions]{Choices and conventions%
  \sectionmark{}}
\sectionmark{}
\label{sec.choices}

Throughout this book, we have chosen to focus on polynomial functors of a single variable on the category of sets.
The motivation for this seemingly narrow scope is twofold: to keep matters as concrete and intuitive as possible, with easy access to elements that we can work with directly; and to demonstrate the immense versatility of even this small corner of the theory. Furthermore, the subject of multi-variate polynomials arises by considering what are called comonoids and comodules in $\poly$ \cite{spivak2023functorial}.

\index{polynomial functor!multivariate}\index{multivariate polynomial|see{polynomial functor, multivariate}}

Below is a list of conventions we adopt; while it is not comprehensive, any unusual choices are justified within the text, often as a footnote.

The natural numbers include $0$, so $\nn\coloneqq\{0,1,2,\ldots\}$. Throughout this book, when referring to finite sets, we will adopt the following convention: $\0\coloneqq\{\}=\varnothing,$ $\1\coloneqq\{1\},$ $\2\coloneqq\{1,2\},$ $\3\coloneqq\{1,2,3\}$, and so on, with $\ord{n}\coloneqq\{1,\ldots,n\}$, an $n$-element set, for each natural number $n$.
For example, in standard font, $5$ represents the usual natural number, while in sans serif font, $\5$ represents the $5$-element set $\{1,2,3,4,5\}$.
When the same variable name appears in both italicized and sans serif fonts, the italicized variable denotes a natural number and the sans serif variable denotes the corresponding set; for example, if we state that $m\in\nn$, then we also understand $\ord{m}$ to mean the set $\{1,\ldots,m\}$.
\index{set!finite}

\index{natural numbers, $\nn$}\index{set!ordinal}

The names of categories will be capitalized.
We will mostly ignore size issues, but roughly speaking small categories will be written in script (e.g.$\ \cat{C}, \cat{D}$), while large categories (usually, but not always, named) will be written in bold (e.g.\ $\poly, \Cat{C}$).
We use $\smset$ to denote the category of (small) sets and functions and $\smcat$ to denote the category of (small) categories and functors.
We use exponential notation $\cat{D}^{\cat{C}}$ to denote the category of functors $\cat{C}\to\cat{D}$ and natural transformations.
\index{set!category of sets}\index{functor!category of functors}

We write either $c\in\Ob\cat{C}$ or $c\in\cat{C}$ to denote an object $c$ of a category $\cat{C}$.
We use $\sum$ rather than $\coprod$ to denote coproducts.
We denote the collection of morphisms $f\colon c\to d$ in a category $\cat{C}$ by using the name of the category itself, followed by the ordered pair of objects: $\cat{C}(c,d)$.
We denote the domain of $f$ by $\dom f$ and the codomain of $f$ by $\cod f$.
We use $\coloneqq$ for definitions and temporary assignments, as opposed to $=$ for identifications that can be observed and proven.
We use $\iso$ to indicate an isomorphism of objects and $=$ to indicate an equality of objects, although the choice of the former does not preclude the possibility of the latter, nor does the latter necessarily imply any significance beyond an arbitrary selection that has been made.
We will freely use the definite article ``the'' to refer to objects that are unique only up to unique isomorphism.

\index{coproduct}\index{product}\index{isomorphism}\index{equality!use of}

We list nullary operations before binary ones: for example, we denote a monoidal category $\cat{C}$ with monoidal unit $I$ and monoidal product $\odot$ by $(\cat{C},I,\odot)$, or say that $(I,\odot)$ is a monoidal structure on $\cat{C}$.

\section[Past, present, and future]{Past, present, and future%
  \sectionmark{}}
\sectionmark{}

The idea for this book began in 2020, originally as part of a joint work with David Jaz Myers on using categories to model dynamical systems.
It soon became clear, however, that our writing and his---while intimately related---would be better off as separate volumes.
His book is nonetheless an excellent companion to ours: see \cite{jaz}.

In the summer of 2021, we taught a course on a draft of this book that was livestreamed from the Topos Institute.
Lecture recordings are freely available at \url{https://topos.site/poly-course/}.
A follow-up workshop with additional recorded lectures and write-ups was held in early 2024; see \url{https://topos.site/events/poly-at-work/}.

The theory and application of polynomial functors comprise an active area of research.
We have laid the foundations here, but work is still ongoing.
Even while writing this, we discovered new results and uses for polynomial functors, which only goes to show how bountiful $\poly$ can be in its rewards---but of course, we had to cut things off somewhere.
We say this in the hope that you will keep the following in mind: where this book ends, the story will have just begun.

\section[Acknowledgments]{Acknowledgments%
  \sectionmark{}}
\sectionmark{}

Special thanks to David Jaz Myers: a brilliant colleague, a wonderful conversation partner, a congenial housemate, a superb chef, and an all-around good guy.

Thanks go to John Baez, Eric Bond, Spencer Breiner, Kris Brown, Matteo Capucci, Valeria de Paiva, Joseph Dorta, Brendan Fong, Richard Garner, Bruno Gavranovi\'c, Neil Ghani, Ben Goertzel, Tim Hosgood, Samantha Jarvis, Max Lieblich, Shaowei Lin, Owen Lynch, Joachim Kock, J\'er\'emie Koenig, Sophie Libkind, Joshua Meyers, Dominic Orchard, Nathaniel Osgood, Evan Patterson, Brandon Shapiro, Juliet Szatko, Tish Tanski, Todd Trimble, Adam Vandervorst, Jonathan Weinberger, and Christian Williams.

This material is based upon work supported by the AFOSR under award numbers FA9550-20-1-0348 and FA9550-23-1-0376.

\clearpage
\tableofcontents*
\clearpage

\mainmatter

\setcounter{chapter}{0}%Just finished 0.

%---------------- Part ----------------%
\part{The category of polynomial functors}\label{part.poly}

\Opensolutionfile{solutions}[solution-file1]

%------------ Chapter ------------%
\chapter{Representable functors from the category of sets}
\chaptermark{Representable functors} \label{ch.poly.rep-sets}
%\dnote{Get rid of ``solution here''}

In this chapter, we lay the categorical groundwork needed to define our category of interest, the category of polynomial functors.
We begin by examining a special kind of polynomial functor that you may already be familiar with---representable functors from the category $\smset$ of sets and functions.%
\index{functor!representable|(}
\index{functor!polynomial|see{polynomial functor}}
\index{polynomial functor}
\index{functor!set-valued}
We highlight the role these representable functors play in what is arguably the fundamental theorem of category theory, the Yoneda lemma.
We will also discuss sums and products of sets and of set-valued functors, which we will need to construct our polynomial functors.

%-------- Section --------%
\section[Representable functors and the Yoneda lemma]{Representable functors and the Yoneda lemma%
  \sectionmark{Representable functors \& the Yoneda lemma}}
\sectionmark{Representable functors \& the Yoneda lemma} \label{sec.poly.rep-sets.yon}

\index{Yoneda lemma}

Representable functors---special functors to the category of sets---provide the foundation for the category $\poly$.
While much of the following theory applies to representable functors from any category, we need only focus on representable functors $\smset\to\smset$.

\begin{definition}[Representable functor] \label{def.representable}
    For a set $S$, we denote by $\yon^S\colon\smset\to\smset$ the functor that sends each set $X$ to the set $X^S=\smset(S,X)$ and each function $h\colon X\to Y$ to the function $h^S\colon X^S\to Y^S$ that sends $g\colon S\to X$ to $g\then h\colon S\to Y$.\tablefootnote{Throughout this text, given morphisms $f \colon A \to B$ and $g \colon B \to C$ in a category, we will denote their composite morphism $A \to C$ interchangeably as $f \then g$ or $g \circ f$ (or even $gf$), whichever is more convenient.}

    We call a functor (isomorphic to one) of this form a \emph{representable functor}, or a \emph{representable}.
    In particular, we call $\yon^S$ the functor \emph{represented by} $S$, and we call $S$ the \emph{representing set of $\yon^S$}.
    As $\yon^S$ denotes raising a variable to the power of $S$, we will also refer to representables as \emph{pure powers}.%
    \index{pure power|see{functor, representable}}
\end{definition}

The symbol $\yon$ stands for Yoneda, for reasons we will explain in \cref{lemma.yoneda} and \cref{exc.finish_proof_yoneda} \cref{exc.finish_proof_yoneda.embed}.

Throughout this book, we will use the notation $\0\coloneqq\{\}=\varnothing,$ $\1\coloneqq\{1\},$ $\2\coloneqq\{1,2\},$ $\3\coloneqq\{1,2,3\}$, and so on, with $\ord{n}\coloneqq\{1,\ldots,n\}$
\begin{example}
    The functor that sends each set $X$ to $X\times X$ and each function $h\colon X\to Y$ to $(h\times h)\colon (X\times X)\to(Y\times Y)$ is representable.
    After all, $X\times X \iso X^\2$, so this functor is the pure power $\yon^\2$.
\end{example}

\begin{exercise}\label{exc.representable_fun}\index{functor!representable}\index{representable polynomial|seealso{functor!representable}}
    For each of the following functors $\smset\to\smset$, say if it is representable or not; if it is, give the set that represents it.
    \begin{enumerate}
        \item The identity functor $X\mapsto X$, which sends each function to itself.
        \item The constant functor $X\mapsto\2$, which sends every function to the identity on $\2$.
        \item The constant functor $X\mapsto\1$, which sends every function to the identity on $\1$.
        \item The constant functor $X\mapsto\0$, which sends every function to the identity on $\0$.
        \item A functor $X\mapsto X^\nn$.
        If it could be representable, where should it send each function?
        \item A functor $X\mapsto \2^X$.
        If it could be representable, where should it send each function? \qedhere
    \end{enumerate}
\index{functor!identity}\index{functor!constant}

    \begin{solution}
        \begin{enumerate}
            \item The identity functor $X\mapsto X$ is represented by the set $\1$: a function $\1 \to X$ can be identified with an element of $X$, so $\smset(\1,X)\iso X$.
            Alternatively, note that $X^\1 \iso X$.
            \item \label{sol.representable_fun.2} The constant functor $X\mapsto\2$ is not representable: it sends $\1$ to $\2$, but $\1^S \iso \1 \not\iso \2$ for any set $S$.
            \item The constant functor $X\mapsto\1$ is represented by $S=\0$: there is exactly one function $\0 \to X$, so $\smset(\0,X) \iso \1$.
            Alternatively, note that $X^0 \iso \1$.
            \item The constant functor $X\mapsto\0$ is not representable for the same reason as in \cref{sol.representable_fun.2}.
            \item The functor $\yon^\nn$ that sends $X\mapsto X^\nn$ is represented by $\nn$, by definition.
            It should send each function $h \colon X \to Y$ to the function $h^\nn \colon X^\nn \to Y^\nn$ that sends each $g \colon \nn \to X$ to $g \then h \colon \nn \to Y$.
            \item No $\smset \to \smset$ functor $X\mapsto \2^X$ is representable, for the same reason as in \cref{sol.representable_fun.2}.
            (There \emph{is}, however, a functor $\smset\op \to \smset$ sending $X \mapsto 2^X$ that is understood to be representable in a more general sense.)
        \end{enumerate}
    \end{solution}
\end{exercise}

Now that we have introduced representable functors, we study the maps between them.
As representables are functors, the maps between them are natural transformations.
\index{natural transformation!between representables}

\begin{proposition}\label{prop.representable_nt}
    For any function $f\colon R\to S$, there is an induced natural transformation $\yon^f\colon\yon^S\to \yon^R$; on any set $X$ its $X$-component $X^f\colon X^S\to X^R$ is given by sending $g\colon S\to X$ to $f\then g\colon R\to X$.
\end{proposition}

\begin{proof}
    See \cref{exc.representable_nt}.
\end{proof}

\begin{exercise} \label{exc.representable_nt}
    To prove \cref{prop.representable_nt}, show that for any function $f\colon R\to S$, the given construction $\yon^f\colon\yon^S\to\yon^R$ really is a natural transformation.
    That is, for any function $h\colon X\to Y$, show that the following naturality square commutes:
    \begin{equation} \label{diag.yon_embed_nat}
        \begin{tikzcd}%[bottom base]
            X^S\ar[r, "h^S"]\ar[d, "X^f"']&
            Y^S\ar[d, "Y^f"]\\
            X^R\ar[r, "h^R"']&
            Y^R\ar[ul, phantom, "?"]
        \end{tikzcd}
    \end{equation}
    \qedhere

    \begin{solution}
        To show that \eqref{diag.yon_embed_nat} commutes, we note that by the construction of the components of $\yon^f$ in the statement of \cref{prop.representable_nt}, both vertical maps in the diagram compose functions from $S$ with $f \colon R \to S$ on the left; and by \cref{def.representable}, both horizontal maps compose functions to $X$ with $h \colon X \to Y$ on the right.
        So by the associativity of composition, the diagram commutes: $(f\then g)\then h=f\then(g\then h)$ for all $g\colon S\to X$.
    \end{solution}
\end{exercise}\index{associativity}

\begin{exercise} \label{exc.representable_nt_components}
    Let $X$ be an arbitrary set. For each of the following sets $R,S$ and functions $f\colon R\to S$, describe the $X$-component $X^f\colon X^S\to X^R$ of the natural transformation $\yon^f\colon\yon^S\to\yon^R$.
    \begin{enumerate}
        \item \label{exc.representable_nt_components.id} $R=\5$, $S=\5$, $f=\id_\5$. (You should describe the function $X^{\id_\5}\colon X^\5\to X^\5$.)
        \item $R=\2$, $S=\1$, $f$ is the unique function.
        \item $R=\1$, $S=\2$, $f(1)=1$.
        \item $R=\1$, $S=\2$, $f(1)=2$.
        \item $R=\0$, $S=\5$, $f$ is the unique function.
        \item $R=\nn$, $S=\nn$, $f(n)=n+1$.
        \qedhere
    \end{enumerate}

    \begin{solution}
        In each case, given $f \colon R \to S$, we can find the $X$-component $X^f \colon X^S \to X^R$ of the natural transformation $\yon^f\colon\yon^S\to\yon^R$ by applying \cref{prop.representable_nt}, which says that $X^f$ sends each $g \colon S \to X$ to $f \then g \colon R \to X$.
        \begin{enumerate}
            \item Here $X^{\id_5}\colon X^\5\to X^\5$ is the identity function.
            \item If $f\colon\2\to\1$ is the unique function, then $X^f\colon X^\1\to X^\2$ sends each $g \in X$ (i.e.\ function $g \colon \1 \to X$) to the function that maps both elements of $\2$ to $g$.
            We can think of $X^f$ as the diagonal $X \to X \times X$.
            \item If $f\colon\1\to\2$ sends $1\mapsto1$, then $X^f\colon X^\2\to X^\1$ sends each $g \colon \2 \to X$ to $g(1)$, viewed as a function $\1 \to X$.
            We can think of $X^f$ as the left projection $X \times X \to X$.
            \item If $f\colon\1\to\2$ sends $1\mapsto2$, then $X^f\colon X^\2\to X^\1$ sends each $g \colon \2 \to X$ to $g(2)$, viewed as a function $\1 \to X$.
            We can think of $X^f$ as the right projection $X \times X \to X$.
            \item Here $X^f\colon X^\5\to X^\0\iso\1$ is the unique function.
            \item If $f\colon\nn\to\nn$ sends $n\mapsto n+1$, then $X^f\colon X^\nn\to X^\nn$ sends each $g \colon \nn \to X$ to the function $h \colon \nn \to X$ defined by $h(n)\coloneqq g(n+1)$ for all $n \in \nn$.
            We can think of $X^f$ as removing the first term of an infinite sequence of elements $(g(0),g(1),g(2),\ldots)$ of $X$ to obtain a new sequence $(g(1),g(2),g(3),\ldots)$.
        \end{enumerate}
    \end{solution}
\end{exercise}

These representable functors and natural transformations live in the larger category $\smset^\smset$, whose objects are functors $\smset\to\smset$ and whose morphisms are the natural transformations between them.

\begin{exercise} \label{exc.representable_nt_functorial}
    Show that the construction $f\mapsto\yon^f$ from \cref{prop.representable_nt} defines a functor
    \begin{equation} \label{eqn.yoneda_embedding}
        \yon^-\colon\smset\op\to\smset^\smset
    \end{equation}
    by verifying functoriality, as follows.
    \begin{enumerate}
        \item Show that for any set $S$, the natural transformation $\yon^{\id_S}\colon\yon^S\to\yon^S$ is the identity.
        \item Show that for functions $f\colon R\to S$ and $g\colon S\to T$, we have $\yon^g\then\yon^f=\yon^{f\then g}$. \qedhere
    \end{enumerate}

    \begin{solution}
        \begin{enumerate}
            \item The fact that $\yon^{\id_S}\colon\yon^S\to\yon^S$ is the identity is just a generalization of \cref{exc.representable_nt_components} \cref{exc.representable_nt_components.id}.
            For any set $X$, the $X$-component $X^{\id_S} \colon X^S \to X^S$ of $\yon^{\id_S}$ sends each $h \colon S \to X$ to $\id_S \then h = h$, so $X^{\id_S}$ is the identity natural transformation on $X^S$.
            Hence $\yon^{\id_S}$ is the identity on $\yon^S$.
            \item Fix $f \colon R \to S$ and $g \colon S \to T$; we wish to show that $\yon^g \then \yon^f = \yon^{f \then g}$.
            It suffices to show componentwise that $X^g \then X^f = X^{f \then g}$ for every set $X$.
            Indeed, $X^g$ sends each $h \colon T \to X$ to $g \then h$; then $X^f$ sends $g \then h$ to $f \then g \then h = X^{f \then g}(h)$.
        \end{enumerate}
    \end{solution}
\end{exercise}

We now have all the ingredients we need to state and prove the Yoneda lemma on the category of sets.
\index{Yoneda lemma}

\begin{lemma}[Yoneda lemma]\label{lemma.yoneda}
    Given a functor $F\colon\smset\to\smset$ and a set $S$, there is an isomorphism
    \begin{equation}\label{eqn.yoneda}
        F(S)\iso\smset^\smset(\yon^S,F)
    \end{equation}
    where the right hand side is the set of natural transformations $\yon^S\to F$.
    Moreover, \eqref{eqn.yoneda} is natural in both $S$ and $F$.
\end{lemma}
\begin{proof}[Proof]
    Given a natural transformation $m\colon\yon^S\to F$, consider its $S$-component $m_S\colon S^S\to F(S)$.
    Applying this function to $\id_S\in S^S$ yields an element $m_S(\id_S)\in F(S)$.

    Conversely, given an element $a\in F(S)$, there is a natural transformation we denote by $m^a\colon\yon^S\to F$ whose $X$-component is the function $X^S\to F(X)$ that sends $g\colon S\to X$ to $F(g)(a)$.
    In \cref{exc.finish_proof_yoneda} we ask you to show that this is indeed natural in $X$; that these two constructions, $m\mapsto m_S(\id_S)$ and $a\mapsto m^a$, are mutually inverse; and that the resulting isomorphism is natural.
\end{proof}

\index{natural transformation!and Yoneda embedding}

\begin{exercise}\label{exc.finish_proof_yoneda}
    In this exercise, we fill in the details of the preceding proof.
    \begin{enumerate}
        \item Show that for any $a\in F(S)$, the maps $X^S\to F(X)$ defined in the proof of \cref{lemma.yoneda} are natural in $X$.
        \item Show that the two mappings given in the proof of \cref{lemma.yoneda} are mutually inverse, thus defining the isomorphism \eqref{eqn.yoneda}.
        \item Show that \eqref{eqn.yoneda} as defined is natural in $F$.
        \item Show that \eqref{eqn.yoneda} as defined is natural in $S$.
        \item \label{exc.finish_proof_yoneda.embed} As a corollary of \cref{lemma.yoneda}, show that the functor $\yon^-\colon\smset\op\to\smset^\smset$ from \eqref{eqn.yoneda_embedding} is fully faithful---in particular, there is an isomorphism $S^T\iso \smset^\smset(\yon^S,\yon^T)$ for sets $S,T$.
        For this reason, we call $\yon^-$ the \emph{Yoneda embedding}.
        \qedhere
    \end{enumerate}

    \begin{solution}
        \begin{enumerate}
            \item To check that $X^S \to F(X)$ is natural in $X$, we verify that the naturality square
            \[
            \begin{tikzcd}[ampersand replacement=\&]
                X^S\ar[r, "h^S"]\ar[d, "F(-)(a)"']\&
                Y^S\ar[d, "F(-)(a)"]\\
                F(X)\ar[r, "F(h)"']\&
                F(Y)
            \end{tikzcd}
            \]
            commutes for all $h \colon X \to Y$.
            The top map $h^S$ sends any $g \colon S \to X$ to $g \then h$ (\cref{def.representable}), which is then sent to $F(g \then h)(a)$ by the right map.
            Meanwhile, the left map sends $g$ to $F(g)(a)$, which is then sent to $F(h)(F(g)(a))$ by the bottom map.
            So by the functoriality of $F$, the square commutes.

            \item We show that the maps $m\mapsto m_S(\id_S)$ and $a\mapsto m^a$ defined in the proof of \cref{lemma.yoneda} are mutually inverse.
            First, we show that for any natural transformation $m \colon \yon^S \to F$, we have $m^{m_S(\id_S)} = m$.
            Given a set $X$, the $X$-component of $m^{m_S(\id_S)}$ sends each $g \colon S \to X$ to $F(g)(m_S(\id_S))$; it suffices to show that this is also where the $X$-component of $m$ sends $g$.
            Indeed, by the naturality of $m$, the square
            \[
            \begin{tikzcd}[ampersand replacement=\&]
                S^S\ar[r, "g^S"]\ar[d, "m_S"']\&
                X^S\ar[d, "m_X"]\\
                F(S)\ar[r, "F(g)"']\&
                F(X)
            \end{tikzcd}
            \]
            commutes, so in particular, following $\id_S\in S^S$ around this diagram, we have
            \begin{equation} \label{eqn.finish_proof_yoneda}
                F(g)(m_S(\id_S)) = m_X(g^S(\id_S)) = m_X(\id_S \then g) = m_X(g).
            \end{equation}
            In the other direction, we show that for any $a \in F(S)$, we have $m^a_S(\id_S) = a$: by construction, $m^a_S \colon S^S \to F(S)$ sends $\id_S$ to $F(\id_S)(a) = a$.

            \item It suffices to show that given functors $F, G\colon\smset\to\smset$ and a natural transformation $\alpha \colon F \to G$, the naturality square
            \[
            \begin{tikzcd}[ampersand replacement=\&]
                \smset^\smset(\yon^S,F)\ar[d, "- \then \alpha"']\ar[r, "\sim"]\&
                F(S)\ar[d, "\alpha_S"]\\
                \smset^\smset(\yon^S,G)\ar[r, "\sim"]\&
                G(S)
            \end{tikzcd}
            \]
            commutes.
            The top map sends any $m \colon \yon^S \to F$ to $m_S(\id_S)$, which in turn is sent by the right map to $\alpha_S(m_S(\id_S)) = (m \then \alpha)_S(\id_S)$.
            This is also where the bottom map sends $m \then \alpha$, so the square commutes.

            \item It suffices to show that given a function $g \colon S \to X$, the naturality square
            \[
            \begin{tikzcd}[ampersand replacement=\&]
                \smset^\smset(\yon^S,F)\ar[d, "\yon^g \then -"']\ar[r, "\sim"]\&
                F(S)\ar[d, "F(g)"]\\
                \smset^\smset(\yon^X,F)\ar[r, "\sim"]\&
                F(X)
            \end{tikzcd}
            \]
            commutes.
            The left map sends any $m \colon \yon^S \to F$ to $\yon^g \then m$, which is sent by the bottom map to $(\yon^g \then m)_X(\id_X) = m_X(X^g(\id_X)) = m_X(g \then \id_X) = m_X(g)$.
            Meanwhile, the top map sends $m$ to $m_S(\id_S)$, which is sent by the right map to $F(g)(m_S(\id_S))$.
            So the square commutes by \eqref{eqn.finish_proof_yoneda}.

            \item To show that $\smset^\smset(\yon^S, \yon^T) \iso S^T$, just take $F\coloneqq\yon^T$ in \cref{lemma.yoneda} so that $F(S)\iso S^T$.
        \end{enumerate}
    \end{solution}
\end{exercise}

\index{functor!representable|)}

How will we go from these representable functors to polynomial ones?%
\index{polynomial functor}
Recall that, in algebra, a polynomial is just a sum of pure powers.
So we will define a \emph{polynomial functor} $\smset\to\smset$ to be a sum of pure power functors---that is, the representable functors $\yon^A$ for each set $A$ we just introduced.%
\footnote{This analogy isn't perfect: in algebra, polynomials are generally finite sums of pure powers, whereas our polynomial functors may be infinite sums of representables.
However, we are not the first to use the term ``polynomial'' this way, and the name stuck.}

All of our polynomials will be in one variable, $\yon$.
Every other letter or number that shows up in our notation for a polynomial will denote a set.
For example, in the polynomial
\begin{equation} \label{eqn.biz_poly}
    \rr\yon^\zz+\3\yon^{\3}+\yon^A+\sum_{i\in I}\yon^{R_i+Q_i^\2},
\end{equation}
$\rr$ denotes the set of real numbers, $\zz$ denotes the set of integers, $\2$ and $\3$ respectively denote the sets $\{1,2\}$ and $\{1,2,3\}$, and $A$, $I$, $Q_i$, and $R_i$ denote arbitrary sets.

To make sense of these polynomials, we need to define functor addition,%
\index{functor!sum of functors}
\index{functor!product of functors}
\index{polynomial functor!sum of polynomials}
\index{polynomial functor!product of polynomials}\index{product|seealso{polynomial functor, product of polynomials}}
 both in the binary case (i.e.\ what is $\yon^A+\yon^B$?) and more generally over arbitrary sets (i.e.\ what is $\sum_{i\in I}\yon^{A_i}$?).
This will allow us to interpret polynomials like \eqref{eqn.biz_poly}.
In particular, just as $3y^3=y^3+y^3+y^3$ in algebra, the summand $\3\yon^\3$ of \eqref{eqn.biz_poly} denotes the sum of representables $\yon^\3+\yon^\3+\yon^\3$, while the summand $\rr\yon^\zz$ denotes the sum over $\rr$ of copies of $\yon^\zz$.

While polynomial functors will be defined as sums, products of polynomials will turn out to be polynomials as well, again mimicking polynomials in algebra.
To make sense of these products, we will also define functor multiplication.
The construction of sums and products of functors $\smset\to\smset$ will rely on the construction of sets and products of sets themselves.

%-------- Section --------%
\section[Sums and products of sets]{Sums and products of sets%
  \sectionmark{Sums \& products of sets}}
\sectionmark{Sums \& products of sets}
\label{sec.poly.rep-sets.sum-prod-set}

\index{set!sum of sets}
\index{set!product of sets}
\index{set!indexing}\index{category!discrete}\index{indexed set}

Let $I$ be a set, and let $X_i$ be a set for each $i\in I$.
Classically, we may denote this \emph{$I$-indexed family of sets} by $(X_i)_{i\in I}$.
Categorically, we may view this data as a specific kind of functor: if we identify the set $I$ with the \emph{discrete category} on $I$, whose objects are the elements of $I$ and whose morphisms are all identities, then $(X_i)_{i\in I}$ can be identified with a functor $X\colon I\to\smset$ with $X(i)\coloneqq X_i$.
\index{indexed family of objects}
To compromise, we will denote an indexed family of sets by $X\colon I\to\smset$ for a set $I$ viewed as a discrete category (although we will occasionally use the classical notation when convenient), but denote the set obtained by evaluating $X$ at each $i\in I$ by $X_i$ rather than $X(i)$.

To pick out an element of one of the sets in the indexed family $X\colon I\to\smset$, we need to specify both an index $i\in I$ and an element $x\in X_i$.
We call the set of such pairs $(i,x)$ the \emph{sum} of this indexed family, as below.

\begin{definition}[Sum of sets] \label{def.sum_sets}
    Let $I$ be a set and $X\colon I\to\smset$ be an $I$-indexed family of sets.
    The \emph{sum $\sum_{i\in I}X_i$ of the indexed family $X$} is the set
    \[
    \sum_{i\in I}X_i\coloneqq\{(i,x)\mid i\in I\text{ and }x\in X_i\}.
    \]
    When $I\coloneqq\{i_1,\ldots,i_n\}$ is finite, we may alternatively denote this sum as
    \[
    X_{i_1}+\cdots+X_{i_n}.
    \]
\end{definition}
\index{element!of a sum of sets}
\index{element!of a product of sets}
\index{dependent function}\index{indexed set}

Say instead we pick an element from \emph{every} set in the indexed family: that is, we construct an assignment $i\mapsto x_i$, where each $x_i\in X_i$.
If every $X_i$ were the same set $X$, then this would just be a function $I\to X$.
More generally, this assignment is what we call a \emph{dependent function}: its codomain $X_i$ \emph{depends} on its input $i$.
\index{dependent function}
We write the signature of such a dependent function as
\[
f \colon (i \in I) \to X_i.
\]
Note that the indexed family of sets $X\colon I\to\smset$ completely determines this signature.
The set of all dependent functions whose signature is determined by a given indexed family of sets is the \emph{product} of that indexed family, as below.

\begin{definition}[Product of sets] \label{def.prod_sets}
    Let $I$ be a set and $X\colon I\to\smset$ be an $I$-indexed family of sets.
    The \emph{product $\prod_{i\in I}X_i$ of the indexed family $X$} is the set of dependent functions
    \[
    \prod_{i\in I}X_i\coloneqq\{f \colon (i \in I) \to X_i\}.
    \]
    When $I\coloneqq\{i_1,\ldots,i_n\}$ is finite, we may alternatively denote this product as
    \[
    X_{i_1}\times\cdots\times X_{i_n} \qqor X_{i_1}\cdots X_{i_n}.
    \]
\end{definition}\index{dependent function}

For a dependent function $f\colon(i\in I)\to X_i$, we may denote the element of $X_i$ that $f$ assigns to $i\in I$ by $f(i), fi,$ or $f_i$.
When $I\coloneqq\{i_1,\ldots,i_n\}$ is finite, we may identify $f$ with the $n$-tuple $(f(i_1),\ldots,f(i_n))$; similarly, when $I\coloneqq\nn$, we may identify $f$ with the infinite sequence $(f_0,f_1,f_2,\ldots)$.

\index{set!cardinality of}
\begin{example}\label{ex.two_sums_and_prods}
    If $I\coloneqq\2=\{1,2\}$, then an $I$-indexed family $X\colon I\to \smset$ consists of two sets---say $X_1\coloneqq\{a,b,c\}$ and $X_2\coloneqq\{c,d\}$.
    Their sum is then the disjoint union
    \[
    \sum_{i\in \2}X_i=X_1+X_2=\{(1,a),(1,b),(1,c),(2,c),(2,d)\}.
    \]
    The cardinality\tablefootnote{The \emph{cardinality} of a set is the number of elements it contains, at least when the set is finite; with care the notion can be extended to infinite sets as well.} of $X_1+X_2$ will always be the sum of the cardinalities of $X_1$ and $X_2$, justifying the use of the word ``sum.''

    Meanwhile, their product is the usual cartesian product
    \index{cartesian product|see{product}}\index{product}
    \[\prod_{i\in \2}X_i \cong X_1\times X_2=\{(a,c),(a,d),(b,c),(b,d),(c,c),(c,d)\}.\]
    The cardinality of $X_1\times X_2$ will always be the product of the cardinalities of $X_1$ and $X_2$, justifying the use of the word ``product.''
\end{example}

\begin{exercise}\label{exc.on_sums_prods_sets}\index{indexed set}
    Let $I$ be a set.
    \begin{enumerate}
        \item \label{exc.on_sums_prods_sets.sum} Show that there is an isomorphism of sets $I\iso\sum_{i\in I}\1$.
        \item \label{exc.on_sums_prods_sets.prod} Show that there is an isomorphism of sets $\1\iso\prod_{i\in I}\1$.
    \end{enumerate}
    As a special case, suppose that $I\coloneqq\0=\varnothing$ and that $X\colon \varnothing\to\smset$ is the unique empty indexed family of sets.
    \begin{enumerate}[resume]
        \item Is it true that $X_i=\1$ for each $i\in I$?
        \item Justify the statement ``the empty sum is $\0$'' by showing that there is an isomorphism of sets $\sum_{i\in\varnothing}X_i\iso\0$.
        \item Justify the statement ``the empty product is $\1$'' by showing that there is an isomorphism of sets $\prod_{i\in\varnothing}X_i\iso\1$.
        \qedhere
    \end{enumerate}

    \begin{solution}
        \begin{enumerate}
            \item \label{sol.on_sums_prods_sets.sum}
            To show that $I\iso\sum_{i \in I}\1$, observe that $x \in \1 = \{1\}$ if and only if $x = 1$, so $\sum_{i \in I} \1 = \{(i, 1) \mid i \in I\}$.
            Then the function $I \to \sum_{i \in I} \1$ that sends each $i \in I$ to $(i, 1)$ is clearly an isomorphism.

            \item \label{sol.on_sums_prods_sets.prod}
            To show that $\1 \iso \prod_{i \in I} \1$, it suffices to show that there is a unique dependent function $f \colon (i \in I) \to \1$.
            As $\1 = \{1\}$, such a function $f$ must always send $i \in I$ to $1$.
            This uniquely characterizes $f$, so there is indeed only one such dependent function.

            \item \label{sol.on_sums_prods_sets.vac} Yes: since $I$ is empty, there are no $i \in I$.
            So it is true that $X_i = 1$ holds whenever $i \in I$ holds, because $i \in I$ never holds.
            We say that this sort of statement is \emph{vacuously true}.

            \item As $I = \0 = \varnothing$, we have $\sum_{i\in\varnothing}X_i = \sum_{i\in I}\1 \iso I = \0$, where the equation on the left follows from \cref{sol.on_sums_prods_sets.vac} and the isomorphism in the middle follows from \cref{sol.on_sums_prods_sets.sum}.

            \item As $I = \varnothing$, we have $\prod_{i\in\varnothing}X_i = \prod_{i\in I}\1 \iso \1$, where the equation on the left follows from \cref{sol.on_sums_prods_sets.vac} and the isomorphism on the right follows from \cref{sol.on_sums_prods_sets.prod}.
        \end{enumerate}
    \end{solution}
\end{exercise}

The following standard fact describes the constructions from \cref{def.sum_sets,def.prod_sets} categorically and further justifies why we call them sums and products.
\index{product}
\index{coproduct}
\index{diagram}

\begin{proposition} \label{prop.set_prod_coprod}\index{functor!set-valued}
    Let $I$ be a set and $X\colon I\to\smset$ be an $I$-indexed family of sets. Then the sum $\sum_{i\in I}X_i$ is the categorical coproduct of these sets in $\smset$ (i.e.\ the colimit of the functor $X\colon I\to\smset$, viewed as a diagram), and the product $\prod_{i\in I}X_i$ is the categorical product of these sets in $\smset$ (i.e.\ the limit of the functor $X\colon I\to\smset$, viewed as a diagram).
\end{proposition}
\index{coproduct!inclusion map into}
\index{coproduct!universal property of}
\index{product!projection map out of}
\index{product!universal property of}\index{colimit}

\begin{proof}
    The sum $\sum_{i\in I}X_i$ comes equipped with an inclusion $\iota_j\colon X_j\to\sum_{i\in I}X_i$ for each $j\in I$ given by $x\mapsto(j,x)$.
    The product $\prod_{i\in I}X_i$ comes equipped with a projection $\pi_j\colon\prod_{i\in I}X_i\to X_j$ for each $j\in I$ sending each $f\colon(i\in I)\to X_i$ to $f(j)$.
    These satisfy the universal properties for categorical coproducts and products, respectively; see \cref{exc.set_prod_coprod}.
\end{proof}

\begin{exercise} \label{exc.set_prod_coprod}
    \begin{enumerate}
        \item Show that $\sum_{i\in I}X_i$ along with the inclusions $\iota_j\colon X_j\to\sum_{i\in I}X_i$ described in the proof of \cref{prop.set_prod_coprod} satisfy the universal property of a categorical coproduct: for any set $Y$ with functions $g_j\colon X_j\to Y$ for each $j\in I$, there exists a unique function $h\colon\sum_{i\in I}X_i\to Y$ for which $\iota_j\then h=g_j$ for all $j\in I$.
        \item Show that $\prod_{i\in I}X_i$ along with the projections $\pi_j\colon\prod_{i\in I}X_i\to X_j$ described in the proof of \cref{prop.set_prod_coprod} satisfy the universal property of a categorical product: for any set $Y$ with functions $g_j\colon Y\to X_j$ for each $j\in I$, there exists a unique function $h\colon Y\to\prod_{i\in I}X_i$ for which $h\then\pi_j=g_j$ for all $j\in I$. \qedhere
    \end{enumerate}
    \begin{solution}
        \begin{enumerate}
            \item Any function $h\colon\sum_{i\in I}X_i\to Y$ for which $\iota_j\then h=g_j$ for all $j\in I$ must satisfy $h(j,x)=h(\iota_j(x))=g_j(x)$ for all $j\in I$ and $x\in X_j$.
            This uniquely characterizes $h$, so if we define $h(j,x)\coloneqq g_j(x)$ we are done.

            \item Any function $h\colon Y\to\prod_{i\in I}X_i$ for which $h\then\pi_j=g_j$ for all $j\in I$ must satisfy $h(y)_j=\pi_j(h(y))=g_j(y)$ for all $y\in Y$ and $j\in I$.
            This uniquely characterizes $h(y)$ and thus $h$, so if we define $h(y)\colon(i\in I)\to X_i$ to be the dependent function given by $i\mapsto g_i(y)$ we are done.
        \end{enumerate}
    \end{solution}
\end{exercise}

Though we proved above explicitly that $\smset$ has all small products and coproducts, from here on out, we will assume the standard categorical fact that $\smset$ is complete (has all small limits) and cocomplete (has all small colimits). % TODO: add ref for this??
\index{limit}
\index{colimit}

We have constructed categorical sums and products of sets, but we can also construct categorical sums and products of the maps between them: functions.

\index{coproduct!of functions}
\index{product!of functions}

\begin{definition}[Categorical sum and product of functions] \label{def.sum-prod-func}
    Let $I$ be a set and $X,Y\colon I\to\smset$ be $I$-indexed families of sets.
    Given a natural transformation $f\colon X\to Y$, i.e.\ an $I$-indexed family of functions $(f_i\colon X_i\to Y_i)_{i\in I}$, its \emph{categorical sum} (or \emph{coproduct}) is the function
    \[
    \sum_{i\in I}f_i\colon\sum_{i\in I}X_i\to\sum_{i\in I}Y_i
    \]
    that, given $i\in I$ and $x\in X_i$, sends $(i,x)\mapsto(i,f_i(x))$; while its \emph{categorical product} is the function
    \[
    \prod_{i\in I}f_i\colon\prod_{i\in I}X_i\to\prod_{i\in I}Y_i
    \]
    that sends each $g\colon(i\in I)\to X_i$ to the \emph{composite dependent function} $(i\in I)\to Y_i$, denoted $g\then f$ or $f\circ g$, which sends $i\in I$ to $f_i(g(i))$.

    When $I\coloneqq\{i_1,\ldots,i_n\}$ is finite, we may alternatively denote this categorical sum and product of functions respectively as\tablefootnote{We will take care to highlight when this notation may clash with a sum (resp.\ product) of functions with common domain and codomain whose codomain has an additive (resp.\ multiplicative) structure.}
    \[
    f_{i_1}+\cdots+f_{i_n} \qqand f_{i_1}\times\cdots\times f_{i_n}.
    \]
\end{definition}

\begin{exercise} \label{exc.sum-prod-func}
    \begin{enumerate}
        \item Show that the categorical sum of functions is the one induced by the universal property of the categorical sum of sets.
        That is, given a set $I$, two $I$-indexed families of sets $X,Y\colon I\to\smset$, and a natural transformation $f\colon X\to Y$, the function $\sum_{i\in I}f_i\colon\sum_{i\in I}X_i\to\sum_{i\in I}Y_i$ that we called the categorical sum is induced by the following composite maps for $j\in I$:
        \[
        X_j\To{f_j}Y_j\To{\iota'_j}\sum_{i\in I}Y_i,
        \]
        where $\iota'_j$ is the inclusion.
        It then follows by a standard categorical argument that the sum is functorial, i.e.\ that the sum of identities is an identity and that the sum of composites is the composite of sums.

        \item Similarly, show that the categorical product of functions is the one induced by the universal property of the categorical product of sets.
        That is, given the same setup as the previous part, the function $\prod_{i\in I}f_i\colon\prod_{i\in I}X_i\to\prod_{i\in I}Y_i$ that we called the categorical product is induced by the following composite maps for $j\in J$:
        \[
        \prod_{i\in I}X_i\To{\pi_j}X_j\To{f_j}Y_j.
        \]
        Again, this implies that the product is functorial. \qedhere
    \end{enumerate}
    \begin{solution}
        \begin{enumerate}
            \item It suffices to show that the following square, where the vertical maps are the inclusions, commutes for all $j\in I$:
            \[
            \begin{tikzcd}[column sep=large]
                X_j \ar[d, "\iota_j"] \ar[r, "f_j"] & Y_j \ar[d, "\iota'_j"] \\
                \sum_{i\in I}X_i \ar[r, "\sum_{i\in I}f_i"] & \sum_{i\in I}Y_i
            \end{tikzcd}
            \]
            Given $x\in X_j$, the left inclusion map sends $x$ to $(j,x)$, which the bottom sum of maps sends to $(j,f_j(x))$.
            Meanwhile, the top map sends $x$ to $f_j(x)$, which the right inclusion map again sends to $(j,f_j(x))$.

            \item It suffices to show that the following square, where the vertical maps are the projections, commutes for all $j\in I$:
            \[
            \begin{tikzcd}[column sep=large]
                \prod_{i\in I}X_i \ar[d, "\pi_j"] \ar[r, "\prod_{i\in I}f_i"] & \prod_{i\in I}Y_i \ar[d, "\pi'_j"] \\
                X_j \ar[r, "f_j"] & Y_j
            \end{tikzcd}
            \]
            Given $g\colon(i\in I)\to X_i$ in $\prod_{i\in I}X_i$, the top product of maps sends $g$ to $f\circ g$, which the right projection map sends to $f_j(g(j))$.
            Meanwhile, the left projection map sends $g$ to $g(j)$, which the bottom map again sends to $f_j(g(j))$.
        \end{enumerate}
    \end{solution}
\end{exercise}

We now highlight some tools and techniques to help us work with sum and product sets.

\begin{exercise}\label{exc.product_as_sections}
    Let $I$ be a set and $X \colon I \to \smset$ be an indexed family.
    There is a
    projection function
    $\pi_1 \colon \sum_{i \in I} X_i \to I$
    defined by $\pi_1(i, x) \coloneqq i$.
    \begin{enumerate}
        \item What is the signature of the second projection $\pi_2(i, x) \coloneqq x$?
        (Hint: it's a dependent function.)
        \item A \emph{section} of a function $r \colon A \to B$ is a function $s \colon B \to A$ such that $s \then r = \id_B$.
        Show that the product of the indexed family is isomorphic to the set of sections of $\pi_1$:
        \[\prod_{i \in I} X_i \cong \left\{s \colon I \to \sum_{i \in I} X_i \,\middle|\, s \then \pi_1 = \id_I\right\}.\]
        \qedhere
    \end{enumerate}
    \begin{solution}
        \begin{enumerate}
            \item The second projection $\pi_2(i, x) = x$ sends each pair $p \coloneqq (i, x) \in \sum_{i \in I} X_i$ to $x$, an element of $X_i$.
            Note that we can write $i$ in terms of $p$ as $\pi_1(p)$.
            This allows us to write the signature of $\pi_2$ as $\pi_2 \colon (p \in \sum_{i \in I} X_i) \to X_{\pi_1(p)}$.

            \item Let $S := \{s \colon I \to \sum_{i \in I} X_i \mid s \then \pi_1 = \id_I\}$ be the set of sections of $\pi_1$. To show that $\prod_{i \in I} X_i \cong S$, we will exhibit maps in either direction and show that they are mutually inverse.
            For each $f \colon (i \in I) \to X_i$ in $\prod_{i \in I} X_i$, we have $f(i) \in X_i$ for $i \in I$, so we can define a function $s_f \colon I \to \sum_{i \in I} X_i$ that sends $i\mapsto(i, f(i))$.
            Then $\pi_1(s_f(i)) = \pi_1(i, f(i)) = i$, so $s_f$ is a section of $\pi_1$.
            Hence $f \mapsto s_f$ is a map $\prod_{i \in I} X_i \to S$.

            In the other direction, for each section $s \colon I \to \sum_{i \in I} X_i$ we have $\pi_1(s(i)) = i$ for $i \in I$, so we can write $s(i)$ as an ordered pair $(i, \pi_2(s(i)))$ with $\pi_2(s(i)) \in X_i$.
            Hence we can define a dependent function $f_s \colon (i \in I) \to X_i$ sending $i\mapsto\pi_2(s(i))$.
            Then $s \mapsto f_s$ is a map $S \to \prod_{i \in I} X_i$.
            By construction $s_{f_s}(i) = (i, f_s(i)) = (\pi_1(s(i)), \pi_2(s(i))) = s(i)$ and $f_{s_f}(i) = \pi_2(s_f(i)) = \pi_2(i, f(i)) = f(i)$, so these maps are mutually inverse.
        \end{enumerate}
    \end{solution}
\end{exercise}\index{dependent function}

A helpful way to think about sum or product sets is to consider what choices must be made to specify an element of such a set.
In the following examples, say that we have a set $I$ and an $I$-indexed family $X \colon I \to \smset$.

Below, we give the instructions for choosing an element of $\sum_{i \in I} X_i$.

\begin{quote}
    To choose an element of $\sum_{i \in I} X_i$:
    \begin{enumerate}
        \item choose an element $i \in I$;
        \item choose an element of $X_i$.
    \end{enumerate}
\end{quote}
\index{element!of a dependent sum}

Then the projection $\pi_1$ from \cref{exc.product_as_sections} sends each element of $\sum_{i \in I} X_i$ to the element of $i \in I$ chosen in step 1, while the projection $\pi_2$ sends each element of $\sum_{i \in I} X_i$ to the element of $X_i$ chosen in step 2.

Next, we give the instructions for choosing an element of $\prod_{i \in I} X_i$.

\begin{quote}
    To choose an element of $\prod_{i \in I} X_i$:
    \begin{enumerate}
        \item for each element $i \in I$:
        \begin{enumerate}[label*=\arabic*.]
            \item choose an element of $X_i$.
        \end{enumerate}
    \end{enumerate}
\end{quote}

\index{element!of a dependent product}
\index{element!of a nested dependent set}

Armed with these interpretations, we can tackle more complicated expressions, including those with nested $\sum$'s and $\prod$'s such as
\begin{equation}\label{eqn.sum_prod_sum}
    A \coloneqq \sum_{i\in I}\prod_{j\in J(i)}\sum_{k\in K(i,j)}X(i,j,k).
\end{equation}
The instructions for choosing an element of $A$ form a nested list, as follows.

\begin{quote}
    To choose an element of $A$:
    \begin{enumerate}
        \item choose an element $i \in I$;
        \item for each element $j \in J(i)$:
        \begin{enumerate}[label*=\arabic*.]
            \item choose an element $k \in K(i,j)$;
            \item choose an element of $X(i,j,k)$.
        \end{enumerate}
    \end{enumerate}
\end{quote}

Here the choice of $k\in K(i,j)$ may depend on $i$ and $j$: different values of $i$ and $j$ may lead to different sets $K(i,j)$.

By describing $A$ like this, we see that each $a \in A$ can be projected to an element $i\coloneqq\pi_1(a) \in I$, chosen in step 1, and a dependent function $\pi_2(a)$, chosen in step 2.
This dependent function in turn sends each $j \in J(i)$ to a pair that can be projected to an element $k\coloneqq\pi_1(\pi_2(a)(j)) \in K(i, j)$ chosen in step 2.1 and an element $\pi_2(\pi_2(a)(j)) \in X(i,j,k)$ chosen in step 2.2.

\begin{example}%[Notation for $\sum\prod$ stuff]
    \label{ex.notation_sum_prod}
    % Here we give notation for the elements of a set involving $\sum$'s and $\prod$'s such as that in \eqref{eqn.sum_prod_sum}.

    Let $I\coloneqq\{1,2\}$; let $J(1)\coloneqq\{j\}$ and $J(2)\coloneqq\{j,j'\}$; let $K(1,j)\coloneqq\{k_1,k_2\}$, $K(2,j)\coloneqq\{k_1\}$, and $K(2,j')\coloneqq\{k'\}$; and let $X(i,j,k)\coloneqq\{x,y\}$ for all $i,j,k$. Now the formula
    \[\sum_{i\in I}\prod_{j\in J(i)}\sum_{k\in K(i,j)}X(i,j,k)\]
    from \eqref{eqn.sum_prod_sum} specifies a fixed set:
    \[
    \left\{
    \begin{gathered}
        \big(1, j\mapsto(k_1,x)\big),
        \big(1, j\mapsto(k_1,y)\big),
        \big(1, j\mapsto(k_2,x)\big),
        \big(1, j\mapsto(k_2,y)\big),\\
        \big(2, j\mapsto(k_1,x), j'\mapsto(k',x)\big),
        \big(2, j\mapsto(k_1,x), j'\mapsto(k',y)\big),\\
        \big(2, j\mapsto(k_1,y), j'\mapsto(k',x)\big),
        \big(2, j\mapsto(k_1,y), j'\mapsto(k',y)\big)
    \end{gathered}
    \right\}.
    \]
    In each case, we first chose an element $i\in I$, either 1 or 2. Then for each $j\in J(i)$ we chose an element $k\in K(i,j)$ and an element of $X(i,j,k)$.
\end{example}

\begin{exercise}
    Consider the set
    \begin{equation}\label{eqn.prod_sum_prod}B \coloneqq \prod_{i\in I}\sum_{j\in J(i)}\prod_{k\in K(i,j)}X(i,j,k).\end{equation}
    \begin{enumerate}
        \item Give the instructions for choosing an element of $B$ as a nested list, like we did for $A$ just below \eqref{eqn.sum_prod_sum}.
        \item With $I$, $J$, $K$, and $X$ as in \cref{ex.notation_sum_prod}, how many elements are in $B$?
        \item Write out three of these elements in the style of \cref{ex.notation_sum_prod}.
        \qedhere
    \end{enumerate}
    \begin{solution}
        \begin{enumerate}
            \item Here are the instructions for choosing an element of $B$ as a nested list.
            \begin{quote}
                To choose an element of $B$:
                \begin{enumerate}[label=\arabic*.]
                    \item for each element $i \in I$:
                    \begin{enumerate}[label*=\arabic*.]
                        \item choose an element $j \in J(i)$;
                        \item for each element $k \in K(i, j)$:
                        \begin{enumerate}[label*=\arabic*.]
                            \item choose an element of $X(i,j,k)$.
                        \end{enumerate}
                    \end{enumerate}
                \end{enumerate}
            \end{quote}
            \item Given $I\coloneqq\{1,2\}$, $J(1)\coloneqq\{j\}$, $J(2)\coloneqq\{j,j'\}$, $K(1,j)\coloneqq\{k_1,k_2\}$, $K(2,j)\coloneqq\{k_1\}$, $K(2,j')\coloneqq\{k'\}$, and $X(i,j,k)\coloneqq\{x,y\}$ for all $i,j,k$, our goal is to count the number of elements in $B$.
            To compute the cardinality of $B$, we can use the fact that the cardinality of a sum (resp.\ product) is the sum (resp.\ product) of the cardinalities of the summands (resp.\ factors).
            So
            \begin{align*}
                |B| &= \prod_{i\in I}\sum_{j\in J(i)}\prod_{k\in K(i,j)}|X(i,j,k)| \\
                &= \prod_{i\in \{1,2\}}\sum_{j\in J(i)}\prod_{k\in K(i,j)}2 \\
                &= \left(\sum_{j\in J(1)} 2^{|K(1,j)|}\right)\left(\sum_{j\in J(2)} 2^{|K(2,j)|}\right) \\
                &= \left(2^2\right)\left(2^1 + 2^1\right) = 16.
            \end{align*}
            \item Here are three of the elements of $B$ (you may have written down others):
            \begin{itemize}
                \item $(1 \mapsto (j, k_1 \mapsto x, k_2 \mapsto y), 2 \mapsto (j', k' \mapsto x))$
                \item $(1 \mapsto (j, k_1 \mapsto y, k_2 \mapsto y), 2 \mapsto (j, k_1 \mapsto y))$
                \item $(1 \mapsto (j, k_1 \mapsto y, k_2 \mapsto x), 2 \mapsto (j', k' \mapsto y))$
            \end{itemize}
            \qedhere
        \end{enumerate}
    \end{solution}
\end{exercise}

%Henceforth, we will omit the full sequence of nested instructions corresponding to every sum or product of sets; we will assume you can read them for yourself.

%-------- Section --------%
\section{Expanding products of sums} \label{sec.poly.rep-sets.expand}

\index{type theoretic axiom of choice}
\index{distributive law}

We will often encounter sums of sets nested within products, as in \eqref{eqn.sum_prod_sum} and \eqref{eqn.prod_sum_prod}.
The following proposition helps us work with these; it is sometimes called the \emph{type-theoretic axiom of choice}, but it is perhaps more familiar as a set-theoretic analogue of the \emph{distributive property} of multiplication over addition.
While the identity may look foreign, it captures for sets the same process that you would use to multiply multi-digit numbers from grade school arithmetic or polynomials from high school algebra.

\begin{proposition}[Pushing $\prod$ past $\sum$]\label{prop.push_prod_sum_set}
    For any sets $I,(J(i))_{i\in I},$ and $(X(i,j))_{i\in I, j\in J(i)}$, we have a natural isomorphism\index{isomorphism!natural}
    \begin{equation}\label{eqn.set_completely_distributive}
        \prod_{i\in I}\sum_{j\in J(i)}X(i,j)
        \iso
        \sum_{\bar{j}\in \prod_{i\in I}J(i)}\;\prod_{i\in I}X(i,\bar{j}(i)).\tablefootnote{We draw a bar over $j$ in $\bar{j}$ to remind ourselves that $\bar{j}$ is no longer just an index but a (dependent) function.}
    \end{equation}
\end{proposition}
\begin{proof}
    First, we construct a map from the left hand set to the right. An element of the set on the left is a dependent function $f \colon (i \in I) \to \sum_{j \in J(i)} X(i, j)$, which we can compose with projections from its codomain to yield $\pi_1(f(i)) \in J(i)$ and $\pi_2(f(i)) \in X(i, \pi_1(f(i)))$ for every $i \in I$.
    We can then form the following pair:\footnote{We omit parentheses for function application here and throughout the text for compactness whenever the meaning is clear.}
    \[
    (i \mapsto \pi_1 fi, i \mapsto \pi_2 fi).
    \]
    This is an element of the right hand set, because $i \mapsto \pi_1 fi$ is a dependent function in $\prod_{i\in I}J(i)$ and $i \mapsto \pi_2 fi$ is a dependent function in $\prod_{i\in I}X(i,\pi_1 fi)$.

    Now we go from right to left.
    An element of the right hand set is a pair of dependent functions, $\bar{j} \colon (i \in I) \to J(i)$ and $g \colon
    (i \in I) \to X(i, \bar{j}i)$.
    We map this pair to the following element of the left hand set, a dependent function $(i\in I)\to\sum_{j\in J(i)}X(i,j)$:
    \[
    i \mapsto (\bar{j}i, gi).
    \]

    Finally, we verify that the maps are mutually inverse.
    An element $(\bar{j}, g)$ of the right hand set is sent by one map and then the other to the pair
    \[
    (i \mapsto \pi_1(\bar{j}i, gi), i \mapsto \pi_2(\bar{j}i, gi))=(i \mapsto \bar{j}i, i \mapsto gi)=(\bar{j},g)
    \]
    Meanwhile, an element $f$ of the left hand set is sent by one map and then the other to the dependent function
    \[
    i \mapsto (\pi_1 fi, \pi_2 fi).
    \]
    As $fi$ is a pair whose components are $\pi_1 fi$ and $\pi_2 fi$, the dependent function above is precisely $f$.
\end{proof}

When $J(i)=J$ does not depend on $i\in I$, we can simplify the formula in \eqref{eqn.set_completely_distributive}.

\begin{corollary} \label{cor.push_prod_sum_set_indep}
    For any sets $I, J,$ and $(X(i, j))_{i \in I, j \in J}$, we have a natural isomorphism\index{isomorphism!natural}
    \begin{equation} \label{eqn.push_prod_sum_set_indep}
        \prod_{i\in I}\sum_{j\in J}X(i,j)\cong\sum_{\bar{j}\colon I\to J}\prod_{i\in I}X(i,\bar{j}i),
    \end{equation}
    where $\bar{j}$ ranges over all (standard, non-dependent) functions $I\to J$.
\end{corollary}
\begin{proof}
    Take $J(i)\coloneqq J$ for all $i \in I$ in \eqref{eqn.set_completely_distributive}.
    Then the set $\prod_{i\in I}J(i)$ becomes $\prod_{i\in I}J$ (which we may denote in exponential form by $J^I$); its elements, dependent functions $\bar{j}\colon(i\in I)\to J(i)=J$, become standard functions $\bar{j}\colon I\to J$.
\end{proof}

It turns out that being able to push $\prod$ past $\sum$ as in \eqref{eqn.set_completely_distributive} is not a property that is unique to sets.
In general, we refer to a category having this property as follows.

\begin{definition}[Completely distributive category]\index{coproduct}\index{product}
    A category $\Cat{C}$ with all small products and coproducts is \emph{completely distributive}\tablefootnote{While our terminology generalizes that of a completely distributive lattice, which has the additional requirement that the category be a poset, it is unfortunately not standard: a completely distributive category refers to a different concept in some categorical literature. We will not use this other concept, so there is no ambiguity.} if products distribute over coproducts as in \eqref{eqn.set_completely_distributive}; that is, for any set $I$, sets $(J(i))_{i\in I}$, and objects $(X(i,j))_{i\in I,j\in J(i)}$ from $\Cat{C}$, we have a natural isomorphism\index{isomorphism!natural}
    \begin{equation}\label{eqn.cat_completely_distributive}
        \prod_{i\in I}\sum_{j\in J(i)}X(i,j)
        \iso
        \sum_{\bar{j}\in \prod_{i\in I}J(i)}\;\prod_{i\in I}X(i,\bar{j}i).
    \end{equation}
\end{definition}

The term ``completely distributive'' comes from lattice theory. As such it is consistent with two different extensions to categories that may not be posets. We use it to mean that the category has all sums and products and that products distribute over sums. Other authors use it to mean that the category has all colimits and limits and a sort of distributivity between them.
\index{completely distributive category}

So \cref{prop.push_prod_sum_set} states that $\smset$ is completely distributive.
Once we define the category of polynomial functors, we will see that it, too, is completely distributive.

\index{completely distributive category!$\poly$ as}
\index{completely distributive category!$\smset$ as}

\cref{cor.push_prod_sum_set_indep} generalizes to all completely distributive categories as well; we state this formally below.

\begin{corollary} \label{cor.push_prod_sum_obj_indep}
    Let $\Cat{C}$ be a completely distributive category.
    For any sets $I$ and $J$ and objects $(X(i, j))_{i \in I, j \in J}$ in $\Cat{C}$, we have a natural isomorphism\index{isomorphism!natural}
    \begin{equation} \label{eqn.push_prod_sum_obj_indep}
        \prod_{i\in I}\sum_{j\in J}X(i,j)\iso\sum_{\bar{j}\colon I\to J}\prod_{i\in I}X(i,\bar{j}i).
    \end{equation}
\end{corollary}
\begin{proof}
    Again, take $J(i)\coloneqq J$ for all $i \in I$ in \eqref{eqn.cat_completely_distributive}.
\end{proof}

\index{completely distributive category!distributive law in|see{distributive law}}
\index{distributive law}

\begin{exercise}
    Let $\Cat{C}$ be a completely distributive category.
    How is the usual distributive law
    \[
    X\times(Y+Z)\iso X\times Y+X\times Z
    \]
    for $X,Y,Z\in\cat{C}$ a special case of \eqref{eqn.cat_completely_distributive}?
    \begin{solution}
        We wish to show that $X\times (Y+Z)\iso X\times Y+X\times Z$ using \eqref{eqn.cat_completely_distributive}.
        On the left hand side, we are taking a $2$-fold product: a single object times a $2$-fold sum.
        So we should let $I\coloneqq\2$ and let $J(1)\coloneqq\1$, with $X(1,1)\coloneqq X$; and $J(2)\coloneqq\2$, with $X(2,1)\coloneqq Y$ and $X(2,2)\coloneqq Z$.
        Then
        \[
        X\times(Y+Z)\iso \prod_{i\in\2}\sum_{j\in J(i)}X(i,j) \iso \sum_{\bar{j}\in\prod_{i\in\2}J(i)}\;\prod_{i\in \2}X(i,\bar{j}(i)) \iso \sum_{\bar{j}\in\prod_{i\in\2}J(i)}X(1,\bar{j}(1))\times X(2,\bar{j}(2)),
        \]
        where the middle isomorphism follows from \eqref{eqn.cat_completely_distributive}.
        The set $\prod_{i\in\2}J(i)$ contains two functions: $(1\mapsto1,2\mapsto1)$ and $(1\mapsto1,2\mapsto2)$.
        So we can rewrite the right hand side as
        \[
        X(1,1)\times X(2,1)+X(1,1)\times X(2,2)\iso X\times Y+X\times Z.
        \]
    \end{solution}
\end{exercise}

Throughout this book, such as in the exercise below, you will see expressions consisting of alternating products and sums.
Using \eqref{eqn.cat_completely_distributive}, you can always rewrite such an expression as a sum of products, in which every $\sum$ appears before every $\prod$.\footnote{When an expression is written so that every $\sum$ appears before every $\prod$, it is said to be in \emph{disjunctive normal form}.}
This is analogous to how products of sums in high school algebra can always be expanded into sums of products via the distributive property.

\index{disjunctive normal form}

\begin{exercise} \label{exc.push_prod_sum_set}
    Let $I, (J(i))_{i\in I},$ and $(K(i,j))_{(i,\,j)\in IJ}$ be sets, and for each $(i,j,k)\in IJK$, let $X(i,j,k)$ be an object in a completely distributive category.
    \begin{enumerate}
        \item Rewrite
        \[
        \sum_{i\in I}\prod_{j\in J(i)}\sum_{k\in K(i,\,j)}X(i,j,k)
        \]
        so that every $\sum$ appears before every $\prod$.
        \item Rewrite
        \[
        \prod_{i\in I}\sum_{j\in J(i)}\prod_{k\in K(i,\,j)}X(i,j,k)
        \]
        so that every $\sum$ appears before every $\prod$.
        \item Rewrite
        \[
        \prod_{i\in I}\prod_{j\in J(i)}\sum_{k\in K(i,\,j)}X(i,j,k)
        \]
        so that every $\sum$ appears before every $\prod$.\qedhere
    \end{enumerate}
    \begin{solution}
        \begin{enumerate}
            \item By applying \eqref{eqn.cat_completely_distributive}, we can rewrite
            \[
            \sum_{i\in I}\prod_{j\in J(i)}\sum_{k\in K(i,j)}X(i,j,k)
            \]
            as
            \[
            \sum_{i\in I}\sum_{\bar{k}\in \prod_{j\in J}K(i,j)}\prod_{j\in J(i)}X(i,j,\bar{k}j).
            \]
            \item By applying \eqref{eqn.cat_completely_distributive}, we can rewrite
            \[
            \prod_{i\in I}\sum_{j\in J(i)}\prod_{k\in K(i,j)}X(i,j,k)
            \]
            as
            \[
            \sum_{\bar{j}\in \prod_{i\in I}J(i)}\;\prod_{i\in I}X(i,\bar{j}i)\prod_{k\in K(i,\bar{j}i)}X(i,\bar{j}i,k).
            \]
            \item By applying \eqref{eqn.cat_completely_distributive}, we can rewrite
            \[
            \prod_{i\in I}\prod_{j\in J(i)}\sum_{k\in K(i,j)}X(i,j,k)
            \]
            once as
            \[
            \prod_{i\in I}\sum_{\bar{k}\in\prod_{j\in J}K(i,j)}\prod_{j\in J(i)}X(i,j,\bar{k}j)
            \]
            and then again as
            \[
            \sum_{\bar{\bar{k}}\in\prod_{i\in I}\prod_{j\in J}K(i,j)}\prod_{i\in I}\prod_{j\in J(i)}X(i,j,\bar{\bar{k}}(i,j)).
            \]
        \end{enumerate}
    \end{solution}
\end{exercise}

Now that we understand sums and products of sets, we are ready to explore sums and products of set-valued functors.\index{functor!set-valued}

%-------- Section --------%
\section[Sums and products of functors $\smset\to\smset$]{Sums and products of functors $\smset\to\smset$%
  \sectionmark{Sums \& products of functors $\smset\to\smset$}}
\sectionmark{Sums \& products of functors $\smset\to\smset$} \label{sec.poly.rep-sets.sum-prod-func}

Recall that our goal is to define polynomial functors such as $\yon^\2+\2\yon+\1$ and the maps between them.
We have defined representable functors such as $\yon^\2$, $\yon$, and $\1$; we just need to interpret sums of functors $\smset\to\smset$.
But we might as well introduce products of functors at the same time, because they will very much come in handy.
Both these concepts generalize to limits and colimits in $\smset^\smset$.

\begin{proposition} \label{prop.presheaf_lim_ptwise}\index{functor!limit of functors}\index{functor!colimit of functors}\index{limit!of functors}\index{colimit!of functors}
    The category $\smset^\smset$ has all small limits and colimits, and they are computed pointwise.
    In particular, on objects, given a small category $\cat{J}$ and a functor $F\colon \cat{J}\to\smset^\smset$, for all $X\in\smset$, the limit and colimit of $F$ satisfy isomorphisms
    \[
      \left(\lim_{j\in\cat{J}} F(j)\right)(X) \iso \lim_{j\in\cat{J}} \left(F(j)(X)\right)
    \]
    and
    \[
      \left(\colim_{j\in\cat{J}} F(j)\right)(X) \iso \colim_{j\in\cat{J}} \left(F(j)(X)\right)
    \]
    natural in $X$.
\end{proposition}
\begin{proof}
    This is a special case of a more general fact when $\smset^\smset$ is replaced by an arbitrary functor category $\Cat{D}^\Cat{C}$, where $\Cat{D}$ is a category that (like $\smset$) has all small limits and colimits; see \cite[pages 22--23, displays (24) and (25)]{macLane1992sheaves}.
\end{proof}

\index{functor!category of functors}\index{limit!of functors $\smset\to\smset$}\index{colimit!of functors $\smset\to\smset$}

Focusing on the case of coproducts and products, the following corollary is immediate.
\index{coproduct}

\begin{corollary}[Sums and products of functors $\smset\to\smset$] \label{cor.sum_prod_set_endofuncs}
    Given functors $F,G\colon\smset\to\smset$, their categorical coproduct or sum in $\smset^\smset$, denoted $F+G$, is the functor $\smset\to\smset$ defined for $X,Y\in\smset$ and $f\colon X\to Y$ by
    \[
    (F+G)(X)\coloneqq F(X)+G(X) \qqand (F+G)(f)\coloneqq F(f)+G(f);
    \]
    while their categorical product in $\smset^\smset$, denoted $F\times G$ or $FG$, is the functor $\smset\to\smset$ defined for $X,Y\in\smset$ and $f\colon X\to Y$ by
    \[
    (F\times G)(X)\coloneqq F(X)\times G(X) \qqand (F\times G)(f)\coloneqq F(f)\times G(f).
    \]

    More generally, given functors $(F_i)_{i\in I}$ indexed over a set $I$, their categorical coproduct or sum and categorical product in $\smset^\smset$, respectively denoted
    \[
    \sum_{i\in I}F_i\colon\smset\to\smset
    \qqand
    \prod_{i\in I}F_i\colon\smset\to\smset,
    \]
    are respectively defined for $X\in\smset$ by
    \[
    \left(\sum_{i\in I}F_i\right)(X)\coloneqq\sum_{i\in I} F_i(X)
    \qqand
    \left(\prod_{i\in I}F_i\right)(X)\coloneqq\prod_{i\in I} F_i(X).
    \]
    and for functions $f\colon X\to Y$ by
    \[
    \left(\sum_{i\in I}F_i\right)(f)\coloneqq\sum_{i\in I} F_i(f)
    \qqand
    \left(\prod_{i\in I}F_i\right)(f)\coloneqq\prod_{i\in I} F_i(f).
    \]
\end{corollary}

% TODO: sum and products of nat trans & functoriality?

We also note the special case of initial and terminal objects.
Given a set $I\in\smset$, we will also use $I$ to denote the constant functor $\smset\to\smset$ that sends every set to $I$.

\index{functor!constant}

\begin{corollary}[Initial and terminal functors $\smset\to\smset$]
    The constant functor $\0\colon\smset\to\smset$ is initial in $\smset^\smset$, while the constant functor $\1\colon\smset\to\smset$ is terminal in $\smset^\smset$.
\end{corollary}
\begin{proof}
    As the set $\0$ is initial in $\smset$ (for every set $X$ there is a unique map $\0\to X$), \cref{prop.presheaf_lim_ptwise} implies that the constant functor $\0$ is initial in $\smset^\smset$.
    Similarly, as the set $\1$ is terminal in $\smset$ (for every set $X$ there is a unique map $X\to\1$), \cref{prop.presheaf_lim_ptwise} implies that the constant functor $\1$ is terminal in $\smset^\smset$.
\end{proof}

\index{functor!initial}\index{functor!terminal}

Finally, we note that $\smset^\smset$ inherits the distributivity of $\smset$.

\begin{proposition}\label{prop.set_endofunc_distrib}\index{completely distributive category!$\smset^\smset$ as}
    The category $\smset^\smset$ is completely distributive.
\end{proposition}
\begin{proof}
    This follows directly from the fact that $\smset$ itself is completely distributive (\cref{prop.push_prod_sum_set}) and the fact that sums and products in $\smset^\smset$ are computed pointwise (\cref{cor.sum_prod_set_endofuncs}).
\end{proof}

The following exercises justify some notational shortcuts we will use when denoting polynomial functors.
First, for any set $A$ and functor $F\colon\smset\to\smset$, we may write an $A$-indexed sum of copies of $F$ as $AF$, the product of $F$ and the constant functor $A$; for instance, $\yon+\yon\iso\2\yon$.

\begin{exercise} \label{exc.repeated_sum_is_product}
    Show that for a set $A\in\smset$ and a functor $F\colon\smset\to\smset$, an $A$-indexed sum of copies of $F$ is isomorphic to the product of the constant functor $A$ and $F$:
    \[
    \sum_{a \in A}F\iso AF.
    \]
    (This is analogous to the fact that adding up $n$ copies of number is equal to multiplying that same number by $n$.)
    \begin{solution}
        It suffices to show that for all $X\in\smset$, there is an isomorphism
        \[
        \sum_{a\in A}F(X)\iso(AF)(X)
        \]
        natural in $X$.
        The left hand side is the set $\{(a,s)\mid a\in A\text{ and }s\in F(X)\} \iso A \times F(X)$, while the right hand side is also naturally isomorphic to the set $A(X)\times F(X)\iso A\times F(X)$.
        Alternatively, since $\smset^\smset$ is completely distributive by \cref{prop.set_endofunc_distrib}, the result also follows from \eqref{eqn.cat_completely_distributive}, with $I\coloneqq\2, J(1)\coloneqq A, J(2)\coloneqq\1, X(1,a)\coloneqq\1$ (the constant functor) for $a\in A$, and $X(2,1)\coloneqq F$:
        \[
        AF \iso
        \left(\sum_{a\in A}\1\right)F \iso
        \prod_{i\in\2}\sum_{j\in J(i)}X(i,j)
        \iso
        \sum_{\ol{j}\in\prod_{i\in\2}J(i)}\prod_{i\in\2}X(i,\ol{j}(i))
        \iso
        \sum_{a\in A}\1\times F
        \iso
        \sum_{a\in A}F.
        \]
        Here we used the fact that $A\iso\sum_{a\in A}\1$ from \cref{exc.on_sums_prods_sets} \cref{exc.on_sums_prods_sets.sum} (there we proved the statement for sets, but the same statement for the corresponding constant set-valued functors follows immediately).
    \end{solution}
\end{exercise}\index{completely distributive category}

Similarly, we may wish to write an $A$-indexed product of copies of $F$ in exponential form as $F^A$.
But since we have already introduced exponential notation for representable functors, this yields two possible interpretations for the functor $\smset\to\smset$ denoted by $\yon^A$: as the functor represented by $A$, or as the $A$-indexed product of copies of the identity functor $\yon\colon\smset\to\smset$.
In fact, the following exercise shows that there is no ambiguity, as the two interpretations are isomorphic.

\begin{exercise}
    \begin{enumerate}
        \item Show that for a set $I\in\smset$, an $I$-indexed product of copies of the identity functor $\yon\colon\smset\to\smset$ is isomorphic to the functor $\yon^I\colon\smset\to\smset$ represented by $I$:
        \[
        \prod_{i\in I}\yon\iso\yon^I.
        \]
        (This is analogous to the fact that multiplying $n$ copies of a number together is equal to raising that same number to the power of $n$.)
        \item Show that the $I$-indexed product of copies of a representable functor $\yon^A\colon\smset\to\smset$ for some $A\in\smset$ is isomorphic to the functor $\yon^{IA}\colon\smset\to\smset$ represented by the product set $IA$:
        \[
        \prod_{i\in I}\yon^A\iso\yon^{IA}.
        \]
        (Hint: You may use the fact that following natural isomorphism holds between sets of functions:\index{isomorphism!natural}
        \[
        \{ f \colon I\times A\to X \} \iso \{ g \colon I\to X^A \}.
        \]
        The process of converting a function $f$ in the left hand set to the corresponding function $i\mapsto(a\mapsto f(i,a))$ in the right is known as \emph{currying}.) \qedhere
    \end{enumerate}
    \begin{solution}
        \begin{enumerate}
            \item It suffices to show that for all $X\in\smset$, there is an isomorphism
            \[
            \prod_{i\in I} \yon(X) \iso \yon^I(X).
            \]
            natural in $X$.
            We have that $\yon(X)\iso X$ and that $\yon^I(X)\iso X^I$.
            So both sides are naturally isomorphic to the set of functions $I\to X$.

            \item It suffices to show that for all $X\in\smset$, there is an isomorphism
            \[
            \prod_{i\in I} \yon^A(X) \iso \yon^{IA}(X).
            \]
            natural in $X$.
            We have that $\yon^A(X)\iso X^A$, so $\prod_{i\in I}\yon^A(X)\iso(X^A)^I$, and that $\yon^{IA}(X)\iso X^{IA}$.
            By currying, both sides are naturally isomorphic to the set of functions $I\times A\to X$.
        \end{enumerate}
    \end{solution}
\end{exercise}

Henceforth, given $A\in\smset$ and a functor $F\colon\smset\to\smset$, we define
\[
F^A\coloneqq\prod_{a\in A}F.
\]
The exercise above shows that this notation does not conflict with the way we write representable functors as powers of the identity functor $\yon$.
The exercise also shows how a power of a representable functor can be simplified to a single representable functor.

With these ingredients, we are finally ready to define what a polynomial functor is.
We will begin with this definition in the next chapter.

%-------- Section --------%
\section[Summary and further reading]{Summary and further reading%
  \sectionmark{Summary \& further reading}}
\sectionmark{Summary \& further reading}

In this chapter, we reviewed the definition of a \emph{representable functor} $\yon^S\colon\smset\to\smset$ for $S\in\smset$ sending $X\mapsto X^S$.
We then stated and proved the \emph{Yoneda lemma}, a foundational result characterizing maps out of these representables: for an arbitrary functor $F\colon\smset\to\smset$, natural transformations $\yon^S\to F$ are in natural correspondence with elements of $F(S)$.

\index{functor!representable}\index{Yoneda lemma}\index{category!discrete}\index{coproduct}\index{product}

We then reviewed other categorical constructions in $\smset$, many of which carry over to the polynomial functors we introduce in the next chapter.
For a set $I$, we can view it as a discrete category and consider a functor $X\colon I\to\smset$ as an \emph{$I$-indexed family of sets} comprised of a set $X_i$ for each $i\in I$.
An $I$-indexed family of sets $X$ has a \emph{sum} (or \emph{coproduct}), the set of pairs $(i,x)$ with $i\in I$ and $x\in X_i$; and a \emph{product}, the set of \emph{dependent functions} $f\colon(i\in I)\to X_i$.
Such a dependent function sends each $i\in I$ to an element of $X_i$.
These constructions satisfy the universal properties of coproducts and products; moreover, products distribute over coproducts, making $\smset$ a completely distributive category.
All these constructions and properties are inherited by $\smset^\smset$, whose limits (including products) and colimits (including coproducts) are computed pointwise: on one object at a time according to limits and colimits in $\smset$.

\index{dependent function}

For other introductions to the Yoneda lemma, the category of sets, or both, take your pick of \cite{Pierce:1991,Borceux:1994a,MacLane:1998a,Leinster:2014a,Riehl:2017a,fong2019seven,cheng_2022}.

%-------- Section --------%
\section{Exercise solutions}
\Closesolutionfile{solutions}
{\footnotesize
    \input{solution-file1}}

\Opensolutionfile{solutions}[solution-file2]

%------------ Chapter ------------%
\chapter{Polynomial functors} \label{ch.poly.obj}

In this chapter, we will formally introduce our objects of study: polynomial functors.
In addition to the set-theoretic perspective, we will present several more concrete ways to think about polynomials to aid intuition that we will use throughout the rest of this book.
We keep the mathematical content of this chapter fairly light, preferring to solidify our conceptual understanding of polynomials, before advancing to deeper categorical content.

%-------- Section --------%
\section{Introducing polynomial functors} \label{sec.poly.obj.intro}

\begin{definition}[Polynomial functor]
    A \emph{polynomial functor} (or simply \emph{polynomial}) is a functor $p\colon\smset\to\smset$ such that there exists a set $I$, an $I$-indexed family of sets $(p[i])_{i\in I}$, and an isomorphism
    \[
    p\iso\sum_{i\in I}\yon^{p[i]}
    \]
    to the corresponding $I$-indexed sum of representables.
\end{definition}

So, up to isomorphism, a polynomial functor is just a sum of representables.

\index{polynomial functor!as sum of representables}

\begin{remark}
    Given sets $I, A \in \smset$, it follows from \cref{exc.repeated_sum_is_product} that we have an isomorphism of polynomials
    \[
    \sum_{i \in I} \yon^A \iso I\yon^A.
    \]
    So when we write down a polynomial, we will often combine identical representable summands $\yon^A$ by writing them in the form $I\yon^A$.
    In particular, the constant functor $\1$ is a representable functor ($\1 \iso \yon^\0$), so every constant functor $I$ is a polynomial functor: $I \iso \sum_{i \in I} \1$.
\end{remark}

\begin{exercise}
    Consider the polynomial $q\coloneqq\yon^\8+\4\yon$.
    \begin{enumerate}
        \item Does the polynomial $q$ have a representable summand $\yon^\2$?
        \item Does the polynomial $q$ have a representable summand $\yon$?
        \item Does the polynomial $q$ have a representable summand $\4\yon$?
        \qedhere
    \end{enumerate}
    \begin{solution}
    \begin{enumerate}
        \item No, $q$ does not have $\yon^\2$ as a representable summand.
        \item Yes, $q$ does have $\yon$ as a representable summand.
        \item No, $q$ does not have $\4\yon$ as a representable summand, because $\4\yon$ is not a representable functor!
        But to make amends, we could say that $\4\yon$ is a \emph{summand}; this means that there is some $q'$ such that $q\iso q'+\4\yon$, namely $q'\coloneqq\yon^8$. So $\3\yon$ is also a summand, but $\yon^\2$ and $\5\yon$ are not.
    \end{enumerate}
    \end{solution}
\end{exercise}

\index{polynomial functor!summand of}

\begin{example}\label{ex.verbose_poly_eval}
    Consider the polynomial $p\coloneqq\yon^\2+\2\yon+\1$.
    It denotes a functor $\smset\to\smset$; where does this functor send the set $X\coloneqq\{a,b\}$?
    To be precise, we will rather verbosely say that $I\coloneqq\4$ and
    \[
    p[1]\coloneqq\2,\quad
    p[2]\coloneqq\1,\quad
    p[3]\coloneqq\1,\qqand
    p[4]\coloneqq\0\
    \]
    so that $p\cong\sum_{i\in I}\yon^{p[i]}$. Now we have
    \begin{align*}
        p(X) &\iso
        \sum_{i\in\4}\{a,b\}^{p[i]} \\ &=
        \{a,b\}^\2 + \{a,b\}^\1 + \{a,b\}^\1 + \{a,b\}^\0 \\ &\iso
        \{(1,(a,a)),(1,(a,b)),(1,(b,a)),(1,(b,b)),(2,(a)),(2,(b)), \\ &
        \qquad (3,(a)),(3,(b)),(4,())\}.
    \end{align*}
    Above, we denote each function $f\colon p[i]\to\{a,b\}$ in the set $\{a,b\}^{p[i]}$ by the $n$-tuple $(f(1),\ldots,f(n))$ whenever $p[i]\coloneqq\ord{n}$.
    For ease of reading, we may drop the parentheses around these $n$-tuples to obtain the equivalent set
    \begin{align*}
    p(X) &\iso
    \{(1,a,a),(1,a,b),(1,b,a),(1,b,b),(2,a),(2,b), \\ & \qquad (3,a),(3,b),(4)\}.
    \end{align*}
    As we might expect, the set $p(X)$ contains $2^2+2+2+1=9$ elements, equal to the value obtained when we plug $|X|=2$ into the original polynomial $p$ when we interpret its coefficients and exponents as numbers instead of sets.
\end{example}

In general, a polynomial $p\coloneqq\sum_{i\in I}\yon^{p[i]}$ applied to a set $X$ expands to
\[
\sum_{i\in I}X^{p[i]}
\]
and can be thought of as the set of all pairs comprised of an element of $I$ and a function $p[i]\to X$ or, equivalently, a $p[i]$-tuple of elements of $X$.

\index{polynomial functor!action on sets}

\begin{exercise}
    In the verbose style of \cref{ex.verbose_poly_eval}, write out all the elements of $p(X)$ for $p$ and $X$ as follows (if there are infinitely many, denote the set $p(X)$ some other way):
    \begin{enumerate}
        \item $p\coloneqq\yon^\3$ and $X\coloneqq\{4,9\}.$
        \item $p\coloneqq\3\yon^\2+\1$ and $X\coloneqq\{a\}$.
        \item $p\coloneqq\0$ and $X\coloneqq\nn$.
        \item $p\coloneqq\4$ and $X\coloneqq\nn$.
        \item $p\coloneqq\yon$ and $X\coloneqq\nn$.
        \qedhere
    \end{enumerate}
    \begin{solution}
        \begin{enumerate}
            \item Let $I\coloneqq\1$ and $p[1]\coloneqq\3$ so that $p\coloneqq\yon^\3\iso\sum_{i\in I}\yon^{p[i]}$.
            Then
            \begin{align*}
            p(X)&\iso\{(1, 4, 4, 4), (1, 4, 4, 9), (1, 4, 9, 4), (1, 4, 9, 9), \\
            &\qquad(1, 9, 4, 4), (1, 9, 4, 9), (1, 9, 9, 4), (1, 9, 9, 9)\}.
            \end{align*}

            \item Let $I\coloneqq\4$, $p[1]\coloneqq p[2]\coloneqq p[3]\coloneqq\2$, and $p[4]\coloneqq\1$, so that $p\coloneqq\3\yon^\2+\1\iso\sum_{i\in I}\yon^{p[i]}$.
            Then $p(X)\iso\{(1,a,a),(2,a,a),(3,a,a),(4)\}$.

            \item Let $I\coloneqq\0$ so that $p\coloneqq\0\iso\sum_{i\in I}\yon^{p[i]}$.
            Then $p(X)\iso\0$.
            Alternatively, note that $\0$ is the constant functor that sends every set to $\0$.

            \item Let $I\coloneqq\4$ and $p[i]\coloneqq\0$ for every $i\in I$ so that $p\coloneqq\4\iso\sum_{i\in I}\yon^{p[i]}$.
            Then $p(X)\iso\{(1), (2), (3), (4)\}\iso\4$.
            Alternatively, note that $\4$ is the constant functor that sends every set to $\4$.

            \item Let $I \coloneqq \1$ and $p[1] \coloneqq \1$ so that $p\coloneqq\yon\iso\sum_{i\in I}\yon^{p[i]}$.
            So $p(X)\iso\{(1,n)\mid n\in\nn\}\iso\nn$.
            Alternatively, note that $\yon$ is the identity functor, so it sends $\nn$ to itself.
        \end{enumerate}
    \end{solution}
\end{exercise}

The following proposition shows how the polynomial functor $p$ itself determines the set $I$ over which we sum up representables to obtain $p$.

\index{polynomial functor!$p(1)$ as summands of}

\begin{proposition}\label{prop.apply1}
    Let $p\coloneqq\sum_{i\in I}\yon^{p[i]}$ be an arbitrary polynomial functor. Then $I\cong p(\1)$, so there is an isomorphism of functors
    \begin{equation}\label{eqn.sum_p1}
        p\iso\sum_{i\in p(\1)}\yon^{p[i]}.
    \end{equation}
\end{proposition}
\begin{proof}
    We need to show that $I\iso p(\1)$; the latter claim follows directly.
    In \cref{exc.on_sums_prods_sets} \cref{exc.on_sums_prods_sets.sum}, we showed that $I\iso\sum_{i\in I}\1$, so it suffices to show that $(\yon^{p[i]})(\1)\iso\1$ for all $i \in I$.
    Indeed, $\1^{p[i]}\iso \1$ because there is a unique function $p[i]\to \1$ for each $p[i]$.
\end{proof}
We can draw an analogy between \cref{prop.apply1} and evaluating $p(1)$ for a polynomial $p$ from high school algebra, which yields the sum of the coefficients of $p$.
The notation in \eqref{eqn.sum_p1} will be how we denote arbitrary polynomials from now on, and we will use the following terms to denote the sets $p(\1)$ and $p[i]$ for $i\in p(\1)$ on which a polynomial $p$ depends.

\begin{definition}[Position and direction]
    Given a polynomial functor
    \[
    p\iso\sum_{i\in p(\1)}\yon^{p[i]},
    \]
    we call an element $i\in p(\1)$ a \emph{position of $p$} or a \emph{$p$-position}, and we call an element $a\in p[i]$ a \emph{direction of $p$ at $i$} or a \emph{$p[i]$-direction}.
    We call $p(\1)$ the \emph{position-set of $p$} and $p[i]$ the \emph{direction-set of $p$ at $i$}.
\end{definition}

\index{polynomial functor!positions and directions}
\index{positions!of polynomial|seealso{polynomial functor!positions and directions}}
\index{directions!of polynomial|seealso{polynomial functor!positions and directions}}

Note that the position-set $p(\1)$ along with the $p(\1)$-indexed family of direction-sets $p[-]\colon p(\1)\to\smset$ uniquely characterize a polynomial $p$ up to isomorphism.
Throughout this book, we will often specify a polynomial by giving its positions and its directions at each position.

\begin{exercise}\label{exc.apply0}
    We saw in \cref{prop.apply1} how to interpret the position-set $p(\1)$ of a polynomial $p$, e.g.\ $p\coloneqq\yon^\3+\3\yon^\2+\4$, as the sum of the coefficients of $p$: here $p(\1)\iso\1+\3+\4\iso\8$.
    How might you interpret $p(\0)$?
    \begin{solution}
        We consider $p(\0)$ for arbitrary polynomials $p$.
        A representable functor $\yon^S$ for $S\in\smset$ sends $\0\mapsto\0$ if $S\neq\0$ (as there are then no functions $S\to\0$), but sends $\0\mapsto\1$ if $S=\0$ (as there is a unique function $\0\to\0$).
        So
        \[
        p(\0)\iso\sum_{i\in p(\1)}\left(\yon^{p[i]}\right)(\0)\iso\sum_{\substack{i\in p(\1),\\ p[i]\neq\0}}\0+\sum_{\substack{i\in p(\1),\\ p[i]=\0}}\1\iso\{i\in p(\1)\mid p[i]=\0\}.
        \]
        That is, $p(\0)$ is the set of \emph{constant} positions of $p$, the positions of $p$ that have no directions.
        For example, if $p\coloneqq\yon^\3+\3\yon^\2+\4$, then $p(\0)=\4$.
        In the language of high school algebra, we might call $p(\0)$ the \emph{constant term} of $p$.
    \end{solution}
\end{exercise}

As a functor $\smset\to\smset$, a polynomial should act on functions as well as on sets.
Below, we explain how.

\index{polynomial functor!action on functions}

\begin{proposition} \label{prop.poly_on_functions}
    Let $p$ be an arbitrary polynomial functor, which our notation lets us write as $p\iso\sum_{i\in p(\1)} \yon^{p[i]}$, and let $f\colon X\to Y$ be an arbitrary function.
    Then $p(f)\colon p(X)\to p(Y)$ sends each $(i, g)\in p(X)$, with $i\in p(\1)$ and $g\colon p[i]\to X$, to $(i, g\then f)$ in $p(Y)$.
\end{proposition}
\begin{proof}
    For each $i\in p(\1)$, by \cref{def.representable}, the functor $\yon^{p[i]}$ sends $f$ to the function $X^{p[i]}\to Y^{p[i]}$ mapping each $g\colon p[i]\to X$ to $g\then f\colon p[i]\to Y$.
    So the sum of these functors over $i\in p(\1)$ sends each $(i, g)\in p(X)$ to $(i, g\then f)\in p(Y)$.
\end{proof}

\begin{example}
    Suppose $p\coloneqq\yon^\2+\2\yon+1$. Let $X\coloneqq\{a_1,a_2,b_1\}$ and $Y\coloneqq\{a,b,c\}$, and let $f\colon X\to Y$ be the function sending $a_1,a_2\mapsto a$ and $b_1\mapsto b$. The induced function $p(f)\colon p(X)\to p(Y)$, according to \cref{prop.poly_on_functions}, is shown below:
    \[
    \begin{tikzcd}[row sep=2pt, column sep=3pt, shorten <=-5pt, shorten >=-5pt, dashed, font=\small]
        \LMO{(1,a_1,a_1)}\ar[rrrr, bend left=25pt]&\LMO{(1,a_1,a_2)}\ar[rrr, bend left=25pt]&\LMO{(1,a_1,b_1)}\ar[rrr, bend left=25pt]&[30pt]&
        \LMO{(1,a,a)}&\LMO{(1,a,b)}&\LMO{(1,a,c)}&&
        \\
        \LMO{(1,a_2,a_1)}\ar[rrrru, bend left=10pt]&\LMO{(1,a_2,a_2)}\ar[rrru, bend left=10pt]&\LMO{(1,a_2,b_1)}\ar[rrru, bend left=10pt]&&
        \LMO{(1,b,a)}&\LMO{(1,b,b)}&\LMO{(1,b,c)}&&
        \\
        \LMO{(1,b_1,a_1)}\ar[rrrru, bend left=10pt]&\LMO{(1,b_1,a_2)}\ar[rrru, bend left=10pt]&\LMO{(1,b_1,b_1)}\ar[rrru, bend left=10pt]&&
        \LMO{(1,c,a)}&\LMO{(1,c,b)}&\LMO{(1,c,c)}&&
        \\
        \LMO{(2,a_1)}\ar[rrrr, bend left=25pt]&\LMO{(2,a_2)}\ar[rrr, bend left=25pt]&\LMO{(2,b_1)}\ar[rrr, bend left=25pt]&&
        \LMO{(2,a)}&\LMO{(2,b)}&\LMO{(2,c)}&&
        \\
        \LMO{(3,a_1)}\ar[rrrr, bend left=25pt]&\LMO{(3,a_2)}\ar[rrr, bend left=25pt]&\LMO{(3,b_1)}\ar[rrr, bend left=25pt]&&
        \LMO{(3,a)}&\LMO{(3,b)}&\LMO{(3,c)}&&
        \\&
        \LMO{(4)}\ar[rrrr]&&&&\LMO{(4)}
    \end{tikzcd}
    \]
\end{example}

\begin{exercise}
    Let $p\coloneqq\yon^\2+\yon$. Choose a function $f\colon\1\to\2$ and write out the induced function $p(f)\colon p(\1)\to p(\2)$.
    \begin{solution}
        We have
        \[
        p(\1) \iso \{(1, 1, 1), (2, 1)\} \qqand p(\2) \iso \{(1, 1, 1), (1, 1, 2), (1, 2, 1), (1, 2, 2), (2, 1), (2, 2)\}.
        \]
        Say we choose the function $f\colon\1\to\2$ that sends $1 \mapsto 1$.
        Then $p(f)$ would send $(1, 1, 1) \mapsto (1, 1, 1)$ and $(2, 1) \mapsto (2, 1)$.
        If we had instead picked $1 \mapsto 2$ as our function $f$, then $p(f)$ would send $(1, 1, 1) \mapsto (1, 2, 2)$ and $(2, 1) \mapsto (2, 2)$.
    \end{solution}
\end{exercise}

%-------- Section --------%
\section{Special classes of polynomial functors} \label{sec.poly.obj.spec}

Here we describe several special classes of polynomials.
We have already defined two such special classes: \emph{representables} and \emph{constants}.
A \emph{representable polynomial} (or simply a \emph{representable}) is a representable functor, i.e.\ a polynomial functor isomorphic to $\yon^A$ for some set $A$.
Meanwhile, a \emph{constant polynomial} (or simply a \emph{constant}) is a constant functor, i.e.\ a polynomial functor isomorphic to $I$, interpreted as a functor, for some set $I$.

\index{polynomial functor!constant|see{constant polynomial}}\index{polynomial functor!representable|see{representable polynomial}}\index{polynomial functor!linear|see{linear polynomial}}
\index{set!as constant polynomial}
\index{constant polynomial}\index{linear polynomial}\index{representable polynomial}

\begin{exercise}
\begin{enumerate}
  \item Characterize when a polynomial $p$ is \textit{representable} in terms of its positions and/or its directions.
  \item Characterize when a polynomial $p$ is \textit{constant} in terms of its positions and/or its directions. \qedhere
\end{enumerate}
\begin{solution}
\begin{enumerate}
  \item A polynomial $p$ is representable when $p\iso\yon^A$ for some set $A$, and $\yon^A$ has exactly $1$ position.
  Conversely, if a polynomial $p$ has exactly $1$ position, then $p(\1)\iso\1$, so we may write $p$ (up to isomorphism) as $p\iso\sum_{1\in\1}\yon^{p[1]}\iso\yon^{p[1]}$, which is representable.
  So a polynomial $p$ is representable if and only if it has exactly $1$ position.

  \item A polynomial $p$ is constant when $p\iso I$ for some set $I$, and $I$ has no directions at any of its positions.
  Conversely, if every direction-set of a polynomial $p$ is empty, then
  \[
    p\iso\sum_{i\in p(\1)}\yon^{p[i]}\iso\sum_{i\in p(\1)}\yon^\0\iso\sum_{i\in p(\1)}\1\iso p(\1),
  \]
  i.e.\ the set $p(\1)$ viewed as a constant functor.
  So a polynomial $p$ is constant if and only if it has exactly $0$ directions at each position.
\end{enumerate}
\end{solution}
\end{exercise}

Like constants, the other two special classes of polynomials we define here will share their names with their algebraic analogues.
Throughout, let $p\iso\sum_{i\in p(\1)}\yon^{p[i]}$ be a polynomial functor.

\begin{definition}[Linear, affine]
We say that $p$ is \emph{linear}\tablefootnote{Unlike linear polynomials from high school algebra (which are really \emph{affine linear functions} rather than necessarily \emph{linear functions}), our linear polynomial functors have no (nonzero) constant terms: they always send $\0$ to $\0$. A polynomial is \emph{affine} if it is of the form $A\yon+B$, though we will not use this concept much in the book.}
if $p\iso I\yon$ for some set $I$.
\end{definition}
\index{linear polynomial}
\index{polynomial functor!affine|see{affine polynomial}}\index{affine polynomial}

\begin{definition}[Monomial]
We say that $p$ is a \emph{monomial} if $p\iso I\yon^A$ for sets $I$ and $A$.
\end{definition}
\index{monomial}

\begin{example}\index{constant polynomial}\index{linear polynomial}\index{representable polynomial}
  Every constant polynomial $I\iso I\yon^\0$ is a monomial, as is every linear polynomial $I\yon\iso I\yon^\1$ and every representable $\yon^A \iso \1\yon^A$.
  On the other hand, there are monomials that are neither constant, linear, nor representable, such as $\2\yon^\2$ or $\nn\yon^\rr$.
  Moreover, there are polynomials that are not monomials, such as $\yon^\4+\3$ or $\sum_{n\in\nn}\yon^{\ord{n}}$.

  There is only one polynomial that is both constant and linear, namely $\0\iso\0\yon$.
  Similarly, there is only one polynomial (up to isomorphism) that is both constant and representable, namely $\1\iso\yon^\0$.
  Finally, there is only one polynomial (up to isomorphism) that is both linear and representable, namely the identity functor $\1\yon\iso\yon\iso\yon^\1$.

  In general, every set $S$ has a corresponding constant polynomial $S$, linear polynomial $S\yon$, and representable polynomial $\yon^S$; and as long as $|S|\geq2$, these are all distinct.
\end{example}

\begin{exercise}
  \begin{enumerate}
    \item Characterize when a polynomial $p$ is \textit{linear} in terms of its positions and/or its directions.
    \item Characterize when a polynomial $p$ is a \textit{monomial} in terms of its positions and/or its directions. \qedhere
  \end{enumerate}
  \begin{solution}
    \begin{enumerate}
      \item A polynomial $p$ is linear when $p\iso I\yon$ for some set $I$, and $I\yon$ has exactly $1$ direction at each position. (Note that this is even true when $p\iso\0\yon\iso\0$, for then it is true vacuously.)
      Conversely, if a polynomial $p$ has exactly $1$ direction at each position, then $p[i]\iso\1$ for all $i\in p(\1)$, so
      \[
        p
          \iso
        \sum_{i\in p(\1)}\yon^{p[i]}
          \iso
        \sum_{i\in p(\1)}\yon^\1
          \iso
        \sum_{i\in p(\1)}\yon
          \iso
        p(\1)\yon,
      \]
      which is linear.
      So a polynomial $p$ is linear if and only if it has exactly $1$ direction at each position.

      \item A polynomial $p$ is a monomial when $p\iso I\yon^A$ for sets $I$ and $A$, implying that there is an isomorphism of direction-sets $p[i]\iso A\iso p[j]$ for all $p$-positions $i$ and $j$ (i.e.\ all the direction-sets of $p$ have the same cardinality).
      Conversely, if all the direction-sets of a polynomial $p$ are isomorphic to each other, then they are all isomorphic to some set $A$, so we have
      \[
      p
      \iso
      \sum_{i\in p(\1)}\yon^{p[i]}
      \iso
      \sum_{i\in p(\1)}\yon^A
      \iso
      p(\1)\yon^A,
      \]
      which is a monomial.
      So a polynomial $p$ is a monomial if and only if all of its direction-sets have the same cardinality.
    \end{enumerate}
  \end{solution}
\end{exercise}

Later on in \cref{sec.poly.bonus.adj}, we will see how all four of these special classes of polynomials arise from various adjunctions.

%-------- Section --------%
\section[Interpreting positions and directions]{Interpreting positions and directions%
  \sectionmark{Interpreting positions \& directions}}
\sectionmark{Interpreting positions \& directions}

\index{polynomial functor!positions and directions|(}
Let us make an informal digression on how we will think about positions and directions of polynomials in this book.
While this section has little mathematical content, the intuition we build here will guide us as we delve into the deeper theory of polynomials and their applications to modeling interaction.

The main idea is that a \emph{position} is some status that may be held, while the \emph{directions} at each position are the options available when holding that status.
While these positions and directions may be imagined abstractly, here we give some concrete examples.

\begin{example}[Directions as menu options] \label{ex.reps-as-menus}
    Consider a representable and thus polynomial functor $\yon^A$ for a set $A$.
    It has $1$ position and the elements of $A$ as its directions.
    We may think of $A$ as a menu of options to choose from.

    The menu may consist of dinner options available at a wedding; then the corresponding representable functor could be
    \[
        \yon^{\{\const{chicken},\,\const{beef},\,\const{vegetarian}\}};
    \]
    or it may be the menu of a text editor, in which case the representable could be
    \[
        \yon^{\{\const{cut},\,\const{copy},\, \const{paste}\}}.
    \]
    In both these cases, there are exactly $3$ directions, so there is an isomorphism of representable functors
    \[
        \yon^{\{\const{chicken},\,\const{beef},\,\const{vegetarian}\}} \iso \yon^{\{\const{cut},\,\const{copy},\, \const{paste}\}}.
    \]

    Similarly, we may interpret the representable $\yon^\2$ as a $2$-option menu.
    Such menus are ubiquitous in life: yes or no, true or false, heads or tails, 0 or 1.
    A $1$-option menu, represented by $\yon^\1\iso\yon$, is also familiar as an unavoidable choice, the only option: ``sorry, ya just gotta go through it.''
    Having no options, represented by $\yon^\0\iso\1$, is when you actually don't get through it: an impossible decision, a ``dead end.''

    In contrast, we may interpret the representable $\yon^{[0,1]}$ as a menu with an infinite range of options: a slider with one end labeled $0$ and the other labeled $1$, able to take on any value in between.
\end{example}\index{decision!impossible}\index{isomorphism!of representable functors}

For consistency, we will favor the term ``direction'' over ``option'' when referring to the elements of $A$ for a summand $\yon^A$ of a polynomial.
Nevertheless, when we think of a polynomial's directions, we will often think of them as options to choose from.

\cref{ex.reps-as-menus} shows how we may interpret the directions of a single representable summand as options in a menu.
By having multiple representable summands---one for each position---a polynomial may capture more general scenarios with a range of possible menus.

\begin{example}[Modeling with a polynomial] \label{ex.coin-jar}
    Consider a coin jar with a slot that may be open or closed.
    When the slot is open, the jar may accept a penny, a nickel, a dime, or a quarter---there are $4$ options to choose from.
    When the slot is closed, the jar may not accept any coins at all---there are $0$ options.
    We may model this scenario with the polynomial
    \[
        \yon^{\{\const{penny},\,\const{nickel},\,\const{dime},\,\const{quarter}\}}+\yon^\0 \iso \yon^\4+\1.
    \]
    This polynomial has $2$ positions, corresponding to the two statuses the slot could take: open or closed.
    To delineate these positions, we could take advantage of the fact that every singleton set is isomorphic to $\1$ and that $\1\yon^A\iso\yon^A$ to rewrite the above polynomial as
    \[
        \{\const{open}\}\yon^{\{\const{penny},\,\const{nickel},\,\const{dime},\,\const{quarter}\}}+\{\const{closed}\}\yon^\0 \iso \yon^\4+\1.
    \]
\end{example}

\begin{exercise}
    Give another example of a real-world scenario that may be modeled by a polynomial with more than $1$ position.
\begin{solution}
    The stopwatch app on my phone has three positions: a $\const{zero}$ position, from which I may tap a single $\const{start}$ button; a $\const{running}$ position, from which I may tap either a $\const{lap}$ button or a $\const{stop}$ button; and a $\const{stopped}$ position, from which I may tap either a $\const{start}$ button or a $\const{reset}$ button.
    Thinking of the buttons available to press as the directions at each position, the corresponding polynomial is
    \[
        \{\const{zero}\}\yon^{\{\const{start}\}}+\{\const{running}\}\yon^{\{\const{lap},\,\const{stop}\}}+\{\const{stopped}\}\yon^{\{\const{start},\, \const{reset}\}}\iso\yon+\2\yon^\2.
    \]
\end{solution}
\end{exercise}

\index{polynomial functor!positions and directions|)}

%-------- Section --------%
\section{Corolla forests}

We would like to have graphical depictions of our polynomials to make them easy to visualize.
These will take the form of special graphs known as \emph{corolla forests}.
We build up to defining them as follows.

\index{corolla forest|see{polynomial functor!associated corolla forest}}\index{polynomial functor!associated corolla forest}\index{corolla forest|(}

Our first definition will be familiar to students of graph theory, although we will add some technical details suited to our purposes.
\begin{definition}[Rooted tree]\index{tree!rooted}\index{tree!corolla and}
    A \emph{rooted tree} is a directed acyclic graph with a distinguished vertex called the \emph{root} such that there exists a unique directed path from the root to each vertex.
    \index{rooted tree|see{tree, rooted}}

    We allow infinitely and even uncountably many vertices and infinitely and even uncountably many edges incident to each vertex; on the other hand, each pair of vertices is connected by a (necessarily unique) path of finitely many adjacent edges.
\end{definition}
Since all our trees will be rooted, we may refer to them simply as \emph{trees}---roots are implied.
We will draw our trees with roots at the bottom and other vertices ``growing'' upward.

The following terminology will be handy when working with our trees; these terms should be familiar, or at the very least they should match your intuition.
\begin{definition}[Rooted path; height]
    A \emph{rooted path} is a (directed) path in a rooted tree from its root to any vertex.

    \index{rooted tree!rooted path in}
    \index{rooted tree!height of vertex}

    Given a vertex of a rooted tree, its \emph{height} is the length of (i.e.\ number of edges in) the rooted path to that vertex.
\end{definition}
In any rooted tree, the root has height $0$, the length of the empty rooted path to the root itself; every neighbor of the root has height $1$; every neighbor of a vertex of height $1$ either is the root or has height $2$; and so forth.

Now we can define a special kind of tree that we will use to depict representable functors.
\begin{definition}[Corolla]
    A \emph{corolla} is a rooted tree in which every vertex aside from the root has height $1$. We call these vertices the \emph{leaves} of the corolla.

    The \emph{corolla associated to a representable functor $\yon^A$} for $A\in\smset$ is the corolla whose leaves are in bijection with $A$.
\end{definition}

\index{corolla}

\begin{example}
    Here are the corollas associated to various representables:
    
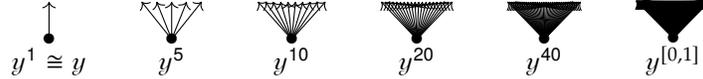
\begin{equation*}
    \begin{tikzpicture}[trees, sibling distance=2mm]
        \node["$\yon^{\1}\iso\yon$" below] (1) {$\bullet$}
        child;
    \end{tikzpicture}
    \qquad
    \begin{tikzpicture}[trees, sibling distance=2mm]
        \node["$\yon^{\5}$" below] (1) {$\bullet$}
        child foreach \i in {1,...,5}
        ;
    \end{tikzpicture}
    \qquad
    \begin{tikzpicture}[trees, sibling distance=1mm]
        \node["$\yon^{\1\0}$" below] (1) {$\bullet$}
        child foreach \i in {1,...,10}
        ;
    \end{tikzpicture}
    \qquad
    \begin{tikzpicture}[trees, sibling distance=0.5mm]
        \node["$\yon^{\2\0}$" below] (1) {$\bullet$}
        child foreach \i in {1,...,20}
        ;
    \end{tikzpicture}
    \qquad
    \begin{tikzpicture}[trees, sibling distance=0.25mm]
        \node["$\yon^{\4\0}$" below] (1) {$\bullet$}
        child foreach \i in {1,...,40}
        ;
    \end{tikzpicture}
    \qquad
    \begin{tikzpicture}[trees, sibling distance=0.0625mm]
        \node["$\yon^{[0,1]}$" below] (1) {$\bullet$}
        child foreach \i in {1,...,160}
        ;
    \end{tikzpicture}
\end{equation*}
\end{example}

In the example above, the roots are indicated by dots ($\bullet$), and the leaves are indicated by arrows ($\uparrow$).
Because the direction-set of the representable is in bijection with the leaves of the associated corolla, we can think of each leaf as a direction, so it makes sense to draw the leaves as arrows pointing in different directions.
Thinking of directions as menu options like in the previous section, we may view these corollas as mini-decision trees, indicating all the possible options we could select.
\index{corolla forest!leaf of}\index{corolla forest!corolla of}
\index{decision tree}

\begin{example}
    The corolla associated to $\yon^\0\iso\1$ has \textit{no} leaves: it is the rooted tree consisting of one vertex---the root---and no edges.
    
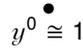
\begin{equation*}
        \begin{tikzpicture}[trees, sibling distance=2mm]
            \node["$\yon^\0\iso\1$" below] (1) {$\bullet$};
        \end{tikzpicture}
    \end{equation*}
    By definition, the root itself is \textit{not} a leaf, so the corolla above does in fact have $0$ leaves.
    With no arrows pointing out, it is the corolla associated to a representable with no directions.
\end{example}

As each representable functor has an associated corolla, each polynomial functor will have an associated disjoint union of corollas that we call a \emph{corolla forest}.

\index{corolla!representable as}

\begin{definition}[Corolla forest]
    A \emph{corolla forest} is a disjoint union of corollas.

    The \emph{corolla forest associated to a polynomial functor} $p\iso\sum_{i\in p(\1)}\yon^{p[i]}$ is the disjoint union of the corollas associated to each representable summand $\yon^{p[i]}$ of $p$.
    When we draw the corolla forest associated to $p$, we may say that we are \emph{drawing $p$ as a (corolla) forest}.
    We call the corollas in this forest corresponding to $p$-positions \emph{$p$-corollas} and the leaves corresponding to $p[i]$-directions \emph{$p[i]$-leaves}.
\end{definition}

\begin{example} \label{ex.corolla-forest}
    We may draw $p\coloneqq\yon^\2+\2\yon+\1$ as a forest like so:
    
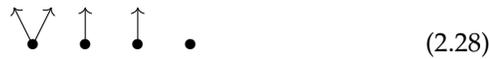
\begin{equation} \label{pic.forest-example}
    \begin{tikzpicture}[trees]
        \node (1) {$\bullet$}
        child {}
        child {};
        \node[right=.5 of 1] (2) {$\bullet$}
        child {};
        \node[right=.5 of 2] (3) {$\bullet$}
        child {};
        \node[right=.5 of 3] (4) {$\bullet$};
    \end{tikzpicture}
    \end{equation}
    Each of the $4$ corollas in \eqref{pic.forest-example} corresponds to one of the $4$ representable summands of $p$.
    The $4$ roots in \eqref{pic.forest-example} correspond to the $4$ positions of $p$, and the leaves connected to each root correspond to the directions at each position.
    Note that $p$ has $1$ position with $2$ directions, $2$ positions with $1$ direction each, and $1$ position with $0$ directions.
    Hence \eqref{pic.forest-example} is the disjoint union of $1$ corolla with $2$ leaves, $2$ corollas with $1$ leaf each, and $1$ corolla with $0$ leaves.

    Since $p(\1)\iso\4$, we could label the positions of $p$ with the elements of $\4=\{1,2,3,4\}$ so that
    \[
        p[1] = \2, \qquad p[2] = \1, \qquad p[3] = \1, \qquad p[4] = \0.
    \]
    Then we could give these same labels to the roots in \eqref{pic.forest-example}:
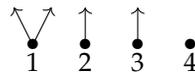
    \[
    \begin{tikzpicture}[trees]
        \node["$1$" below] (1) {$\bullet$}
        child {}
        child {};
        \node["$2$" below, right=.5 of 1] (2) {$\bullet$}
        child {};
        \node["$3$" below, right=.5 of 2] (3) {$\bullet$}
        child {};
        \node["$4$" below, right=.5 of 3] (4) {$\bullet$};
    \end{tikzpicture}
    \]
    Similarly, we could label the directions and their corresponding leaves, but we will reserve leaf labels for another purpose.
\end{example}

\begin{exercise}
    Consider the polynomial $p\coloneqq\2\yon^\3+\2\yon+\1$.
    \begin{enumerate}
        \item Draw $p$ as a corolla forest.
        \item How many roots does this forest have?
        \item How many positions of $p$ do these roots represent?
        \item For each $p$-corolla, say how many leaves it has.
        \item For each $p$-position, say how many directions it has. \qedhere
    \end{enumerate}
    \begin{solution}
        \begin{enumerate}
            \item Here is the corolla forest associated to $p\coloneqq\2\yon^\3+\2\yon+\1$ (note that the order in which the corollas are drawn does not matter):
            \[
            \begin{tikzpicture}[trees, sibling distance=3mm]
                \node (1) {$\bullet$}
                child {}
                child {}
                child {};
                \node[right=.7 of 1] (2) {$\bullet$}
                child {}
                child {}
                child {};
                \node[right=.5 of 2] (3) {$\bullet$}
                child {};
                \node[right=.3 of 3] (4) {$\bullet$}
                child {};
                \node[right=.3 of 4] (5) {$\bullet$};
            \end{tikzpicture}
            \]
            \item The forest has $5$ roots.
            \item The roots represent the $5$ positions, one position per root.
            \item \label{sol.forest.leaves} The first and second corollas have $3$ leaves each, the third and fourth corollas have $1$ leaf each, and the fifth corolla has $0$ leaves.
            \item The directions at each position correspond to the leaves in each corolla, so just copy the answer from \cref{sol.forest.leaves}, replacing ``corolla'' with ``position'' and ``leaf'' with ``direction'': the first and second positions have $3$ directions each, the third and fourth positions have $1$ direction each, and the fifth position has $0$ directions.
        \end{enumerate}
    \end{solution}
\end{exercise}

The position-set or any of the direction-sets of a polynomial may be infinite.
This makes their associated corolla forests impossible to draw precisely, but they may be approximated.
We sketch the polynomial $\yon^\3+\nn\yon^{[0,1]}$ as a forest below.
\[%\label{eqn.represented_interval}
\begin{tikzpicture}[trees, sibling distance=0.0625mm]
\node (1) {$\bullet$}
child[sibling distance=3mm] foreach \i in {1,2,3}
;
\node[right=1 of 1] (2) {$\bullet$}
child foreach \i in {1,...,160}
;
\node[right=1 of 2] (3) {$\bullet$}
child foreach \i in {1,...,160}
;
\node[right=1 of 3] (4) {$\bullet$}
child foreach \i in {1,...,160}
;
\node[right=.7 of 4] (5) {$\cdots$};
\end{tikzpicture}
\]

\begin{exercise}
If you were a suitor choosing the corolla forest you love, aesthetically speaking, which would strike your interest? Answer by selecting the associated polynomial:
\begin{enumerate}
    \item $\yon^\2+\yon+\1$
    \item $\yon^\3+\3\yon^\2+\3\yon+\1$
    \item $\yon^\2$
    \item $\yon+\1$
    \item $(\nn\yon)^\nn$
    \item $S\yon^S$ for some set $S$
    \item $\yon^{\1\0\0}+\yon^\2+\3\yon$
    \item $\yon + \2\yon^\4 + \3\yon^\9 + \4\yon^{\1\6} + \cdots$
    \item Your polynomial's name $p$ here.
\end{enumerate}
Any reason for your choice? Draw a sketch of your forest.
\begin{solution}
    Aesthetically speaking, here is a polynomial that may be drawn as a beautiful corolla forest:
    \[
        p\coloneqq\yon^\0+\yon^\1+\yon^\2+\yon^\3+\cdots
    \]
    It is reminiscent (and formally related) to the notion of lists: if $A$ is any set, then $p(A)\iso A^\0+A^\1+A^\2+\cdots$ is the set $\lst(A)$ of lists (i.e.\ finite ordered sequences) with entries in $A$.
    Here is a picture of the lovely forest associated to $p$:
    \[
    \begin{tikzpicture}[trees, sibling distance=3mm]
        \node (1) {$\bullet$};
        \node[right=.3 of 1] (2) {$\bullet$}
        child {};
        \node[right=.4 of 2] (3) {$\bullet$}
        child {}
        child {};
        \node[right=.6 of 3] (4) {$\bullet$}
        child {}
        child {}
        child {};
        \node[right=.6 of 4] {$\cdots$};
    \end{tikzpicture}
    \]
\end{solution}
\end{exercise}

Corolla forests help us visualize the positions and directions of polynomials, and they will especially come in handy in the next chapter, when we describe the morphisms between our polynomials and how they interact with positions and directions.
They may also depict the elements of a polynomial functor applied to a given set, as follows.
We have seen that for a polynomial $p$ and a set $X$,
\[
    p(X) \iso \sum_{i\in p(\1)}X^{p[i]} \iso \{(i,f)\mid i\in p(\1), f\colon p[i]\to X\}.
\]
So an element of $p(X)$ is a $p$-position $i$ along with a function $f$ that maps each direction at $i$ to an element of $X$.
Equivalently, it is a $p$-corolla along with a function that maps each of its leaves to an element of $X$.
Then to draw an element $(i,f)\in p(X)$, we simply need to draw the $p$-corolla corresponding to $i$ and label its leaves with elements of $X$ according to $f$.

\begin{example} \label{ex.corolla-apply-poly}
    In \cref{ex.corolla-forest}, we drew $p\coloneqq\yon^\2+\2\yon+\1$ as a corolla forest like so:
    
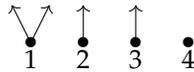
\[
    \begin{tikzpicture}[trees]
        \node["$1$" below] (1) {$\bullet$}
        child {}
        child {};
        \node["$2$" below, right=.5 of 1] (2) {$\bullet$}
        child {};
        \node["$3$" below, right=.5 of 2] (3) {$\bullet$}
        child {};
        \node["$4$" below, right=.5 of 3] (4) {$\bullet$};
    \end{tikzpicture}
    \]
    Previously, in \cref{ex.verbose_poly_eval}, we wrote out all $9$ elements of $p$ applied to the set $X\coloneqq\{a,b\}$ as tuples.
    We could draw them out instead---an element of $p(X)$ may be depicted as one of the four corollas above with each of its leaves labeled with an element of $X$:
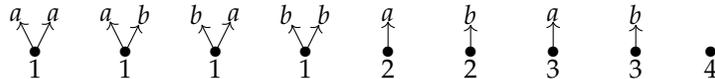
    \[
    \begin{tikzpicture}[trees, sibling distance=5mm]
        \node["$1$" below] (1) {$\bullet$}
            child {node {$a$}}
            child {node {$a$}};
        \node["$1$" below, right=1 of 1] (2) {$\bullet$}
            child {node {$a$}}
            child {node {$b$}};
        \node["$1$" below, right=1 of 2] (3) {$\bullet$}
            child {node {$b$}}
            child {node {$a$}};
        \node["$1$" below, right=1 of 3] (4) {$\bullet$}
            child {node {$b$}}
            child {node {$b$}};
        \node["$2$" below, right=.9 of 4] (5) {$\bullet$}
            child {node {$a$}};
        \node["$2$" below, right=.9 of 5] (6) {$\bullet$}
            child {node {$b$}};
        \node["$3$" below, right=.9 of 6] (7) {$\bullet$}
            child {node {$a$}};
        \node["$3$" below, right=.9 of 7] (8) {$\bullet$}
            child {node {$b$}};
        \node["$4$" below, right=.8 of 8] (9) {$\bullet$};
    \end{tikzpicture}
    \]
\end{example}

Throughout this book, we will generally use corolla forests to depict polynomials with relatively small numbers of positions or directions, where drawing out entire corolla forests is manageable.
Later, we will study how building larger rooted trees out of these corollas corresponds to conducting various categorical operations on our polynomials.

\index{corolla forest|)}

%% TODO: how much later?

%-------- Section --------%
\section{Polyboxes}

Before we conclude this chapter, we introduce one more tool for visualizing polynomials whose full power will not be evident until later.

\index{polybox|(}

Throughout this book, we may depict a polynomial $p$ as a pair of boxes stacked on top of each other, like so:
\[
  \begin{tikzpicture}[polybox, tos]
    \node[poly, "$p$" below] (p) {};
    \node[left=0pt of p_pos] {$p(\1)$};
    \node[left=0pt of p_dir] {$p[-]$};
  \end{tikzpicture}
\]
We call this picture the \emph{polyboxes for $p$}.
Think of these boxes as cells in a spreadsheet.
The bottom cell, or the \emph{position box}, is restricted to values in the set $p(\1)$ (as indicated by the label to its left)---it must be filled with a $p$-position, say $i\in p(\1)$:
\[
  \begin{tikzpicture}[polybox, tos]
    \node[poly, "$p$" below] (p) {\at$i$};
    \node[left=0pt of p_pos] {$p(\1)$};
    \node[left=0pt of p_dir] {$p[-]$};
  \end{tikzpicture}
\]
The top cell, or the \emph{direction box}, cannot be filled until the position box below it is.
Once the position box contains a $p$-position $i$, the direction box must be filled with a $p[i]$-direction, say $a\in p[i]$:
\[
  \begin{tikzpicture}[polybox, tos]
    \node[poly, "$p$" below] (p) {$a$\at$i$};
    \node[left=0pt of p_pos] {$p(\1)$};
    \node[left=0pt of p_dir] {$p[-]$};
  \end{tikzpicture}
\]
The $p[-]$ label to the left of the direction box reminds us that the $a$ within it is an element of $p[i]$, where $i$ is the entry in the position box.
Once we are accustomed to polyboxes, we will often drop these reminder labels, so that
\[
  \begin{tikzpicture}[polybox, tos]
    \node[poly, "$p$" below] (p) {$a$\at$i$};
  \end{tikzpicture}
\]
serves as a graphical shorthand for the statement ``consider a polynomial functor $p$ with position $i\in p(\1)$ and direction $a\in p[i]$.''

\index{polybox!spreadsheet}\index{polybox!position box}\index{polybox!direction box}

Viewing polynomials as these restricted two-cell spreadsheets reinforces the idea that directions are like menu options: imagine a dropdown menu for the direction box above a filled position box that lists all the directions to choose from at the given position.
Polyboxes also help us conceptualize the possible pairs of positions and directions of a polynomial whose corolla forest is impractical to draw, as suggested by the following example.

\begin{example}
  Consider the polynomial
  \[
    p\coloneqq\sum_{r\in\rr}\yon^{[-|r|,|r|]},
  \]
  whose positions are the real numbers and whose directions at position $r$ are the real numbers with magnitude at most $|r|$.
  There is no clear way to draw $p$ as a corolla forest, but we could draw its polyboxes
  \[
  \begin{tikzpicture}[polybox, tos]
    \node[poly, "$p$" below] (p) {$s$\at$r$};
    \node[left=0pt of p_pos] {$p(\1)$};
    \node[left=0pt of p_dir] {$p[-]$};
  \end{tikzpicture}
  \]
  with the condition that $r$ and $s$ are real numbers satisfying $|s|\leq|r|$.
\end{example}

We may also use polyboxes to highlight our special classes of polynomials.
When a position box may only be filled with one possible entry, we shade it in like so:
\[
\begin{tikzpicture}[polybox, tos]
  \node[poly, "$p\iso\yon^A$" below, pure] (p) {};
  \node[left=0pt of p_pos] {$p(\1)\iso\1$};
  \node[left=0pt of p_dir] {$A$};
\end{tikzpicture}
\]
The idea is that if there is only one entry that could fill a given box, then it should come pre-filled---no further choice needs to be made to fill it.
Here $p(\1)\iso\1$, so $p$ is \textit{representable}; indeed, $p\iso\yon^A$, where $A$ is the set of possible entries for the unfilled direction box.

Similarly, the polyboxes for a \textit{linear} polynomial $I\yon$, whose direction-set at each position is a singleton, can be drawn like so:
\[
\begin{tikzpicture}[polybox, tos]
  \node[poly, "$I\yon$" below, linear] (p) {};
  \node[left=0pt of p_pos] {$I$};
\end{tikzpicture}
\]
No matter what fills the position box, there is exactly one entry that could fill the direction box, so it comes pre-filled.
The identity polynomial functor $\yon$, which is both representable and linear, therefore has the following polyboxes:
\[
\begin{tikzpicture}[polybox, tos]
  \node[poly, "$\yon$" below, identity] (p) {};
\end{tikzpicture}
\]
It has exactly one position and exactly one direction, so both its boxes come pre-filled.

\index{polybox!for linear polynomials}\index{polybox!for representable polynomials}

Finally, a \textit{constant} polynomial $I$ for some set $I$ has empty direction-sets.
We indicate this by coloring its direction box red:
\[
\begin{tikzpicture}[polybox, tos]
  \node[poly, "$I$" below, constant] (p) {};
  \node[left=0pt of p_pos] {$I$};
\end{tikzpicture}
\]
Because every direction-set is empty, there is nothing that may be written in the direction box.
The red suggests an error---the direction box cannot be filled.
The polynomial functor $\1$, which is both representable and constant, therefore has the following polyboxes:
\[
\begin{tikzpicture}[polybox, tos]
  \node[poly, "$\1$" below, terminal] (p) {};
\end{tikzpicture}
\]

\index{polybox!for constant polynomials}\index{constant polynomial}

In the next chapter, we will introduce the morphisms between polynomial functors and see how their behavior may be depicted using polyboxes.

\index{polybox|)}

%-------- Section --------%
\section[Summary and further reading]{Summary and further reading%
  \sectionmark{Summary \& further reading}}
\sectionmark{Summary \& further reading}

In this chapter, we introduced the main objects of study in this book: \emph{polynomial functors}, which are sums of representable functors $\smset\to\smset$.
We write a polynomial $p$ as $\sum_{i\in p(\1)}\yon^{p[i]}$, calling the elements of $p(\1)$ the \emph{positions} of $p$ and the elements of $p[i]$ the \emph{directions} of $p$ at position $i$.
As a polynomial is determined up to isomorphism by its position-set and direction-sets, we can think of the data of a polynomial as an indexed family of sets $(p[i])_{i\in p(\1)}$.

We highlighted four special classes of polynomials (here $I$ and $A$ are sets):
\begin{itemize}
  \item constants $I$, whose direction-sets are all empty;
  \item linear polynomials $I\yon$, whose direction-sets are all singletons;
  \item representables $\yon^A$, whose position-sets are singletons;
  \item monomials $I\yon^A$, whose direction-sets all have the same cardinality.
\end{itemize}

\index{constant polynomial}\index{linear polynomial}\index{representable polynomial}\index{polynomial functor!monomial|see{monomial}}\index{monomial}

Throughout this book, we will use polynomials to model decision-making agents that hold positions and take directions from those positions.
We can draw a polynomial $p$ graphically as a \emph{corolla forest}, with a \emph{corolla} (a rooted tree whose non-root vertices are all leaves) for every $p$-position $i$ that has a leaf for every $p[i]$-direction.
We can also depict a polynomial as a \emph{polybox picture}, resembling two stacked cells in a spreadsheet, to be filled in with an element of $i$ below and an element of $p[i]$ above.

\index{corolla forest}\index{polybox}

There are many fine sources on polynomial functors. Some of the computer science literature is more relaxed about what a polynomial is. For example, the ``coalgebra community'' often defines a polynomial to include finite power sets (see e.g.\ \cite{jacobs2017introduction}). Other computer science communities use the same definition of polynomial, but refer to it as a \emph{container} and use different words for its positions (they call them ``shapes'') and directions (they call them, rather unfortunately, ``positions''). See e.g.\ \cite{abbot2003categoriesthesis,abbott2005containers}.\index{coalgebra!community}

\index{polynomial functor!term usage by ``coalgebra community''}

But the notion of polynomial functors seems to have originated from Andr\'{e} Joyal. A good introduction to polynomial functors---including an extensive bibliography of references---can be found in \cite{kock2012polynomial} and more extensive notes in \cite{kock2016}; in particular the related work section on page~3 provides a nice survey of the field. A reader may also be interested in the Workshops on Polynomial Functors organized by the Topos Institute: \url{https://topos.site/p-func-workshop/}.

\index{Joyal, Andr\'{e}}\index{Kock, Joachim}

%-------- Section --------%
\section{Exercise solutions}
\Closesolutionfile{solutions}
{\footnotesize
    \input{solution-file2}}

\Opensolutionfile{solutions}[solution-file3]

%------------ Chapter ------------%
\chapter{The category of polynomial functors} \label{ch.poly.cat}
\index{category!of polynomial functors|(}

% TODO: Should we name it "Morphisms between polynomial functors" or maybe "Dependent lenses" instead?? If so, split off symmetric monoidal structures except little bit about coproducts

In this chapter, we will define $\poly$, our main category of interest, so that we have a firm foundation from which to speak about interactive systems.
The objects of $\poly$ are the polynomial functors that we defined in the previous chapter.
Here we will examine the morphisms of $\poly$: natural transformations between polynomial functors.
Along the way, we will present some of $\poly$'s most versatile categorical properties.

\index{natural transformation!between polynomials|see{lens}}

%-------- Section --------%
\section[Dependent lenses between polynomial functors]{Dependent lenses between polynomial functors%
  \sectionmark{Dependent lenses}}
\sectionmark{Dependent lenses}
\label{sec.poly.cat.morph}

Before we define the category $\poly$ of polynomial functors, we note that polynomial functors live inside a category already: the category $\smset^\smset$ of functors $\smset\to\smset$, whose morphisms are natural transformations.
This leads to a natural definition of morphisms between polynomial functors, from which we can derive a category of polynomial functors for free.
We call such a morphism a \emph{dependent lens}, or a \emph{lens} for short.
If you are familiar with lenses from functional programming, we'll see in \cref{subsec.poly.cat.morph.bimorphic-lens} how our notion of a dependent lens is related.

\index{functional programming!lenses in}
\index{polynomial functor!morphism of|see{lens}}
\index{dependent lens|see{lens}}
\index{lens|(}

\begin{definition}[Dependent lens, $\poly$] \label{def.poly_cat}
Given polynomial functors $p$ and $q$, a \emph{dependent lens} (or simply \emph{lens}) \emph{from $p$ to $q$} is a natural transformation $p\to q$.
Then $\poly$ is the category whose objects are polynomial functors and whose morphisms are dependent lenses.
\end{definition}

In other words, $\poly$ is the full subcategory of $\smset^\smset$ spanned by the polynomial functors: we take the category $\smset^\smset$, throw out all the objects that are not (isomorphic to) polynomials, but keep all the same morphisms between the objects that remain.

\index{polynomial functor!category of polynomials}
\index{polynomial functor!full subcategory of $\smset^\smset$}

Unraveling the familiar definition of a natural transformation, a dependent lens between polynomial functors $p \to q$ thus consists of a function $p(X) \to q(X)$ for every set $X$ such that naturality squares commute.
That is a lot of data to keep track of!
Fortunately, there is a much simpler way to think about these lenses, which we will discover using the Yoneda lemma.

\index{Yoneda lemma}\index{isomorphism!natural}

\begin{exercise} \label{exc.poly_morph_yoneda}
Given a set $S$ and a polynomial $q$, show that a lens $\yon^S \to q$ can be naturally identified with an element of the set $q(S)$.
That is, show that there is an isomorphism
\[
    \poly(\yon^S, q) \iso q(S).
\]
natural in both $S$ and $q$.
Hint: Use the Yoneda lemma (\cref{lemma.yoneda}).
\begin{solution}
We know that $\poly$ is the full subcategory of $\smset^\smset$ spanned by polynomial functors, including $\yon^S$ and $q$.
So $\poly(\yon^S, q)=\smset^\smset(\yon^S, q)$.
Hence the natural isomorphism $\poly(\yon^S, q)\iso q(S)$ follows directly from the Yoneda lemma (\cref{lemma.yoneda}) with $F\coloneqq q$.
\end{solution}
\end{exercise}

The above exercise gives us an alternative characterization for lenses out of representable functors.
But before we can characterize lenses out of polynomial functors in general, we need to describe how coproducts work in $\poly$.
Fortunately, since polynomial functors are defined as coproducts of representables, coproducts in $\poly$ are easy to understand.

\index{coproduct}\index{polynomial functor!coproduct of polynomials|see{polynomial functors, sum of polynomials}}
\index{polynomial functor!product of polynomials|see{polynomial functors, product of polynomials}}

\begin{proposition} \label{prop.poly_coprods}
  The category $\poly$ has all small coproducts, coinciding with coproducts in $\smset^\smset$ given by the operation $\sum_{i \in I}$ for each set $I$.
\end{proposition}
\begin{proof}
  By \cref{cor.sum_prod_set_endofuncs}, the category $\smset^\smset$ has all small coproducts given by $\sum_{i\in I}$.
  The full subcategory inclusion $\poly\to \smset^\smset$ reflects these coproducts, and by definition $\poly$ is closed under the operation $\sum_{i \in I}$.
\end{proof}

Explicitly, given an $I$-indexed family of polynomials $(p_i)_{i \in I}$, its coproduct is
\begin{equation} \label{eqn.poly_coprod}\index{coproduct!of indexed family}
  \sum_{i \in I} p_i \iso \sum_{i \in I} \sum_{j \in p_i(\1)} \yon^{p_i[j]} \iso \sum_{(i,\,j) \in \sum_{i \in I} p_i(\1)} \yon^{p_i[j]}
\end{equation}
by \cref{cor.sum_prod_set_endofuncs}.
This coincides with our notion of polynomial addition from high school algebra: just add all the terms together, combining like terms to simplify.
Binary coproducts are given by binary sums of functors, appropriately denoted by $+$, while the initial object of $\poly$ is the constant polynomial $\0$.

In particular, \eqref{eqn.poly_coprod} implies that for any polynomials $p$ and $q$, their coproduct $p+q$ is given as follows.
The position-set of $p+q$ is the coproduct of sets $p(\1) + q(\1)$.
At position $(1,i) \in p(\1) + q(\1)$ with $i \in p(\1)$, the directions of $p+q$ are just the $p[i]$-directions; at position $(2,j) \in p(\1) + q(\1)$ with $j \in q(\1)$, the directions of $p+q$ are just the $q[j]$-directions.

Crucially, we have the following corollary.

\begin{corollary} \label{cor.poly-coprod-repr}
  In the category $\poly$, every polynomial $p$ is the coproduct of its representable summands $(\yon^{p[i]})_{i\in p(\1)}$.
\end{corollary}

In other words, writing $p$ as the sum $\sum_{i\in p(\1)}\yon^{p[i]}$ is not just a coproduct in $\smset^\smset$; it is also a coproduct in $\poly$ itself.\index{coproduct!in $\poly$}

We are now ready to give our alternative characterization of dependent lenses.
Recall that a polynomial $p\iso\sum_{i\in p(\1)}\yon^{p[i]}$ can be uniquely identified with an indexed family $p[-]\colon p(\1)\to\smset$, a functor from the set $p(\1)$ viewed as a discrete category.\index{category!discrete}

\begin{proposition}\label{prop.lens-prod-sum}
Given polynomials $p$ and $q$, there is an isomorphism
\begin{equation}\label{eqn.main_formula}
\poly(p,q)\cong\prod_{i\in p(\1)}\sum_{j\in q(\1)}{p[i]}^{q[j]}
\end{equation}
natural in $p$ and $q$.
In particular, a lens $f\colon p\to q$ can be identified with a pair $(f_\1,f^\sharp)$
\begin{equation}\label{eqn.colax_poly_map}
\begin{tikzcd}[column sep=small]
	p(\1)\ar[dr, "p{[-]}"']\ar[rr, "f_\1"]&~&
	q(\1)\ar[dl, "q{[-]}"]\\&
	\smset\ar[u, phantom, near end, "\overset{f^\sharp}{\Leftarrow}"]
\end{tikzcd}
\end{equation}
where $f_\1\colon p(\1)\to q(\1)$ is a function (equivalently, a functor between discrete categories) and $f^\sharp \colon q[f_\1(-)] \to p[-]$ is a natural transformation: a function $f^\sharp_i\colon q[f_\1i]\to p[i]$ for each $i\in p(\1)$.
\end{proposition}
\begin{proof}
We have $p\iso\sum_{i\in p(\1)}\yon^{p[i]}$.
Then by \cref{cor.poly-coprod-repr} and the universal property of the coproduct, we have a natural isomorphism\index{coproduct!universal property of}\index{isomorphism!natural}
\[
    \poly\left(\sum_{i\in p(\1)}\yon^{p[i]}, q\right) \iso \prod_{i\in p(\1)}\poly(\yon^{p[i]},q).
\]
Applying \cref{exc.poly_morph_yoneda} (i.e.\ the Yoneda lemma) and the fact that $q\iso\sum_{j\in q(\1)}\yon^{q[j]}$ yields the natural isomorphism
\[
  \prod_{i\in p(\1)}\poly(\yon^{p[i]},q) \iso \prod_{i\in p(\1)}q(p[i]) \iso \prod_{i\in p(\1)}\sum_{j\in q(\1)}p[i]^{q[j]},
\]
so \eqref{eqn.main_formula} follows.

\index{lens!positions and directions}\index{Yoneda lemma}

The right hand side of \eqref{eqn.main_formula} is the set of dependent functions $f\colon(i\in p(\1))\to\sum_{j\in q(\1)}p[i]^{q[j]}$.
Each such dependent function is uniquely determined by its two projections $\pi_1f\colon(i\in p(\1))\to q(\1)$ and $\pi_2f\colon(i\in p(\1))\to p[i]^{q[\pi_1fi]}$.
These can be identified respectively with a (non-dependent) function $f_\1\coloneqq\pi_1f$ with signature $p(\1)\to q(\1)$ and a natural transformation $f^\sharp\colon q[f_\1(-)]\to p[-]$ whose $i$-component for $i\in p(\1)$ is $f^\sharp_i\coloneqq\pi_2 fi\in p[i]^{q[f_\1i]}$.
\end{proof}

We have now greatly simplified our characterization of a dependent lens $f\colon p\to q$: rather than infinitely many functions satisfying infinitely many naturality conditions, $f$ may simply be specified by a function $f_\1\colon p(\1)\to q(\1)$ and, for each $i\in p(\1)$, a function $f^\sharp_i\colon q[f_\1i]\to p[i]$, without any additional restrictions.
This characterization can be expressed entirely in the language of positions and directions: $f_\1$ is a function from $p$-positions to $q$-positions, while $f^\sharp_i$ for a $p$-position $i$ is a function from $q[f_\1i]$-directions to $p[i]$-directions.
This leads to the following definition.

\begin{definition}[On-positions function, on-directions map and function]
  Given a lens $f\colon p\to q$, let $(f_\1, f^\sharp)$ denote the pair identified with $f$ via \cref{prop.lens-prod-sum}.
  Then we call the function $f_\1\colon p(\1)\to q(\1)$ the \emph{(forward) on-positions function of $f$}, while we call the natural transformation $f^\sharp\colon q[f_\1(-)]\to p[-]$ the \emph{(backward) on-directions map of $f$}.
  For $i\in p(\1)$, we call the $i$-component $f^\sharp_i\colon q[f_\1i]\to p[i]$ of $f^\sharp$ the \emph{(backward) on-directions function of $f$ at $i$}.
\end{definition}

The above definition highlights the bidirectional nature of a lens $f\colon p\to q$: it consists of a function going \textit{forward} on positions, following the direction of $f$ from $p$ to $q$, as well as functions going \textit{backward} on directions, opposing the direction of $f$ from $q$ to $p$.
This forward-backward interaction is what drives the applications of $\poly$ we will study.

\index{lens!bidirectionality of}

We prefer to call a morphism between polynomial functors a ``lens'' rather than a ``natural transformation'' because we wish to emphasize this concrete on-positions and on-directions perspective.
Whenever we do need to view a morphism in $\poly$ as a natural transformation, we will refer to them as such.

In the next several sections, we will give some examples of lenses and intuition for thinking about them in terms of interaction protocols, corolla forests, and polyboxes.

%-------- Section --------%
\section{Dependent lenses as interaction protocols}

Here is our first example of a dependent lens and a real-world interaction it might model.

\begin{example}[Modeling an interaction protocol with a lens]
  Recall our coin jar polynomial from \cref{ex.coin-jar}:
  \[
    q\coloneqq\{\const{open}\}\yon^{\{\const{penny},\,\const{nickel},\,\const{dime},\,\const{quarter}\}}+\{\const{closed}\}\yon^\0.
  \]
  It has $2$ positions: its $\const{open}$ position has $4$ directions representing the $4$ denominations of coins it may take, while its $\const{closed}$ position has $0$ directions to indicate that it cannot take anything.

  \index{lens!as interaction protocol|(}

  Now imagine that we model the owner of this coin jar with the following polynomial:
  \begin{align*}
    p\coloneqq\:
    &\{\const{needy}\}\yon^{\{\const{save},\, \const{spend}\}}
      \\
    +\:
    &\{\const{greedy}\}\yon^{\{\const{accept},\, \const{reject},\,\const{ask for more}\}}
      \\
    +\:
    &\{\const{content}\}\yon^{\{\const{count},\,\const{rest}\}}.
  \end{align*}
  Each of its $3$ positions represents a possible mood of the owner, and the directions at each position represent the options available to an owner in the corresponding mood.
  We will construct a lens $f\colon p\to q$ to model the interaction between the owner and their coin jar.

  Say that a $\const{needy}$ or $\const{greedy}$ owner will keep their coin jar $\const{open}$, while a $\const{content}$ owner will keep their coin jar $\const{closed}$.
  We can express this with an on-positions function $f_\1$ from the set of $p$-positions (on the left) to the set of $q$-positions (on the right), as follows (the dashed arrows indicate the function assigments):
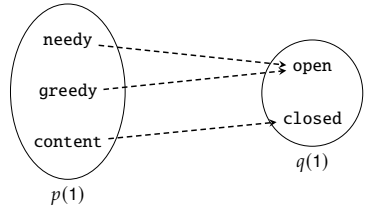
  \begin{center}
  \scalebox{.7}{
  \begin{tikzpicture}
    % q-positions (right)
      % from top to bottom: open, closed
    \node (open) {\const{open}};
    \node[below=.5 of open] (closed) {\const{closed}};

    \node[draw, ellipse, inner sep=0pt, fit=(open)(closed),
      label={[anchor=south,below]270:$q(\1)$}] (qpos) {};

    % p-positions (left)
      % from top to bottom: needy, greedy, content
      % greedy aligned with middle of q-positions
    \node[left=3 of qpos] (greedy) {\const{greedy}};
    \node[above=.5 of greedy] (needy) {\const{needy}};
    \node[below=.5 of greedy] (content) {\const{content}};

    \node[draw, ellipse, inner sep=0pt, fit=(greedy)(needy)(content),
      label={[anchor=south,below]270:$p(\1)$}] (ppos) {};

    \draw[mapsto] (needy) -- (open);
    \draw[mapsto] (greedy) -- (open);
    \draw[mapsto] (content) -- (closed);
  \end{tikzpicture}
  }
  \end{center}

  From there, say that a $\const{needy}$ owner whose coin jar receives a $\const{nickel}$ or higher will choose to $\const{save}$ it, but one whose coin jar receives a $\const{penny}$ will choose to $\const{spend}$ it.
  Meanwhile, a $\const{greedy}$ owner whose coin jar receives a $\const{penny}$ or $\const{nickel}$ will ask for more, but one whose coin jar receives a $\const{dime}$ or $\const{quarter}$ will accept it.
  We can express this behavior with an on-directions map $f^\sharp\colon q[f_\1(-)]\to p[-]$.
  Its $\const{needy}$-component is the on-directions function $f^\sharp_{\const{needy}}\colon q[f_\1(\const{needy})]\to p[\const{needy}]$ drawn as follows:
  \begin{center}
  \scalebox{.7}
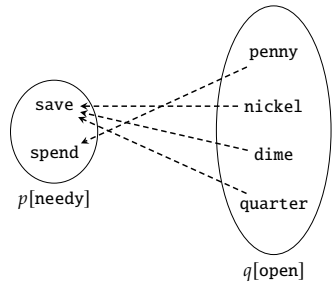
{
  \begin{tikzpicture}
    % p-directions (left)
      % from top to bottom: save, spend
    \node (save) {\const{save}};
    \node[below=.5 of save] (spend) {\const{spend}};

    \node[draw, ellipse, inner sep=0pt, fit=(save)(spend),
      label={[anchor=south,below]270:$p[\const{needy}]$}] (pdir) {};

    % q-directions (right)
      % from top to bottom: penny, nickel, dime, quarter
      % nickel aligned with save
    \node[right=3 of save] (nickel) {\const{nickel}};
    \node[above=.5 of nickel] (penny) {\const{penny}};
    \node[below=.5 of nickel] (dime) {\const{dime}};
    \node[below=.5 of dime] (quarter) {\const{quarter}};

    \node[draw, ellipse, inner sep=0pt,
      fit=(penny)(nickel)(dime)(quarter),
      label={[anchor=south,below]270:$q[\const{open}]$}] (qdir) {};

    \draw[mapsto] (penny) -- (spend);
    \draw[mapsto] (nickel) -- (save);
    \draw[mapsto] (dime) -- (save);
    \draw[mapsto] (quarter) -- (save);
  \end{tikzpicture}
  }
  \end{center}
  Notice that, by keeping the positions and directions of $p$ on the left and those of $q$ on the right, the on-positions function is drawn from left to right, while the on-directions functions must be drawn right to left.
  The $\const{greedy}$-component of $f^\sharp$ is the on-directions function $f^\sharp_{\const{greedy}}\colon q[f_\1(\const{greedy})]\to p[\const{greedy}]$ drawn like so:
  \begin{center}
  \scalebox{.7}
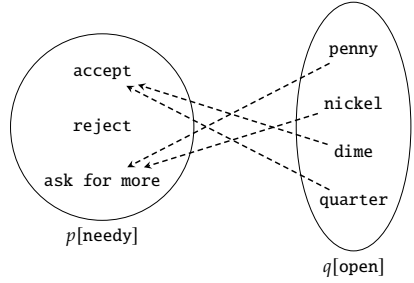
{
    \begin{tikzpicture}
      % q-directions (right)
        % from top to bottom: penny, nickel, dime, quarter
      \node (nickel) {\const{nickel}};
      \node[above=.5 of nickel] (penny) {\const{penny}};
      \node[below=.5 of nickel] (dime) {\const{dime}};
      \node[below=.5 of dime] (quarter) {\const{quarter}};

      \node[draw, ellipse, inner sep=0pt,
        fit=(penny)(nickel)(dime)(quarter),
        label={[anchor=south,below]270:$q[\const{open}]$}] (qdir) {};

      % p-directions (left)
        % from top to bottom: accept, reject, ask for more
        % reject aligned with middle of q-directions
      \node[left=3 of qdir] (reject) {\const{reject}};
      \node[above=.5 of reject] (accept) {\const{accept}};
      \node[below=.5 of reject] (ask) {\const{ask for more}};

      \node[draw, ellipse, inner sep=0pt, fit=(accept)(reject)(ask),
        label={[anchor=south,below]270:$p[\const{needy}]$}] (pdir) {};

      \draw[mapsto] (penny) -- (ask);
      \draw[mapsto] (nickel) -- (ask);
      \draw[mapsto] (dime) -- (accept);
      \draw[mapsto] (quarter) -- (accept);
    \end{tikzpicture}
  }
  \end{center}
  Finally, as $q[f_\1(\const{content})]=q[\const{closed}]=\0$, the $\const{content}$-component of $f^\sharp$ is the vacuously-defined on-directions function $f^\sharp_{\const{content}}\colon \0\to p[\const{content}]$.
  Together, the on-positions function $f_\1$ and the on-directions map $f^\sharp$ defined above completely characterize a lens $f\colon p\to q$ depicting the interaction between the coin jar and its owner.
\end{example}

\index{lens!bidirectionality of}

More generally, a lens depicts what we call an \emph{interaction protocol}, a kind of dialogue between two agents regarding their positions and directions.
Say that one agent is represented by a polynomial $p$ and another by a polynomial $q$.
Then a lens $f\colon p\to q$ is an interaction protocol that prescribes how the positions of $p$ influence the positions of $q$ and how the directions of $q$ influence the directions of $p$.
Each $p$-position $i\in p(\1)$ is passed forward via the on-positions function of $f$ to a $q$-position $f_\1i\in q(\1)$.
Then each $q[f_\1i]$-direction $b$ is passed backward via the on-directions function of $f$ at $i$ to a $p[i]$-direction $f^\sharp_ib$.

To visualize these lenses, we may use either our corolla forests or our polyboxes.

\index{lens!as interaction protocol|)}

%-------- Section --------%
\section{Corolla forest pictures of dependent lenses}

The corolla forest associated to a polynomial concretely depicts its positions and directions, making it easy to extend our corolla forest pictures to depict the dependencies between the positions and directions of two polynomials that a lens between them prescribes.

\index{lens!corolla forest depiction}

\begin{example}\label{ex.practice_with_poly_morphisms}
Let $p\coloneqq \yon^\3+\2\yon$ and $q\coloneqq\yon^\4+\yon^\2+\2$.
We draw them as corolla forests with their positions labeled (we use different shapes to distinguish the positions of different polynomials):
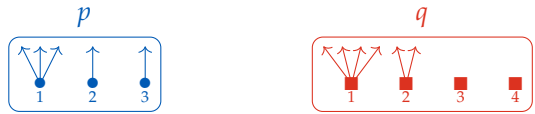
\[
\begin{tikzpicture}[rounded corners]
	\node (p1) [draw, my-blue, "\color{my-blue} $p$" above] {
	\begin{tikzpicture}[trees, sibling distance=2.5mm]
    \node["\tiny 1" below] (1) {$\bullet$}
      child {}
      child {}
      child {};
    \node[right=.5 of 1,"\tiny 2" below] (2) {$\bullet$}
      child {};
    \node[right=.5 of 2,"\tiny 3" below] (3) {$\bullet$}
      child {};
  \end{tikzpicture}
  };
	\node (p2) [draw, my-red, right=2 of p1, "\color{my-red} $q$" above] {
	\begin{tikzpicture}[trees, sibling distance=2.5mm]
    \node["\tiny 1" below] (1) {\small $\blacksquare$}
      child {}
      child {}
      child {}
      child {};
    \node[right=.5 of 1,"\tiny 2" below] (2) {\small $\blacksquare$}
      child {}
      child {};
    \node[right=.5 of 2,"\tiny 3" below=.1] (3) {\small $\blacksquare$}
    ;
    \node[right=.5 of 3,"\tiny 4" below=.2] (4) {\small $\blacksquare$}
    ;
  \end{tikzpicture}
  };
\end{tikzpicture}
\]
To give a lens $p\to q$, we must send each $p$-position $i\in p(\1)$ to a $q$-position $j\in q(\1)$, then send each direction in $q[j]$ back to one in $p[i]$.
We can draw such a lens as follows.
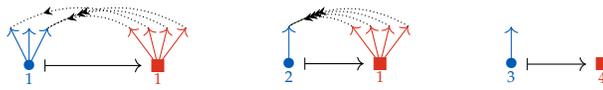
\[
\begin{tikzpicture}
	\node (p1) {
	\begin{tikzpicture}[trees, sibling distance=2.5mm]
    \node[my-blue, "{\color{my-blue}\tiny 1}" below] (1) {$\bullet$}
      child[my-blue] {coordinate (11)}
      child[my-blue] {coordinate (12)}
      child[my-blue] {coordinate (13)};
    \node[right=1.5 of 1, my-red, "{\color{my-red}\tiny 1}" below] (2) {\small $\blacksquare$}
      child[my-red] {coordinate (21)}
      child[my-red] {coordinate (22)}
      child[my-red] {coordinate (23)}
      child[my-red] {coordinate (24)};
    \draw[|->, shorten <= 3pt, shorten >= 3pt] (1) -- (2);
    \begin{scope}[densely dotted, bend right, decoration={markings, mark=at position 0.75 with \arrow{stealth}}]
      \draw[postaction={decorate}] (21) to (13);
      \draw[postaction={decorate}] (22) to (11);
      \draw[postaction={decorate}] (23) to (13);
      \draw[postaction={decorate}] (24) to (13);
    \end{scope}
  \end{tikzpicture}
	};
	\node (p2) [right=1 of p1] {
	\begin{tikzpicture}[trees, sibling distance=2.5mm]
    \node[my-blue, "{\color{my-blue}\tiny 2}" below] (1) {$\bullet$}
      child[my-blue] {coordinate (11)};
    \node[right=of 1, my-red, "{\color{my-red}\tiny 1}" below] (2) {\small $\blacksquare$}
      child[my-red] {coordinate (21)}
      child[my-red] {coordinate (22)}
      child[my-red] {coordinate (23)}
      child[my-red] {coordinate (24)};
    \draw[|->, shorten <= 3pt, shorten >= 3pt] (1) -- (2);
    \begin{scope}[densely dotted, bend right, decoration={markings, mark=at position 0.75 with \arrow{stealth}}]
      \draw[postaction={decorate}] (21) to (11);
      \draw[postaction={decorate}] (22) to (11);
      \draw[postaction={decorate}] (23) to (11);
      \draw[postaction={decorate}] (24) to (11);
    \end{scope}
  \end{tikzpicture}
	};
	\node (p3) [below right=-10mm and 1 of p2] {
	\begin{tikzpicture}[trees, sibling distance=2.5mm]
    \node[my-blue, "{\color{my-blue}\tiny 3}" below] (1) {$\bullet$}
      child[my-blue] {};
    \node[right=of 1, my-red, "{\color{my-red}\tiny 4}" below] (2) {\small $\blacksquare$}
		;
    \draw[|->, shorten <= 3pt, shorten >= 3pt] (1) -- (2);
  \end{tikzpicture}
	};
\end{tikzpicture}
\]
This represents one possible lens $f\colon p\to q$.
The horizontal solid arrows pointing rightward in the picture above tell us that the on-positions function $f_\1\colon p(\1)\to q(\1)$ is given by
\[
  f_\1(1) \coloneqq 1, \qquad f_\1(2) \coloneqq 1, \qqand f_\1(3) \coloneqq 4.
\]
Then the curved dashed arrows pointing leftward in the picture above describe the on-directions map $f^\sharp\colon q[f_\1(-)]\to p[-]$.
On the left, the arrows depict a possible on-directions function $f^\sharp_1\colon q[1]\to p[1]$ from the $4$ directions in $q[1]$ to the $3$ directions in $p[1]$.
In the middle, the arrows depict the only possible on-directions function $f^\sharp_2\colon q[1]\to p[2]$ because $|p[2]|=1$.
Finally, on the right, there are no curved arrows, depicting the only possible on-directions function $f^\sharp_3\colon q[4]\to p[3]$ because $|q[4]|=0$.
\end{example}

\begin{exercise}\label{exc.practice_poly_maps}
\begin{enumerate}
	\item Draw the corolla forests associated to $p\coloneqq\yon^\3+\yon+\1$, $q\coloneqq \yon^\2+\yon^\2+\2$, and $r\coloneqq\yon^\3$.
	\item Pick an example of a dependent lens $p\to q$ and draw it as we did in \cref{ex.practice_with_poly_morphisms}.
	\item Explain the behavior of your lens as an interaction protocol in terms of positions and directions.
	\item Explain in those terms why there can't be any lenses $p\to r$.
\qedhere
\end{enumerate}
\begin{solution}
\begin{enumerate}
    \item Here are the corolla forests associated to $p\coloneqq\yon^\3+\yon+\1$, $q\coloneqq \yon^\2+\yon^\2+\2$, and $r\coloneqq\yon^\3$ (with each root labeled for convenience).
    \[
    \begin{tikzpicture}[rounded corners]
    	\node (p) [draw, my-blue, "\color{my-blue} $p$" above] {
    	\begin{tikzpicture}[trees, sibling distance=2.5mm]
            \node["\tiny 1" below] (1) {$\bullet$}
              child {}
              child {}
              child {};
            \node[right=.5 of 1,"\tiny 2" below] (2) {$\bullet$}
              child {};
            \node[right=.5 of 2,"\tiny 3" below] (3) {$\bullet$};
        \end{tikzpicture}
        };
    	\node (q) [draw, my-red, right=2 of p, "\color{my-red} $q$" above] {
    	\begin{tikzpicture}[trees, sibling distance=2.5mm]
            \node["\tiny 1" below] (1) {\tiny$\blacksquare$}
              child {}
              child {};
            \node[right=.5 of 1,"\tiny 2" below] (2) {\tiny$\blacksquare$}
              child {}
              child {};
            \node[right=.5 of 2,"\tiny 3" below] (3) {\tiny$\blacksquare$};
            \node[right=.5 of 3,"\tiny 4" below] (4) {\tiny$\blacksquare$};
        \end{tikzpicture}
        };
    	\node (r) [draw, right=2 of q, "$r$" above] {
    	\begin{tikzpicture}[trees, sibling distance=2.5mm]
        \node["\tiny 1" below] (1) {\tiny$\blacktriangle$}
          child {}
          child {}
          child {};
        \end{tikzpicture}
      };
    \end{tikzpicture}
    \]
	\item Here is one possible lens $p\to q$ (you may have drawn others).

	\[
    \begin{tikzpicture}
    	\node (p1) {
        	\begin{tikzpicture}[trees, sibling distance=2.5mm]
                \node[my-blue, "\color{my-blue}\tiny 1" below] (1) {$\bullet$}
                  child[my-blue] {coordinate (11)}
                  child[my-blue] {coordinate (12)}
                  child[my-blue] {coordinate (13)};
                \node[right=1.5 of 1, my-red, "\color{my-red}\tiny 2" below] (2) {\tiny$\blacksquare$}
                  child[my-red] {coordinate (21)}
                  child[my-red] {coordinate (22)};
                \draw[|->, shorten <= 3pt, shorten >= 3pt] (1) -- (2);
                \begin{scope}[densely dotted, bend right, decoration={markings, mark=at position 0.75 with \arrow{stealth}}]
                  \draw[postaction={decorate}] (21) to (13);
                  \draw[postaction={decorate}] (22) to (11);
                \end{scope}
            \end{tikzpicture}
    	};
    	\node (p2) [right=1 of p1, yshift=-2mm] {
        	\begin{tikzpicture}[trees, sibling distance=2.5mm]
                \node[my-blue, "\color{my-blue}\tiny 2" below] (1) {$\bullet$}
                  child[my-blue] {coordinate (11)};
                \node[right=of 1, my-red, "\color{my-red}\tiny 4" below] (2) {\tiny$\blacksquare$};
                \draw[|->, shorten <= 3pt, shorten >= 3pt] (1) -- (2);
            \end{tikzpicture}
    	};
    	\node (p3) [right=1 of p2, yshift=-2mm] {
        	\begin{tikzpicture}[trees, sibling distance=2.5mm]
                \node[my-blue, "\color{my-blue}\tiny 3" below] (1) {$\bullet$};
                \node[right=of 1, my-red, "\color{my-red}\tiny 3" below] (2) {\tiny$\blacksquare$};
                \draw[|->, shorten <= 3pt, shorten >= 3pt] (1) -- (2);
            \end{tikzpicture}
    	};
    \end{tikzpicture}
    \]
    \item As depicted, our lens assigns to the first position of $p$ the second position of $q$, whose first and second directions are passed back to the third and first directions, respectively, of the first position of $p$.
    Then the second position of $p$ is assigned the fourth position of $q$, which has no directions; effectively, the choice of direction of the second position of $p$ has been canceled.
    Finally, the third position of $p$ is assigned the third position of $q$; here neither position has any directions.
    \item There cannot be a lens $p \to r$ for the following reason: if we send the third position of $p$, which has no directions, to the sole position of $r$, which has $3$ directions, then there is no way to pass a choice of one of those $3$ directions back to any of the options on the third menu of $p$, as there are no such options.
\end{enumerate}
\end{solution}
\end{exercise}

%-------- Section --------%
\section{Polybox pictures of dependent lenses}

\index{lens!polybox depiction|(}

Another way to visualize a dependent lens $f\colon p\to q$ is to draw the polyboxes for $p$ and $q$.
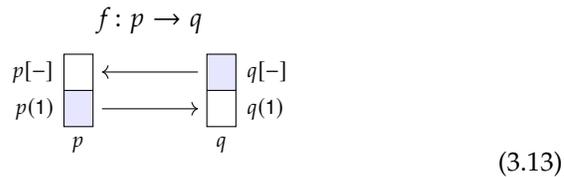
\begin{equation} \label{eqn.polybox_lens}
  \begin{tikzpicture}
    \node (f) ["$f\colon p\to q$" above] {
      \begin{tikzpicture}[polybox, tos]
        \node[poly, dom, "$p$" below] (p) {};
        \node[left=0pt of p_pos] {$p(\1)$};
        \node[left=0pt of p_dir] {$p[-]$};

        \node[poly, cod, right=of p, "$q$" below] (q) {};
        \node[right=0pt of q_pos] {$q(\1)$};
        \node[right=0pt of q_dir] {$q[-]$};

        \draw (p_pos) -- node[below] {} (q_pos);
        \draw (q_dir) -- node[above] {} (p_dir);
      \end{tikzpicture}
    };
  \end{tikzpicture}
\end{equation}
Thinking of the polyboxes as cells in a spreadsheet, the lens prescribes how the values of the cells are computed.
The boxes colored blue accept user input, while the other boxes are computed from that input according to the spreadsheet's rules.

\index{lens!spreadsheet depiction}

The arrows track the flow of information, starting from the lower left.
When the user fills the blue position box of $p$ with some $i\in p(\1)$, the arrow pointing right indicates that the lens fills the position box of $q$ with some $j\in q(\1)$ based on $i$.
The corresponding assignment $i\mapsto j$ is the on-position function $f_\1$ of the lens.

Then when the user fills the blue direction box of $q$ with some $b\in q[j]$, the arrow pointing left indicates that the lens fills the direction box of $p$ with some $a \in p[i]$ based on $i$ and $b$.
Fixing $i\in p(\1)$, the assignment $b\mapsto a$ is an on-directions function $f^\sharp_i$ of the lens.

Once all the boxes are filled, we obtain the following:
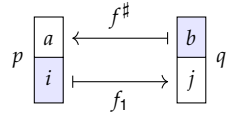
\[ \label{eqn.polybox_lens_filled}
\begin{tikzpicture}[polybox, mapstos]
  \node[poly, dom, "$p$" left] (p) {$a\vphantom{b}$\at$i\vphantom{j}$};
  \node[poly, cod, "$q$" right, right=of p] (q) {$b$\at$j$};
  \draw (p_pos) -- node[below] {$f_\1$} (q_pos);
  \draw (q_dir) -- node[above] {$f^\sharp$} (p_dir);
\end{tikzpicture}
\]
Here $j\coloneqq f_\1i$ and $a\coloneqq f^\sharp_ib$.
So a lens is any protocol that will fill the remaining boxes once the user fills the blue boxes, following the directions of the arrows drawn.
Be careful: although the arrow $f^\sharp$ is drawn from the codomain's direction box, it also takes into account what is entered into the domain's position box previously.
After all, the on-directions map of a lens is dependent on both a position of the domain and a direction of the codomain.

Here is an example of a lens depicted with polyboxes that would be difficult to draw with corolla forests.

\begin{example}[Modeling with a lens in polyboxes] \label{ex.lend-return}
  Caroline asks each of her parents for $20$ dollars. Each parent gives Caroline a positive amount of money not exceeding $20$ dollars. Caroline spends some of the money she receives before returning the remainder to each parent proportionally according to the amount she received from each.

  To model this interaction as a lens $f\colon p\to q$, we first define the polynomials $p$ and $q$.
  Let
  \[
    p\coloneqq\sum_{(i,\,j)\in(0,20]\times(0,20]}\yon^{[0,i]\times[0,j]}
  \]
  be the polynomial that models the parents: its position $(i,j)\in(0,20]\times(0,20]$ consists of the quantities of money that each parent gives to Caroline, and its direction $(i',j')\in[0,i]\times[0,j]$ at that position consists of the quantities of money that Caroline returns to each parent.
  Then let
  \[
    q\coloneqq\sum_{k\in(0,\,\infty)}\yon^{[0,k]}
  \]
  be the polynomial that models Caroline: its position $k\in(0,\infty)$ is the total money that Caroline receives (perhaps Caroline is prepared to receive more money than she is asking for, even if her parents are not prepared to give it), while its direction $r\in[0,k]$ is the money Caroline has remaining after spending some of it.
  Then we draw the lens $f$ that models their interaction in polyboxes as
  
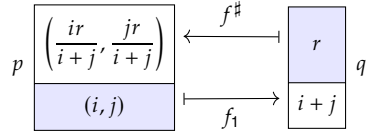
\[
  \begin{tikzpicture}[polybox, mapstos]
    \node[poly, dom, "$p$" left] (p) {$\left(\dfrac{ir}{i+j}, \dfrac{jr}{i+j}\right)$\at$(i,j)$};
    \node[poly, cod, "$q$" right, right=of p] (q) {$r\vphantom{\dfrac{jr}{i+j}}$\at$i+j$};
    \draw (p_pos) -- node[below] {$f_\1$} (q_pos);
    \draw (q_dir) -- node[above] {$f^\sharp$} (p_dir);
  \end{tikzpicture}
  \]
  We interpret this as saying that the on-positions function $f_\1$ from $p(\1)=(0,20]\times(0,20]$ to $q(\1)=(0,\infty)$ is defined on $(i,j)\in(0,20]\times(0,20]$ to be
  \[
    f_\1(i,j)\coloneqq i+j,
  \]
  while the on-directions function $f^\sharp_{(i,\,j)}$ from $q[i+j]=[0,i+j]$ to $p[(i,j)]=[0,i]\times[0,j]$ is defined on $r\in[0,i+j]$ to be
  \[
    f^\sharp_{(i,\,j)}(r) \coloneqq \left(\dfrac{ir}{i+j}, \dfrac{jr}{i+j}\right);
  \]
  the polybox picture expresses this more compactly.
  Notice that the position box of $q$ depends only on the position box of $p$, while the direction box of $p$ depends on the position box of $p$ as well as the direction box of $q$.
\end{example}

The above example illustrates how we can use polyboxes to specify a particular lens---or, equivalently, how we can use polyboxes to \emph{define} a lens, the same way we might define a function by writing it as a formula in a dummy variable.
Later on we will see how polyboxes help depict how lenses compose.
\index{lens!polybox depiction|)}

%-------- Section --------%
\section{Computations with dependent lenses}

Our concrete characterization of dependent lenses allows us to enumerate them.

\begin{example}[Enumerating lenses]
  Let $p\coloneqq \yon^\3+\2\yon$ and $q\coloneqq\yon^\4+\yon^\2+\2$, as in \cref{ex.practice_with_poly_morphisms}.
  Again, we draw them as corolla forests with their positions labeled:
  
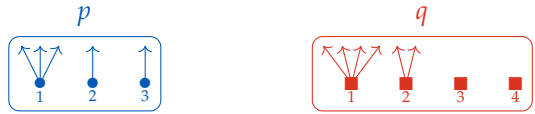
\[
  \begin{tikzpicture}[rounded corners]
    \node (p1) [draw, my-blue, "\color{my-blue} $p$" above] {
      \begin{tikzpicture}[trees, sibling distance=2.5mm]
        \node["\tiny 1" below] (1) {$\bullet$}
        child {}
        child {}
        child {};
        \node[right=.5 of 1,"\tiny 2" below] (2) {$\bullet$}
        child {};
        \node[right=.5 of 2,"\tiny 3" below] (3) {$\bullet$}
        child {};
      \end{tikzpicture}
    };
    \node (p2) [draw, my-red, right=2 of p1, "\color{my-red} $q$" above] {
      \begin{tikzpicture}[trees, sibling distance=2.5mm]
        \node["\tiny 1" below] (1) {\small$\blacksquare$}
        child {}
        child {}
        child {}
        child {};
        \node[right=.5 of 1,"\tiny 2" below] (2) {\small$\blacksquare$}
        child {}
        child {};
        \node[right=.5 of 2,"\tiny 3" below] (3) {\small$\blacksquare$}
        ;
        \node[right=.5 of 3,"\tiny 4" below] (4) {\small$\blacksquare$}
        ;
      \end{tikzpicture}
    };
  \end{tikzpicture}
  \]

  \index{lens!enumeration}

  How many lenses are there from $p$ to $q$? The first $p$-position must be sent to any $q$-position: 1, 2, 3, or 4. Sending it to 1 would require choosing an on-directions function $q[1]\to p[1]$, or $\4\to\3$; there are $3^4$ of these.
  Similarly, there are $3^2$ possible on-directions functions $q[2]\to p[1]$, as well as $3^0$ on-directions functions $q[3]\to p[1]$ and $3^0$ on-directions functions $q[4]\to p[1]$.
  Hence there are a total of $3^4+3^2+3^0+3^0=92$ ways to choose $f_\1(1)$ and $f^\sharp_1$.

  The second $p$-position must also be sent to 1, 2, 3, or 4 before selecting $f^\sharp_2$; there are $1^4+1^2+1^0+1^0=4$ ways to do this.
  Identically, there are $4$ ways to choose $f_\1(3)$ and $f^\sharp_3$.

  In total, there are $92 \cdot 4 \cdot 4=1472$ lenses $p\to q$.
  This coincides with what we obtain by taking the cardinality of both sides of \eqref{eqn.main_formula} and plugging in our values for $p$ and $q$:
  \begin{align*}
    |\poly(p, q)| &= \prod_{i \in p(\1)} |q(p[i])| \\
    &= \prod_{i\in\3}\left(|p[i]|^4 + |p[i]|^2 + 2\right) \\
    &= (3^4 + 3^2 + 2)(1^4 + 1^2 + 2)^2 \\
    &= 92 \cdot 4^2 = 1472.
  \end{align*}
\end{example}

\index{Yoneda lemma}
\begin{exercise}
For any polynomial $p$ and set $A$, e.g.\ $A\coloneqq\2$, the Yoneda lemma gives an isomorphism $\poly(\yon^A,p)\iso p(A)$, so the number of lenses $\yon^A\to p$ should be equal to the cardinality of $p(A)$.
\begin{enumerate}
	\item Choose a polynomial $p$ with finitely many positions and directions and draw both $\yon^\2$ and $p$ as corolla forests.
	\item Count all the lenses $\yon^\2\to p$. How many are there?
	\item Compute the cardinality of $p(\2)$. Is this the same as the previous answer?
\qedhere
\end{enumerate}
\begin{solution}
\begin{enumerate}
    \item We let $p\coloneqq\yon^\3+\1$ (you could have selected others) and draw both $\yon^\2$ and $p$ as corolla forests, labeling each root for convenience.
    \[
    \begin{tikzpicture}[rounded corners]
    	\node (y2) [draw, my-blue, "\color{my-blue}$\yon^\2$" above] {
    	\begin{tikzpicture}[trees, sibling distance=2.5mm]
            \node["\tiny 1" below] (1) {$\bullet$}
              child {}
              child {};
        \end{tikzpicture}
        };
    	\node (p) [draw, my-red, right=2 of y2, "\color{my-red}$p$" above] {
    	\begin{tikzpicture}[trees, sibling distance=2.5mm]
            \node["\tiny 1" below] (1) {$\bullet$}
              child {}
              child {}
              child {};
            \node[right=.5 of 1,"\tiny 2" below] (2) {$\bullet$};
        \end{tikzpicture}
        };
    \end{tikzpicture}
    \]

    \item When constructing a lens $\yon^2\to p$, the unique position of $\yon^\2$ can be sent to either $p$-position.
    If it is sent to the first $p$-position, then each of the $3$ directions in $p[1]$ must be sent to one of the $2$ directions of $\yon^\2$, for a total of $2^3 = 8$ lenses.
    Otherwise, the unique position of $\yon^\2$ is sent to the second $p$-position, at which there are no directions; so there is exactly $1$ lens like this.
    Hence there are $8+1=9$ lenses $\yon^\2\to p$.

    \item The cardinality of $p(\2)$ is $|\2^\3+\1|=9$, which agrees with previous answer, as predicted by the Yoneda lemma.
\end{enumerate}
\end{solution}
\end{exercise}

\begin{exercise}
For each of the following polynomials $p,q$, compute the number of lenses $p\to q$.
\begin{enumerate}
	\item $p\coloneqq\yon^\3$,\quad $q\coloneqq\yon^\4$.
	\item $p\coloneqq\yon^\3+\1$,\quad $q\coloneqq\yon^\4$.
	\item $p\coloneqq\yon^\3+\1$,\quad $q\coloneqq\yon^\4+\1$.
	\item $p\coloneqq\4\yon^\3+\3\yon^\2+\yon$,\quad $q\coloneqq\yon$.
	\item $p\coloneqq\4\yon^\3$,\quad $q\coloneqq\3\yon$.
\qedhere
\end{enumerate}
\begin{solution}
By \eqref{eqn.main_formula}, we have for all $p,q\in\poly$ that
\[
    |\poly(p,q)| = \prod_{i\in p(\1)}|q(p[i])|.
\]
\begin{enumerate}
	\item If $p\coloneqq\yon^\3$ and $q\coloneqq\yon^\4$, then
	\[
	    |\poly(p,q)| = \prod_{i\in\1}|p[i]|^4 = 3^4 = 81.
	\]
	\item If $p\coloneqq\yon^\3+\1$ and $q\coloneqq\yon^\4$, then
	\[
	    |\poly(p,q)| = \prod_{i\in\2}|p[i]|^4 = 3^4\cdot0^4 = 0.
	\]
	\item If $p\coloneqq\yon^\3+\1$ and $q\coloneqq\yon^\4+\1$, then
	\[
	    |\poly(p,q)| = \prod_{i\in\2}\left(|p[i]|^4+1\right) = (3^4 + 1)(0^4 + 1) = 82.
	\]
	\item If $p\coloneqq\4\yon^\3+\3\yon^\2+\yon$ and $q\coloneqq\yon$, then
	\[
	    |\poly(p,q)| = \prod_{i\in\8}|p[i]| = 3^4\cdot2^3\cdot1 = 648.
	\]
	\item If $p\coloneqq\4\yon^\3$ and $q\coloneqq\3\yon$, then
	\[
	    |\poly(p,q)| = \prod_{i\in\4}3|p[i]| = (3\cdot3)^4 = 6561.
	\]
\qedhere
\end{enumerate}
\end{solution}
\end{exercise}

The following exercises provide alternative formulas for the set of lenses between two polynomials.

\begin{exercise}\label{exc.practice_sum_prod}
\begin{enumerate}
\item Show that the following are isomorphic:
\begin{equation}\label{eqn.poly_p_q}
  \poly(p,q)
  \iso
  \prod_{i\in p(\1)}\sum_{j\in q(\1)}\prod_{b\in q[j]}\sum_{a\in p[i]}\1.
\end{equation}
\item \label{exc.practice_sum_prod.useful} Show that the following are isomorphic:
	\begin{equation}\label{eqn.useful_misc472}
	\poly(p,q)\iso\sum_{f_\1\colon p(\1)\to q(\1)}\,\prod_{j\in q(\1)}\smset\Bigg(q[j],\prod_{\substack{i \in p(\1), \\ f_\1i = j}}p[i]\Bigg).
	\end{equation}
\item Using the language of positions and directions, describe how an element of the right hand side of \eqref{eqn.useful_misc472} corresponds to a lens $p\to q$.
\qedhere
\end{enumerate}
\begin{solution}
\begin{longenum}
\item By \eqref{eqn.main_formula}, it suffices to show that for all $i\in p(\1)$ and $j\in q(\1)$, we have
\[
    p[i]^{q[j]} \iso \prod_{b\in q[j]}\sum_{a\in p[i]}\1.
\]
Indeed, by \eqref{eqn.push_prod_sum_set_indep} and \cref{exc.on_sums_prods_sets}, we have
\begin{align*}
    \prod_{b \in q[j]} \sum_{a \in p[i]} \1 &\iso \sum_{\bar{a} \colon q[j] \to p[i]} \, \prod_{b \in q[j]} \1 \tag*{\eqref{eqn.push_prod_sum_set_indep}} \\
    &\iso \sum_{\bar{a} \colon q[j] \to p[i]} \1 \tag{\cref{exc.on_sums_prods_sets} \cref{exc.on_sums_prods_sets.prod}} \\
    &\iso \smset(q[j], p[i]) \tag{\cref{exc.on_sums_prods_sets} \cref{exc.on_sums_prods_sets.sum}} \\
    &\iso p[i]^{q[j]}.
\end{align*}

\item By \eqref{eqn.main_formula} and \eqref{eqn.push_prod_sum_set_indep}, we have
\begin{align*}
    \poly(p, q) &\iso \prod_{i \in p(\1)} \sum_{j \in q(\1)} p[i]^{q[j]} \tag*{\eqref{eqn.main_formula}} \\
    &\iso \sum_{f_\1 \colon p(\1) \to q(\1)} \; \prod_{i \in p(\1)} p[i]^{q[f_\1i]} \tag*{\eqref{eqn.push_prod_sum_set_indep}} \\
    &\iso \sum_{f_\1 \colon p(\1) \to q(\1)} \; \prod_{j \in q(\1)} \; \prod_{\substack{i \in p(\1), \\ f_\1i = j}} p[i]^{q[j]} \tag{$\ast$} \\
    &\iso \sum_{f_\1\colon p(\1)\to q(\1)} \; \prod_{j\in q(\1)}\smset\Bigg(q[j],\prod_{\substack{i \in p(\1), \\ f_\1i = j}}p[i]\Bigg) \tag{Universal property of products}
\end{align*}
where ($\ast$) follows from the fact that for any function $f_\1\colon p(\1)\to q(\1)$, its domain $p(\1)$ can be written as the disjoint union of preimages $f_\1\inv(j) = \{i \in p(1) \mid f_\1i = j\}$ for each $j \in q(\1)$.

\item To explain how an element of the set
\[
	D_{p,q} \coloneqq \sum_{f_\1\colon p(\1)\to q(\1)}\prod_{j\in q(\1)}\smset\Bigg(q[j],\prod_{\substack{i \in p(\1), \\ f_\1i = j}}p[i]\Bigg)
\]
corresponds to a lens $p\to q$, we first give the instructions for choosing an element of $D_{p,q}$ as a nested list.
\begin{quote}
To choose an element of $D_{p,q}$:
\begin{longenum}
    \item choose a function $f_\1 \colon p(\1) \to q(\1)$;
    \item for each element $j \in q(\1)$:
    \begin{longenum}
        \item for each element of $q[j]$:
        \begin{longenum}
            \item for each element $i \in p(\1)$ satisfying $f_\1i = j$:
            \begin{longenum}
                \item choose an element of $p[i]$.
            \end{longenum}
        \end{longenum}
    \end{longenum}
\end{longenum}
\end{quote}
So $f_\1$ sends each $p$-position to a $q$-position, as an on-positions function should.
Then for each $q$-position $j$, each $q[j]$-direction $b$, and each $p$-position $i$ that $f_\1$ sends to $j$, we choose an element of $p[i]$ that $f^\sharp_i\colon q[j]\to p[i]$ assigns to $b$.
As every $p$-position $i$ is sent to some $q$-position $j$ by $f_\1$, this completely characterizes $f^\sharp_i$ for every $p$-position $i$.
\end{longenum}
\end{solution}
\end{exercise}

In \eqref{eqn.poly_coprod}, we gave an explicit formula for coproducts in $\poly$ inherited from coproducts in $\smset^\smset$, as justified by \cref{prop.poly_coprods}.
We can now use \eqref{eqn.main_formula} to directly verify that the expression on the right hand side of \eqref{eqn.poly_coprod} satisfies the universal property for the coproduct of polynomials $(p_i)_{i\in I}$ in $\poly$.

\index{polynomial functor!sum of polynomials}

\begin{exercise}%\label{exc.poly_coprod}
Use \eqref{eqn.main_formula} to verify that
\[
    \poly\left(\sum_{(i,\,j)\in\sum_{i\in I}p_i(\1)}\yon^{p_i[j]}, q\right) \iso \prod_{i \in I}\poly(p_i, q)
\]
for all polynomials $(p_i)_{i\in I}$ and $q$.
\begin{solution}
We have
\begin{align*}
    \poly\left(\sum_{(i,\,j) \in \sum_{i \in I} p_i(\1)} \yon^{p_i[j]}, q\right) &\iso \prod_{(i,\,j) \in \sum_{i \in I} p_i(\1)} q(p_i[j])
    \tag*{\eqref{eqn.main_formula}} \\
    &\iso \prod_{i \in I} \prod_{j \in p_i(\1)} q(p_i[j]) \\
    &\iso \prod_{i \in I} \poly(p_i, q).
    \tag*{\eqref{eqn.main_formula}}
\end{align*}
\end{solution}
\end{exercise}

For the remaining exercises in this section, we introduce the concept of the \emph{derivative} $\dot{p}$ of a polynomial $p$.

\begin{example}[Derivatives]\label{ex.derivatives}
The \emph{derivative} of a polynomial $p$, denoted $\dot{p}$, is defined as
\[
  \dot{p}\coloneqq\sum_{i\in p(\1)}\sum_{a\in p[i]}\yon^{p[i]-\{a\}}.
\]
For example, if $p\coloneq\yon^{\{U,V,W\}}+\{A\}\yon^{\{X\}}+\{B\}\yon^{\{X\}}$ then
\[\dot{p}=\{U\}\yon^{\{V,W\}}+\{V\}\yon^{\{U,W\}}+\{W\}\yon^{\{U,V\}}+\{(A,X)\}\yon^\0+\{(B,X)\}\yon^\0.\]
Up to isomorphism $p\cong\yon^\3+\2\yon$ and $\dot{p}\cong\3\yon^\2+\2$.
Indeed, this coincides with the familiar notion of derivatives of polynomials from calculus.

Since
\[
  \dot{p}\yon\iso\sum_{i\in p(\1)}\sum_{a\in p[i]}\yon^{p[i]-\{a\}}\yon\iso\sum_{i\in p(\1)}\sum_{a\in p[i]}\yon^{p[i]-\{a\}+\1}\iso\sum_{i\in p(\1)}\sum_{a\in p[i]}\yon^{p[i]},
\]
there exists a canonical lens $\dot{p}\yon\to p$; you will define this lens in \cref{exc.deriv-ops}.
The lens arises in computer science in the context of ``plugging in to one-hole contexts''; we will not explore that here, but see \cite{mcbride2001derivative} and \cite{abbott2003derivatives} for details.%The Derivative of a Regular Type is its Type of One-Hole Contexts. http://citeseerx.ist.psu.edu/viewdoc/download?doi=10.1.1.22.8611&rep=rep1&type=pdf.
%https://www.cs.nott.ac.uk/~psztxa/publ/tlca03.pdf

A lens $f'\colon p\to \dot{q}$ is similar to a lens $f\colon p\to q$, except that each $p$-position explicitly selects a direction of $q$ to remain unassigned.
More precisely, while the on-positions function of $f$ sends each $p$-position to a $q$-position, the on-positions function of $f'$ sends each $p$-position $i$ to some $(j,a)\in \sum_{j\in q(\1)}q[j]$, picking out not only a $q$-position $j$, but also a $q[j]$-direction $a$.
Then the on-directions function of $f'$ at $i$ sends every $q[j]$-direction \emph{other than $a$} back to a $p[i]$-direction.
\index{polynomial functor!derivative of|see{derivative}}\index{derivative}
\end{example}

\begin{exercise} \label{exc.deriv-ops}
The derivative is not very well-behaved categorically, but it is nevertheless intriguing.
Take $p,q\in\poly$.
\begin{enumerate}
	\item Give an explicit construction for the canonical lens $\dot{p}\yon\to p$ from \cref{ex.derivatives}.
	\item Is there always a lens $p\to \dot{p}$?
  If so, prove it; if not, give a counterexample.
	\item Is there always a lens $\dot{p}\to p$?
  If so, prove it; if not, give a counterexample.
	\item Given a lens $p\to q$, is there always a lens $\dot{p}\to\dot{q}$?
  If so, prove it; if not, give a counterexample.
	\item We will define the binary operations $\otimes$ and $\ihom{-,-}$ on $\poly$ later on in \eqref{eqn.parallel_def} and \eqref{eqn.par_hom}; and in \cref{exc.dir_hom_p_yon_dir_p}, you will be able to use \cref{exc.par_hom_sum} to deduce that
	\begin{equation} \label{eqn.dir_hom_p_yon_dir_p}
	    \ihom{p, \yon} \otimes p \iso \sum_{f \in \prod_{i \in p(\1)} p[i]} \; \sum_{i \in p(\1)} \yon^{p(\1) \times p[i]},
	\end{equation}
	Construct a canonical lens $\ihom{p,\yon}\otimes p\to \dot{p}$.
	\item In \cref{ex.derivatives}, we described a lens $p\to\dot{q}$ in terms of ``unassigned'' directions. Describe a lens $p\yon\to q$ in terms of ``unassigned'' directions as well.
	\qedhere
\end{enumerate}
\begin{solution}
\begin{enumerate}
  \item We construct a lens $\dot{p}\yon\to p$, or a lens
  \[
    \sum_{i\in p(\1)}\sum_{a\in p[i]}\yon^{p[i]}\to\sum_{i\in p(\1)}\yon^{p[i]},
  \]
  as follows.
  It sends each position $(i,a)\in\sum_{i\in p(\1)}p[i]$ of $\dot{p}\yon$ to its first projection $i\in p(\1)$, and it is the identity on directions.

	\item There is not always a lens $p\to \dot{p}$: if $p\coloneqq\1$, then $\dot{p}\coloneqq\0$, and there is no lens $\1\to\0$.

	\item There is not always a lens $\dot{p}\to p$: take $p\coloneqq\yon$, so $\dot{p}\coloneqq\1$.
	A lens $\1\to\yon$ must have an on-directions function $\1\to\0$, but there is no such function.

	\item We show that even when there is a lens $p \to q$, there is not necessarily a lens $\dot{p}\to\dot{q}$.
	Take $p \coloneqq \yon$ and $q \coloneqq \1$.
	Then there is a lens $p \to q$ that sends the unique position of $\yon$ to the unique position of $\1$ and is the empty function on directions.
	But $\dot{p} = \1$ and $\dot{q} = \0$, and there is no lens $\1 \to \0$.

	\item We construct a lens $g\colon\ihom{p,\yon}\otimes p\to \dot{p}$, where $\ihom{p,\yon} \otimes p$ is given by \eqref{eqn.dir_hom_p_yon_dir_p}, as follows.
	The on-positions function $g_\1$ takes $f\in\prod_{i\in p(\1)}p[i]$ and $i\in p(\1)$ and sends the pair of them to the $\dot{p}$-position corresponding to $i\in p(\1)$ and $fi\in p[i]$.
	Then $\dot{p}[(i, fi)]=p[i]-\{fi\}$ and $(\ihom{p,\yon}\otimes p)[(f, i)]\iso p(\1)\times p[i]$, so the on-directions function $g^\sharp_{(f,\,i)}$ can send each $a\in p[i]-\{fi\}$ to $(i,a)\in p(\1)\times p[i]$.

	\item We describe a lens $p\yon\to q$ in terms of ``unassigned'' directions.
	Observe that $p\yon$ has the same positions as $p$ but has one more direction than $p$ does at each position.
	Given a position $i\in p(\1)$, we denote this extra $(p\yon)[i]$-direction by $\ast_i$, identifying $(p\yon)[i]$ with $p[i]+\{\ast_i\}$.
	So a lens $f\colon p\yon\to q$ sends each $p$-position $i$ to a $q$-position $j$, but every $q[j]$-direction could be sent back to either an original $p[i]$-direction or the extra $(p\yon)[i]$-direction $\ast_i$.
	So a lens $p\yon\to q$ is like a lens $p\to q$, except that some of the directions of $q$ may remain ``unassigned'' to any direction of $p$, which we signify by assigning them to $\ast_i$ instead.
  In other words, a lens $p\yon\to q$ could be interpreted as a lens $p\to q$ whose on-directions functions may only be partially defined.
\end{enumerate}
\end{solution}
\index{derivative!not well-behaved in $\poly$}
\end{exercise}

We will not use derivatives very much in the rest of this text, except as shorthand to denote the set of all directions of a polynomial: given a polynomial $p$, its directions comprise the set $\dot{p}(\1)$.

\begin{exercise} \label{exc.deriv-directions}
  Show that $\dot{p}(\1)$ is isomorphic to the set of all directions of $p$ (i.e.\ the sum of all direction-sets of $p$).
  \begin{solution}
    We have
    \[
      \dot{p}(\1)\iso\sum_{i\in p(\1)}\sum_{a\in p[i]}\1^{p[i]-\{a\}} \iso \sum_{i \in p(\1)} p[i],
    \]
    which is precisely the set of all directions of $p$.
  \end{solution}
  \index{derivative! $\dot{p}(1)$ as directions of $p$}
\end{exercise}

%-------- Section --------%
\section{Dependent lenses between special polynomials}

In \cref{sec.poly.obj.spec}, we considered four special classes of polynomials (here $I$ and $A$ are sets): constant polynomials $I$, linear polynomials $I\yon$, representable polynomials $\yon^A$, and monomials $I\yon^A$.
Of special note are the constant and linear polynomial $\0$, the constant and representable polynomial $\1$, and the linear and representable polynomial $\yon$.
We now consider lenses with these special polynomials as domains or codomains, highlighting some important examples using polyboxes and leaving most of the rest as exercises.
Let $p$ be a polynomial throughout.

\begin{example}[Lenses from linear polynomials] \label{ex.lens-from-lin}\index{lens!from linear polynomial}
  A lens $f\colon I\yon\to p$ can be drawn in polyboxes as follows:
  
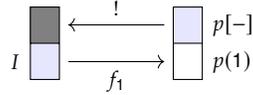
\[
    \begin{tikzpicture}[polybox, tos]
      \node[poly, linear dom] (p) {};
      \node[left=0pt of p_pos] {$I$};

      \node[poly, cod, right=of p] (q) {};
      \node[right=0pt of q_pos] {$p(\1)$};
      \node[right=0pt of q_dir] {$p[-]$};

      \draw (p_pos) -- node[below] {$f_\1$} (q_pos);
      \draw (q_dir) -- node[above] {$!$} (p_dir);
    \end{tikzpicture}
  \]
  Recall that we shade in the direction box of a linear polynomial to indicate that it can only be filled with one entry.
  Hence the on-directions map of $f$ is uniquely determined, so $f$ is completely characterized by its on-positions function $f_\1$.
  We conclude that lenses $I\yon\to p$ can be identified with functions $I\to p(\1)$.
\end{example}

\begin{exercise}[Lenses from $\0$ and $\yon$] \label{exc.lens-from-0-or-yon}
\begin{enumerate}
  \item Use \cref{ex.lens-from-lin} to verify that $\0$ is the initial object of $\poly$.
  \item Use \cref{ex.lens-from-lin} to show that lenses $\yon\to p$ can be identified with $p$-positions. \qedhere
\end{enumerate}
\begin{solution}
  \begin{enumerate}
    \item Since $\0\iso\0\yon$ is a linear polynomial, \cref{ex.lens-from-lin} tells us that lenses $\0\to p$ can be identified with functions $\0\to p(\1)$.
    There is exactly one function $\0\to p(\1)$, so there is exactly one lens $\0\to p$.
    \item Since $\yon\iso\1\yon$ is a linear polynomial, \cref{ex.lens-from-lin} tells us that lenses $\yon\to p$ can be identified with functions $\1\to p(\1)$, which in turn can be identified with elements of $p(\1)$.
  \end{enumerate}
\end{solution}
\end{exercise}

\begin{exercise}[Lenses to linear polynomials]\index{lens!to linear polynomial}
  Characterize lenses $p\to I\yon$ in terms of $I$ and the positions and directions of $p$.
\begin{solution}
  A lens $f\colon p\to I\yon$ consists of an on-positions function $f_\1\colon p(\1)\to I$ and, for each $j\in p(\1)$, an on-directions function $f^\sharp_j\colon\1\to p[j]$.
  Equivalently, this is a function $p(\1)\to I$ and a choice of direction at every position of $p$, i.e.\ a dependent function $(j\in p(\1))\to p[j]$.
\end{solution}
\end{exercise}

\begin{example}[Lenses to constants] \label{ex.lens-to-const}\index{lens!to constant polynomial}
  A lens $f\colon p\to I$ can be drawn in polyboxes as follows:
  
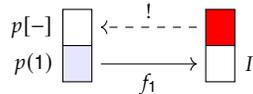
\[
  \begin{tikzpicture}[polybox, tos]
    \node[poly, dom] (p) {};
    \node[left=0pt of p_pos] {$p(\1)$};
    \node[left=0pt of p_dir] {$p[-]$};

    \node[poly, constant, right=of p] (q) {};
    \node[right=0pt of q_pos] {$I$};

    \draw (p_pos) -- node[below] {$f_\1$} (q_pos);
    \draw[dashed] (q_dir) -- node[above] {$!$} (p_dir);
  \end{tikzpicture}
  \]
  Recall that we color the direction box of a constant polynomial red to indicate that it cannot be filled with any entry.
  Hence the on-directions map of $f$ is again uniquely determined, so $f$ is completely characterized by its on-positions function $f_\1$.
  (While the on-directions map of $f$ does exist---it is vacuous---it can never produce an element to fill the direction box of $p$, so we draw it with a dashed line.)
  We conclude that lenses $p\to I$ can be identified with functions $p(\1)\to I$.
\end{example}

\begin{exercise}[Lenses to $\1$ and $\0$]
  \begin{enumerate}
    \item Use \cref{ex.lens-to-const} to show that $\1$ is the terminal object of $\poly$.
    \item Show that there is exactly one lens whose codomain is $\0$.
    What is its domain? \qedhere
  \end{enumerate}
  \begin{solution}
    \begin{enumerate}
      \item Since $\1$ is a constant, \cref{ex.lens-to-const} tells us that lenses $p\to\1$ can be identified with functions $p(\1)\to\1$.
      There is exactly one function $p(\1)\to\1$, so there is exactly one lens $p\to\1$.
      \item Since $\0$ is a constant, \cref{ex.lens-to-const} tells us that lenses $p\to\0$ can be identified with functions $p(\1)\to\0$, which only exist when $p(\1)=\0$.
      The only polynomial with an empty position-set is $\0$ itself, and there is a unique function from the set $\0$ to itself, so there is a unique lens from the constant polynomial $\0$ to itself as well.
      If $p$ is not the constant $\0$, then there are no lenses $p\to\0$.
    \end{enumerate}
  \end{solution}
\end{exercise}

\begin{exercise}[Lenses from constants such as $\1$]\index{lens!from constant polynomial}
  \begin{enumerate}
    \item Characterize lenses $I\to p$ in terms of $I$ and the positions and directions of $p$.
    You may find it helpful to refer to $p(\0)$; see \cref{exc.apply0}.
    \item Use the previous part to characterize lenses $\1\to p$. \qedhere
  \end{enumerate}
  \begin{solution}
  \begin{enumerate}
    \item A lens $f\colon I\to p$ consists of an on-positions function $f_\1\colon I\to p(\1)$ and, for each $i\in I$, an on-directions function $f^\sharp_i\colon p[i]\to\0$.
    There is exactly one such on-directions function when $p[i]=\0$ and no such on-directions function otherwise.
    It follows that a lens $f\colon I\to p$ can be identified with a function $f_\1\colon I\to p(\1)$ whose image is contained in the set of $p$-positions with no directions.
    By \cref{exc.apply0}, this set of $p$-positions can be identified with the set $p(\0)$ (the \emph{constant} term of $p$); so a lens $I\to p$ is equivalent to a function $I\to p(\0)$.
    \item From the previous part, a lens $\1\to p$ may be identified with a function $\1\to p(\0)$ and thus an element of $p(\0)$.
  \end{enumerate}
  \end{solution}
\end{exercise}

\index{Yoneda lemma}

We already know from the Yoneda lemma (see \cref{exc.poly_morph_yoneda}) that lenses $\yon^A\to p$ correspond to elements of $p(A)$, so we understand lenses from representables. Thus we turn our attention to lenses $p\to\yon^A$.\index{lens!from representable polynomial}

\begin{example}[Lenses to representables]\index{lens!to representable polynomial}
  A lens $f\colon p\to\yon^A$ can be drawn in polyboxes as follows:
  
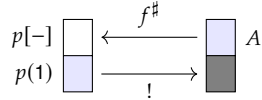
\[
  \begin{tikzpicture}
    \node (f) {
      \begin{tikzpicture}[polybox, tos]
        \node[poly, dom] (p) {};
        \node[left=0pt of p_pos] {$p(\1)$};
        \node[left=0pt of p_dir] {$p[-]$};

        \node[poly, pure cod, right=of p] (q) {};
        \node[right=0pt of q_dir] {$A$};

        \draw (p_pos) -- node[below] {$!$} (q_pos);
        \draw (q_dir) -- node[above] {$f^\sharp$} (p_dir);
      \end{tikzpicture}
    };
  \end{tikzpicture}
  \]
  Recall that we shade in the position box of a representable to indicate that it can only be filled with one entry.
  Hence the on-positions function of $f$ is uniquely determined, so $f$ is completely characterized by its on-directions map $f^\sharp$, which takes a $p$-position $i$ and a direction $a\in A$ and sends them to a direction $b\in p[i]$.
  We conclude that lenses $p\to\yon^A$ can be identified with dependent functions $((i,a)\in p(\1)\times A)\to p[i]$.
\end{example}

\begin{example}[Lenses to $\yon$]
  As a special case of the previous example, a lens $\gamma\colon p\to\yon$ can be drawn in polyboxes as follows:
  
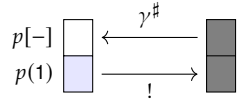
\[
  \begin{tikzpicture}
    \node (f) {
      \begin{tikzpicture}[polybox, tos]
        \node[poly, dom] (p) {};
        \node[left=0pt of p_pos] {$p(\1)$};
        \node[left=0pt of p_dir] {$p[-]$};

        \node[poly, identity, right=of p] (q) {};

        \draw (p_pos) -- node[below] {$!$} (q_pos);
        \draw (q_dir) -- node[above] {$\gamma^\sharp$} (p_dir);
      \end{tikzpicture}
    };
  \end{tikzpicture}
  \]\index{dependent function!lens to $\yon$ as}
  Such lenses can be identified with dependent functions $(i\in p(\1))\to p[i]$, which, abusing notation, we also denote by $\gamma$.
  For each $p$-position $i$ in the blue position box, $\gamma$ picks out a $p[i]$-direction to fill the unshaded direction box.
  (Remember that the arrow labeled $\gamma^\sharp$ depends not only on the direction box to its right, but also on the position box of $p$.)
  So it makes sense to abbreviate the polybox picture of $\gamma$ like so:
  
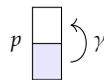
\begin{equation}\label{eqn.map_to_0ary_composite}
    \begin{tikzpicture}[polybox, tos]
      \node[poly, dom, "$p$" left] (p) {};
      \draw (p_pos) to[climb'] node[right] {$\gamma$} (p_dir);
    \end{tikzpicture}
  \end{equation}
\end{example}

\index{dependent function!lens to $\yon$ as}

The correspondence between lenses $p\to\yon$ and dependent functions $(i\in p(\1))\to p[i]$ exhibited in the previous example also follows directly from \eqref{eqn.main_formula}: taking $q\coloneqq\yon$, we have
\begin{equation} \label{eqn.gamma_prod}
\poly(p,\yon)\iso\prod_{i\in p(\1)}\sum_{j\in\1}p[i]^\1\iso\prod_{i\in p(\1)}p[i],
\end{equation}
where the right hand side is precisely the set of dependent functions $(i\in p(\1))\to p[i]$.
By \cref{exc.product_as_sections}, such functions may be identified with the \emph{sections} of the projection function from $\dot{p}(\1)\iso\sum_{i\in p(\1)}p[i]$, the set of all directions of $p$ (see \cref{exc.deriv-directions}), to $p(\1)$, the set of all positions of $p$, sending each $(i,a)\in\dot{p}(\1)$ to $i\in p(\1)$.
The fact that this projection determines $p$ (up to isomorphism) motivates the following definition.

\begin{definition}[Section; bundle] \label{def.sec-bun}\index{section|(}\index{lens!to $\yon$ as section}\index{bundle}\index{section|(}
  For $p\in\poly$, a \emph{section} of $p$ is a lens $p\to\yon$.
  We denote the set of all sections of $p$ by $\Gamma(p)$; that is,
  \begin{equation} \label{eqn.gamma_def}
    \Gamma(p)\coloneqq\poly(p,\yon).
  \end{equation}
  The \emph{bundle} of $p$, denoted $\pi_p$, is the projection function
  \[
    \dot{p}(\1)\iso\sum_{i\in p(\1)}p[i]\to p(\1)
  \]
  sending $(i,a)\mapsto i$.
\end{definition}

With this terminology, we can say that $p$ is determined (up to isomorphism) by its bundle, and that the sections of $p$ can be identified with the sections of its bundle.

To visualize the bundle of $p$, simply draw it as a corolla forest: a \emph{bundle} of arrows.
The bundle projects each leaf down to its root.
To visualize a section of $p$, picture its corollas piled atop each other; a section $\gamma\colon p\to\yon$ is then a \emph{cross-section} of this pile of $p$-corollas, picking out an arrow from each one---a direction at each position.

Alternatively, you could think of the arrow curving back to the polyboxes for $p$ in our picture  \eqref{eqn.map_to_0ary_composite} of a section $\gamma\colon p\to\yon$ as \emph{sectioning} off the polyboxes for $p$ from any polyboxes that may otherwise appear to its right.
We clarify this intuition by returning to a previous example of a polynomial and considering its sections.

\begin{example}[Modeling with sections] \label{ex.spend-section}
  Recall from \cref{ex.lend-return} the polynomial
  \[
    q\coloneqq\sum_{k\in(0,\infty)}\yon^{[0,k]}
  \]
  whose positions $k\in(0,\infty)$ are the possible quantities of money that Caroline receives from her parents and whose directions $r\in[0,k]$ are the possible quantities of money that Caroline has remaining after spending some of it.

  One component that was missing from our model was how Caroline spends her money.
  A section for $q$ fills in this gap by closing the loop from the money Caroline receives to the money she has remaining.
  Explicitly, a section $\gamma\colon q\to\yon$ corresponds to a dependent function $(k\in(0,\infty))\to[0,k]$ that uses the amount of money that Caroline receives to determine the amount of money that she will have remaining.

  For instance, if Caroline always spends half the money she receives, then the polyboxes for the section $\gamma\colon q\to\yon$ that models this behavior can be drawn as follows:
  
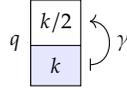
\[
  \begin{tikzpicture}[polybox,mapstos]
    \node[poly, dom, "$q$" left] (p) {$k/2$\at$k$};
    \draw (p_pos) to[climb'] node[right] {$\gamma$} (p_dir);
  \end{tikzpicture}
  \]
  Without $\gamma$, we do not know how much money Caroline will decide to spend; having $\gamma$ makes her decision deterministic and sections this decision off from unknown variables.
  Of course, the position box of $q$, the amount of money Caroline receives, is still open to outside influence, as determined by a lens to $q$ such as the one from \cref{ex.lend-return}.
\end{example}

The definition of $\Gamma$ given in \eqref{eqn.gamma_def} makes $\Gamma$ a functor $\poly\to\smset\op$ satisfying the following.

\begin{proposition}\label{prop.gamma_pres_coproduct}
  The sections functor $\Gamma\colon\poly\to\smset\op$ sends $(0,+)$ to $(1,\times)$:
  \[
  \Gamma(\0)\iso\1
  \qqand
  \Gamma(p+q)\iso\Gamma(p)\times\Gamma(q).
  \]
\end{proposition}

\index{section|)}

\begin{exercise}\index{colimit!in $\poly$}
  Prove \cref{prop.gamma_pres_coproduct}.
  \begin{solution}
    The functor $\Gamma$ is defined as the hom-functor $\poly(-,\yon)\colon\poly\to\smset\op$, which exhibits the universal property of colimits by sending colimits in $\poly$ to limits in $\smset$.
    Hence \cref{prop.gamma_pres_coproduct} follows from \cref{prop.poly_coprods}.
    More explicitly, $\Gamma(\0)=\poly(\0,\yon)\iso\1$ since $\0$ is initial in $\poly$, and $\Gamma(p+q)=\poly(p+q,\yon)\iso\poly(p+q,\yon)=\Gamma(p)\times\Gamma(q)$ since $+$ gives coproducts in $\poly$.
  \end{solution}
\end{exercise}

\index{functor!of sections, $\Gamma$}
\index{section|)}
%\begin{remark}
%  The sections functor $\Gamma\colon\poly\to\smset\op$ is also normal lax monoidal in the sense that there are canonical functions
%  \[
%  \1\cong\Gamma(\yon)
%  \qqand
%  \Gamma(p)\times\Gamma(q)\to\Gamma(p\otimes q)
%  \]
%  satisfying certain well-known laws. But we won't need this, so we omit its proof.
%\end{remark}

We conclude this section by discussing lenses between monomials, which arise in functional programming.

\begin{example}[Lenses between monomials are bimorphic lenses] \label{subsec.poly.cat.morph.bimorphic-lens}
  Lenses whose domains and codomains are both monomials are especially simple to write down, because they can be characterized as a pair of (standard, not dependent) functions that are independent of each other, as follows.

  \index{lens!bimorphic|see{lens, between monomials}}\index{lens!between monomials}

  Given $I,J,A,B\in\smset$, a lens $f\colon I\yon^A\to J\yon^B$ is determined by an on-positions function $f_\1\colon I\to J$ and an on-directions map: for each $i\in I$, an on-directions function $f^\sharp_i\colon B\to A$.
  But the data of such an on-directions map may be repackaged as a single function $f^\sharp\colon I\times B\to A$.
  We can do this because every position in $I$ has the same direction-set $A$, and every position in $J$ has the same direction-set $B$.

  In functional programming, such a pair of functions is called a \emph{bimorphic lens}, or a \emph{lens} for short.
  In categorical terms, we may say that the monomials in $\poly$ span a full subcategory of $\poly$ equivalent to \emph{the category of bimorphic lenses}, defined in \cite{hedges2018limits} (here the category is named after its morphisms rather than its objects).
  When such a lens arises in functional programming, the two functions that comprise it are given special names:
  \begin{equation}\label{eqn.bimorphic_lens}
    \begin{aligned}
      \lensget\coloneqq f_\1 &\colon J\to I\\
      \lensput\coloneqq f^\sharp &\colon I\times B\to A
    \end{aligned}
  \end{equation}
  Each position $i\in I$ \emph{gets} a position $f_\1i\in J$ and \emph{puts} each direction $b\in B$ back to a direction $f^\sharp(i,b)\in A$.

  So a natural transformation between two monomial functors is a bimorphic lens.
  Then a natural transformation between two polynomial functors is a more general kind of lens: a \emph{dependent} lens, where the direction-sets depend on the positions.
  Favoring the dependent version, we call these natural transformations \emph{lenses}.
\end{example}

\begin{example}[Very well-behaved lenses] \label{ex.lens_get_put}\index{lens!very well-behaved}
  Consider the monomial $S\yon^S$.
  Its position-set is $S$, and its direction-set at each position $s\in S$ is again just $S$.
  We could think of each direction as pointing to the `next' position to move to.
  We will start to formalize this idea in \cref{ex.do_nothing} and continue this work throughout the following chapters.

  Then here is one way we can think of a lens $f\colon S\yon^S\to T\yon^T$.
  Say that Otto takes positions in $S$, while Tim takes positions in $T$.
  Tim will act as Otto's proxy as follows.
  Tim will model Otto's position via the on-positions function $S\to T$ of $f$: if Otto is at position $s\in S$, then Tim will be at position $f_\1 s\in T$.
  On the other hand, Otto will take his directions from Tim via the on-directions map $S\times T\to S$ of $f$: if Tim follows the direction $t'\in T$, then Otto will head from his current position $s\in S$ in the direction $f^\sharp(s,t')\in S$.
  We interpret these directions as new positions for Otto and Tim to move to.
  So as Otto moves through the positions in $S$, he is both modeled and directed by Tim moving through the positions in $T$.

  With this setup, there are three conditions that we might expect the lens $f\colon S\yon^S\to T\yon^T$ to satisfy:
  \begin{enumerate}
    \item With Otto at $s\in S$, if Tim stays put at $f_\1s$ (i.e.\ the direction he selects at $f_\1s$ is still $f_\1s$), then Otto should stay put at $s$ (i.e.\ the direction he selects at $s$ is still $s$):
    \[
      f^\sharp(s,f_\1s)=s.
    \]

    \item Once Tim moves to $t\in T$ and Otto moves from $s\in S$ accordingly, Tim's new position should model Otto's new position:
    \[
      f_\1(f^\sharp(s,t))=t.
    \]

    \item If Tim moves to $t$, then to $t'$, Otto should end up at the same position as where he would have ended up if Tim had moved directly to $t'$ in the first place:
    \[
      f^\sharp(f^\sharp(s,t),t')=f^\sharp(s,t')
    \]
  \end{enumerate}
  Such a lens is known to functional programmers as a \emph{very well-behaved lens}; the three conditions above are its \emph{lens laws}.
  We will see these conditions emerge from more general theory in \cref{ex.very_well_behaved_lenses}.\index{lens!laws}
\end{example}

%-------- Section --------%
\section[Translating between natural transformations and lenses]{Translating between natural transformations and lenses%
  \sectionmark{Translating lenses}}
\sectionmark{Translating lenses}
\label{subsec.poly.cat.morph.translate}

\index{lens!as natural transformation|(}

We now know that we can specify a morphism $p\to q$ in $\poly$ in one of two ways:
\begin{itemize}
    \item in the language of functors, by specifying a natural transformation $p \to q$, i.e.\ for each $X\in\smset$, a function $p(X)\to q(X)$ such that naturality squares commute; or
    \item in the language of positions and directions, by specifying a lens $p\to q$, i.e.\ a function $f_\1 \colon p(\1) \to q(\1)$ and, for each $i \in p(\1)$, a function $f^\sharp_i \colon q[f_\1i] \to p[i]$.
\end{itemize}
But how are these two formulations related?
Given the data of a lens and that of a natural transformation between polynomials, how can we tell if they correspond to the same morphism?
We want to be able to translate between these two languages.

Our Rosetta Stone turns out to be the proof of the Yoneda lemma.
The lemma itself is the crux of the proof of \cref{prop.lens-prod-sum}, which states that these two formulations of morphisms between polynomials are equivalent; so unraveling these proofs reveals the translation we seek.

\index{Yoneda lemma}

\begin{proposition} \label{prop.morph_arena_to_func}
Given $p,q\in\poly$, let $f_\1\colon p(\1)\to q(\1)$ be a function between their position-sets (like an on-positions function) and $f^\sharp\colon q[f_\1(-)]\to p[-]$ be a natural transformation whose components are functions between their direction-sets (like an on-directions map).
Then the isomorphism in \eqref{eqn.main_formula} identifies $(f_\1,f^\sharp)$ with the natural transformation $f\colon p\to q$ whose $X$-component $f_X\colon p(X)\to q(X)$ for $X\in\smset$ sends each
\[
    (i,g)\in\sum_{i\in p(\1)} X^{p[i]}\iso p(X)
\]
with $i\in p(\1)$ and $g\colon p[i]\to X$ to
\[
    (f_\1i,f^\sharp_i\then g)\in\sum_{j\in q(\1)}X^{q[j]}\iso q(X).
\]
\end{proposition}
\begin{proof}\index{Yoneda lemma}
As an element of the product over $I$ on the right hand side of \eqref{eqn.main_formula}, the pair $(f_\1,f^\sharp)$ is equivalently an $I$-indexed family of pairs $((f_\1i,f^\sharp_i))_{i\in I}$, where each pair $(f_\1i,f^\sharp_i)$ is an element of
\[
    \sum_{j\in q(\1)}p[i]^{q[j]}\iso q(p[i]).
\]
By the Yoneda lemma (\cref{lemma.yoneda}), we have an isomorphism $q(p[i])\iso\poly(\yon^{p[i]}, q)$; and by the proof of the Yoneda lemma, this isomorphism sends $(f_\1i, f^\sharp_i)$ to the natural transformation $f^i\colon\yon^{p[i]} \to q$ whose $X$-component is the function $f^i_X\colon X^{p[i]}\to q(X)$ given by sending $g\colon p[i]\to X$ to
\begin{align*}
    q(g)(f_\1i, f^\sharp_i) &= \left(\sum_{j \in q(\1)} g^{q[j]}\right)(f_\1i, f^\sharp_i) \tag{\cref{cor.sum_prod_set_endofuncs}} \\
    &= \left(f_\1i, g^{q[f_\1i]}(f^\sharp_i)\right) \tag{\cref{def.sum-prod-func} and \cref{exc.sum-prod-func}} \\
    &= (f_\1i, f^\sharp_i\then g) \tag{\cref{def.representable}}.
\end{align*}
Then the $p(\1)$-indexed family of natural transformations $(f^i)_{i\in p(\1)}$ is an element of
\[
  \prod_{i\in p(\1)}\poly(\yon^{p[i]}, q)\iso\poly\left(\sum_{i\in p(\1)}\yon^{p[i]},q\right),
\]
where the isomorphism follows from the universal property of coproducts, as in the proof of \cref{prop.lens-prod-sum}.
Unwinding this isomorphism, we find that $(f^i)_{i\in I}$ corresponds to the natural transformation $f$ from $\sum_{i\in p(\1)}\yon^{p[i]}\iso p$ to $q$ that we desire.
\end{proof}

\begin{example} \label{ex.morph-corolla-with-labels}
Let us return once more to the polynomials $p \coloneqq \yon^\3 + \2\yon$ and $q \coloneqq \yon^\4 + \yon^\2 + \2$ from \cref{ex.practice_with_poly_morphisms} and the lens $f \colon p \to q$ depicted below:
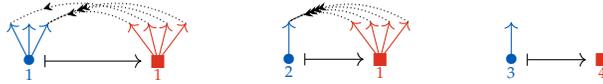
\[
\begin{tikzpicture}
	\node (p1) {
	\begin{tikzpicture}[trees, sibling distance=2.5mm]
    \node[my-blue, "\color{my-blue}\tiny 1" below] (1) {$\bullet$}
      child[my-blue] {coordinate (11)}
      child[my-blue] {coordinate (12)}
      child[my-blue] {coordinate (13)};
    \node[right=1.5 of 1, my-red, "\color{my-red}\tiny 1" below] (2) {\small$\blacksquare$}
      child[my-red] {coordinate (21)}
      child[my-red] {coordinate (22)}
      child[my-red] {coordinate (23)}
      child[my-red] {coordinate (24)};
    \draw[|->, shorten <= 3pt, shorten >= 3pt] (1) -- (2);
    \begin{scope}[densely dotted, bend right, decoration={markings, mark=at position 0.75 with \arrow{stealth}}]
      \draw[postaction={decorate}] (21) to (13);
      \draw[postaction={decorate}] (22) to (11);
      \draw[postaction={decorate}] (23) to (13);
      \draw[postaction={decorate}] (24) to (13);
    \end{scope}
  \end{tikzpicture}
	};
	\node (p2) [right=1 of p1] {
	\begin{tikzpicture}[trees, sibling distance=2.5mm]
    \node[my-blue, "\color{my-blue}\tiny 2" below] (1) {$\bullet$}
      child[my-blue] {coordinate (11)};
    \node[right=of 1, my-red, "\color{my-red}\tiny 1" below] (2) {\small$\blacksquare$}
      child[my-red] {coordinate (21)}
      child[my-red] {coordinate (22)}
      child[my-red] {coordinate (23)}
      child[my-red] {coordinate (24)};
    \draw[|->, shorten <= 3pt, shorten >= 3pt] (1) -- (2);
    \begin{scope}[densely dotted, bend right, decoration={markings, mark=at position 0.75 with \arrow{stealth}}]
      \draw[postaction={decorate}] (21) to (11);
      \draw[postaction={decorate}] (22) to (11);
      \draw[postaction={decorate}] (23) to (11);
      \draw[postaction={decorate}] (24) to (11);
    \end{scope}
  \end{tikzpicture}
	};
	\node (p3) [below right=-10mm and 1 of p2] {
	\begin{tikzpicture}[trees, sibling distance=2.5mm]
    \node[my-blue, "\color{my-blue}\tiny 3" below] (1) {$\bullet$}
      child[my-blue] {};
    \node[right=of 1, my-red, "\color{my-red}\tiny 4" below] (2) {\small$\blacksquare$}
		;
    \draw[|->, shorten <= 3pt, shorten >= 3pt] (1) -- (2);
  \end{tikzpicture}
	};
\end{tikzpicture}
\]
Fix a set $X \coloneqq \{\const{a},\const{b},\const{c},\const{d},\const{e}\}$.
When viewed as a natural transformation, $f$ has as its $X$-component a function $f_X \colon p(X) \to q(X)$.
In other words, for each element of $p(X)$, the lens $f$ should tell us how to obtain an element of $q(X)$.

We saw in \cref{ex.corolla-apply-poly} that each $(i,g)\in p(X)$ may be drawn as a $p$-corolla (corresponding to $i$) whose leaves are labeled with elements of $X$ (according to $g$).
For example, here we draw $(1,g)\in p(X)$, where $g\colon p[1]\to X$ is given by $1\mapsto \const{c}, 2 \mapsto \const{e},$ and $3 \mapsto \const{a}$:
\begin{equation} \label{diag.corolla-apply-elt}
\begin{tikzpicture}[trees, sibling distance=2.5mm]
    \node[my-blue, "\color{my-blue}\tiny 1" below] (1) {$\bullet$}
      child[my-blue] {node {$\const{c}$}}
      child[my-blue] {node {$\const{e}$}}
      child[my-blue] {node {$\const{a}$}};
\end{tikzpicture}
\end{equation}
Similarly, each element of $q(X)$ can be drawn as a $q$-corolla whose leaves are labeled with elements of $X$.
So what element of $q(X)$ is $f_X(1, g)$?

\cref{prop.morph_arena_to_func} tells us that $f_X(1, g)=(f_\1(1), f^\sharp_1\then g)$, so we need only focus on the behavior of $f$ at $p$-position $1$:
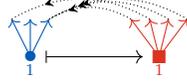
\[
\begin{tikzpicture}[trees, sibling distance=2.5mm]
    \node[my-blue, "\color{my-blue}\tiny 1" below] (1) {$\bullet$}
      child[my-blue] {coordinate (11)}
      child[my-blue] {coordinate (12)}
      child[my-blue] {coordinate (13)};
    \node[right=1.5 of 1, my-red, "\color{my-red}\tiny 1" below] (2) {\small$\blacksquare$}
      child[my-red] {coordinate (21)}
      child[my-red] {coordinate (22)}
      child[my-red] {coordinate (23)}
      child[my-red] {coordinate (24)};
    \draw[|->, shorten <= 3pt, shorten >= 3pt] (1) -- (2);
    \begin{scope}[densely dotted, bend right, decoration={markings, mark=at position 0.75 with \arrow{stealth}}]
      \draw[postaction={decorate}] (21) to (13);
      \draw[postaction={decorate}] (22) to (11);
      \draw[postaction={decorate}] (23) to (13);
      \draw[postaction={decorate}] (24) to (13);
    \end{scope}
\end{tikzpicture}
\]
To draw $(f_\1(1), f^\sharp_1\then g)$, we first draw the $q$-corolla corresponding to $f_\1(1)$, the corolla on the right hand side above.
Then we label each leaf of that corolla by following the arrow from that leaf to a $p[1]$-leaf, and use the label there from \eqref{diag.corolla-apply-elt} (as prescribed by $g$).
So $f_X(1, g)$ looks like
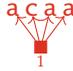
\[
\begin{tikzpicture}[trees, sibling distance=2.5mm]
    \node[my-red, "\color{my-red}\tiny 1" below] (1) {\small$\blacksquare$}
      child[my-red] {node {$\const{a}$}}
      child[my-red] {node {$\const{c}$}}
      child[my-red] {node {$\const{a}$}}
      child[my-red] {node {$\const{a}$}};
\end{tikzpicture}
\]

\end{example}

\cref{prop.morph_arena_to_func} lets us translate from lenses to natural transformations.
The following corollary tells us how to go in the other direction.
In particular, it justifies the notation $f_\1$ for the on-positions function of $f$: it is the $\1$-component of $f$ as a natural transformation.

\begin{corollary} \label{cor.morph_func_to_arena}
Let $p$ and $q$ be polynomial functors, and let $f \colon p \to q$ be a natural transformation between them.
Then the isomorphism in \eqref{eqn.main_formula} sends $f$ to the lens whose on-positions function $f_\1\colon p(\1)\to q(\1)$ is the $\1$-component of $f$ and whose on-directions map $f^\sharp\colon q[f_\1(-)]\to p[-]$ satisfies, for all $i\in p(\1)$,
\[
    (f_\1i, f^\sharp_i) = f_{p[i]}(i, \id_{p[i]}).
\]
\end{corollary}
\begin{proof}
By \cref{prop.morph_arena_to_func}, the $\1$-component of $f$ is a function $p(\1)\to q(\1)$ sending every $i \in p(\1)$ to $f_\1i \in q(\1)$, so the on-positions function $f_\1$ is indeed equal to the $\1$-component of $f$.
Moreover, for each $i\in p(\1)$, the $p[i]$-component $f_{p[i]} \colon p(p[i]) \to q(p[i])$ of $f$ sends $(i,\id_{p[i]})\in p(p[i])$ to $(f_\1i, f^\sharp_i \then \id_{p[i]}) = (f_\1i, f^\sharp_i)$.
\end{proof}

\index{lens!as natural transformation|)}

%-------- Section --------%
\section[Identity lenses and lens composition]{Identity lenses and lens composition%
  \sectionmark{Identity lenses \& lens composition}}
\sectionmark{Identity lenses \& lens composition} \label{subsec.poly.cat.morph.id-comp}

Thus far, we have seen how the category $\poly$ of polynomial functors and natural transformations can be identified with the category of indexed families of sets and lenses.
But in order to actually discuss the latter category, we need to be able to give identity lenses and describe how lenses compose.
We can do so by translating between lenses and natural transformations.

For instance, given a polynomial $p$, its identity lens should correspond to the identity natural transformation of $p$ as a functor.

\index{lens!identity}

\begin{exercise}[Identity lenses] \label{exc.arena_morph_id}
For $p\in\poly$, let $\id_p \colon p \to p$ be its identity natural transformation, whose $X$-component $(\id_p)_X\colon p(X)\to p(X)$ for $X\in\smset$ is the identity function on $p(X)$; that is, $(\id_p)_X=\id_{p(X)}$.

Use \cref{cor.morph_func_to_arena} to demonstrate that the on-positions function $(\id_p)_\1\colon p(\1)\to p(\1)$ as well as, for each $i\in p(\1)$, the on-directions function $(\id_p)^\sharp_i\colon p[(\id_p)_\1i]\to p[i]$ are all identity functions.
\begin{solution}
Fix $i\in p(\1)$.
Since the $p[i]$-component $(\id_p)_{p[i]}$ of $\id_p$ is the identity function on $p(p[i])$, by \cref{cor.morph_func_to_arena},
\[
    ((\id_p)_\1i, (\id_p)^\sharp_i) = (\id_p)_{p[i]}(i,\id_{p[i]}) = (i,\id_{p[i]}).
\]
Hence the on-positions function $(\id_p)_\1\colon p(\1)\to p(\1)$ maps every $i\in p(\1)$ to itself, so it is an identity function; and each on-directions function $(\id_p)^\sharp_i\colon p[i]\to p[i]$ is equal to $\id_{p[i]}$.
\end{solution}
\end{exercise}

Similarly, we may infer how two lenses compose by translating them to natural transformations, composing those, then translating back to lenses.

\index{lens!composition of}

\begin{exercise}[Composing lenses] \label{exc.arena_morph_comp}
For $p,q,r\in\poly$, let $f\colon p\to q$ and $g\colon q\to r$ be natural transformations, and let $h\coloneqq f\then g$ be their composite, whose $X$-component $h_X\colon p(X)\to r(X)$ for $X\in\smset$ is the composite of the $X$-components of $f$ and $g$; that is, $h_X=f_X\then g_X$.

Use \cref{prop.morph_arena_to_func} and \cref{cor.morph_func_to_arena} to show that the on-positions function $h_\1\colon p(\1)\to r(\1)$ of $h$ is given by $h_\1=f_\1\then g_\1$, while the on-directions function $h^\sharp_i$ of $h$ for $i\in p(\1)$ is given by $h^\sharp_i=g^\sharp_{f_\1i}\then f^\sharp_i$.
\begin{solution}
Fix $i\in p(\1)$.
By \cref{prop.morph_arena_to_func} and \cref{cor.morph_func_to_arena},
\begin{align*}
    (h_\1i, h^\sharp_i) &= h_{p[i]}(i, \id_{p[i]}) \tag{\cref{cor.morph_func_to_arena}} \\
    &= g_{p[i]}(f_{p[i]}(i, \id_{p[i]})) \tag{$h = f \then g$} \\
    &= g_{p[i]}(f_\1i, f^\sharp_i) \tag{\cref{cor.morph_func_to_arena}} \\
    &= (g_\1f_\1i, g^\sharp_{f_\1i} \then f^\sharp_i). \tag{\cref{prop.morph_arena_to_func}}
\end{align*}
\end{solution}
\end{exercise}

% TODO: turn following to prop?
% The following corollary about interpreting commutative diagrams in $\poly$ is immediate from the preceding exercise.

The following proposition, a restatement of the previous exercise, allows us to interpret commutative diagrams of polynomials in $\poly$ in terms of commutative diagrams of their position- and direction-sets in $\smset$.

\begin{proposition} \label{prop.comm_poly}
Given $p,q,r\in\poly$ and lenses $f\colon p\to q, g\colon q\to r,$ and $h\colon p\to r$, the diagram
\[
\begin{tikzcd}
    p \ar[r, "f"] \ar[dr, "h"'] & q \ar[d, "g"] \\
    & r
\end{tikzcd}
\]
commutes in $\poly$ if and only if the forward on-positions diagram
\[
\begin{tikzcd}
    p(\1) \ar[r, "f_\1"] \ar[dr, "h_\1"'] & q(\1) \ar[d, "g_\1"] \\
    & r(\1)
\end{tikzcd}
\]
commutes in $\smset$ and, for each $i \in p(\1)$, the backward on-directions diagram
\[
\begin{tikzcd}
    p[i] & q[f_\1i] \ar[l, "f^\sharp_i"'] \\
    & r[h_\1i] \ar[u, "g^\sharp_{f_\1i}"'] \ar[ul, "h^\sharp_i"]
\end{tikzcd}
\]
commutes in $\smset$.
\end{proposition}

We can use this fact to determine whether a given diagram in $\poly$ commutes, as in the following exercise.

\begin{exercise}
Using \cref{prop.comm_poly}, verify explicitly that, for $p, q \in \poly$, the polynomial $p+q$ given by the binary sum of $p$ and $q$ satisfies the universal property of the coproduct of $p$ and $q$.
That is, provide lenses $\iota \colon p \to p + q$ and $\kappa \colon q \to p + q$, then show that for any other polynomial $r$ equipped with lenses $f \colon p \to r$ and $g \colon q \to r$, there exists a unique lens $h\colon p+q\to r$ (shown dashed) making the following diagram commute:
\begin{equation} \label{eqn.coprod_univ_prop}
\begin{tikzcd}
	p \ar[r, "\iota"] \ar[dr, "f"'] &
	p + q \ar[d, "h", dashed] &
	q \ar[l, "\kappa"'] \ar[dl, "g"] \\
	& r
\end{tikzcd}
\end{equation}
\begin{solution}
We provide lenses $\iota\colon p\to p+q$ and $\kappa\colon q\to p+q$ as follows.
On positions, they are the canonical inclusions $\iota_\1\colon p(\1)\to p(\1)+q(\1)$ and $\kappa_\1\colon q(\1)\to p(\1)+q(\1)$; on directions, they are identities.
To show that $p+q$ equipped with $\iota$ and $\kappa$ satisfies the universal property of the coproduct, we apply \cref{prop.comm_poly}.
In order for \eqref{eqn.coprod_univ_prop} to commute, it must commute on positions---that is, the following diagram of sets must commute:
\begin{equation} \label{eqn.coprod_univ_prop_pos}
\begin{tikzcd}
	p(\1) \ar[r, "\iota_\1"] \ar[dr, "f_\1"'] &
	p(\1) + q(\1) \ar[d, "h_\1", dashed] &
	q(\1) \ar[l, "\kappa_\1"'] \ar[dl, "g_\1"] \\
	& r(\1)
\end{tikzcd}
\end{equation}
But since $p(1)+q(\1)$ along with the inclusions $\iota_\1$ and $\kappa_\1$ form the coproduct of $p(\1)$ and $q(\1)$ in $\smset$, there exists a unique $h_\1$ for which \eqref{eqn.coprod_univ_prop_pos} commutes.
Hence $h$ is uniquely characterized on positions.
In particular, it must send each $(1,i) \in p(\1)+q(\1)$ with $i \in p(\1)$ to $f_\1i$ and each $(2,j) \in p(\1)+q(\1)$ with $j \in q(\1)$ to $g_\1j$.

Moreover, if \eqref{eqn.coprod_univ_prop} is to commute on directions, then for every $i\in p(\1)$ and $j \in q(\1)$, the following diagrams of sets must commute:
\begin{equation} \label{eqn.coprod_univ_prop_dir}
\begin{tikzcd}[sep=large]
	p[i] & (p+q)[(1,i)] \ar[l, "\iota^\sharp_i"'] & (p+q)[(2,j)] \ar[r, "\kappa^\sharp_j"] & q[j] \\
	& r[f_\1i] \ar[ul, "f^\sharp_i"] \ar[u, "h^\sharp_{(1,\,i)}"', dashed] & r[g_\1j] \ar[u, "h^\sharp_{(2,\,j)}", dashed] \ar[ur, "g^\sharp_j"']
\end{tikzcd}
\end{equation}
But $(p+q)[(1,i)] \iso p[i]$ and $\iota^\sharp_i$ is the identity, so we must have $h^\sharp_{(1,\,i)} = f^\sharp_i$.
Similarly, $(p+q)[(2,j)] \iso q[j]$ and $\kappa^\sharp_j$ is the identity, so we must have $h^\sharp_{(2,\,j)} = g^\sharp_j$.
Hence $h$ is also uniquely characterized on directions, so it is unique overall.
Moreover, we have shown that we can define $h$ on positions so that \eqref{eqn.coprod_univ_prop_pos} commutes, and that we can define $h$ on directions such that the diagrams in \eqref{eqn.coprod_univ_prop_dir} commute.
As the commutativity of the diagrams in \eqref{eqn.coprod_univ_prop_pos} and \eqref{eqn.coprod_univ_prop_dir} together imply the commutativity of \eqref{eqn.coprod_univ_prop}, it follows that there exists $h$ for which \eqref{eqn.coprod_univ_prop} commutes.
\end{solution}
\end{exercise}

Now that we know how lens composition works in $\poly$, we have a better handle on how it behaves categorically.
For instance, we can verify functoriality in $\poly$, as in the following exercise.

\begin{exercise}[A functor $\Cat{Top}\to\poly$] \label{exc.top_poly_func}\index{topological space}
This exercise is for those who know what topological spaces and continuous maps are. It will not be used again in this book.
\begin{enumerate}
	\item Given a topological space $X$, define a polynomial $p_X$ whose positions are the points in $X$ and whose directions at $x\in X$ are the open neighborhoods of $x$.
  That is,
  \[
    p_X\coloneqq\sum_{x \in X}\yon^{\{ U\ss X \mid x\in U, \, U\text{ open} \}}.
  \]
  Given a continuous map $f\colon X\to Y$, define a lens $p_X\to p_Y$ either by writing down its formula or drawing it in polyboxes.
	\item Show that the assignment above defines a functor $\Cat{Top}\to\poly$.
	\item Is this functor full? Is it faithful?
\qedhere
\end{enumerate}
\begin{solution}
\begin{enumerate}
	\item \label{exc.top_poly_func.morphs} Given a continuous map $f \colon X \to Y$, we define a lens $p_f \colon p_X \to p_Y$ as follows.
	The on-positions function is just $f$; then for each $p_X$-position $x\in X$, the on-directions function $(p_f)^\sharp_x\colon p_Y[f(x)]\to p_X[x]$ sends each open neighborhood $U$ of $f(x)$ to $f\inv(U)$, which we know is an open neighborhood of $x$ because $f$ is continuous.
  The polybox picture for $p_f$ is as follows:
  \[
  \begin{tikzpicture}[polybox, mapstos]
    \node[poly, dom, "$p$" below] (p) {$f\inv(U)$\at$x\vphantom{f}$};

    \node[poly, cod, right=of p, "$q$" below] (q) {$U\vphantom{f\inv}$\at$f(x)$};

    \draw (p_pos) -- node[below] {$f$} (q_pos);
    \draw (q_dir) -- node[above] {$f\inv$} (p_dir);
  \end{tikzpicture}
  \]

	\item To show that $p_X$ is functorial in $X$, it suffices to show that sending continuous maps $f\colon X\to Y$ to their induced lenses $p_f\colon p_X\to p_Y$ preserves identities and composition.
	First, we show for $X\in\Cat{Top}$ that the lens $p_{\id_X}$ is an identity.
	By \cref{exc.top_poly_func.morphs}, the on-positions function of $p_{\id_X}$ is $\id_X$, and for each $x\in X$ the on-directions function $(p_f)^\sharp_x\colon p_X[x]\to p_X[x]$ sends $U\in p_X[x]$ to $(\id_X)\inv(U) = U$.
	Hence $p_{\id_X}$ is the identity on both positions and directions; it follows from \cref{exc.arena_morph_id} that $p_{\id_X}$ is the identity lens.

	We now show for $X,Y,Z\in\Cat{Top}$ and continuous maps $f\colon X\to Y$ and $g\colon Y\to Z$ that $p_f\then p_g = p_{f\then g}$.
	By \cref{exc.top_poly_func.morphs} and \cref{exc.arena_morph_comp}, the on-positions functions of both $p_f\then p_g$ and $p_{f\then g}$ are equal to $f\then g$, so it suffices to show for each $x\in X$ that
	\[
	    (p_{f \then g})^\sharp_x = (p_g)^\sharp_{f(x)} \then (p_f)^\sharp_x.
	\]
	By \cref{exc.top_poly_func.morphs}, the left hand side sends each $U\in p_Z[g(f(x))]$ to $(f \then g)\inv(U)$, while the right hand side sends $U$ to $f\inv(g\inv(U))$; by elementary set theory, these sets are equal.

	\item The functor is not full.
	Consider the spaces $X\coloneqq\2$ with the indiscrete topology (i.e.\ the only open sets are the empty set and $X$) and $Y\coloneqq\2$ with the discrete topology (i.e.\ all subsets are open).
	Then $p_X\iso\2\yon$ (each point in $X$ has exactly one open neighborhood: the entire space $X$) and $p_Y\iso\2\yon^\2$ (each point in $Y$ has exactly two open neighborhoods: a singleton set and $Y$ itself), so our functor induces a map from the set of continuous functions $X\to Y$ to the set of lenses $\2\yon\to\2\yon^\2$.
	We claim this map is not surjective: in particular, consider the lens $h\colon\2\yon\to\2\yon^\2$ that is the identity on positions (and uniquely defined on directions).
	Then a continuous function $f\colon X\to Y$ that our functor sends to $h$ must also be the identity on the underlying sets of $X$ and $Y$.
	But such a function cannot be continuous: a singleton subset of $Y$ is open, but its preimage under $f$ is a singleton subset of $X$ and therefore not open.
	So our functor sends no continuous function $X\to Y$ to $h$ and therefore is not full.
	The functor is, however, faithful: given spaces $X$ and $Y$ and continuous function $f\colon X\to Y$, we can uniquely recover $f$ from $p_f$ by taking its on-positions function $(p_f)_\1=f$.
\end{enumerate}
\end{solution}
\end{exercise}

%-------- Section --------%
\section{Polybox pictures of lens composition} \label{sec.poly.cat.}

\index{polybox!for lens composition|(}

Given lenses $f\colon p\to q$ and $g\colon q\to r$, we can piece their polyboxes together to form polyboxes for their composite, $f\then g\colon p\to r$:
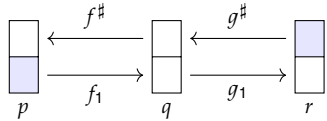
\[
\begin{tikzpicture}[polybox, tos]
  \node[poly, dom, "$p$" below] (p) {};

  \node[poly, right=of p, "$q$" below] (q) {};

  \node[poly, cod, right=of q, "$r$" below] (r) {};

  \draw (p_pos) -- node[below] {$f_\1$} (q_pos);
  \draw (q_dir) -- node[above] {$f^\sharp$} (p_dir);

  \draw (q_pos) -- node[below] {$g_\1$} (r_pos);
  \draw (r_dir) -- node[above] {$g^\sharp$} (q_dir);
\end{tikzpicture}
\]
The position box for $q$, which would be blue as part of the polyboxes for $g\colon q\to r$ alone, is instead filled in via $f_\1$; similarly, the direction box for $q$, which would be blue as part of the polyboxes for just $f\colon p\to q$, is filled in via $g^\sharp$.
This forms a spreadsheet-filling protocol that acts as the polyboxes for $f\then g$.

\index{lens!spreadsheet depiction}

As we follow the arrows from left to right and up and left again, take care to note that the arrow $g^\sharp$ depends not only on the direction box of $r$, but also the position box of $q$ that came before it.
Similarly, $f^\sharp$ depends on both the position box of $p$ and the direction box of $q$.
On the other hand, the arrow $g_\1$ depends only on the position box of $q$, and not the position box of $p$ that came before it: $g_\1$ is the on-positions function for a lens $q\to r$ and therefore depends only on its domain $q$. (Of course, changing the position box of $p$ may change the position box of $q$ via $f_\1$, thus indirectly affecting what $g_\1$ enters in the position box for $r$; we mean that if the position box of $p$ changes but the position box of $q$ does not, $g_\1$ will not change the position box of $r$.)
Similarly, $g^\sharp$ does not depend on the position box of $p$, and $f^\sharp$ does not depend on either box of $r$.
The key is to let each arrow depend on exactly the boxes that come before it in the domain and codomain of the lens that the arrow is a part of.

If we have another lens $h\colon p\to r$, we can interpret the equation $f\then g = h$ by filling in their polyboxes and comparing them:
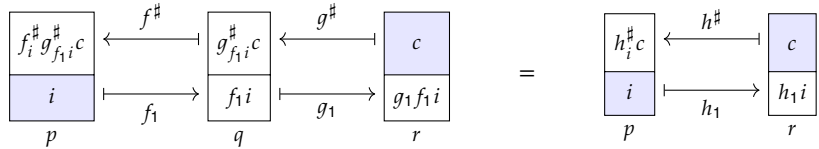
\[
\begin{tikzpicture}
  \node (1) {
    \begin{tikzpicture}[polybox, mapstos]
      \node[poly, dom, "$p$" below] (p) {$f^\sharp_ig^\sharp_{f_\1i}c$\at$i\vphantom{f_\1}$};

      \node[poly, right=of p, "$q$" below] (q) {$g^\sharp_{f_\1i}c$\at$f_\1i$};

      \node[poly, cod, right=of q, "$r$" below] (r) {$c\vphantom{f^\sharp_ig^\sharp_{f_\1i}c}$\at$g_\1f_\1i$};

      \draw (p_pos) -- node[below] {$f_\1$} (q_pos);
      \draw (q_dir) -- node[above] {$f^\sharp$} (p_dir);

      \draw (q_pos) -- node[below] {$g_\1$} (r_pos);
      \draw (r_dir) -- node[above] {$g^\sharp$} (q_dir);
    \end{tikzpicture}
  };
  \node[right=1.8 of 1] (2) {
    \begin{tikzpicture}[polybox, mapstos]
      \node[poly, dom, "$p$" below] (p) {$h^\sharp_ic$\at$i\vphantom{h_\1}$};
      \node[poly, cod, "$r$" below, right=of p] (q) {$\vphantom{h^\sharp_i}c$\at$h_\1i$};
      \draw (p_pos) -- node[below] {$h_\1$} (q_pos);
      \draw (q_dir) -- node[above] {$h^\sharp$} (p_dir);
    \end{tikzpicture}
  };
  \node at ($(1.east)!.5!(2.west)$) {=};
\end{tikzpicture}
\]
Here we have filled the blue boxes on either side with the same entries.
Then if we match up the uncolored boxes in the domain and codomain on either side, we can read off the equations
\[
  g_\1f_\1i = h_\1i \qqand f^\sharp_ig^\sharp_{f_\1i}c = h^\sharp_ic
\]
for every $p$-position $i$ and $r[h_\1i]$-direction $c$, which agrees with \cref{exc.arena_morph_comp} and \cref{prop.comm_poly}.
Throughout this book, we will often read off equalities of positions and directions from polybox pictures of lens equations in this way.

Note that there is redundancy in the above polybox picture: we have filled in all the boxes for clarity, but their entries are determined by the entries in the blue boxes and the labels on the arrows.
So we may omit the entries in the uncolored boxes without losing information, leaving the reader to fill in the blanks:
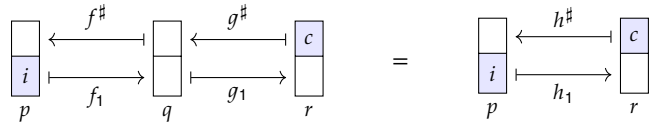
\[
\begin{tikzpicture}
  \node (1) {
    \begin{tikzpicture}[polybox, mapstos]
      \node[poly, dom, "$p$" below] (p) {$\phantom{c}$\at$i$};

      \node[poly, right=of p, "$q$" below] (q) {$\phantom{c}$\at$\phantom{i}$};

      \node[poly, cod, right=of q, "$r$" below] (r) {$c$\at$\phantom{i}$};

      \draw (p_pos) -- node[below] {$f_\1$} (q_pos);
      \draw (q_dir) -- node[above] {$f^\sharp$} (p_dir);

      \draw (q_pos) -- node[below] {$g_\1$} (r_pos);
      \draw (r_dir) -- node[above] {$g^\sharp$} (q_dir);
    \end{tikzpicture}
  };
  \node[right=1.8 of 1] (2) {
    \begin{tikzpicture}[polybox, mapstos]
      \node[poly, dom, "$p$" below] (p) {$\phantom{c}$\at$i$};
      \node[poly, cod, "$r$" below, right=of p] (q) {$c$\at$\phantom{i}$};
      \draw (p_pos) -- node[below] {$h_\1$} (q_pos);
      \draw (q_dir) -- node[above] {$h^\sharp$} (p_dir);
    \end{tikzpicture}
  };
  \node at ($(1.east)!.5!(2.west)$) {=};
\end{tikzpicture}
\]

\begin{remark}
  The reader may be concerned that when working with polyboxes, we refer to ``spreadsheets'' and ``protocols'' without being rigorous about what they are or what it means to set them equal.
  We choose to elide this issue to highlight the graphical intuition rather than grinding through the details.
  This is not to say our work with polyboxes will lack rigor moving forward---if you are particularly worried, you should think of polyboxes as an alternate way to present information about indexed families of sets, dependent functions, and sum and product sets that can be systematically translated---via elementary steps, though perhaps with some laborious bookkeeping---into the more standard $\in$ and $\sum$ and $\prod$ notation we have been using thus far.

  For example, given lenses $f\colon p\to q$ and $g\colon q\to r$, the polyboxes on the left hand side of the equation above should be interpreted as the element of the set
  \[
    \prod_{i\in p(\1)}\sum_{k\in r(\1)}p[i]^{r[k]}\iso\poly(p,r)
  \]
  corresponding to the lens $p\to r$ whose on-positions function $p(\1)\to r(\1)$ is the composite of the on-positions functions $f_\1$ and $g_\1$ and whose on-directions function $r[g_\1f_\1i]\to p[i]$ at $i\in p(\1)$ is equal to the composite of the on-directions functions $g^\sharp_{f_\1i}$ and $f^\sharp_i$.
  In other words, the polyboxes represent the composite lens $f\then g$.
  But the polyboxes show how lenses pass positions and directions back and forth far more legibly than the last two sentences can.
  Throughout the rest of this book, we will see how this polybox notation provides immediate, reader-friendly computations and justifications; but all these results can be translated back into more grounded mathematical language as desired.
\end{remark}

\begin{example}[Modeling with a composite lens in polyboxes]
  By composing the lens $f\colon p\to q$ from \cref{ex.lend-return} that models the exchange of money between Caroline (modeled by $q$) and her parents (modeled by $p$) with the lens $\gamma\colon q\to\yon$ from \cref{ex.spend-section} that models how Caroline spends her money, we obtain a lens $f\then\gamma\colon p\to\yon$ that models how Caroline's parents spend their money through Caroline.
  The polybox picture of the composite lens $f\then\gamma$ is given by merging the polybox pictures of $f$ and $\gamma$:
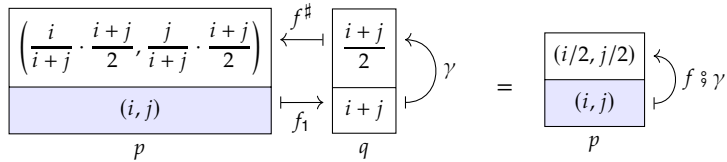
  \[
  \begin{tikzpicture}
    \node (1) {
      \begin{tikzpicture}[polybox, mapstos]
        \node[poly, dom, "$p$" below] (p) {$\left(\dfrac{i}{i+j}\cdot\dfrac{i+j}2,\dfrac{j}{i+j}\cdot\dfrac{i+j}2\right)$\at$(i,j)$};

        \node[poly, right=.8 of p, "$q$" below] (q) {$\vphantom{\left(\dfrac{j}{i+j}\dfrac{i+j}2\right)}\dfrac{i+j}2$\at$i+j$};

        \draw (p_pos) -- node[below] {$f_\1$} (q_pos);
        \draw (q_dir) -- node[above] {$f^\sharp$} (p_dir);

        \draw (q_pos) to[climb'] node[right] {$\gamma$} (q_dir);
      \end{tikzpicture}
    };
    \node[right=.8 of 1, yshift=-2mm] (2) {
      \begin{tikzpicture}[polybox, mapstos]
        \node[poly, dom, "$p$" below] (p) {$(i/2,j/2)$\at$(i,j)$};

        \draw (p_pos) to[climb'] node[right] {$f\then\gamma$} (p_dir);
      \end{tikzpicture}
    };
    \node at ($(1.east)!.5!(2.west)$) {=};
  \end{tikzpicture}
  \]
  Here $(i,j)\in p(\1)=(0,20]\times(0,20]$.
  The right hand side summarizes what happens to the parents: if the first parent gives away $i$ dollars and the second parent gives away $j$ dollars, eventually the first parent will receive $i/2$ dollars and the second parent will receive $j/2$ dollars.
  The factored left hand side describes how this happens: the parents give $i$ and $j$ dollars respectively to Caroline, who takes the $i+j$ dollars total and spends half of it.
  She then returns the remaining half to her parents, splitting the money proportionately according to the amount each parent contributed.
\end{example}

\index{polybox!for lens composition|)}
\index{lens|)}

%-------- Section --------%
\section[Symmetric monoidal products of polynomial functors]{Symmetric monoidal products of polynomial functors%
  \sectionmark{Symmetric monoidal products in $\poly$}}
\sectionmark{Symmetric monoidal products in $\poly$}
\label{sec.poly.cat.monoidal}

\index{monoidal structure}\index{polynomial functor!monoidal structure on|see{monoidal structure}}

One of the reasons $\poly$ is so versatile is that there is an abundance of monoidal structures on it.
Monoidal structures are the key ingredient to many applications of categories to real-world settings, and $\poly$ is no different in that regard.
As a bonus, if you know how to add and multiply polynomials from high school algebra, then you already know how to compute two of the monoidal products on $\poly$.

We have already seen one of these monoidal structures on $\poly$: the cocartesian monoidal structure, which gives $\poly$ its finite coproducts.
In fact, we know from \cref{prop.poly_coprods} that $\poly$ has all coproducts: they are given by an operation that looks just like addition.
It turns out $\poly$ has all products as well, giving it a cartesian monoidal structure that looks just like multiplication.

\index{monoidal structure!cartesian (products)}\index{monoidal structure!cocartesian (sums)}\index{coproduct}

\begin{proposition}\label{prop.poly_prods}\index{polynomial functor!product of polynomials}
  The category $\poly$ has arbitrary products, coinciding with products in $\smset^\smset$ given by the operation $\prod_{i \in I}$.
\end{proposition}
\begin{proof}
  Unsurprisingly, the proof is very similar to that of \cref{prop.poly_coprods}.

  By \cref{cor.sum_prod_set_endofuncs}, the category $\smset^\smset$ has arbitrary products given by $\prod_{i \in I}$.
  The full subcategory inclusion $\poly \to \smset^\smset$ reflects these products.
  It remains to show that $\poly$ is closed under the operation $\prod_{i \in I}$.\index{polynomial functor!sum of polynomials}

  By \cref{prop.set_endofunc_distrib}, $\smset^\smset$ is completely distributive.
  Hence, given polynomials $(p_i)_{i \in I}$, we can use \eqref{eqn.cat_completely_distributive} to write their product in $\smset^\smset$ as
  \begin{equation} \label{eqn.poly_prod}
    \prod_{i \in I} p_i \iso \prod_{i \in I} \sum_{j \in p_i(\1)} \yon^{p_i[j]} \iso \sum_{\bar{j} \in \prod_{i \in I} p_i(\1)} \prod_{i \in I} \yon^{p_i[\bar{j}i]} \iso \sum_{\bar{j} \in \prod_{i \in I} p_i(\1)} \yon^{\sum_{i \in I} p_i[\bar{j}i]},
  \end{equation}
  which, as a coproduct of representables, is in $\poly$.
  % We will see that $\1$ is a terminal object and that the product of $p$ and $q$ in $\poly$ is the usual product of $p$ and $q$ as polynomials. That is, if $p\coloneqq\sum_{i\in p(\1)}\yon^{p[i]}$ and $q\coloneqq\sum_{j\in q(\1)}\yon^{q[j]}$ are in standard notation, then
  % \begin{equation}\label{eqn.poly_times}
    % p\times q\cong\sum_{i\in p(\1)}\sum_{j\in q(\1)}\yon^{p[i]+q[j]}.
    % \end{equation}
  % We leave the proof as an exercise; see \cref{exc.poly_times}.
\end{proof}

\index{polynomial functor!product of polynomials}

\begin{corollary} \label{prop.poly_completely_distributive}\index{completely distributive category}
  The category $\poly$ is completely distributive.
\end{corollary}
\begin{proof}
  This is a direct consequence of the fact that $\poly$ has arbitrary (co)products coinciding with (co)products in $\smset^\smset$ (\cref{prop.poly_coprods,prop.poly_prods}) and the fact that $\smset^\smset$ itself is completely distributive (\cref{prop.set_endofunc_distrib}).
\end{proof}

The result above will allow us to apply \eqref{eqn.cat_completely_distributive}, or sometimes specifically \eqref{eqn.push_prod_sum_obj_indep}, to push $\prod$'s past $\sum$'s of polynomials whenever we so desire.

\begin{exercise}%\label{exc.poly_prod}
    Use \eqref{eqn.main_formula}
    to verify that
    \[
    \poly\!\left(q, \prod_{i \in I} p_i\right) \iso \prod_{i \in I} \poly(q, p_i)
    \]
    for all polynomials $(p_i)_{i \in I}$ and $q$, as one would expect from the universal property of products.
    \qedhere
  \begin{solution}
    Given $q \in \poly$ and $p_i \in \poly$ for each $i \in I$ for some set $I$, we use \eqref{eqn.main_formula} to verify that
      \begin{align*}
        \poly\left(q, \prod_{i \in I} p_i\right) &\iso \prod_{k \in q(\1)} \left(\prod_{i \in I} p_i\right)(q[k])
        \tag*{\eqref{eqn.main_formula}} \\
        &\iso \prod_{k \in q(\1)} \prod_{i \in I} p_i(q[k]) \\
        &\iso \prod_{i \in I} \prod_{k \in q(\1)} p_i(q[k]) \\
        &\iso \prod_{i \in I} \poly(q, p_i).
        \tag*{\eqref{eqn.main_formula}}
      \end{align*}
      % \end{enumerate}
  \end{solution}
\end{exercise}

\begin{exercise}
  Let $p_1\coloneqq\yon+\1, p_2\coloneqq\yon+\2,$ and $p_3\coloneqq\yon^\2$.
  What is $\prod_{i\in\3}p_i$ according to \eqref{eqn.poly_prod}? Is the answer what you would expect?
  \begin{solution}
    Given $p_1\coloneqq\yon+\1,p_2\coloneqq\yon+\2,$ and $p_3\coloneqq\yon^\2$, we compute $\prod_{i\in\3}p_i$ via \eqref{eqn.poly_prod} as follows:
    \begin{align*}
      \prod_{i\in\3} p_i
      &\iso
      \sum_{\bar{j} \in \prod_{i\in\3} p_i(\1)} \yon^{\sum_{i\in\3} p_i[\bar{j}(i)]}
      \tag*{\eqref{eqn.poly_prod}} \\
      &\iso
      \sum_{\bar{j} \colon (i\in\3) \to p_i(\1)} \yon^{p_1[\bar{j}(1)] + p_2[\bar{j}(2)] + p_3[\bar{j}(3)]} \\
      &\iso
      \yon^{p_1[1] + p_2[1] + p_3[1]}
      + \yon^{p_1[1] + p_2[2] + p_3[1]}
      + \yon^{p_1[1] + p_2[3] + p_3[1]} \\
      &+ \yon^{p_1[2] + p_2[1] + p_3[1]}
      + \yon^{p_1[2] + p_2[2] + p_3[1]}
      + \yon^{p_1[2] + p_2[3] + p_3[1]} \\
      &\iso
      \yon^{\1 + \1 + \2}
      + \yon^{\1 + \0 + \2}
      + \yon^{\1 + \0 + \2} \\
      &+ \yon^{\0 + \1 + \2}
      + \yon^{\0 + \0 + \2}
      + \yon^{\0 + \0 + \2} \\
      % &\iso
      % \yon^\4 + \yon^\3 + \yon^\3 + \yon^\3 + \yon^\2 + \yon^\2 \\
      &\iso
      \yon^\4 + \3\yon^\3 + \2\yon^\2,
    \end{align*}
    as we might expect from standard polynomial multiplication.
  \end{solution}
\end{exercise}

It follows from \eqref{eqn.poly_prod} that the terminal object of $\poly$ is $\1$, and that binary products are given by
\begin{equation}\label{eqn.poly_times}
  p \times q \iso \sum_{i \in p(\1)} \sum_{j \in q(\1)} \yon^{p[i] + q[j]}.
\end{equation}

We will sometimes write $pq$ rather than $p\times q$:
\[
  pq\coloneqq p\times q.
\]

\index{corolla forest!for product of polynomials}
\begin{example}
  We can draw the product of two polynomials in terms of their associated forests. Let $p\coloneqq\yon^\3+\yon$ and $q\coloneqq\yon^\4+\yon^\2+\1$.
  
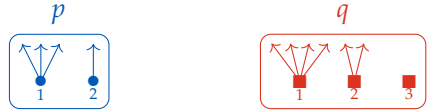
\[
  \begin{tikzpicture}[rounded corners]
    \node (p1) [draw, my-blue, "\color{my-blue} $p$" above] {
      \begin{tikzpicture}[trees, sibling distance=2.5mm]
        \node["\tiny 1" below] (1) {$\bullet$}
        child {}
        child {}
        child {};
        \node[right=.5 of 1,"\tiny 2" below] (2) {$\bullet$}
        child {};
      \end{tikzpicture}
    };
    \node (p2) [draw, my-red, right=2 of p1, "\color{my-red} $q$" above] {
      \begin{tikzpicture}[trees, sibling distance=2.5mm]
        \node["\tiny 1" below] (1) {\small$\blacksquare$}
        child {}
        child {}
        child {}
        child {};
        \node[right=.5 of 1,"\tiny 2" below] (2) {\small$\blacksquare$}
        child {}
        child {};
        \node[right=.5 of 2,"\tiny 3" below] (3) {\small$\blacksquare$}
        ;
      \end{tikzpicture}
    };
  \end{tikzpicture}
  \]
  Then $pq\cong\yon^\7+\2\yon^\5+\2\yon^\3+\yon$.
  We take all pairs of positions, and for each pair we take the disjoint union of the directions.
  
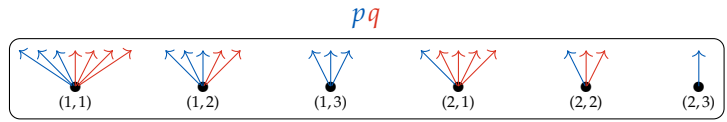
\[
  \begin{tikzpicture}[rounded corners]
    \node (p1) [draw, "${\color{my-blue}p}{\color{my-red}q}$" above] {
      \begin{tikzpicture}[trees, sibling distance=2.5mm]
        \node["\tiny {(1,\,1)}" below] (11) {$\bullet$}
        child[my-blue] {}
        child[my-blue] {}
        child[my-blue] {}
        child[my-red] {}
        child[my-red] {}
        child[my-red] {}
        child[my-red] {};
        \node[right=1.5 of 11, "\tiny {(1,\,2)}" below] (12) {$\bullet$}
        child[my-blue] {}
        child[my-blue] {}
        child[my-blue] {}
        child[my-red] {}
        child[my-red] {};
        \node[right=1.5 of 12, "\tiny {(1,\,3)}" below] (13) {$\bullet$}
        child[my-blue] {}
        child[my-blue] {}
        child[my-blue] {};
        \node[right=1.5 of 13, "\tiny {(2,\,1)}" below] (21) {$\bullet$}
        child[my-blue] {}
        child[my-red] {}
        child[my-red] {}
        child[my-red] {}
        child[my-red] {};
        \node[right=1.5 of 21, "\tiny {(2,\,2)}" below] (22) {$\bullet$}
        child[my-blue] {}
        child[my-red] {}
        child[my-red] {};
        \node[right=1.3 of 22, "\tiny {(2,\,3)}" below] (23) {$\bullet$}
        child[my-blue] {};
      \end{tikzpicture}
    };
  \end{tikzpicture}
  \]
\end{example}

In practice, we can multiply polynomial functors the same way we would multiply two polynomials in high school algebra.

\begin{exercise} \label{exc.general_poly_times}
  \begin{enumerate}
    \item \label{exc.general_poly_times.monomial} Show that for sets $A_1, B_1, A_2, B_2$, we have
    \[
    B_1\yon^{A_1} \times B_2\yon^{A_2} \iso B_1 B_2\yon^{A_1 + A_2}.
    \]
    \item \label{exc.general_poly_times.polynomial} Show that for sets $(A_i)_{i \in I},(A_j)_{j \in J},(B_i)_{i \in I},$ and $(B_j)_{j \in J}$, we have
    \[
    \left(\sum_{i \in I} B_i\yon^{A_i}\right) \times \left(\sum_{j \in J} B_j\yon^{A_j}\right) \iso \sum_{i \in I} \sum_{j \in J} B_i B_j \yon^{A_i + A_j}.
    \]
  \end{enumerate}
  \begin{solution}
    \begin{enumerate}
      \item We compute the product using \eqref{eqn.poly_times}:
      \begin{align*}
        B_1\yon^{A_1} \times B_2\yon^{A_2} &\iso \left(\sum_{i \in B_1} \yon^{A_1}\right) \times \left(\sum_{j \in B_2} \yon^{A_2}\right) \\
        &\iso \sum_{i \in B_1} \sum_{j \in B_2} \yon^{A_1 + A_2} \\
        &\iso B_1 B_2\yon^{A_1 + A_2}.
      \end{align*}
      \item We expand the product by applying \eqref{eqn.set_completely_distributive}, with $I_1 \coloneqq I$ and $I_2 \coloneqq J$:
      \begin{align*}
        \left(\sum_{i \in I} B_i\yon^{A_i}\right) \times \left(\sum_{j \in J} B_j\yon^{A_j}\right) &\iso \prod_{k \in \2} \sum_{i \in I_k} B_i\yon^{A_i} \\
        &\iso \sum_{\bar{i} \in \prod_{k \in \2} I_k} \prod_{k \in \2} B_{\bar{i}(k)}\yon^{A_{\bar{i}(k)}} \\
        &\iso \sum_{(i,\,j) \in IJ} B_i\yon^{A_i} \times B_j\yon^{A_j} \\
        &\iso \sum_{i \in I} \sum_{j \in J} B_i B_j \yon^{A_i + A_j}
      \end{align*}
      where the last isomorphism follows from \cref{exc.general_poly_times.monomial}.
    \end{enumerate}
  \end{solution}
\end{exercise}

As lenses, the canonical projections $\pi \colon pq \to p$ and $\varphi \colon pq \to q$ behave quite naturally: on positions, they are the projections from $(pq)(\1) \iso p(\1) \times q(\1)$ to $p(\1)$ and $q(\1)$, respectively; on directions, they are the inclusions $p[i] \to p[i] + q[j]$ and $q[j] \to p[i] + q[j]$ for each position $(i, j)$ of $pq$.

\begin{exercise} \label{exc.poly_prod}
  Verify that, for $p, q \in \poly$, the polynomial $pq$ given by \eqref{eqn.poly_times} along with the lenses $\pi \colon pq \to p$ and $\varphi \colon pq \to q$ described above satisfy the universal property of the product of $p$ and $q$.
  \begin{solution}
    We wish to show that, for $p, q \in \poly$, the polynomial $pq$ along with the lenses $\pi \colon pq \to p$ and $\varphi \colon pq \to q$ as described in the text satisfy the universal property of the product.
    That is, we must show that for any $r \in \poly$ and lenses $f \colon r \to p$ and $g \colon r \to q$, there exists a unique lens $h \colon r \to pq$ for which the following diagram commutes:
    \begin{equation} \label{eqn.prod_univ_prop}
      \begin{tikzcd}
        r \ar[d, "f"'] \ar[r, "g"] \ar[dr, "h", dashed] & q \\
        p & pq. \ar[l, "\pi"] \ar[u, "\varphi"']
      \end{tikzcd}
    \end{equation}
    We apply \cref{prop.comm_poly}.
    In order for \eqref{eqn.prod_univ_prop} to commute, it must commute on positions---that is, the following diagram of sets must commute:
    \begin{equation} \label{eqn.prod_univ_prop_pos}
      \begin{tikzcd}
        r(\1) \ar[d, "f_\1"'] \ar[r, "g_\1"] \ar[dr, "h_\1", dashed] & q(\1) \\
        p(\1) & (pq)(\1). \ar[l, "\pi_\1"] \ar[u, "\varphi_\1"']
      \end{tikzcd}
    \end{equation}
    But since $(pq)(\1) \iso p(1) \times q(\1)$ along with the projections $\pi_\1$ and $\varphi_\1$ form the product of $p(\1)$ and $q(\1)$ in $\smset$, there exists a unique $h_\1$ for which \eqref{eqn.prod_univ_prop_pos} commutes.
    Hence $h$ is uniquely characterized on positions.
    In particular, it must send each $k \in r(\1)$ to the pair $(f_\1k, g_\1k) \in (pq)(\1)$.

    Moreover, if \eqref{eqn.coprod_univ_prop} is to commute on directions, then for every $k \in r(\1)$, the following diagram of sets must commute:
    \begin{equation} \label{eqn.prod_univ_prop_dir}
      \begin{tikzcd}[sep=large]
        r[k] & q[g_\1k] \ar[l, "g^\sharp_k"'] \ar[d, "\varphi^\sharp_{(f_\1k,\, g_\1k)}"] \\
        p[f_\1k] \ar[u, "f^\sharp_k"] \ar[r, "\pi^\sharp_{(f_\1k,\, g_\1k)}"'] & (pq)[(f_\1k, g_\1k)]. \ar[ul, "h^\sharp_k"', dashed]
      \end{tikzcd}
    \end{equation}
    As $(pq)[(f_\1k, g_\1k)] \iso p[f_\1k] + q[g_\1k]$ along with the inclusions $\pi^\sharp_{(f_\1k,\,g_\1k)}$ and $\varphi^\sharp_{(f_\1k,\,g_\1k)}$ form the coproduct of $p[f_\1k]$ and $q[g_\1k]$ in $\smset$, there exists a unique $h^\sharp_k$ for which \eqref{eqn.prod_univ_prop_dir} commutes.
    Hence $h$ is also uniquely characterized on directions, so it is unique overall.
    Moreover, we have shown that we can define $h$ on positions so that \eqref{eqn.prod_univ_prop_pos} commutes, and that we can define $h$ on directions such that \eqref{eqn.prod_univ_prop_dir} commutes.
    As the commutativity of \eqref{eqn.prod_univ_prop_pos} and \eqref{eqn.prod_univ_prop_dir} together imply the commutativity of \eqref{eqn.prod_univ_prop}, it follows that there exists $h$ for which \eqref{eqn.prod_univ_prop} commutes.
  \end{solution}
\end{exercise}

Much of \cref{part.comon} will focus on the features of another monoidal structure, an asymmetric one, whose definition we will postpone---we will save its surprises for when we can better savor them.
But here we introduce a third symmetric monoidal structure, given by an operation you were not allowed to do to polynomials back in high school.
% TODO: explain coprods, prods (both univ prop and monoidal), and parallel products in terms of interaction

\index{monoidal structure!parallel|see{Parallel product}}
\index{Dirichlet product|see{parallel product}}
\index{parallel product|(}

\begin{definition}[Parallel product of polynomials] \label{def.parallel}
Let $p$ and $q$ be polynomials. Their \emph{parallel product} (also called \emph{Dirichlet product}), denoted $p\otimes q$, is given by the formula
\begin{equation}\label{eqn.parallel_def}
p\otimes q\coloneqq\sum_{i\in p(\1)}\sum_{j\in q(\1)}\yon^{p[i]\times q[j]}.
\end{equation}
\end{definition}

One should compare this with the formula for the product of polynomials shown in \eqref{eqn.poly_times}. The difference is that the parallel product multiplies exponents where the categorical product adds them.

\begin{exercise} \label{exc.general_poly_parallel_times}
  \begin{enumerate}
    \item \label{exc.general_poly_parallel_times.monomial} Show that for sets $A_1, B_1, A_2, B_2$, we have
    \[
    B_1\yon^{A_1} \otimes B_2\yon^{A_2} \iso B_1 B_2\yon^{A_1 A_2}.
    \]
    \item \label{exc.general_poly_parallel_times.polynomial} Show that for sets $(A_i)_{i \in I},(A_j)_{j \in J},(B_i)_{i \in I},$ and $(B_j)_{j \in J}$, we have
    \[
    \left(\sum_{i \in I} B_i\yon^{A_i}\right) \otimes \left(\sum_{j \in J} B_j\yon^{A_j}\right) \iso \sum_{i \in I} \sum_{j \in J} B_i B_j \yon^{A_i A_j}.
    \]
  \end{enumerate}
  \begin{solution}
    \begin{enumerate}
      \item We compute the parallel product using \eqref{eqn.parallel_def}:
      \begin{align*}
        B_1\yon^{A_1} \otimes B_2\yon^{A_2} &\iso \left(\sum_{i \in B_1} \yon^{A_1}\right) \otimes \left(\sum_{j \in B_2} \yon^{A_2}\right) \\
        &\iso \sum_{i \in B_1} \sum_{j \in B_2} \yon^{A_1 \times A_2} \\
        &\iso B_1 B_2\yon^{A_1 A_2}.
      \end{align*}

      \item We expand the parallel product as follows:
      \begin{align*}
        \left(\sum_{i \in I} B_i\yon^{A_i}\right) \otimes \left(\sum_{j \in J} B_j\yon^{A_j}\right) &\iso \left(\sum_{i \in I} \sum_{i' \in B_i} \yon^{A_i}\right) \otimes \left(\sum_{j \in J} \sum_{j' \in B_j} \yon^{A_j}\right) \\
        &\iso \sum_{i \in I} \sum_{i' \in B_i} \sum_{j \in J} \sum_{j' \in B_j} \yon^{A_i \times A_j} \\
        &\iso \sum_{i \in I} \sum_{j \in J} \sum_{i' \in B_i} \sum_{j' \in B_j} \yon^{A_i A_j} \\
        &\iso \sum_{i \in I} \sum_{j \in J} B_i B_j \yon^{A_i A_j}.
      \end{align*}
    \end{enumerate}
  \end{solution}
\end{exercise}

\begin{exercise} \label{exc.some_parallel_prods}
  \begin{enumerate}
    \item \label{exc.some_parallel_prods.const} If $p\coloneqq A$ and $q\coloneqq B$ are constant polynomials, what is $p\otimes q$?
    \item \label{exc.some_parallel_prods.lin} If $p\coloneqq A\yon$ and $q\coloneqq B\yon$ are linear polynomials, what is $p\otimes q$?
    \item \label{exc.some_parallel_prods.pos_prod} For arbitrary $p,q\in\poly$, show that the sets $(p\otimes q)(\1)$ and $p(\1)\times q(\1)$ are isomorphic.
    \qedhere
  \end{enumerate}
  \begin{solution}
    \begin{enumerate}
      \item By \cref{exc.general_poly_parallel_times} \cref{exc.general_poly_parallel_times.monomial}, we have $A\otimes B \iso A\yon^\0\otimes B\yon^\0 \iso AB\yon^{\0} \iso AB$.
      \item By \cref{exc.general_poly_parallel_times} \cref{exc.general_poly_parallel_times.monomial}, we have $A\yon\otimes B\yon \iso A\yon^\1\otimes B\yon^\1 \iso AB\yon^\1 \iso AB\yon$.
      \item By \eqref{eqn.parallel_def},
      \[
      (p \otimes q)(\1) \iso \sum_{i \in p(\1)} \sum_{j \in q(\1)} \1^{p[i] \times q[j]} \iso p(\1) \times q(\1).
      \]
    \end{enumerate}
  \end{solution}
\end{exercise}

\begin{exercise}
  Consider the polynomials $p\coloneqq \2\yon^\2+\3\yon$ and $q\coloneqq\yon^\4+\3\yon^\3$.
  \begin{enumerate}
    \item What is $p\times q$?
    \item What is $p\otimes q$?
    \item Expand the following expression in the variable $y$ according to the ordinary laws of arithmetic.
    \[
    (2\cdot2^y+3\cdot 1^y) \cdot
    (1\cdot4^y+3\cdot 3^y)
    \]
    The factors of the above product are called \emph{Dirichlet series}.\index{Dirichlet series}
    \item Describe the connection between the last two parts. (This is why the parallel product $\otimes$ is also known as the \emph{Dirichlet product}.) \qedhere
  \end{enumerate}
  \begin{solution}
    \begin{enumerate}
      \item We compute $p \times q$ using \cref{exc.general_poly_times} \cref{exc.general_poly_times.polynomial}:
      \begin{align*}
        p \times q &\iso \2\yon^{\2 + \4} + (\2 \times \3)\yon^{\2 + \3} + \3\yon^{\1 + \4} + (\3 \times \3)\yon^{\1 + \3} \\
        &\iso \2\yon^\6 + \6\yon^\5 + \3\yon^\5 + \9\yon^\4 \\
        &\iso \2\yon^\6 + \9\yon^\5 + \9\yon^\4.
      \end{align*}

      \item We compute $p \otimes q$ using \cref{exc.general_poly_parallel_times} \cref{exc.general_poly_parallel_times.polynomial}:
      \begin{align*}
        p \otimes q &\iso \2\yon^{\2 \times \4} + (\2 \times \3)\yon^{\2 \times \3} + \3\yon^\4 + (\3 \times \3)\yon^\3 \\
        &\iso \2\yon^\8 + \6\yon^\6 + \3\yon^\4 + \9\yon^\3.
      \end{align*}

      \item We evaluate $(2\cdot2^y+3\cdot1^y+1) \cdot
      (1\cdot4^y+3\cdot3^y+2)$ using ordinary laws of arithmetic:
      \begin{align*}
        (2\cdot2^y+3\cdot 1^y) \cdot (1\cdot4^y+3\cdot 3^y)
          &=
        (2\cdot2^y+3)(4^y+3\cdot 3^y) \\
        &= 2\cdot2^y\cdot4^y + 2\cdot3\cdot2^y\cdot3^y + 3\cdot4^y + 3\cdot3\cdot3^y \\
        &= 2\cdot8^y + 6\cdot6^y + 3\cdot4^y + 9\cdot3^y.
      \end{align*}

      \item We describe the connection between the last two parts as follows.
      Given a polynomial $p$, we let $d(p)$ denote the Dirichlet series $\sum_{i \in p(\1)} |p[i]|^y$.
      Then by \eqref{eqn.parallel_def},
      \begin{align*}
        d(p \otimes q) &= \sum_{i \in p(\1)} \sum_{j \in q(\1)} |p[i] \times q[j]|^y \\
        &= \sum_{i \in p(\1)} |p[i]|^y \sum_{j \in q(\1)} |q[j]|^y \\
        &= d(p) \cdot d(q).
      \end{align*}
      The last two parts are simply an example of this identity for a specific choice of $p$ and $q$.
    \end{enumerate}
  \end{solution}
\end{exercise}

\begin{example}
We can draw the parallel product of two polynomials in terms of their associated forests. Let $p\coloneqq\yon^\3+\yon$ and $q\coloneqq\yon^\4+\yon^\2+\1$.
\[
\begin{tikzpicture}[rounded corners]
	\node (p1) [draw, my-blue, "\color{my-blue}$p$" above] {
	\begin{tikzpicture}[trees, sibling distance=2.5mm]
    \node["\tiny 1" below] (1) {$\bullet$}
      child {}
      child {}
      child {};
    \node[right=.5 of 1,"\tiny 2" below] (2) {$\bullet$}
      child {};
  \end{tikzpicture}
  };
	\node (p2) [draw, my-red, right=2 of p1, "\color{my-red}$q$" above] {
	\begin{tikzpicture}[trees, sibling distance=2.5mm]
    \node["\tiny 1" below] (1) {\small$\blacksquare$}
      child {}
      child {}
      child {}
      child {};
    \node[right=.5 of 1,"\tiny 2" below] (2) {\small$\blacksquare$}
      child {}
      child {};
    \node[right=.5 of 2,"\tiny 3" below] (3) {\small$\blacksquare$}
    ;
  \end{tikzpicture}
  };
\end{tikzpicture}
\]
Then $p\otimes q\cong\yon^{\1\2}+\yon^\6+\yon^\4+\yon^\2+\2$.
We take all pairs of positions, and for each pair we take the product of the directions.
\[
\begin{tikzpicture}[rounded corners]
	\node (p1) [draw, "${\color{my-blue}p}\otimes{\color{my-red}q}$" above] {
	\begin{tikzpicture}[trees, sibling distance=2mm]
    \node["\tiny {(1,\,1)}" below] (11) {$\bullet$}
      child {}
      child {}
      child {}
      child {}
      child {}
      child {}
      child {}
      child {}
      child {}
      child {}
      child {}
      child {}
    ;
    \node[right=2 of 11, "\tiny {(1,\,2)}" below] (12) {$\bullet$}
      child {}
      child {}
      child {}
      child {}
      child {}
      child {}
    ;
    \node[right=1.2 of 12, "\tiny {(1,\,3)}" below] (13) {$\bullet$}
    ;
   \node[right=1.2 of 13, "\tiny {(2,\,1)}" below] (21) {$\bullet$}
      child {}
      child {}
      child {}
      child {}
 		;
		\node[right=1.2 of 21, "\tiny {(2,\,2)}" below] (22) {$\bullet$}
      child {}
      child {}
 		;
    \node[right=1 of 22, "\tiny {(2,\,3)}" below] (23) {$\bullet$}
 		;
	\end{tikzpicture}
	};
\end{tikzpicture}
\]
\end{example}

\index{monoidal structure!corolla forest depiction}

\begin{exercise}
  Let $p\coloneqq\yon^\2+\yon$ and $q\coloneqq\2\yon^\4$.
  \begin{enumerate}
    \item Draw $p$ and $q$ as corolla forests.
    \item Draw $pq=p\times q$ as a corolla forest.
    \item Draw $p\otimes q$ as a corolla forest.
    \qedhere
  \end{enumerate}
  \begin{solution}
    \begin{enumerate}
      \item Here are $p$ and $q$ drawn as corolla forests:
      \[
      \begin{tikzpicture}[rounded corners]
        \node (p1) [draw, my-blue, "\color{my-blue}$p$" above] {
          \begin{tikzpicture}[trees, sibling distance=2.5mm]
            \node["\tiny 1" below] (1) {$\bullet$}
            child {}
            child {};
            \node[right=.5 of 1,"\tiny 2" below] (2) {$\bullet$}
            child {};
          \end{tikzpicture}
        };
        \node (p2) [draw, my-red, right=2 of p1, "\color{my-red}$q$" above] {
          \begin{tikzpicture}[trees, sibling distance=2.5mm]
            \node["\tiny 1" below] (1) {\small$\blacksquare$}
            child {}
            child {}
            child {}
            child {};
            \node[right=1 of 1,"\tiny 2" below] (2) {\small$\blacksquare$}
            child {}
            child {}
            child {}
            child {};
          \end{tikzpicture}
        };
      \end{tikzpicture}
      \]

      \item Here is $pq$ drawn as a corolla forest:
      \[
      \begin{tikzpicture}[rounded corners]
        \node (p1) [draw, "$pq$" above] {
          \begin{tikzpicture}[trees, sibling distance=2.5mm]
            \node["\tiny {(1,\,1)}" below] (11) {$\bullet$}
            child[my-blue] {}
            child[my-blue] {}
            child[my-red] {}
            child[my-red] {}
            child[my-red] {}
            child[my-red] {};
            \node[right=1.5 of 11, "\tiny {(1,\,2)}" below] (12) {$\bullet$}
            child[my-blue] {}
            child[my-blue] {}
            child[my-red] {}
            child[my-red] {}
            child[my-red] {}
            child[my-red] {};
            \node[right=1.5 of 12, "\tiny {(2,\,1)}" below] (21) {$\bullet$}
            child[my-blue] {}
            child[my-red] {}
            child[my-red] {}
            child[my-red] {}
            child[my-red] {};
            \node[right=1.5 of 21, "\tiny {(2,\,2)}" below] (22) {$\bullet$}
            child[my-blue] {}
            child[my-red] {}
            child[my-red] {}
            child[my-red] {}
            child[my-red] {};
          \end{tikzpicture}
        };
      \end{tikzpicture}
      \]
      \item Here is $p \otimes q$ drawn as a corolla forest:
      \[
      \begin{tikzpicture}[rounded corners]
        \node (p1) [draw, "$p\otimes q$" above] {
          \begin{tikzpicture}[trees, sibling distance=2mm]
            \node["\tiny {(1,\,1)}" below] (11) {$\bullet$}
            child {}
            child {}
            child {}
            child {}
            child {}
            child {}
            child {}
            child {}
            ;
            \node[right=2 of 11, "\tiny {(1,\,2)}" below] (12) {$\bullet$}
            child {}
            child {}
            child {}
            child {}
            child {}
            child {}
            child {}
            child {}
            ;
            \node[right=1.5 of 12, "\tiny {(2,\,1)}" below] (21) {$\bullet$}
            child {}
            child {}
            child {}
            child {}
            ;
            \node[right=1.5 of 21, "\tiny {(2,\,2)}" below] (22) {$\bullet$}
            child {}
            child {}
            child {}
            child {}
            ;
          \end{tikzpicture}
        };
      \end{tikzpicture}
      \]
    \end{enumerate}
  \end{solution}
\end{exercise}

\begin{exercise}\label{exc.prepare_poly_smc}
Let $p,q,r\in\poly$ be any polynomials.
\begin{enumerate}
  \item Show that there is an isomorphism $p\otimes\yon\iso p$.
  \item Show that there is an isomorphism $(p\otimes q)\otimes r\iso p\otimes (q\otimes r)$.
  \item Show that there is an isomorphism $p\otimes q \iso q\otimes p$.
 \qedhere
\end{enumerate}
\begin{solution}
\begin{enumerate}
  \item We show that $p\otimes\yon\iso p$:
  \begin{align*}
      p \otimes y &\iso \sum_{i \in p(\1)} \sum_{j \in \1} \yon^{p[i] \times \1} \tag*{\eqref{eqn.parallel_def}} \\
      &\iso \sum_{i \in p(\1)} \yon^{p[i]} \iso p.
  \end{align*}

  \item We show that $(p\otimes q)\otimes r\cong p\otimes (q\otimes r)$:
  \begin{align*}
      (p \otimes q) \otimes r &\iso \left(\sum_{i \in p(\1)} \sum_{j \in q(\1)} \yon^{p[i] \times q[j]}\right) \otimes r \tag*{\eqref{eqn.parallel_def}} \\
      &\iso \sum_{i \in p(\1)} \sum_{j \in q(\1)} \left(\sum_{k \in r(\1)} \yon^{(p[i] \times q[j]) \times r[k]}\right) \tag*{\eqref{eqn.parallel_def}} \\
      &\iso \sum_{i \in p(\1)} \left(\sum_{j \in q(\1)} \sum_{k \in r(\1)} \yon^{p[i] \times (q[j] \times r[k])}\right) \tag{Associativity of $\sum$ and $\times$} \\
      &\iso p \otimes \left(\sum_{j \in q(\1)} \sum_{k \in r(\1)} \yon^{q[j] \times r[k]}\right) \tag*{\eqref{eqn.parallel_def}} \\
      &\iso p \otimes (q \otimes r). \tag*{\eqref{eqn.parallel_def}} \\
  \end{align*}

  \item We show that $(p\otimes q)\cong(q\otimes p)$:
  \begin{align*}
      p \otimes q &\iso \sum_{i \in p(\1)} \sum_{j \in q(\1)} \yon^{p[i] \times q[j]} \tag*{\eqref{eqn.parallel_def}} \\
      &\iso \sum_{j \in q(\1)} \sum_{i \in p(\1)} \yon^{q[j] \times p[i]} \tag{Commutativity of $\sum$ and $\times$} \\
      &\iso q \otimes p. \tag*{\eqref{eqn.parallel_def}} \\
  \end{align*}
\end{enumerate}
\end{solution}
\end{exercise}

In \cref{exc.prepare_poly_smc}, we have gone most of the way to proving that $(\poly,\yon,\otimes)$ is a symmetric monoidal category.
We sketch the rest of the proof as follows.

\begin{proposition}\label{prop.parallel_monoidal}
The category $\poly$ has a symmetric monoidal structure $(\yon,\otimes)$ where $\otimes$ is the parallel product from \cref{def.parallel}.
\end{proposition}
\begin{proof}[Sketch of proof]
Given two lenses $f\colon p\to p'$ and $g\colon q\to q'$, we need to define a lens $(f\otimes g)\colon (p\otimes q)\to (p'\otimes q')$. This is easiest to define using polyboxes, keeping in mind that the positions and directions of a parallel product are pairs of positions and directions of its constituent factors:
\[
\begin{tikzpicture}[polybox, mapstos]
  \node[poly, dom, "$p\otimes q$" left] (p) {$(f^\sharp_i a,g^\sharp_j b)$\at$(i,j)\vphantom{(f_\1i,g_\1j)}$};
  \node[poly, cod, "$p'\otimes q'$" right, right=of p] (q) {$(a,b)\vphantom{(f^\sharp_i a,g^\sharp_j b)}$\at$(f_\1i,g_\1j)$};
  \draw (p_pos) -- node[below] {$(f\otimes g)_\1$} (q_pos);
  \draw (q_dir) -- node[above] {$(f\otimes g)^\sharp$} (p_dir);
\end{tikzpicture}
\]
Here $i\in p(\1), j\in q(\1), a\in p'[f_\1 i],$ and $b\in q'[g_\1 j]$.

Then \cref{exc.prepare_poly_smc} gives us the unitors, associator, and braiding.
We have not proven the functoriality of $\otimes$, the naturality of the isomorphisms from \cref{exc.prepare_poly_smc}, or all the coherences between these isomorphisms, but we ask the reader to take them on trust or to check them for themselves.
Alternatively, we may invoke the Day convolution to obtain the monoidal structure $(\yon, \otimes)$ directly: see \cref{prop.day}.
\end{proof}\index{Day convolution}

\begin{exercise}
  \begin{enumerate}
    \item What is $(\3\yon^\5+\6\yon^\2)\otimes\4$? Hint: $\4=\4\yon^\0$.
    \item Is the class of constant polynomials a \emph{$\otimes$-ideal}; that is, is the parallel product of a polynomial and a constant polynomial always a constant? \qedhere
  \end{enumerate}
  \begin{solution}
    \begin{enumerate}
      \item We compute $(\3\yon^\5+\6\yon^\2)\otimes\4$ using \cref{exc.general_poly_parallel_times} \cref{exc.general_poly_parallel_times.polynomial} and the fact that $\4=\4\yon^\0$:
      \begin{align*}
        (\3\yon^\5+\6\yon^\2)\otimes\4\yon^\0 &\iso (\3\times\4)\yon^{\5\times\0} + (\6\times\4)\yon^{\2\times\0} \\
        &\iso \1\2\yon^\0 + \2\4\yon^\0 \\
        &\iso \3\6.
      \end{align*}
      \item Given a polynomial $p$ and a set $J$ viewed as a constant polynomial, we have
      \begin{align*}
        p\otimes J &\iso \left(\sum_{i\in p(\1)}\yon^{p[i]}\right)\otimes\left(\sum_{j\in J}\yon^\0\right) \\
        &\iso \sum_{i\in p(\1)}\sum_{j\in J}\yon^{p[i]\times\0} \\
        &\iso p(\1)J\yon^\0 \\
        &\iso p(\1)J,
      \end{align*}
      itself a constant polynomial; so the class of constant polynomials is indeed a $\otimes$-ideal.
    \end{enumerate}
  \end{solution}
\end{exercise}

\index{parallel product!of special polynomials}

\begin{exercise}\label{exc.dir_closed_classes}\index{affine polynomial}\index{linear polynomial}\index{monomial}
Which of the following special classes of polynomials are closed under $\otimes$? Note also whether they contain $\yon$.
\begin{enumerate}
	\item The class $\{A\yon^\0\mid A\in\smset\}$ of constant polynomials.
	\item The class $\{A\yon\mid A\in\smset\}$ of linear polynomials.
	\item The class $\{A\yon+B\mid A,B\in\smset\}$ of affine polynomials.
	\item The class $\{A\yon^\2+B\yon+C\mid A,B,C\in\smset\}$ of quadratic polynomials.
	\item The class $\{A\yon^B\mid A,B\in\smset\}$ of monomials.
	\item The class $\{S\yon^S\mid S\in\smset\}$.
	\item The class $\{p\in\poly\mid p(\1)\text{ is finite}\}$. \qedhere
\end{enumerate}
\begin{solution}
For each of the following classes of polynomials, we determine whether they are closed under $\otimes$ and whether they contain $\yon$.
\begin{enumerate}
	\item The set $\{A\yon^\0\mid A\in\smset\}$ of constant polynomials is closed under $\otimes$ by the solution to \cref{exc.some_parallel_prods} \cref{exc.some_parallel_prods.const}.
	But the set does not contain $\yon$, as $\yon$ is not a constant polynomial.
	\item The set $\{A\yon\mid A\in\smset\}$ of linear polynomials is closed under $\otimes$ by the solution to \cref{exc.some_parallel_prods} \cref{exc.some_parallel_prods.lin} and does contain $\yon$, as $\yon \iso \1\yon$.
	\item The set $\{A\yon+B\mid A,B\in\smset\}$ of affine polynomials is closed under $\otimes$, for \cref{exc.general_poly_parallel_times} \cref{exc.general_poly_parallel_times.polynomial} yields
	\[
	    (A\yon + B) \otimes (A'\yon + B') \iso AA'\yon + AB' + BA' + BB'.
	\]
	The set contains $\yon$, as $\yon \iso \1\yon + \0 $.
	\item The set $\{A\yon^\2+B\yon+C\mid A,B,C\in\smset\}$ of quadratic polynomials is not closed under $\otimes$, for even though $\yon^\2 \iso \1\yon^\2 + \0\yon + \0$ is a quadratic polynomial, \cref{exc.general_poly_parallel_times} \cref{exc.general_poly_parallel_times.monomial} implies that
	\[
	    \yon^\2 \otimes \yon^\2 \iso \yon^\4,
	\]
	which is not quadratic.
	The set contains $\yon$, as $\yon \iso \0\yon^\2 + \1\yon + \0$.
	\item The set $\{A\yon^B\mid A,B\in\smset\}$ of monomials is closed under $\otimes$ by \cref{exc.general_poly_parallel_times} \cref{exc.general_poly_parallel_times.monomial} and does contain $\yon$, as $\yon \iso \1\yon^\1$.
	\item The set $\{S\yon^S\mid S\in\smset\}$ is closed under $\otimes$, for \cref{exc.general_poly_parallel_times} \cref{exc.general_poly_parallel_times.monomial} returns
	\[
	    S\yon^S \otimes T\yon^T \iso ST\yon^{ST}.
	\]
	The set contains $\yon$, as $\yon \iso \1\yon^\1$.
	\item The set $\{p\in\poly\mid p(\1)\text{ is finite}\}$ is closed under $\otimes$ by the solution to \cref{exc.some_parallel_prods} \cref{exc.some_parallel_prods.pos_prod} and the fact that the product of two finite sets is itself finite.
	The set contains $\yon$, as $\yon(\1) \iso \1$ is finite.
\end{enumerate}
\end{solution}
\end{exercise}

\begin{exercise}
What is the smallest class of polynomials that is closed under $\otimes$ and contains $\yon$?
\begin{solution}
The smallest class of polynomials that is closed under $\otimes$ and contains $\yon$ is just $\{\yon\}$.
This is because by \cref{exc.general_poly_parallel_times} \cref{exc.general_poly_parallel_times.monomial}, we have $\yon \otimes \yon \iso \yon$.
\end{solution}
\end{exercise}

\index{distributive law!for parallel product}

\begin{exercise}
Show that for any $p_1,p_2,q\in\poly$ there is an isomorphism
\[
(p_1+p_2)\otimes q\iso (p_1\otimes q)+(p_2\otimes q).
\]
\begin{solution}
We show that $(p_1 + p_2) \otimes q \iso (p_1 \otimes q) + (p_2 \otimes q)$ using \eqref{eqn.parallel_def}:
\begin{align*}
    (p_1 + p_2) \otimes q &\iso \sum_{k \in \2} \sum_{i \in p_k(\1)} \sum_{j \in q(\1)} \yon^{p_k[i] \times q[j]} \\
    &\iso \sum_{i \in p_1(\1)} \sum_{j \in q(\1)} \yon^{p_1[i] \times q[j]} + \sum_{i \in p_2(\1)} \sum_{j \in q(\1)} \yon^{p_2[i] \times q[j]} \\
    &\iso (p_1 \otimes q) + (p_2 \otimes q).
\end{align*}
\end{solution}
\end{exercise}

\begin{remark}
  Monoids in $\poly$ with respect to the parallel product $\otimes$ are particularly interesting---they have a kind of collective semantics, letting agents aggregate their contributions and distribute returns on those contributions in a coherent way.
  We leave discussion of them to future work, so as not to distract us from our main story.

  % TODO: cite collectives paper
\end{remark}

\index{monoid!for parallel product}

There is a more general way to obtain monoidal structures on $\poly$ like $\times$ and $\otimes$ using a construction known as the \emph{Day convolution}, defined by a special kind of colimit known as a \emph{coend}. If you have not seen the Day convolution or coends before, do not fret: we will not use them elsewhere in the book, and rest assured that the fact about coends known as the \emph{co-Yoneda lemma} employed in the following proof is a standard and purely formal result.

\index{Day convolution}\index{distributive law!for Day monoidal structures}
\index{monoidal structure!from monoidal structure on $\smset$}

\begin{proposition} \label{prop.day}
For any monoidal structure $(I,\star)$ on $\smset$, there is a corresponding monoidal structure $(\yon^I, \odot)$ on $\poly$, where $\odot$ is the Day convolution.
Moreover, $\odot$ distributes over coproducts.\index{coproduct!distributing over}

In the case of $(\0,+)$ and $(\1,\times)$, this procedure returns the $(\1,\times)$ and $(\yon,\otimes)$ monoidal structures respectively.
\end{proposition}
\begin{proof}
Any monoidal structure $(I,\star)$ on $\smset$ induces a monoidal structure on $\smset^\smset$ with the Day convolution $\odot$ as the tensor product and $\yon^I$ as the unit.
To prove that this monoidal structure restricts to $\poly$, it suffices to show that $\poly$ is closed under the Day convolution.

Given polynomials $p$ and $q$, their Day convolution in $\smset^\smset$ is given by the coend
\begin{equation} \label{eqn.day_conv.coend}
    p \odot q \iso \int^{(A,B)\in\smset^\2} \yon^{A \star B} \times p(A) \times q(B).
\end{equation}
We can rewrite the product $p(A) \times q(B)$ as
\[
    p(A) \times q(B) \iso \left(\sum_{i \in p(\1)} A^{p[i]}\right) \times \left(\sum_{j \in q(\1)} B^{q[i]}\right) \iso \sum_{(i,\,j) \in p(\1) \times q(\1)} A^{p[i]} \times B^{q[i]}
\]
So because products distribute over coproducts in $\smset^\smset$ and coends always commute with coproducts (as they are both colimits), we can rewrite \eqref{eqn.day_conv.coend} as
\begin{align*}
    p \odot q &\iso \sum_{(i,\,j) \in p(\1) \times q(\1)} \int^{(A,B)\in\smset^\2} \yon^{A \star B} \times A^{p[i]} \times B^{q[i]} \\
    &\iso \sum_{(i,\,j) \in p(\1) \times q(\1)} \int^{(A,B)\in\smset^\2} \yon^{A \star B} \times \smset^\2((p[i],q[j]),(A,B))
\end{align*}
which, by the co-Yoneda lemma, can be rewritten as
\begin{equation} \label{eqn.day_conv.poly}
    p \odot q \iso \sum_{(i,\,j) \in p(\1) \times q(\1)} \yon^{p[i] \star q[j]},
\end{equation}
which is in $\poly$.
That the Day convolution distributes over coproducts also follows from the fact that products distribute over coproducts in $\smset^\smset$ and that coends commute with coproducts; or, alternatively, directly from \eqref{eqn.day_conv.poly}.

We observe that \eqref{eqn.day_conv.poly} gives $(\yon^I, \odot) = (\1, \times)$ when $(I, \star) \coloneqq (\0, +)$ and $(\yon^I, \odot) = (\yon, \otimes)$ when $(I, \star) \coloneqq (\1, \times)$.
\end{proof}\index{Day convolution}

\begin{exercise}
There is a monoidal structure on $\smset$ whose unit is $\0$ and whose product is given by $(A, B)\mapsto A+AB+B$.
\begin{enumerate}
	\item Verify that the operation $(A, B)\mapsto A+AB+B$ on $\smset$ is associative.
	\item Verify that $\0$ is the unit for the above operation.
	\item Let $(\1,\odot)$ denote the corresponding monoidal structure on $\poly$ obtained via \cref{prop.day}. Compute the monoidal product $(\yon^\3+\yon)\odot(\2\yon^\2+\2)$.
\qedhere
\end{enumerate}
\begin{solution}
\begin{enumerate}
    \item To show that the operation $(A,B)\mapsto A+AB+B$ on $\smset$ is associative, observe that
    \begin{align*}
        (A + AB + B) + (A + AB + B)C + C &\iso A + AB + B + AC + ABC + BC + C \\
        &\iso A + AB + ABC + AC + B + BC + C \\
        &\iso A + A(B + BC + C) + (B + BC + C).
    \end{align*}
    \item To show that $\0$ is the unit for this operation, observe that
    \[
        (A, \0) \mapsto A + A\0 + \0 \iso A
    \]
    and
    \[
        (\0, B) \mapsto \0 + \0B + B \iso B.
    \]
    \item Taking $A \star B \coloneqq A + AB + B$ in \cref{prop.day} to obtain a monoidal product $\odot$ on $\poly$, we can use \eqref{eqn.day_conv.poly} to compute that
    \begin{align*}
        (\yon^\3+\yon)\odot(\2\yon^\2+\2) &\iso (\yon^\3+\yon^\1)\odot(\yon^\2+\yon^\2+\yon^\0+\yon^\0) \\
        &\iso \yon^{\3\star\2} + \yon^{\3\star\2} + \yon^{\3\star\0} + \yon^{\3\star\0} + \yon^{\1\star\2} + \yon^{\1\star\2} + \yon^{\1\star\0} + \yon^{\1\star\0} \\
        &\iso \2\yon^{\1\1} + \2\yon^\3 + \2\yon^\5 + \2\yon^\1.
    \end{align*}
\end{enumerate}
\end{solution}
\end{exercise}

\index{parallel product|)}

%-------- Section --------%
\section[Summary and further reading]{Summary and further reading%
  \sectionmark{Summary \& further reading}}
\sectionmark{Summary \& further reading}

In this chapter, we introduced the category $\poly$, whose objects are polynomial functors and whose morphisms are the natural transformations between them.
We call these natural transformations \emph{dependent lenses}, or \emph{lenses} for short.
We also proved our first categorical property of $\poly$: that it has all small coproducts.
%In the next chapter, we will review the coproduct construction before turning to other operations on our polynomials.

The main result of this chapter was a concrete characterization of our dependent lenses between polynomial functors.
A dependent lens $f\colon p\to q$ is characterized by its
\begin{itemize}
  \item \emph{on-positions function}, a function that we denote by $f_\1\colon p(\1)\to q(\1)$ sending $p$-positions forward to $q$-positions; and its
  \item \emph{on-directions functions}, one for each $p$-position $i$ that we denote by $f^\sharp_i\colon q[f_\1 i]\to p[i]$ sending $q[f_\1 i]$-directions backward to $p[i]$-directions.
\end{itemize}
This forward-backward relationship is what makes dependent lenses so well-suited for modeling \emph{interaction protocols}.
Given two agents with positions and directions, a dependent lens between them defines an interaction protocol that describes how the position of the first agent determines the position of the second agent and how the direction of the second agent determines the direction of the first.
This perspective is exhibited by our corolla and polybox pictures for lenses.
We studied examples of lenses between special polynomials: in particular, lenses between monomials are known as \emph{bimorphic lenses} in functional programming literature.

\index{lens!between monomials}

We then unwound our interpretation of natural transformations between polynomials as dependent lenses with on-positions and on-directions functions to describe what happens to these functions when lenses compose.
This gave us an accessible way to interpret commutative diagrams in $\poly$ that is particularly convenient to express using polyboxes.

Finally, we considered various categorical structures on $\poly$, e.g.\ that it has all products and coproducts, and that these distribute.
%$\prod\sum\to\sum\prod$.
%\begin{align*}
%	\sum_{a\in A}p_a&\coloneqq\sum_{(a,i)\in\sum_{a\in A}p_a(1)}\yon^{p_a[i]}&
%	\prod_{a\in A}p_a&\coloneqq\sum_{i\in\prod_{a\in A}p_a(1)}\yon^{\sum_{a\in A}p_a[i a]}
%\\
%	p_1+p_2&\coloneqq\sum_{(a,i)\in\{(1,i_1)\mid i_1\in p_1(\1)\}+\{(2,i_2)\mid i_2\in p_2(\1)\}}\yon^{p_a[i]}&
%	p_1\times p_2&\coloneqq\sum_{(i_1,i_2)\in p_1(1)\times p_2(1)}\yon^{p_1[i_1]+p_2[i_2]}
%\end{align*}
We also discussed how one can take any monoidal product $\star$ from $\smset$ and lift it to a monoidal product $\cdot$ on $\poly$:
\[
	p_1\odot p_2\coloneqq\sum_{(i_1,i_2)\in p_1(1)\times p_2(1)}\yon^{p_1[i_1]\cdot p_2[i_2]}
\]
A special case of this is the product structure $\times$ on $\poly$, which emerges from the coproduct structure $+$ on $\smset$. The other case of interest is the parallel (or Dirichlet) product structure $\otimes$ on $\poly$, which emerges from the product structure $\times$ on $\smset$:
\[
	p_1\otimes p_2\coloneqq\sum_{(i_1,i_2)\in p_1(1)\times p_2(1)}\yon^{p_1[i_1]\times p_2[i_2]}
\]

Variants of lenses are studied in compositional game theory \cite{hedges2016compositionality,hedges2017coherence,hedges2018limits,hedges2018morphisms}, in categorical database theory \cite{johnson2012lenses}, in functional programming and programming language theory \cite{bohannon2006relational,oconnor2011functor,abou2016reflections}, and in more generalized categorical settings \cite{gibbons2012relating,spivak2019generalized}.

\index{compositional game theory}
\index{database}
\index{functional programming!lenses in}

%-------- Section --------%
\section{Exercise solutions}
\Closesolutionfile{solutions}
{\footnotesize
\input{solution-file3}}

\Opensolutionfile{solutions}[solution-file4]
\index{category!of polynomial functors|)}

%------------ Chapter ------------%
\chapter{Dynamical systems as dependent lenses} \label{ch.poly.dyn_sys}

One of the main goals of this book is to use dependent lenses in $\poly$ to model dynamical systems and automata.
In this chapter, we will begin to see how to do this through an array of examples.

%-------- Section --------%
\section{Moore machines}\label{sec.poly.dyn_sys.moore}

\index{Moore machine|(}

We start with our simplest example of a dynamical system: a deterministic state machine with a fixed range of states, inputs, and outputs.
At any point in time, this machine will inhabit one of its possible states and return output according to that current state. It can also update its current state according to the input it receives.

\begin{definition}[Moore machine]\label{def.moore_machine}
  A \emph{Moore machine} consists of the following data: three sets,
  \begin{itemize}
    \item a set $S$, called the \emph{state-set}, whose elements are \emph{states};
    \item a set $I$, called the \emph{position-set} (or \emph{output-set}), whose elements are \emph{positions} (or \emph{outputs});
    \item a set $A$, called the \emph{direction-set} (or \emph{input-set}), whose elements are \emph{directions} (or \emph{inputs});
  \end{itemize}
  and two functions,
  \begin{itemize}
    \item $\text{return}\colon S\to I$;
    \item $\text{update}\colon S\times A\to S$.
  \end{itemize}
  To emphasize the role that the three sets play, we can specify that this is an $(A,I)$\emph{-Moore machine with states} $S$.
\end{definition}

\index{Moore machine!states of}
\index{Moore machine!update}\index{Moore machine!return}

The input and output terminology is standard, while the position and direction terminology is our own: we will soon see how the positions and directions of a Moore machine relate to that of a polynomial.

We should interpret an $(A,I)$-Moore machine as follows.
At any time, the machine inhabits one of the states in its state-set $S$.
Say its current state is $s\in S$.
We can ask the machine to perform one of the following two tasks.
\begin{itemize}
  \item We can ask the machine to \emph{return its position}: it should then produce the position $\text{return}(s) \in I$.
  \item We can feed the machine one of its \emph{directions} $a\in A$ and ask it to \emph{update its state}: it should then replace its current state with the new state $\text{update}(s,a)\in S$.
  Note that the new state depends not only on the direction the machine receives but also on the state the machine inhabits when it receives that direction.
\end{itemize}

We may visualize a Moore machine with a \emph{transition diagram} as follows.

\index{Moore machine!transition diagram}
\index{interface!monomial|(}

\begin{example}[A Moore machine's transition diagram]\label{ex.Moore_three}
Given $A\coloneqq\{\Red,\Blue\}$ and $I\coloneqq\{0,1\}$, we can draw a transition diagram for an $(A,I)$-Moore machine with $S\coloneqq\3$ states as follows:
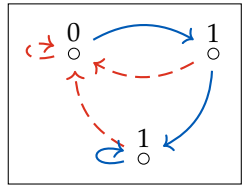
\begin{equation} \label{eqn.trans_diag}
\begin{tikzpicture}
	\node[draw] {
  \begin{tikzcd}[column sep=small]
  	\LMO{0}\ar[rr, my-blue, thick, bend left]\ar[loop left, thick, my-red, dashed]&&
  	\LMO{1}\ar[ll, thick, my-red, dashed, bend left]\ar[dl, bend left, thick, my-blue]\\&
  	\LMO{1} \ar[ul, thick, my-red, dashed, bend left] \ar[loop left, thick, my-blue]
  \end{tikzcd}
  };
\end{tikzpicture}
\end{equation}
Each state is labeled by the position it returns according to the machine's return function.
Additionally, each state has two outgoing arrows, one $\Red$ (dashed) and one $\Blue$ (solid), corresponding to the two possible directions.
The targets of the arrows indicate the updated state according to the machine's update function.

Say the machine starts at the bottom state.
By feeding it a sequence of directions---say $(\Red,\Red,\Blue,\Red,\ldots)$---we can send the machine through its states via its update function and return the position at each state:
\begin{enumerate}
    \item Starting at the bottom state, the machine returns the position $1$.
    \item Following the $\Red$ arrow from the bottom state, the machine updates its state to the left state.
    \item At the left state, the machine returns the position $0$.
    \item Following the $\Red$ arrow from the left state, the machine updates its state to---once again---the left state.
    \item At the left state, the machine returns the position $0$.
    \item Following the $\Blue$ arrow from the left state, the machine updates its state to the right state.
    \item At the right state, the machine returns the position $1$.
    \item Following the $\Red$ arrow from the right state, the machine updates its state to the left state.
    \item At the left state, the machine returns the position $0$.

    \ldots
\end{enumerate}
In summary, starting from the bottom state, this Moore machine sends the sequence $(\Red,\Red,\Blue,\Red,\ldots)$ of directions in $A$ to the sequence $(1,0,0,1,0,\ldots)$ of positions in $I$.
\end{example}

In general, given an initial state $s_0\in S$, an $(A,I)$-Moore machine with states $S$ sends every sequence $(a_1,a_2,a_3,\ldots)$ of directions in $A$ to a sequence $(i_0,i_1,i_2,i_3,\ldots)$ of positions in $I$, defined inductively as follows, via an intermediary sequence $(s_0,s_1,s_2,s_3,\ldots)$ of states in $S$:
\[
    i_k\coloneqq \text{return}(s_k) \qqand s_{k+1}\coloneqq \text{update}(s_k, a_{k+1})
\]
for all $k\in\nn$.
We will see that $\poly$ gives us a more concise way to express this in \cref{ex.input_output}.

Comparing \cref{def.moore_machine} with \cref{subsec.poly.cat.morph.bimorphic-lens}, we find that an $(A,I)$-Moore machine with states $S$ is precisely a lens between monomials $\varphi\colon S\yon^S\to I\yon^A$ with on-positions function $\varphi_\1\coloneqq\text{return}\colon S\to I$ and on-directions map $\varphi^\sharp\coloneqq S\times A\to S$.
The positions and directions of the Moore machine are the positions and directions of the codomain of the corresponding lens, while the domain of the lens has the states of the Moore machine as both its positions and its directions.
So we can repackage \cref{def.moore_machine} as follows.

\begin{definition}[Moore machine, version 2]\label{def.moore_machine2}
  For $S,I,A\in\smset$, an $(A,I)$-\emph{Moore machine with states} $S$ is a lens \[\varphi\colon S\yon^S\to I\yon^A\] in $\poly$.
  We call
  \begin{itemize}
    \item the domain monomial $S\yon^S$ the machine's \emph{state system}: its position-set (equivalently, its direction-set) is the machine's \emph{state-set}, and its positions (equivalently, its directions) are the machine's \emph{states};
    \item the codomain monomial $I\yon^A$ the machine's \emph{interface}: its position-set and direction-set are the machine's \emph{position-set} and \emph{direction-set}, and its positions and directions are the machine's \emph{positions} and \emph{directions};
    \item the on-positions function $\varphi_\1\colon S\to I$ the machine's \emph{return function};
    \item the on-directions map $\varphi^\sharp\colon S\times A\to S$ the machine's \emph{update function}.
  \end{itemize}
\end{definition}
\index{state system}\index{interface}\index{Moore machine}

We call the codomain of a Moore machine its \emph{interface} because it encodes how an outsider interacts with the machine: an outsider observes the positions of the interface that the machine returns and feeds the directions of the interface to the machine to update it.
Rather than directly observing and altering the machine's states, an outsider must interact with the machine via its interface.

\index{Moore machine!interface of|seealso{interface}}

\begin{exercise}
In this exercise, we will write the Moore machine from \cref{ex.Moore_three} as a lens $\varphi$ between monomials.
\begin{enumerate}
    \item What is the machine's state system, the domain of $\varphi$?
    \item What is the machine's interface, the codomain of $\varphi$?
\end{enumerate}
Call the left state $\const{L}$, the right state $\const{R}$, and the bottom state $\const{B}$.
\begin{enumerate}[resume]
    \item What is the machine's return function, the on-positions function of $\varphi$?
    \item What is the machine's update function, the on-directions map of $\varphi$?
    \item Draw the first two steps listed in \cref{ex.Moore_three} of the machine's operation (starting at the bottom state and receiving the direction $\Red$) using polyboxes.
    \qedhere
\end{enumerate}
\begin{solution}
\begin{enumerate}
    \item As $S=\3$, the state system is $S\yon^S=\3\yon^\3$.
    \item As $I=\{0,1\}$ and $A=\{\Red,\Blue\}$, the interface is $I\yon^A=\{0,1\}\yon^{\{\Red,\,\Blue\}}$.
    \item The return function $S\to I$ sends $\const{L}\mapsto 0, \const{R}\mapsto 1,$ and $\const{B}\mapsto 1$.
    \item The update function $S\times A\to S$ sends
    \begin{align*}
        (\const{L}, \Red)\mapsto \const{L}&,\quad(\const{L}, \Blue)\mapsto \const{R}, \\
        (\const{R}, \Red)\mapsto \const{L}&,\quad(\const{R}, \Blue)\mapsto \const{B}, \\
        (\const{B}, \Red)\mapsto \const{L}&,\quad(\const{B}, \Blue)\mapsto \const{B}.
    \end{align*}
    \item In the first two steps, the machine sends its bottom state $\const{B}$ to its position $1$ via its return function, then sends its bottom state $\const{B}$ and its direction $\Red$ to its left state $\const{L}$ via its update function.
    We can interpret this in terms of the lens $\varphi$ and depict the steps in polyboxes as
    \[
    \begin{tikzpicture}[polybox, mapstos]
      \node[poly, dom] (p) {$\const{L}$\at$\const{B}\vphantom{1}$};
        \node[left=0pt of p_pos] {$S$};
        \node[left=0pt of p_dir] {$S$};

      \node[poly, cod, right=of p] (q) {$\vphantom{\const{L}}\Red$\at$1\vphantom{\const{B}}$};
        \node[right=0pt of q_pos] {$I$};
        \node[right=0pt of q_dir] {$A$};

      \draw (p_pos) -- node[below] {$f$} (q_pos);
      \draw (q_dir) -- node[above] {$f\inv$} (p_dir);
    \end{tikzpicture}
    \]
\end{enumerate}
\end{solution}
\end{exercise}

Here are some more examples of Moore machines.

\begin{example}[Counter]\label{ex.counting_trajectory}
There is a $(\1,\nn)$-Moore machine with states $\nn$ that, with initial state $0\in\nn$, returns the sequence of natural numbers $(0,1,2,3,\ldots)$.
The machine is given by the lens $\nn\yon^\nn\to\nn\yon$ whose on-positions function is the identity $\nn\to\nn$ and whose on-directions map $\nn\times\1\iso\nn\to\nn$ sends $n\mapsto n+1$.
Here it is in polyboxes (recall that the shaded direction box indicates that the direction-set is a singleton, i.e.\ there is no choice to be made in filling it in):
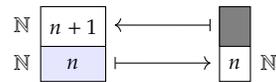
\[
\begin{tikzpicture}[polybox, mapstos]
  \node[poly, dom] (s) {$n+1$\at$n$};
    \node[left=0pt of s_pos] {$\nn$};
    \node[left=0pt of s_dir] {$\nn$};

 	\node[poly, cod, linear, right=of s] (p) {$\vphantom{1}$\at$n$};
    \node[right=0pt of p_pos] {$\nn$};

 	\draw (s_pos) to[first] (p_pos);
 	\draw (p_dir) to[last] (s_dir);
\end{tikzpicture}
\]
The picture tells us that if the current state (the left position box) is $n\in\nn$, the next state (the left direction box) is $n+1\in\nn$.
Since the machine just returns its current state as a position, the sequence of positions returned will always be an increasing sequence of consecutive natural numbers starting at the initial state.
\end{example}

\begin{example}[Moving in the plane]\label{ex.R2_moore}
  Let us construct a Moore machine with positions in $\rr^2$, which we may think of as locations in the coordinate plane, and directions in $[0,\infty)\times[0,2\pi)$, which we may think of as commands to move a certain distance $r\in[0,\infty)$ at a certain angle $\theta\in[0,2\pi)$.
  We will let the machine's state-set be $\rr^\2$ as well, so the machine is a lens
  \[
    \rr^\2\yon^{\rr^\2}\to\rr^\2\yon^{[0,1]\times[0,2\pi)}.
  \]
  We can define such a lens using polyboxes:
  \[
  \begin{tikzpicture}[polybox, mapstos]
    \node[poly, dom] (s) {$(x+r\cos\theta,y+r\sin\theta)$\at$(x,y)$};
    \node[left=0pt of s_pos] {$\rr^\2$};
    \node[left=0pt of s_dir] {$\rr^\2$};

    \node[poly, cod, right=of s] (p) {$(r,\theta)\vphantom{y}$\at$(x,y)$};
    \node[right=0pt of p_pos] {$\rr^\2$};
    \node[right=0pt of p_dir] {$[0,\infty)\times[0,2\pi)$};

    \draw (s_pos) to[first] (p_pos);
    \draw (p_dir) to[last] (s_dir);
  \end{tikzpicture}
  \]
\end{example}

\begin{exercise}
Explain in words what the Moore machine in \cref{ex.R2_moore} does.
\begin{solution}
At any time, the Moore machine in \cref{ex.R2_moore} is located at a point on the coordinate plane, say $(x, y) \in \rr^\2$.
This location is its current state.
When we ask the machine to return its position, it will tell us those coordinates, since the return function is the identity.
Then if we give the machine a direction $(r, \theta)$ for some distance $r \in [0,1]$ and angle $\theta \in [0, 2\pi)$, the machine will move by that distance, at that angle counterclockwise from the positive $x$-axis, from $(x, y)$ to
\[
    (x+r\cos\theta, y+r\sin\theta) = (x,y) + r(\cos\theta, \sin\theta)
\]
(here we treat $\rr^\2$ as a vector space, so that $r(\cos\theta, \sin\theta)$ is a vector of length $r$ at the angle $\theta$).
\end{solution}
\end{exercise}

\index{Moore machine!memoryless}
\index{Moore machine!induced by function}

\begin{example}[Functions as memoryless Moore machines]\label{ex.funs_to_moore}
Given a function $f\colon A\to I$, there is a corresponding $(A,I)$-Moore machine with states $I$ that takes in an element of $A$ and returns the element of $I$ obtained by applying $f$.

It is given by the lens $I\yon^I\to I\yon^A$ defined as follows:
\[
\begin{tikzpicture}[polybox, mapstos]
  \node[poly, dom] (s) {$f(a)$\at$i$};
  \node[left=0pt of s_pos] {$I$};
  \node[left=0pt of s_dir] {$I$};

  \node[poly, cod, right=of s] (p) {$\vphantom{f}a$\at$i$};
  \node[right=0pt of p_pos] {$I$};
  \node[right=0pt of p_dir] {$A$};

  \draw (s_pos) to[first] (p_pos);
  \draw (p_dir) to[last] (s_dir);
\end{tikzpicture}
\]
That is, this lens is the identity on positions, returning the state directly as its position, and on directions it is the function $I\times A\To{\pi_2}A\To{f} I$, which ignores the current state and applies $f$ to the direction received to compute the new state.

If the machine starts in state $i_0$ and is given a sequence of directions $(a_1,a_2,\ldots)$ from $A$, the machine will return the positions $(i_0,f(a_1),f(a_2),\ldots)$. We say this machine is \emph{memoryless}, because at no point does the state of the machine actually depend on any previous states; instead, its state depends only on the last direction it received.
\end{example}

\begin{exercise}\label{exc.funs_to_moore}
Suppose we have a function $f\colon A\times I\to I$.
\begin{enumerate}
	\item Find a corresponding $(A,I)$-Moore machine $I\yon^I\to I\yon^A$.
  You may draw it out in polyboxes.
	\item Would you say the machine is memoryless?
\qedhere
\end{enumerate}
\begin{solution}
\begin{enumerate}
    \item We seek an $(A,I)$-Moore machine $I\yon^I\to I\yon^A$ corresponding to the function $f\colon A\times I\to I$.
    We know that an $(A,I)$-Moore machine $I\yon^I \to I\yon^A$ consists of a return function $I \to I$ and an update function $I \times A \to I$.
    So we can simply let the return function be the identity on $I$ and the update function be $I \times A \iso A \times I \To{f} B$, i.e.\ the function $f$ with its inputs swapped.
    In polyboxes, the machine looks like
    \[
    \begin{tikzpicture}[polybox, mapstos]
      \node[poly, dom] (s) {$f(a,i)$\at$i$};
      \node[left=0pt of s_pos] {$I$};
      \node[left=0pt of s_dir] {$I$};

      \node[poly, cod, right=of s] (p) {$\vphantom{f}a$\at$i$};
      \node[right=0pt of p_pos] {$I$};
      \node[right=0pt of p_dir] {$A$};

      \draw (s_pos) to[first] (p_pos);
      \draw (p_dir) to[last] (s_dir);
    \end{tikzpicture}
    \]

    \item Generally, such a machine is not memoryless.
    Unlike in \cref{ex.funs_to_moore}, the update function $I \times A \iso A \times I \To{f} I$ does appear to depend on its first input, namely the previous state, which $f$ takes as its second input.
    We can see this from out polybox picture above: the left direction box, which contains the new state $f(a,i)$, depends on the current state $i\in I$ in the left position box.

    However, if $f$ factors through the projection $\pi_1 \colon A \times I \to A$, i.e.\ if $f$ can be written as a composite $A \times I \To{\pi_1} A \To{f'} B$ for some $f' \colon A \to B$, then the resulting machine \emph{is} memoryless: it is the memoryless Moore machine from \cref{ex.funs_to_moore} corresponding to $f'$.
\end{enumerate}
\end{solution}
\end{exercise}

\begin{exercise}
Find $A,I\in\smset$ such that the following can be identified with a lens $S\yon^S\to I\yon^A$, and explain in words what the corresponding $(A,I)$-Moore machine does (there may be multiple possible solutions):
\begin{enumerate}
	\item a \emph{discrete dynamical system}, i.e.\ a set of states $S$ and a transition function $S\to S$ that describes how to transfer from state to state.
	\item a \emph{magma}, i.e.\ a set $S$ and a function $S\times S\to S$.
	\item a set $S$ and a subset $S'\ss S$.\qedhere
\end{enumerate}
\begin{solution}
For each of the following constructs, we find $A,I\in\smset$ such that the construct can be identified with a lens $\varphi\colon S\yon^S\to I\yon^A$, i.e.\ a return function $\varphi_\1\colon S\to B$ and an update function $\varphi^\sharp\colon S\times A\to S$.
\begin{enumerate}
  \item Given a discrete dynamical system with states $S$ and transition funtion $n\colon S\to S$, we can set $A\coloneqq I\coloneqq\1$.
  Then $\varphi_\1\colon S\to\1$ is unique, while $\varphi^\sharp\colon S\times\1\to S$ is given by $S\times\1\iso S\To{n}S$.
  The corresponding Moore machine can only be fed one direction (you could think of that direction as a button that simply says ``advance to the next state'') and can only return one position (which tells us no information).
  So it is just a set of states and a deterministic way to move from state to state.

  We could have also set $A\coloneqq\0$ and $I\coloneqq S$, so that $\varphi_\1\coloneqq n$ and $\varphi^\sharp\colon S\times\0\to S$ is unique, but this formulation is somewhat less satisfying: this is a Moore machine that never moves between its states, effectively functioning as a lookup table between whatever state the machine happens to be in and its position, which also happens to refer to some state.

  \item Given a magma consisting of a set $S$ and a function $m\colon S\times S\to S$, we can set $A \coloneqq S$ and $I\coloneqq\1$.
  Then $\varphi_\1\colon S\to\1$ is unique, while $\varphi^\sharp\colon S\times S\to S$ is equal to $m$.
  The corresponding Moore machine always returns the same position.
  It uses the binary operation $m$ to combine the current state with a given direction---which also refers to a state---to obtain the new state.

  Alternatively, we could have set the update function to be $m$ with its inputs swapped.
  The difference here is that the new state is given by applying $m$ with the direction on the left and the current state on the right, rather than the other way around.
  If $m$ is noncommutative, this would yield a different Moore machine.

  We could have also set $A\coloneqq\0$ and $I\coloneqq S^S$, so that $\varphi^\sharp\colon S\times\0\to S$ is unique, while currying $m$ gives $\varphi_\1$, so that $\varphi_\1 s$ is the function $S \to S$ given by $s' \mapsto m(s, s')$.
  Alternatively, $\varphi_\1 s$ could be the function $s' \mapsto m(s', s)$.
  Either way, this is again a Moore machine that never moves between its states, functioning as a lookup table between the machine's current state and the function $m$ partially applied to that state on one side or the other.

  \item Given a set $S$ and a subset $S' \ss S$, we can set $A\coloneqq\0$ and $I\coloneqq\2$.
  Then $\varphi^\sharp\colon S\times\0\to S$ is unique, while we define $\varphi_\1\colon S\to\2$ by
  \[
      \varphi_\1 s =
      \begin{cases}
          1 & \text{if } s \in S' \\
          2 & \text{if } s \notin S'
      \end{cases}
  \]
  so that $S'$ can be recovered from $\varphi_\1$ as its fiber over $1$.
  The corresponing Moore machine never moves between its states, but returns one of two positions indicating whether or not the current state is in the subset $S'$.
\end{enumerate}
\end{solution}
\end{exercise}

The previous examples of Moore machines mostly had identities as return functions.
In the following exercises, we will build examples of Moore machines that do not return their entire states as positions.

\begin{exercise}[Robot with health]
Think of the Moore machine in \cref{ex.R2_moore} as a robot and modify it as follows.

Add to its state a ``health meter,'' which takes a real value between 0 and 1 representing the robot's health.
Have the robot lose half its health each time it moves to a location whose $x$-coordinate is negative.
Do not return the robot's health; instead, use its health $h$ as a multiplier, allowing it to move a distance of $hr$ given an input of $r$.
\begin{solution}
The original Moore machine had states $\rr^\2$, so to add a health meter that takes values in $[0,1]$, we take the cartesian product to obtain a new set of states $\rr^\2\times[0,1]$.
The position-set and direction-set are unchanged, so the Moore machine is a lens
\[
    \rr^\2\times[0,1]\yon^{\rr^\2\times[0,1]} \to \rr^\2\yon^{[0,\infty)\times[0,2\pi)}.
\]\index{product}
Its return function $\rr^\2\times[0,1]\to\rr^\2$ is the canonical projection, as the machine returns only its location in $\rr^\2$ and not its health; while its update function
\[
    \rr^\2 \times [0,1] \times [0,\infty) \times [0,2\pi) \to \rr^\2 \times [0,1]
\]
sends $(x, y, h, r, \theta)$ to
\[
    (x + hr\cos\theta, y + hr\sin\theta, h'),
\]
where $h' = h/2$ if the machine's new $x$-coordinate $x + hr\cos\theta < 0$ and $h' = h$ otherwise.
As polyboxes, the lens is
\[
\begin{tikzpicture}[polybox, mapstos]
  \node[poly, dom] (s) {$(x + hr\cos\theta, y + hr\sin\theta, h')$\at$(x,y,h)$};
  \node[left=0pt of s_pos] {$\rr^\2\times[0,1]$};
  \node[left=0pt of s_dir] {$\rr^\2\times[0,1]$};

  \node[poly, cod, right=of s] (p) {$(r,\theta)\vphantom{hy}$\at$(x,y)$};
  \node[right=0pt of p_pos] {$\rr^\2$};
  \node[right=0pt of p_dir] {$[0,\infty)\times[0,2\pi)$};

  \draw (s_pos) to[first] (p_pos);
  \draw (p_dir) to[last] (s_dir);
\end{tikzpicture}
\]
\end{solution}
\end{exercise}

\index{Moore machine!tape of Turing machine}

\begin{exercise}[Tape of a Turing machine]
A Turing machine has a tape consisting of a cell for each integer.
Each cell bears a value $v\in V\coloneqq\{0,1,-\}$, and one of the cells $c\in\zz$ is distinguished as the ``current'' cell.
So the set of states of the tape is $V^\zz\times\zz$.

The Turing machine interacts with the tape by asking for the value of the current cell, an element of $V$; and by changing the value of the current cell before moving left (i.e.\ replacing the current cell $c\in\zz$ with the new cell $c-1$) or right (i.e.\ replacing $c$ with $c+1$).
Hence the tape's position-set is $V$ and its direction-set is $V\times\{\const{left},\,\const{right}\}$.

\begin{enumerate}
	\item If we model the tape as a Moore machine $t\colon S\yon^S\to I\yon^A$, what are $S,I,$ and $A$?
	\item Write down the specific $t$ that makes it act like a tape as specified above.
\qedhere
\end{enumerate}
\begin{solution}
\begin{enumerate}
    \item The tape has states $S\coloneqq V^\zz \times \zz$, positions $I\coloneqq V$, and directions $A\coloneqq V\times\{\const{left},\,\const{right}\}$; as a Moore machine, it is a lens
    \[
        t \colon (V^\zz \times \zz)\yon^{V^\zz \times \zz} \to V\yon^{V \times \{\const{left},\,\const{right}\}}.
    \]
    \item The return function of $t$ should give the value in the current cell of the tape.
    So $t_\1\colon V^\zz\times\zz\to V$ is the evaluation map: it sends $(f,c)$ with $f\colon\zz\to V$ and $c\in\zz$ to $f(c)\in V$.
    Then on a given direction $(v,d)\in V\times\{\const{left},\,\const{right}\}$, the update function of $t$ writes $v$ in the tape's current cell before shifting the current cell number up or down by one according to whether $d$ is $\const{right}$ or $\const{left}$.
    More precisely,
    \[
        t^\sharp \colon (V^\zz \times \zz) \times (V \times \{\const{left},\,\const{right}\}) \to V^\zz \times \zz
    \]
    sends current tape $f\colon\zz\to V$, current cell number $c\in\zz$, new value $v \in V$, and $d \in \{\const{left},\,\const{right}\}$ to the new tape $f' \colon \zz \to V$ defined by
    \[
        f'(n)\coloneqq
        \begin{cases}
            v & \text{if } n = c \\
            f(n) & \text{if } n \neq c
        \end{cases}
    \]
    and the new cell number $c-1$ if $d=\const{left}$ and $c+1$ if $d=\const{right}$.
\end{enumerate}
\end{solution}
\end{exercise}

\index{Moore machine!file reader}

\begin{exercise}[File-reader]\label{exc.file_reader}
Say that a \emph{file} of length $n$ is a function $f\colon\ord{n}\to\Set{ascii}$, where $\Set{ascii}\coloneqq\2\5\6$.
We refer to elements of $\ord{n}=\{1,\ldots,n\}$ as \emph{entries} in the file and, for each entry $i\in\ord{n}$, the value $f(i)\in\Set{ascii}$ as the \emph{character} at entry $i$.

Given a file $f$, design a file-reading Moore machine whose position-set is $\Set{ascii} + \{\const{done}\}$
and whose direction-set is
\[
\{(s,t)\mid 1\leq s\leq t\leq n\}+\{\const{continue}\}.
\]
Given a direction $(s,t)$, the file-reader should go to entry $s$ in the file and return the character at that entry.
If the given direction is instead $\const{continue}$, the file-reader should move to the next entry (i.e.\ from $s$ to $s+1$) and read that character---unless the new entry would be greater than $t$, in which case the file-reader should return $\const{done}$ until it receives another $(s,t)$ pair.
\begin{solution}
There are many options for the machine's state-set; we choose to use pairs of entries $(i,t)\in\ord{n}^\2$, where $i$ is the entry where the file-reader is currently located and $t$ is the entry where the file-reader should stop.
We will also include a $\const{stopped}$ state for when the file-reader has already stopped.
So our Moore machine is a lens
\[
    (\ord{n}^\2+\{\const{stopped}\})\yon^{\ord{n}^\2+\{\const{stopped}\}} \to (\Set{ascii} + \{\const{done}\})\yon^{\{(s,t)\mid 1\leq s\leq t\leq n\}+\{\const{continue}\}}.
\]
If the file-reader's current state is $\const{stopped}$, then the file-reader should return the position ``done.''
Otherwise, the file-reader should return the character at the entry where the file-reader is currently located.
So its return function $\ord{n}^\2+\{\const{stopped}\} \to \Set{ascii} + \{\const{done}\}$ sends $(i, t)$ to $f(i)$ and $\const{stopped}$ to $\const{done}$.
Meanwhile, the update function
\[
    (\ord{n}^\2+\{\const{stopped}\}) \times (\{(s,t)\mid 1\leq s\leq t\leq n\}+\{\const{continue}\}) \to \ord{n}^\2+\{\const{stopped}\}
\]
behaves as follows on the current state and given direction: regardless of the current state, if the given direction is a pair $(s,t)$, the new state will also be $(s,t)$.
On the other hand, if the given direction is $\const{continue}$ and the current state is a pair $(i,t)$, the new state should be the pair $(i+1,t)$ if $i+1\leq t$ and $\const{done}$ otherwise.
Finally, if the given direction is $\const{continue}$ and the current state is $\const{stopped}$, the new state should still be $\const{stopped}$.
\end{solution}
\end{exercise}

While \cref{exc.file_reader} gives us a functioning file-reader, it is rather awkward that we are still able to give the direction $\const{continue}$ even when the position is $\const{done}$, or provide a new range of entries before the file-reader has finished reading from the previous range.
In \cref{sec.poly.dyn_sys.depend_sys}, we will introduce a generalization of Moore machines to handle cases like these, where the array of directions the machine can receive changes depending on its current position.
In particular, we will be able to let the file-reader ``close its port,'' so that it cannot receive signals while it is busy reading, but open its port once it is \const{done}; see \cref{ex.generalized_file_reader}.

\subsection{Deterministic state automata}

The diagram in \cref{ex.Moore_three} may look familiar to those who have studied automata theory; in fact, a deterministic state automaton can be expressed as a Moore machine with a distinguished initial state.

\index{deterministic state automaton}\index{Moore machine!deterministic state automaton|see{deterministic state automaton}}

\begin{definition}[Deterministic state automaton, language]\label{def.dfa}
A \emph{deterministic state automaton} consists of
\begin{itemize}
	\item a set $S$ of \emph{states};
	\item a set $A$ of \emph{symbols};
	\item an \emph{update function} $u\colon S\times A\to S$;
	\item an \emph{initial state} $s_0\in S$;
	\item a subset $F\ss S$ of \emph{accept states}.
\end{itemize}
Let\index{list}
\[
  \lst(A)=\sum_{n\in\nn}A^\ord{n}
\]
denote the set of finite sequences $(a_1,\ldots,a_n)$ of symbols in $A$; we call such a sequence a \emph{word}.
We say that the automaton \emph{accepts} the word $(a_1,\ldots,a_n)$ if starting at the initial state and following the symbols in the word leads us to an accept state---or, more formally, if the sequence $(s_0,s_1,\ldots,s_n)$ defined inductively by
\[
  s_{k+1}\coloneqq u(s_k,a_{k+1})
\]
is such that $s_n$ is an accept state: $s_n\in F$.

We call a subset of $\lst(A)$ a \emph{language}, and we say that the set of all words in $\lst(A)$ that the automaton accepts is the language \emph{recognized} by the automaton.
\end{definition}

\index{deterministic state automaton!language of}

\begin{remark}
When we study a deterministic state automaton, we are usually interested in which words the automaton accepts and, more generally, what language the automaton recognizes.
While intuitive, the condition we provided for when an automaton accepts a word can be cumbersome to work with.
In \cref{ex.dsa_lang_recog}, we will give a more compact way of describing whether an automaton accepts a word and specifying the language the automaton recognizes.
Better yet, we will find that this alternative formulation arises naturally from the theory of $\poly$.
\end{remark}

\begin{proposition} \label{prop.dsa}
A deterministic state automaton with a set of states $S$ and a set of symbols $A$ can be identified with a pair of lenses
\[
  \yon\to S\yon^S\to \2\yon^A.
\]
\end{proposition}
\begin{proof}
By \cref{exc.lens-from-0-or-yon}, a lens $\yon\to S\yon^S$ can be identified with an initial state $s_0\in S$.
Then a lens $S\yon^S\to\2\yon^A$ consists of a return function $f\colon S\to\2$, which can be identified with a subset of accept states $F \ss S$, together with an update function $u \colon S\times A\to S$.
\end{proof}

In other words, we can think of a deterministic state automaton as a Moore machine with position set $\2$ along with a distinguished initial state; the Moore machine has the same states and update function as the automaton and the automaton's symbols as its directions.

Now imagine if we wanted to construct a version of this automaton that stops reading symbols (i.e.\ directions) whenever the machine enters an accept state (i.e.\ returns one position instead of the other).
To do this would require a machine whose set of possible directions is dependent on its current position.
Instead of an update function $u\colon S\times A\to S$, we would need an update function that takes a direction $a\in A$ if the state $s\in S$ is \emph{not} an accept state (say, if $f(s)=1$) but takes a direction in $\0$ (i.e.\ no direction) if the state $s$ \emph{is} an accept state (if $f(s)=2$).
So there would be one update function $u_s\colon A\to S$ if $f(s)=1$ and a different update function $u_s\colon\0\to S$ if $f(s)=2$.
But these are exactly the on-directions functions of a lens $S\yon^S\to\yon^A+\1$!
Indeed, replacing our interface monomial with a general polynomial is exactly how we will obtain our generalized dependent Moore machines.

\index{interface!monomial|)}
\index{Moore machine|)}

%-------- Section --------%
\section{Dependent dynamical systems}\label{sec.poly.dyn_sys.depend_sys}

\index{Moore machine!dependent|see{dynamical system, dependent}}
\index{dynamical system|(}\index{interface!polynomial}

Each of our Moore machines above has a monomial $I\yon^A$ as an interface.
Every representable summand of such an interface has the same representing set $A$, so the set of directions that can be fed into the machine is always $A$.
But by replacing $I\yon^A$ with an arbitrary polynomial $p$, which may have a different direction-set at each position, we can model a broader class of machines.

\begin{definition}[Dependent dynamical system]\label{def.gen_moore}
  A \emph{dependent dynamical system} (or a \emph{dependent Moore machine}, or simply a \emph{dynamical system}) is a lens \[\varphi\colon S\yon^S\to p\] for some $S\in\smset$ and $p\in\poly$.
  We call
  \begin{itemize}
    \item the domain $S\yon^S$ the machine's \emph{state system}: its position-set (equivalently, its direction-set) is the machine's \emph{state-set}, and its positions (equivalently, its directions) are the machine's \emph{states};
    \item the codomain $p$ the machine's \emph{interface}: its position-set and direction-sets are the same as the machine's \emph{position-set} and \emph{direction-sets}, and its positions and directions are the machine's \emph{positions} and \emph{directions};
    \item the on-positions function $\varphi_\1\colon S\to p(\1)$ the machine's \emph{return function};
    \item the on-directions map $\varphi^\sharp\colon p[\varphi_\1(-)]\to S$ the machine's \emph{update map}, and the on-directions function $\varphi^\sharp_s\colon p[\varphi_\1s]\to S$ at $s\in S$ the machine's \emph{update function} at $s$.
  \end{itemize}
\end{definition}
\index{state system}

\begin{example}[Dynamical systems as polyboxes]
  We can express a dynamical system $\varphi\colon S\yon^S\to p$ in polyboxes as
  
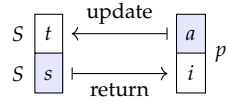
\[
    \begin{tikzpicture}[polybox, mapstos]
      \node[poly, dom] (S) {$t$\at$s\vphantom{i}$};
        \node[left=0pt of S_pos] {$S$};
        \node[left=0pt of S_dir] {$S$};

      \node[poly, cod, right=of S, "$p$" right] (p) {$a\vphantom{t}$\at$i$};

      \draw (S_pos) -- node[below] {return} (p_pos);
      \draw (p_dir) -- node[above] {update} (S_dir);
    \end{tikzpicture}
  \]
  We can visualize $\varphi$ as a channel between the internal state system on the left and the external interface on the right.
  The state system enters its current state $s\in S$ into the left position box, and the return function converts this state to a position $i\in p(\1)$ of the interface.
  At $i$, the interface has a direction-set $p[i]$; an interacting agent selects one of these directions $a\in p[i]$ to enter into the right direction box.
  Finally, the update map uses the current state $s$ and the position $i$ to fill the left direction box with the new state $t\in S$.
  Then the process repeats with $t$ in place of $s$.
\end{example}
\index{interface}

\begin{remark}
It may seem limiting that the set of possible directions a dependent dynamical system can receive should depend on the current \emph{position} rather than the current \emph{state}; but this makes sense philosophically if we accept that the system's interface should capture \emph{everything} about how it interacts with the outside world.
In particular, the system's position should capture everything an external observer could possibly perceive about the system, while the direction-set should capture all the ways in which an external agent can choose to interact with the system.
But if the set of directions available to an external agent \emph{changes}, the external agent should be able to detect this fact---the system's position must have changed as well!
On the other hand, if the internal state changes, but the external position remains the same, the agent wouldn't see any difference---they wouldn't know to interact with the system any differently, so the directions available to them would have to stay the same, too.
\end{remark}

Here are some examples of dependent dynamical systems.
We begin by finishing the example at the end of the last section.

\index{dependent Moore machine|see{dynamical system}}
\index{deterministic state automaton!halting}

\begin{example}[Halting deterministic state automata]\label{ex.regular_lang_stop}
Recall deterministic state automata from \cref{def.dfa}.
Say we want such an automaton to halt after reaching an accept state and read no more symbols.
Then rather than a lens $S\yon^S\to \2\yon^A$, we could use a lens
\[
  \varphi\colon S\yon^S\to \yon^A+\1\iso\{\const{reject}\}\yon^A+\{\const{accept}\}.
\]
To give such a lens, we first need to provide a return function $\varphi_\1\colon S\to\{\const{reject},\,\const{accept}\}$.
We let $\varphi$ send the accept states to $\const{accept}$ and every other state to $\const{reject}$.

If we reach an accept state, we want the machine to halt.
So at the position $\const{accept}$, corresponding to the summand $\1$, there are no directions available.
This makes the update function $\varphi^\sharp_s$ vacuous when $\varphi_\1s=\const{accept}$.

On the other hand, when $\varphi_\1s=\const{reject}$, the update functions $\varphi^\sharp_s\colon A\to S$ specify how the machine updates its state for each direction in $A$ if the current state is $s$.
This corresponds to the automaton's update function.

When equipped with an initial state $s_0\in S$ specified by a lens $\yon\to S\yon^S$, we call these dependent dynamical systems \emph{halting deterministic state automata}.
Given a word $(a_1,\ldots,a_n)\in\lst(A)$, we say that the automaton \emph{accepts} this word if starting at the initial state and following the elements in the sequence leads us to an accept state, \emph{without reaching an accept state any earlier}---or, more formally, if the sequence $(s_0,s_1,\ldots,s_n)$ defined inductively by
\[
  s_{k+1}\coloneqq \varphi^\sharp_{s_k}a_{k+1}
\]
is such that $s_n$ is the sequence's first accept state:
\[
  \varphi_\1s_k=\begin{cases}
    \const{reject} & \text{if } k<n \\
    \const{accept} & \text{if } k=n
  \end{cases}
\]
We call the set of all words accepted by the automaton the language \emph{recognized} by the automaton.
\end{example}

\index{deterministic state automaton!language of}

\begin{remark}
Again, the conditions for when such an automaton accepts a word are rather awkward to formally state.
We will see in \cref{ex.halt_dsa_accept} an alternative way of saying whether a word is accepted by a halting deterministic state automaton.
\end{remark}

\begin{exercise}\label{exc.halt_dsa}
Consider the halting deterministic state automaton shown below:
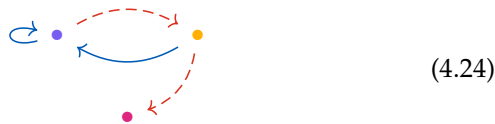
\begin{equation} \label{eqn.halt_dsa}
\begin{tikzcd}[column sep=small]
	\bul[my-lavender]\ar[rr, bend left, my-red, dashed]\ar[loop left, my-blue]&&
	\bul[my-yellow]\ar[dl, bend left, my-red, dashed]\ar[ll, my-blue, bend left]\\&
	\bul[my-magenta]
\end{tikzcd}
\end{equation}
% \begin{equation} \label{eqn.halt_dsa}
% \begin{tikzpicture}
% \node[draw] {
% \begin{tikzcd}[column sep=small]
% 	\LMO{}\ar[rr, bend left, orange]\ar[loop left, dgreen]&&
% 	\LMO{}\ar[dl, bend left, orange]\ar[ll, dgreen, bend left]\\&
% 	\LMO{}
% \end{tikzcd}
% };
% \end{tikzpicture}
% \end{equation}
Let the left state $\bul[my-lavender]$ be $1$, the right state $\bul[my-yellow]$ be $2$, and the bottom state $\bul[my-magenta]$ be $3$.
We designate $\bul[my-lavender]$, state $1$, as the initial state.
We can also call the dashed red arrows $\Red$ and the solid blue arrows $\Blue$.
Answer the following questions, in keeping with the notation from \cref{ex.regular_lang_stop}.

\begin{enumerate}
	\item What is $S$?
	\item What is $A$?
	\item Based on the labeled transition diagram, which states are accept states, and which are not?
	\item Specify the corresponding lens $S\yon^S\to\yon^A+\1$.
	\item Name a word that is accepted by this automaton.
	\item Name a word that is not accepted by this automaton.
	Why not?
	Can you find another word that is not accepted by this automaton for a different reason?
\qedhere
\end{enumerate}
\begin{solution}
\begin{enumerate}
    \item The set of states is $S\coloneqq\{1,2,3\}=\3$ (or equivalently $S\coloneqq\{\bul[my-lavender],\bul[my-yellow],\bul[my-magenta]\}\iso\3$).
    \item The set of input symbols is $A\coloneqq\{\Red,\Blue\}$.
    \item The automaton should halt at the accept states, so the accept states are exactly the states that have no arrows coming out of them---in this case, only state 3.
    States 1 and 2 are not accept states.
    \item Let the corresponding lens be $\varphi\colon S\yon^S\to\yon^A+\1$, or $\varphi\colon \3\yon^\3\to\yon^{\{\Red,\,\Blue\}}+\1$.
    According to the previous part, $\varphi$ has a return function $\varphi_\1\colon S\to\2$ sending states 1 and 2, as non-accept states, to 1; and sending state 3, as an accept state, to 2.
    Then $\varphi^\sharp_3$ is vacuous, while the other two update functions are given by the the targets of the arrows in \eqref{eqn.halt_dsa} as follows:
    \begin{align*}
        \varphi^\sharp_1(\Red)\coloneqq2&,\varphi^\sharp_1(\Blue)\coloneqq1;\\
        \varphi^\sharp_2(\Red)\coloneqq3&,\varphi^\sharp_2(\Blue)\coloneqq1,
    \end{align*}
    \item Some words accepted by this automaton include the word $(\Red,\Red),$ the word $(\Red,\Blue,\Red,\Red),$ and the word $(\Blue,\Red,\Blue,\Blue,\Blue,\Red,\Red)$.
    \item Some words are not accepted by the automaton because they lead you to a non-accept state (1 or 2); others are not accepted by the automaton because they lead you to an accept state (3) too early.
    Some examples of the former possibility include the words $(\Blue,\Blue)$ and $(\Red,\Blue,\Red,\Blue)$, while some examples of the latter possibility include the words $(\Blue,\Red,\Red,\Blue)$ and $(\Red,\Red,\Red,\Red,\Red,\Red,\Red)$.
\end{enumerate}
\end{solution}
\end{exercise}

For further examples, every graph gives rise to a dynamical system; but to ensure that we are discussing the same concept, let us fix the definition of a graph.

\begin{definition}[Graph] \label{def.graph}\index{graph}
A \emph{graph} $G \coloneqq (E \tto V)$ consists of
\begin{itemize}
  \item a set $E$ of \emph{edges};
  \item a set $V$ of \emph{vertices};
  \item a \emph{source function} $s\colon E\to V$ that assigns each edge a source vertex;
  \item a \emph{target function} $t\colon E\to V$ that assigns each edge a target vertex.
\end{itemize}
\end{definition}

So when we say ``graph,'' we mean a \emph{directed} graph, and we allow multiple edges between the same pair of vertices as well as self-loops.

\begin{example}[Graphs as dynamical systems] \label{ex.graph_dyn}\index{graph}
Given a graph $G\coloneqq(E\tto V)$ with source and target functions $s,t\colon E\to V$, there is an associated polynomial
\[
    g\coloneqq\sum_{v \in V} \yon^{s\inv(v)}.
\]
Its positions are the vertices of the graph, and its directions at $v\in V$ are the edges coming out of $v$.
We call this the \emph{emanation polynomial} of $G$.

\index{graph}

The graph itself induces a dynamical system $\varphi\colon V\yon^V\to g$, where $\varphi_\1 = \id_V$ and $\varphi^\sharp_v e = t(e)$.
So its states as well as its positions are the vertices of the graph, and a direction at a vertex $v\in V$ is an edge $e\in E$ coming out of $v$ that takes us from $v=s(e)$ along the edge $e$ to its target vertex $\varphi^\sharp_v e=t(e)$.
\end{example}

\index{graph!as dynamical system}

\begin{exercise}
Pick your favorite graph $G$, and consider the associated dynamical system as in \cref{ex.graph_dyn}.
Draw its labeled transition diagram as in \eqref{eqn.trans_diag} or \eqref{eqn.halt_dsa}.
\begin{solution}
No matter what graph you chose, \cref{ex.graph_dyn} tells us that if you were to draw the labeled transition diagram of its associated dynamical system, you would just end up with a picture of your graph!
The vertices of your graph are the states, and the edges of your graph are the possible transitions between them.
\end{solution}
\end{exercise}

% \begin{example}[Inputting an initial state]
% Suppose you have a closed system $f^\sharp\colon S\yon^S\to\yon$. The modeler can choose an initial state $\yon\to S\yon^S$, but what if we want some other system to choose the initial state? We haven't gotten to wiring diagrams yet, but the idea is to create a system that starts as not-closed---accepting as input a state $s\in S$---and then dives into its closed loop with that initial state.

% Let $S'\coloneqq S+\1$, so that the initial state $\yon\to S'\yon^{S'}$ now is canonical: it's the new $\1$. We also have a canonical inclusion $S\To{i}S'$. We will give a lens
% \[
% S'\yon^{S'}\to\yon+\yon^S
% \]
% that starts out with its outer box in the mode $\yon^S$ of accepting an $S$-input, and then moves to the mode $\yon$ so that it is a closed system forever after.

% To give a lens $S'\yon^{S'}\to\yon+\yon^S$, it is sufficient to give two morphisms: $S\yon^{S'}\to\yon$ and $\yon^{S'}\to\yon^S$. The first is equivalent to a function $S\to S'$ and we take the map $S\To{f^\sharp}S\To{i}S'$; this means that whenever we want to update the state from a state in $S$ we'll just do whatever our original closed system did. The second is also equivalent to a function $S\to S'$ and we use $i$; this means that whatever state is input at the beginning will be what we take as our first noncanonical state.
% \end{example}

\begin{example}\label{ex.generalized_file_reader}
In \cref{exc.file_reader}, we built a file-reader as a Moore machine, where a file is a function $f\colon\ord{n}\to\Set{ascii}$ from entries to characters.
Now we turn that file-reader into a dependent dynamical system $\varphi\colon S\yon^S\to p$ with only one direction while reading.

We let $S \coloneqq \{(s,t)\mid 1\leq s\leq t\leq n\}$, so that each state consists of a current entry $s$ and a terminal entry $t$.
Meanwhile, our interface $p$ will have two labeled copies of $\Set{ascii}$ as positions:
\[
    p(\1)\coloneqq\{\const{ready},\const{busy}\}\times\Set{ascii}.
\]
So each $p$-position is a pair $(m,c)$, where $c\in\Set{ascii}$ and $m$ is one of two modes: $\const{ready}$ or $\const{busy}$.
Then we define the direction-sets of $p$ for each $c\in\Set{ascii}$ as follows:
\[
    p[(\const{ready}, c)]\coloneqq S \qqand p[(\const{busy}, c)]\coloneqq\{\const{advance}\}\iso\1.
\]
That way, our file-reader can receive as its direction any pair of entries in $S$ when it is $\const{ready}$ but can only be told to $\const{advance}$ when it is $\const{busy}$.

We want our file-reader to be $\const{ready}$ if its current entry is the terminal entry; otherwise, it will be $\const{busy}$.
In either case, it will return the character at the current entry.
So we define the return function $\varphi_\1$ such that, for all $(s,t)\in S$,
\begin{align*}
  \varphi_\1(s, t) =
  \begin{cases}
    (\text{\const{ready}}, f(s)) &\mbox{if $s = t$}\\
    (\text{\const{busy}}, f(s)) &\mbox{otherwise}
  \end{cases}
\end{align*}

While the file-reader is $\const{ready}$, we want to set its new current and terminal entries to equal the given direction.
So for each $(s,s)\in S$, define the update function $\varphi^\sharp_{(s,\,s)}\colon S\to S$ to be the identity on $S$.

On the other hand, while the file-reader is $\const{busy}$, we want it to step forward through the file each time it receives an input.
So for each $(s,t)\in S$ for which $s<t$, we let the update function $\varphi^\sharp_{(s,\,t)}\colon\1\to S$ specify the element $(s+1,t)\in S$, thus shifting its current entry up by $1$.
\end{example}

\index{file reader|see{Moore machine, file reader}}
\index{Moore machine!file reader}
\index{file searcher|see{Moore machine, file searcher}}
\index{Moore machine!file searcher}

\begin{exercise} \label{exc.file_searcher}
Say instead of a file-reader, we wanted a file-searcher, which acts just like the file-reader from \cref{ex.generalized_file_reader} except that it only returns $c\in\Set{ascii}$ in its position when $c$ is a specific character; say $c=100$.
Otherwise, it returns the placeholder character $\_$.
Give the lens for this file-searcher by explicitly defining its return (on-positions) and update (on-directions) functions.
Hint: You should be able to use the same state system.
\begin{solution}
We give a file-searcher $\psi\colon S\yon^S\to q$ according to the specification as follows.
Its possible positions should form the set
\[
    q(\1)\coloneqq\{\const{ready}, \const{busy}\}\times\{100,\_\}.
\]
The direction-sets of $q$ can be defined in the same way we defined the direction-sets of $p$: for each $c\in\{100,\_\}$, we have
\[
    q[(\const{ready}, c)]\coloneqq S \qqand q[(\const{busy}, c)]\coloneqq\1.
\]
Then we set the return function $\psi_\1$ to behave like $\varphi_\1$, but with characters not equal to $100$ replaced with $\_$: so for all $(s,t)\in S$,
\begin{align*}
  \psi_\1(s,t) =
  \begin{cases}
    (\const{ready}, 100) &\mbox{if $s=t$ and $f(s)=100$}\\
    (\const{ready}, \_) &\mbox{if $s=t$ and $f(s)\neq100$}\\
    (\const{busy}, 100) &\mbox{if $s\neq t$ and $f(s)=100$}\\
    (\const{busy}, \_) &\mbox{otherwise}\\
  \end{cases}
\end{align*}
Then the update functions of $\psi$ behave just like those of $\varphi$.
For each $(s,t)\in S$ for which $s=t$, we define the update function $\psi^\sharp_{(s,\,t)}\colon S\to S$ to be the identity on $S$.
On the other hand, for each $(s,t)\in S$ for which $s\neq t$, we let the update function $\psi^\sharp_{(s,\,t)}\colon \1\to S$ specify the element $(s+1, t)\in S$, thus shifting its current entry up by $1$.
\end{solution}
\end{exercise}

In the previous exercise, we manually constructed a file-searcher that acted very much like a file-reader.
In \cref{exc.file_searcher_wrap}, we will see a simpler way to construct a file-searcher by leveraging the file-reader we have already defined.
Moreover, this construction will highlight precisely how our file-searcher is related to our file-reader.
This will be possible using \emph{wrapper interfaces}, which we will introduce in \cref{subsec.poly.dyn_sys.new.wrap}.

\index{interface!wrapper}
\index{robot|(}

\begin{example}\label{ex.grid_robot}
Choose $n\in\nn$, a \emph{grid size}, and for each $i\in\ord{n}$, let $D_i$ be the set
\[
	D_i\coloneqq
	\begin{cases}
		\{0,+1\}&\tn{ if }i=1\\
		\{-1,0,+1\}&\tn{ if } 1<i<n\\
		\{-1,0\}&\tn{ if }i=n
	\end{cases}
\]
We can think of $D_i$ as the set of ways a robot could move from location $i$.
If $1<i<n$, a robot may shift its location by $-1$ (move left/down), $0$ (remain still), or $+1$ (move right/up).
But a robot already at $i=1$ cannot shift its location by $-1$; likewise, a robot already at $i=n$ cannot shift its location by $+1$.

Then we can model a robot told to move in an $\ord{n}\times\ord{n}$ grid as a dependent dynamical system
\[
    \varphi\colon (\ord{n}\times\ord{n})\yon^{\ord{n}\times\ord{n}}\to\sum_{(i,\,j)\in\ord{n}\times\ord{n}}\yon^{D_i\times D_j}.
\]
The robot's state is a location $(i,j)\in\ord{n}\times\ord{n}$ in the grid.
We let $\varphi_\1\coloneqq\id_{\ord{n}\times\ord{n}}$ so that the dynamical system returns its state as its position: the robot expresses its position by moving to that location in the grid.

For each $(i,j)\in\ord{n}\times\ord{n}$, we let $\varphi^\sharp_{(i,\,j)}$ send each pair $(d,e)\in D_i\times D_j$ to the grid location $(i+d,j+e)\in(n,n)$.
Concretely, this says that if a robot located at $(i,j)$ receives the pair $(d,e)$ as its direction, its new position will be $(i+d,j+e)$.
As polyboxes, the dynamical system is given by
  
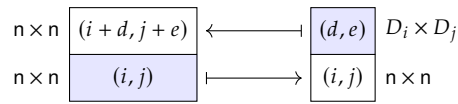
\[
\begin{tikzpicture}[polybox, mapstos]
  \node[poly, dom] (s) {$(i+d,j+e)$\at$(i,j)$};
  \node[left=0pt of s_pos] {$\ord{n}\times\ord{n}$};
  \node[left=0pt of s_dir] {$\ord{n}\times\ord{n}$};

  \node[poly, cod, right=of s] (p) {$(d,e)\vphantom{j}$\at$(i,j)$};
  \node[right=0pt of p_pos] {$\ord{n}\times\ord{n}$};
  \node[right=0pt of p_dir] {$D_i\times D_j$};

  \draw (s_pos) to[first] (p_pos);
  \draw (p_dir) to[last] (s_dir);
\end{tikzpicture}
\]

Our definition of $D_i$ for each $i\in\ord{n}$ guarantees that this position is still inside the grid; with this setup, the robot has fewer ways to move on the sides or corners of the grid than anywhere else:
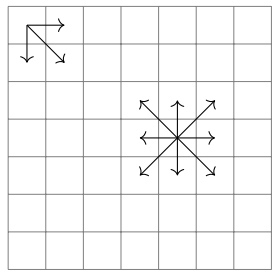
\[
\begin{tikzpicture}[scale=.5]
  \draw[step=1cm,gray,very thin] (-3,-3) grid (4,4);
	\draw[->] (-2.5,3.5) -- (-1.5, 3.5);
	\draw[->] (-2.5,3.5) -- (-1.5, 2.5);
	\draw[->] (-2.5,3.5) -- (-2.5, 2.5);
	\draw[->] (1.5, 0.5) -- (0.5, 0.5);
	\draw[->] (1.5, 0.5) -- (2.5, 0.5);
	\draw[->] (1.5, 0.5) -- (1.5, 1.5);
	\draw[->] (1.5, 0.5) -- (1.5,-0.5);
	\draw[->] (1.5, 0.5) -- (0.5, 1.5);
	\draw[->] (1.5, 0.5) -- (0.5,-0.5);
	\draw[->] (1.5, 0.5) -- (2.5, 1.5);
	\draw[->] (1.5, 0.5) -- (2.5,-0.5);
\end{tikzpicture}
\]
In this picture, $\ord{n}\coloneqq\7$, and the $4$ directions at position $(1,7)$ and $9$ directions at position $(5,4)$ are shown (recall that remaining still is an option in either case).

Note that in this example, the positions are literally the positions in the grid where the robot could be, and the directions at each position are literally the directions in which the robot can move!
\end{example}

\begin{exercise} \label{exc.grid_reward}
Modify the dynamical system from \cref{ex.grid_robot} as follows.
\begin{enumerate}
	\item Replace the interface with a new polynomial $p$ so that at each grid location, the robot can receive not only the direction it should move in but also a ``reward value'' $r\in\rr$.
	\item Replace the state-set with a new set $S$ so that an element $s\in S$ may include both the robot's position and a list of all reward values so far.
	\item With your new $p$ and $S$, define a new dynamical system $\varphi'\colon S\yon^S\to p$ that preserves the behavior of the dynamical system from \cref{ex.grid_robot} while updating the robot's reward list without returning it externally.
\qedhere
\end{enumerate}
\begin{solution}
We modify the dynamical system from \cref{ex.grid_robot}.
\begin{enumerate}
    \item Previously, the direction-set at position $(i,j)\in\ord{n}\times\ord{n}$ of our interface was $D_i\times D_j$.
    But now we also want to give the robot a ``reward value'' $r\in\rr$.
    So our new direction-set should be $D_i\times D_j\times\rr$:
    \[
        p\coloneqq\sum_{(i,\,j)\in\ord{n}\times\ord{n}} \yon^{D_i\times D_j\times\rr}.
    \]
    \item Previously, a state was just a location in the grid: an element of $\ord{n}\times\ord{n}$.
    But now we want to be able to record a list of reward values as well.
    Since each reward value is a real number, we define the state-set to be $S\coloneqq\ord{n}\times\ord{n}\times\lst(\rr)$.\index{list}
    \item The former return function $\varphi_\1$ was the identity on $\ord{n}\times\ord{n}$.
    The new return function $\varphi'_\1$ should still just return the robot's current grid position; but since it is now a function from $S=\ord{n}\times\ord{n}\times\lst(\rr)$, it should instead be the canonical projection $\varphi'_\1\colon \ord{n}\times\ord{n}\times\lst(\rr)\to\ord{n}\times\ord{n}$.

    For each former state $(i,j)\in\ord{n}\times\ord{n}$, the former update function $\varphi^\sharp_{(i,\,j)}\colon D_i\times D_j\to\ord{n}\times\ord{n}$ sent $(d,e)\mapsto(i+d,j+e)$.
    With an extra component $(r_1,\ldots,r_k)\in\lst(\rr)$ of the state, the new update function $(\varphi')^\sharp_{(i,\,j,\,(r_1,\,\ldots,\,r_k))}\colon D_i\times D_j\times\rr\to\ord{n}\times\ord{n}\times\lst(\rr)$ sends $(d,e,r)\mapsto(i+d,j+e,(r_1,\ldots,r_k,r))$, also updating the list of rewards.
    As polyboxes, the new dynamical system is given by
    \[
    \begin{tikzpicture}[polybox, mapstos]
      \node[poly, dom] (s) {$(i+d,j+e,(r_1,\ldots,r_k,r))$\at$(i,j,(r_1,\ldots,r_k))$};
      \node[left=0pt of s_pos] {$\ord{n}\times\ord{n}\times\lst(\rr)$};
      \node[left=0pt of s_dir] {$\ord{n}\times\ord{n}\times\lst(\rr)$};

      \node[poly, cod, right=of s] (p) {$(d,e,r)\vphantom{j}$\at$(i,j)$};
      \node[right=0pt of p_pos] {$\ord{n}\times\ord{n}$};
      \node[right=0pt of p_dir] {$D_i\times D_j\times\rr$};

      \draw (s_pos) to[first] (p_pos);
      \draw (p_dir) to[last] (s_dir);
    \end{tikzpicture}
    \]
\end{enumerate}
\end{solution}
\end{exercise}

In the previous exercise, we added a reward system to the robot on the grid by manually redefining the associated lens.
But there is a simpler way to think about the new system: it is the juxtaposition of two systems, a robot system and a reward system, in parallel.
We will see how to express this in terms of lenses in \cref{exc.grid_reward_par}, once we explain how to juxtapose systems like this in general in \cref{subsec.poly.dyn_sys.new.par}.
In fact, we will see in \cref{exc.grid_robot_par} that the robot-on-a-grid system itself can be viewed as the juxtaposition of two systems, and this perspective will provide a structured way to generalize \cref{ex.grid_robot} to more than two dimensions.

\index{robot|)}
\index{dynamical system|)}

%-------- Section --------%
\section{New dynamical systems from old}\label{sec.poly.dyn_sys.new}

\index{dynamical system!constructing new from old|(}

We have seen how dependent dynamical systems can be modeled as lenses in $\poly$ of the form $S\yon^S\to p$.
But we have yet to take full advantage of the categorical structure that $\poly$ provides.
In particular, based only on what we know of $\poly$ so far from \cref{ch.poly.cat}, we have three rather different ways of constructing new dynamical systems from old ones:
\begin{enumerate}
    \item Given dynamical systems $S\yon^S\to p$ and $S\yon^S\to q$, we can use the universal property of the \emph{categorical product} to obtain a dynamical system $S\yon^S\to p\times q$; see \cref{subsec.poly.dyn_sys.new.prod}.
    \item Given dynamical systems $\varphi\colon S\yon^S\to p$ and $\psi\colon T\yon^T\to q$, we can take their parallel product to obtain a dynamical system $\varphi\otimes\psi\colon ST\yon^{ST}\to p\otimes q$; see \cref{subsec.poly.dyn_sys.new.par}.
    \item Given a dynamical system $\varphi\colon S\yon^S\to p$ and a lens $f\colon p\to q$, we can compose them to obtain a dynamical system $\varphi\then f\colon S\yon^S\to q$; see \cref{subsec.poly.dyn_sys.new.wrap}.
\end{enumerate}
Each of these operations has a concrete interpretation in terms of the systems' behavior.
In this section, we will review each of them in turn.

%---- Subsection ----%
\subsection[Categorical products of interfaces]{Categorical products: multiple interfaces operating on the same states}
\label{subsec.poly.dyn_sys.new.prod}

\index{interface!multiple}

\index{dynamical system!constructing from product}

Let $I$ be a set, and say that we have an $I$-indexed family of dependent dynamical systems $(\varphi_i\colon S\yon^S\to p_i)_{i\in I}$ that all share the same state-set $S$.
Then since $\poly$ has all small products, the universal property of products induces a lens \[\varphi\colon S\yon^S\to\prod_{i\in I}p_i,\] which is itself a dynamical system with state-set $S$.
By \eqref{eqn.poly_prod}, the interface of $\varphi$ (i.e.\ the product that is its codomain) has position-set
\[
  \left(\prod_{i\in I}p_i\right)\!(\1)\iso\prod_{i\in I}\left(p_i(\1)\right)
\]
and, at each position $\bar{j}\colon(i\in I)\to p_i(\1)$, direction-set
\[
  \sum_{i\in I}p_i[\bar{j}i].
\]
We then characterize the dynamics of $\varphi$ in terms of each $\varphi_i$ by leveraging the universal property of products in $\poly$, detailed in the binary case in the solution to \cref{exc.poly_prod}, as follows.\index{dynamics}
The return function
\[
  \varphi_\1\colon S\to\prod_{i\in I}\left(p_i(\1)\right)
\]
sends each state $s\in S$ to the dependent function $(i\in I)\to p_i(\1)$ sending $i\in I$ to the position $(\varphi_i)_\1s$ returned by the corresponding dynamical system $\varphi_i$ at the state $s$.
Then the update function
\[
  \varphi^\sharp_s\colon\sum_{i\in I}p_i[(\varphi_i)_\1s]\to S
\]
at $s\in S$ sends each pair $(i,d)$ in its domain, with $i\in I$ and $d\in p_i[(\varphi_i)_\1s]$, to where the update function of $\varphi_i$ at $s$ sends the direction $d$: namely $(\varphi_i)^\sharp_s d$.
We can write $\varphi$ using polyboxes as follows:
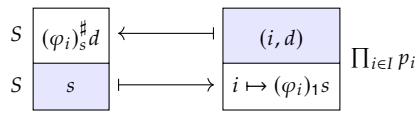
\[
\begin{tikzpicture}[polybox, mapstos]
  \node[poly, dom] (S) {$(\varphi_i)^\sharp_sd$\at$s\vphantom{(\varphi_i)_\1s}$};
  \node[left=0pt of S_pos] {$S$};
  \node[left=0pt of S_dir] {$S$};

  \node[poly, cod, right=of S, "$\prod_{i\in I}p_i$" right] (p) {$(i,d)\vphantom{(\varphi_i)^\sharp_sd}$\at$i\mapsto(\varphi_i)_\1s$};

  \draw (S_pos) to[first] (p_pos);
  \draw (p_dir) to[last] (S_dir);
\end{tikzpicture}
\]

In other words, if there are multiple interfaces that drive the same set of states, we may view them as a single product interface that drives those states.
This single dynamical system returns positions in all of the original systems at once; then it can receive a direction from any one of the original systems' direction-sets and update its state accordingly.
It is as though all of the dynamical systems can detect the current state of the combined system, but only one of them can change the state at a time.
So products give us a universal way to combine multiple polynomial interfaces into one.

\begin{exercise}
Given a set $I$, suppose we have an $(A_i,B_i)$-Moore machine with state-set $S$ for each $i\in I$.
Show that there is an induced $\left(\sum_{i\in I}A_i, \prod_{i\in I}B_i\right)$-Moore machine, again with state-set $S$.
\begin{solution}
We are given a lens $S\yon^S\to B_i\yon^{A_i}$ for each $i\in I$.
The universal property of products in $\poly$ then induces a lens
\[
    S\yon^S\to\prod_{i\in I}B_i\yon^{A_i}.
\]
By \eqref{eqn.poly_prod}, its codomain is the product of monomials
\[
  \prod_{i\in I}B_i\yon^{A_i}\iso\left(\prod_{i\in I}B_i\right)\yon^{\sum_{i\in I}A_i},
\]
which, in particular, is still a monomial.
Hence the induced lens is a $\left(\sum_{i\in I}A_i,\prod_{i\in I}B_i\right)$-Moore machine with state-set $S$.
\end{solution}
\end{exercise}\index{Moore machine}

\begin{example} \label{ex.prod_diagrams}
Consider two four-state dependent dynamical systems $\varphi\colon\4\yon^\4\to\rr\yon^{\{\Red,\,\Blue\}}$ and $\psi\colon\4\yon^\4\to \zz_{\geq0}\yon^{\{\Black,\,\Yellow\}}+\zz_{<0}\yon^{\{\Black\}}$, drawn below as labeled transition diagrams, with $\Red$ solid arrows and $\Blue$ dashed arrows on the left and $\Black$ dotted arrows and $\Yellow$ solid arrows:
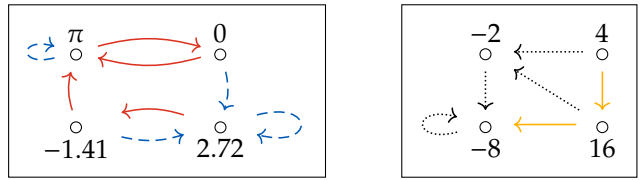
\[
\begin{tikzpicture}
	\node[draw] (1) {
  \begin{tikzcd}[row sep=15pt]
  	\LMO{\pi}\ar[r, bend left=15pt, my-red]\ar[loop left=15pt, my-blue, dashed]&
  	\LMO{0}\ar[l, bend left=15pt, my-red]\ar[d, bend left=15pt, my-blue, dashed]\\
  	\LMO[under]{-1.41}\ar[u,bend left=15pt, my-red]\ar[r, bend right=15pt, my-blue, dashed]&
  	\LMO[under]{2.72}\ar[l, bend right=15pt, my-red]\ar[loop right=15pt, my-blue, dashed]
  \end{tikzcd}
	};
	\node[draw, right=of 1] {
  \begin{tikzcd}[row sep=15pt]
  	\LMO{-2}\ar[d, densely dotted]&
  	\LMO{4}\ar[l, densely dotted]\ar[d, my-yellow]\\
  	\LMO[under]{-8}\ar[loop left, densely dotted]&
  	\LMO[under]{16}\ar[ul, densely dotted]\ar[l, my-yellow]
  \end{tikzcd}
  };
 \end{tikzpicture}
\]

The universal property of products provides a unique way to put these systems together to obtain a dynamical system $\4\yon^\4\to\rr\zz_{\geq0}\yon^{\{\Red,\,\Blue,\,\Black,\,\Yellow\}}+\rr\zz_{<0}\yon^{\{\Red,\,\Blue,\,\Black\}}$ that looks like this:
\[
\begin{tikzpicture}
	\node[draw] (1) {
  \begin{tikzcd}
  	\LMO{(\pi,-2)}\ar[r, bend left=15pt, my-red]\ar[loop left, my-blue, dashed]\ar[d, bend left=15pt, densely dotted]&
  	\LMO{(0,4)}\ar[l, bend left=15pt, my-red]\ar[d, bend left=15pt, my-blue, dashed]\ar[d, bend right=15pt, my-yellow]\ar[l, densely dotted]\\
  	\LMO[under]{(-1.41,-8)}\ar[u,bend left=15pt, my-red]\ar[r, bend right=15pt, my-blue, dashed]\ar[loop left, densely dotted]&
  	\LMO[under]{(2.72,16)}\ar[l, bend right=15pt, my-red]\ar[l, my-yellow]\ar[loop right=15pt, my-blue, dashed]\ar[ul, densely dotted]
  \end{tikzcd}
  };
\end{tikzpicture}
\]
Each state now returns two positions: one according to the return function of $\varphi$, and another according to the return function of $\psi$.
As for the possible directions, we can now choose either a direction of $\varphi$ (either $\Red$ or $\Blue$), in which case the dynamical system will update its state according to the update map of $\varphi$; or a direction of $\psi$ (either $\Black$ or sometimes $\Yellow$), in which case the dynamical system will update its state according to the update map of $\psi$.
\end{example}

\index{event-based systems}

\begin{exercise}[Toward event-based systems]
Let $\varphi\colon S\yon^S\to p$ be a dynamical system.
We can think of it as requiring a direction at each time step to update its state.

Suppose we want to change $\varphi$ into an \emph{event-based system}: one that does not always receive a direction and changes state only when it does.
That is, we want every position of $p$ to have an extra direction that, when selected, never changes the state.
We want a new system $\varphi'\colon S\yon^S\to p'$ that has this behavior; what should $p'$ and $\varphi'$ be?
\begin{solution}
Given a dynamical system $\varphi\colon S\yon^S\to p$, we seek a new dynamical system $\varphi'\colon S\yon^S\to p'$ with an additional direction that does not change the state.
We can think of this as having two different interfaces acting on the same system: the original interface $p$ of $\varphi$ and a new interface with only one possible direction that does not change the state.
This latter interface also need not distinguish between its positions; it should have a single position that provides no additional information.
So the second interface we want acting on $S\yon^S$ is $\yon$.

If $\yon$ were the only interface acting on the state system, we would have a Moore machine $\epsilon\colon S\yon^S\to\yon$ whose return function is the unique function $S\to\1$ and whose update function should be the identity function on $S$, since the direction never changes the state.
Then $p'$ is the product of the two interfaces $p$ and $\yon$, while $\varphi'\colon S\yon^S\to p'$ is the unique lens induced by $\varphi\colon S\yon^S\to p$ and $\epsilon\colon S\yon^S\to\yon$.
In particular, $p'\iso p\yon\iso\sum_{i\in p(\1)}\yon^{p[i]+\1}$, and if we let $\ast$ denote the additional direction not in $p$ at each position of $p\yon$, we can write $\varphi'$ in polyboxes as
\[
\begin{tikzpicture}[polybox, mapstos]
  \node[poly, dom] (S) {$\varphi^\sharp_sa$\at$s\vphantom{\varphi_\1s}$};
  \node[left=0pt of S_pos] {$S$};
  \node[left=0pt of S_dir] {$S$};

  \node[poly, cod, right=of S, "$p\yon$" right] (p) {$a\vphantom{\varphi^\sharp_sa}$\at$\varphi_\1s$};

  \draw (S_pos) to[first] (p_pos);
  \draw (p_dir) to[last] (S_dir);
\end{tikzpicture}
\]
when the chosen direction is from the original $p$ (i.e.\ $a\in p[\varphi_\1s]$ above), coinciding with the behavior of the original system $\varphi$; or as
\[
\begin{tikzpicture}[polybox, mapstos]
  \node[poly, dom] (S) {$s$\at$s\vphantom{\varphi_\1s}$};
  \node[left=0pt of S_pos] {$S$};
  \node[left=0pt of S_dir] {$S$};

  \node[poly, cod, right=of S, "$p\yon$" right] (p) {$\ast$\at$\varphi_\1s$};

  \draw (S_pos) to[first] (p_pos);
  \draw (p_dir) to[last] (S_dir);
\end{tikzpicture}
\]
when the chosen direction is the additional direction $\ast$ that does not change the state.
It turns out this construction is universal; it is known as \emph{copointing}.
\end{solution}
\end{exercise}

\index{polynomial functor!copointed}

%---- Subsection ----%
\subsection{Parallel products: juxtaposing dynamical systems}\label{subsec.poly.dyn_sys.new.par}

\index{dynamical system!juxtaposing}
\index{parallel product!and dynamical systems|(}

Another way to combine two polynomials---and indeed two lenses---is by taking their parallel product, as in \cref{def.parallel}.
In particular, the parallel product of two state systems is still a state system.
So parallel products give us another way to create new dynamical systems from old ones.
The procedure is straightforward: take the product of the state-sets as the new state-set, the product of the position-sets as the new position-set, and the product of the direction-sets at each of position in a tuple of positions as the new direction-set at that tuple.

For $n\in\nn$, say that we have $n$ dynamical systems: a lens $\varphi_i\colon S_i\yon^{S_i}\to p_i$ for each $i\in\ord{n}$.
Then we can take the parallel product of all of them to obtain a lens
\[
  \varphi\colon\bigotimes_{i\in\ord{n}}\left(S_i\yon^{S_i}\right)\to\bigotimes_{i\in\ord{n}}p_i.
\]
By inductively applying \cref{exc.general_poly_parallel_times}, we find that the domain of $\varphi$ is
\[
  \bigotimes_{i\in\ord{n}}\left(S_i\yon^{S_i}\right)\iso\left(\prod_{i\in\ord{n}}S_i\right)\yon^{\prod_{i\in\ord{n}}S_i},
\]
so $\varphi$ is itself a dependent dynamical system with state-set $\prod_{i\in\ord{n}}S_i$.
Meanwhile, inductively applying \eqref{eqn.parallel_def} yields
\[
  \left(\bigotimes_{i\in\ord{n}}p_i\right)\!(\1)\iso\prod_{i\in\ord{n}}\left(p_i(\1)\right)
\]
as the position-set and, at each position $(j_i)_{i\in\ord{n}}$,
\[
  \prod_{i\in\ord{n}}p_i[j_i]
\]
as the direction-set of the interface of $\varphi$.

We can characterize the dynamics of $\varphi$ in terms of each constituent dynamical system $\varphi_i$ as follows.
By our proof sketch of \cref{prop.parallel_monoidal}, the return function
\[
  \varphi_\1\colon \prod_{i\in\ord{n}}S_i\to\prod_{i\in\ord{n}}\left(p_i(\1)\right)
\]
sends each $n$-tuple of states $(s_i)_{i\in\ord{n}}$ in its domain, with each $s_i\in S_i$, to the $n$-tuple of positions $((\varphi_i)_\1s_i)_{i\in\ord{n}}$ returned by each of the constituent dynamical systems at each state.
Then at the $n$-tuple of states $(s_i)_{i\in\ord{n}}\in\prod_{i\in\ord{n}}S_i$, the update function
\[
  \varphi^\sharp_{(s_i)_{i\in\ord{n}}}\colon\prod_{i\in\ord{n}}p_i[(\varphi_i)_\1s_i]\to\prod_{i\in\ord{n}}S_i
\]
sends each $n$-tuple of directions $(d_i)_{i\in\ord{n}}$ in its domain, with each $d_i\in p_i[(\varphi_i)_\1s_i]$, to the $n$-tuple $((\varphi_i)^\sharp_{s_i}d_i)_{i\in\ord{n}}$ consisting of states to which the update function of each $\varphi_i$ at $s_i$ sends $d_i$.
We can write $\varphi$ using polyboxes as follows:
\[
\begin{tikzpicture}[polybox, mapstos]
  \node[poly, dom] (S) {$((\varphi_i)^\sharp_{s_i}d_i)_{i\in\ord{n}}$\at$(s_i)_{i\in\ord{n}}\vphantom{((\varphi_i)_\1s_i)_{i\in\ord{n}}}$};
  \node[left=0pt of S_pos] {$\prod_{i\in\ord{n}}S_i$};
  \node[left=0pt of S_dir] {$\prod_{i\in\ord{n}}S_i$};

  \node[poly, cod, right=of S, "$\bigotimes_{i\in\ord{n}}p_i$" right] (p) {$(d_i)_{i\in\ord{n}}\vphantom{((\varphi_i)^\sharp_{s_i}d_i)_{i\in\ord{n}}}$\at$((\varphi_i)_\1s_i)_{i\in\ord{n}}$};

  \draw (S_pos) to[first] (p_pos);
  \draw (p_dir) to[last] (S_dir);
\end{tikzpicture}
\]

In other words, multiple dynamical systems running in parallel can be thought of as a single dynamical system.
This system stores the states of all the constituent systems at once and returns positions from all of them together; then it can receive directions from all of the constituent systems' direction-sets at those positions at once and update each constituent state accordingly.
So parallel products give us a way to juxtapose multiple dynamical systems in parallel to form a single system.

\begin{exercise}
Given $n\in\nn$, suppose we have an $(A_i,B_i)$-Moore machine with state-set $S_i$ for every $i \in \ord{n}$.
Show that there is an induced $\left(\prod_{i \in \ord{n}} A_i, \prod_{i \in \ord{n}} B_i\right)$-Moore machine with state-set $\prod_{i \in \ord{n}} S_i$.
\begin{solution}
We are given a lens $S_i\yon^{S_i}\to B_i\yon^{A_i}$ for each $i\in\ord{n}$.
By inductively applying \cref{exc.general_poly_parallel_times}, we find that their parallel product in $\poly$ is a lens
\[
  \left(\prod_{i\in\ord{n}}S_i\right)\yon^{\prod_{i\in\ord{n}}S_i}\iso\bigotimes_{i\in\ord{n}}S_i\yon^{S_i}\to\bigotimes_{i\in\ord{n}}B_i\yon^{A_i}\iso\left(\prod_{i\in\ord{n}}B_i\right)\yon^{\prod_{i\in\ord{n}}A_i},
\]
which is a $\left(\prod_{i \in \ord{n}} A_i, \prod_{i \in \ord{n}} B_i\right)$-Moore machine with state-set $\prod_{i\in\ord{n}}S_i$.
This works because the parallel product of monomials is still a monomial.
\end{solution}
\end{exercise}\index{Moore machine}

\begin{example} %%CHECK \label{ex.par_diagrams}
Consider dependent dynamical systems $\varphi\colon\2\yon^\2\to \rr_{<0}\yon^{\{\const{b},\,\const{r}\}}+\rr_{\geq0}\yon^{\{\const{b}\}}$ and $\psi\colon\3\yon^\3\to\zz_{<0}\yon^{\{\const{r}\}}+\{0\}\yon^{\{\const{r},\,\const{y}\}}+\zz_{>0}\yon^{\{\const{y}\}}$, drawn below as labeled transition diagrams (we think of $\const{b},\const{r},$ and $\const{y}$ as $\const{b}$lue and dotted, $\const{r}$ed and solid, and $\const{y}$ellow and dashed, respectively):
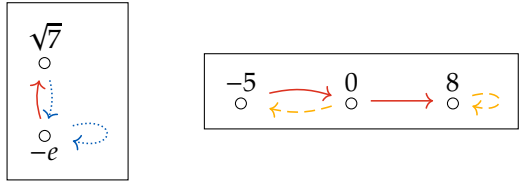
\[
\begin{tikzpicture}
	\node[draw] (1) {
  \begin{tikzcd}[row sep=15pt]
  	\LMO{\sqrt{7}}\ar[d, bend left=15pt, my-blue, densely dotted]\\
  	\LMO[under]{-e}\ar[u, bend left=15pt, my-red]\ar[loop right=15pt, my-blue, densely dotted]
  \end{tikzcd}
	};
	\node[draw, right=of 1] {
  \begin{tikzcd}[row sep=15pt]
  	\LMO{-5}\ar[r, bend left=15pt, my-red]&
  	\LMO{0}\ar[l, bend left=15pt, my-yellow, dashed]\ar[r, my-red]&
  	\LMO{8}\ar[loop right=15pt, my-yellow, dashed]
  \end{tikzcd}
  };
 \end{tikzpicture}
\]
Taking their parallel product, we obtain a dynamical system with state system $\6\yon^\6$ and interface
\begin{align*}
    &\rr_{<0}\zz_{<0}\yon^{\{(\const{b},\,\const{r}),\,(\const{r},\,\const{r})\}}+\rr_{<0}\{0\}\yon^{\{(\const{b},\,\const{r}),\,(\const{b},\,\const{y}),\,(\const{r},\,\const{r}),\,(\const{r},\,\const{y})\}}\\
    +\:&\rr_{<0}\zz_{>0}\yon^{\{(\const{b},\,\const{y}),\,(\const{r},\,\const{y})\}}+\rr_{\geq0}\zz_{<0}\yon^{\{(\const{b},\,\const{r})\}}\\
    +\:&\rr_{\geq0}\{0\}\yon^{\{(\const{b},\,\const{r}),\,(\const{b},\,\const{y})\}}+\rr_{\geq0}\zz_{>0}\yon^{\{(\const{b},\,\const{y})\}}
\end{align*}
that looks like this (we use purple and dotted to indicate $(\const{b},\const{r})$; red and solid to indicate $(\const{r},\const{r})$; green, dotted, and dashed to indicate $(\const{b},\const{y})$; and orange and dashed to indicate $(\const{r},\const{y})$):
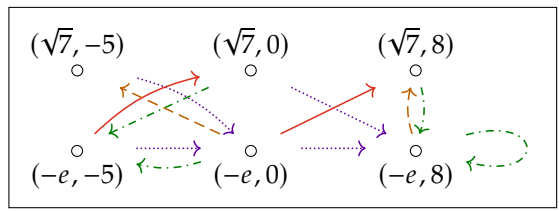
\[
\begin{tikzpicture}
	\node[draw] (1) {
  \begin{tikzcd}
    \LMO{(\sqrt{7},-5)}\ar[dr, bend left=15pt, blue!50!purple, densely dotted] &
    \LMO{(\sqrt{7},0)}\ar[dl, green!50!black, dash dot]\ar[dr, blue!50!purple, densely dotted] &
    \LMO{(\sqrt{7},8)}\ar[d, bend left=15pt, green!50!black, dash dot] \\
    \LMO[under]{(-e,-5)}\ar[ur, bend left=15pt, my-red]\ar[r, blue!50!purple, densely dotted] &
    \LMO[under]{(-e,0)}\ar[ul, orange!75!black, dashed]\ar[ur, my-red]\ar[l, bend left=15pt, green!50!black, dash dot]\ar[r, blue!50!purple, densely dotted] &
    \LMO[under]{(-e,8)}\ar[u, bend left=15pt, orange!75!black, dashed]\ar[loop right=15pt, green!50!black, dash dot]
  \end{tikzcd}
  };
\end{tikzpicture}
\]
Each state---actually a pair of states from the constituent state-sets---returns two positions, one according to the return function of $\varphi$ and another according to the return function of $\psi$.
Then every direction must be a pair of directions from the constituent interfaces at those positions, with the update function updating each state in the pair according to each direction in the pair via the constituent update functions of $\varphi$ and $\psi$.
\end{example}

\begin{exercise} \label{exc.grid_reward_par}
Explain how the dynamical system $\varphi'\colon S'\yon^{S'}\to p'$ you built in \cref{exc.grid_reward} can be expressed as the parallel product of the robot-on-a-grid dynamical system $\varphi\colon S\yon^S\to p$ from \cref{ex.grid_robot} with another dynamical system, $\psi\colon T\yon^T\to q$.
Be sure to specify $T, q,$ and $\psi$.
\begin{solution}
We will show that taking the parallel product of the robot-on-a-grid dynamical system $\varphi\colon S\yon^S\to p$ from \cref{ex.grid_robot} and a reward-tracking dynamical system $\psi\colon T\yon^T\to q$ we will define shortly yields the dynamical system $\varphi'\colon S'\yon^{S'}\to p'$ from \cref{exc.grid_reward}.

The reward-tracking dynamical system should have states in $\lst(\rr)$ to record a list of reward values, unchanging position, and directions in $\rr$ to give new reward values.
So it is the lens $\lst(\rr)\yon^{\lst(\rr)}\to\yon^\rr$ that has a uniquely defined return function, while its update map sends each state $(r_1,\ldots,r_k)\in\lst(\rr)$ and each direction $r\in\rr$ to the new state $(r_1,\ldots,r_k,r)$.

Then the dynamical system from \cref{exc.grid_reward} is the parallel product of the robot-on-a-grid dynamical system from \cref{ex.grid_robot} with the reward-tracking dynamical system $\lst(\rr)\yon^{\lst(\rr)}\to\yon^\rr$, as can be seen in the solution to \cref{exc.grid_reward}.
\end{solution}
\end{exercise}

\begin{exercise} \label{exc.grid_robot_par}
\begin{enumerate}
    \item Explain how the robot-on-a-grid dynamical system $\varphi\colon S\yon^S\to p$ from \cref{ex.grid_robot} can be written as the parallel product of some dynamical system with itself.
    \item Use $k$-fold parallel products to generalize \cref{ex.grid_robot} to robots on $k$-dimensional grids.\qedhere
\end{enumerate}
\begin{solution}
\begin{enumerate}
    \item The robot-on-a-grid dynamical system from \cref{ex.grid_robot} can be written as the parallel product of two robot-on-a-line dynamical systems of the form $\lambda\colon\ord{n}\yon^{\ord{n}}\to\sum_{i\in\ord{n}}\yon^{D_i}$, where $\lambda_\1\coloneqq\id_{\ord{n}}$ and $\lambda^\sharp_i$ for each $i\in\ord{n}$ sends each direction $d\in D_i$ to the position on the line given by $i+d$.
    This yields a robot that can move along a single axis, and the parallel product of this robot with itself yields a robot that can move along two different axes at once, which is precisely our robot-on-a-grid dynamical system.
    \item To create a dynamical system consisting of a robot moving in a $k$-dimensional grid of size $n$ along every dimension, we just take the $k$-fold parallel product of the dynamical system $\lambda\colon\ord{n}\yon^{\ord{n}}\to\sum_{i\in\ord{n}}\yon^{D_i}$ we just defined to obtain a dynamical system \[\lambda^{\otimes k}\colon\ord{n}^\ord{k}\yon^{\ord{n}^\ord{k}}\to\sum_{(i_1,\,\ldots,\,i_k)\in\ord{n}^\ord{k}}\yon^{\prod_{j\in\ord{k}}D_{i_j}}.\]
    In fact, we could have used a different $n_j$ for each $j\in\ord{k}$ instead of $n$ to obtain a robot moving in an arbitrary $k$-dimensional grid of size $n_1\times\cdots\times n_k$ as a $k$-fold parallel product.
\end{enumerate}
\end{solution}
\end{exercise}
\index{robot}

Intuitively, the parallel product takes two dynamical systems and puts them in the same room together so that they can be run at the same time.
But it does not allow for any interaction \emph{between} the two systems.
For that, we will need to use what we call a wrapper interface.
We will introduce wrapper interfaces in the next section before describing how they can be used in conjunction with parallel products to model general interaction in \cref{sec.poly.dyn_sys.interact}.

\index{parallel product!and dynamical systems|)}

%---- Subsection ----%
\subsection{Composing lenses: wrapper interfaces}\label{subsec.poly.dyn_sys.new.wrap}

\index{dynamical system!wrapper interface|see{interface, wrapper}}\index{interface!wrapper|(}

Given a dynamical system $\varphi\colon S\yon^S\to p$, say that we want to interact with its state system using a new interface $q$ rather than $p$.
We can do this whenever we have a lens $f\colon p\to q$, which we could compose with our original dynamical system to obtain a new system $S\yon^S\To{\varphi}p\To{f}q$.
We call the lens $f$ the \emph{wrapper} and its codomain $q$ the \emph{wrapper interface}, which we \emph{wrap} around $\varphi$ (or sometimes just $p$, if a dynamical system $\varphi$ has yet to be specified) using $f$.

How does this new composite system $\varphi\then f$ relate to the original dynamical system $\varphi$?
The lens $f$ converts a position $i$ from $p$ to a position $f_\1i$ from $q$; at the same time, it allows the choice of direction from $p[i]$ to depend on a choice of direction from $q[f_\1i]$, converting directions of the wrapper interface $q$ to directions of the original interface $p$.
Precomposing $f$ with a dynamical system yields a new dynamical system that lets an agent interact with the original system using only this new interface wrapped around it.

\begin{example} \label{ex.wrap_diagrams} %%CHECK
Consider a dependent dynamical system $\varphi\colon\6\yon^\6\to p$ with
\[
    p\coloneqq\{1\}\yon^{\{b,y,r\}}+\{2\}\yon^{\{b,r\}}+\{3\}\yon^{\{b\}}+\{4\}\yon^{\{r\}},
\]
drawn below as a labeled transition diagram (we think of $b,y,$ and $r$ as blue, yellow, and red, respectively):
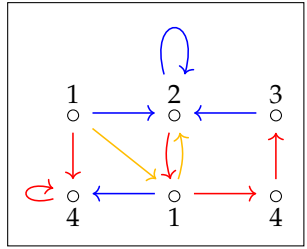
\[
\begin{tikzpicture}
	\node[draw] (1) {
  \begin{tikzcd}
    \LMO{1}\ar[r, blue]\ar[dr, dyellow]\ar[d, red] &
    \LMO{2}\ar[loop above=5pt, blue]\ar[d, bend right=15pt, red] &
    \LMO{3}\ar[l, blue] \\
    \LMO[under]{4}\ar[loop left=15pt, red] &
    \LMO[under]{1}\ar[l, blue]\ar[u, bend right=15pt, dyellow]\ar[r, red] &
    \LMO[under]{4}\ar[u, red]
  \end{tikzcd}
  };
\end{tikzpicture}
\]
We will wrap the interface \[q\coloneqq\{a\}\yon^{\{g,p,o\}}+\{b\}\yon^{\{g,p\}}+\{c\}\] around $\varphi$ using the following lens $f\colon p\to q$ (we think of $g,p,$ and $o$ as green, purple, and orange, respectively):
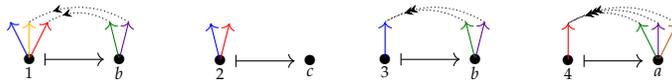
\[
\begin{tikzpicture}
	\node (p1) {
	\begin{tikzpicture}[trees, sibling distance=2.5mm]
    \node["\tiny 1" below] (1) {$\bullet$}
      child[blue] {coordinate (11)}
      child[dyellow] {coordinate (12)}
      child[red] {coordinate (13)};
    \node[right=1 of 1, "\tiny $b$" below] (2) {$\bullet$}
      child[green!50!black] {coordinate (21)}
      child[blue!50!purple] {coordinate (22)};
    \draw[|->, shorten <= 3pt, shorten >= 3pt] (1) -- (2);
    \begin{scope}[densely dotted, bend right, decoration={markings, mark=at position 0.75 with \arrow{stealth}}]
      \draw[postaction={decorate}] (21) to (13);
      \draw[postaction={decorate}] (22) to (12);
    \end{scope}
  \end{tikzpicture}
	};
	\node (p2) [below right=-1 and .7 of p1] {
	\begin{tikzpicture}[trees, sibling distance=2.5mm]
    \node["\tiny 2" below] (1) {$\bullet$}
      child[blue] {coordinate (11)}
      child[red] {coordinate (12)};
    \node[right=1 of 1, "\tiny $c$" below] (2) {$\bullet$};
    \draw[|->, shorten <= 3pt, shorten >= 3pt] (1) -- (2);
  \end{tikzpicture}
	};
	\node (p3) [right=3 of p1] {
	\begin{tikzpicture}[trees, sibling distance=2.5mm]
    \node["\tiny 3" below] (1) {$\bullet$}
      child[blue] {coordinate (11)};
    \node[right=1 of 1, "\tiny $b$" below] (2) {$\bullet$}
      child[green!50!black] {coordinate (21)}
      child[blue!50!purple] {coordinate (22)};
    \draw[|->, shorten <= 3pt, shorten >= 3pt] (1) -- (2);
    \begin{scope}[densely dotted, bend right, decoration={markings, mark=at position 0.75 with \arrow{stealth}}]
      \draw[postaction={decorate}] (21) to (11);
      \draw[postaction={decorate}] (22) to (11);
    \end{scope}
  \end{tikzpicture}
	};
	\node (p4) [right=.7 of p3] {
	\begin{tikzpicture}[trees, sibling distance=2.5mm]
    \node["\tiny 4" below] (1) {$\bullet$}
      child[red] {coordinate (11)};
    \node[right=1 of 1, "\tiny $a$" below] (2) {$\bullet$}
      child[green!50!black] {coordinate (21)}
      child[blue!50!purple] {coordinate (22)}
      child[orange!75!black] {coordinate (23)};
    \draw[|->, shorten <= 3pt, shorten >= 3pt] (1) -- (2);
    \begin{scope}[densely dotted, bend right, decoration={markings, mark=at position 0.75 with \arrow{stealth}}]
      \draw[postaction={decorate}] (21) to (11);
      \draw[postaction={decorate}] (22) to (11);
      \draw[postaction={decorate}] (23) to (11);
    \end{scope}
  \end{tikzpicture}
	};
\end{tikzpicture}
\]
Composing $\varphi$ with $f$, we obtain a dynamical system $\6\yon^\6\To{\varphi}p\To{f}q$ that looks like this:
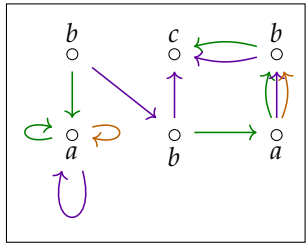
\[
\begin{tikzpicture}
	\node[draw] (1) {
  \begin{tikzcd}
    \LMO{b}\ar[d, green!50!black]\ar[dr, blue!50!purple] &
    \LMO{c} &
    \LMO{b}\ar[l, bend right=15pt, green!50!black]\ar[l, bend left=15pt, blue!50!purple] \\
    \LMO[under]{a}\ar[loop left=15pt, green!50!black]\ar[loop below=5pt, blue!50!purple]\ar[loop right=15pt, orange!75!black] &
    \LMO[under]{b}\ar[r, green!50!black]\ar[u, blue!50!purple] &
    \LMO[under]{a}\ar[u, bend left=20pt, green!50!black]\ar[u, blue!50!purple]\ar[u, bend right=20pt, orange!75!black]
  \end{tikzcd}
  };
\end{tikzpicture}
\]
Each state returns a $q$-position $j$ according to where the on-positions function of $f$ sends the $p$-position $i$ that the state returns.
Then each $q[j]$-direction is sent to a $p[i]$-direction via the on-directions function of $f$ at $i$, and the update function of $\varphi$ uses this $p[i]$-direction to compute the new state.
So $f$ allows us to operate $\varphi$ with the wrapper interface $q$ instead of the original interface $p$.
\end{example}

\begin{exercise} \label{exc.file_searcher_wrap}
In \cref{exc.file_searcher}, we built a file-searcher $\psi\colon S\yon^S\to q$ by taking the file-reader $\varphi\colon S\yon^S\to p$ from \cref{ex.generalized_file_reader} and replacing its interface $p$ with a new interface $q$ while keeping its state system $S\yon^S$ the same.
Express this construction as wrapping $q$ around $\varphi$ by giving a lens $f\colon p\to q$ for which composing $\varphi$ with $f$ yields $\psi$.
\begin{solution}
We give a lens $f\colon p\to q$ for which composing the file-reader $\varphi\colon S\yon^S\to p$ from \cref{ex.generalized_file_reader} with $f$ yields the file-searcher $\psi\colon S\yon^S\to q$ from \cref{exc.file_searcher}.
The file-searcher returns the same position as the file-reader when the second coordinate of that position is $100$, but replaces the second coordinate with a blank $\_$ otherwise.
So the on-positions function of $f$ should send each $(m,c)\in p(\1)$ to
\[
    f_\1(m,c) \coloneqq
        \begin{cases}
            (m,c) & \text{if } c = 100 \\
            (m,\_) & \text{otherwise}.
        \end{cases}
\]
Then the file-searcher acts just like the file-reader does on inputs, so every on-directions function of $f$ should be an identity function.
\end{solution}
\end{exercise}

\begin{example}[Polybox pictures of wrapper interfaces]
  In polyboxes, composing a dynamical system $\varphi\colon S\yon^S\to p$ with a wrapper $f\colon p\to q$ looks like this:
  
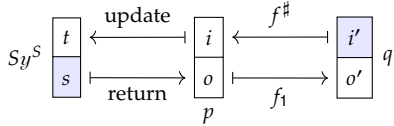
\begin{equation*}
    \begin{tikzpicture}[polybox, mapstos]
      \node[poly, dom, "$S\yon^S$" left] (S) {$t$\at$s\vphantom{o'}$};

      \node[poly, right=of S, "$p$" below] (p) {$i$\at$o\vphantom{o'}$};

      \node[poly, cod, right=of p, "$q$" right] (q) {$i'$\at$o'$};

      \draw (S_pos) -- node[below] {return} (p_pos);
      \draw (p_dir) -- node[above] {update} (S_dir);
      \draw (p_pos) -- node[below] {$f_\1$} (q_pos);
      \draw (q_dir) -- node[above] {$f^\sharp$} (p_dir);
    \end{tikzpicture}
  \end{equation*}
  The position $o$ displayed by the intermediary interface $p$ is instead exposed as a position $f_\1(o)=o'$ of the wrapper interface $q$ in the rightmost position box.
  Moreover, the direction box of $p$ is no longer blue: an agent who wishes to interact with the middle interface $p$ can only do so via the rightmost interface $q$.
  The on-directions function of the wrapper at $o$ converts a direction $i'\in q[o']$ from the rightmost direction box into a direction $i\in p[o]$.

  Picture the agent standing to the right of all the polyboxes (i.e.\ ``outside'' of the system) with their attention directed leftward (i.e.\ ``inward''), receiving positions from the white position box and feeding directions into the blue direction box.
  To an agent who is unaware of its inner workings, the composite dynamical system $\varphi\then f$ might as well look like this:
  
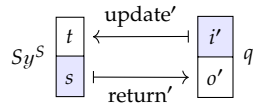
\begin{equation*}
    \begin{tikzpicture}[polybox, mapstos]
      \node[poly, dom, "$S\yon^S$" left] (S) {$t$\at$s\vphantom{o'}$};

      \node[poly, cod, right=of S, "$q$" right] (q) {$i'$\at$o'$};

      \draw (S_pos) -- node[below] {return$'$} (q_pos);
      \draw (q_dir) -- node[above] {update$'$} (S_dir);
    \end{tikzpicture}
  \end{equation*}
\end{example}

In the next section, we describe a special kind of wrapper.

%---- Subsection ----%
\subsection{Sections as wrappers}\label{subsec.poly.dyn_sys.new.sit_encl}

Say we wanted to model a dynamical system $\varphi\colon S\yon^S\to p$ within a closed system, for which an external agent can perceive no change in position and effect no change in direction.
We can think of this as wrapping $\yon$, the interface with one position and one direction, around $\varphi$.
To do so, we must specify a wrapper $\gamma\colon p\to\yon$.
In the language of \cref{def.sec-bun}, this is precisely a \emph{section} of $p$.
As we noted then, this name is appropriate, since $\gamma$ a way of sectioning off the interface $p$ from the outside world.

Recall that a section $\gamma\colon p\to\yon$ can be identified with a dependent function of the form $(i\in p(\1))\to p[i]$ that sends each $p$-position $i$ to a $p[i]$-direction, fixing a direction at every position of $p$.
So a section for an interface dictates the direction it receives given any position it inhabits; there is no need for any further outside interference.

\index{dependent function}

\begin{exercise} \label{exc.enclosures_as_functions}
Let $\varphi\colon S\yon^S\to B\yon^A$ be an $(A,B)$-Moore machine.
\begin{enumerate}
	\item Is it true that a section $\gamma\colon B\yon^A\to\yon$ can be identified with a function $A\to B$?
	\item Describe how to interpret a section $\gamma\colon B\yon^A\to\yon$ as a wrapper around an interface $B\yon^A$.
	\item Given a section $\gamma$, describe the dynamics of the composite Moore machine \[S\yon^S\To{\varphi}B\yon^A\To{\gamma}\yon\] obtained by wrapping $\yon$ around $\varphi$ using $\gamma$.
\qedhere
\end{enumerate}\index{dynamics}
\begin{solution}
\begin{enumerate}
	\item No, it represents a function $B\to A$!
	A section sends each position $b\in B$ to a direction $a\in A$.
	\item As a wrapper around an interface $B\yon^A$, a section $\gamma\colon B\yon^A\to\yon$ corresponds to a function $g\colon B\to A$ that feeds the direction $g(b)\in A$ into the system whenever it returns the position $b\in B$.
	\item Composing our original Moore machine $S\yon^S\to B\yon^A$ with a section $\gamma$ yields a Moore machine $S\yon^S\To{\varphi}B\yon^A\To{\gamma}\yon$ that returns unchanging output and receives unchanging input.
	If we identify the Moore machine with its return function $S\to B$ and its update function $S\times A\to S$, and if we identify the section $\gamma$ with a function $g\colon B\to A$, then the composite Moore machine $S\yon^S\To{\varphi}B\yon^A\To{\gamma}\yon$ can be identified with a function $S\to S$, equal to the composite
	\[
	    S\To{\Delta}S\times S\To{\id_S\times\text{return}}S\times B\To{\id_S\times g}S\times A\To{\text{update}}S,
	\]
	where $\Delta$ is the diagonal map $s\mapsto(s,s)$.
	This composite map $S\to S$ sends every state to the next according to the position the original state returns, the direction that the section gives in response to that position, and the update function that sends the original state and the selected direction to the new state.
\end{enumerate}
\end{solution}
\end{exercise}

\begin{example}[The do-nothing section] \label{ex.do_nothing}
There is something rather off-putting about the way we model dynamical systems as lenses $\varphi\colon S\yon^S\to p$.
We know that $\varphi$ sends states-as-positions to positions of the interface and, at each state-as-position, sends directions of the interface to states-as-directions.
But we rely only on the labels of elements in $S$ to tell us which positions and directions refer to the same states!

Nothing inherent in the language of $\poly$ makes these associations between states-as-positions and states-as-directions for us; we have to rely on the position-set and direction-sets of the state system being the same set for the machine to work properly.
Put another way, the monomials $\{4,6\}\yon^{\{4,6\}},\{4,6\}\yon^{\{4,8\}},$ and $\{3,5\}\yon^{\{6,7\}}$ are all isomorphic in $\poly$, but the first can be a state system while the other two cannot!

To address this issue, we need a way to connect the positions of a polynomial to its own directions in the language of $\poly$.
This is where sections can help: a lens $S\yon^S\to\yon$ is just a way of assigning to each position in $S$ a direction in $S$.
So we can define $\epsilon\colon S\yon^S\to\yon$ to be the section that sends each position $s\in S$ to the direction $s$ at $s$ with the same name, corresponding to the same state.
Note that $\epsilon$ can be identified with the identity function on $S$ (see \cref{exc.enclosures_as_functions}).
Now $\poly$ knows for each position which direction is associated with the same state the position is.

In this way, we can generalize our notion of state systems to monomials $S\yon^{S'}$ equipped with a bijection $S\to S'$, which we can then translate to a section $\epsilon\colon S\yon^{S'}\to\yon$.
But for convenience of notation, we will continue to identify the position-set of a state system with each of its direction-sets.

More concretely, the section $\epsilon\colon S\yon^S\to\yon$ acts as a very special (if rather unexciting) dynamical system: it is the \emph{do-nothing section}, with only one possible position and one possible direction that always keeps the current state the same.
While the system literally does nothing, we do know one key fact about it: given any state system, regardless of the state-set, we can always define a do-nothing section on it.\tablefootnote{It may have bothered you that we call $S\yon^S$, which is a single polynomial, a state \emph{system}, when we also use the word ``system'' to refer to dependent dynamical systems, which are lenses $S\yon^S\to p$.
The existence of the do-nothing section explains why our terminology does not clash: every polynomial $S\yon^S$ comes equipped with a dependent dynamical system $\epsilon\colon S\yon^S\to\yon$, so it really is a state \emph{system}.
%Later on we will see other systems that come with $S\yon^S$ for free.
}

Yet this is not the whole story.
The do-nothing section knows, at each position, the direction that keeps the system at the same state; but it does not know which of the other directions at that position correspond to which of the other states of the system.
We are still relying on the labels of direction-sets being the same for that: for instance, the polynomials $\{1',2',3'\}\yon^{\{1,2,3\}}$ and $\{1'\}\yon^{\{0,1,4\}}+\{2'\}\yon^{\{2,5,6\}}+\{3'\}\yon^{\{-8,-1,3\}}$ are isomorphic, but even though each has a do-nothing section matching $1'\mapsto1,2'\mapsto2,3'\mapsto3$ that makes the first one into a state system, we do not have a way to tell $\poly$ how to make the second one a state system yet.

From another perspective, $\epsilon\colon S\yon^S\to\yon$ does nothing, while $\varphi\colon S\yon^S\to p$ does ``one thing'': it steps through the system once, producing the current state's position with the return function and taking in directions with the update function.
It is ready to take another step, but how does $\poly$ know which state to visit next?
Is there a lens that does ``two things,'' ``$n$ things,'' or ``arbitrarily many things''?
Can we actually $\emph{run}$ a dynamical system in $\poly$?
We will develop the machinery to answer these questions over the course of \cref{part.comon}, starting in \cref{subsec.comon.comp.def.dyn_sys}.
\end{example}

\begin{example}[Polybox pictures of sections as wrappers]
  In polyboxes, composing a system $S\yon^S\to p$ with a section $\gamma$ of $p$ can be depicted as
  
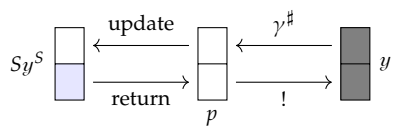
\begin{equation*}
    \begin{tikzpicture}[polybox, tos]
      \node[poly, dom, "$S\yon^S$" left] (S) {};

      \node[poly, right=of S, "$p$" below] (p) {};

      \node[poly, identity, right=of p, "$\yon$" right] (yon) {};

      \draw (S_pos) -- node[below] {return} (p_pos);
      \draw (p_dir) -- node[above] {update} (S_dir);
      \draw (p_pos) -- node[below] {$!$} (yon_pos);
      \draw (yon_dir) -- node[above] {$\gamma^\sharp$} (p_dir);
    \end{tikzpicture}
  \end{equation*}
  or, equivalently, as
  
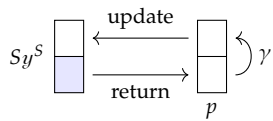
\begin{equation*}
    \begin{tikzpicture}[polybox, tos]
      \node[poly, dom, "$S\yon^S$" left] (S) {};

      \node[poly, right=of S, "$p$" below] (p) {};

      \draw (S_pos) -- node[below] {return} (p_pos);
      \draw (p_dir) -- node[above] {update} (S_dir);
      \draw (p_pos) to[climb'] node[right] {$\gamma$} (p_dir);
    \end{tikzpicture}
  \end{equation*}
  Remember: in a polybox depiction of a dynamical system, the world outside the system exists to the right of all the boxes.
  So the first picture represents $\yon$ as a gray wall, cutting off any interaction between the system to its left and the world to its right.
  Meanwhile, the second picture illustrates how a sectioned-off system independently selects directions of the intermediary interface $p$ via $\gamma$ according to the $p$-positions that the inner (leftward) system $S\yon^S\to p$ returns.
  While the second picture shows us why the closed system neither seeks nor requires external directions, the first picture helps remind us that any returned $p$-positions never reach the outside world either.
  The composite system is therefore equivalent to the section drawn as follows:
  
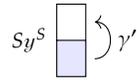
\begin{equation*}
    \begin{tikzpicture}[polybox, tos]
      \node[poly, dom, "$S\yon^S$" left] (S) {};

      \draw (S_pos) to[climb'] node[right] {$\gamma'$} (S_dir);
    \end{tikzpicture}
  \end{equation*}
  In lens parlance, $\gamma'\colon S\yon^S\to\yon$ is the original system $S\yon^S\to p$ composed with $\gamma\colon p\to\yon$; in the language of dependent functions, $\gamma'\colon S\to S$ is given by
  \[
  \gamma'(s)=\oper{update}(s,\gamma(\oper{return}(s))) \qquad \text{for all }s\in S,
  \]
  where we interpret $\gamma$ as a dependent function $(i\in p(\1))\to p[i]$.
  We can deduce this equation by equating the previous two polybox pictures, knowing they represent the same lens:
  
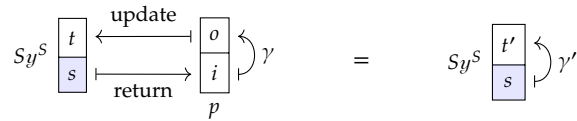
\[
  \begin{tikzpicture}
    \node (1) {
      \begin{tikzpicture}[polybox, mapstos]
        \node[poly, dom, "$S\yon^S$" left] (S) {$t$\at$s$};

        \node[poly, right=of S, "$p$" below] (p) {$o$\at$i$};

        \draw (S_pos) -- node[below] {return} (p_pos);
        \draw (p_dir) -- node[above] {update} (S_dir);
        \draw (p_pos) to[climb'] node[right] {$\gamma$} (p_dir);
      \end{tikzpicture}
    };
    \node[right=1.8 of 1] (2) {
      \begin{tikzpicture}[polybox, mapstos]
        \node[poly, dom, "$S\yon^S$" left] (S) {$t'$\at$s$};

        \draw (S_pos) to[climb'] node[right] {$\gamma'$} (S_dir);
      \end{tikzpicture}
    };
    \node at ($(1.east)!.5!(2.west)$) {=};
  \end{tikzpicture}
  \]
  Equating the directions boxes of the domain on either side, we have that $t=t'$, so
  \[
  \gamma'(s)=\oper{update}(s,o)=\oper{update}(s,\gamma(i))=\oper{update}(s,\gamma(\oper{return}(s))).
  \]
  Later on we will read more intricate equations off of polyboxes in this manner, although we will not spell out the procedure in so much detail; we encourage you to trace through the arrows on your own.
\end{example}

\begin{example}[The do-nothing section in polyboxes]  \label{ex.do_nothing_polybox}
  In \cref{ex.do_nothing}, we saw that every state system $S\yon^S$ can be equipped with a section $\epsilon\colon S\yon^S\to\yon$ called the do-nothing section, which assigns each state-as-position to its corresponding state-as-direction, thus leaving the state unchanged.
  That is, it is the section whose polyboxes can be drawn as follows:
  
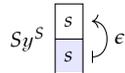
\begin{equation*}
    \begin{tikzpicture}[polybox, mapstos]
      \node[poly, dom, "$S\yon^S$" left] (S) {$s$\at$s$};

      \draw (S_pos) to[climb'] node[right] {$\epsilon$} (S_dir);
    \end{tikzpicture}
  \end{equation*}
\end{example}

\index{dynamical system!constructing new from old|)}

%-------- Section --------%
\section{General interaction}\label{sec.poly.dyn_sys.interact}

We are now ready to use $\poly$ to model interactions between dependent dynamical systems that can change their interfaces and interaction patterns.
\index{interaction|(}
\index{interaction pattern|see{interaction}}

\subsection[Wrapping juxtaposed systems together]{Wrapping juxtaposed dynamical systems together}
\index{parallel product!and interaction}

When wrapper interfaces are used in conjunction with parallel products, they may encode multiple interacting dynamical systems as a single system.
Explicitly, given $n\in\nn$ and a dynamical system $\varphi_i\colon S_i\yon^{S_i}\to p_i$ for each $i\in\ord{n}$, we can first juxtapose them to form a single dynamical system $\varphi\then f$:
\[\varphi\colon \left(\prod_{i\in\ord{n}}S_i\right)\yon^{\prod_{i\in\ord{n}}S_i}\to\bigotimes_{i\in\ord{n}}p_i\]
by taking their parallel product.
Then we can wrap an interface $q$ around $\varphi$ using a wrapper $f\colon\bigotimes_{i\in\ord{n}}p_i\to q$, yielding a new composite dynamical system
\[\left(\prod_{i\in\ord{n}}S_i\right)\yon^{\prod_{i\in\ord{n}}S_i}\To{\varphi}\bigotimes_{i\in\ord{n}}p_i\To{f}q.\]
On positions, $f$ gives a way of combining all the positions of the constituent interfaces into a single position of the wrapper interface.
On directions, $f$ takes into account the current positions of each constituent interface along with a direction for the wrapper interface to specify a direction for each of the constituent interfaces.
In particular, a judiciously chosen on-directions function could feed positions of some interfaces as directions to others.
When $f$ is a wrapper around a parallel product of interfaces, we call $f$ the \emph{interaction pattern} between those interfaces.

\begin{example}[Repeater]\index{dynamical system!repeater}
Suppose we have a dynamical system $\varphi\colon S\yon^S\to A\yon+\yon$ that takes unchanging directions and sometimes returns elements of $A$ while other times returning only silence (the position associated with the right hand summand $\yon$).
What if we wanted to construct a system $\psi$ that operates just like $\varphi$, but \emph{always} returns elements of $A$ as output?
Where $\varphi$ would have returned silence, we want $\psi$ to instead \emph{repeat} the last element of $A$ that was returned. (We allow $\psi$ to repeat an arbitrary element of $A$ if $\varphi$ returns silence before it has returned any elements of $A$ yet.)

What we need is a way to store an element of $A$ and retrieve it as needed.
So whenever $\varphi$ returns an element of $A$, we store it; then when $\varphi$ returns silence, we retrieve the last element of $A$ we stored and return that instead.

What should this storage-retrieval dynamical system look like?
It needs to take elements of $A$ as directions, return elements of $A$ as positions, and store elements of $A$ as states.
In fact, the identity lens $\iota\colon A\yon^A\to A\yon^A$ works perfectly: it returns the element of $A$ currently stored as its position and updates its state to the direction it receives.

Now we can juxtapose our original system $\varphi$ with the storage-retrieval system $\iota$ by taking their parallel product, yielding a dynamical system
\[
    \varphi\otimes\iota\colon S\yon^S\otimes A\yon^A\iso SA\yon^{SA}\to (A\yon+\yon)\otimes A\yon^A
\]
that runs both systems in parallel: simultaneously and independently.
But what we want is for $\varphi$ and $\iota$ to interact with each other, and for the resulting system to only return elements of $A$.
To do so, we need to wrap an interface $A\yon$ around $\varphi\otimes\iota$ by composing it with some lens
\[
    f\colon(A\yon+\yon)\otimes A\yon^A\to A\yon,
\]
the interaction pattern between the interfaces $A\yon+\yon$ and $A\yon^A$ that we must define.
Then we will be able to define $\psi\coloneqq(\varphi\otimes\iota)\then f$.

Since $\otimes$ distributes over $+$, by the universal property of the coproduct, it suffices to give lenses
\[
    g\colon(A\times A)\yon^A\to A\yon \qqand h\colon A\yon^A\to A\yon.
\]
The former corresponds to the case where $\varphi$ returns an element of $A$, while the latter corresponds to the case where $\varphi$ is silent.

When $\varphi$ returns an element of $A$, we want the composite system $\psi$ to return that same element, but we also want to give that position as a direction to $\iota$ so that it can be stored.
We do not need to do anything with the position returned by $\iota$; we can simply discard it.
So $g$ should map $(a,a')\mapsto a$ on positions, yielding the position returned by $\varphi$ and discarding the position returned by $\iota$; and the on-directions function $g^\sharp_{(a,\,a')}\colon\1\to A$ should pick out the direction $a\in A$, feeding the position returned by $\varphi$ as a direction for $\iota$.

On the other hand, when $\varphi$ returns silence, we want $\psi$ to return the position of $\iota$ instead.
We also need to feed this position back to $\iota$ as a direction so that it can continue to be stored.
So $h$ should be the identity on positions as well as the identity on directions.
\end{example}

\begin{example}[Paddling]\label{ex.paddler}
Say we wanted to build a Moore machine with interface $\nn\yon$; we will interpret the natural number position it returns as describing the machine's current location.
Suppose we want to be very strict about what how far the machine can move and what can make it move.

To model this, we introduce two intermediary systems, which we call the \emph{paddler} and the \emph{tracker}:%
\tablefootnote{Perhaps one could refer to the tracker as the \emph{demiurge}; it is responsible for maintaining the material universe.}
\[
  \text{paddler}\colon S\yon^S\to\2\yon
  \qqand
  \text{tracker}\colon T\yon^T\to\nn\yon^\2
\]
The paddler has interface $\2\yon$ because it is blind (i.e.\ takes no directions) and can only move its paddle (i.e.\ return a position) $\const{left}$ or $\const{right}$: its position-set is $\2\iso\{\const{left},\,\const{right}\}$. The tracker has interface $\nn\yon^\2$ because it will return the location of the machine (as an element $n\in\nn$) as its position and take in the position of the paddler (as an element of $\2$) as its direction.
We can wrap an interface $\nn\yon$ around them both using an interaction pattern
\[
    \2\yon\otimes\nn\yon^\2\iso\2\nn\yon^\2\to\nn\yon
\]
whose on-positions function is the canonical projection $\2\nn\to\nn$, returning the location returned by the tracker, and whose on-directions map is the projection $\2\nn\to\2$, passing the position of the paddler as a direction to the tracker.

\index{dynamics}

Let us leave the paddler's dynamics alone---how you the paddler may behave is arbitrary---and instead focus on the dynamics of the tracker.
We want it to watch for when the paddle switches from $\const{left}$ to $\const{right}$ or from $\const{right}$ to $\const{left}$; at that moment it should push the machine forward one unit. Thus the states of the tracker are given by $T\coloneqq\2\nn$, storing what side the paddler is on and the machine's current location.
The on-positions function of the tracker is the canonical projection $\2\nn\to\nn$ that returns the current location; then at each $(d,i)\in\2\nn$, the on-directions function of the tracker $\2\to\2\nn$ sends
\[
  d'\mapsto
	\begin{cases}
		(d',i)&\tn{if }d=d'\\
		(d',i+1)&\tn{if }d\neq d',
	\end{cases}
\]
storing the new position of the paddler as well as moving the machine forward one unit if the paddle switches while keeping the machine still if the paddle stays still.
\end{example}

\index{dynamical system!paddler}
\begin{exercise}
Change the dynamics and state system of the tracker in \cref{ex.paddler} so that it exhibits the following behavior.

When the paddle switches once and stops, the tracker increases the location by one unit and stops, as before in \cref{ex.paddler}. But when the paddle switches twice in a row, the tracker increases the location by two units on the second switch! So if the paddler is stable for a while, then switches three times in a row, the tracker will increase the location by one, then two, then two again.
\begin{solution}
We define a new tracker $T'\yon^{T'}\to\nn\yon^\2$ based on the one from \cref{ex.paddler} to watch for when the paddle switches sides once, at which point the tracker should increase its location by one, and watch for when the paddle switches sides twice in a row, at which point the tracker should increase its location by two.
To do this, we need the tracker to remember not just the current side the paddle is on, but the previous side the paddle was on as well.
The tracker should still remember the current location.
Thus the states of the tracker are given by $T\coloneqq\2\times\2\times\nn$, storing the previous side the paddler was on, the current side the paddler is on, and the current location.
The on-positions function of the tracker is the canonical projection $\2\times\2\times\nn\to\nn$ that returns the current location; then at each $(d,d',i)\in\2\times\2\times\nn$, the on-directions function of the tracker $\2\to\2\times\2\times\nn$ sends
\[
  d''\mapsto
	\begin{cases}
		(d',d'',i)&\tn{if }d'=d''\\
		(d',d'',i+1)&\tn{if }d'\neq d''\tn{ and }d=d'\\
		(d',d'',i+2)&\tn{if }d'\neq d''\tn{ and }d\neq d'
	\end{cases}
\]
storing both the last side the paddle was on and the new side the paddle is on as well as moving the machine forward one unit if the paddle switches after not switching and two units if the paddle switches after just switching.
\end{solution}
\end{exercise}

\begin{example}
Suppose you have two systems with the same interface $p\coloneqq q\coloneqq\rr^\2\yon^{\rr^\2-\{(0,0)\}}$.
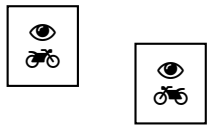
\[
\begin{tikzpicture}
	\node (m1) {\faMotorcycle};
	\node[above=-.15 of m1] (e1) {\faEye};
	\node[draw, thick, fit = (m1) (e1)] {};
	\node[below right=0 and 1 of m1] (m2) {\scalebox{-1}[1]{\faMotorcycle}};
	\node[above=-.15 of m2] (e2) {\faEye};
	\node[draw, thick, fit = (m2) (e2)] {};
\end{tikzpicture}
\]
The ordered pair comprising the position of each interface indicates the location of the corresponding system, while the range of possible directions indicate the locations that the system could observe, relative to the location of the system itself.
Taking all pairs of reals except $(0,0)$ corresponds to the fact that the eye cannot see anything at the same location as the eye itself.

Let us have the two systems approach each other accelerating at a rate equal to the reciprocal of the squared distance between them, modeled with discrete time units.
If they finally collide, let us have both systems halt.
To do this, we want the wrapper interface to be $\{\const{go}\}\yon+\{\const{stop}\}$, so that if the system returns $\const{go}$, it can still advance to the next state; but if it returns $\const{stop}$, it halts.
The wrapper $\rr^\2\yon^{\rr^\2-\{(0,0)\}}\otimes\rr^\2\yon^{\rr^\2-\{(0,0)\}}\to\{\const{go}\}\yon+\{\const{stop}\}$ is given on positions by
\[
  \big((x_\1,y_1),(x_2,y_2)\big)\mapsto
	\begin{cases}
		\const{stop}&\mbox{ if $x_1=x_2$ and $y_1=y_2$}\\
		\const{go}&\mbox{ otherwise}.
	\end{cases}
\]
On directions, we use the function
\[
  \big((x_1,y_1),(x_2,y_2)\big)\mapsto \big((x_2-x_1,y_2-y_1),(x_1-x_2,y_1-y_2)\big),
\]
so that each system is able to see the location of the other system relative to its own, i.e.\ the vector pointing from itself to the other system (unless that vector is zero, in which case the system should have returned the $\const{stop}$ position and halted).

\index{dynamics}

We can use these vectors to define the internal dynamics of each system so that they move the way we want them to.
Each system will hold as its internal state its current location and velocity as vectors, i.e.\ $S\coloneqq\rr^\2\times\rr^\2$.
To define a lens $S\yon^S\to\rr^\2\yon^{\rr^\2-\{(0,0)\}}$ we simply return the current location, update the current location by adding the current velocity vector, and update the current velocity vector by adding an acceleration vector with appropriate magnitude pointing to the other system:
\[
\begin{tikzpicture}[polybox, mapstos]
  \node[poly, dom] (S) {$\left(x+v_x,y+v_y,v_x+\frac{a}{(a^2+b^2)^{3/2}},v_y+\frac{b}{(a^2+b^2)^{3/2}}\right)$\at$((x,y),(v_x,v_y))$};
    \node[below=0pt of S_pos] {$\rr^2\times\rr^2$};
    \node[above=0pt of S_dir] {$\rr^2\times\rr^2$};

  \node[poly, cod, right=of S] (p) {$(a,b)\vphantom{\left(x+v_x,y+v_y,v_x+\frac{a}{(a^2+b^2)^{3/2}},v_y+\frac{b}{(a^2+b^2)^{3/2}}\right)}$\at$(x,y)\vphantom{((x,y),(v_x,v_y))}$};
    \node[below=0pt of p_pos] {$\rr^2$};
    \node[above=0pt of p_dir] {$\rr^2-\{(0,0)\}$};

  \draw (S_pos) -- node[below] {return} (p_pos);
  \draw (p_dir) -- node[above] {update} (S_dir);
\end{tikzpicture}
\]
\end{example}

\begin{exercise}\index{robot}
Suppose $(X,d)$ is a \emph{metric space}, i.e.\ $X$ is a set of \emph{points} and $d\colon X\times X\to\rr_{\geq0}$ is a function called the \emph{distance function} for which $d(x,y)=d(y,x),d(x,y)=0$ if and only if $x=y$, and $d(x,y)+d(y,z)\geq d(x,z)$ for all $x,y,z\in X$.
Let us have robots interact in this space.

Let $A,A'$ be sets, each thought of as a set of signals, and let $a_0\in A$ and $a_0'\in A'$ be elements, each thought of as a default value. Let $p\coloneqq AX\yon^{A'X}$ and $p'\coloneqq A'X\yon^{AX}$, and imagine there are two robots, one with interface $p$, returning a signal as an element of $A$ and its location as a point in $X$, and one with interface $p'$, returning a signal as an element of $A'$ and also its location as a point in $X$.
\begin{enumerate}
	\item Write down an interaction pattern $p\otimes p'\to\yon$ such that each robot receives the other's location but only receives the other's signal when their locations $x,x'$ are sufficiently close, namely when $d(x,x')<1$.
	Otherwise, it receives the default signal.
	\item Modify the previous interaction pattern to specify a new interaction pattern $p\otimes p'\to\yon^{[0,5]}$ where the value $s\in [0,5]$ is a scalar that changes the distance threshold for the signal to $s$.
	\item Suppose that each robot has a set $S,S'$ of possible private states in addition to their locations.
	What functions are involved in providing a dynamical system $\varphi\colon SX\yon^{SX}\to AX\yon^{A'X}$, if the location state $x\in X$ is directly returned without modification?
	\item Change the setup in any way so that each robot only extends a port to hear the other's signal when the distance between them is less than $s$. Otherwise, they can only detect the position (element of $X$) that the other currently inhabits.
	(Don't worry too much about timing---one missed signal when the robots first get close or one extra signal when the robots first get far is okay.)
\qedhere
\end{enumerate}
\begin{solution}
\begin{enumerate}
    \item An interaction pattern $p\otimes p'=AX\yon^{A'X}\otimes A'X\yon^{AX}\to\yon$ consists of a trivial on-positions function $AX\times A'X\to\1$ and an on-directions map $AX\times A'X\to A'X\times AX$ indicating what directions the robots should receive according to the positions they return.
    To model the robots receiving each others' locations but only receiving each others' signals when the distance between their locations is less than $1$, this on-directions function should send
    \[
        ((a,x),(a',x'))\mapsto
          \begin{cases}
          	((a',x'),(a,x))&\tn{ if }d(x,x')<1\\
          	((a'_0,x'),(a_0,x))&\tn{ otherwise}.
          \end{cases}
    \]
    \item An interaction pattern $p\otimes p'\to\yon^{[0,5]}$ that changes the distance threshold for the signal to $s\in[0,5]$ consists of a still trivial on-positions function and an on-directions map $AX\times A'X\times[0,5]\to A'X\times AX$ indicating what directions the robots should receive according to the external direction $s\in[0,5]$ and the positions they return.
    To model the fact that the robots receive each others' locations, but only receive each others' signals when the distance between their locations is less than $s$, this on-directions function should send
    \[
        ((a,x),(a',x'),s)\mapsto
          \begin{cases}
          	((a',x'),(a,x))&\tn{ if }d(x,x')<s\\
          	((a'_0,x'),(a_0,x))&\tn{ otherwise}.
          \end{cases}
    \]
    \item To provide a dynamical system $\varphi\colon SX\yon^{SX}\to AX\yon^{A'X}$ under the condition that the on-positions function preserves the second coordinate $x\in X$, we must provide the first projection $SX\to A$ of an on-positions function that turns the robot's private state and current location into the signal it returns, as well as an on-directions function $SX\times A'X\to SX$ that provides a new private state and location for the robot given its old private state, old location, and the signal and location it receives from the other robot.

    \item To have the robots listen for each others' signals only when they are sufficiently close, we must move away from monomial interfaces and Moore machines to leverage dependency.
    There are several ways of doing this; we give just one method below.
    With $D\coloneqq\{\text{`close'},\text{`far'}\}$, let the robots' new interfaces be
    \[
        p\coloneqq \{\text{`close'}\}AX\yon^{DA'X}+\{\text{`far'}\}AX\yon^{DX} \qqand p'\coloneqq \{\text{`close'}\}A'X\yon^{DAX}+\{\text{`far'}\}A'X\yon^{DX},
    \]
    so that they may receive input telling them whether they are close or far, but cannot receive signals in $A$ or $A'$ when they are `far.'

    Then by the distributivity of $\otimes$ over $+$, their new interaction pattern $p\otimes p'\to\yon^{[0,5]}$ can be specified by four lenses, all trivial on positions: the lens
    \[
        \{\text{`close'}\}AX\yon^{DA'X}\otimes\{\text{`close'}\}A'X\yon^{DAX}\to\yon^{[0,5]},
    \]
    given by the on-directions function
    \[
        ((\text{`close'},a,x),(\text{`close'},a',x'),s)\mapsto
          \begin{cases}
          	((\text{`close'},a',x'),(\text{`close'},a,x))&\tn{ if }d(x,x')<s\\
          	((\text{`far'},a'_0,x'),(\text{`far'},a_0,x))&\tn{ otherwise};
          \end{cases}
    \]
    the lens
    \[
        \{\text{`far'}\}AX\yon^X\otimes\{\text{`far'}\}A'X\yon^X\to\yon^{[0,5]},
    \]
    given by the on-directions function
    \[
        ((\text{`far'},a,x),(\text{`far'},a',x'),s)\mapsto
          \begin{cases}
          	((\text{`close'},x'),(\text{`close'},x))&\tn{ if }d(x,x')<s\\
          	((\text{`far'},x'),(\text{`far'},x))&\tn{ otherwise};
          \end{cases}
    \]
    and two other lenses that can be defined arbitrarily, as they should never come up in practice.

    Finally, in order for each robot to properly remember whether the other is close or far, we record an element of $D$ in its state that is returned and updated: one robot is a lens
    \[
        \varphi\colon DSX\yon^{DSX}\to \{\text{`close'}\}AX\yon^{DA'X}+\{\text{`far'}\}AX\yon^{DX}
    \]
    whose on-positions function preserves not just the third coordinate $x\in X$ but also the first coordinate $d\in D$, while the on-directions function also preserves the first coordinate $d\in D$; and the other robot is constructed similarly.
\end{enumerate}
\end{solution}
\end{exercise}

\index{interface!wrapper|)}

\subsection[Sectioning off juxtaposed systems]{Sectioning off juxtaposed dynamical systems}

We saw in \cref{subsec.poly.dyn_sys.new.sit_encl} that a section (i.e.\ lens to $\yon$) for the interface of a dynamical system sections that dynamical system off as a closed system.
So it should not come as a surprise that a section for a parallel product of interfaces yields an interaction pattern between the interfaces that only allows the interfaces to interact with each other, cutting off any other interaction with the outside world.

\index{parallel product!section for}

\index{interaction!picking up the chalk|(}\index{interface!closed}

\begin{example}[Picking up the chalk]\label{ex.pickup_chalk}
Imagine that you see some chalk and you pinch it between your thumb and forefinger.
An amazing thing about reality is that you will then have the chalk, in the sense that you can move it around.
How might we model this in $\poly$?
We will construct a closed dynamical system---one with interface $\yon$---consisting of only you and the chalk.
To do so, we will provide an interface for you, and interface for the chalk, and a section for your juxtaposition.

Say that your hand can be at one of two heights, $\const{down}$ or $\const{up}$, and that your fingers can either be $\const{pressed}$ (with pressure between your thumb and forefinger) or $\const{unpressed}$. Say too that you take in information about the chalk's height, which can be $\const{down}$ or $\const{up}$ as well. Here are the two sets we will be using:
\[
	H\coloneqq\{\const{down},\,\const{up}\}
	\qqand
	P\coloneqq\{\const{pressed},\,\const{unpressed}\}.
\]
Your interface is $HP\yon^H$: your position is your own height and pressure, and your possible directions are the chalk's possible heights.
As for the chalk, it is either $\const{in}$ your possession or $\const{out}$ of it as well as either $\const{down}$ or $\const{up}$.
The direction the chalk receives includes whether it is $\const{pressed}$ or $\const{unpressed}$.
When it's $\const{out}$ of your possession, that is the entire direction, but when it is $\const{in}$ your possession, its direction also comprises your hand's height.
In summary, here are the two interfaces:
\[
	\text{You}\coloneqq HP\yon^H
	\qqand
	\text{Chalk}\coloneqq \{\const{out}\}H\yon^P + \{\const{in}\}H\yon^{HP}.
\]

Now we want to give the interaction pattern between you and the chalk.
As we said before, you see the chalk's height.
If your hand is not at the height of the chalk, the chalk remains $\const{unpressed}$.
Otherwise, your hand is at the height of the chalk, so the chalk receives your pressure (or lack thereof).
Furthermore, if the chalk is in your possession, it also receives your hand's height.

To provide a lens $\gamma\colon\text{You}\otimes\text{Chalk}\to\yon$, we use the fact that $\text{Chalk}$ is a coproduct and that $\otimes$ distributes over coproduct to write $\text{You}\otimes\text{Chalk}$ as a coproduct itself:
\[
  \text{You}\otimes\text{Chalk}
    \iso
  HP\yon^H\otimes\{\const{out}\}H\yon^P+HP\yon^H\otimes\{\const{in}\}H\yon^{HP}.
\]
Then by the universal property of coproducts, to define $\gamma$, it suffices to define two lenses
\[
	\alpha\colon HP\yon^H\otimes\{\const{out}\}H\yon^P\to\yon
	\qqand
	\beta\colon HP\yon^H\otimes \{\const{in}\}H\yon^{HP}\to\yon
\]
The lens $\beta$, corresponding to when the chalk is $\const{in}$ your possession, is easy to describe: it can be identified with a function $HPH\to HHP$, and we take it to be the obvious function sending your height and pressure to the chalk and the chalk's height to you; see \cref{exc.pickup_chalk}.
But $\alpha$, corresponding to when the chalk is $\const{out}$ of your possession, is more semantically interesting: it can be identified with a function $HPH\to HP$ given by
\[
  (h_\text{You},p_\text{You},h_\text{Chalk})\mapsto
  \begin{cases}
  	(h_\text{Chalk},\const{unpressed}) & \tn{ if } h_\text{You} \neq h_\text{Chalk} \\
  	(h_\text{Chalk},p_\text{You}) & \tn{ if } h_\text{You} = h_\text{Chalk}.
  \end{cases}
\]
In words, this says that if you and the chalk are at different heights, then regardless of your pressure, the chalk remains unpressed; but if you are at the same height as the chalk, the chalk receives your pressure.

Now that you and the chalk are sectioned off together by $\gamma$, we are ready to add some dynamics.
Your dynamics can be whatever you want, so let us focus on giving dynamics to the chalk (you will be able to give yourself dynamics in \cref{exc.pickup_chalk}).
The chalk's state is comprised of its height and whether or not it is in your possession, so we give a dynamical system with state-set $C\coloneqq\{\const{out},\,\const{in}\} \times H$ and interface $\text{Chalk}$: that is, a lens
\begin{equation}\label{eqn.chalk_dynamics}
	\{\const{out},\,\const{in}\}\times H\yon^{\{\const{out},\,\const{in}\}\times H}\to\{\const{out}\}H\yon^P + \{\const{in}\}H\yon^{HP}.
\end{equation}
On positions, the chalk returns its height and whether it is in your possession directly: in other words, the on-positions function is the identity.
On directions, we have two cases.
If the chalk is $\const{out}$ of your possession, it falls $\const{down}$ unless you catch it, making it $\const{pressed}$ so that it becomes $\const{in}$ your possession and retains its current height.
So we can express the on-directions function of \eqref{eqn.chalk_dynamics} at $(\const{out}, h_\text{Chalk})$ as
\begin{align*}
	\const{unpressed} &\mapsto (\const{out},\,\const{down}) \\
	\const{pressed} &\mapsto (\text{`in'}, h_\text{Chalk})
\end{align*}
On the other hand, if the chalk is $\const{in}$ your possession, it takes whatever height it receives from you, remaining $\const{in}$ your possession if $\const{pressed}$ but coming $\const{out}$ of your possession if $\const{unpressed}$.
So the on-directions function of \eqref{eqn.chalk_dynamics} at $(\const{in}, h_\text{Chalk})$ is given by
\begin{align*}
	(h_\text{You}, \const{unpressed}) &\mapsto (\const{out}, h_\text{You}) \\
	(h_\text{You}, \const{pressed}) &\mapsto (\const{in}, h_\text{You}).
\end{align*}

We have thus defined an interaction pattern that allows one system to engage with or disengage from another system and control the behavior of the other system only when the two are engaged.
\end{example}\index{control}

\begin{exercise}\label{exc.pickup_chalk}
\begin{enumerate}
	\item In \cref{ex.pickup_chalk}, we said that $\beta\colon HP\yon^H\otimes \{\const{in}\}H\yon^{HP}\to\yon$ was easy to describe and given by a function $HPH\to HHP$. Explain what is being said, and provide the function.
	\item Provide dynamics to the $\text{You}$ interface (i.e.\ specify a dynamical system with interface $\text{You}=HP\yon^H$) so that you repeatedly reach down and grab the chalk, lift it with your hand, and drop it.
\qedhere
\end{enumerate}
\begin{solution}
\begin{enumerate}
    \item A lens $HP\yon^H\otimes \{\const{in}\}H\yon^{HP}\to\yon$ is a section of $HP\yon^H\otimes \{\const{in}\}H\yon^{HP}\iso HPH\yon^{HHP}$ and can thus be identified with a function from its position-set to its direction-set: $HPH\to HHP$.
    Here the first factor of the domain refers to your height, the second factor to the your pressure, and the third factor to the chalk's height; while in the codomain, the first factor refers to the chalk's height that you receive, the second factor refers to your height that the chalk receives, and the third factor refers to your pressure that the chalk receives.
    We should therefore define the function $HPH\to HHP$ to be the isomorphism that sends $(h_\text{You},p_\text{You},h_\text{Chalk})\mapsto(h_\text{Chalk},h_\text{You},p_\text{You})$.
    \item If you alwayss cycle through three possible actions---reaching $\const{down}$ and grabbing the chalk so that it is $\const{pressed}$, moving your hand $\const{up}$ while keeping it $\const{pressed}$ to lift the chalk, and dropping the chalk by leaving it $\const{unpressed}$ while your hand is still $\const{up}$---you only need $3$ possible states.
    So we can provide your dynamics using a lens $\3\yon^\3\to\text{You}=HP\yon^H$.
    The return function $\3\to HP$ indicates what happens at each state, sending $1\mapsto(\const{down},\,\const{pressed}),2\mapsto(\const{up},\,\const{pressed})$, and $3\mapsto(\const{up},\,\const{unpressed})$.
    Meanwhile the update function $\3H \to \3$ always changes the state to the next one, regardless of direction: it ignores the $H$ coordinate and sends $1\mapsto2$, $2\mapsto3$, and $3\mapsto1$.
\end{enumerate}
\end{solution}
\end{exercise}

\index{interaction!picking up the chalk|)}
\index{section}

Given $n\in\nn$ and polynomials $p_1,\ldots,p_n$ as interfaces, a section $p_1\otimes\cdots\otimes p_n\to\yon$ sections off these $n$ interfaces together.
The following proposition provides an alternative perspective on such sections.

\begin{proposition}\label{prop.situations2}
Given polynomials $p,q\in\poly$, there is a bijection
\begin{equation} \label{eqn.situations2}
\Gamma(p\otimes q)\cong\smset\big(q(\1),\Gamma(p)\big)\;\times\;\smset\big(p(\1),\Gamma(q)\big).
\end{equation}
\end{proposition}
The idea is that specifying a section for the interfaces $p$ and $q$ together is equivalent to specifying a section for $p$ for every output $q$ might return and specifying a section for $q$ for every output $p$ might return.
\begin{proof}[Proof of \cref{prop.situations2}]
This is a direct calculation:
\begin{align*}
	\Gamma(p\otimes q) &\iso
	\prod_{i\in p(\1)}\prod_{j\in q(\1)}(p[i]\times q[j]) \\
	&\iso
	\left(\prod_{j\in q(\1)}\prod_{i\in p(\1)}p[i]\right)\times
		 \left(\prod_{i\in p(\1)}\prod_{j\in q(\1)}q[j]\right) \\
	&\iso
	\smset(q(\1),\Gamma(p))\times\smset(p(\1),\Gamma(q)).
\end{align*}
\end{proof}

\begin{example}
A section $f\colon I\yon^A\otimes I'\yon^{A'}\to\yon$ corresponds to a function $I\times I'\to A\times A'$. In other words, for every pair of positions $(i,i')\in I\times I'$, the section $f$ specifies a pair of directions $(a,a')\in A\times A'$.

Let us think of the positions in $I$ and $I'$ as locations that two machines may occupy.
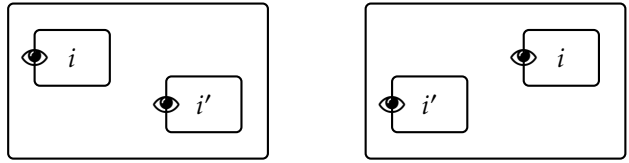
\[
\begin{tikzpicture}[oriented WD, bb port length=0]
	\node[bb={1}{0}] (p) {$i$};
	\node[bb={1}{0}, below right=-0.5 and 0.5 of p] (q) {$i'$};
	\node[bb={0}{0}, inner sep=10pt, fit=(p) (q)] {};
	\node at (p_in1) {\faEye};
	\node at (q_in1) {\faEye};
\end{tikzpicture}
\hspace{.5in}
\begin{tikzpicture}[oriented WD, bb port length=0]
	\node[bb={1}{0}] (p) {$i$};
	\node[bb={1}{0}, below left=-0.5 and 0.5 of p] (q) {$i'$};
	\node[bb={0}{0}, inner sep=10pt, fit=(p) (q)] {};
	\node at (p_in1) {\faEye};
	\node at (q_in1) {\faEye};
\end{tikzpicture}
\]
Then given a pair of locations $(i,i')\in I\times I'$, the interaction pattern $f$ as a function $I\times I'\to A\times A'$ tells us the directions the machines observe, i.e.\ the ordered pair $(a,a')\in A\times A'$ comprised of what each machine receives.
Equivalently, \eqref{eqn.situations2} says that the interaction pattern tells us what direction the first machine observes at each location when the second machine's location is fixed at $i'$, along with the direction the second machine observes at each location when the first machine's location is fixed at $i$.

Here we see that \eqref{eqn.situations2} provides two ways to interpret the interaction pattern between two interfaces in a closed system: either as a section around each interface parametrized by the other's position, or as a single section around them both.
\end{example}

\begin{exercise}
Let $p\coloneqq q\coloneqq\nn\yon^\nn$.
We wish to specify a section around their juxtaposition.
\begin{enumerate}
  \item Say we wanted to feed the position of $q$ as a direction for $p$.
  What function $f\colon q(\1)\to\Gamma(p)$ captures this behavior?
  \item Say we wanted to feed the sum of the positions of $p$ and $q$ as a direction for $q$.
  What function $g\colon p(\1)\to\Gamma(q)$ captures this behavior?
  \item What section $\gamma\colon p\otimes q\to\yon$ does the pair of functions $(f,g)$ correspond to via \eqref{eqn.situations2}?
	\item Let dynamical systems $\varphi\colon\nn\yon^\nn\to p$ and $\psi\colon\nn\yon^\nn\to q$ both be the identity on $\nn\yon^\nn$.
	Suppose $\varphi$ starts in the state $0\in\nn$ and $\psi$ starts in the state $1\in\nn$.
	Describe the behavior of the system obtained by sectioning $\varphi$ and $\psi$ off together with $\gamma$, i.e.\ the system $(\varphi\otimes\psi)\then\gamma$.
\qedhere
\end{enumerate}
\begin{solution}
\begin{enumerate}
    \item The function $f\colon q(\1)\to\Gamma(p)$ should send each $q$-position $b$ to the section of $p$ corresponding to the constant function $\nn\to\nn$ that sends every $p$-position to $b$ itself.
    That is, $f$ is the function $b\mapsto(\_\mapsto b)$.
    \item The function $g\colon p(\1)\to\Gamma(q)$ should send each position $p$-position $a$ to the section of $q$ corresponding to the function $\nn\to\nn$ that sends each $q$-position $b$ to the sum $a+b$.
    That is, $g$ is the function $a\mapsto(b\mapsto a+b)$.
    \item Together, $f$ and $g$ form a function $\nn\times\nn\to\nn\times\nn$ mapping $(a,b)\mapsto((fb)a,(ga)b)=(b,a+b)$.
    Then $\gamma\colon p\otimes q\to\yon$ is the section with this function as its on-directions function.
    \item As $\varphi$ and $\psi$ are both the identity, their parallel product is the identity as well, so $(\varphi\otimes\psi)\then\gamma=\gamma$.
    From the previous part, if the current state of this system is $(a,b)$, its next state will be $(b,a+b)$.
    So if its initial state is $(0,1)$, its following states will be $(1,1),(1,2),(2,3),(3,5),(5,8),(8,13),\ldots$, forming the familiar Fibonacci sequence.
\end{enumerate}
\end{solution}
\end{exercise}

\index{interaction!picking up the chalk}

\begin{exercise}
We will use \eqref{eqn.situations2} to consider the interaction pattern $\gamma$ between \text{You} and \text{Chalk} from \cref{ex.pickup_chalk} as a pair of functions $\oper{You}(\1)\to\Gamma(\text{Chalk})$ and $\oper{Chalk}(\1)\to\Gamma(\text{You})$.
\begin{enumerate}
	\item How does the chalk's position specify a section for you? That is, describe the function $\oper{Chalk}(\1)\to\Gamma(\text{You})$.
	\item How does your position specify a section for the chalk? That is, describe the function $\oper{You}(\1)\to\Gamma(\text{Chalk})$.
\qedhere
\end{enumerate}
\begin{solution}
\begin{enumerate}
    \item Fix a position $(s_\text{Chalk}, h_\text{chalk})\in\oper{Chalk}(\1)=\{\const{out},\const{in}\}H$ of the chalk.
    If $s_\text{Chalk}=\const{out}$, then the corresponding section $\text{You}\iso HP\yon^H\to\yon$ given by $f$ via \eqref{eqn.situations2} can be thought of as the function $HP\to H$ sending
    \[
        (h_\text{You},p_\text{You})\mapsto h_\text{Chalk},
    \]
    according to the behavior of $\alpha\colon HP\yon^H\otimes H\yon^P\to\yon$ when we fix the position of the domain's right factor to be $h_\text{chalk}\in H$ and focus on the direction in $H$ of the domain's left factor.
    Meanwhile, if $s_\text{Chalk}=\const{in}$, then the corresponding section $\text{You}\iso HP\yon^H\to\yon$ can also be thought of as the function $HP\to H$ sending
    \[
        (h_\text{You},p_\text{You})\mapsto h_\text{Chalk},
    \]
    according to the behavior of $\beta\colon HP\yon^H\otimes H\yon^{HP}\to\yon$ when we again fix the position of the domain's right factor to be $h_\text{chalk}\in H$ and focus on the direction in $H$ of the domain's left factor.
    So overall, the desired function $\oper{Chalk}(\1)\to\Gamma(\text{You})$ is given by
    \[
        (\_,h_\text{chalk})\mapsto((\_,\_)\mapsto h_\text{Chalk}).
    \]

    \item Fix a position $(h_\text{You},p_\text{You})\in\oper{You}(\1)=HP$ of the chalk.
    Then the corresponding section $\text{Chalk}\iso\{\const{out}\}H\yon^P + \{\const{in}\}H\yon^{HP}\to\yon$ can be thought of as a pair of functions: one $\{\const{out}\}H\to P$ sending
    \[
        (\const{out},h_\text{Chalk})\mapsto
            \begin{cases}
  	            \const{unpressed} & \tn{ if } h_\text{You} \neq h_\text{Chalk} \\
  	            p_\text{You} & \tn{ if } h_\text{You} = h_\text{Chalk},
            \end{cases}
    \]
    according to the behavior of $\alpha\colon HP\yon^H\otimes H\yon^P\to\yon$ when we fix the position of the domain's left factor to be $(h_\text{You},p_\text{You})\in HP$ and focus on the direction in $P$ of the domain's right factor; and another $\{\const{in}\}H\to HP$ sending
    \[
        (\const{in},h_\text{Chalk})\mapsto(h_\text{You},p_\text{You})
    \]
    according to the behavior of $\beta\colon HP\yon^H\otimes H\yon^{HP}\to\yon$ when we again fix the position of the domain's left factor to be $(h_\text{You},p_\text{You})\in HP$ and focus on the direction in $HP$ of the domain's right factor.
    So overall, the desired function $\oper{You}(\1)\to\Gamma(\text{Chalk})$ is given by
    \[
        (h_\text{You},p_\text{You})\mapsto\left(
            \begin{aligned}
                (\const{out},h_\text{Chalk})&\mapsto
                    \begin{cases}
          	            \const{unpressed} & \tn{ if } h_\text{You} \neq h_\text{Chalk} \\
          	            p_\text{You} & \tn{ if } h_\text{You} = h_\text{Chalk}
                    \end{cases} \\
                (\const{in},h_\text{Chalk})&\mapsto(h_\text{You},p_\text{You})
            \end{aligned}
        \right).
    \]
\end{enumerate}
\end{solution}
\end{exercise}

\begin{exercise}
\begin{enumerate}
	\item State and prove a generalization of \eqref{eqn.situations2} from \cref{prop.situations2} for $n$-many polynomials $p_1,\ldots,p_n\in\poly$.
	\item Generalize the ``idea'' statement between \cref{prop.situations2} and its proof.
\qedhere
\end{enumerate}
\begin{solution}
\begin{enumerate}
    \item We generalize \eqref{eqn.situations2} for $n$ polynomials as follows.
    Given polynomials $p_1,\ldots,p_n\in\poly$, we claim there is a bijection
    \[
        \Gamma\left(\bigotimes_{i=1}^n p_i \right) \iso \prod_{i=1}^n \smset\left(\prod_{\substack{1 \leq j \leq n, \\ j \neq i}} p_j(\1), \Gamma(p_i)\right).
    \]
    The $n=1$ case is tautological, and the $n=2$ case is given by \eqref{eqn.situations2}.
    Then by induction on $n$, we have
    \begin{align*}
        \Gamma\left(\bigotimes_{i=1}^n p_i \right) &\iso \smset\left(p_n(\1), \Gamma\left(\bigotimes_{i=1}^{n-1} p_i \right)\right) \times \smset\left(\prod_{i=1}^{n-1}\left(p_i(\1)\right), \Gamma(p_n)\right) \tag*{\eqref{eqn.situations2}} \\
        &\iso \smset\left(p_n(\1), \prod_{i=1}^{n-1} \smset\left(\prod_{\substack{1 \leq j \leq n-1, \\ j \neq i}} p_j(\1), \Gamma(p_i)\right)\right) \times \smset\left(\prod_{i=1}^{n-1} p_i(\1), \Gamma(p_n)\right) \tag{Inductive hypothesis} \\
        &\iso \prod_{i=1}^{n-1}\left( \smset\left(\prod_{\substack{1 \leq j \leq n, \\ j \neq i}} p_j(\1), \Gamma(p_i)\right)\right) \times \smset\left(\prod_{i=1}^{n-1} p_i(\1), \Gamma(p_n)\right) \tag{Universal properties of products and internal homs},
    \end{align*}
    and the result follows.
    \item The general idea is that specifying a section for interfaces $p_1,\ldots,p_n$ together is equivalent to specifying a section for $p_i$ for every combination of positions that all the other interfaces might return together, for each $i\in\ord{n}$.
\end{enumerate}
\end{solution}
\end{exercise}

\index{section}

\subsection{Wiring diagrams as interaction patterns}

\index{wiring diagram|see{interaction, wiring diagram}}
\index{interaction!wiring diagram|(}

A \emph{wiring diagram} is a graphical depiction of interactions between systems.
Wiring diagrams depict systems as boxes, showing how they send signals to each other through the wires between them, as well as how multiple systems can combine to form a larger system whenever smaller boxes are nested within a larger box.

\index{interface!monomial}

Formally, and more precisely, we can think of each box in a wiring diagram as an interface given by some monomial.
The box itself is not a dynamical system, but it becomes a dynamical system once we equip it with a lens from a state system to the interface the box represents.
Then the entire wiring diagram---specifying how these boxes nest within a larger box---is just an interaction pattern between the corresponding interfaces, with the larger box playing the role of the wrapper interface.
Once every nested box is equipped with a lens from a state system, we obtain a dynamical system whose interface is the larger box.

In the examples to come, we follow the convention that the signals emitted by a box, i.e.\ positions returned by the corresponding interface, travel along wires out of the right side of that box; while the signals received, i.e.\ directions observed by the corresponding interface, by a box travel along wires into the left side of that box.
A wire may optionally be labeled by the name of the set of elements that may travel as signals along that wire.

\begin{example}\index{control}
Here is a simple wiring diagram.
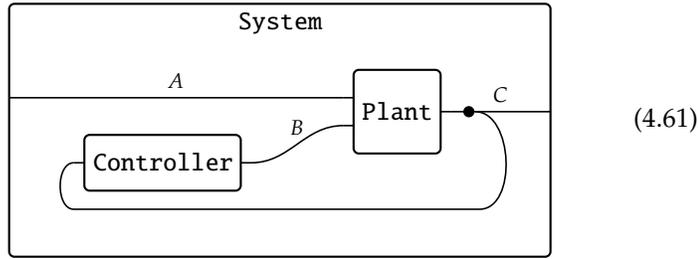
\begin{equation}\label{eqn.control_diag}
\begin{tikzpicture}[oriented WD, baseline=(B)]
	\node[bb={2}{1}] (plant) {\texttt{Plant}};
	\node[bb={1}{1}, below left=-1 and 1 of plant]  (cont) {\texttt{Controller}};
	\node[circle, inner sep=1.5pt, fill=black, right=.1] at (plant_out1) (pdot) {};
	\node[bb={0}{0}, inner ysep=25pt, inner xsep=1cm, fit=(plant) (pdot) (cont)] (outer) {};
	\coordinate (outer_out1) at (outer.east|-plant_out1);
	\coordinate (outer_in1) at (outer.west|-plant_in1);
	\begin{scope}[above, font=\footnotesize]
  	\draw (outer_in1) -- node {$A$} (plant_in1);
  	\draw (cont_out1) to node (B) {$B$} (plant_in2);
  	\draw (plant_out1) to node {$C$} (outer_out1);
  	\draw
  		let
  			\p1 = (cont.south west-| pdot),
  			\p2 = (cont.south west),
  			\n1 = \bby,
  			\n2 = \bbportlen
  		in
  			(pdot) to[out=0, in=0]
  			(\x1+\n2, \y1-\n1) --
  			(\x2-\n2, \y2-\n1) to[out=180, in=180]
  			(cont_in1);
		\end{scope}
	\node[below=0of outer.north] {\texttt{System}};
\end{tikzpicture}
\end{equation}
The $\const{Plant}$ is receiving information from the world outside the $\const{System}$ along the wire labeled $A$ as well as from the $\const{Controller}$ along the wire labeled $B$.
It is also producing information for the outside world along the wire labeled $C$ which is also being monitored by the $\const{Controller}$.\index{control}

There are three boxes shown in \eqref{eqn.control_diag}: the $\const{Controller}$, the $\const{Plant}$, and the $\const{System}$.
Each has a fixed set of positions corresponding to the wire(s) connected to its right and a fixed set of directions corresponding to the wire(s) connected to its left, so we can consider each box as a monomial interface, as follows:\index{interface!monomial}
\begin{equation}\label{eqn.basic_diagram}
	\const{Controller}\coloneqq B\yon^C
	\qquad\quad
  \const{Plant}\coloneqq C\yon^{AB}
	\qquad\quad
	\const{System}\coloneqq C\yon^A.
\end{equation}
Note that in the case of the $\const{Plant}$, two wires labeled $A$ and $B$ enter the box from the left, so we take their cartesian product to be the direction-set of the $\const{Plant}$.

The wiring diagram itself is a wrapper
\[
	w\colon\const{Controller}\otimes\const{Plant}\to\const{System},
\]
specifying an interaction pattern between the $\const{Controller}$ and the $\const{Plant}$ with the $\const{System}$ as the wrapper interface.
Concretely, $w$ is a lens $BC\yon^{CAB}\to C\yon^A$ that prescribes how wires feed positions to directions.
As a lens between monomials, $w$ consists of an on-positions function $BC\to C$ and an on-directions map $BCA\to CAB$.

The wiring diagram is a picture that tells us what the on-positions function and on-directions map to use.
In particular, the on-positions function sends positions of the inner interfaces to positions of the outer interface, so it is depicted by how the wires coming from the right sides of the inner boxes connect to the right side of the outer box.
Given inner boxes that return positions in $B$ and $C$, the outer box must return a position in $C$.
Here the wire labeled $B$ is not connected to the outer box, but the wire labeled $C$ is, so the on-positions function $BC\to C$ sends $(b,c)\mapsto c$.

Meanwhile, the on-directions function sends positions of the inner interfaces and directions of the outer interface to directions of the inner interfaces, so it is depicted by how the wires coming from both the right sides of the inner boxes and the left side of the outer box connect to the left sides of the inner boxes.
Given inner boxes that return positions in $B$ and $C$ and an outer box that receives directions in $A$, the inner boxes must receive directions in $C$ and $B$.
Again, we can read the on-directions map $BCA\to CAB$ off the wiring diagram: it sends $(b,c,a)\mapsto (c,a,b)$.

Note that neither the wiring diagram nor any of the boxes within it represent dynamical systems on their own.
Rather, each box is a monomial that could be the interface of a dynamical system.
When we assign to a box a dynamical system having that box as its interface, we say that we \emph{give dynamics} to the box.
So the entire wiring diagram is a wrapper that tells us how, once we give dynamics for each inner box,
\[
\varphi\colon S\yon^S\to\const{Controller}
\qqand
\psi\colon T\yon^T\to\const{Plant},
\]
we have given dynamics for the entire outer box:
\[
ST\yon^{ST}\To{\varphi\:\otimes\:\psi}\const{Controller}\otimes\const{Plant}\To{w}\const{System}.
\]
\end{example}\index{dynamics}

\begin{exercise}
\begin{enumerate}\index{control}
	\item Make a new wiring diagram like \eqref{eqn.control_diag} except where the controller also receives information from the outside world as an element of a set $A'$.
	\item What are the monomials represented by the boxes in your diagram (replacing \eqref{eqn.basic_diagram})?
	\item What is the interaction pattern represented by this wiring diagram?
	Give the corresponding lens, including its on-positions and on-directions functions.
\qedhere
\end{enumerate}
\begin{solution}
\begin{enumerate}
  \item Here is the wiring diagram \eqref{eqn.control_diag} modified so that the controller also receives information from the outside world as an element of $A'$.
  \[
  \begin{tikzpicture}[oriented WD, baseline=(B)]
  	\node[bb={2}{1}] (plant) {\texttt{Plant}};
  	\node[bb={2}{1}, below left=-1 and 1 of plant]  (cont) {\texttt{Controller}};
  	\node[circle, inner sep=1.5pt, fill=black, right=.1] at (plant_out1) (pdot) {};
  	\node[bb={0}{0}, inner ysep=25pt, inner xsep=1cm, fit=(plant) (pdot) (cont)] (outer) {};
  	\coordinate (outer_out1) at (outer.east|-plant_out1);
  	\coordinate (outer_in1) at (outer.west|-plant_in1);
  	\coordinate (outer_in2) at (outer.west|-cont_in1);
  	\begin{scope}[above, font=\footnotesize]
    	\draw (outer_in1) -- node {$A$} (plant_in1);
    	\draw (outer_in2) -- node {$A'$} (cont_in1);
    	\draw (cont_out1) to node (B) {$B$} (plant_in2);
    	\draw (plant_out1) to node {$C$} (outer_out1);
    	\draw
    		let
    			\p1 = (cont.south west-| pdot),
    			\p2 = (cont.south west),
    			\n1 = \bby,
    			\n2 = \bbportlen
    		in
    			(pdot) to[out=0, in=0]
    			(\x1+\n2, \y1-\n1) --
    			(\x2-\n2, \y2-\n1) to[out=180, in=180]
    			(cont_in2);
  		\end{scope}
  	\node[below=0of outer.north] {\texttt{System}};
  \end{tikzpicture}
  \]
  \item The monomials represented by the boxes in this diagram are the same, except that the $\const{Controller}$ and the $\const{System}$ both have extra $A'$ factors in their exponent:
  \[
	\const{Controller}\coloneqq B\yon^{A'C}
  	\qquad\quad
  \const{Plant}\coloneqq C\yon^{AB}
  	\qquad\quad
	\const{System}\coloneqq C\yon^{AA'}.
  \]
  \item The interaction pattern represented by this wiring diagram is the lens
  \[
   	w'\colon \const{Controller}\otimes\const{Plant}\to\const{System}
  \]
  consisting of an on-positions function $BC\to C$ given by $(b,c)\mapsto c$ and an on-directions function $BCAA'\to A'CAB$ given by $(b,c,a,a')\mapsto(a',c,a,b)$.
\end{enumerate}
\end{solution}
\end{exercise}

\begin{exercise}
Consider the following wiring diagram.
\[
\begin{tikzpicture}[oriented WD, font=\footnotesize, bb port sep=1, bb port length=2.5pt, bb min width=.4cm, bby=.2cm, inner xsep=.2cm, x=.5cm, y=.3cm, text height=1.5ex, text depth=.5ex]
  	\node[bb={2}{1}] (Trf) {$\const{Alice}$};
  	\node[bb={1}{2}, below=1 of Trf] (Trg) {$\const{Bob}$};
		\node[bb={2}{2}] at ($(Trf)!.5!(Trg)+(1.5,0)$) (Trh) {$\const{Carl}$};
  	\node[bb={0}{0}, fit={($(Trf.north west)+(-.25,4)$) (Trg) ($(Trh.north east)+(.25,0)$)}] (Tr) {};
		\node[below] at (Tr.north) {$\const{Team}$};
  	\node[coordinate] at (Tr.west|-Trf_in2) (Tr_in1) {};
  	\node[coordinate] at (Tr.west|-Trg_in1) (Tr_in2) {};
  	\node[coordinate] at (Tr.east|-Trh_out2) (Tr_out1) {};
  	\node at ($(Trg_out2)+(5pt,0)$) (dot) {$\bullet$};
\begin{scope}[font=\tiny]
  	\draw[shorten <=-2pt] (Tr_in1) -- node[below=-3pt] {$A$} (Trf_in2);
  	\draw[shorten <=-2pt] (Tr_in2) -- node[below=-3pt] {$B$} (Trg_in1);
		\draw (Trf_out1) to node[above=-3pt] {$D$} (Trh_in1);
		\draw (Trg_out1) to node[above=-3pt] {$E$} (Trh_in2);
  	\draw (Trg_out2) -- node[below=-3pt] {$F$} (dot.center);
  	\draw[shorten >=-2pt] (Trh_out2) -- node[below=-3pt] {$G$} (Tr_out1);
  	\draw let \p1=(Trh.east), \p2=(Trf.north west), \n1=\bbportlen, \n2=\bby in
  		(Trh_out1) to[in=0] (\x1+\n1,\y2+\n2) -- node[pos=.3, below=-3pt] {$H$} (\x2-\n1,\y2+\n2) to[out=180] (Trf_in1);
	\end{scope}
\end{tikzpicture}
\]
\begin{enumerate}
	\item Write out the monomials for $\const{Alice}$, $\const{Bob}$, and $\const{Carl}$.
	\item Write out the monomial for the outer box, $\const{Team}$.
	\item The wiring diagram constitutes a lens $f$ in $\poly$; what is its domain and codomain?
	\item What lens is $f$?
	\item Say we have dynamical systems $\alpha\colon A\yon^A\to\const{Alice}$, $\beta\colon B\yon^B\to\const{Bob}$, and $\gamma\colon C\yon^C\to\const{Carl}$. What is the induced dynamical system with interface $\const{Team}$?
\qedhere
\end{enumerate}
\begin{solution}
\begin{enumerate}
    \item According to the wiring diagram, we have that $\const{Alice}\coloneqq D\yon^{HA},$ that $\const{Bob}\coloneqq EF\yon^B,$ and that $\const{Carl}\coloneqq HG\yon^{DE}.$
    \item According to the wiring diagram, we have that $\const{Team}\coloneqq G\yon^{AB}.$
    \item The wiring diagram constitutes a wrapper
    \[
        f\colon\const{Alice}\otimes\const{Bob}\otimes\const{Carl}\to\const{Team}.
    \]
    Its domain is $\const{Alice}\otimes\const{Bob}\otimes\const{Carl}\iso DEFHG\yon^{HABDE}$, while its codomain is $\const{Team}=G\yon^{AB}$.
    \item On positions, the lens $f$ is a function $DEFHG\to G$ that sends $(d,e,f,h,g)\mapsto g$.
    On directions, $f$ is a function $DEFHGAB\to HABDE$ that sends $(d,e,f,h,g,a,b)\mapsto(h,a,b,d,e)$.
    \item Given dynamical systems $\alpha\colon A\yon^A\to\const{Alice}$, $\beta\colon B\yon^B\to\const{Bob}$, and $\gamma\colon C\yon^C\to\const{Carl}$, the dynamical system induced by the wiring diagram is given by the composite lens
    \[
        ABC\yon^{ABC}\To{\alpha\:\otimes\:\beta\:\otimes\:\gamma}\const{Alice}\otimes\const{Bob}\otimes\const{Carl}\To{f}\const{Team}.
    \]
\end{enumerate}
\end{solution}
\end{exercise}

\begin{exercise}[Long division] \label{exc.long_div}\index{long division}
\begin{enumerate}
	\item Let $\fun{divmod}\colon\nn\times\nn_{\geq1}\to\nn\times\nn$ send $(a,b)\mapsto(a\bdiv b, a\bmod b)$; for example, it sends $(10,7)\mapsto(1,3)$ and $(30,7)\mapsto(4,2)$.
	Use \cref{ex.funs_to_moore} to turn it into a dynamical system.
	\item
  In the following wiring diagram, we have already given dynamics to each box, as follows.
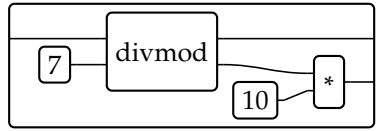
\[
\begin{tikzpicture}[oriented WD, bb small]
	\node[bb port sep=3, bb={2}{2}] (divmod) {divmod};
	\node[bb={0}{1}, left=of divmod_in2] (7) {$7$};
	\node[bb port sep=2, bb={2}{1}, below right=-1 and 3 of divmod_out2] (times) {$*$};
	\node[bb={0}{1}, below left=-1 and 1 of times_in2] (10) {$10$};
	\node[bb={0}{0}, inner xsep=\bbx, fit=(divmod) (times)(7) (10)] (outer) {};
	\coordinate (outer_in1) at (outer.west|-divmod_in1);
	\coordinate (outer_out1) at (outer.east|-divmod_out1);
	\coordinate (outer_out2) at (outer.east|-times_out1);
	\draw (outer_in1) -- (divmod_in1);
	\draw (7_out1) -- (divmod_in2);
	\draw (10_out1) -- (times_in2);
	\draw (divmod_out1) -- (outer_out1);
	\draw (divmod_out2) to (times_in1);
	\draw (times_out1) -- (outer_out2);
\end{tikzpicture}
\]
  The dynamical system corresponding to the box $\fun{divmod}$, with the box as its interface, is the dynamical system from the previous part (the upper wires correspond to the left hand factors of the domain and codomain of the $\fun{divmod}$ function, while the lower wires correspond to the right hand factors).
  Similarly, the box labeled $\ast$ corresponds to the  to the dynamical system arising from the multiplication function $\nn\times\nn\to\nn$ sending $(m,n)\mapsto mn$ following \cref{ex.funs_to_moore}.
  Meanwhile the boxes labeled $7$ and $10$ correspond to dynamical systems with $1$ state each that always return $7\in\nn_{\geq1}$ and $10\in\nn$, respecively.
  Describe the behavior of the dynamical system corresponding to the entire outer box.
	\item Using the outer box from the wiring diagram above as the inner box of the wiring diagram below, pick an initial state so that the resulting dynamical system alternates between returning $0$'s and the base-$10$ digits of $1/7$ after the decimal point, like so:
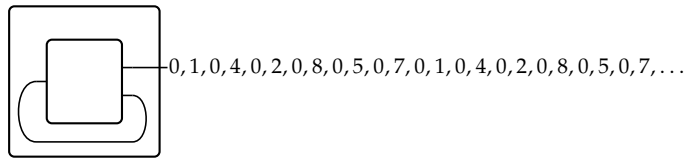
\[
\begin{tikzpicture}[oriented WD]
	\node[bb={1}{2}] (inner) {};
	\node[bb={0}{0}, inner xsep=1cm, inner ysep=1cm] (outer) {};
	\coordinate (outer_out1) at (outer.east|-inner_out1);
	\draw[shorten >=-3pt] (inner_out1) -- (outer_out1);
	\draw
		let \p1=(inner.south east), \p2=(inner.south west), \n1=\bbportlen, \n2=\bby in
		(inner_out2) to[in=0] (\x1+\n1,\y1-\n2) -- (\x2-\n1,\y1-\n2) to[out=180] (inner_in1);
		\node[right, font=\footnotesize] at (outer_out1) {$0,1,0,4,0,2,0,8,0,5,0,7,0,1,0,4,0,2,0,8,0,5,0,7,\ldots$};
\end{tikzpicture}
\]
We will see in \cref{subsec.comon.sharp.state.run} how to make a dynamical system run twice as fast, then apply this to the above system in \cref{ex.long_div_skip} so that it skips the $0$'s.
\qedhere
\end{enumerate}
\begin{solution}
\begin{enumerate}
    \item Using \cref{ex.funs_to_moore}, we can turn $\fun{divmod}$ into the dynamical system $\fun{divmod}\colon\nn\times\nn\yon^{\nn\times\nn}\to\nn\times\nn\yon^{\nn\times\nn_{\geq1}}$ whose return function is the identity on $\nn\times\nn$ and whose update map  $\nn\times\nn\times\nn\times\nn_{\geq1}\to\nn\times\nn$ sends $(\_,\_,a,b)\mapsto(a\bdiv b,a\bmod b)$.
    \item From left to right, the inner boxes represent monomial interfaces $\nn_{\geq1}\yon, \nn\times\nn\yon^{\nn\times\nn_{\geq1}}, \nn\yon,$ and $\nn\yon^{\nn\times\nn}$.\index{interface!monomial}
    The box labeled $7$ is given dynamics $7\colon\yon\to\nn_{\geq1}\yon$ so that it always returns the position $7$; similarly, the box labeled $10$ is given dynamics $10\colon\yon\to\nn\yon$ so that it always returns the position $10$.
    Meanwhile, the box labeled $\fun{divmod}$ is given dynamics $\fun{divmod}\colon\nn\times\nn\yon^{\nn\times\nn}\to\nn\times\nn\yon^{\nn\times\nn_{\geq1}}$ from the previous part; and applying \cref{exc.funs_to_moore} to the multiplication function $\ast\colon\nn\times\nn\to\nn$ yields the dynamics for the box labeled $\ast$: a dynamical system $\ast\colon\nn\yon^\nn\to\nn\yon^{\nn\times\nn}$ whose return function is the identity on $\nn$ and whose update map $\nn\times\nn\times\nn\to\nn$ sends $(\_,m,n)\mapsto mn$.

    Then the outer box is the monomial interface $\nn\times\nn\yon^\nn$, and the wiring diagram is the interaction pattern
    \[
        w\colon\nn_{\geq1}\yon\otimes\left(\nn\times\nn\yon^{\nn\times\nn_{\geq1}}\right)\otimes\nn\yon\otimes\nn\yon^{\nn\times\nn}\to\nn\times\nn\yon^\nn
    \]
    with on-positions function $(s,q,r,t,p)\mapsto(q,p)$ and on-directions map $(s,q,r,t,p,a)\mapsto(a,s,r,t)$.
    So the dynamical system induced by the wiring diagram is the composite lens $\varphi$ given by
    \[
        \yon\otimes\left(\nn\times\nn\yon^{\nn\times\nn}\right)\otimes\yon\otimes\nn\yon^\nn\To{7\:\otimes\:\fun{divmod}\:\otimes\:10\:\otimes\:\ast}\nn_{\geq1}\yon\otimes\left(\nn\times\nn\yon^{\nn\times\nn_{\geq1}}\right)\otimes\nn\yon\otimes\nn\yon^{\nn\times\nn}\To{w}\nn\times\nn\yon^\nn,
    \]
    whose return function is given by the composite map $(q,r,p)\mapsto(7,q,r,10,p)\mapsto(q,p)$ and whose update function at state $(q,r,p)$ is given by the composite map $a\mapsto(a,7,r,10)\mapsto(a\bdiv7,a\bmod7,10r)$.

    In other words, the dynamical system $\varphi$ behaves as follows: its state consists of a quotient $q$, a remainder $r$, and a product $p$, of which it returns the quotient and the product.
    Then it is fed a dividend $a$ and evaluates $a\bdiv7$ to obtain the new quotient and $a\bmod7$ to obtain the new remainder.
    Meanwhile, the new product is given by the previous remainder multiplied by $10$.

    \item This second wiring diagram specifies an interaction pattern
    \[
        w'\colon\nn\times\nn\yon^\nn\to\nn\yon
    \]
    with on-positions function $(q,p)\mapsto q$ and on-directions function $(q,p)\mapsto p$.
    So the dynamical system induced by nesting the first wiring diagram within the inner box of the second wiring diagram is the composite lens
    \[
        \left(\nn\times\nn\yon^{\nn\times\nn}\right)\otimes\nn\yon^\nn\To{\varphi}\nn\times\nn\yon^\nn\To{w'}\nn\yon
    \]
    whose return function is given by the composite map $(q,r,p)\mapsto(q,p)\mapsto q$ and whose update function at state $(q,r,p)$ specifies the new state $(p\bdiv7,p\bmod7,10r)$.

    In other words, the dynamical system $\varphi$ behaves as follows: its state consists of a quotient $q$, a remainder $r$, and a product $p$, of which it returns just the quotient.
    Then it advances to a new state by evaluating $p\bdiv7$ to obtain the new quotient and $p\bmod7$ to obtain the new remainder.
    Meanwhile, the new product is given by the previous remainder multiplied by $10$.

    If we set the initial state to be $(q,r,p)\coloneqq(0,0,10)$, then the subsequent states will be as follows, with the values of $q$ in the left column giving us the positions returned:
    \begin{table}[hbt!]
        \centering
        \footnotesize
        \begin{tabular}{c|c|c}
            $q$ $(p\bdiv7)$ & $r$ $(p\bmod7)$ & $p$ $(10r)$ \\
            \hline
            0 & 0 & 10 \\
            1 & 3 & 0 \\
            0 & 0 & 30 \\
            4 & 2 & 0 \\
            0 & 0 & 20 \\
            2 & 6 & 0 \\
            0 & 0 & 60 \\
            8 & 4 & 0 \\
            0 & 0 & 40 \\
            5 & 5 & 0 \\
            0 & 0 & 50 \\
            7 & 1 & 0 \\
            0 & 0 & 10 \\
            $\vdots$ & $\vdots$ & $\vdots$
        \end{tabular}
    \end{table}
\end{enumerate}
\end{solution}
\end{exercise}

\begin{example}[Graphs as wiring diagrams and cellular automata]\label{ex.graph_interaction}\index{cellular automata}\index{graph} % TODO: check if this still works with multiple parallel edges
Suppose we have a graph $G=(E\tto V)$ as in \cref{def.graph} and a set $\tau(v)$ associated with each vertex $v\in V$:
\[
\begin{tikzcd}
	E\ar[r, shift left=3pt, "s"]\ar[r, shift right=3pt, "t"']&
	V\ar[r, "\tau"]&
	\smset
\end{tikzcd}
\]
We can think of $G$ as an alternative representation of a specific kind of wiring diagram, one in which each inner box has exactly one position wire coming out and the outer box is closed (i.e.\ no wires in or out, representing the interface $\yon$ with exactly $1$ position and $1$ direction).
The vertices $v\in V$ are the inner boxes, the set $\tau(v)$ is the position-set associated with the wire coming out of $v$, and each edge $e$ is a wire connecting the position wire of its target $t(e)$ to a direction wire of its source $s(e)$.
An edge from a vertex $v_0$ to a vertex $v_1$ indicates that the directions received by $v_0$ depend on the positions returned by $v_1$.\tablefootnote{We could swap the roles of the sources and the targets, so that edges point in the direction of data flow rather than in the direction of data dependencies; this is an arbitrary choice.}

In other words, we can associate each vertex $v\in V$ with the monomial
\[
	p_v\coloneqq\tau(v)\yon^{\prod_{e\in E_v}\tau(t(e))}
\]
specifying its positions and directions, where $E_v\coloneqq s\inv(v)\ss E$ denotes the set of edges emanating from $v$.
The graph then determines a section
\[
\gamma\colon\bigotimes_{v\in V}p_v\to\yon
\]
given by a function
\[
    \prod_{v\in V}\tau(v)\too\prod_{e\in E}\tau(t(e))
\]
that sends each dependent function $i\colon(v\in V)\to\tau(v)$ to the dependent function $(e\in E)\to\tau(t(e))$ sending $e\mapsto i(t(e))$.
In other words, given the $p_v$-position $i(v)\in\tau(v)$ returned by each vertex $v\in V$, we know for each edge $e\in E$ that the direction that its source vertex $s(e)$ receives is the position $i(t(e))$ returned by its target vertex $t(e)$.
\index{dependent function}\index{dynamics}\index{graph}

So once we give dynamics to each $p_v$, namely by specifying a dynamical system $S_v\yon^{S_v}\to p_v$ with positions in $\tau(v)$ and directions in $\prod_{e\in E_v}\tau(t(e))$, we will obtain a closed dynamical system that updates the state of each vertex according to the information that they observe from each other along their edges.

Effectively, by interpreting a graph as a wiring diagram and giving each vertex dynamics, we have created what is known as a \emph{cellular automaton}---a network of vertices (or \emph{cells}), each with an internal state, in which each vertex $v\in V$ broadcasts a signal (i.e.\ returns a position) according to its current state, then updates its state according to the signals broadcasted by its \emph{neighbors} in $t(E_v)$ (i.e.\ positions returned by its neighbors it receives as directions).

For example, many cellular automata have their cells on a 2-dimensional integer lattice.
The corresponding graph has vertices $V\coloneqq\zz\times\zz$ and edges given by
\[
  E\coloneqq(\{-1,0,1\}\times\{-1,0,1\}-\{(0,0)\})\times V,
\]
with $s((i,j),(m,n))=(m,n)$ and $t((i,j),(m,n))=(m+i,n+j)$, so that the neighbors of each vertex are the eight vertices that surround it.
\end{example}

\begin{exercise}[Conway's Game of Life]\label{exc.conway}\index{graph}\index{interaction!cellular automata and}
Conway's Game of Life is a cellular automaton taking place on a 2-dimensional integer lattice as follows.
Each lattice point is either \emph{live} or \emph{dead}, and each point observes its eight \emph{neighbors} to which it is horizontally, vertically, or diagonally adjacent.
The following occurs at every time step:
\begin{itemize}
    \item Any live point with 2 or 3 live neighbors remains live.
    \item Any dead point with 3 live neighbors becomes live.
    \item All other points either become or remain dead.
\end{itemize}
We can use \cref{ex.graph_interaction} to model Conway's Game of Life as a closed dynamical system.
\begin{enumerate}
	\item What is the appropriate graph $E\tto V$?
	\item What is the appropriate assignment of sets $\tau\colon V\to\smset$?
	\item What are the monomials $p_v$ from \cref{ex.graph_interaction}?
	\item What is the appropriate state-set $S_v$ for each interface $p_v$?
	\item What is the appropriate dynamical system lens $S_v\yon^{S_v}\to p_v$?
\qedhere
\end{enumerate}
\begin{solution}
\begin{enumerate}
    \item Following the suggestion from the end of \cref{ex.graph_interaction}, we can use a graph with $V\coloneq\zz\times\zz$ and $E\coloneqq(\{-1,0,1\}\times\{-1,0,1\}-\{(0,0)\})\times V$ with $s((i,j),(m,n))=(m,n)$ and $t((i,j),(m,n))=(m+i,n+j)$ to model cellular automata like Conway's Game of Life on a 2-dimensional integer lattice in which each point observes only its eight immediate neighbors.
    \item Each vertex only needs to return whether it is live or dead, so we assign $\tau(v)\coloneqq\{\const{live},\,\const{dead}\}$ for every $v\in V$.
    \item For each $v\in V$, the monomial represented by $v$ from \cref{ex.graph_interaction} can be written as
    \[
        p_v\iso\{\const{live},\,\const{dead}\}\yon^{\smset(\{-1,0,1\}\times\{-1,0,1\}-\{(0,0)\},\, \{\const{live},\,\const{dead}\})}.
    \]
    Every vertex returns either $\const{live}$ or $\const{dead}$ as its position and receives as its direction whether each of its eight neighbors is $\const{live}$ or $\const{dead}$.
    \item Each vertex $v\in V$ only needs to record whether it is $\const{live}$ or $\const{dead}$, so $S_v\coloneqq\{\const{live},\,\const{dead}\}$.
    \item The appropriate dynamical system lens $S_v\yon^{S_v}\to p_v$ for each vertex $v\in V$ should have the identity function on $\{\const{live},\,\const{dead}\}$ as its return function, while its update map should be a function
    \[
        \{\const{live},\,\const{dead}\}\times\smset(\{-1,0,1\}\times\{-1,0,1\}-\{(0,0)\},\, \{\const{live},\,\const{dead}\})\to\{\const{live},\,\const{dead}\}
    \]
    that takes whether $v$ is $\const{live}$ or $\const{dead}$ as its first coordinate and a function from $\{-1,0,1\}\times\{-1,0,1\}-\{(0,0)\}$ to $\{\const{live},\,\const{dead}\}$ that says whether each of its eight neighbors is $\const{live}$ or $\const{dead}$ as a second coordinate, then executes the rules from Conway's Game of Life to determine whether it should be live or dead in the next time step.
\end{enumerate}
\end{solution}
\end{exercise}

\index{interaction!wiring diagram|)}

\subsection{More examples of general interaction}

While wiring diagrams are a handy visualization tool for certain simple interaction patterns, there are more general interaction patterns that cannot be captured by such a static diagram.
For example, here we generalize our previous cellular automata example.

\begin{example}[Generalized cellular automata: voting on who your neighbors are]\label{ex.cell_auto_vote_interaction}\index{cellular automata}\index{graph}
Recall from \cref{ex.graph_interaction} how we constructed a cellular automaton on a graph $G=(E\tto V)$.
For each $v\in V$, the graph specifies the set $N(v)\coloneqq t(E_v)$ of target vertices of edges emanating from $v$.
These vertices are the \emph{neighbors} of $v$, or the vertices that $v$ can ``observe.''
We call the function $N\colon V\to\2^V$ from each vertex to the set of its neighbors the \emph{neighbor function}.
For simplicitly, we let each vertex store and return one of two states, so $S_v\coloneqq\tau(v)\coloneqq\2$.

Now consider only the vertices of our graph and forget the edges.
Suppose we are then given a function $n\colon V\to\nn$ that we can think of as specifying the number $n(v)$ of neighbors each $v\in V$ could potentially have.
Let $\ord{n}(v)\coloneqq\{1,2,\ldots,n(v)\}$.
Then the monomial each vertex represents is
\[
    p_v\iso\2\yon^{\2^{\ord{n}(v)}},
\]
with its own possible states as positions and its potential neighbors' possible states as directions.

Say that a neighbor function $N\colon V\to\2^V$ \emph{respects} $n$ if we have an isomorphism $N(v)\iso\ord{n}(v)$ for each $v\in V$.
Now suppose we have a function $N'_-\colon \2^V\to (\2^V)^V$ that sends each set of vertices $S\in\2^V$ to a neighbor function $N'_S\colon V\to \2^V$ that respects $n$.
In other words, each possible state configuration $S$ of all the vertices in $V$ determines a neighbor function $N'_S$.
In the case of \cref{ex.graph_interaction}, when we had a graph, it told us what the neighbor function should always be.
Now we could imagine all the vertices returning their states to vote, via $N'$, on what neighbor function to use to determine which vertices are observing which others.

We can put this all together by providing a section for all the vertices,
\begin{equation}\label{eqn.polymap_misc9237}
    \bigotimes_{v\in V}p_v\cong\2^V\yon^{\2^{\sum_{v\in V}\ord{n}(v)}}\too\yon.
\end{equation}
Such a section is equivalent to a function $g\colon \2^V\to\2^{\sum_{v\in V}\ord{n}(v)}$ that sends each possible state configuration $S\in\2^V$ of all the vertices in $V$ to a function $g(S)\colon\sum_{v\in V}\ord{n}(v)\to\2$ specifying the states every vertex observes.
But we already have a neighbor function assigned to $S$ that respects $\ord{n}$, namely $N'_S$: we have $N'_S(v)\iso\ord{n}(v)$ for all $v\in V$.
So we can think of $g(S)$ equivalently as a function $g(S)\colon\sum_{v\in V}N'_S(v)\to\2$ that says for each $v\in V$ what signal in $\2$ it should receive from its neighbor $w\in N'_S(v)$.
We will let it receive the current state of that neighbor, as given by $S$:
\[
    g(S)(v,w)\coloneqq S(w).
\]

We have accomplished our goal: the vertices ``vote'' on how they should be connected in that their states together determine the neighbor function.
We do not mean to imply that this vote needs to be democratic or fair in any way: it is an arbitrary function $N'_-\colon \2^V\to(\2^V)^V$.
For instance, the state of a given vertex $v_0\in V$ may completely determine the neighbor function $V\to\2^V$; this would be expressed by saying that $N'_-$ factors as $\2^V\to\2^{\{v_0\}}\iso\2\to(\2^V)^V$, where the left function is the evaluation map at $v_0\in V$.
\end{example}

Here are some more examples of dependent dynamical systems with interaction patterns in which the way the constituent components are wired together may change.

\begin{example}\label{ex.bonds_break}\index{interaction!breaking bonds}
In the picture below, forces are applied to the connected boxes on the left; we would like to model how too much force could cause the connection between the boxes to sever, as depicted on the right.
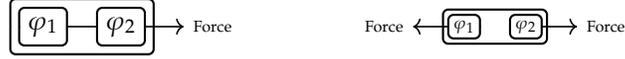
\[
\begin{tikzpicture}[oriented WD, bb small, bb port length=0]
	\node[bb={1}{1}] (x1) {$\varphi_1$};
	\node[bb={1}{1}, right=of x1] (x2) {$\varphi_2$};
	\node[bb={1}{1}, fit= (x1) (x2)] (outer) {};
5mm, font=\tiny] {Force};
	\draw (x1_out1) -- (x2_in1);
	\draw[->, shorten >= -4mm] (x2_out1) -- (outer_out1) node[right=4.5mm, font=\tiny] (L) {Force};
	\node[bb={1}{1}, right=7 of L] (y1) {$\varphi_1$};
	\node[bb={1}{1}, right=of y1] (y2) {$\varphi_2$};
	\node[bb={1}{1}, fit= (y1) (y2)] (outer) {};
	\draw[->, shorten >= -4mm] (y1_in1) -- (outer_in1) node[left=4.5mm, font=\tiny] (R){Force};
	\draw[->, shorten >= -4mm] (y2_out1) -- (outer_out1) node[right=4.5mm, font=\tiny] {Force};
\end{tikzpicture}
\]
We will imagine the dependent dynamical systems $\varphi_1\colon S\yon^S\to p_1$ and $\varphi_2\colon S\yon^S\to p_2$ as initially connected in space.
They experience forces from the outside world, and---for as long as they are connected---they experience forces from each other.
More precisely, their interfaces are given by
\[
	p_1\coloneqq p_2\coloneqq F\yon^{F\times F}+\{\const{snapped}\}\yon^F,
\]
where $F$ is our set of \emph{forces}.
We need to be able to add and compare forces, so we need $F$ to be an ordered monoid; let us say $F\coloneqq\nn$ for simplicity.
Here each interface has two kinds of positions it can return: either a force $f\in F$ that will be applied to the other interface (i.e.\ sent to the other interface as a direction) or $\const{snapped}$, indicating that the interfaces are no longer connected.
The interface always receives a force from the outside world as part of its direction, but when the position of the interface is not $\const{snapped}$, it receives a force from the other interface as part of its direction as well.
So its direction set is $F\times F$ at positions in $F$, when the two interfaces are connected; and just $F$ at the position $\const{snapped}$, when the two interfaces are not.
We define the wrapper interface to be
\[
    p\coloneqq\yon^{F\times F};
\]
it can return only 1 position, while its directions are ordered pairs of forces $(f_L, f_R)$ indicating the two external forces acting on the composite system.

\index{interaction pattern}
Though the systems $\varphi_1$ and $\varphi_2$ may be initially connected, if the forces on either one surpass a threshold, that system stops communicating with the other system.
The connection is broken and neither system ever receives forces from the other again. To implement this explicitly, we need to define an interaction pattern $\kappa\colon p_1\otimes p_2\to p$ that wraps $p$ around $\varphi_1$ and $\varphi_2$.
That is, we need to give a lens
\[
    \kappa\colon (F\yon^{F\times F}+\{\const{snapped}\}\yon^F)\otimes (F\yon^{F\times F}+\{\const{snapped}\}\yon^F)\to\yon^{F\times F}.
\]
Distributing and leveraging coproducts as usual, we find that it suffices to give four lenses:
\begin{equation}\label{eqn.snapped_maps}
\arraycolsep=1.4pt
\begin{array}{lll}
	\kappa_{11}\colon&~ F\times F\yon^{(F\times F)(F\times F)}&\to\yon^{F\times F}\\
	\kappa_{12}\colon&~ F\{\const{snapped}\}\yon^{(F\times F)F}&\to\yon^{F\times F}\\
	\kappa_{21}\colon&~ \{\const{snapped}\}F\yon^{F(F\times F)}&\to\yon^{F\times F}\\
	\kappa_{22}\colon&~ \{\const{snapped}\}\{\const{snapped}\}\yon^{F\times F}&\to\yon^{F\times F}
\end{array}
\end{equation}
The lenses $\kappa_{12}$ and $\kappa_{21}$ will not actually occur in our dynamics (when one interface returns $\const{snapped}$, both should), so we take them to be arbitrary.
We take the lens $\kappa_{22}$ to be the obvious isomorphism, passing the external forces to the two internal interfaces.
Finally, the lens $\kappa_{11}$ is equivalent to a function $(F\times F)(F\times F)\to (F\times F)(F\times F)$ which, taking care to remember what each $F$ refers to, we find should send $((f_1,f_2),(f_L,f_R))\mapsto((f_L,f_2),(f_1,f_R))$.
While the multiple $F$'s may be a little hard to keep track of, what this map says is that if $\varphi_1$ returns the force $f_1$ on $\varphi_2$ as output and $\varphi_2$ returns the force $f_2$ on $\varphi_1$ as output, then $\varphi_1$ receives the force $f_2$ from the right as input and $\varphi_2$ receives the force $f_1$ from the left as input; and in the meantime the left external force $f_L$ is given to $\varphi_1$ on the left, while the right external force is given to $\varphi_2$ on the right.

Now that we have the interfaces wrapped together, it remains to specify each dynamical system.
The state-sets for the two systems will be the same, namely $S\coloneqq F+\{\const{snapped}\}$: each system is either applying a force to the other system or not.
The dynamical systems themselves will be the same as well, up to a symmetry swapping left and right; we will define only the left system.
It is given by a lens
\[\varphi_1\colon (F+\{\const{snapped}\})\yon^{F+\{\const{snapped}\}}\to F\yon^{F\times F}+\{\const{snapped}\}\yon^F\]
which we write as the sum of two lenses
\[F\yon^{F+\{\const{snapped}\}}\to F\yon^{F\times F} \text{ and } \{\const{snapped}\}\yon^{F+\{\const{snapped}\}}\to\{\const{snapped}\}\yon^F.\]
Both lenses are identities on positions, directly returning their current states.
The second lens corresponds to when the connection is broken, after which the connection should remain broken: so its on-directions function is constant, sending any direction to $\const{snapped}$.
Meanwhile, the first lens corresponds to the case where the systems are still connected; in this state, the system can receive a pair of forces as its direction and must update its state---either the force it applies or $\const{snapped}$---accordingly.
We let the on-directions function $F(F\times F)\to F+\{\const{snapped}\}$ send
\[
(f_1,(f_L,f_2))\mapsto
\begin{cases}
	\const{snapped}&\tn{ if }f_1+f_2\geq100\\
	f_L&\tn{ otherwise}
\end{cases}
\]
Thus, when the sum of forces is above a certain threshold (arbitrarily chosen here to be $100$), the internal state is updated to the $\const{snapped}$ state; otherwise, the internal state is set to the external force applied to the system, which it is now ready to transfer to the other system.
\end{example}

\begin{example}\label{ex.supplier_change}\index{interaction!supplier change}
Consider the case of a company that may change its supplier based on its internal state. The company returns two possible positions, corresponding to whether it wants to receive gizmos in $G$ from the first supplier or widgets in $W$ from the second:
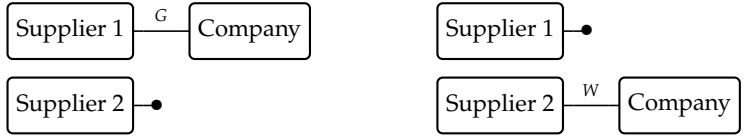
\[
\begin{tikzpicture}[oriented WD]
	\node[bb={0}{1}, font=\small] (s1) {Supplier 1};
	\node[bb={0}{1}, below=of s1, font=\small] (s2) {Supplier 2};
	\node[bb={1}{0}, right=0.5 of s1, font=\small] (c) {Company};
	\draw (s1_out1) to node[above, fill=none, font=\tiny] {$G$} (c_in1);
	\draw (s2_out1) to +(5pt,0) node[fill=none] {$\bullet$};
\begin{scope}[xshift=2.25in]
	\node[bb={0}{1}, font=\small] (s1') {Supplier 1};
	\node[bb={0}{1}, below=of s1', font=\small] (s2') {Supplier 2};
	\node[bb={1}{0}, right=0.5 of s2', font=\small] (c') {Company};
	\draw (s2'_out1) to node[above, fill=none, font=\tiny] {$W$} (c'_in1);
	\draw (s1'_out1) to +(5pt,0) node[fill=none] {$\bullet$};
\end{scope}
\end{tikzpicture}
\]
So the company has interface $\{1\}\yon^G+\{2\}\yon^W$, the first supplier has interface $G\yon$, and the second supplier has interface $W\yon$.
Then a section for the company and the suppliers is a lens
\[
  \left(\{1\}\yon^G+\{2\}\yon^W\right)\otimes G\yon\otimes W\yon\to\yon,
\]
corresponding to a pair of functions $\{1\}\times GW\iso GW\to G$ and $\{2\}\times GW\iso GW\to W$ given by canonical projections.
In other words, the company's position determines its supplier and what it receives.
\end{example}

\begin{example}\label{ex.assemble_machine}\index{interaction!assembling}
When someone assembles a machine, their own positions dictate the interaction pattern of the machine's components.
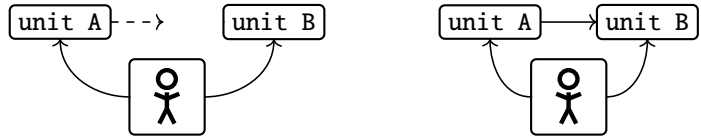
\begin{equation*}%\label{eqn.someone2}
\begin{tikzpicture}[oriented WD, font=\ttfamily, bb port length=0, baseline=(someone.north)]
	\node[bb port sep=.5, bb={0}{1}] (A) {unit A};
	\node[bb port sep=.5, bb={1}{0}, right=of A] (B) {unit B};
	\coordinate (helper) at ($(A)!.5!(B)$);
	\node[bb={1}{1}, below=2 of helper] (someone) {\tikzsymStrichmaxerl[3]};
	\draw[->] (someone_in1) to[out=180, in=270] (A.270);
	\draw[->] (someone_out1) to[out=0, in=270] (B.270);
	\draw[->, dashed] (A_out1) -- +(20pt,0);
\begin{scope}[xshift=2.25in]
	\node[bb port sep=.5, bb={0}{1}] (A') {unit A};
	\node[bb port sep=.5, bb={1}{0}, right=.5of A'] (B') {unit B};
	\coordinate (helper') at ($(A')!.5!(B')$);
	\node[bb={1}{1}, below=2 of helper'] (someone') {\tikzsymStrichmaxerl[3]};
	\draw[->] (someone'_in1) to[out=180, in=270] (A'.270);
	\draw[->] (someone'_out1) to[out=0, in=270] (B'.270);
	\draw[->] (A'_out1) -- (B'_in1);
\end{scope}
\end{tikzpicture}
\end{equation*}
Define $S\coloneqq\{\const{attach},\,\const{separate}\}$.
We say that \texttt{unit A} has interface
\[
  \left(\{\const{attached}\}\times X+\{\const{separated}\}\right)\yon^S.
\]
It can always receive either the direction to $\const{attach}$ or to $\const{separate}$ from $S$, while its position is either $\const{separated}$ or $\const{attached}$ and returning some element of a fixed set $X$.
Meanwhile, \texttt{unit B} has interface
\[
  \{\const{attached}\}\yon^{X\times S}+\{\const{separated}\}\yon^S.
\]
It can also always receive either direction from $S$, but when it is $\const{attached}$ it can further receive an element of $X$.
Finally, the role of the person is simply to return whether the units should $\const{attach}$ or $\const{separate}$, so we give it the interface $S\yon$.

Then a section for the person and the units is a lens
\begin{align*}
  &\left(\{\const{attached}\}\times X+\{\const{separated}\}\right)\yon^S \\
    \otimes
  &\left(\{\const{attached}\}\yon^{X\times S}+\{\const{separated}\}\yon^S\right)
    \otimes
  S\yon
    \to
  \yon.
\end{align*}
Such a lens corresponds to four functions, two of which can be arbitrary because our dynamics should never return them (either both units are $\const{attached}$ or both are $\const{separated}$).
The other two functions consist of one function
\[
  \{(\const{attached},\const{attached})\}\times X \times S\iso X\times S\to S\times X\times S
\]
that sends $(x,s)\mapsto(s,x,s)$ and another function
\[
  \{(\const{separated},\const{separated})\}\times S\iso S\to S\times S
\]
that sends $s\mapsto(s,s)$.

In words: the person's position tells the units whether they should $\const{attach}$ or $\const{separate}$.
If, and only if, the units are $\const{attached}$, one unit sends elements of $X$ to the other.
\end{example}

We can easily generalize \cref{ex.assemble_machine}.
Indeed, we will see in the next section that there is an interface $\ihom{q_1\otimes\cdots\otimes q_k\,,\,r}$ that represents all the interaction patterns between $q_1,\ldots,q_k$ with wrapper interface $r$, and that wrapping it around $p$ can be interpreted as a larger interaction pattern with wrapper interface $r$:
\[
\poly(p,[q_1\otimes\cdots\otimes q_k\,,\,r])\cong\poly(p\otimes q_1\otimes\cdots\otimes q_k\,,\,r).
\]
In other words, if the positions $p$ returns is deciding the interaction pattern between $q_1,\ldots,q_k$ with wrapper interface $r$, and the directions $p$ receives is from the subsequent behavior of that interaction pattern itself, then this is equivalent to an interaction pattern with wrapper interface $r$ that $p$ is part of alongside $q_1,\ldots,q_k$.

What it also means is that a dynamical system with interface $[q_1\otimes\cdots\otimes q_k,r]$ is simply selecting interaction patterns $q_1\otimes\cdots\otimes q_k\to r$.
Let us see how this works.

\index{interaction|)}
%-------- Section --------%
\section{Closure of $\otimes$}\label{sec.closure}%[-,-]
\index{closed monoidal structure}
\index{monoidal closed structure|see{closed monoidal structure}}
\index{parallel product!closure for}

The parallel monoidal product is closed---we have a closed monoidal structure on $\poly$---meaning that there is a closure operation, which we denote $\ihom{-,-}\colon\poly\op\times\poly\to\poly$, such that there is an isomorphism
\begin{equation}\label{eqn.monoidal_closure}
  \poly(p\otimes q,r) \iso \poly(p,\ihom{q,r})
\end{equation}
natural in $p,q,r$.
The closure operation is defined on $q,r$ as follows:
\begin{equation}\label{eqn.par_hom}
	\ihom{q,r} \coloneqq \prod_{j\in q(\1)}r\circ(q[j]\yon)
\end{equation}
Here $\circ$ denotes standard functor composition; informally, $r \circ (q[j]\yon)$ is the polynomial obtained by replacing each appearance of $\yon$ in $r$ by $q[j]\yon$.
Composition, together with the unit $\yon$, is in fact yet another monoidal structure, as we will cover in greater depth in \cref{part.comon}.

Before we prove that the isomorphism \eqref{eqn.monoidal_closure} holds naturally, let us investigate the properties of the closure operation, starting with some simple examples.

\begin{exercise}
Calculate $\ihom{q,r}$ for $q,r\in\poly$ given as follows.
\begin{enumerate}
	\item $q\coloneqq \0$ and $r$ arbitrary.
	\item $q\coloneqq \1$ and $r$ arbitrary.
	\item $q\coloneqq\yon$ and $r$ arbitrary.
	\item $q\coloneqq A$ for $A\in\smset$ (constant) and $r$ arbitrary.
	\item $q\coloneqq A\yon$ for $A\in\smset$ (linear) and $r$ arbitrary.
	\item $q\coloneqq\yon^\2+\2\yon$ and $r\coloneqq\2\yon^\3+\3$.
\qedhere
\end{enumerate}
\begin{solution}
We compute $\ihom{q,r}$ for various values of $q,r\in\poly$ using \eqref{eqn.par_hom}.
\begin{enumerate}
  \item If $q\coloneqq\0$, then $q(\1)\iso\0$, so $\ihom{q,r}$ is an empty product.
    Hence $\ihom{q,r}\iso\1$.
  \item If $q\coloneqq\1$, then $q(\1)\iso\1$ and $q[1]\iso\0$, so $\ihom{q,r}\iso r\circ(\0\yon)\iso r(\0)$.
  \item If $q\coloneqq\yon$, then $q(\1)\iso\1$ and $q[1]\iso\1$, so $\ihom{q,r}\iso r\circ(\1\yon)\iso r$.
	\item If $q\coloneqq A$ for $A\in\smset$, then $q(\1)\iso A$ and $q[j]\iso\0$ for every $j\in A$, so $\ihom{q,r}\iso\prod_{j\in A}(r\circ(\0\yon))\iso r(\0)^A$.
	\item If $q\coloneqq A\yon$ for $A\in\smset$, then $q(\1)\iso A$ and $q[j]\iso\1$ for every $j\in A$, so $\ihom{q,r}\iso\prod_{j\in A}(r\circ(\1\yon))\iso r^A$.
	\item If $q\coloneqq\yon^\2+\2\yon$ and $r\coloneqq\2\yon^\3+\3$, then
	\begin{align*}
	    \ihom{q,r} &\iso (r \circ (\2\yon))(r \circ (\1\yon))^\2 \\
	    &\iso \left(\2(\2\yon)^\3 + \3\right)\left(\2\yon^\3 + \3\right)^\2 \\
	    &\iso \6\4\yon^\9 + \2\0\4\yon^\6 + \1\8\0\yon^\3 + \2\7.
	\end{align*}
\end{enumerate}
\end{solution}
\end{exercise}

\begin{exercise}\label{exc.sum_times_closure}
Show that for any polynomials $p_1,p_2,q$, we have an isomorphism
\[
\ihom{p_1 + p_2, q} \iso \ihom{p_1, q} \times \ihom{p_2, q}.
\]
\begin{solution}
We wish to show that for all $p_1, p_2, q \in \poly$, we have $\ihom{p_1 + p_2, q} \iso \ihom{p_1, q} \times \ihom{p_2, q}$.
By \eqref{eqn.par_hom},
\[
    \ihom{p_1 + p_2, q} \iso \left(\prod_{i \in p_1(\1)} q \circ (p_1[i]\yon)\right) \left(\prod_{i \in p_2(\1)} q \circ (p_2[i]\yon)\right) \iso \ihom{p_1, q} \times \ihom{p_2, q}.
\]
\end{solution}
\end{exercise}

\begin{exercise} \label{exc.par_hom_sum}
Show that there is an isomorphism
\begin{equation} \label{eqn.par_hom_sum}
\scalebox{1.3}{$\displaystyle
\ihom{q,r} \iso \sum_{f\colon q\to r}\yon^{\sum_{j\in q(\1)}r[f_\1j]}
$}
\end{equation}
where the sum is indexed over $f\in\poly(q,r)$.
\begin{solution}
We may compute
\begin{align*}
    \ihom{q, r} &\iso \prod_{j \in q(\1)} r \circ (q[j]\yon) \tag*{\eqref{eqn.par_hom}} \\
    &\iso \prod_{j \in q(\1)} \, \sum_{k \in r(\1)} (q[j]\yon)^{r[k]} \tag{Replacing each $\yon$ in $r$ by $q[j]\yon$} \\
    &\iso \sum_{f_\1 \colon q(\1) \to r(\1)} \, \prod_{j \in q(\1)} (q[j]\yon)^{r[f_\1(j)]} \tag*{\eqref{eqn.push_prod_sum_set_indep}} \\
    &\iso \sum_{f_\1 \colon q(\1) \to r(\1)} \, \left(\prod_{j \in q(\1)} q[j]^{r[f_\1(j)]} \right)\left(\prod_{j \in q(\1)} \yon^{r[f_\1(j)]} \right) \\
    &\iso \sum_{f_\1 \colon q(\1) \to r(\1)} \; \sum_{f^\sharp \in \prod_{j \in q(\1)} q[j]^{r[f_\1(j)]}} \yon^{\sum_{j \in q(\1)} r[f_\1(j)]} \\
    &\iso \sum_{f \colon q \to r} \yon^{\sum_{j \in q(\1)} r[f_\1(j)]}. \tag*{\eqref{eqn.main_formula}}
\end{align*}
\end{solution}
\end{exercise}

\begin{exercise} \label{exc.dir_hom_p_yon_dir_p}
Verify that \eqref{eqn.dir_hom_p_yon_dir_p} holds.
\begin{solution}
We verify \eqref{eqn.dir_hom_p_yon_dir_p} as follows:
\begin{align*}
    \ihom{p, \yon} \otimes p
    &\iso
    \left(\sum_{f \colon p \to \yon} \yon^{\sum_{i \in p(\1)} \yon[f_\1i]}\right) \otimes p
    \tag*{\eqref{eqn.par_hom_sum}} \\
    &\iso
    \sum_{f \in \Gamma(p)} \yon^{p(\1)} \otimes \sum_{i \in p(\1)} \yon^{p[i]} \\
    &\iso
    \sum_{f \in \Gamma(p)} \; \sum_{i \in p(\1)} \yon^{p(\1) \times p[i]}
    \tag*{\eqref{eqn.parallel_def}} \\
    &\iso
    \sum_{f \in \prod_{i \in p(\1)} p[i]} \; \sum_{i \in p(\1)} \yon^{p(\1) \times p[i]}.
    \tag*{\eqref{eqn.gamma_prod}}
\end{align*}
\end{solution}
\end{exercise}

\begin{example}\label{ex.parallel_dual}
For any $A\in\smset$ we have
\[
  \ihom{\yon^A,\yon} \iso A\yon
  \qqand
  \ihom{A\yon,\yon} \iso \yon^A.
\]
More generally, for any polynomial $p\in\poly$ we have
\begin{equation}\label{eqn.dir_dual}
  \ihom{p,\yon} \iso \Gamma(p)\yon^{p(\1)}.
\end{equation}
All these facts follow directly from \eqref{eqn.par_hom}.
\end{example}

\begin{exercise}
Verify the three facts above.
\begin{solution}
We have that
\[
    \ihom{\yon^A, \yon} \iso \prod_{j \in \yon^A(\1)} \yon \circ (\yon^A[j]\yon) \iso \prod_{j \in \1} A\yon \iso A\yon,
\]
that
\[
    \ihom{A\yon, \yon} \iso \prod_{j \in A\yon(\1)} \yon \circ ((A\yon)[j]\yon) \iso \prod_{j \in A} \yon \iso \yon^A,
\]
and that
\begin{align*}
    \ihom{p, \yon} &\iso \sum_{f \colon p \to \yon} \yon^{\sum_{i \in p(\1)} \yon[f_\1i]} \tag*{\eqref{eqn.par_hom_sum}} \\
    &\iso \sum_{f \in \Gamma(p)} \yon^{\sum_{i \in p(\1)} \1} \\
    &\iso \Gamma(p)\yon^{p(\1)}.
\end{align*}
\end{solution}
\end{exercise}

\begin{exercise}
Show that for any $p\in\poly$, if there is an isomorphism $\ihom{\ihom{p,\yon},\yon} \iso p$, then $p$ is either linear $A\yon$ or representable $\yon^A$ for some $A$. Hint: first show that $p$ must be a monomial.
\begin{solution}
Given $p \in \poly$ and an isomorphism $\ihom{\ihom{p,\yon},\yon} \iso p$, we wish to show that $p$ is either linear or representable.
Applying \eqref{eqn.dir_dual} twice, we have that
\[
    \ihom{\ihom{p,\yon},\yon} \iso \Gamma\left(\Gamma(p)\yon^{p(\1)}\right)\yon^{\Gamma(p)}.
\]
By \eqref{eqn.gamma_prod},
\[
    \Gamma\left(\Gamma(p)\yon^{p(\1)}\right) \iso \prod_{\gamma \in \Gamma(p)} p(\1) \iso p(\1)^{\Gamma(p)}.
\]
Hence taking $\ihom{\ihom{p,\yon},\yon} \iso p$ and rewriting the left hand side using the isomorphisms above yields
\begin{equation} \label{eqn.p_as_gamma_monomial}
    p(\1)^{\Gamma(p)}\yon^{\Gamma(p)}\iso p.
\end{equation}
In particular, $p$ is a monomial, so we can write $p \coloneqq I\yon^A$ for some $I,A \in \smset$.
Then $p(\1) \iso I$ and \eqref{eqn.gamma_prod} tells us that $\Gamma(p) \iso A^I$.
Equating the direction-sets on either side of \eqref{eqn.p_as_gamma_monomial} yields $A^I\iso A$; then equating position-sets gives $I^A\iso I^{A^I}\iso I$.

We conclude with some elementary set theory.
If either one of $I$ or $A$ were (isomorphic to) $\1$, then $p$ would be either linear or representable, and we would be done.
Meanwhile, if either one of $I$ or $A$ were $\0$, then the other would be $\1$, and we would again be done.
Otherwise, $|A|,|B| \geq 2$.
But by Cantor's theorem,
\[
    |I| < \big|\2^I\big| \leq \big|A^I\big| = |A| \qqand |A| < \big|\2^A\big| \leq \big|I^A\big| = |I|,
\]
a contradiction.
\end{solution}
\end{exercise}

\begin{proposition}\label{prop.parallel_closure}\index{isomorphism!natural}
With $\ihom{-,-}$ as defined in \eqref{eqn.par_hom}, there is a natural isomorphism
\begin{equation}\label{eqn.poly_closure_brackets}
	\poly(p\otimes q,r)\cong\poly(p,\ihom{q,r}).
\end{equation}
\end{proposition}
\begin{proof}
We have the following chain of natural isomorphisms:
\begin{align*}
	\poly(p\otimes q,r)
	&\iso
	\poly\Big(\sum_{i\in p(\1)}\sum_{j\in q(\1)}\yon^{p[i]q[j]},r\Big) \\
	&\iso
	\prod_{i\in p(\1)}\prod_{j\in q(\1)}\poly(\yon^{p[i]q[j]},r)
	\tag{Universal property of coproducts} \\
	&\iso
	\prod_{i\in p(\1)}\prod_{j\in q(\1)}r(p[i]q[j])
	\tag{Yoneda lemma} \\
	&\iso
	\prod_{i\in p(\1)}\prod_{j\in q(\1)}\poly(\yon^{p[i]},r\circ(q[j]\yon))
	\tag{Yoneda lemma} \\
	&\iso
	\poly\Big(\sum_{i\in p(\1)}\yon^{p[i]},\prod_{j\in q(\1)}r\circ(q[j]\yon)\Big)
	\tag{Universal property of (co)products} \\
	&\iso
	\poly(p,\ihom{q,r}).
\end{align*}
\end{proof}\index{Yoneda lemma}

\begin{exercise}\label{exc.poly_plug_1}
Show that for any $p,q$ we have an isomorphism of sets
\[
\poly(p,q) \iso \ihom{p,q}(\1).
\]
Hint: you can either use the formula \eqref{eqn.par_hom}, or just use
\eqref{eqn.poly_closure_brackets} with the Yoneda lemma and the fact that $\yon\otimes p \iso p$.
\begin{solution}
The isomorphism $\poly(p,q) \iso \ihom{p,q}(\1)$ follows directly from \cref{exc.par_hom_sum} when both sides are applied to $\1$.
Alternatively, we can apply \eqref{eqn.poly_closure_brackets}.
Since $p \iso \yon \otimes p$, we have that
\begin{align*}
    \poly(p, q) &\iso \poly(\yon \otimes p, q) \\
    &\iso \poly(\yon, \ihom{p,q}) \tag*{\eqref{eqn.poly_closure_brackets}} \\
    &\iso \ihom{p,q}(\1). \tag{Yoneda lemma}
\end{align*}
\end{solution}
\end{exercise}\index{Yoneda lemma}

The closure of $\otimes$ implies that for any $q,r\in\poly$, there is a canonical \emph{evaluation} lens
\begin{equation}\label{eqn.eval_parallel}
  \fun{eval}\colon \ihom{q,r}\otimes q\to r
\end{equation}
given by sending the identity lens on $\ihom{q,r}$ leftward through the natural isomorphism
\[
  \poly(\ihom{q,r}\otimes q,r)\iso\poly(\ihom{q,r},\ihom{q,r})
\]
As in any closed monoidal category, such an evaluation lens has the universal property that for any $p\in\poly$ and lens $f\colon p\otimes q\to r$, there is a unique lens $f'\colon p\to\ihom{q,r}$ such that the following diagram commutes:
\[
    \begin{tikzcd}
    	p\otimes q\ar[r, "f'\:\otimes\:q"]\ar[rr, bend right, "f"']&
    	{\ihom{q,r}}\otimes q\ar[r, "\fun{eval}"]&
    	r
    \end{tikzcd}
\]

\begin{exercise} \label{exc.eval_parallel}
Describe the behavior of the evaluation lens $\fun{eval}\colon \ihom{p,q}\otimes p\to q$ from \eqref{eqn.eval_parallel}.
\begin{solution}
To obtain the evaluation lens $\fun{eval}\colon \ihom{q,r}\otimes q\to r$, we need to send the identity lens on $\ihom{q,r}$ leftward through the natural isomorphism
\[
\poly(\ihom{q,r}\otimes q,r)\iso\poly(\ihom{q,r},\ihom{q,r})
\]
To do so, we can start from the identity lens on $\ihom{q,r}$ and work our way along a chain of natural isomorphisms from $\poly(\ihom{q,r},\ihom{q,r})$ until we get to $\poly(\ihom{q,r}\otimes q,r)$.
To start, \cref{exc.par_hom_sum} implies that
\begin{align*}\index{isomorphism!natural}
    \poly(\ihom{q,r},\ihom{q,r})
    &\iso
    \poly\left(\sum_{f\colon q\to r} \, \prod_{i'\in q(\1)} \yon^{r[f_\1i']}, \prod_{i \in q(\1)} \, \sum_{j \in r(\1)} (q[i]\yon)^{r[j]}\right) \\
    &\iso
    \prod_{f\colon q\to r} \, \prod_{i\in q(\1)} \poly\left(\prod_{i' \in q(\1)} \yon^{r[f_\1i']}, \sum_{j \in r(\1)} (q[i]\yon)^{r[j]}\right),
\end{align*}
where the second isomorphism follows from the universal properties of products and coproducts.
In particular, under this isomorphism, the identity lens on $\ihom{q,r}$ corresponds to a collection of lenses, namely for each $f\colon q\to r$ and each $i\in q(\1)$ the composite
\[
    \prod_{i' \in q(\1)} \yon^{r[f_\1i']} \to \yon^{r[f_\1i]} \to \sum_{g \colon r[f_\1i] \to q[i]} \yon^{r[f_\1i]} \iso (q[i]\yon)^{r[f_\1i]} \to \sum_{j \in r(\1)} (q[i]\yon)^{r[j]}
\]
of the canonical projection with index $i'\coloneqq i$, the canonical inclusion with index $g\coloneqq f^\sharp_i$, and the canonical inclusion with index $j\coloneqq f_\1i$.
On positions, this lens picks out the position of $\sum_{j \in r(\1)} (q[i]\yon)^{r[j]}$ corresponding to $j = f_\1i \in r(\1)$ and $f^\sharp_i \colon r[f_\1i] \to q[i]$; on directions, the lens is the canonical inclusion $r[f_\1i] \to \sum_{i' \in q(\1)} r[f_\1i']$ with index $i' = i$.
By the Yoneda lemma, we can reinterpret each of these lenses as a lens
\[
    \yon^{q[i] \times \sum_{i' \in q(\1)} r[f_\1i']} \to \sum_{j \in r(\1)} \yon^{r[j]} \iso r
\]
that, on positions, picks out the position $f_\1i \in r(\1)$ of $r$ and, on directions, is the map $r[f_\1i] \to q[i] \times \sum_{i'\in q(\1)} r[f_\1i']$ induced by the universal property of products applied to the map $f^\sharp_i \colon r[f_\1i] \to q[i]$ and the inclusion $r[f_\1i] \to \sum_{i'\in q(\1)} r[f_\1i']$.
Then by the universal property of coproducts, this collection of lenses induces a single lens $\fun{eval}\colon\ihom{q,r}\otimes q\to r$ that sends each position $f\colon q\to r$ of $\ihom{q,r}$ and position $i\in q(\1)$ of $q$ to the position $f_\1i$ of $r$, with the same behavior on directions as the corresponding lens described previously.
\end{solution}
\end{exercise}

\begin{exercise}
\begin{enumerate}
	\item For any set $S$, obtain the do-nothing section $S\yon^S\to\yon$ from \cref{ex.do_nothing} whose on-directions is the identity on $S$ using eval and \cref{ex.parallel_dual}.
	\item Show that four lenses in \eqref{eqn.snapped_maps} from \cref{ex.bonds_break}, written equivalently as
	\begin{equation} \label{eqn.snapped_maps2}
	\arraycolsep=1.4pt
    \begin{array}{lll}
    	\kappa_{11}\colon&~ F\yon^{FF}\otimes F\yon^{FF}&\to\yon^F\otimes\yon^F\\
    	\kappa_{12}\colon&~ F\yon^{FF}\otimes \yon^F&\to\yon^F\otimes\yon^F\\
    	\kappa_{21}\colon&~ \yon^F\otimes F\yon^{FF}&\to\yon^F\otimes\yon^F\\
    	\kappa_{22}\colon&~ \yon^F\otimes\yon^F&\to\yon^F\otimes\yon^F,
    \end{array}
	\end{equation}
	can be obtained by taking the parallel product of identity lenses and evaluation lenses.
\qedhere
\end{enumerate}
\begin{solution}
\begin{enumerate}
    \item Given a set $S$, \cref{ex.parallel_dual} shows that
    \[
        \ihom{S\yon, \yon} \otimes (S\yon) \iso \yon^S \otimes (S\yon) \iso S\yon^S,
    \]
    so by setting $q\coloneqq S\yon$ and $r\coloneqq\yon$ in \eqref{eqn.eval_parallel}, we obtain an evaluation lens $\fun{eval} \colon S\yon^S \to \yon$.
    By the solution to \cref{exc.eval_parallel}, given a position $s \in S$ of $S\yon^S$, the evaluation lens on directions is the map $\1 \to S$ that picks out $s$.
    In other words, it is indeed the identity on directions.
    \item We wish to write the four lenses in \eqref{eqn.snapped_maps2} from \cref{ex.bonds_break} as the parallel product of identity lenses and evaluation lenses.
    By the solution to \cref{exc.eval_parallel}, the evaluation lens $\ihom{F\yon,\yon^F}\otimes F\yon\to\yon^F$
    is a lens from
    \[
        \ihom{F\yon,\yon^F}\otimes F\yon\iso F\left(\sum_{f\colon F\yon\to\yon^F}\,\prod_{i\in F}\yon^F\right)\iso F\yon^{FF}
    \]
    to $\yon^F$ that is uniquely determined on positions and has the on-directions map $FF\to FF$ given by the identity.
    Then we can verify that $\kappa_{11}$ is equivalent to the parallel product of this evaluation lens with itself.
    We can define $\kappa_{12}$ and $\kappa_{21}$ to be the parallel product of this evaluation lens with the identity on $\yon^F$, while $\kappa_{22}$ is the parallel product of the identity on $\yon^F$ with itself.
\end{enumerate}
\end{solution}
\end{exercise}

\begin{example}[Modeling your environment without knowing what it is]\index{environment!universal}\index{interaction pattern}
Imagine a robot whose interface is an arbitrary polynomial $q$ that is part of an interaction pattern $f\colon p\otimes q\to r$ alongside its environment with interface $p$. Then $f$ induces a lens $f'\colon p\to \ihom{q,r}$ such that the following diagram commutes:
\[
\begin{tikzcd}
  p\otimes q\ar[r, "f'\:\otimes\:q"]\ar[rr, bend right, "f"']&
  {\ihom{q,r}}\otimes q\ar[r, "\fun{eval}"]&
  r
\end{tikzcd}
\]

In other words, $\ihom{q,r}$ holds within it all of the possible ways $q$ can interact with other systems when they are wrapped in $r$ together.
For example, in the case of $r\coloneqq\yon$, note that $\ihom{q,\yon}\iso\prod_{i\in q(\1)}q[i]\yon$.
That is, a position of $\ihom{q,\yon}$ sends each $q$-position to a direction at that position, which is exactly what we need to specify to put $q$ in a closed system.

\index{dynamics}\index{environment}

Now suppose we were to give dynamics to $q$ by specifying a lens $S\yon^S\to q$. One could aim to choose a set $S$ along with an interesting map $g\colon S\to\poly(q,r)$. Then each state $s$ would include a guess $g(s)$ about the state of its environment. This is not the real environment $p$, but just the environment as it affects $q$, namely $\ihom{q,r}$. The robot's states model its environmental conditions.
\end{example}

\begin{example}[Chu $\&$]\index{Chu space}
Say we have $p_1,p_2,q_1,q_2,r\in\poly$ and lenses
\[
	\varphi_1\colon p_1\otimes q_1\to r
	\qqand
	\varphi_2\colon p_2\otimes q_2\to r.
\]
One might call these ``$r$-Chu spaces.'' One operation you can do with these as Chu spaces is to return something denoted $\varphi_1\&\varphi_2$, or ``$\varphi_1$ \emph{with} $\varphi_2$,'' of the following type:
\[
\varphi_1\&\varphi_2\colon (p_1\times p_2)\otimes (q_1+q_2)\to r
\]
Suppose we are given a position in $p_1$ and a position in $p_2$. Then given a position in either $q_1$ or $q_2$, one evaluates either $\varphi_1$ or $\varphi_2$ respectively to get a position in $r$; given a direction there, one returns the corresponding direction in $q_1$ or $q_2$ respectively, as well as a direction in $p_1\times p_2$ which is either a direction in $p_1$ or in $p_2$.

This sounds complicated, but it can be constructed formally via the monoidal closure.
We use the closure to rewrite $\varphi_1$ and $\varphi_2$:
\[
\psi_1\colon p_1\to \ihom{q_1,r}
\qqand
\psi_2\colon p_2\to \ihom{q_2,r}
\]
Now we take their categorical product to obtain $\psi_1\times\psi_2\colon p_1\times p_2\to\ihom{q_1,r}\times\ihom{q_2,r}$.
Then we apply \cref{exc.sum_times_closure} to find that $\ihom{q_1,r}\times\ihom{q_2,r}\cong\ihom{q_1+q_2,r}$, and finally leverage the monoidal closure again to obtain $(p_1\times p_2)\otimes(q_1+q_2)\to r$ as desired.
\end{example}

%-------- Section --------%
\section[Summary and further reading]{Summary and further reading%
  \sectionmark{Summary \& further reading}}
\sectionmark{Summary \& further reading}

In this chapter we explained how discrete dynamical systems can be expressed as certain lenses between polynomial functors. For example, a Moore machine has an input set $A$, an output set $B$, a set of states $S$, a return function $S\to B$, and an update function $A\times S\to S$. All this is captured in a depedent lens
\[S\yon^S\to B\yon^A.\]
We discussed a generalization $S\yon^S\to p$, where the output is an arbitrary polynomial $p\in\poly$. We also talked about how to wire machines in parallel by using the parallel product $\otimes$ and how to add wrapper interfaces by composing with lenses $p\to q$.

Throughout the chapter we gave quite a few different examples. For example, we discussed how every function $A\to B$ counts as a memoryless dynamical system. In fact, it was shown in \cite{beurier2019memoryless} that every dynamical system can be obtained by wiring together memoryless ones. We discussed examples such as file-readers, moving robots, colliding particles, companies that change their suppliers, materials that break when too much force is applied, etc.

For further reading on the mathematics of Moore machines, see \cite{conway2012regular}. For more on mode-dependent interaction, see \cite{spivak2017nesting}. For a similar and complementary categorical approach to dynamical systems, we recommend David Jaz Myers' \emph{Categorical Systems Theory} book, found at \url{http://davidjaz.com/Papers/DynamicalBook.pdf}.

\index{interaction!mode dependent}

%-------- Section --------%
\section{Exercise solutions}
\Closesolutionfile{solutions}
{\footnotesize
	\input{solution-file4}}

\Opensolutionfile{solutions}[solution-file5]

%------------ Chapter ------------%
\chapter{More categorical properties of polynomials}
\chaptermark{More categorical properties of $\poly$}
\label{ch.poly.bonus}

The category $\poly$ has very useful formal properties, including many adjunctions with $\smset$, factorization systems, cartesian closure, completion under colimits and limits, and so on.
We detail some of these properties here.
Most of the following material is not necessary for the development of our main story, but we collect it here for reference.
The reader may skip directly to \cref{part.comon} if so inclined and check back here when needed; or they may skim the chapter to see how well-behaved and versatile $\poly$ is.

%-------- Section --------%
\section{Special polynomials and adjunctions} \label{sec.poly.bonus.adj}

In \cref{sec.poly.obj.spec}, we identified four special classes of polynomials:
\begin{enumerate}[label=\alph*)]
	\item constant polynomials $I$ for $I\in\smset$;
	\item linear polynomials $I\yon$ for $I\in\smset$;
	\item representable polynomials $\yon^A$ for $A\in\smset$; and
	\item monomials $I\yon^A$ for $I,A\in\smset$.
\end{enumerate}
The first two classes, constant and linear polynomials, are interesting because they may each be viewed as a copy of $\smset$ inside $\poly$, as we will see in \cref{prop.ff_const_set_to_poly,prop.ff_lin_set_to_poly}.
The third class may be viewed as a copy of $\smset\op$ inside $\poly$ via the Yoneda embedding, as we saw in \cref{exc.finish_proof_yoneda}.
Finally, the fourth may be viewed as a copy of the category of bimorphic lenses in $\poly$, as we saw in \cref{subsec.poly.cat.morph.bimorphic-lens}.

\index{Yoneda embedding|see{Yoneda lemma}}

\begin{exercise}
Which of the four classes above are closed under
\begin{enumerate}
	\item coproducts?
	\item cartesian products?
	\item parallel products?
	\item composition? Recall that the composite of $p,q\in\poly$ is $p\circ q$, given by replacing each appearance of $\yon$ in $p$ by $q$.
\qedhere
\end{enumerate}
\begin{solution}
Here $I,J,A,B\in\smset$.
\begin{enumerate}
    \item We determine whether each class of polynomials is closed under coproducts.
    \begin{enumerate}
        \item Constants are closed under coproducts: given constants $I,J$, their coproduct $I+J$ is also a constant.
        \item Linear polynomials are closed under coproducts: given linear polynomials $I\yon, J\yon$, their coproduct $I\yon + J\yon \iso (I+J)\yon$ is also linear.
        \item Representables are \emph{not} closed under coproducts: for example, $\yon$ is a representable, but the coproduct of $\yon$ with itself, $\2\yon$, is not.
        \item Monomials are \emph{not} closed under addition: for example, $\yon$ and $\2\yon^\3$ are monomials, but their coproduct $\yon+\2\yon^\3$ is not.
    \end{enumerate}
    \item We determine whether each class of polynomials is closed under products.
%    The results below follow from \cref{exc.general_poly_times} \cref{exc.general_poly_times.monomial}.
    \begin{enumerate}
        \item Constants are closed under products: given constants $I, J$, their product $IJ$ is also a constant.
        \item Linear polynomials are \emph{not} closed under products: for example, $\yon$ and $\2\yon$ are linear polynomials, but their product $\2\yon^\2$ is not.
        \item Representables are closed under products: given representables $\yon^A, \yon^B$, their product $\yon^{A+B}$ is also a representable.
        \item Monomials are closed under products: given monomials $I\yon^A, J\yon^B$, their product $IJ\yon^{A+B}$ is also a monomial.
    \end{enumerate}
    \item We determine whether each class of polynomials is closed under parallel products.
%    The results below follow from \cref{exc.general_poly_parallel_times} \cref{exc.general_poly_parallel_times.monomial}.
    \begin{enumerate}
        \item Constants are closed under parallel products: given constants $I, J$, their parallel product $IJ$ is also a constant.
        \item Linear polynomials are closed under parallel products: given linear polynomials $I\yon, J\yon$, their parallel product $IJ\yon$ is also linear.
        \item Representables are closed under parallel products: given representables $\yon^A, \yon^B$, their parallel product $\yon^{AB}$ is also a representable.
        \item Monomials are closed under parallel products: given monomials $I\yon^A, J\yon^B$, their parallel product $IJ\yon^{AB}$ is also a monomial.
    \end{enumerate}
    \item We determine whether each class of polynomials is closed under composition.
    \begin{enumerate}
        \item Constants are closed under composition: given constants $I,J$, their composite $I\circ J\iso I$ is also a constant.
        \item Linear polynomials are closed under composition: given linear polynomials $I\yon,J\yon$, their composite $I\yon\circ J\yon\iso I(J\yon)\iso IJ\yon$ is also linear.
        \item Representables are closed under composition: given representables $\yon^A, \yon^B$, their composite $\yon^A\circ\yon^B\iso(\yon^B)^A\iso\yon^{BA}$ is also a representable.
        \item Monomials are closed under composition: given monomials $I\yon^A, J\yon^B$, their composite $I\yon^A\circ J\yon^B\iso I(J\yon^B)^A \iso IJ^A\yon^{BA}$ is also a monomial.
    \end{enumerate}
\end{enumerate}
\end{solution}
\end{exercise}

\begin{proposition}\label{prop.ff_const_set_to_poly}
There is a fully faithful functor $\smset\to\poly$ sending sets to constants: $I\mapsto I\yon^\0=I$.
\end{proposition}
\begin{proof}
A lens $f\colon I\yon^\0\to J\yon^\0$ consists of a function $f_\1\colon I\to J$ and, for each $i\in I$, a function $\0\to\0$. There is exactly one function $\0\to\0$, so $f$ can be identified with just a function between sets $I\to J$; in particular $\smset(I,J)\iso\poly(I\yon^\0,J\yon^\0)$.
This identification respects identities and composition, so it defines a fully faithful functor.
\end{proof}

\begin{proposition}\label{prop.ff_lin_set_to_poly}
There is a fully faithful functor $\smset\to\poly$ sending $I\mapsto I\yon$.
\end{proposition}
\begin{proof}
A lens $f\colon I\yon\to J\yon$ consists of a function $f_\1\colon I\to J$ and, for each $i\in I$, a function $\1\to\1$. There is exactly one function $\1\to\1$, so $f$ can be identified with just a function between sets $I\to J$; in particular $\smset(I,J)\iso\poly(I\yon,J\yon)$.
This identification respects identities and composition, so it defines a fully faithful functor.
\end{proof}

\index{adjunction!between $\smset$ and $\poly$}
\index{set!as constant polynomial}

\begin{theorem}\label{thm.adjoint_quadruple}
$\poly$ has an adjoint quadruple with $\smset$:%
\tablefootnote{We denote an adjunction with a double arrow $\imp$ pointing in the direction of the left adjoint.}
\begin{equation}\label{eqn.adjoints_galore}
\begin{tikzcd}[column sep=60pt, background color=theoremcolor]
  \smset
  	\ar[r, shift left=7pt, "I" description]
		\ar[r, shift left=-21pt, "I\yon"']&
  \poly
  	\ar[l, shift right=21pt, "q(\0)"']
  	\ar[l, shift right=-7pt, "q(\1)" description]
	\ar[l, phantom, "\scriptstyle\Leftarrow"]
	\ar[l, phantom, shift left=14pt, "\scriptstyle\Rightarrow"]
	\ar[l, phantom, shift right=14pt, "\scriptstyle\Rightarrow"]
\end{tikzcd}
\end{equation}
where the functors have been labeled by where they send $I\in\smset$ and $q\in\poly$.
\end{theorem}
\begin{proof}\index{Yoneda lemma}
For each $A\in\smset$, there is a representable functor $\poly\to\smset$ sending $q\mapsto\poly(\yon^A,q)$.
By the Yoneda lemma, this functor is isomorphic to the functor sending $q\mapsto q(A)$.
This, along with \cref{prop.ff_const_set_to_poly,prop.ff_lin_set_to_poly}, gives us the four functors in the quadruple.
It remains to provide the following three natural isomorphisms:\index{isomorphism!natural}
\[
\poly(I,q)\iso\smset(I,q(\0)),\qquad
\poly(q,I)\iso\smset(q(\1),I),
\]
\[
\text{and}\quad\poly(I\yon,q)\iso\smset(I,q(\1)).
\]
All three come from our formula \eqref{eqn.main_formula} for computing hom-sets in $\poly$; we leave the details to the reader in \cref{exc.adjoint_quadruple}.
\end{proof}

\begin{exercise}\label{exc.adjoint_quadruple}
Here we prove the remainder of \cref{thm.adjoint_quadruple} using \eqref{eqn.main_formula}:
\begin{enumerate}\index{isomorphism!natural}
	\item Exhibit a natural isomorphism $\poly(I,q)\iso\smset(I,q(\0))$.
	\item \label{exc.adjoint_quadruple.pos_const} Exhibit a natural isomorphism $\poly(q,I)\iso\smset(q(\1),I)$.
	\item \label{exc.adjoint_quadruple.linear_pos} Exhibit a natural isomorphism $\poly(I\yon,q)\iso\smset(I,q(\1))$.
\qedhere
\end{enumerate}
\begin{solution}
\begin{enumerate}
    \item By \eqref{eqn.main_formula}, we have the natural isomorphism
    \begin{align*}
      \poly(I,q)
        &\iso
      \prod_{i\in I}\sum_{j\in q(\1)}\0^{q[j]} \\
        &\iso
      \prod_{i\in I}q(\0) \\
        &\iso
      \smset(I,q(\0)).
    \end{align*}
    \item By \eqref{eqn.main_formula}, we have the natural isomorphism
    \begin{align*}
      \poly(q,I)
        &\iso
      \prod_{j\in q(\1)}\sum_{i\in I}q[j]^\0 \\
        &\iso
      \prod_{j\in q(\1)}\sum_{i\in I}\1 \\
        &\iso
      \prod_{j\in q(\1)}I \\
        &\iso
      \smset(q(\1),I).
    \end{align*}
    \item By \eqref{eqn.main_formula}, we have the natural isomorphism
    \begin{align*}
        \poly(I\yon,q)
          &\iso
        \prod_{i\in I}\sum_{j\in q(\1)}\1^{q[j]} \\
          &\iso
        \prod_{i\in I}q(\1) \\
          &\iso
        \smset(I,q(\1)).
    \end{align*}
\end{enumerate}
\end{solution}
\end{exercise}

\begin{exercise}\label{exc.positions_maps_yon}
Show that for any polynomial $q$, its set $q(\1)$ of positions is in bijection with the set of functions $\yon\to q$.
\begin{solution}
Given $q\in\poly$, we wish to show that $q(\1)$ is in bijection with the set of functions $\yon\to q$.
In fact, this follows directly from the Yoneda lemma, but we can also invoke the isomorphism from \cref{exc.adjoint_quadruple} \cref{exc.adjoint_quadruple.linear_pos} with $I\coloneqq\1$ to observe that
\[
    \poly(\yon,q)\iso\smset(\1,q(\1))\iso q(\1).
\]
\end{solution}
\end{exercise}\index{Yoneda lemma}

In \cref{thm.adjoint_quadruple} we see that $q\mapsto q(\0)$ and $q\mapsto q(\1)$ have left adjoints.
It turns out this is true more generally for any set $A$ in place of $\0$ and $\1$, as we will show in \cref{cor.substituting_adj}.
However, the fact that $q\mapsto q(\1)$ is itself the left adjoint of the left adjoint of $q\mapsto q(\0)$---and hence that we have the \emph{quadruple} of adjunctions in \eqref{eqn.adjoints_galore}---is special to $A=\0,\1$.

We also have a copower-hom-power two-variable adjunction between $\poly,\smset,$ and $\poly$.

\index{adjunction!two-variable}

\begin{proposition}\label{prop.two_var_adj}
There is a two-variable adjunction between $\poly$, $\smset$, and $\poly$:
\begin{equation}\label{eqn.two_var_adj}
\poly(Ip,q) \iso \smset(I,\poly(p,q)) \iso \poly(p,q^I)
\end{equation}
for $I\in\smset$ and $p,q\in\poly$.
\end{proposition}
\begin{proof}
Since $Ip$ is the $I$-fold coproduct of $p$ and $q^I$ is the $I$-fold product of $q$, the universal properties of coproducts and products give natural isomorphisms\index{isomorphism!natural}
\[\poly(Ip,q)\iso\prod_{i\in I}\poly(p,q)\iso\poly(p,q^I).\]
The middle set is then naturally isomorphic to $\smset(A,\poly(p,q))$.
\end{proof}\index{coproduct!indexed}

Taking $p\coloneqq\yon^A$ in \eqref{eqn.two_var_adj} yields the following via the Yoneda lemma.

\begin{corollary}\label{cor.substituting_adj}
For any set $A$ there is an adjunction
\[
\adj{\smset}{I\yon^A}{q(A)}{\poly}
\]
where the functors are labeled by where they send $I\in\smset$ and $q\in\poly$.
\end{corollary}\index{Yoneda lemma}

\begin{exercise}
Prove \cref{cor.substituting_adj} from \cref{prop.two_var_adj}.
\begin{solution}
To prove \cref{cor.substituting_adj}, it suffices to exhibit a natural isomorphism\index{isomorphism!natural}
\[
    \poly(I\yon^A,q)\iso\smset(I,q(A)).
\]
Taking $p\coloneqq\yon^A$ in \eqref{eqn.two_var_adj} yields the natural isomorphism
\[
    \poly(I\yon^A,q)\iso\smset(I, \poly(\yon^A,q)).
\]
By the Yoneda lemma, $\poly(\yon^A,q)$ is naturally isomorphic to $q(A)$, yielding the desired result.
\end{solution}
\end{exercise}

\index{Yoneda lemma}
\index{adjunction!between $\smset\op$ and $\poly$}\index{global sections}

\begin{proposition}\label{prop.yoneda_left_adjoint}
The Yoneda embedding $A\mapsto \yon^A$ has a left adjoint
\[
\adjr{\smset\op}{\yon^-}{\Gamma}{\poly}
\]
where $\Gamma(p) \coloneqq \poly(p, \yon) \iso \prod_{i\in p(\1)}p[i]$, as in \eqref{eqn.gamma_prod} and \eqref{eqn.gamma_def}.
That is, there is a natural isomorphism\index{isomorphism!natural}
\begin{equation} \label{eqn.yoneda_left_adjoint}
    \poly(p, \yon^A) \iso \smset(A, \Gamma(p)).
\end{equation}
\end{proposition}
\begin{proof}
By \eqref{eqn.main_formula}, we have the natural isomorphism
\[
    \poly(p, \yon^A) \iso \prod_{i \in p(\1)} p[i]^A \iso \left(\prod_{i\in p(\1)}p[i]\right)^A,
\]
which in turn is naturally isomorphic to $\smset(A, \Gamma(p))$ by \eqref{eqn.gamma_prod}.
\end{proof}

\begin{exercise}\index{global sections}
  Give an alternate proof for \cref{prop.yoneda_left_adjoint} using \cref{prop.two_var_adj}.
  \begin{solution}
    Taking $I\coloneqq A$ and $q\coloneqq\yon$ in the second isomorphism in \eqref{eqn.two_var_adj} yields the natural isomorphism\index{isomorphism!natural}
    \[
    \smset(A,\poly(p,\yon)) \iso \poly(p,\yon^A).
    \]
    As $\Gamma(p)=\poly(p,\yon)$, this yields the desired result.
  \end{solution}
\end{exercise}

% Previously an exercise; but already basically done
% Show that $\Gamma(p)\cong\ihom{p,\yon}(\1)$ where $\ihom{-,-}$ is as in \cref{prop.parallel_closure}.

\index{monomial!principal}

\begin{corollary}[Principal monomial]\label{cor.principal_monomial}
There is an adjunction
\[
    \adj{\poly}{{(p(\1),\,\Gamma(p))}}{I\yon^A}{\smset\times\smset\op}
\]
where the functors are labeled by where they send $p\in\poly$ and $(I,A)\in\smset\times\smset\op$.
That is, there is a natural isomorphism\index{isomorphism!natural}
\begin{equation} \label{eqn.principal_monomial}
    \poly(p,I\yon^A) \iso \smset(p(1),I) \times \smset(A,\Gamma(p)).
\end{equation}
\end{corollary}
\begin{proof}
By the universal property of the product of $I$ and $\yon^A$, we have a natural isomorphism
\[
    \poly(p, I\yon^A) \iso \poly(p, I) \times \poly(p, \yon^A).
\]
Then the desired natural isomorphism follows from \cref{exc.adjoint_quadruple} \cref{exc.adjoint_quadruple.pos_const} and \eqref{eqn.yoneda_left_adjoint}.
\end{proof}

\begin{exercise}\index{global sections}
Use \eqref{eqn.principal_monomial} together with \eqref{eqn.dir_dual} and \eqref{eqn.poly_closure_brackets} to find an alternate proof for \cref{prop.situations2}, i.e.\ that there is an isomorphism
\[
    \Gamma(p\otimes q) \iso \smset\big(q(\1),\Gamma(p)\big) \times \smset\big(p(\1),\Gamma(q)\big).
\]
for any $p,q\in\poly$.
\begin{solution}
% The universal properties of the adjunctions
% \[
%       \adj[40pt]{\poly}{{(-(\1),\Gamma(-))}}{-\yon^-}{\smset\times\smset\op}
%       \qqand
%       \adj{\poly}{-\otimes q}{{\ihom{q,-}}}{\poly}
% \]
% from \cref{cor.principal_monomial,prop.parallel_closure} together with the isomorphism $\ihom{p,\yon} \iso \Gamma(p)\yon^{p(\1)}$ from \eqref{eqn.dir_dual} give
We have the following chain of natural isomorphisms involving global sections: \index{global sections}\index{isomorphism!natural}
\begin{align*}
	\Gamma(p \otimes q) &=
	\poly(p \otimes q,\yon)
	\tag*{\eqref{eqn.gamma_def}} \\
	&\iso
	\poly(p, \ihom{q,\yon})
	\tag*{\eqref{eqn.poly_closure_brackets}} \\
	&\iso
	\poly(p, \Gamma(q)\yon^{q(\1)})
	\tag*{\eqref{eqn.dir_dual}} \\
% 	&\iso
% 	(\smset\times\smset\op)\Big(\big((p(\1),\Gamma(p)\big),\big(\Gamma(q),q(\1)\big)\Big) \\
	&\iso
	\smset\big(p(\1),\Gamma(q)\big) \times \smset\big(q(\1),\Gamma(p)\big).
	\tag*{\eqref{eqn.principal_monomial}}
\qedhere
\end{align*}
\end{solution}
\end{exercise}

%-------- Section --------%
\section{Epi-mono factorization of lenses}

\index{factorization system!epi-mono}

We begin by characterizing monomorphisms and epimorphisms in $\poly$.

\index{lens!monomorphism}

\begin{proposition}\label{prop.monics_in_poly}
Let $f \colon p \to q$ be a lens in $\poly$. It is a monomorphism if and only if the on-positions function $f_\1 \colon p(\1) \to q(\1)$ is a monomorphism in $\smset$ and, for each $i \in p(\1)$, the on-directions function $f^\sharp_i \colon q[f_\1i]\to p[i]$ is an epimorphism in $\smset$.
\end{proposition}
\begin{proof}
To prove the forward direction, suppose that $f$ is a monomorphism.
Since $p\mapsto p(\1)$ is a right adjoint (\cref{thm.adjoint_quadruple}), it preserves monomorphisms, so the on-positions function $f_\1$ is also a monomorphism.

We now need to show that for any $i\in p(\1)$, the on-directions function $f^\sharp_i \colon q[f_\1i] \to p[i]$ is an epimorphism.
Suppose we are given a set $A$ and a pair of functions $g^\sharp,h^\sharp\colon p[i]\tto A$ with $f^\sharp_i \then g^\sharp = f^\sharp_i \then h^\sharp$.
Then there exist lenses $g,h \colon \yon^A \tto p$ whose on-positions functions both pick out $i$ and whose on-directions functions are $g^\sharp$ and $h^\sharp$, so that $g \then f = h \then f$.
As $f$ is a monomorphism, $g = h$; in particular, their on-directions functions $g^\sharp$ and $h^\sharp$ are equal, as desired.

Conversely, suppose that $f_\1$ is a monomorphism and that, for each $i\in p(\1)$, the function $f^\sharp_i$ is an epimorphism.
Let $r$ be a polynomial and $g,h\colon r\tto p$ be two lenses such that $g \then f = h \then f$.
Then $g_\1 \then f_\1 = h_\1 \then f_\1$, which implies $g_\1 = h_\1$.
We also have that $f^\sharp_{g_\1k} \then g^\sharp_k = f^\sharp_{g_\1k} \then h^\sharp_k$ for any $k \in r(\1)$. But $f^\sharp_{g_\1k}$ is an epimorphism, so in fact $g^\sharp_k = h^\sharp_k$, as desired.
\end{proof}

\begin{example}\label{ex.clock_in_N}\index{clock}
Choose a finite nonempty set $\ord{k}$ for $1\leq k\in\nn$, e.g.\ $\ord{k}\coloneqq\1\2$. There is a monomorphism
\[
f \colon\ord{k}\yon^{\ord{k}}\to\nn\yon^\nn
\]
defined as follows.
On positions, we have $f_\1i\coloneqq i$ for all $i \in \ord{k}$. On directions, for any $i \in \ord{k}$, we have $f^\sharp_i(n)\coloneqq n \bmod k$ for all $n \in \nn$.
\end{example}

\begin{exercise}
In \cref{ex.clock_in_N}, we gave a lens $\1\2\yon^{\1\2}\to\nn\yon^\nn$. This allows us to turn any dynamical system with state-set $\nn$ into a dynamical system with 12 states, while keeping the same interface---say, $p$---by composing with this lens.

Explain how the behavior of the new system $\1\2\yon^{\1\2}\to p$ relates to the behavior of the old system $\nn\yon^\nn\to p$.
\begin{solution}
We are given a monomorphism $f \colon \1\2\yon^{\1\2} \to \nn\yon^\nn$ from \cref{ex.clock_in_N}.
Let $g \colon \nn\yon^\nn \to p$ be a dynamical system with return function $g_\1 \colon \nn \to p(\1)$ and update functions $g^\sharp_n \colon p[g_\1(n)] \to \nn$ for each state $n \in \nn$.
Then the new composite dynamical system $h \coloneqq f \then g$ has a return function $h_\1 \colon \1\2 \to p(\1)$ which sends each state $i \in \1\2$ to the position $h_\1i = g_\1f_\1i = g_\1i$, the same position that the original system returns in the state $i \in \nn$.
Meanwhile, the update function for each state $i \in \1\2$ is a function $h^\sharp_i \colon p[g_\1i] \to \1\2$ which, given a direction $a \in p[g_\1i]$, updates the state from $i$ to $h^\sharp_ia = f^\sharp_{g_\1i}(g^\sharp_ia) = g^\sharp_ia \mod 12$, which is where the original system would have taken the same state, but reduced modulo 12.
In other words, the new system behaves like the old system but with only the states in $\1\2 \ss \nn$ retained, and on any input that would have caused the old system to move to a state outside of $\1\2$, the new system moves to the equivalent state (modulo 12) within $\1\2$ instead.
\end{solution}
\end{exercise}

\index{lens!epimorphism}

\begin{proposition}\label{prop.epis_in_poly}
Let $f \colon p \to q$ be a lens in $\poly$. It is an epimorphism if and only if the function $f_\1 \colon p(\1) \to q(\1)$ is an epimorphism in $\smset$ and, for each $j\in q(\1)$, the induced function
\[
    f^\flat_j \colon q[j] \to \prod_{\substack{i\in p(\1), \\ f_\1i=j}} p[i]
\]
from \eqref{eqn.useful_misc472} is a monomorphism.
\end{proposition}
\begin{proof}
To prove the forward direction, suppose that $f$ is an epimorphism. Since $p \mapsto p(\1)$ is a left adjoint (\cref{thm.adjoint_quadruple}), it preserves epimorphisms, so the on-positions function $f_\1$ is also a epimorphism.

We now need to show that for any $j\in q(\1)$, the induced function $f^\flat_j$ is a monomorphism.
Suppose we are given a set $A$ and a pair of functions $g',h'\colon A\tto q[j]$ with $g' \then f^\flat_j = h' \then f^\flat_j$.
Then we can construct lenses $g,h\colon q\tto\{0\}\yon^A+\{1\}$ that send $j\mapsto0$ and every other position to $1$ on positions; on directions, $g^\sharp_j\coloneqq g'$ and $h^\sharp_j\coloneqq h'$, while every other on-direction function is uniquely determined.
Then $f \then g = f \then h$.
As $f$ is an epimorphism, $g = h$; in particular, their on-directions functions are equal, so $g' = h'$, as desired.

Conversely, suppose that $f_\1$ is an epimorphism and that, for each $j\in q(\1)$, the function $f^\flat_j$ is a monomorphism.
Let $r$ be a polynomial and $g,h\colon q\tto r$ be two lenses such that $f \then g = f \then h$.
Then $f_\1 \then g_\1 = f_\1 \then h_\1$, which implies $g_\1=h_\1$.
We also have that $g^\sharp_{f_\1i} \then f^\sharp_i = h^\sharp_{f_\1i} \then f^\sharp_i$ for any $i\in p(\1)$.
It follows that, for any $j \in q(\1)$, the two composites
\[
\begin{tikzcd}
	r[g_\1j] \ar[r, shift left, "g^\sharp_j"] \ar[r, shift right, "h^\sharp_j"'] & q[j] \ar[r, "f^\flat_j"] & \displaystyle\prod_{\substack{i\in p(\1), \\ f_\1i=j}} p[i]
\end{tikzcd}
\]
are equal, which implies that $g^\sharp_j=h^\sharp_j$ as desired.
\end{proof}

% Insert exercise exploring the difference between the epi proposition and the one about monos

\begin{exercise}
Show that the only way for a lens $p\to\yon$ to \emph{not} be an epimorphism is when $p=0$.
\begin{solution}
Given $p \in \poly$ and a lens $f \colon p \to \yon$, we will use \cref{prop.epis_in_poly} to show that either $f$ is an epimorphism or $p = \0$.
First, note that $f_\1 \colon p(\1) \to \1$ must be an epimorphism unless $p(\1) \iso \0$, in which case $p = \0$.
Next, note that the induced function
\[
    f^\flat \colon \1 \to \prod_{i\in p(\1)} p[i]
\]
from \eqref{eqn.useful_misc472} must be a monomorphism.
So it follows from \cref{prop.epis_in_poly} that either $f$ is an epimorphism or $p = \0$.
\end{solution}
\end{exercise}

\begin{exercise}
Let $A$ and $B$ be sets and $AB$ their product. Find an epimorphism $\yon^A+\yon^B\surj\yon^{AB}$.
\begin{solution}
Given sets $A$ and $B$, by \cref{prop.epis_in_poly}, a lens $f \colon \yon^A + \yon^B \to \yon^{AB}$ is an epimorphism if its on-positions function $f_\1 \colon \2 \to \1$ is an epimorphism (which must be true) and if the induced function
\[
    f^\flat \colon AB \to \prod_{i \in \2} (\yon^A + \yon^B)[i] \iso AB
\]
is a monomorphism.
If we take the on-directions functions $AB \to A$ and $AB \to B$ of $f$ to be the canonical projections, then the induced function $f^\flat \colon AB \to AB$ would be the identity, which is indeed a monomorphism.
So $f$ would be an epimorphism.
\end{solution}
\end{exercise}

\begin{exercise}
Suppose a lens $f\colon p\to q$ is both a monomorphism and an epimorphism; it is then an isomorphism? (That is, is $\poly$ \emph{balanced}?)

Hint: You may use the following facts.
\begin{enumerate}
    \item A function that is both a monomorphism and an epimorphism in $\smset$ is an isomorphism.
    \item A lens is an isomorphism if and only if the on-positions function is an isomorphism and every on-directions function is an isomorphism. \qedhere
\end{enumerate}
\begin{solution}
Let $f \colon p \to q$ be a lens in $\poly$ that is both a monomorphism and an epimorphism.
We claim that $f$ is an isomorphism.
By \cref{prop.monics_in_poly} and \cref{prop.epis_in_poly}, the on-positions function $f_\1 \colon p(\1) \to q(\1)$ is both a monomorphism and an epimorphism, so it is an isomorphism.
Meanwhile, \cref{prop.epis_in_poly} says that, for each $j \in q(\1)$, the induced function
\[
    f^\flat_j \colon q[j] \to \prod_{\substack{i \in p(\1), \\ f_\1i = j}} p[i]
\]
is a monomorphism.
As $f_\1$ is an isomorphism, it follows that for each $i \in p(\1)$, the function
\[
    f^\flat_{f_\1i} \colon q[f_\1i] \to p[i]
\]
is a monomorphism.
But this is just the on-directions function $f^\sharp_i$ of $f$.
From \cref{prop.monics_in_poly}, we also know that $f^\sharp_i$ is an epimorphism.
It follows that every on-directions function of $f$ is an isomorphism.
Hence $f$ itself is an isomorphism.
\end{solution}
\end{exercise}

We are often interested in whether epimorphisms and monomorphisms form what is called a \emph{factorization system} in a given category, which we define below.

\begin{definition}[Factorization system] \label{def.factor}
Given a category $\cat{C}$ and two classes of morphisms $E$ and $M$ in $\cat{C}$, we say that $(E, M)$ is a \emph{factorization system} of $\cat{C}$ if:
\begin{enumerate}
    \item every morphism $f$ in $\cat{C}$ factors uniquely (up to unique isomorphism) as a morphism $e \in E$ composed with a morphism $m\in M$, so that $f = e\then m$;
    \item $E$ and $M$ each contain every isomorphism; and
    \item $E$ and $M$ are each closed under composition.
\end{enumerate}
If $E$ is the class of epimorphisms and $M$ is the class of monomorphisms (in which case conditions 2 and 3 are automatically satisfied), we say that $\cat{C}$ has \emph{epi-mono factorization}.
\end{definition}

\begin{example}[Epi-mono factorization in $\smset$] \label{ex.epi_mono_set}
The category $\smset$ has epi-mono factorization: a function $f\colon X\to Y$ can be uniquely factored into an epimorphism (surjection) $e$ followed by a monomorphism (injection) $i$, as follows.
The epimorphism $e\colon X\to f(X)$ is given by restricting the codomain of $f$ to its image (also known as \emph{corestricting} $f$), so $e$ sends $x\mapsto f(x)$ for all $x\in X$.
The monomorphism $i\colon f(X)\to Y$ is then given by including the image into the codomain, so $i$ sends $y\mapsto y$ for all $y\in f(X)\ss Y$.
\end{example}

\begin{proposition}
$\poly$ has epi-mono factorization.
\end{proposition}
\begin{proof}
Take an arbitrary lens $\varphi\colon p\to q$.
It suffices to show that there exists a unique polynomial $r$ equipped with an epimorphism $\epsilon\colon p\to r$ and a monomorphism $\mu\colon r\to q$ such that $\varphi=\epsilon\then\mu$.

On positions, we must have $\varphi_\1 = \epsilon_\1\then\mu_\1$, with $\mu_\1$ a monomorphism and $\epsilon_\1$ an epimorphism per \cref{prop.monics_in_poly,prop.epis_in_poly}.
By \cref{ex.epi_mono_set}, since $\smset$ has epi-mono factorization, such $r(\1), \epsilon_\1,$ and $\mu_\1$ uniquely exist.
In particular, we must have that $r(\1)\iso \varphi_\1(p(\1))$, that $\epsilon_\1\colon p(\1)\to\varphi_\1(p(\1))$ is the corestriction of $\varphi_\1$ sending $i\mapsto\varphi_\1(i)$ for each $p$-position $i$, and that $\mu_\1\colon\varphi_\1(p(\1))\to q(\1)$ is the inclusion sending $j\mapsto j$ for each $r$-position $j$.

Then on directions, for any $i\in p(\1)$, we must have that
\[
\begin{tikzcd}
    q[\varphi_\1(i)] \ar[r, "\mu^\sharp_{\varphi_\1(i)}"] \ar[dr, "\varphi^\sharp_i"'] & r[\varphi_\1(i)] \ar[d, "\epsilon^\sharp_i"] \\
    & p[i]
\end{tikzcd}
\]
commutes---or, equivalently, for every $j\in r(\1)\iso\varphi_\1(p(\1))$,
\[
\begin{tikzcd}
    q[j] \ar[r, "\mu^\sharp_j"] \ar[dr, "\varphi^\flat_j"'] & r[j] \ar[d, "\epsilon^\flat_j"] \\
    & \prod\limits_{\substack{i\in p(\1), \\ \varphi_\1(i)=j}} p[i]
\end{tikzcd}
\]
commutes (here $\varphi^\flat_j$ and $\epsilon^\flat_j$ are the induced functions from \eqref{eqn.useful_misc472}), with $\mu_j^\sharp$ an epimorphism and $\epsilon^\flat_j$ a monomorphism per
\cref{prop.monics_in_poly,prop.epis_in_poly}.
So again since $\smset$ has epi-mono factorization, such $r[j], \mu^\sharp_j,$ and $\epsilon^\flat_j$ uniquely exist.
Hence such $p \To{\epsilon} r \To{\mu} q$ uniquely exists overall.
\end{proof}

\index{factorization system!epi-mono}

%-------- Section --------%
\section{Cartesian closure}

\index{polynomial functor!cartesian closed structure}
\index{polynomial functor!exponentials of polynomials}
\index{closed monoidal structure}
\index{exponential|(}

We have already seen in \cref{sec.closure} the closure operation $\ihom{-,-}$ for one monoidal structure on $\poly$, namely $(\yon,\otimes)$.
But this is not the only closed monoidal structure on $\poly$: in fact, we will show that $\poly$ is cartesian closed as well.

For any two polynomials $q,r$, define $r^q\in\poly$ by the formula
\begin{equation}\label{eqn.exponential}
  r^q\coloneqq\prod_{j\in q(\1)}r\circ(\yon+q[j])
\end{equation}
where $\circ$ denotes composition.

Before proving that this really is an exponential in $\poly$, which we do in \cref{thm.poly_cart_closed}, we first get some practice with it.

\begin{example}
Let $A$ be a set. We've been writing the polynomial $A\yon^\0$ simply as $A$, so it better be true that there is an isomorphism
\[
    \yon^A \iso \yon^{A\yon^\0}
\]
in order for the notation to be consistent.
Luckily, this is true.
By \eqref{eqn.exponential}, we have
\[
    \yon^{A\yon^\0} = \prod_{a\in A} \yon \circ (\yon+\0) \iso \yon^A
\]
\end{example}

\begin{exercise}
Compute the following exponentials in $\poly$ using \eqref{eqn.exponential}:
\begin{enumerate}
	\item $p^\0$ for an arbitrary $p\in\poly$.
	\item $p^\1$ for an arbitrary $p\in\poly$.
	\item $\1^p$ for an arbitrary $p\in\poly$.
	\item $A^p$ for an arbitrary $p\in\poly$ and $A\in\smset$.
	\item $\yon^\yon$.
	\item $\yon^{\4\yon}$.
	\item $(\yon^A)^{\yon^B}$ for arbitrary sets $A,B\in\smset$.
\qedhere
\end{enumerate}
\begin{solution}
We use \eqref{eqn.exponential} to compute various exponentials.
Here $p \in \poly$ and $A, B \in \smset$.
\begin{enumerate}
    \item We have that $p^\0$ is an empty product, so $p^\0 \iso \1$ as expected.
	\item We have that $p^\1 \iso p \circ (\yon + \0) \iso p$, as expected.
	\item We have that $\1^p \iso \prod_{i \in p(\1)} \1 \circ (\yon + p[i]) \iso \1$, as expected.
	\item We have that $A^p \iso \prod_{i \in p(\1)} A \circ (\yon + p[i]) \iso A^{p(\1)}$.
	\item We have that $\yon^\yon \iso \yon \circ (\yon + \1) \iso \yon + \1$.
	\item We have that $\yon^{\4\yon} \iso \prod_{j \in \4} \yon \circ (\yon + \1) \iso (\yon + \1)^\4 \iso \yon^\4 + \4\yon^\3 + \6\yon^\2 + \4\yon + \1$.
	\item We have that $(\yon^A)^{\yon^B} \iso (\yon^A) \circ (\yon + B) \iso (\yon + B)^A \iso \sum_{f \colon A \to \2} B^{f\inv(1)} \yon^{f\inv(2)}$.
\end{enumerate}
\end{solution}
\end{exercise}

% \begin{exercise} % Immediate from previous exercise
% Using \eqref{eqn.exponential}, show that the functor $\smset\to\poly$ that sends each set $A$ to the constant polynomial $A$ preserves exponentials.
% That is, given sets $A, B \in \smset$, the set $B^A$ as a constant polynomial coincides with the exponential in $\poly$ that is the constant polynomial $B$ raised to the constant polynomial $A$.
% \begin{solution}
% By \eqref{eqn.exponential}, the exponential that is the constant polynomial $B$ raised to the constant polynomial $A$ can be written as
% \[
%     \prod_{a \in A} B \circ (\yon + A[a]) \iso \prod_{a \in A} B \iso B^A.
% \]
% \end{solution}
% \end{exercise}

\begin{theorem}\label{thm.poly_cart_closed}
The category $\poly$ is cartesian closed. That is, we have a natural isomorphism\index{isomorphism!natural}
\[
    \poly(p,r^q) \iso \poly(p\times q,r),
\]
where $r^q$ is the polynomial defined in \eqref{eqn.exponential}.
\end{theorem}
\begin{proof}
We have the following chain of natural isomorphisms:
\begin{align*}
	\poly(p, r^q) &\iso
	\poly\Big(p, \prod_{j \in q(\1)} r \circ (\yon+q[j])\Big)
	\tag*{\eqref{eqn.exponential}} \\
% 	&\iso
% 	\prod_{j\in q(\1)}\poly(p,r\circ(\yon+q[j]))
% 	\tag{Universal property of products} \\
	&\iso
	\prod_{i\in p(\1)}\prod_{j\in q(\1)}\poly\big(\yon^{p[i]},r\circ(\yon+q[j])\big)
	\tag{Universal property of (co)products} \\
	&\iso
	\prod_{i\in p(\1)}\prod_{j\in q(\1)}r\circ(p[i]+q[j])
	\tag{Yoneda lemma} \\
	&\iso
	\prod_{i\in p(\1)}\prod_{j\in q(\1)}\sum_{k\in r(\1)}(p[i]+q[j])^{r[k]}
	\\
	&\iso
	\prod_{(i,\,j) \in (p \times q)(\1)} \; \sum_{k\in r(\1)}(p \times q)[(i, j)]^{r[k]}
	\tag*{\eqref{eqn.poly_times}} \\
	&\iso
	\poly(p \times q,r).
	\tag*{\eqref{eqn.main_formula}}
\end{align*}
\end{proof}

\begin{exercise}
Use \cref{thm.poly_cart_closed} to show that for any polynomials $p,q$, there is a canonical evaluation lens
\begin{equation*}%\label{eqn.eval_times}
	\text{eval}\colon p^q \times q \to p.
\end{equation*}
\begin{solution}
By \cref{thm.poly_cart_closed}, there is a natural isomorphism\index{isomorphism!natural}
\[
    \poly(p^q, p^q) \iso \poly(p^q \times q, p).
\]
Under this isomorphism, there exists a lens $\text{eval} \colon p^q \times q \to p$ corresponding to the identity lens on $p^q$.
The lens $\text{eval}$ is the canonical evaluation lens.
\end{solution}
\end{exercise}

\index{exponential|)}\index{exponential|seealso{closed monoidal structure}}\index{Cartesian closed category}

%-------- Section --------%
\section[Limits and colimits of polynomials]{Limits and colimits of polynomials%
  \sectionmark{Limits \& colimits of polynomials}}
\sectionmark{Limits \& colimits of polynomials}

\index{polynomial functor!limits of polynomials|(}

We have already seen that $\poly$ has all coproducts (\cref{prop.poly_coprods}) and products (\cref{prop.poly_prods}).
We will now see that $\poly$ has all small limits and colimits.

\begin{theorem}\label{thm.poly_limits}\index{limit!in $\poly$}
The category $\poly$ has all small limits.
\end{theorem}
\begin{proof}
A category has all small limits if and only if it has products and equalizers, so by \cref{prop.poly_prods}, it suffices to show that $\poly$ has equalizers.

We claim that equalizers in $\poly$ are simply equalizers on positions and coequalizers on directions.
More precisely, let $f,g \colon p \tto q$ be two lenses.
We construct the equalizer $p'$ of $f$ and $g$ as follows.\footnote{If we're being precise, a ``(co)equalizer'' is an object equipped with a morphism, but we will use the term to refer to either just the object or just the morphism when the context is clear.}
We define its position-set $p'(\1)$ to be the equalizer of $f_\1,g_\1 \colon p(\1) \tto q(\1)$ in $\smset$; that is,
\[
    p'(\1) \coloneqq \{i \in p(\1) \mid f_\1i = g_\1i\}.
\]
Then for each $i \in p'(\1)$, we can define the direction-set $p'[i]$ to be the coequalizer of $f^\sharp_i, g^\sharp_i \colon q[f_\1i] \tto p[i]$.
In this way, we obtain a polynomial $p'$ that comes equipped with a lens $e \colon p' \to p$.
One can check that $p'$ together with $e$ satisfies the universal property of the equalizer of $f$ and $g$; see \cref{exc.poly_limits}.
\end{proof}

\begin{exercise}\label{exc.poly_limits}
Complete the proof of \cref{thm.poly_limits} as follows:
\begin{enumerate}
	\item We said that $p'$ comes equipped with a lens $e \colon p' \to p$; what is it?
	\item Show that $e \then f = e \then g$.
	\item Show that $e$ is the equalizer of the pair $f,g$.
\qedhere
\end{enumerate}
\begin{solution}
\begin{enumerate}
    \item The lens $e \colon p' \to p$ can be characterized as follows.
    The on-positions function $e_\1 \colon p'(\1) \to p(\1)$ is the equalizer of $f_\1, g_\1 \colon p(\1) \tto q(\1)$ in $\smset$.
    In particular, $e_\1$ is the canonical inclusion that sends each element of $p'(\1)$ to the same element in $p(\1)$.
    Then for each $i \in p'(\1)$, the on-directions function $e^\sharp_i \colon p[i] \to p'[i]$ is the coequalizer of $f^\sharp_i, g^\sharp_i \colon q[f_\1i] \tto p[i]$ in $\smset$.

    \item To show that $e \then f = e \then g$, it suffices to show that both sides are equal on positions and on directions.
    On positions, $e_\1$ is defined to be the equalizer of $f_\1$ and $g_\1$, so $e_\1 \then f_\1 = e_\1 \then g_\1$.
    Then for each $i \in p'(\1)$, the on-directions function $e^\sharp_i$ is defined to be the coequalizer of $f^\sharp_i$ and $g^\sharp_i$, so $f^\sharp_i \then e^\sharp_i = g^\sharp_i \then e^\sharp_i$.

    \item To show that $e$ is the equalizer of $f$ and $g$, it suffices to show that for any $r \in \poly$ and lens $a \colon r \to p$ satisfying $a \then f = a \then g$, there exists a unique lens $h \colon r \to p'$ for which $a = h \then e$, so that the following diagram commutes.
    \begin{equation*} %\label{eqn.eq_univ_prop}
    \begin{tikzcd}
        p' \ar[r, "e"] & p \ar[r, "f", shift left] \ar[r, "g"', shift right] & q \\
        r \ar[u, "h", dashed] \ar[ur, "a"']
    \end{tikzcd}
    \end{equation*}
    In order for $a = h \then e$ to hold, we must have $a_\1 = h_\1 \then e_\1$ on positions.
    But we have that $a_\1 \then f_\1 = a_\1 \then g_\1$, so by the universal property of $p'(\1)$ and the map $e_\1$ as the equalizer of $f_\1$ and $g_\1$ in $\smset$, there exists a unique $h_\1$ for which $a_\1 = h_\1 \then e_\1$.
    Hence $h$ is uniquely characterized on positions.
    In particular, it must send each $k \in r(\1)$ to $a_\1(k) \in p'(\1)$.

    Then for $a = h \then e$ to hold on directions, we must have that $a^\sharp_k = e^\sharp_{a_\1(k)} \then h^\sharp_k$ for each $k \in r(\1)$.
    But we have that $f^\sharp_{a_\1(k)} \then a^\sharp_{a_\1(k)} = g^\sharp_{a_\1(k)} \then a^\sharp_{a_\1(k)}$, so by the universal property of $p'[a_\1(k)]$ and the map $e^\sharp_{a_\1(k)}$ as the coequalizer of $f^\sharp_{a_\1(k)}$ and $g^\sharp_{a_\1(k)}$ in $\smset$, there exists a unique $h^\sharp_k$ for which $a^\sharp_k = e^\sharp_{a_\1(k)} \then h^\sharp_k$, so that the diagram below commutes.
    \begin{equation*} %\label{eqn.eq_univ_prop_dir}
    \begin{tikzcd}[sep=large]
        p'[a_\1(k)] \ar[d, "h^\sharp_k"', dashed] & p[a_\1(k)] \ar[l, "e^\sharp_{a_\1(k)}"'] \ar[dl, "a^\sharp_k"] & q[f_\1(a_\1(k))] \ar[l, "f^\sharp_{a_\1(k)}"', shift right] \ar[l, "g^\sharp_{a_\1(k)}", shift left] \\
        r[k]
    \end{tikzcd}
    \end{equation*}
    Hence $h$ is also uniquely characterized on directions, so it is unique overall.
    Moreover, we have shown that we can define $h$ on positions so that $a_\1 = h_\1 \then e_\1$, and that we can define $h$ on directions such that $a^\sharp_k = e^\sharp_{a_\1(k)} \then h^\sharp_k$ for all $k \in r(\1)$.
    It follows that there exists $h$ for which $a = h \then e$.
\end{enumerate}
\end{solution}
\end{exercise}

\begin{example}[Computing general limits in $\poly$] \label{ex.compute_limits}\index{limit!positions and directions}
The proof of \cref{thm.poly_limits} justifies the following mnemonic for limits in $\poly$:
\slogan{The positions of a limit are the limit of the positions. \\ The directions of a limit are the colimit of the directions.}
We can make this precise as follows: the limit of a functor $p_-\colon\cat{J}\to\poly$ is the polynomial whose position-set is
\begin{equation} \label{eqn.lim_pos}
    \left(\lim_{j\in\cat{J}} p_j\right)(\1) \iso \lim_{j\in\cat{J}} p_j(\1),
\end{equation}
equipped with a canonical projection $\pi_j$ to each $p_j(\1)$, and whose direction-set for each position $i$ is
\begin{equation} \label{eqn.lim_dir}
    \left(\lim_{j\in\cat{J}} p_j\right)[i] \iso \colim_{j\in\cat{J}\op} p_j[\pi_j(i)].
\end{equation}
This notation obscures what is occuring on lenses, but in particular, each lens $\varphi\colon p_j\to p_{j'}$ in the diagram $p_-$ induces an on-positions function $\varphi_\1\colon p_j(\1)\to p_{j'}(\1)$ in the diagram whose limit we take in \eqref{eqn.lim_pos} and, for every position $i$ of the limit, an on-directions function $\varphi^\sharp_{\pi_j(i)}\colon p_{j'}[\pi_{j'}(i)]\to p_j[\pi_j(i)]$ in the diagram whose colimit we take in \eqref{eqn.lim_dir}.
(Note that, by the definition of a limit, $\varphi_\1(\pi_j(i)) = \pi_{j'}(i)$.)

We have seen \eqref{eqn.lim_pos} and \eqref{eqn.lim_dir} to be true for products: the position-set of the product is just the product of the original position-sets, while the direction-set at a tuple of the original positions is just the coproduct of the direction-sets at every position in the tuple.
We have also just shown \eqref{eqn.lim_pos} and \eqref{eqn.lim_dir} to be true for equalizers in the proof of \cref{thm.poly_limits}.
It follows from the construction of any limit as an equalizer of products that it is true for arbitrary limits.
\end{example}

\index{polynomial functor!pullback of polynomials}\index{pullback|seealso{polynomial functor, pullback of polynomials}}

\begin{example}[Pullbacks in $\poly$]\label{ex.pullbacks_in_poly}
Given $q,q',r \in \poly$ and lenses $q\To{f} r\From{f'} q'$, the pullback
\[
\begin{tikzcd}
	p\ar[r, "g'"]\ar[d, "g"']&
	q'\ar[d, "f'"]\\
	q\ar[r, "f"']&
	r\ar[ul, phantom, very near end, "\lrcorner"]
\end{tikzcd}
\]
is given as follows.
The position-set of $p$ is the pullback of the position-sets of $q$ and $q'$ over that of $r$ in $\smset$.
Then at each position $(i, i') \in p(\1) \ss q(\1) \times q'(\1)$ with $f_\1i=f'_\1i'$, we take the direction-set $p[(i, i')]$ to be the pushout of the direction-sets $q[i]$ and $q'[i']$ over $r[f_\1i]=r[f_\1'i']$ in $\smset$.
These pullback and pushout squares also give the lenses $g$ and $g'$ on positions and on directions:
\begin{equation}\label{eqn.pullback_poly}
\begin{tikzcd}
	p(\1)\ar[r, "g'_\1"]\ar[d, "g_\1"']&
	q'(\1)\ar[d, "f_\1'"]\\
	q(\1)\ar[r, "f_\1"']&
	r(\1)\ar[ul, phantom, very near end, "\lrcorner"]
\end{tikzcd}
\qqand
\begin{tikzcd}
	p[(i,i')]\ar[from=r, "(g')^\sharp_{(i,\,i')}"']\ar[from=d, "g^\sharp_{(i,\,i')}"]&
	q'[i']\ar[from=d, "(f')^\sharp_{i'}"']\\
	q[i]\ar[from=r, "f^\sharp_i"]&
	r[f_1(i)]\ar[ul, phantom, very near end, "\lrcorner"]
\end{tikzcd}
\end{equation}
\end{example}

\begin{exercise}\index{polynomial functor!pullback of polynomials}
Let $p$ be any polynomial.
\begin{enumerate}
	\item There is a canonical choice of lens $\eta\colon p\to p(\1)$; what is it?
	\item Given an element $i\in p(\1)$, i.e.\ a function (or lens between constant polynomials) $i\colon\1\to p(\1)$, let $p_i$ be the pullback
	\[
	\begin{tikzcd}
	p_i\ar[r, "g"]\ar[d, "f"']&
	p\ar[d, "\eta"]\\
	\1\ar[r, "i"']&
	p(\1)\ar[ul, phantom, very near end, "\lrcorner"]
	\end{tikzcd}
	\]
	What is $p_i$? What are the lenses $f \colon p_i \to \1$ and $g \colon p_i \to p$? \qedhere
\end{enumerate}
\begin{solution}
Here $p \in \poly$.
\begin{enumerate}
    \item The canonical lens $\eta \colon p \to p(\1)$ is the identity $\eta_\1 \colon p(\1) \to p(\1)$ on positions and the empty function on directions.

    \item On positions, we have that $p_i(\1)$ along with $f_\1$ and $g_\1$ form the following pullback square in $\smset$:
    \[
	\begin{tikzcd}
    	p_i(\1) \ar[r, "g_\1"] \ar[d, "f_\1"'] &
    	p(\1) \ar[d, equals] \\
    	\1 \ar[r, "i"'] &
    	p(\1) \ar[ul, phantom, very near end, "\lrcorner"]
	\end{tikzcd}
	\]
	So $p_i(\1) \coloneqq \{(a, i') \in \1 \times p(\1) \mid i = i' \} = \{(1, i)\}$, with $f_\1$ uniquely determined and $g_1$ picking out $i \in p(\1)$.
	Then on directions, we have that $p_i[(1,i)]$ along with $f^\sharp_{(1,\,i)}$ and $g^\sharp_{(1,\,i)}$ form the following pushout square in $\smset$:
	\[
	\begin{tikzcd}
    	p_i[(1,i)] \ar[from=r, "g^\sharp_{(1,\,i)}"'] \ar[from=d, "f^\sharp_{(1,\,i)}"] &
    	p[i] \ar[from=d, "!"'] \\
    	\0 \ar[from=r, "!"] &
    	\0 \ar[ul, phantom, very near end, "\lrcorner"]
    \end{tikzcd}
    \]
    So $p_i[(1,i)] \coloneqq p[i]$, with $f^\sharp_{(1,\,i)}$ uniquely determined and $g^\sharp_{(1,\,i)}$ as the identity.
    It follows that $p_i \coloneqq \{(1,i)\}\yon^{p[i]} \iso \yon^{p[i]}$, where $f \colon p_i \to \1$ is uniquely determined and $g \colon p_i \to p$ picks out $i \in p(\1)$ on positions and is the identity on $p[i]$ on directions.
\end{enumerate}
\end{solution}
\end{exercise}

\begin{exercise}
Let $q\coloneqq \yon^\2+\yon$, $q'\coloneqq\2\yon^\3+\yon^\2$, and $r\coloneq\yon+\1$.
\begin{enumerate}
	\item Choose lenses $f\colon q\to r$ and $f'\colon q'\to r$ and write them down.
	\item Find the pullback of $q\To{f} r\From{f'} q'$.
\qedhere
\end{enumerate}
\begin{solution}\index{polynomial functor!pullback of polynomials}
\begin{enumerate}
    \item There are many possible answers, but one lens $f \colon q \to r$, on positions, sends $1 \in q(\1)$ (corresponding to $\yon^\2$) to $2 \in r(\1)$ (corresponding to $\1$) and $2 \in q(\1)$ (corresponding to $\yon$) to $1 \in r(\1)$ (corresponding to $\yon$).
    Then the on-directions functions $f^\sharp_\1 \colon \0 \to \2$ and $f^\sharp_2 \colon \1 \to \1$ are uniquely determined.
    Another morphism $f' \colon q' \to r$, on positions, sends $1 \in q'(\1)$ (corresponding to one of the $\yon^\3$ terms) to $2 \in r(\1)$ and both $2 \in q'(\1)$ (corresponding to the other $\yon^\3$ term) and $3 \in q'(\1)$ (corresponding to the $\yon^\2$ term) to $1 \in r(\1)$.
    Then the on-directions function $(f')^\sharp_1 \colon \0 \to \3$ is uniquely determined, while we can let $(f')^\sharp_2 \colon \1 \to \3$ pick out $3$ and $(f')^\sharp_3 \colon \1 \to \2$ pick out $1$.

    \item We compute the pullback $p$ along with the lenses $g \colon p \to q$ and $g' \colon p \to q'$ of $q\To{f} r\From{f'} q'$ by following \cref{ex.pullbacks_in_poly}.
    We can compute $p(\1)$ by taking the pullback in $\smset$:
    \[
        p(\1) \coloneqq \{(i, i') \in \2 \times \3 \mid f_\1i = f'_\1(i)\} = \{(1,1), (2,2), (2,3)\}.
    \]
    Moreover, the on-positions functions $g_\1$ and $g'_\1$ send each pair in $p(\1)$ to its left component and its right component, respectively.

    To compute the direction-set at each $p$-position, we must compute a pushout.
    At $(1,1)$, we have $r[f_\1(1)] = r[f'_\1(1)] = r[2] = \0$, so the pushout $p[(1,1)]$ is just the sum $q[1] + q'[1] = \2 + \3 \iso \5$.
    Moreover, the on-directions functions $g^\sharp_{(1,\,1)}$ and $(g')^\sharp_{(1,\,1)}$ are the canonical inclusions $\2 \to \2 + \3$ and $\3 \to \2 + \3$.

    At $(2,2)$, we have $r[f_\1(2)] = r[f'_\1(2)] = r[1] = \1$, with $f^\sharp_2$ picking out $1 \in \1 = q[2]$ and $(f')^\sharp_2$ picking out $3 \in \3 = q'[2]$.
    So the pushout $p[(2,2)]$ is the set $\1 + \3 = \{(1,1), (2,1), (2,2), (2,3)\}$ but with $(1,1)$ identified with $(2,3)$; we can think of it as the set of equivalence classes $p[(2,2)] \iso \{\{(1,1), (2,3)\}, \{(2,1)\}, \{(2,2)\}\} \iso \3$.
    Moreover, the on-directions function $g^\sharp_{(2,\,2)}$ maps $1 \mapsto \{(1,1), (2,3)\}$, while the on-directions function $(g')^\sharp_{(2,\,2)}$ maps $1 \mapsto \{(2,1)\}, 2 \mapsto \{(2,2)\},$ and $3 \mapsto \{(1,1), (2,3)\}$.

    Finally, at $(2,3)$, we have $r[f_\1(2)] = r[f'_\1(3)] = r[1] = \1$, with $f^\sharp_2$ still picking out $1 \in \1 = q[2]$ and $(f')^\sharp_3$ picking out $1 \in \2 = q'[3]$.
    So the pushout $p[(2,3)]$ is the set $\1 + \2 = \{(1,1), (2,1), (2,2)\}$ but with $(1,1)$ identified with $(2,1)$; we can think of it as the set of equivalence classes $p[(2,3)] \iso \{\{(1,1), (2,1)\}, \{(2,2)\}\} \iso \2$.
    Moreover, the on-directions function $g^\sharp_{(2,\,3)}$ maps $1 \mapsto \{(1,1), (2,1)\}$, while the on-directions function $(g')^\sharp_{(2,\,3)}$ maps $1 \mapsto \{(1,1), (2,1)\}$ and $2 \mapsto \{(2,2)\}$.

    It follows that $p \iso \yon^\5 + \yon^\3 + \yon^\2$, with $g$ and $g'$ as described.
\end{enumerate}
\end{solution}
\end{exercise}

\index{polynomial functor!equalizer of polynomials}

\begin{exercise} \label{exc.refl_limits}
An alternative way to prove \cref{thm.poly_limits} would have been to show that the equalizer of two natural transformations between polynomial functors in $\smset^\smset$ is still a polynomial functor---since the full subcategory inclusion $\poly \to \smset^\smset$ reflects these equalizers, it would follow that $\poly$ has equalizers.
But we already know what polynomial the equalizer should be from the proof of \cref{thm.poly_limits}.
So in this exercise, we will show that the equalizer of polynomials we found in $\poly$ is also the equalizer of those same functors in $\smset^\smset$.

Let $f,g\colon p\tto q$ be a pair of natural transformations $f,g\colon p\tto q$ between polynomial functors $p$ and $q$, and let $e\colon p'\to p$ be their equalizer in $\poly$ that we computed in the proof of \cref{thm.poly_limits}.
\begin{enumerate}
    \item Given a set $X$, show that $e_X\colon p'(X)\to p(X)$ is the equalizer of the $X$-components $f_X,g_X\colon p(X)\tto q(X)$ in $\smset$.
    \item Deduce that equalizers in $\poly$ coincide with equalizers in $\smset^\smset$.
    \item Conclude that limits in $\poly$ coincide with limits in $\smset^\smset$. \qedhere
\end{enumerate}
\begin{solution}
\begin{enumerate}
    \item By \cref{prop.morph_arena_to_func}, $f_X$ (resp.\ $g_X$) sends each $(i,h)\in p(X)$ with $i\in p(\1)$ and $h\colon p[i]\to X$ to $(f_\1i, f^\sharp_i\then h)$ (resp.\ $(g_\1i, g^\sharp_i\then h)$) in $q(X)$.
    So the equalizer of $f_X$ and $g_X$ is the set of all $(i,h)\in p(X)$ for which both $f_\1i = g_\1i$ and $f^\sharp_i\then h = g^\sharp_i\then h$.

    Indeed, by our construction of $p'$, the set $p'(X)$ consists of all pairs $(i,h')$ with $i\in p(\1)$ such that $f_\1i = g_\1i$ and $h'\colon p'[i]\to X$, where $p'[i]$ is the coequalizer of $f^\sharp_i,g^\sharp_i\colon q[f_\1i]\tto p[i]$.
    By the universal property of the coequalizer, functions $h'\colon p'[i]\to X$ precisely correspond to functions $h\colon p[i]\to X$ for which $f^\sharp_i\then h = g^\sharp_i\then h$.
    So $p'(X)$ is indeed the equalizer of $f_X$ and $g_X$.

    The equalizer natural transformation $e'\colon p'\to p$ has the inclusion $e'_X\colon p'(X)\to p(X)$ as its $X$-component, so by \cref{cor.morph_func_to_arena}, it is the lens whose on-positions function is the canonical equalizer inclusion $e'_\1\colon p'(\1)\to p(\1)$, while its on-directions function at $i\in p'(\1)$ is the map $p[i]\to p'[i]$ corresponding to the identity on $p'[i]$ given by the universal property of the coequalizer---which is just the canonical coequalizer map $p[i]\to p'[i]$.
    But this is exactly the lens $e\colon p'\to p$ constructed in the proof of \cref{prop.poly_prods}, as desired.
    \item By \cref{prop.presheaf_lim_ptwise}, limits---including equalizers---in $\smset^\smset$ are computed pointwise.
    So if $e_X\colon p'(X)\to p(X)$ is the equalizer of $f_X,g_X\colon p(X)\tto q(X)$ for every $X\in\smset$, then $e\colon p'\to p$ is the equalizer of $f,g\colon p(X)\tto q(X)$.
    \item We have just shown that equalizers in $\poly$ coincide with equalizers in $\smset^\smset$.
    We saw in the proof of \cref{prop.poly_prods} that products in $\poly$ also coincide with products in $\smset^\smset$.
    Since every limit can be computed as an equalizer of products, we can conclude that limits in $\poly$ coincide with limits in $\smset^\smset$.
\end{enumerate}
\end{solution}
\end{exercise}\index{limits!computed pointwise}

\index{polynomial functor!limits of polynomials|)}
\index{polynomial functor!colimits of polynomials|(}\index{colimit!in $\poly$}

\begin{theorem}\label{thm.poly_colimits}
The category $\poly$ has all small colimits.
\end{theorem}
\begin{proof}
A category has all small colimits if and only if it has coproducts and coequalizers, so by \cref{prop.poly_coprods}, it suffices to show that $\poly$ has coequalizers.

Let $s,t \colon p \tto q$ be two lenses.
We construct the coequalizer $q'$ of $s$ and $t$ as follows.
The pair of functions $s_\1, t_\1 \colon p(\1) \tto q(\1)$ define a graph $G \colon \fbox{$\bullet\tto\bullet$} \to \smset$ with vertices in $q(\1)$, edges in $p(\1)$, sources indicated by $s_\1$, and targets indicated by $t_\1$.
Then the set $C$ of connected components of $G$ is given by the coequalizer $g_\1 \colon q(\1) \to C$ of $s_\1$ and $t_\1$.
We define the position-set of $q'$ to be $C$.
Each direction-set of $q'$ will be a limit of a diagram of direction-sets of $p$ and $q$, but expressing this limit, as we proceed to do, is a bit involved.

\index{graph}

For each connected component $c \in C$, we have a connected subgraph $G_c \ss G$ with vertices $V_c \coloneqq g_\1\inv(c)$ and edges $E_c \coloneqq s_\1\inv(g_\1\inv(c)) = t_\1\inv(g_\1\inv(c))$.
Note that $E_c\ss p(\1)$ and $V_c\ss q(\1)$, so to each $e\in E_c$ (resp.\ to each $v\in V_c$) we have an associated direction-set $p[e]$ (resp.\ $q[v]$).

\index{element!category of elements}\index{category!of elements}

The category of elements $\int G_c$ has objects $E_c+V_c$ and two kinds of (non-identity) morphisms, $e \to s_\1(e)$ and $e \to t_\1(e)$, associated to each $e \in E_c$, all pointing from an object in $E_c$ to an object in $V_c$.
There is a functor $F \colon (\int G_c)\op \to \smset$ sending every $v \mapsto q[v]$, every $e \mapsto p[e]$, and every morphism to a function between them, namely either $s^\sharp_e \colon q[s_\1(e)] \to p[e]$ or $t^\sharp_e \colon q[t_\1(e)] \to p[e]$.
So we can define $q'[c]$ to be the limit of $F$ in $\smset$.

We claim that $q'\coloneqq\sum_{c\in C}\yon^{q'[c]}$ is the coequalizer of $s$ and $t$. We leave the complete proof to the interested reader in \cref{exc.poly_colimits}.
\end{proof}

\begin{exercise}\label{exc.poly_colimits}
Complete the proof of \cref{thm.poly_colimits} as follows:
\begin{enumerate}
	\item Provide a lens $g \colon q \to q'$.
	\item Show that $s \then g = t \then g$.
	\item Show that $g$ is a coequalizer of the pair $s, t$.
\qedhere
\end{enumerate}
\begin{solution}
\begin{enumerate}
    \item We define a lens $g \colon q \to q'$ as follows.
    The on-positions function $g_\1 \colon q(\1) \to q'(\1)$ is the coequalizer of $s_\1, t_\1 \colon p(\1) \tto q(\1)$.
    In particular, $g_\1$ sends each vertex in $q(\1)$ to its corresponding connected component in $q'(\1) = C$.
    Then for each $v \in q(\1)$, if we let its corresponding connected component be $c \coloneqq g_\1(v)$, we can define the on-directions function $g^\sharp_v \colon q'[c] \to q[v]$ to be the projection from the limit $q'[c]$ to its component $q[v]$.

    \item To show that $s \then g = t \then g$, we must show that both sides are equal on positions and on directions.
    The on-positions function $g_\1$ is defined to be the coequalizer of $s_\1$ and $t_\1$, so $s_\1 \then g_\1 = t_\1 \then g_\1$.
    So it suffices to show that for all $e \in p(\1)$, if we let its corresponding connected component be $c \coloneqq g_\1(s_\1(e)) = g_\1(t_\1(e))$, then the following diagram of on-directions functions commutes:
    \[
    \begin{tikzcd}[sep=small]
        & q[s_\1(e)] \ar[dl, "s^\sharp_e"'] \\
        p[e] & & q'[c] \ar[ul, "g^\sharp_{s_\1(e)}"'] \ar[dl, "g^\sharp_{t_\1(e)}"] \\
        & q[t_\1(e)] \ar[ul, "t^\sharp_e"]
    \end{tikzcd}
    \]
    But this is automatically true by the definition of $q'[c]$ as a limit---specifically the limit of a functor with $s^\sharp_e$ and $t^\sharp_e$ in its image---and the definitions of $g^\sharp_{s_\1(e)}$ and $g^\sharp_{t_\1(e)}$ as projections from this limit.

    \item To show that $g$ is the coequalizer of $s$ and $t$, it suffices to show that for any $r \in \poly$ and lens $f \colon q \to r$ satisfying $s \then f = t \then f$, there exists a unique lens $h \colon q' \to r$ for which $f = g \then h$, so that the following diagram commutes.
    \begin{equation*} %\label{eqn.eq_univ_prop}
    \begin{tikzcd}
        p \ar[r, "s", shift left] \ar[r, "t"', shift right] & q \ar[r, "g"] \ar[dr, "f"'] & q' \ar[d, "h", dashed] \\
        & & r
    \end{tikzcd}
    \end{equation*}
    In order for $f = g \then h$ to hold, we must have $f_\1 = g_\1 \then h_\1$ on positions.
    But we have that $s_\1 \then f_\1 = t_\1 \then f_\1$, so by the universal property of $q'(\1)$ and the map $g_\1$ as the coequalizer of $s_\1$ and $t_\1$ in $\smset$, there exists a unique $h_\1$ for which $f_\1 = g_\1 \then h_\1$.
    Hence $h$ is uniquely characterized on positions.
    In particular, it must send each connected component $c \in q'(\1)$ to the element in $r(\1)$ to which $f_\1$ sends every vertex $v \in V_c = g_\1\inv(c)$ that lies in the connected component $c$.

    Then for $f = g \then h$ to hold on directions, we must have that $f^\sharp_v = h^\sharp_{g_\1(v)} \then g^\sharp_v$ for each $v \in q(\1)$.
    Put another way, given $c \in q'(\1)$, we must have that $f^\sharp_v = h^\sharp_c \then g^\sharp_v$ for every $v \in V_c$.
    But $s \then f = t \then f$ implies that for each $e \in E_c = s_\1\inv(g_\1\inv(c)) = t_\1\inv(g_\1\inv(c)) \ss p(\1)$, the following diagram of on-directions functions commutes:
    \[
    \begin{tikzcd}[sep=small]
        & q[s_\1(e)] \ar[dl, "s^\sharp_e"'] \\
        p[e] & & r[f_\1(v)] \ar[ul, "f^\sharp_{s_\1(e)}"'] \ar[dl, "f^\sharp_{t_\1(e)}"] \\
        & q[t_\1(e)] \ar[ul, "t^\sharp_e"]
    \end{tikzcd}
    \]
    It follows that $r[f_\1(v)]$ together with the maps $(f^\sharp_v)_{v \in V_c}$ form a cone over the functor $F$.
    So by the universal property of the limit $q'[c]$ of $F$ with projection maps $(g^\sharp_v)_{v \in V_c}$, there exists a unique $h^\sharp_c \colon r[f_\1(v)] \to q'[c]$ for which $f^\sharp_v = h^\sharp_c \then g^\sharp_v$ for every $v \in V_c$.
    Hence $h$ is also uniquely characterized on directions, so it is unique overall.
    Moreover, we have shown that we can define $h$ on positions so that $f_\1 = g_\1 \then h_\1$, and that we can define $h$ on directions such that $f^\sharp_v = h^\sharp_c \then g^\sharp_v$ for all $c \in q'(\1)$ and $v \in V_c$.
    It follows that there exists $h$ for which $f = g \then h$.
\end{enumerate}
\end{solution}
\end{exercise}

\begin{example}
Given a diagram in $\poly$, one could either take its (co)limit as a diagram of \emph{polynomial} functors (i.e.\ its (co)limit in $\poly)$ or its (co)limit simply as a diagram of functors (i.e.\ its (co)limit in $\smset^\smset$).
We saw in \cref{exc.refl_limits} that in the case of limits, these yield the same result.
So, too, in the case of coproducts, per \cref{prop.poly_coprods}.

But in the case of general colimits, there are diagrams that yield different results: by the co-Yoneda lemma, \emph{every} functor $\smset \to \smset$---even those that are not polynomials---can be written as the colimit of representable functors in $\smset^\smset$, yet the colimit of the same representables in $\poly$ can only be another polynomial.

For a concrete example, consider the two distinct projections $\yon^\2\to\yon$, which form the diagram
\begin{equation} \label{eqn.2_projs}
    \yon^\2\tto\yon.
\end{equation}
According to \cref{thm.poly_colimits}, the colimit of \eqref{eqn.2_projs} in $\poly$ has the coequalizer of $\1 \tto \1$, namely $\1$, as its position-set, and the limit of the diagram $\1 \tto \2$ consisting of the two inclusions as its sole direction-set.
But this latter limit is just $\0$, so in fact the colimit of \eqref{eqn.2_projs} in $\poly$ is the constant functor $\1\yon^\0\iso\1$.

\index{colimit!$\poly$ disagrees with $\smset^\smset$}

But as functors, by \cref{prop.presheaf_lim_ptwise}, the colimit of \eqref{eqn.2_projs} can be computed pointwise: it is the (nonconstant!) functor
\[
  X\mapsto
  \begin{cases}
  	\0&\tn{ if }X=\0\\
  	\1&\tn{ if }X\neq\0
  \end{cases}
\]
\end{example}

\begin{exercise}\index{global sections}
By \cref{thm.adjoint_quadruple}, for any polynomial $p$, there are canonical lenses involving positions and global sections:
\[
	\epsilon \colon p(\1)\yon\to p
	\qqand
	\eta \colon p\to \yon^{\Gamma(p)}.
\]
\begin{enumerate}
	\item Characterize the behavior of the canonical lens $\epsilon \colon p(\1)\yon\to p$.
	\item Characterize the behavior of the canonical lens $\eta \colon p\to \yon^{\Gamma(p)}$.
	\item Show that the following is a pushout in $\poly$:
    \begin{equation} \label{eqn.pushout_adjoint}
    \begin{tikzcd}
    	p(\1)\yon\ar[r, "!"]\ar[d, "\epsilon"']&
    	\yon\ar[d, "!"]\\
    	p\ar[r, "\eta"']&
    	\yon^{\Gamma(p)}\ar[ul, phantom, very near start, "\ulcorner"]
    \end{tikzcd}
    \end{equation}
\end{enumerate}
\begin{solution}
\begin{enumerate}
    \item We characterize the lens $\epsilon \colon p(\1)\yon \to p$ as follows.
    On positions, it is the identity on $p(\1)$.
    Then for each $i \in p(\1)$, on directions, it is the unique map $p[i] \to \1$.

    \item We characterize the lens $\eta \colon p \to \yon^{\Gamma(p)}$ as follows.\index{global sections}
    On positions, it is the unique map $p(\1) \to \1$.
    Then for each $i \in p(\1)$, on directions, it is the canonical projection $\Gamma(p) \iso \prod_{i' \in p(\1)} p[i'] \to p[i]$.

    \item Showing that \eqref{eqn.pushout_adjoint} is a pushout square is equivalent to showing that, in the diagram
    \begin{equation} \label{eqn.coeq_adjoint}
    \begin{tikzcd}[sep=large]
        & \yon \ar[d, "\iota"] \ar[dr, "!"] \\
        p(\1)\yon \ar[ur, "!"] \ar[dr, "\epsilon"'] \ar[r, "s", shift left] \ar[r, "t"', shift right] & \yon + p \ar[r, "g"] & \yon^{\Gamma(p)} \\
        & p \ar[u, "\iota'"'] \ar[ur, "\eta"']
    \end{tikzcd}
    \end{equation}
    in which $\iota, \iota'$ are the canonical inclusions and the four triangles commute, $\yon^{\Gamma(p)}$ equipped with the lens $g$ is the coequalizer of $s$ and $t$.
    To do so, we apply \cref{thm.poly_colimits} to compute the coequalizer $q'$ of $s$ and $t$.
    The position-set of $q'$ is the coequalizer of $s_\1 = (! \then \iota)_\1$, which sends every $i \in p(\1)$ to the position of $\yon + p$ corresponding to the summand $\yon$, and $t_\1 = (\epsilon \then \iota')_\1$, which sends each $i \in p(\1)$ to the corresponding position in the summand $p$ of $\yon + p$.
    It follows that the coequalizer of $s_\1$ and $t_\1$ is $\1$, so $q'(\1) \iso \1$.

    Then the direction-set of $q'$ at its sole position is the limit of the functor $F$ whose image consists of lenses of the form $\1 \to \1$ or $p[i] \to \1$ for every $i \in p(\1)$.
    It follows that the limit of $F$ is just a product, namely $\prod_{i \in p(\1)} p[i] \iso \Gamma(p)$.
    Hence $q' \iso \yon^{\Gamma(p)}$, as desired.

    It remains to check that the upper right and lower right triangles in \eqref{eqn.coeq_adjoint} commute.
    The upper right triangle must commute by the uniqueness of morphisms $\yon \to \yon^{\Gamma(p)}$; and the lower right triangle must commute on positions.
    Moreover, the on-directions function of the coequalizer morphism $g$ at each position $i \in p(\1) \ss (\yon+p)(\1)$ must be the canonical projection $\Gamma(p) \to p[i]$, which matches the behavior of the corresponding on-directions function of $\eta$; hence the lower right triangle also commutes on directions.
\end{enumerate}
\end{solution}
\end{exercise}

\index{polynomial functor!pushout of polynomials}

\begin{proposition}\label{prop.tensor_as_pushout}
For polynomials $p,q$, the following is a pushout:
\[
\begin{tikzcd}
	p(\1)\yon\otimes q(\1)\yon\ar[r]\ar[d]&
	p(\1)\yon\otimes q\ar[d]\\
	p\otimes q(\1)\yon\ar[r]&
	p\otimes q\ar[ul, phantom,very near start, "\ulcorner"]
\end{tikzcd}
\]
\end{proposition}
\begin{proof}
All the lenses shown are identities on positions, so the displayed diagram is the coproduct over all $(i,j)\in p(\1)\times q(\1)$ of the diagram shown left
\[
\begin{tikzcd}
	\yon\ar[r]\ar[d]&
	\yon^{q[j]}\ar[d]\\
	\yon^{p[i]}\ar[r]&
	\yon^{p[i]\times q[j]}\ar[ul, phantom, very near start, "\ulcorner"]
\end{tikzcd}
\begin{tikzcd}
	\1\ar[dr, phantom, very near end, "\ulcorner"]&
	q[j]\ar[l]\\
	p[i]\ar[u]&
	p[i]\times q[j]\ar[l]\ar[u]
\end{tikzcd}
\]
where we used $(p\otimes q)[(i,j)]\cong p[i]\times q[j]$. This is the image under the Yoneda embedding of the diagram of sets shown right, which is clearly a pullback. The result follows by \cref{prop.yoneda_left_adjoint}.
\end{proof}\index{polynomial functor!pullback of polynomials}

This means that to give a lens $\varphi\colon p\otimes q\to r$, it suffices to give two lenses $\varphi_p\colon p\otimes q(\1)\yon\to r$ and $\varphi_q\colon p(\1)\yon\otimes q\to r$ that agree on positions. The lens $\varphi_p$ says how information about $q$'s position is transferred to $p$, and the lens $\varphi_q$ says how information about $p$'s position is transferred to $q$.

\index{monoidal structure!parallel product as pushout}
\index{parallel product!as pushout}

\begin{corollary}\label{cor.tensor_as_pushout}
Suppose we have polynomials $p_1,\ldots,p_n\in\poly$. Then $p_1\otimes\cdots\otimes p_n$ is isomorphic to the wide pushout of
\[
  \begin{tikzcd}[column sep=-15pt]
  	&
		\bigotimes_{i=1}^n \left(p_i(\1)\yon\right)\ar[dl]\ar[dr]\\
		p_1\otimes p_2(\1)\yon\otimes\cdots\otimes p_n(\1)\yon&
		\cdots&
		p_1(\1)\yon\otimes\cdots\otimes p_{n-1}(\1)\yon \otimes p_n
  \end{tikzcd}
\]
\end{corollary}
\begin{proof}
We proceed by induction on $n\in\nn$.
When $n=0$, the wide pushout has no legs and the empty parallel product is $\yon$, so the result holds.
If the result holds for $n$, then it holds for $n+1$ by \cref{prop.tensor_as_pushout}.
\end{proof}

\index{polynomial functor!colimits of polynomials|)}

%-------- Section --------%
\section{Vertical-cartesian factorization of lenses}

\index{factorization system!vertical-cartesian|(}

Aside from epi-mono factorization, there is another factorization system on $\poly$ that will show up frequently.

\begin{definition}[Vertical and cartesian lenses] \label{def.vert_cart}
Let $f\colon p\to q$ be a lens.
It is called \emph{vertical} if $f_\1\colon p(\1)\to q(\1)$ is an isomorphism.
It is called \emph{cartesian} if, for each $i\in p(\1)$, the function $f^\sharp_i\colon q[f(i)]\to p[i]$ is an isomorphism.
\end{definition}

\index{lens!vertical}
\index{lens!cartesian}

\begin{proposition}\label{prop.vert_cart_factorization}
Vertical and cartesian lenses form a factorization system of $\poly$.
\end{proposition}
\begin{proof}
It is easy to check that isomorphisms are both vertical and cartesian, and that vertical and cartesian lenses are each closed under composition.
It remains to show that every lens in $\poly$ can be uniquely (up to unique isomorphism) factored as a vertical lens composed with a cartesian lens.

Recall from \eqref{eqn.colax_poly_map} that a lens in $\poly$ can be written as to the left; we can thus rewrite it as to the right:
\[
\begin{tikzcd}[column sep=small]
	p(\1)\ar[dr, bend right, "{p[-]}"']\ar[rr, "f_\1"]&~&
	q(\1)\ar[dl, bend left, "{q[-]}"]\\&
	\smset\ar[u, phantom, near end, "\overset{f^\sharp}{\Leftarrow}"]
\end{tikzcd}
\hspace{1in}
\begin{tikzcd}
	p(\1)\ar[dr, bend right, "{p[-]}"']\ar[r, equal, ""' name=equal]&
	p(\1)\ar[d, "{q[f_\1(-)]}"]\ar[r, "f_\1"]&
	q(\1)\ar[dl, bend left, "{q[-]}"]\\&
	|[alias=set]|\smset\ar[from=equal, to=set, pos=.3, phantom, "\overset{f^\sharp}{\Leftarrow}"]
\end{tikzcd}
\]
We can see that the intermediary object $\sum_{i\in p(\1)} \yon^{q[f_\1i]}$ is unique up to unique isomorphism.
\end{proof}

\index{2-out-of-3}

\begin{proposition}
Vertical lenses satisfy 2-out-of-3: given $p\To{f}q\To{g}r$ with $h = f \then g$, if any two of $f,g,h$ are vertical, then so is the third.

If $g$ is cartesian, then $h$ is cartesian if and only if $f$ is cartesian.
\end{proposition}
\begin{proof}
Given $h = f \then g$, we have that $h_\1 = f_\1 \then g_\1$.
Since isomorphisms satisfy 2-out-of-3, it follows that vertical lenses satisfy 2-out-of-3 as well.

Now assume $g$ is cartesian.
On directions, $h = f \then g$ implies that for every $i \in p(\1)$, we have $h^\sharp_i = g^\sharp_{f_\1i} \then f^\sharp_i$.
Since $g^\sharp_{f_\1i}$ is an isomorphism, it follows that every $h^\sharp_i$ is an isomorphism if and only if every $f^\sharp_i$ is an isomorphism, so $h$ is cartesian if and only if $f$ is cartesian.
\end{proof}

\begin{exercise}
Give an example of polynomials $p,q,r$ and lenses $p\To{f}q\To{g}r$ such that $f$ and $f \then g$ are cartesian but $g$ is not.
\begin{solution}
Consider the lenses $\yon \To{f} \yon^\2 + \yon \To{g} \yon$ where $f$ is the canonical inclusion and $g$ is uniquely determined on positions and picks out $1 \in \2$ and $1 \in \1$ on directions.
Then the only on-directions function of $f$ is a function $\1 \to \1$, an isomorphism, so $f$ is cartesian.
Meanwhile, one of the on-directions functions of $g$ is a function $\1 \to \2$, which is not an isomorphism, so $g$ is not cartesian.
Finally, $f \then g$ can only be the unique lens $\yon \to \yon$, namely the identity, which is cartesian.
\end{solution}
\end{exercise}

Here is an alternative characterization of a cartesian lens in $\poly$.
Recall from \cref{exc.deriv-directions} that for any polynomial $p$, there is a corresponding function $\pi_p\colon\dot{p}(\1)\to p(\1)$, i.e.\ the set of all directions mapping to the set of positions.
A lens $(f_\1,f^\sharp)\colon p\to q$ can then be described as a function $f_\1\colon p(\1)\to q(\1)$ along with a function $f^\sharp$ that makes the following diagram in $\smset$ commute:
\begin{equation}\label{eqn.poly_map_usu}
\begin{tikzcd}
	\dot{p}(\1)\ar[d, "\pi_p"']&
	\bullet\ar[l, "f^\sharp"']\ar[r]\ar[d]&
	\dot{q}(\1)\ar[d, "\pi_q"]\\
	p(\1)\ar[r, equal]&
	p(\1)\ar[r, "f_\1"']&
	q(\1)\ar[ul, phantom, very near end, "\lrcorner"]
\end{tikzcd}
\end{equation}
Here, the pullback denoted by the dot $\bullet$ is the set of pairs comprised of a $p$-position $i$ and a $q[f_\1i]$-direction $e$.
The function $f^\sharp$ sends each such pair to a direction $f^\sharp_i(e)$ of $p$, and the commutativity of the left square implies that $f^\sharp_i(e)$ is specifically a $p[i]$-direction.
So $f^\sharp_i$ is indeed our familiar on-directions function $q[f_\1i]\to p[i]$, and $f^\sharp$ is just the sum of all these on-directions functions over $i\in p(\1)$.

\index{polynomial functor!pullback of polynomials}

\begin{exercise} \label{exc.cart_pullbacks}
Show that a lens $f\colon p\to q$ in $\poly$ is cartesian if and only if the square on the left hand side of \eqref{eqn.poly_map_usu} is also a pullback:
\[
\begin{tikzcd}
	\dot{p}(\1)\ar[d, "\pi_p"']&
	\bullet\ar[l, "f^\sharp"']\ar[r]\ar[d]&
	\dot{q}(\1)\ar[d, "\pi_q"]\\
	p(\1)\ar[r, equal]\ar[ur, phantom, very near end, "\llcorner"]&
	p(\1)\ar[r, "f_\1"']&
	q(\1)\ar[ul, phantom, very near end, "\lrcorner"]
\end{tikzcd}
\]
\begin{solution}
We wish to show that a lens $f\colon p\to q$ in $\poly$ is cartesian if and only if the square on the left hand side of \eqref{eqn.poly_map_usu} is a pullback.
We already know that that square commutes, so it is a pullback if and only if $f^\sharp$ is an isomorphism.
The right pullback square tells us that the $\bullet$ is $\sum_{i\in p(\1)}q[f_\1i]$.
So $f^\sharp_i\colon q[f_\1i]\to p[i]$ is an isomorphism for every $i\in p(\1)$ if and only if their sum $f^\sharp\colon\sum_{i\in p(\1)}q[f_\1i]\to\sum_{i\in p(\1)}p[i]\iso\dot{p}(\1)$ is an isomorphism as well.
Hence $f$ is cartesian if and only if $f^\sharp$ is an isomorphism, as desired.
\end{solution}
\end{exercise}

\begin{exercise}\index{global sections}
Is the pushout of a cartesian lens always cartesian?
\begin{solution}
The pushout of a cartesian lens is \emph{not} necessarily cartesian.
Take the pushout square \eqref{eqn.pushout_adjoint}.
The lens $!\colon p(\1)\yon\to\yon$ has $\1\to\1$ as every on-directions function, so it is cartesian, but its pushout $\eta\colon p\to\yon^{\Gamma(p)}$ is not going to be cartesian as long as there is some $i\in p(\1)$ for which $\Gamma(p)\not\iso p[i]$.
For instance, when $p\coloneqq\yon+\1$, we have that $\Gamma(p)\iso\0\not\iso\1\iso p[1]$, so $\eta$ is not cartesian.
\end{solution}
\end{exercise}

Why do we use the word \emph{cartesian} to describe cartesian morphisms? It turns out that, as natural transformations, cartesian morphisms are precisely what are known as cartesian natural transformations.

\index{natural transformation!cartesian|seealso{lens!cartesian}}\index{polynomial functor!pullback of polynomials}\index{pullback!cartesian natural transformation and}

\begin{definition}[Cartesian natural transformation] \label{def.cart_nat_trans}
A \emph{cartesian natural transformation} is a natural transformation whose naturality squares are all pullbacks.
That is, given categories $\cat{C},\cat{D}$, functors $F,G$, and natural transformation $\alpha$, we say that $\alpha$ is \emph{cartesian} if for all morphisms $h\colon c\to c'$ in $\cat{C}$,
\[
\begin{tikzcd}
    Fc \ar[d, "Fh"'] \ar[r, "\alpha_c"] & Gc \ar[d, "Gh"] \\
    Fd \ar[r, "\alpha_d"'] & Gd \ar[ul, phantom, very near end, "\lrcorner"]
\end{tikzcd}
\]
is a pullback.
\end{definition}

\begin{proposition}\label{prop.cart_as_nt}
Let $f\colon p\to q$ be a morphism in $\poly$. The following are equivalent:
	\begin{enumerate}
		\item viewed as a lens, $f$ is cartesian in the sense of \cref{def.vert_cart}: for each $i\in p(\1)$, the on-directions function $f^\sharp_i$ is a bijection;
		\item the square on the left hand side of \eqref{eqn.poly_map_usu} is also a pullback:
\[
\begin{tikzcd}
	\dot{p}(\1)\ar[d, "\pi_p"']&
	\bullet\ar[l, "f^\sharp"']\ar[r]\ar[d]&
	\dot{q}(\1)\ar[d, "\pi_q"]\\
	p(\1)\ar[r, equal]\ar[ur, phantom, very near end, "\llcorner"]&
	p(\1)\ar[r, "f_\1"']&
	q(\1)\ar[ul, phantom, very near end, "\lrcorner"]
\end{tikzcd}
\]
		\item viewed as a natural transformation, $f$ is cartesian in the sense of \cref{def.cart_nat_trans}: for any sets $A,B$ and function $h\colon A\to B$, the naturality square
\begin{equation} \label{eqn.cart_nt_pullback}
\begin{tikzcd}
	p(A)\ar[r, "f_A"]\ar[d, "p(h)"']&
	q(A)\ar[d, "q(h)"]\\
	p(B)\ar[r, "f_B"']&
	q(B)\ar[ul, phantom, very near end, "\lrcorner"]
\end{tikzcd}
\end{equation}
is a pullback.
  \end{enumerate}
\end{proposition}
\begin{proof}
We already showed that the first two are equivalent in \cref{exc.cart_pullbacks}, and we will complete this proof in \cref{exc.cart_as_nt}.
\end{proof}

\begin{exercise} \label{exc.cart_as_nt}
In this exercise, you will complete the proof of \cref{prop.cart_as_nt}.

First, we will show that $1\Rightarrow3$.
In the following, let $f\colon p\to q$ be a cartesian lens in $\poly$ and $h\colon A\to B$ be a function.
\begin{enumerate}
    \item Using \cref{prop.morph_arena_to_func} to translate $f$ from a lens in $\poly$ to a natural transformation and \cref{prop.poly_on_functions} to interpret $q(h)$, characterize the pullback of $p(B)\To{f_B}q(B)\From{q(h)}q(A)$ in $\smset$.
    \item Show that this pullback coincides with the naturality square \eqref{eqn.cart_nt_pullback}, hence proving $1\Rightarrow3$.
\end{enumerate}
Next, we show that $3\Rightarrow1$.
In the following, let $f\colon p\to q$ be a lens in $\poly$ that is cartesian when viewed as a natural transformation, so that \eqref{eqn.cart_nt_pullback} is a pullback for any function $h\colon A\to B$.
Also fix $i\in p(\1)$.
\begin{enumerate}[resume]
    \item Show that the diagram
    \begin{equation} \label{eqn.cart_nt_pullback_cone}
        \begin{tikzcd}[column sep=50pt]
        	\1 \ar[r, "{(f_\1i,\,\id_{q[f_\1i]})}"]\ar[d, "{(i,\,\id_{p[i]})}"']&
        	q(q[f_\1i])\ar[d, "q(f^\sharp_i)"]\\
        	p(p[i])\ar[r, "f_{p[i]}"']&
        	q(p[i])\ar[ul, phantom, very near end, "\lrcorner"]
        \end{tikzcd}
    \end{equation}
    commutes.
    Hint: Use \cref{prop.poly_on_functions}, \cref{prop.morph_arena_to_func}, and/or \cref{cor.morph_func_to_arena}.
    \item Apply the universal property of the pullback \eqref{eqn.cart_nt_pullback} to the diagram \eqref{eqn.cart_nt_pullback_cone} above to exhibit an element of $p(q[f_\1i])$.
    Conclude from the existence of this element that $f^\sharp_i$ is an isomorphism, hence proving $3\Rightarrow1$.\qedhere
\end{enumerate}
\begin{solution}
First, we will show that $1\Rightarrow3$ in \cref{prop.cart_as_nt}.
Here $f\colon p\to q$ is a cartesian lens in $\poly$ and $h\colon A\to B$ is a function.
\begin{enumerate}
    \item An element of $p(B)$ is a pair comprised of a $p$-position $i$ and a function $k\colon p[i]\to B$, and \cref{prop.morph_arena_to_func} tells us that $f_B\colon p(B)\to q(B)$ sends $(i,k)\mapsto(f_\1i,f^\sharp_i\then k)$.
    Meanwhile, an element of $q(A)$ is a pair comprised of a $q$-position $j$ and a function $\ell\colon q[j]\to A$, and \cref{prop.poly_on_functions} tells us that $q(h)$ sends $(j,\ell)\mapsto(j,\ell\then h)$.
    So $((i,k),(j,\ell))$ is in the pullback of $p(B)\To{f_B}q(B)\From{q(h)}q(A)$ if and only if $f_\1i=j$ and $f^\sharp_i\then k=\ell\then h$.

    As $f$ is cartesian, $f^\sharp_i$ is an isomorphism, so we can rewrite the latter equation as $k=g_i\then\ell\then h$, where $g_i$ is the inverse of $f^\sharp_i$.
    In fact, if we let $\ell'\coloneqq g_i\then\ell$, we observe that the values of $j,k,$ and $\ell$ are all already determined by the values of $i$ and $\ell'$: we have that $j=f_\1i$, that $k=\ell'\then h$, and that $\ell=f^\sharp_i\then\ell'$
    It follows that the pullback is equivalently the set of pairs $(i,\ell')$ comprised of a $p$-position $i$ and a function $\ell'\colon p[i]\to A$ (with no other restrictions on $i$ and $\ell'$).
    The projection from the pullback to $p(B)$ sends $(i,\ell')\mapsto(i,\ell'\then h)$, and the projection from the pullback to $q(A)$ sends $(i,\ell')\mapsto(f_\1i,f^\sharp_i\then\ell')$.
    \item The pullback described above---the set of pairs $(i,\ell')$ comprised of a $p$-position $i$ and a function $\ell'\colon p[i]\to A$---is exactly the set $p(A)$.
    Moreover, the projection to $p(B)$ sending $(i,\ell')\mapsto(i,\ell'\then h)$ is $p(h)$, and the projection to $q(A)$ sending $(i,\ell')\mapsto(f_\1i,f^\sharp_i\then\ell')$ is $f_A$ by \cref{prop.morph_arena_to_func}.
    So \eqref{eqn.cart_nt_pullback} is a pullback, as desired.
\end{enumerate}
Next, we will show that $3\Rightarrow1$ in \cref{prop.cart_as_nt}, with $f\colon p\to q$ as a lens in $\poly$ that is a cartesian natural transformation and $i\in p(\1)$.
\begin{enumerate}[resume]
    \item By \cref{cor.morph_func_to_arena}, we have that $f_{p[i]}$ sends $(i,\id_{p[i]})\mapsto(f_\1i,f^\sharp_i)$, and by \cref{prop.poly_on_functions}, we have that $q(f^\sharp_i)$ sends $(f_\1i,\id_{q[f_\1i]})\mapsto(f_\1i,f^\sharp_i)$ as well.
    Hence \eqref{eqn.cart_nt_pullback_cone} commutes.

    \item Taking $A\coloneqq q[f_\1i], B\coloneqq p[i],$ and $h\coloneqq f^\sharp_i$ in \eqref{eqn.cart_nt_pullback} and applying its universal property to \eqref{eqn.cart_nt_pullback_cone} induces an element $(i',g)$ of $p(q[f_\1i])$, with $i'\in p(\1)$ and $g\colon p[i']\to q[f_\1i]$, such that $p(f^\sharp_i)$ sends $(i',g)\mapsto(i,\id_{p[i]})$ and $f_{q[f_\1i]}$ sends $(i',g)\mapsto(f_\1i,\id_{q[f_\1i]})$.
    It follows from the behavior of $p(f^\sharp_i)$ (by \cref{prop.poly_on_functions}) that $i'=i$ and $g\then f^\sharp_i=\id_{p[i]}$, and it follows from the behavior of $f_{q[f_\1i]}$ (by \cref{prop.morph_arena_to_func}) that $f^\sharp_i\then g=\id_{q[f_\1i}$.
    So $g$ is the inverse of $f^\sharp_i$, proving that $f^\sharp_i$ is an isomorphism, as desired.
\end{enumerate}
\end{solution}
\end{exercise}

\index{monoidal structure!preservation of vertical and cartesian maps}

\begin{proposition}\label{prop.monoidal_pres_vert_cart}
The monoidal structures $+$, $\times$, and $\otimes$ preserve both vertical and cartesian morphisms.
\end{proposition}
\begin{proof}
Suppose that $f\colon p\to p'$ and $g\colon q\to q'$ are vertical, so that the on-positions functions $f_\1$ and $g_\1$ are isomorphisms.

We can obtain the on-positions function of a lens by passing it through the functor $\poly\To{p(\1)}\smset$ from \cref{thm.adjoint_quadruple}.
As this functor is both a left adjoint and a right adjoint, it preserves both sums and products, so $(f+g)_\1 = f_\1+g_\1$ and $(f\times g)_1 = f_1\times g_1$.
Hence $f+g$ and $f\times g$ are both vertical.
On-positions, the behavior of $\otimes$ is identical to the behavior of $\times$, so $f\otimes g$ must be vertical as well.

Now suppose that $f\colon p\to p'$ and $g\colon q\to q'$ are cartesian.

A position of $p+q$ is a position $i\in p(\1)$ or a position $j\in q(\1)$, and the map $(f+g)^\sharp$ at that position is either $f^\sharp_i$ or $g^\sharp_j$; either way it is an isomorphism, so $f+g$ is cartesian.

A position of $p\times q$ (resp.\ of $p\otimes q$) is a pair $(i,j)\in p(\1)\times q(\1)$. The lens $(f\times g)^\sharp_{(i,\,j)}$ (resp.\ $(f\otimes g)^\sharp_{(i,\,j)}$) is $f^\sharp_i+g^\sharp_j$ (resp.\ $f^\sharp_i\times g^\sharp_j$) which is again an isomorphism if $f^\sharp_i$ and $g^\sharp_j$ are. Hence $f\times g$ (resp.\ $f\otimes g$) is cartesian, completing the proof.
\end{proof}

\begin{proposition}\label{prop.pullback_vert_cart}
Pullbacks preserve vertical (resp.\ cartesian) lenses.
In other words, if $f\colon p\to q$ is a lens and $g\colon q'\to q$ a vertical (resp.\ cartesian) lens, then the pullback $g'$ of $g$ along $p$
\[
\begin{tikzcd}
	p\times_qq'\ar[r]\ar[d, "g'"']&
	q'\ar[d, "g"]\\
	p\ar[r, "f"']&
	q\ar[ul, phantom, very near end, "\lrcorner"]
\end{tikzcd}
\]
is vertical (resp.\ cartesian).
\end{proposition}
\begin{proof}
This follows from \cref{ex.pullbacks_in_poly}, since the pullback (resp.\ pushout) of an isomorphism is an isomorphism.
\end{proof}

\index{factorization system!vertical-cartesian|)}

%-------- Section --------%
\section{Monoidal $*$-bifibration over $\smset$}

\index{monoidal $*$-bifibraion}

We conclude this chapter by showing that the functor $p\mapsto p(\1)$ has special properties that make it what \cite{shulman2008framed} refers to as a \emph{monoidal $*$-bifibration}.
Roughly speaking, this means that $\smset$ acts as a sort of remote controller on the category $\poly$, grabbing every polynomial by its positions and pushing or pulling it this way and that.
The material in this section is even more technical than the rest of this chapter, and we won't use it again in the book, so the reader may wish to skip to \cref{part.comon}.

As an example, suppose one has a set $A$ and a function $f\colon A\to p(\1)$, which we can also think of as a cartesian lens between constant polynomials.
From $f$, we can obtain a new polynomial $f^*p$ with position-set $A$ via a pullback
\begin{equation}\label{eqn.f^*_defined}
\begin{tikzcd}
	f^*p\ar[r, "\fun{cart}"]\ar[d]&
	p\ar[d, "\eta_p"]\\
	A\ar[r, "f"']&
	p(\1)\ar[ul, phantom, very near end, "\lrcorner"]
\end{tikzcd}
\end{equation}
Here $\eta_p$ is the unit of the adjunction $\adjr{\smset}{A}{p(\1)}{\poly}$; it is a vertical lens.
We could evaluate this pullback using \cref{ex.pullbacks_in_poly}.
Alternatively, we can use \cref{prop.pullback_vert_cart} to deduce that the top lens $f^*p\to p$ (which we presciently labeled $\fun{cart}$) is cartesian like $f$ and that the left lens $f^*p\to A$ is vertical like $\eta_p$. Furthermore, $\fun{cart}_1 = f$.
Hence
\[
    f^*p \iso \sum_{a \in A} \yon^{p[f(a)]}.
\]
We'll see this as part of a bigger picture in \cref{prop.basechange,thm.triple_adjoint_basechange}, but first we need the following definitions and a result about cartesian lenses.

\begin{definition}[Slice category] \label{def.slice}
Given an object $c$ in a category $\cat{C}$, the \emph{slice category} of $\cat{C}$ over $c$, denoted $\cat{C}/c$, is the category whose objects are morphisms in $\cat{C}$ with codomain $c$ and whose morphisms are commutative triangles in $\cat{C}$.
\end{definition}

\begin{definition}[Exponentiable morphism]
Given a category $\cat{C}$ with objects $c, d$ and morphism $f \colon c \to d$ such that all pullbacks along $f$ exist in $\cat{C}$, we say that $f$ is \emph{exponentiable} if the functor $f^* \colon \cat{C}/d \to \cat{C}/c$ given by pulling back along $f$ is a left adjoint.
\end{definition}

\index{lens!cartesian}\index{lens!exponentiable}\index{exponentiable lens|see{lens, cartesian}}

\begin{theorem}\label{thm.cart_exponentiable}
Cartesian lenses in $\poly$ are exponentiable.
That is, if $f\colon p\to q$ is cartesian, then the functor $f^*\colon\poly/q\to\poly/p$ given by pulling back along $f$ is a left adjoint:
\[
\begin{tikzcd}[column sep=50pt, background color=theoremcolor]
	\poly/p\ar[r, shift right=7pt, "f_*"']&
	\poly/q\ar[l, shift right=7pt, "f^*"']\ar[l, phantom, "\Leftarrow"]
\end{tikzcd}
\]
\end{theorem}
\begin{proof}
Fix $e\colon p'\to p$ and $g\colon q'\to q$.
\[
\begin{tikzcd}
	p'\ar[d, "e"']&q'\ar[d, "g"]\\
	p\ar[r, "f"']&q
\end{tikzcd}
\]
We need to define a functor $f_*\colon\poly/p\to\poly/q$ and prove the analogous isomorphism establishing it as right adjoint to $f^*$. We first establish some notation. Given a set $Q$ and sets $(P'_i)_{i\in I}$, each equipped with a map $Q\to P'_i$, let $Q/\sum_{i\in I}P'_i$ denote the coproduct in $Q/\smset$, or equivalently the wide pushout of sets $P'_i$ with apex $Q$. Then we give the following formula for $f_*p'$, which we write in larger font for clarity:
\begin{equation}\label{eqn.cart_exp}
f_*p'\coloneqq
\scalebox{1.3}{$\displaystyle
\sum_{j\in q(\1)}\;\sum_{i'\in\prod\limits_{i\in f_\1\inv(j)}e_\1\inv(i)}\;\yon^{q[j]/\sum_{i\in f_\1\inv(j)}p'[i'(i)]}
$}
\end{equation}
Again, $q[j]/\sum_{i\in f_\1\inv(j)}p'[i'(i)]$ is the coproduct of the $p'[i'(i)]$, taken in $q[j]/\smset$. Since $p[i]\cong q[f(i)]$ for any $i\in p(\1)$ by the cartesian assumption on $f$, we have the following chain of natural  isomorphisms:\index{isomorphism!natural}
\begin{align*}
	&\:(\poly/p)(f^*q', p')\\
    \iso&
	\prod_{i\in p(\1)}\;\prod_{j'\in g_\1\inv(f_\1 i)}\;\sum_{i'\in e_\1\inv(i)}(p[i]/\smset)(p'[i'],p[i]+_{q[f_\1i]}q'[j'])\\
    \iso&
	\prod_{i\in p(\1)}\;\prod_{j'\in g_\1\inv(f_\1 i)}\;\sum_{i'\in e_\1\inv(i)}(q[f_\1i]/\smset)(p'[i'],q'[j'])\\
    \iso&
	\prod_{j\in q(\1)}\;\prod_{j'\in g_\1\inv(j)}\;\prod_{i\in f_\1\inv(j)}\;\sum_{i'\in e_\1\inv(i)}(q[j]/\smset)(p'[i'],q'[j'])\\
    \iso&
	\prod_{j\in q(\1)}\;\prod_{j'\in g_\1\inv(j)}\;\sum_{i'\in\prod_{i\in f_\1\inv(j)}e_\1\inv(i)}\;\prod_{i\in f_\1\inv(j)}(q[j]/\smset)(p'[i'(i)],q'[j'])\\
    \iso&
	\prod_{j\in q(\1)}\;\prod_{j'\in g_\1\inv(j)}\;\sum_{i'\in\prod_{i\in f_\1\inv(j)}e_\1\inv(i)}(q[j]/\smset)\!\left(\sum_{i\in f_\1\inv(j)}p'[i'(i)],q'[j']\right)\\
    \iso&\:
	(\poly/q)(q',f_*p')
\end{align*}
\end{proof}

\begin{example}
Let $p\coloneqq\2\yon^\2$, $q\coloneqq\yon^\2+\yon^\0$, and $f\colon p\to q$ the unique cartesian lens between them.
Then for any $e\colon p'\to p$ over $p$, \eqref{eqn.cart_exp} provides the following description for the pushforward $f_*p'$.
%We use the isomorphisms $p(\1)\cong\2$ and $q(\1)\cong\2$ to talk about the positions of $p$ and $q$.

Over the $j=2$ position, $f_\1\inv(2)=\0$ and $q[2]=\0$, so $\prod_{i \in f_\1\inv(2)} e_\1\inv(i)$ is an empty product and $q[2]/\sum_{i\in f_\1\inv(2)} p'[i'(i)]$ is an empty pushout.
Hence the corresponding summand of \eqref{eqn.cart_exp} is simply $\yon^\0\cong\1$.

Over the $j=1$ position, $f_\1\inv(1)=\2$ and $q[1]=p[1]=p[2]=\2$, so $\prod_{i'\in f_\1\inv(1)} e_\1\inv(i) \iso e_\1\inv(1)\times e_\1\inv(2)$.
For $i' \in e_\1\inv(1) \times e_\1\inv(2)$, we have that $q[1]/\sum_{i\in f_\1\inv(2)} p'[i'(i)] \iso X_{i'}$ in the following pushout square:
\[
\begin{tikzcd}
	X_{i'} \ar[from=r] \ar[from=d] &
	p'[i'(2)] \ar[from=d, "e^\sharp_{i'(2)}"'] \\
	p'[i'(1)] \ar[from=r, "e^\sharp_{i'(1)}"] &
	\2 \ar[ul, phantom, very near end, "\lrcorner"]
\end{tikzcd}
\]
Then in sum we have
\[
    f_*p' \iso \left(\sum_{i' \in e_\1\inv(1) \times e_2\inv(2)} \yon^{X_{i'}}\right) + \1.
\]
\end{example}

\begin{exercise}
Prove that the unique lens $f\colon\yon\to\1$ is exponentiable.
\begin{solution}
Choose $p\in\poly$ and $q'\in\poly/\yon$. Then there is $q\in\poly$ such that $q'\cong q\yon$, equipped with the projection $q\yon\to\yon$. The pushforward is given by the exponential
\[f_*(q\yon)\coloneqq q^\yon\]
from the cartesian closure; see \eqref{eqn.exponential}. Indeed, we have
\begin{align*}
	\poly/\yon(f^*p,q\yon)&\cong
	\poly/\yon(p\yon,q\yon)\\&\cong
	\poly(p\yon,q)\\&\cong
	\poly(p,q^\yon).
\end{align*}
\end{solution}
\end{exercise}

For any set $A$, let $\poly[A.]$ denote the category whose objects are polynomials $p$ equipped with an isomorphism $A\cong p(\1)$, and whose morphisms are lenses respecting the isomorphisms with $A$.

\begin{proposition}[Base change]\label{prop.basechange}
For any function $f\colon A\to B$, pullback $f^*$ along $f$ induces a functor $\poly[B.]\to\poly[A.]$, which we also denote $f^*$.
\end{proposition}
\begin{proof}
This follows from \eqref{eqn.pullback_poly} with $q\coloneqq A$ and $r\coloneqq B$, since pullback of an iso is an iso.
\end{proof}

\begin{theorem}\label{thm.triple_adjoint_basechange}
For any function $f\colon A\to B$, the pullback functor $f^*$ has both a left and a right adjoint
\begin{equation}\label{eqn.adjoint_triple_monoidal_fib}
\begin{tikzcd}[column sep=large, background color=theoremcolor]
	\poly[A.]\ar[r, shift left=16pt, "f_!"]\ar[r, shift right=16pt, "f_*"']
	\ar[r, phantom, shift left=9pt, "\Rightarrow"]\ar[r, phantom, shift right=9pt, "\Leftarrow"]
&
	\poly[B.]\ar[l, "f^*" description]
\end{tikzcd}
\end{equation}
Moreover $\otimes$ preserves the op-cartesian arrows, which makes this a monoidal $*$-bifibration in the sense of \cite[Definition 12.1]{shulman2008framed}.
\end{theorem}
\begin{proof}
Let $p$ be a polynomial with $p(\1)\cong A$. Then the formula for $f_!p$ and $f_*p$ are given as follows:
\begin{equation}\label{eqn.f_!andf_*}
f_!p\cong\scalebox{1.3}{$\displaystyle\sum_{b\in B}\yon^{\;\prod\limits_{a\mapsto b}p[a]}$}
\qqand
f_*p\cong\scalebox{1.3}{$\displaystyle\sum_{b\in B}\yon^{\;\sum\limits_{a\mapsto b}p[a]}$}
\end{equation}
It may at first be counterintuitive that the left adjoint $f_!$ involves a product and the right adjoint $f_*$ involves a sum. The reason for this comes from the fact that $\poly$ is equivalent to the Grothendieck construction applied to the functor $\smset\op\to\smcat$ sending each set $A$ to the category $(\smset/A)\op$. The fact that functions $f\colon A\to B$ induces an adjoint triple between $\smset/A$ and $\smset/B$, and hence between $(\smset/A)\op$ and $(\smset/B)\op$ explains the variance in \eqref{eqn.f_!andf_*} and simultaneously establishes the adjoint triple \eqref{eqn.adjoint_triple_monoidal_fib}.

The functor $p\mapsto p(\1)$ is strong monoidal with respect to $\otimes$ and strict monoidal if we choose the lens construction as our model of $\poly$. By \cref{prop.monoidal_pres_vert_cart}, the monoidal product $\otimes$ preserves cartesian lenses; thus we will have established the desired monoidal $*$-bifibration in the sense of \cite[Definition 12.1]{shulman2008framed} as soon as we know that $\otimes$ preserves op-cartesian lenses.

Given $f$ and $p$ as above, the op-cartesian lens is the lens $p\to f_!p$ obtained as the composite $p\to f^*f_!p\to f_!p$ where the first lens is the unit of the $(f_!,f^*)$ adjunction and the second is the cartesian lens for $f_!p$. On positions $p\to f_!p$ acts as $f$, and on directions it is given by projection.

If $f\colon p(\1)\to B$ and $f'\colon p'(\1)\to B'$ are functions then we have
\begin{align*}
	f_!(p)\otimes f'_!(p')&\cong
	\sum_{b\in B}\sum_{b'\in B'}\yon^{\big(\prod_{a\mapsto b}p[a]\big)\times\big(\prod_{a'\mapsto b'}p'[a']\big)}\\&\cong
	\sum_{(b,\,b')\in B\times B'}\yon^{\big(\prod_{(a,\,a')\mapsto(b,\,b')}p[a]\times p[b]\big)}\\&
	\cong (f_!\otimes f'_!)(p\otimes p')
\end{align*}
and the op-cartesian lenses are clearly preserved since projections in the second line match with projections in the first.
\end{proof}

%-------- Section --------%
\section[Summary and further reading]{Summary and further reading%
  \sectionmark{Summary \& further reading}}
\sectionmark{Summary \& further reading}

In this chapter we discussed several of the nice properties of the category $\poly$: it has various adjunctions to $\smset$ and $\smset\op$, is Cartesian closed, has limits and colimits, has an epi-mono factorization system, has a vertical-cartesian factorization system, and comes with a monoidal $*$-bifibration to $\smset$.

The principal monomial functor $p\mapsto p(1)\yon^{\Gamma p}$ discussed in \cref{cor.principal_monomial} is in fact distributive monoidal, and this comes up in work on entropy \cite{spivak2022polynomial} and on noncooperative strategic games \cite{capucci2022diegetic}.

\index{monomial!principal}\index{principal monomial|see{monomial, principal}}

%-------- Section --------%
\section{Exercise solutions}
\Closesolutionfile{solutions}
{\footnotesize
\input{solution-file5}}

\setcounter{chapter}{5}%Just finished 5.

\part{A different category of categories}\label{part.comon}

\Opensolutionfile{solutions}[solution-file6]

%------------ Chapter ------------%
\chapter{The composition product}\label{ch.comon.comp}

\index{monoidal structure!substitution|see{composition product}}
\index{composition product|(}\index{interface!changing}

We have seen that the category $\poly$ of polynomial functors has plenty of interoperating mathematical structure. Further, it is an expressive way to talk about dynamical systems that can change their interfaces and wiring patterns based on their internal states.

But we touched upon one aspect of the theory---what in some sense is the most interesting part of the story---only briefly. That aspect is quite simple to state, and yet has profound consequences. Namely, polynomials can be composed:
\[
\yon^\2\circ(\yon+\1)=(\yon+\1)^\2\iso\yon^\2+\2\yon+\1.
\]
In other words, $(\yon+\1)$ is \emph{substituted in} for the variable $\yon$ in $\yon^\2$. What could be simpler?
\index{polynomial functor, substitution of polynomials|see{polynomial functor!composition of polynomials}}
\index{polynomial functor!composition of polynomials}

It turns out that this operation, which we will soon see is a monoidal product, has a lot to do with time.
\index{time!composition product and}
There is a strong sense---made precise in \cref{prop.poly_closed_comp}---in which the polynomial $p\circ q$ represents ``starting at a position $i$ in $p$, choosing a direction in $p[i]$, landing at a position $j$ in $q$, choosing a direction in $q[j]$, and then landing... somewhere.''
This is exactly what we need to run through multiple steps of a dynamical system, the very thing we didn't know how to do in \cref{ex.do_nothing}.
We'll continue that story in \cref{subsec.comon.comp.def.dyn_sys}.

The composition product has many surprises up its sleeve, as we'll see throughout the rest of the book.%
\footnote{Some authors refer to $\tri$ as the \emph{substitution} product, rather than the composition product. We elected to use the composition product terminology because it provides a good noun form ``the composite'' for $p\tri q$, whereas ``the substitute'' is somehow strange in English.}

%-------- Section --------%
\section{Defining the composition product}\label{sec.comon.comp.def}
We begin with the definition of the composition product in terms of polynomials as functors.

\subsection{Composite functors}\label{subsec.comon.comp.def.functor}

\begin{definition}[Composition product] \label{def.comp}
Given polynomial functors $p, q$, we let $p \circ q$ denote their \emph{composition product}, or their composite as functors.
That is, $p \circ q \colon \smset \to \smset$ sends each set $X$ to the set $p(q(X))$.
\end{definition}
\index{composition product!as composition of functors}

Functor composition gives a monoidal structure on the category $\smset^\smset$ of functors $\smset\to\smset$, but to check that the full subcategory $\poly$ of $\smset^\smset$ inherits this monoidal structure, we need to verify that the composite of two functors in $\poly$ is still a functor in $\poly$.

\begin{proposition}\label{prop.poly_closed_comp}
Suppose $p,q\in\poly$ are polynomial functors. Then their composite $p\circ q$ is again a polynomial functor, and we have the following isomorphism:
\begin{equation} \label{eqn.composite_formula_circ}
p\circ q\iso\sum_{i\in p(\1)}\prod_{a\in p[i]}\sum_{j\in q(\1)}\prod_{b\in q[j]}\yon.
\end{equation}
\end{proposition}
\begin{proof}
We can rewrite $p$ and $q$ as
\[
p\iso\sum_{i\in p(\1)}\yon^{p[i]}\iso\sum_{i\in p(\1)}\prod_{a\in p[i]}\yon
\qqand
q\iso\sum_{j\in q(\1)}\yon^{q[j]}\iso\sum_{j\in q(\1)}\prod_{b\in q[j]}\yon.
\]
For any set $X$ we have
\[
  (p\circ q)(X)=p(q(X))\iso p\left(\sum_{j\in q(\1)}\prod_{b\in q[j]}X\right)\iso\sum_{i\in p(\1)}\prod_{a\in p[i]}\sum_{j\in q(\1)}\prod_{b\in q[j]}X,
\]
so \eqref{eqn.composite_formula_circ} is indeed the formula for the composite $p \circ q$.
To see this is a polynomial, we use \eqref{eqn.push_prod_sum_set_indep}, which says we can rewrite the $\prod\sum$ in \eqref{eqn.composite_formula_circ} as a $\sum\prod$ to obtain
\begin{align}\label{eqn.composite_formula_sums_first_circ}
  p\circ q\iso
  \scalebox{1.3}{$\displaystyle
  \sum_{i\in p(\1)} \; \sum_{\ol{j}\colon p[i]\to q(\1)}\yon^{\sum_{a\in p[i]}q[\ol{j}(a)]}$}
\end{align}
(written larger for clarity), which is a polynomial.
\end{proof}

\index{composition product!formula for}

\begin{corollary} \label{cor.comp_monoidal}
The category $\poly$ has a monoidal structure $(\yon,\circ)$, where $\yon$ is the identity functor and $\circ$ is given by composition.
\end{corollary}

\index{composition product!unit of}

Because we may wish to use $\circ$ to denote composition in arbitrary categories, we use a special symbol for polynomial composition, namely
\[
p\tri q\coloneqq p\circ q.
\]
The symbol $\tri$ looks a bit like the composition symbol in that it is an open shape, and when writing quickly by hand, it's okay if it morphs into a $\circ$.
But $\tri$ highlights the asymmetry of composition, in contrast with the other monoidal structures on $\poly$ we've encountered.
Moreover, we'll soon see that $\tri$ is quite evocative in terms of trees.
For each $n\in\nn$, we'll also use $p\tripow{n}$ to denote the $n$-fold composition product of $p$, i.e.\ $n$ copies of $p$ all composed with each other.\footnote{When we say ``the $n$-fold composition product of $p$,'' we mean $n$ copies of $p$ all composed with each other; but when we discuss an ``$n$-fold composition product'' in general, we refer to an arbitrary composition product of $n$ polynomials that may or may not all be equal to each other. This will apply to composition products of lenses as well, once we define those.}
In particular, $p\tripow0=\yon$ and $p\tripow1=p$.

We repeat the important formulas from \cref{prop.poly_closed_comp} and its proof in the new notation:
\begin{equation}\label{eqn.composite_formula}
p\tri q\iso\sum_{i\in p(\1)}\prod_{a\in p[i]}\sum_{j\in q(\1)}\prod_{b\in q[j]}\yon.
\end{equation}

\begin{align}\label{eqn.composite_formula_sums_first}
  p\tri q\iso
  \scalebox{1.3}{$\displaystyle
  \sum_{i\in p(\1)} \; \sum_{\ol{j}\colon p[i]\to q(\1)}\yon^{\sum_{a\in p[i]}q[\ol{j}(a)]}$}
\end{align}

% \[
% \begin{tikzpicture}[polybox, baseline=(helper)]
% 	\node[poly] (p) {$a:p[i]$\at$i:p(\1)$};
% 	\node[poly, above=of p] (q) {$b:q[j]$\at$j:q(\1)$};
% 	\coordinate (helper) at ($(p.north)!.5!(q.south)$);
% \end{tikzpicture}
% \quad\cong\quad
% \begin{tikzpicture}[polybox, baseline=(p.east)]
% 	\node[poly] (p) {$(a:p[i], b:q[j(a)]$)\at$(i:p(1), j: p[i]\to q(\1))$};
% \end{tikzpicture}
% \]

\begin{exercise}
Let's consider \eqref{eqn.composite_formula_sums_first} piece by piece, with concrete polynomials $p\coloneqq\yon^\2+\yon^\1$ and $q\coloneqq \yon^\3+\1$.
\begin{enumerate}
	\item What is $\yon^\2\tri q$?
	\item What is $\yon^\1\tri q$?
	\item What is $(\yon^\2+\yon^\1)\tri q$? This is what $p\tri q$ ``should be.''
	\item How many functions $\ol{j_1}\colon p[1]\to q(\1)$ are there?
	\item For each function $\ol{j_1}$ as above, what is $\sum_{a\in p[1]} q[\ol{j_1}(a)]$?
	\item How many functions $\ol{j_2}\colon p[2]\to q(\1)$ are there?
	\item For each function $\ol{j_2}$ as above, what is $\sum_{a\in p[2]} q[\ol{j_2}(a)]$?
	\item Write out \[\sum_{i\in p(\1)}\;\sum_{\ol{j}\colon p[i]\to q(\1)}\yon^{\sum_{a\in p[i]}q[\ol{j}(a)]}.\]
	Does the result agree with what $p\tri q$ should be?
\qedhere
\end{enumerate}
\begin{solution}
We are given $p\coloneqq\yon^\2+\yon^\1$ and $q\coloneqq \yon^\3+\1$.
\begin{enumerate}
    \item By standard polynomial multiplication, we have that $\yon^\2 \tri q \iso q \times q \iso \yon^\6 + \2\yon^\3 + \1$.
    \item We have that $\yon^\1 \tri q \iso q \iso \yon^\3 + \1$.
    \item Combining the previous parts, we have that $(\yon^\2 + \yon^\1) \tri q \iso q \times q + q \iso \yon^\6 + \3\yon^\3 + \2$.
    \item Since $p[1] \iso \2$ and $q(\1) \iso \2$, there are $2^2 = 4$ functions $p[1] \to q(\1)$.
    \item When $\ol{j_1} \colon p[1] \to q(\1)$ is one of the two possible bijections, we have that
    \[
        \sum_{a \in p[1]} q[\ol{j_1}(a)] \iso q[1] + q[2] \iso \3 + \0 \iso \3.
    \]
    When $\ol{j_1} \colon p[1] \to q(\1)$ sends everything to $1 \in q(\1)$, we have that
    \[
        \sum_{a \in p[1]} q[\ol{j_1}(a)] \iso q[1] + q[1] \iso \3 + \3 \iso \6.
    \]
    Finally, when $\ol{j_1} \colon p[1] \to q(\1)$ sends everything to $2 \in q(\1)$, we have that
    \[
        \sum_{a \in p[1]} q[\ol{j_1}(a)] \iso q[2] + q[2] \iso \0 + \0 \iso \0.
    \]
    \item Since $p[2] \iso \1$ and $q(\1) \iso \2$, there are $2^1 = 2$ functions $p[2] \to q(\1)$.
    \item When $j_2 \colon p[2] \to q(\1)$ maps to $1 \in q(\1)$, we have that $\sum_{a \in p[2]} q[\ol{j_2}(a)] \iso q[1] \iso \3$, and when $\ol{j_2} \colon p[2] \to q(\1)$ maps to $2 \in q(\1)$, we have that $\sum_{a \in p[2]} q[\ol{j_2}(a)] \iso q[2] \iso \0$.
    \item From the previous parts, we have that
    \[
        \sum_{i\in p(\1)}\;\sum_{\ol{j}\colon p[i]\to q(\1)}\yon^{\sum_{a\in p[i]}q[j_i(a)]} \iso (\2\yon^\3 + \yon^\6 + \yon^\0) + (\yon^\3 + \yon^\0) \iso \yon^\6 + \3\yon^\3 + \2,
    \]
    which agrees with what $p \tri q$ should be.
\end{enumerate}
\end{solution}
\end{exercise}

\index{composition product!of special polynomials}

\begin{exercise}\label{exc.composites_of_specials}
\begin{enumerate}
	\item If $p$ and $q$ are representable, show that $p\tri q$ is too. Give a formula for it.
	\item If $p$ and $q$ are linear, show that $p\tri q$ is too. Give a formula for it.
	\item If $p$ and $q$ are constant, show that $p\tri q$ is too. Give a formula for it.
\qedhere
\end{enumerate}
\begin{solution}
\begin{enumerate}
	\item Given representable polynomials $p \coloneqq \yon^A$ and $q \coloneqq \yon^B$, we have that $p \tri q \iso \left(\yon^B\right)^A \iso \yon^{AB}$, which is also representable.
	\item Given linear polynomials $p \coloneqq A\yon$ and $q \coloneqq B\yon$, we have that $p \tri q \iso A(B\yon) \iso AB\yon$, which is also linear.
	\item Given constant polynomials $p \coloneqq A$ and $q \coloneqq B$, we have that $p \tri q \iso A$, which is also constant (see also \cref{exc.composing_with_constants}).\index{constant polynomial}
\end{enumerate}
\end{solution}
\end{exercise}

\begin{exercise}
Recall the closure operation that we denoted $\ihom{-,-}\colon\poly\op\times\poly\to\poly$ for $\otimes$ from \eqref{eqn.par_hom}.
Show that for all $A\in\smset$ and $q\in\poly$, there is an isomorphism
\[
    \yon^A\tri q\iso\ihom{A\yon,q}.
\]
\begin{solution}
Given $A\in\smset$ and $q\in\poly$, we have
\begin{align*}
    \yon^A\tri q    &\iso
	\sum_{\ol{j}\colon A\to q(\1)}\yon^{\sum_{a\in A}q[\ol{j}(a)]}
	\tag*{\eqref{eqn.composite_formula_sums_first}} \\ &\iso
	\sum_{\varphi\colon A\yon\to q}\yon^{\sum_{a\in A}q[\varphi_\1(a)]}   \\&\iso
	\ihom{A\yon,q},
	\tag*{\eqref{eqn.par_hom_sum}}
\end{align*}
for a lens $\varphi\colon A\yon\to q$ has an on-positions function $A\to q(\1)$ and uniquely determined on-directions functions.
\end{solution}
\end{exercise}

We know how $\tri$ acts on the objects in $\poly$, but what does it do to the morphisms between them?
For any pair of natural transformations $f\colon p\to p'$ and $g\colon q\to q'$ between polynomial functors, their composite $f\tri g\colon p\tri q\to p'\tri q'$ is given by \emph{horizontal composition}.

\index{composition product!on morphisms|see{composition product, on lenses}}
\index{composition product!on lenses|(}

\begin{definition}[Horizontal composition of natural transformations]\label{def.horiz_comp_nat_trans}
Let $f\colon p\to p'$ and $g\colon q\to q'$ be two natural transformations between (polynomial) functors $p,p',q,q'\colon\smset\to\smset$.
Then the \emph{horizontal composite} of $f$ and $g$, denoted $f\tri g$, is the natural transformation $p\tri q\to p'\tri q'$ whose $X$-component for each $X\in\smset$ is the function
\begin{equation} \label{eqn.horiz_comp_nat_trans_comp}
    p(q(X)) \To{f_{q(X)}} p'(q(X)) \To{p'(g_X)} p'(q'(X))
\end{equation}
obtained by composing the $q(X)$-component of $f$ with the functor $p'$ applied to the $X$-component of $g$.
\end{definition}

\begin{exercise}
Show that we could have replaced the composite function \eqref{eqn.horiz_comp_nat_trans_comp} in \cref{def.horiz_comp_nat_trans} with the function
\begin{equation} \label{eqn.horiz_comp_nat_trans_comp2}
    p(q(X)) \To{p(g_X)} p(q'(X)) \To{f_{q'(X)}} p'(q'(X))
\end{equation}
obtained by composing $p$ applied to the $X$-component of $g$ with the $q'(X)$-component of $f$, without altering the definition.
\begin{solution}
We wish to show that \eqref{eqn.horiz_comp_nat_trans_comp2} could replace \eqref{eqn.horiz_comp_nat_trans_comp} in \cref{def.horiz_comp_nat_trans}.
We claim that \eqref{eqn.horiz_comp_nat_trans_comp} and \eqref{eqn.horiz_comp_nat_trans_comp2} are in fact the same function; that is, that the following square commutes:
\[
\begin{tikzcd}
    p(q(X)) \ar[r, "f_{q(X)}"]\ar[d, "p(g_X)"'] & p'(q(X)) \ar[d, "p'(g_X)"] \\
    p(q'(X)) \ar[r, "f_{q'(X)}"'] & p'(q'(X))
\end{tikzcd}
\]
Indeed it does, by the naturality of $f$.
\end{solution}
\end{exercise}

\index{natural transformation!horizontal composition}
\index{natural transformation!vertical composition}

\begin{remark}
There are two very different notions of lens composition floating around, so we'll try to mitigate confusion by standardizing terminology here.
We'll reserve the term \emph{composite lens} for lenses $h\then j\colon r\to t$ obtained by composing a lens $h\colon r\to s$ with a lens $j\colon s\to t$, according to the composition rule of the category $\poly$.
This corresponds to \emph{vertical composition} of natural transformations.
This is also the kind of composition we will mean whenever we use the verb ``\emph{compose},'' if the objects of that verb are lenses.

Meanwhile, we'll use the term \emph{composition product (of lenses)} for lenses $f\tri g\colon p\tri q\to p'\tri q'$ obtained by applying the monoidal product functor $\tri\colon\poly\times\poly\to\poly$ on the lenses $f\colon p\to p'$ and $g\colon q\to q'$.
This corresponds to \emph{horizontal composition} of natural transformations.
In this case, we'll use the verb phrase ``\emph{taking the monoidal product}.''

On the other hand, we'll use the terms ``composite'' and ``composition product'' interchangeably to refer to polynomials $p\tri q$, obtained by composing $p,q\in\poly$ as functors or, equivalently, applying the monoidal product functor $\tri$ on them---as there is no risk of confusion here.

This is another reason we tend to avoid the symbol $\circ$, preferring to use $\then$ for vertical composition and $\tri$ for horizontal composition.
Of course, if you're ever confused, you can always check whether the codomain of the first lens matches up with the domain of the second.
If they don't, we must be taking their monoidal product.
\end{remark}

The composition product of polynomials and lenses will be extremely important in the story that follows.
However, we only sometimes think of it as the composition of functors and the horizontal composition of natural transformations; more often we think of it as certain operations on positions and directions or on corolla forests.

\index{composition product!on lenses|)}\index{natural transformation!horizontal composition}

\subsection{Composite positions and directions}\label{subsec.comon.comp.def.arena}

\index{composition product!positions and directions}

Let us interpret our formula \eqref{eqn.composite_formula_sums_first} for the composition product of two polynomials in terms of positions and directions.
The position-set of $p \tri q$ is
\begin{equation} \label{eqn.comp_pos}
    (p \tri q)(\1) \iso \sum_{i \in p(\1)} \; \sum_{\ol{j} \colon p[i] \to q(1)} \1 \iso \sum_{i \in p(\1)} \smset(p[i], q(\1)).
\end{equation}
In other words, specifying a position of $p \tri q$ amounts to first specifying a $p$-position $i$, then specifying a function $\ol{j} \colon p[i] \to q(\1)$, i.e.\ a $q$-position $\ol{j}(a)$ for each $p[i]$-direction $a$.

Given such a position $(i, \ol{j})$ of $p \tri q$, the direction-set of $p \tri q$ at $(i, \ol{j})$ is
\begin{equation} \label{eqn.comp_dir}
    (p \tri q)[(i, \ol{j})] \iso \sum_{a \in p[i]} q[\ol{j}(a)].
\end{equation}
So a direction of $p \tri q$ at $(i, \ol{j})$ consists of a $p[i]$-direction $a$ and a $q[\ol{j}(a)]$-direction.

While this description completely characterizes $p \tri q$, it may be a bit tricky to wrap your head around.
Here is an alternative perspective that can help us get a better intuition for what's going on with the composition product of polynomials.

Back in \cref{sec.poly.rep-sets.sum-prod-set}, we saw how to write the instructions for choosing an element of a sum or product of sets.
For instance, given a polynomial $p$ and a set $X$, the instructions for choosing an element of
\[
    p\tri X=p(X)\iso\sum_{i\in p(\1)}\prod_{a\in p[i]}X
\]
would be written as follows.
\begin{quote}
To choose an element of $p(X)$:
\begin{enumerate}
    \item choose an element $i\in p(\1)$;
    \item for each element $a\in p[i]$:
    \begin{enumerate}[label*=\arabic*.]
        \item choose an element of $X$.
    \end{enumerate}
\end{enumerate}
\end{quote}
But say we hadn't picked a set $X$ yet; in fact, say we might replace $X$ with a general polynomial instead.
We'll replace ``an element of $X$'' with a placeholder---the words ``a future''---that indicates that we don't yet know what will go there. It depends on what has come before.%
\footnote{In a significant sense, the composition product should be thought of as being about dependency.}
Furthermore, to highlight that these instructions are associated with some polynomial $p$, we will use our familiar positions and directions terminology.
\begin{quote}
The instructions associated with a polynomial $p$ are:
\begin{enumerate}
    \item choose a $p$-position $i$;
    \item for each $p[i]$-direction $a$:
    \begin{enumerate}[label*=\arabic*.]
        \item choose a future.
    \end{enumerate}
\end{enumerate}
\end{quote}

\index{future!as placeholder for dependency}
\index{composition product!as dependency}

If we think of polynomials in terms of their instructions, then \eqref{eqn.composite_formula} tells us that the composition product simply nests one set of instructions within another, as follows.
\begin{quote}
The instructions associated with a polynomial $p\tri q$ are:
\begin{enumerate}
    \item choose a $p$-position $i$;
    \item for each $p[i]$-direction $a$:
    \begin{enumerate}[label*=\arabic*.]
        \item choose a $q$-position $j$;
        \item for each $q[j]$-direction $b$:
        \begin{enumerate}[label*=\arabic*.]
            \item choose a future.
        \end{enumerate}
    \end{enumerate}
\end{enumerate}
\end{quote}
Similarly, we could write down the instructions associated with any $n$-fold composition product by nesting even further.
We might think of such instructions as specifying some sort of length-$n$ \emph{strategy}, in the sense of game theory, for picking positions given any directions---except that the opponent is somehow abstract, having no positions of its own.

\index{composition product!$n$-fold}
\index{composition product!positions as strategies}

When we rewrite \eqref{eqn.composite_formula} \eqref{eqn.composite_formula_sums_first}, we are collapsing the instructions down into the following, highlighting the positions and directions of $p\tri q$.
\begin{quote}
The instructions associated with a polynomial $p\tri q$ are:
\begin{enumerate}
    \item choose a $p$-position $i$ and, for each $p[i]$-direction $a$, a $q$-position $\ol{j}_i(a)$;
    \item for each $p[i]$-direction $a$ and each $q[\ol{j}_i(a)]$-direction $b$:
    \begin{enumerate}[label*=\arabic*.]
        \item choose a future.
    \end{enumerate}
\end{enumerate}
\end{quote}
We will see in \cref{subsec.comon.comp.def.corolla} that these instructions have a very natural interpretation when we depict these polynomials as corolla forests.
\index{future!as placeholder for dependency}

\begin{exercise}
\begin{enumerate}
	\item Let $p$ be an arbitrary polynomial. Write out the (uncollapsed) instructions associated with $p\tripow3=p\tri p\tri p$.
	\item Write out the (uncollapsed) instructions for choosing an element of $p\tri p\tri\1$, but where you would normally write ``choose an element of $\1$,'' just write ``done.'' \qedhere
\end{enumerate}
\begin{solution}
\begin{longenum}
    \item The instructions associated with a polynomial $p\tri p\tri p$ are:
    \begin{enumerate}
        \item choose a $p$-position $i$;
        \item for each $p[i]$-direction $a$:
        \begin{enumerate}[label*=\arabic*.]
            \item choose a $p$-position $i'$;
            \item for each $p[i']$-direction $a'$:
            \begin{enumerate}[label*=\arabic*.]
                \item choose a $p$-position $i''$;
                \item for each $p[i'']$-direction $a''$:
                \begin{enumerate}[label*=\arabic*.]
                    \item choose a future.
                \end{enumerate}
            \end{enumerate}
        \end{enumerate}
    \end{enumerate}
    \item To choose an element of $p\tri p\tri\1$:
    \begin{enumerate}
        \item choose a $p$-position $i$;
        \item for each $p[i]$-direction $a$:
        \begin{enumerate}[label*=\arabic*.]
            \item choose a $p$-position $i'$;
            \item for each $p[i']$-direction $a'$:
            \begin{enumerate}[label*=\arabic*.]
                \item done.
            \end{enumerate}
        \end{enumerate}
    \end{enumerate}
\end{longenum}

\end{solution}
\end{exercise}

\index{future!as placeholder for dependency}

But how does the composition product act on lenses?
Given lenses $f\colon p\to p'$ and $g\colon q\to q'$, we can translate them to natural transformations, take their horizontal composite, then translate this back to a lens.
The following exercise guides us through this process.

\index{composition product!on lenses}

\begin{exercise}[The composition product of lenses] \label{exc.comp_prod_lens}
Fix lenses $f\colon p\to p'$ and $g\colon q\to q'$.
We seek to characterize their composition product $f\tri g\colon p\tri q\to p'\tri q'$.
\begin{enumerate}
    \item\label{exc.comp_prod_lens.1} Use \cref{prop.morph_arena_to_func} to compute the $q(X)$-component of $f$ as a natural transformation.
    \item\label{exc.comp_prod_lens.2} Use \cref{prop.poly_on_functions,prop.morph_arena_to_func} to compute $p'$ applied to the $X$-component of $g$ as a natural transformation.
    \item\label{exc.comp_prod_lens.3} Combine \cref{exc.comp_prod_lens.1} and \cref{exc.comp_prod_lens.2} using \cref{def.horiz_comp_nat_trans} to compute the horizontal composite $f\tri g$ of $f$ and $g$ as natural transformations.
    \item Use \cref{cor.morph_func_to_arena} to translate the natural transformation $f\tri g$ obtained in \cref{exc.comp_prod_lens.3} to a lens $p\tri q\to p'\tri q'$.
    Verify that for each $(i,\ol{j}_i)$ in $(p\tri q)(\1)$ (see \eqref{eqn.comp_pos}), its on-positions function sends
    \begin{equation} \label{eqn.comp_lens_pos}
        (i,\ol{j}_i)\Mapsto{(f\:\tri\:g)_\1}\left(f_\1(i), f^\sharp_i\then\ol{j}_i\then g_\1\right);
    \end{equation}
    while for each $(a',b')$ in $(p'\tri q')[(f_\1(i), f^\sharp_i\then\ol{j}_i\then g_\1)]$ (see \eqref{eqn.comp_dir}), its on-directions function sends
    \begin{equation} \label{eqn.comp_lens_dir}
        (a',b')\Mapsto{(f\:\tri\:g)^\sharp_{(i,\ol{j}_i)}}\left(f^\sharp_i(a'), g^\sharp_{\ol{j}_i(f^\sharp_i(a'))}(b')\right).
    \end{equation}
    \qedhere
\end{enumerate}
\begin{solution}
We have lenses $f\colon p\to p'$ and $g\colon q\to q'$.
\begin{enumerate}
    \item By \cref{prop.morph_arena_to_func}, the $q(X)$-component of $f$ is a function $f_{q(X)}\colon p(q(X))\to p'(q(X))$ that sends every $(i,h)$ with $i\in p(\1)$ and $h\colon p[i]\to q(X)$ to $(f_\1(i),f^\sharp_i\then h)$.
    We can think of the function $h\colon p[i]\to q(X)$ equivalently as a function $\ol{j}_i\colon p[i]\to q(\1)$ and, for each $a\in p[i]$, a function $h_a\colon q[\ol{j}_i(a)]\to X$.
    So $f_{q(X)}\colon (p\tri q)(X)\to (p'\tri q)(X)$ sends \[(i,\ol{j}_i,(h_a)_{a\in p[i]})\mapsto\left(f_\1(i),f^\sharp_i\then\ol{j}_i,\left(h_{f^\sharp_i(a')}\right)_{a'\in p'[f_\1(i)]}\right).\]

    \item By \cref{prop.morph_arena_to_func}, the $X$-component of $g$ is a function $g_X\colon q(X)\to q'(X)$ that sends every $(j,k)$ with $j\in q(\1)$ and $k\colon q[j]\to X$ to $(g_\1(j),g^\sharp_j\then k)$ in $q'(X)$.
    Then by \cref{prop.poly_on_functions}, applying $p'$ to this $X$-component yields a function $p'(q(X))\to p'(q'(X))$ that sends every $(i',\ol{j'}_{i'},(h'_{a'})_{a'\in p'[i']})$ with $i'\in p'(\1)$ as well as $\ol{j'}_{i'}\colon p'[i']\to q(\1)$ and $h'_{a'}\colon q[\ol{j'}_{i'}(a')]\to X$ to \[\left(i',\ol{j'}_{i'}\then g_\1,\left(g^\sharp_{\ol{j'}_{i'}(a')}\then h'_{a'}\right)_{a'\in p'[i']}\right).\]

    \item By \cref{def.horiz_comp_nat_trans}, the horizontal composite of $f$ and $g$ is the natural transformation $f\tri g\colon p\tri p'\to q\tri q'$ whose $X$-component is the composite of the answers to \cref{exc.comp_prod_lens.1} and \cref{exc.comp_prod_lens.2}, sending
    \begin{align*}
        (i,\ol{j}_i,(h_a)_{a\in p[i]})&\mapsto\left(f_\1(i),f^\sharp_i\then\ol{j}_i,\left(h_{f^\sharp_i(a')}\right)_{a'\in p'[f_\1(i)]}\right)\\
        &\mapsto\left(f_\1(i),f^\sharp_i\then\ol{j}_i\then g_\1, \left(g^\sharp_{\ol{j}_{i}(f^\sharp_i(a'))}\then h_{f^\sharp_i(a')}\right)_{a'\in p'[f_\1(i)]}\right).
    \end{align*}

    \item We use \cref{cor.morph_func_to_arena} to translate the answer to \cref{exc.comp_prod_lens.3} into a lens $f\tri g\colon p\tri q\to p'\tri q'$, as follows.
    Its on-positions function is the $\1$-component $(f\tri g)_\1$, which sends every $(i,\ol{j}_i)$ with $i\in p(\1)$ and $\ol{j}_i\colon p[i]\to q(\1)$ to
    \[
        (f_\1(i),f^\sharp_i\then\ol{j}_i\then g_\1).
    \]
    Then for each such $(i,\ol{j}_i)$, if we apply the $(p\tri q)[(i,\ol{j}_i)]$-component of $f\tri g$ to the element $(i,\ol{j}_i,(\iota_d)_{a\in p[i]})$, where $\iota_d\colon q[\ol{j}_i(a)]\to(p\tri q)[(i,\ol{j}_i)]\iso\sum_{a\in p[i]}q[\ol{j}_i(a)]$ is the canonical inclusion, then take the last coordinate of the result, we obtain for each $a'\in p'[f_\1(i)]$ the function
    \[
        q'[g_\1(\ol{j}_i(f^\sharp_i(a')))] \To{g^\sharp_{\ol{j}_{i}(f^\sharp_i(a'))}} q[\ol{j}_i(f^\sharp_i(a'))] \To{\iota_{f^\sharp_i(a')}} \sum_{a\in p[i]}q[\ol{j}_i(a)] \iso (p\tri q)[(i,\ol{j}_i)].
    \]
    These can equivalently be thought of as a single function from
    \[
        \sum_{a'\in p'[f_\1(i)]} q'[g_\1(\ol{j}_i(f^\sharp_i(a')))] \iso (p'\tri q')[(f\tri g)_\1(i,\ol{j}_i)]
    \]
    which \cref{cor.morph_func_to_arena} tells us is the on-directions function of $f\tri g$ at $(i,\ol{j}_i)$, that sends every $(a',b')$ with $a'\in p'[f_\1(i)]$ and $b'\in q'[g_\1(\ol{j}_i(f^\sharp_i(a')))]$ to
    \[
        \left(f^\sharp_i(a'), g^\sharp_{\ol{j}_i(f^\sharp_i(a'))}(b')\right).
    \]
\end{enumerate}
\end{solution}
\end{exercise}

So what does \cref{exc.comp_prod_lens} tell us about the behavior of $f\tri g\colon p\tri q\to p'\tri q'$?
By \eqref{eqn.comp_lens_pos}, on positions, $f\tri g$ takes a $p$-position $i$ and sends it to the $p'$-position $f_\1(i)$; then for each direction $a'$ at this position, the associated $q'$-position is obtained by sending $a'$ back to a $p[i]$-direction via $f^\sharp_i$, checking what $q$-position is associated to that $p[i]$-direction via some $\ol{j}_i$, then sending that $q$-position forward again to a $q'$-position via $g_\1$.

Then by \eqref{eqn.comp_lens_pos}, on directions, $f\tri g$ sends a direction of $p'$ back to a direction of $p$ via an on-directions function of $f$, then sends a direction of $q'$ back to a direction of $q$ via an on-directions funtion of $g$.
We'll get a better sense of what's happening when we see this drawn out as corolla forests in \cref{ex.comp_prod_trees}.

\subsection{Composition product on corolla forests} \label{subsec.comon.comp.def.corolla}

\index{composition product!grafting corollas}

It turns out that the forest of $p\tri q$ is given by grafting $q$-corollas onto the leaves of $p$-corollas in every possible way.
We will demonstrate this using an example.

Let $p\coloneqq\yon^\2+\yon$ and $q\coloneqq\yon^\3+\1$, whose corolla forests we draw as follows:\index{tree!depiction of}
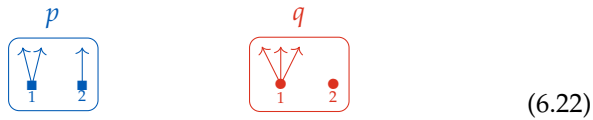
\begin{equation}\label{eqn.pq_misc39}
\begin{tikzpicture}[rounded corners]
	\node (p1) [draw, my-blue, "\color{my-blue} $p$" above] {
	\begin{tikzpicture}[trees, sibling distance=2.5mm]
    \node["\tiny 1" below] (1) {\tiny$\blacksquare$}
      child {}
      child {};
    \node[right=.5 of 1,"\tiny 2" below] (2) {\tiny$\blacksquare$}
      child {};
  \end{tikzpicture}
  };
	\node (p2) [draw, my-red, right=2 of p1, "\color{my-red} $q$" above] {
	\begin{tikzpicture}[trees, sibling distance=2.5mm]
    \node["\tiny 1" below] (1) {$\bullet$}
      child {}
      child {}
      child {};
    \node[right=.5 of 1,"\tiny 2" below] (4) {$\bullet$}
    ;
  \end{tikzpicture}
  };
\end{tikzpicture}
\end{equation}
By \eqref{eqn.comp_pos}, choosing a position of $p \tri q$ amounts to first choosing a $p$-root $i$, then choosing a $q$-root for every $p[i]$-leaf.
So we may depict $(p \tri q)(\1)$ by grafting roots from the corolla forest of $q$ to leaves in the corolla forest of $p$ in every possible way, as follows:
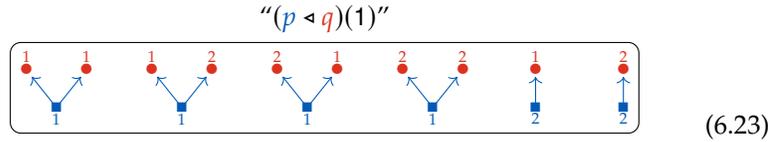
\begin{equation}\label{eqn.comp_pos_forest}
\begin{tikzpicture}[rounded corners]
	\node (p1) [draw, "``$(${\color{my-blue} $p$}$\:\tri\:${\color{my-red}$q$}$)(\1)$''" above] {
	\begin{tikzpicture}[trees,
		level 1/.style={sibling distance=8mm},
	  my-blue]
    \node[my-blue, "\tiny 1" below] (1) {\tiny$\blacksquare$}
      child[my-blue] {node[my-red, "\color{my-red} \tiny 1" above] {$\bullet$}}
      child[my-blue] {node[my-red, "\color{my-red} \tiny 1" above] {$\bullet$}};
    \node[my-blue, right=1.5 of 1, "\tiny 1" below] (2) {\tiny$\blacksquare$}
      child[my-blue] {node[my-red, "\color{my-red} \tiny 1" above] {$\bullet$}}
      child[my-blue] {node[my-red, "\color{my-red} \tiny 2" above] {$\bullet$}};
    \node[my-blue, right=1.5 of 2, "\tiny 1" below] (3) {\tiny$\blacksquare$}
      child[my-blue] {node[my-red, "\color{my-red} \tiny 2" above] {$\bullet$}}
      child[my-blue] {node[my-red, "\color{my-red} \tiny 1" above] {$\bullet$}};
    \node[my-blue, right=1.5 of 3, "\tiny 1" below] (4) {\tiny$\blacksquare$}
      child[my-blue] {node[my-red, "\color{my-red} \tiny 2" above] {$\bullet$}}
      child[my-blue] {node[my-red, "\color{my-red} \tiny 2" above] {$\bullet$}};
    \node[my-blue, right=1.2 of 4, "\tiny 2" below] (5) {\tiny$\blacksquare$}
      child[my-blue] {node[my-red, "\color{my-red} \tiny 1" above] {$\bullet$}};
    \node[my-blue, right=1 of 5, "\tiny 2" below] (6) {\tiny$\blacksquare$}
      child[my-blue] {node[my-red, "\color{my-red} \tiny 2" above] {$\bullet$}};
  \end{tikzpicture}
  };
\end{tikzpicture}
\end{equation}
Now fix one of the positions of $p \tri q$ drawn above: a $p$-root $i$ and a $q$-root grafted to every $p[i]$-leaf.
By \eqref{eqn.comp_dir}, a direction of $p \tri q$ at that position consists of a $p[i]$-leaf $a$ and a second leaf emanating from the $q$-root that has been grafted on	to $a$.
In other words, in the following picture, where we have grafted not just $q$-roots but entire $q$-corollas onto leaves in $p$, the directions of $p \tri q$ at the position corresponding to each tree are the rooted paths\footnote{A \emph{rooted path} of a rooted tree is a path up the tree that starts from the root.} of that tree of length $2$ (we omit the labels):
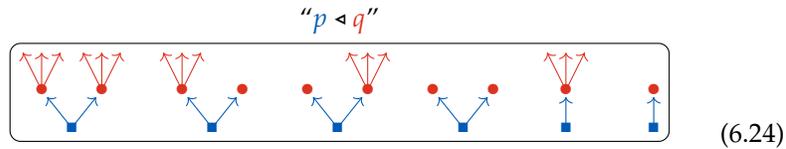
\begin{equation}\label{eqn.prefered_composite}
\begin{tikzpicture}[rounded corners]
	\node (p1) [draw, "``{\color{my-blue} $p$}$\:\tri\:${\color{my-red}$q$}''" above] {
	\begin{tikzpicture}[trees,
		level 1/.style={sibling distance=8mm},
	  level 2/.style={sibling distance=2.5mm},
	  my-blue]
    \node[my-blue] (1) {\tiny$\blacksquare$}
      child {node[my-red] {$\bullet$}
      	child[my-red]
				child[my-red]
				child[my-red]
			}
      child {node[my-red] {$\bullet$}
      	child[my-red]
				child[my-red]
				child[my-red]
			};
    \node[my-blue, right=1.7 of 1] (2) {\tiny$\blacksquare$}
      child {node[my-red] {$\bullet$}
      	child[my-red]
				child[my-red]
				child[my-red]
			}
      child {node[my-red] {$\bullet$}
			};
    \node[my-blue, right=1.5 of 2] (3) {\tiny$\blacksquare$}
      child {node[my-red] {$\bullet$}
			}
      child {node[my-red] {$\bullet$}
      	child[my-red]
				child[my-red]
				child[my-red]
			};
    \node[my-blue, right=1.5 of 3] (4) {\tiny$\blacksquare$}
      child {node[my-red] {$\bullet$}
			}
      child {node[my-red] {$\bullet$}
			};
    \node[my-blue, right=1.2 of 4] (5) {\tiny$\blacksquare$}
      child {node[my-red] {$\bullet$}
      	child[my-red]
				child[my-red]
				child[my-red]
			};
    \node[my-blue, right=1 of 5] (6) {\tiny$\blacksquare$}
      child {node[my-red] {$\bullet$}
			};
  \end{tikzpicture}
  };
\end{tikzpicture}
\end{equation}
Equivalently, we can think of the directions in the picture above as the leaves at the second level of each tree.
So $p \tri q$ has six positions; the first has six directions, the second, third, and fifth have three directions, and the fourth and sixth have no directions.
In total, we can read off that $p\tri q$ is isomorphic to $\yon^\6+\3\yon^\3+\2$.

\index{tree!multi-level}

We put the $p\tri q$ in scare quotes above \eqref{eqn.prefered_composite} because, to be pedantic, the corolla forest of $p \tri q$ has the two levels smashed together as follows:
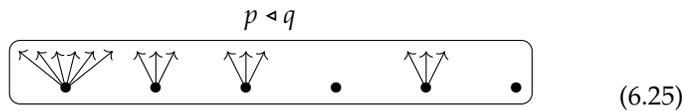
\begin{equation}\label{eqn.actual_composite}
\begin{tikzpicture}[rounded corners]
	\node (p1) [draw, "$p\tri q$" above] {
	\begin{tikzpicture}[trees, sibling distance=2.5mm]
    \node (1) {$\bullet$}
      child {}
      child {}
      child {}
      child {}
      child {}
      child {};
    \node[right=1 of 1] (2) {$\bullet$}
      child {}
      child {}
      child {};
    \node[right=1 of 2] (3) {$\bullet$}
      child {}
      child {}
      child {};
    \node[right=1 of 3] (4) {$\bullet$};
    \node[right=1 of 4] (5) {$\bullet$}
      child {}
      child {}
      child {};
    \node[right=1 of 5] (6) {$\bullet$};
  \end{tikzpicture}
  };
\end{tikzpicture}
\end{equation}
Usually, we will prefer the style of \eqref{eqn.prefered_composite} rather than the more pedantic style of \eqref{eqn.actual_composite}.

We have now seen how to draw a single polynomial as a corolla forest, with height-$1$ leaves as directions; as well as how to draw a two-fold composite of polynomials as a forest of trees, with height-$2$ leaves as directions.
Note that drawing a corolla of $p$ or a tree of $p\tri q$ is just a graphical way of following the instructions associated with the polynomial $p$ or $p\tri q$ that we saw in \cref{subsec.comon.comp.def.arena}, where the arrows---the top-level leaves---are where the ``futures'' would go.
Similarly, we could depict any $n$-fold composite as a forest with height-$n$ leaves as directions.
You will have an opportunity to try this in the following exercise.

\index{future!as placeholder for dependency}

\begin{exercise}
Use $p,q$ as in \eqref{eqn.pq_misc39} and $r\coloneqq \2\yon+\1$ in the following.
\begin{enumerate}
	\item Draw $q\tri p$.
	\item Draw $p\tri p$.
	\item Draw $p\tri p\tri \1$.
	\item Draw $r\tri r$.
	\item Draw $r\tri r\tri r$.
\qedhere
\end{enumerate}
\begin{solution}
We have $p \coloneqq \yon^\2 + \yon$ and $q \coloneqq \yon^\3 + \1$ as in \eqref{eqn.pq_misc39}.
\begin{enumerate}
    \item Here is a picture of $q\tri p$, where each tree is obtained by taking a $q$-corolla and grafting $p$-corollas to every leaf:
\[
\begin{tikzpicture}[rounded corners]
	\node (p1) [draw] {
	\begin{tikzpicture}[trees,
		level 1/.style={sibling distance=4mm},
	  level 2/.style={sibling distance=2.5mm},
	  my-red]
    \node (1) {$\bullet$}
      child {node[my-blue] {\tiny$\blacksquare$}
      	child[my-blue]
				child[my-blue]
			}
      child {node[my-blue] {\tiny$\blacksquare$}
      	child[my-blue]
				child[my-blue]
			}
      child {node[my-blue] {\tiny$\blacksquare$}
      	child[my-blue]
				child[my-blue]
			};
    \node[right=1.5 of 1] (2) {$\bullet$}
      child {node[my-blue] {\tiny$\blacksquare$}
      	child[my-blue]
				child[my-blue]
			}
      child {node[my-blue] {\tiny$\blacksquare$}
      	child[my-blue]
				child[my-blue]
			}
      child {node[my-blue] {\tiny$\blacksquare$}
      	child[my-blue]
			};
    \node[right=1.5 of 2] (3) {$\bullet$}
      child {node[my-blue] {\tiny$\blacksquare$}
      	child[my-blue]
				child[my-blue]
			}
      child {node[my-blue] {\tiny$\blacksquare$}
      	child[my-blue]
			}
      child {node[my-blue] {\tiny$\blacksquare$}
      	child[my-blue]
				child[my-blue]
			};
    \node[right=1.5 of 3] (4) {$\bullet$}
      child {node[my-blue] {\tiny$\blacksquare$}
      	child[my-blue]
				child[my-blue]
			}
      child {node[my-blue] {\tiny$\blacksquare$}
      	child[my-blue]
			}
      child {node[my-blue] {\tiny$\blacksquare$}
      	child[my-blue]
			};

    \node[right=1.5 of 4] (5) {$\bullet$}
      child {node[my-blue] {\tiny$\blacksquare$}
      	child[my-blue]
			}
      child {node[my-blue] {\tiny$\blacksquare$}
      	child[my-blue]
				child[my-blue]
			}
      child {node[my-blue] {\tiny$\blacksquare$}
      	child[my-blue]
				child[my-blue]
			};

    \node[below=1.2 of 1] (6) {$\bullet$}
      child {node[my-blue] {\tiny$\blacksquare$}
      	child[my-blue]
			}
      child {node[my-blue] {\tiny$\blacksquare$}
      	child[my-blue]
				child[my-blue]
			}
      child {node[my-blue] {\tiny$\blacksquare$}
      	child[my-blue]
			};

    \node[below=1.2 of 2] (7) {$\bullet$}
      child {node[my-blue] {\tiny$\blacksquare$}
      	child[my-blue]
			}
      child {node[my-blue] {\tiny$\blacksquare$}
      	child[my-blue]
			}
      child {node[my-blue] {\tiny$\blacksquare$}
      	child[my-blue]
				child[my-blue]
			};
    \node[below=1.2 of 3] (8) {$\bullet$}
      child {node[my-blue] {\tiny$\blacksquare$}
      	child[my-blue]
			}
      child {node[my-blue] {\tiny$\blacksquare$}
      	child[my-blue]
			}
      child {node[my-blue] {\tiny$\blacksquare$}
      	child[my-blue]
			};
    \node[below=1.2 of 4] (9) {$\bullet$};
  \end{tikzpicture}
  };
\end{tikzpicture}
\]

    \item Here is a picture of $p\tri p$:
\[
\begin{tikzpicture}[rounded corners]
	\node (p1) [draw] {
	\begin{tikzpicture}[trees,
		level 1/.style={sibling distance=8mm},
	  level 2/.style={sibling distance=3mm},
	  my-blue]
    \node[my-blue] (1) {\tiny$\blacksquare$}
      child {node {\tiny$\blacksquare$}
      	child
				child
			}
      child {node {\tiny$\blacksquare$}
      	child
				child
			};
    \node[my-blue, right=1.5 of 1] (2) {\tiny$\blacksquare$}
      child {node {\tiny$\blacksquare$}
      	child
				child
			}
      child {node {\tiny$\blacksquare$}
        child
			};
    \node[my-blue, right=1.5 of 2] (3) {\tiny$\blacksquare$}
      child {node {\tiny$\blacksquare$}
        child
			}
      child {node {\tiny$\blacksquare$}
      	child
				child
			};
    \node[my-blue, right=1.5 of 3] (4) {\tiny$\blacksquare$}
      child {node {\tiny$\blacksquare$}
        child
			}
      child {node {\tiny$\blacksquare$}
        child
			};
    \node[my-blue, right=1.2 of 4] (5) {\tiny$\blacksquare$}
      child {node {\tiny$\blacksquare$}
      	child
				child
			};
    \node[my-blue, right=1 of 5] (6) {\tiny$\blacksquare$}
      child {node {\tiny$\blacksquare$}
        child
			};
  \end{tikzpicture}
  };
\end{tikzpicture}
\]

    \item To obtain a picture of $p\tri p\tri\1$, we take our picture of $p\tri p$ and graft the single, leafless $\1$-root onto every (height-$2$) leaf:
\[
\begin{tikzpicture}[rounded corners]
	\node (p1) [draw] {
	\begin{tikzpicture}[trees,
		level 1/.style={sibling distance=8mm},
	  level 2/.style={sibling distance=3mm},
	  my-blue]
    \node[my-blue] (1) {\tiny$\blacksquare$}
      child {node {\tiny$\blacksquare$}
      	child {node[black] {$\bullet$}}
				child {node[black] {$\bullet$}}
			}
      child {node {\tiny$\blacksquare$}
      	child {node[black] {$\bullet$}}
				child {node[black] {$\bullet$}}
			};
    \node[my-blue, right=1.5 of 1] (2) {\tiny$\blacksquare$}
      child {node {\tiny$\blacksquare$}
      	child {node[black] {$\bullet$}}
				child {node[black] {$\bullet$}}
			}
      child {node {\tiny$\blacksquare$}
        child {node[black] {$\bullet$}}
			};
    \node[my-blue, right=1.5 of 2] (3) {\tiny$\blacksquare$}
      child {node {\tiny$\blacksquare$}
        child {node[black] {$\bullet$}}
			}
      child {node {\tiny$\blacksquare$}
      	child {node[black] {$\bullet$}}
				child {node[black] {$\bullet$}}
			};
    \node[my-blue, right=1.5 of 3] (4) {\tiny$\blacksquare$}
      child {node {\tiny$\blacksquare$}
        child {node[black] {$\bullet$}}
			}
      child {node {\tiny$\blacksquare$}
        child {node[black] {$\bullet$}}
			};
    \node[my-blue, right=1.2 of 4] (5) {\tiny$\blacksquare$}
      child {node {\tiny$\blacksquare$}
      	child {node[black] {$\bullet$}}
				child {node[black] {$\bullet$}}
			};
    \node[my-blue, right=1 of 5] (6) {\tiny$\blacksquare$}
      child {node {\tiny$\blacksquare$}
        child {node[black] {$\bullet$}}
			};
  \end{tikzpicture}
  };
\end{tikzpicture}
\]
\end{enumerate}

Now $r\coloneqq \2\yon+\1$. Before we draw the composites, here's a picture of $r$ itself, with different colors and shapes to distinguish the different positions:

\[
\begin{tikzpicture}[rounded corners]
	\node (p1) [draw] {
	\begin{tikzpicture}[trees, sibling distance=2.5mm]
    \node[my-blue] (1) {$\bullet$}
      child[my-blue];
    \node[right=.5 of 1, my-red] (2) {\tiny$\blacksquare$}
      child[my-red];
    \node[right=.5 of 2] (3) {$\diamond$};
  \end{tikzpicture}
  };
\end{tikzpicture}
\]

\begin{enumerate}[resume]
    \item Here is a picture of $r\tri r$:
\[
\begin{tikzpicture}[rounded corners]
	\node (p1) [draw] {
	\begin{tikzpicture}[trees,
		level 1/.style={sibling distance=8mm},
	  level 2/.style={sibling distance=2.5mm},
	  my-blue]
    \node (1) {$\bullet$}
      child {node {$\bullet$}
      	child
			};
    \node[right=.5 of 1] (2) {$\bullet$}
      child {node[my-red] {\tiny$\blacksquare$}
      	child[my-red]
			};
    \node[right=.5 of 2] (3) {$\bullet$}
      child {node[black] {$\diamond$}};
    \node[right=.5 of 3, my-red] (4) {\tiny$\blacksquare$}
      child[my-red] {node {\tiny$\blacksquare$}
        child};
    \node[right=.5 of 4, my-red] (5) {\tiny$\blacksquare$}
      child[my-red] {node[my-red] {\tiny$\blacksquare$}
      	child[my-red]
			};
    \node[right=.5 of 5, my-red] (6) {\tiny$\blacksquare$}
      child[my-red] {node[black] {$\diamond$}};
    \node[right=.5 of 6, black] (7) {$\diamond$};
  \end{tikzpicture}
  };
\end{tikzpicture}
\]

    \item Here is a picture of $r\tri r\tri r$:
\[
\begin{tikzpicture}[rounded corners]
	\node (p1) [draw] {
	\begin{tikzpicture}[trees,
		level 1/.style={sibling distance=8mm},
	  level 2/.style={sibling distance=2.5mm},
	  my-blue]
    \node (1) {$\bullet$}
      child {node {$\bullet$}
      	child {node {$\bullet$}
      	  child
      	      }
			};
    \node[right=.5 of 1] (2) {$\bullet$}
      child {node {$\bullet$}
      	child {node[my-red] {\tiny$\blacksquare$}
      	  child[my-red]
      	      }
			};
    \node[right=.5 of 2] (3) {$\bullet$}
      child {node {$\bullet$}
      	child {node[black] {$\diamond$}}
			};
    \node[right=.5 of 3] (4) {$\bullet$}
      child {node[my-red] {\tiny$\blacksquare$}
      	child[my-red] {node[my-blue] {$\bullet$}
      	  child[my-blue]
      	      }
			};
    \node[right=.5 of 4] (5) {$\bullet$}
      child {node[my-red] {\tiny$\blacksquare$}
      	child[my-red] {node[my-red] {\tiny$\blacksquare$}
      	  child[my-red]
      	      }
			};
    \node[right=.5 of 5] (6) {$\bullet$}
      child {node[my-red] {\tiny$\blacksquare$}
      	child[my-red] {node[black] {$\diamond$}}
			};
    \node[right=.5 of 6] (7) {$\bullet$}
      child {node[black] {$\diamond$}};
    \node[right=.5 of 7, my-red] (8) {\tiny$\blacksquare$}
      child[my-red] {node[my-blue] {$\bullet$}
      	child[my-blue] {node[my-blue] {$\bullet$}
      	  child[my-blue]
      	      }
			};
    \node[right=.5 of 8, my-red] (9) {\tiny$\blacksquare$}
      child[my-red] {node[my-blue] {$\bullet$}
      	child[my-blue] {node[my-red] {\tiny$\blacksquare$}
      	  child[my-red]
      	      }
			};
    \node[right=.5 of 9, my-red] (10) {\tiny$\blacksquare$}
      child[my-red] {node[my-blue] {$\bullet$}
      	child[my-blue] {node[black] {$\diamond$}}
			};
    \node[right=.5 of 10, my-red] (11) {\tiny$\blacksquare$}
      child[my-red] {node[my-red] {\tiny$\blacksquare$}
      	child[my-red] {node[my-blue] {$\bullet$}
      	  child[my-blue]
      	      }
			};
    \node[right=.5 of 11, my-red] (12) {\tiny$\blacksquare$}
      child[my-red] {node[my-red] {\tiny$\blacksquare$}
      	child[my-red] {node[my-red] {\tiny$\blacksquare$}
      	  child[my-red]
      	      }
			};
    \node[right=.5 of 12, my-red] (13) {\tiny$\blacksquare$}
      child[my-red] {node[my-red] {\tiny$\blacksquare$}
      	child[my-red] {node[black] {$\diamond$}}
			};
    \node[right=.5 of 13, my-red] (14) {\tiny$\blacksquare$}
      child[my-red] {node[black] {$\diamond$}};
    \node[right=.5 of 14, black] (15) {$\diamond$};
  \end{tikzpicture}
  };
\end{tikzpicture}
\]
\end{enumerate}
\end{solution}
\end{exercise}

\index{composition product!involving constant polynomials}
\index{polynomial functor!application as composition}

\begin{example}[Composing polynomials with constants] \label{ex.apply_2}
For any set $X$ and polynomial $p$, we can take $p(X)\in\smset$; indeed $p\colon\smset\to\smset$ is a functor! In particular, by this point you've seen us write $p(\1)$ hundreds of times. But we've also seen that $X$ is itself a polynomial, namely a constant one.

It's not hard to see that $p(X)\iso p\tri X$. Here's a picture, where $p\coloneqq\yon^\3+\yon+\1$ and $X\coloneqq\2$.
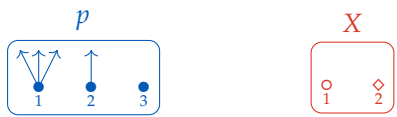
\[
\begin{tikzpicture}[rounded corners]
	\node (p1) [draw, my-blue, "\color{my-blue} $p$" above] {
	\begin{tikzpicture}[trees, sibling distance=2.5mm]
    \node["\tiny 1" below] (1) {$\bullet$}
      child {}
      child {}
      child {};
    \node[right=.5 of 1,"\tiny 2" below] (2) {$\bullet$}
      child {};
      ;
    \node[right=.5 of 2,"\tiny 3" below] (3) {$\bullet$}
      ;
  \end{tikzpicture}
  };
	\node (p2) [draw, my-red, right=2 of p1, "\color{my-red}$X$" above] {
	\begin{tikzpicture}[trees, sibling distance=2.5mm]
    \node["\tiny 1" below] (1) {$\circ$};
    \node[right=.5 of 1,"\tiny 2" below] (4) {$\diamond$};
    \node[above=10pt of 4] {};
    ;
  \end{tikzpicture}
  };
\end{tikzpicture}
\]
Let's see how $(\yon^\3+\yon+\1)\tri\2$ looks.
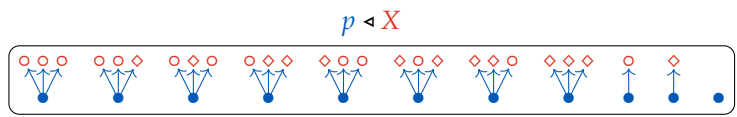
\[
\begin{tikzpicture}[rounded corners]
	\node (p1) [draw, "{\color{my-blue} $p$}$\:\tri\:${\color{my-red}$X$}" above] {
	\begin{tikzpicture}[trees, sibling distance=2.5mm, my-blue]
    \node[my-blue] (1) {$\bullet$}
      child {node[my-red] {$\circ$}}
      child {node[my-red] {$\circ$}}
      child {node[my-red] {$\circ$}};
    \node[my-blue, right=.8 of 1] (2) {$\bullet$}
      child {node[my-red] {$\circ$}}
      child {node[my-red] {$\circ$}}
      child {node[my-red] {$\diamond$}};
    \node[my-blue, right=.8 of 2] (3) {$\bullet$}
      child {node[my-red] {$\circ$}}
      child {node[my-red] {$\diamond$}}
      child {node[my-red] {$\circ$}};
    \node[my-blue, right=.8 of 3] (4) {$\bullet$}
      child {node[my-red] {$\circ$}}
      child {node[my-red] {$\diamond$}}
      child {node[my-red] {$\diamond$}};
    \node[my-blue, right=.8 of 4] (5) {$\bullet$}
      child {node[my-red] {$\diamond$}}
      child {node[my-red] {$\circ$}}
      child {node[my-red] {$\circ$}};
    \node[my-blue, right=.8 of 5] (6) {$\bullet$}
      child {node[my-red] {$\diamond$}}
      child {node[my-red] {$\circ$}}
      child {node[my-red] {$\diamond$}};
    \node[my-blue, right=.8 of 6] (7) {$\bullet$}
      child {node[my-red] {$\diamond$}}
      child {node[my-red] {$\diamond$}}
      child {node[my-red] {$\circ$}};
    \node[my-blue, right=.8 of 7] (8) {$\bullet$}
      child {node[my-red] {$\diamond$}}
      child {node[my-red] {$\diamond$}}
      child {node[my-red] {$\diamond$}};
    \node[my-blue, right=.6 of 8] (9) {$\bullet$}
      child {node[my-red] {$\circ$}};
    \node[my-blue, right=.4 of 9] (10) {$\bullet$}
      child {node[my-red] {$\diamond$}};
    \node[my-blue, right=.4 of 10] (11) {$\bullet$};
	\end{tikzpicture}
	};
\end{tikzpicture}
\]
It has $11$ positions and no height-$2$ leaves, which means it's a set (constant polynomial, with no directions), namely $p\tri X\iso \1\1$.

We could also draw $X\tri p$, since both are perfectly valid polynomials. Here it is:
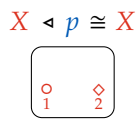
\[
\begin{tikzpicture}[rounded corners]
	\node (p2) [draw, "${\color{my-red}X}\:\tri\:{\color{my-blue}p}\iso{\color{my-red}X}$" above] {
	\begin{tikzpicture}[trees, sibling distance=2.5mm, my-red]
    \node["\tiny 1" below] (1) {$\circ$};
    \node[right=.5 of 1,"\tiny 2" below] (4) {$\diamond$};
    \node[above=10pt of 4] {};
    ;
  \end{tikzpicture}
  };
\end{tikzpicture}
\]
Each of the leaves in $X$, of which there are none, is given a $p$-corolla.
\end{example}

\begin{exercise}\label{exc.composing_with_constants}
\begin{enumerate}
	\item Choose a polynomial $p$ and draw $p\tri\1$ in the style of \cref{ex.apply_2}.
	\item Show that if $X$ is a set (considered as a constant polynomial) and $p$ is any polynomial, then $X\tri p\iso X$.
	\item \label{exc.composing_with_constants.appl} Show that if $X$ is a set and $p$ is a polynomial, then $p\tri X\iso p(X)$, where $p(X)$ is the set given by applying $p$ as a functor to $X$.
\qedhere
\end{enumerate}
\begin{solution}
\begin{enumerate}
    \item We pick the list polynomial, $p\coloneqq\1+\yon+\yon^\2+\yon^\3+\cdots$, drawn as follows:
\[
\begin{tikzpicture}[rounded corners]
\node (p1) [draw] {
  \begin{tikzpicture}[trees, sibling distance=3mm]
    \node (1) {$\bullet$};
    \node[right=.3 of 1] (2) {$\bullet$}
      child {};
    \node[right=.4 of 2] (3) {$\bullet$}
      child {}
      child {};
    \node[right=.6 of 3] (4) {$\bullet$}
      child {}
      child {}
      child {};
    \node[right=.6 of 4] {$\cdots$};
  \end{tikzpicture}
};
\end{tikzpicture}
\]
Then here is a picture of $p\tri\1$:
\[
\begin{tikzpicture}[rounded corners]
\node (p1) [draw] {
  \begin{tikzpicture}[trees, sibling distance=3mm]
    \node (1) {$\bullet$};
    \node[right=.3 of 1] (2) {$\bullet$}
      child {node {$\bullet$}};
    \node[right=.4 of 2] (3) {$\bullet$}
      child {node {$\bullet$}}
      child {node {$\bullet$}};
    \node[right=.6 of 3] (4) {$\bullet$}
      child {node {$\bullet$}}
      child {node {$\bullet$}}
      child {node {$\bullet$}};
    \node[right=.6 of 4] {$\cdots$};
  \end{tikzpicture}
};
\end{tikzpicture}
\]
\end{enumerate}
Below, $X$ is a set and $p$ is a polynomial.
\begin{enumerate}[resume]
    \item A constant functor composed with any functor is still the same constant functor, so $X \tri p \iso X$.
    We can also verify this using \eqref{eqn.composite_formula}:
    \[
        X \tri p \iso \sum_{i \in X} \prod_{a \in \varnothing} \sum_{j \in p(\1)} \prod_{b \in p[j]} \yon \iso \sum_{i \in X} \1 \iso X.
    \]
    \item When viewed as functors, it is easy to see that $p \tri X \iso p(X)$.
    We can also verify this using \eqref{eqn.composite_formula}:
    \[
        p \tri X \iso \sum_{i \in p(\1)} \prod_{a \in p[i]} \sum_{j \in X} \prod_{b \in \varnothing} \yon \iso \sum_{i \in p(\1)} \prod_{a \in p[i]} \sum_{j \in X} \1 \iso \sum_{i \in p(\1)} \prod_{a \in p[i]} X \iso \sum_{i \in p(\1)} X^{p[i]} \iso p(X).
    \]
\end{enumerate}
\end{solution}
\end{exercise}

In particular, this means we could write the position-set of a polynomial $p$ interchangeably as $p(\1)$ or as $p\tri\1$.
We'll generally write $p(\1)$ when we want to emphasize the position-set as a set, and $p\tri\1$ when we want to emphasize the position-set as a polynomial (albeit a constant one, with no directions).

\index{polynomial functor!positions as $p(\1)$}\index{positions|see{polynomial functor, positions and directions}}\index{directions|see{polynomial functor, positions and directions}}

\begin{exercise}\index{natural transformation!positions map as $\1$-component}
Let $\varphi\colon p\to q$ be a lens and $X$ be a set viewed as a constant polynomial.
Consider the lens $\varphi\tri X\colon p\tri X\to q\tri X$, given by taking the composition product of $\varphi$ with the identity lens on $X$.
Show that $\varphi\tri X$, when viewed as a function $p(X)\to q(X)$ between sets, is exactly the $X$-component of $\varphi$ viewed as a natural transformation.
\begin{solution}
Let $\varphi\colon p\to q$ be a lens and $X$ be a set viewed as a constant polynomial.
Note that every component of the identity natural transformation on $X$ as a constant functor is just the identity function $\id_X\colon X\to X$ on $X$ as a set.
Then by \cref{def.horiz_comp_nat_trans}, any component of the composition product $\varphi\tri X$ viewed as a natural transformation is given by the composite function
\[
    p(X) \To{\varphi_X} q(X) \To{q(\id_X)} q(X).
\]
By functoriality, $q(\id_X)$ is itself an identity function, so every component of $\varphi\tri X$ is the $X$-component of $\varphi$.
Therefore $\varphi\tri X$ as a function can be identified with the $X$-component of $\varphi$, as desired.
\end{solution}
\end{exercise}

\begin{exercise}\label{ex.compose_yon}
For any $p\in\poly$ there are natural isomorphisms $p\iso p\tri \yon$ and $p\iso\yon\tri p$.\index{isomorphism!natural}
\begin{enumerate}
	\item Thinking of polynomials as functors $\smset\to\smset$, what functor does $\yon$ represent?
	\item Why are $p\tri\yon$ and $\yon\tri p$ isomorphic to $p$?
	\item Let $p\coloneqq\yon^\3+\yon+\1$.
	In terms of tree pictures, draw $p\tri\yon$ and $\yon\tri p$, and explain pictorially how to see the isomorphisms $p\tri\yon \iso p \iso \yon\tri p$.
\qedhere
\end{enumerate}
\begin{solution}
\begin{enumerate}
    \item The polynomial $\yon$ is the identity functor on $\smset$.
    \item Composing any functor with the identity functor yields the original functor, so $p\tri\yon \iso p \iso \yon\tri p$.
    \item Before we draw $\yon\tri p$ and $p\tri\yon$, here are pictures of $p$ and $\yon$ individually as corolla forests:
\[
\begin{tikzpicture}[rounded corners]
	\node (p1) [draw, my-blue, "\color{my-blue} $p$" above] {
	\begin{tikzpicture}[trees, sibling distance=2.5mm]
    \node (1) {$\bullet$}
      child {}
      child {}
      child {};
    \node[right=.5 of 1] (2) {$\bullet$}
      child {};
    \node[right=.5 of 2] (3) {$\bullet$};
  \end{tikzpicture}
  };
	\node (p2) [draw, right=2 of p1, "$\yon$" above] {
	\begin{tikzpicture}[trees, sibling distance=2.5mm]
    \node (1) {\tiny$\blacksquare$}
      child {};
  \end{tikzpicture}
  };
\end{tikzpicture}
\]
Now here is a picture of $p\tri\yon$, obtained by grafting the one-leaf $\yon$-corolla to all the leaves of each $p$-corolla in turn:
\[
\begin{tikzpicture}[rounded corners]
	\node (p1) [draw, my-blue, "{\color{my-blue}$p$}$\:\tri\:\yon$" above] {
	\begin{tikzpicture}[trees, sibling distance=5mm]
    \node (1) {$\bullet$}
      child {node[black] {\tiny$\blacksquare$} child[black]}
      child {node[black] {\tiny$\blacksquare$} child[black]}
      child {node[black] {\tiny$\blacksquare$} child[black]};
    \node[right=1 of 1] (2) {$\bullet$}
      child {node[black] {\tiny$\blacksquare$} child[black]};
    \node[right=.5 of 2] (3) {$\bullet$};
  \end{tikzpicture}
  };
\end{tikzpicture}
\]
This is just $p$ with every direction extended up one level, so it is still a picture of $p$.
And here is a picture of $\yon\tri p$, obtained by grafting each $p$-corolla to the single leaf of $\yon$:
\[
\begin{tikzpicture}[rounded corners]
	\node (p2) [draw, "$\yon\:\tri\:${\color{my-blue}$p$}" above] {
	\begin{tikzpicture}[trees, sibling distance=2.5mm]
    \node (1) {\tiny$\blacksquare$}
      child {
        node[my-blue] {$\bullet$}
            child[my-blue]
            child[my-blue]
            child[my-blue]
      };
    \node[right=.5 of 1] (2) {\tiny$\blacksquare$}
      child {
        node[my-blue] {$\bullet$}
            child[my-blue]
      };
    \node[right=.5 of 2] (3) {\tiny$\blacksquare$}
      child {
        node[my-blue] {$\bullet$}
      };
  \end{tikzpicture}
  };
\end{tikzpicture}
\]
\end{enumerate}
\end{solution}
This is just $p$ with every position propped up one level, so it is also still a picture of $p$.
\end{exercise}

How shall we think about taking the composition product of lenses in terms of our tree pictures?
We can interpret the results of \cref{exc.comp_prod_lens} as follows.

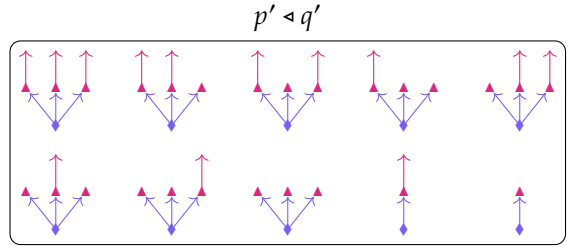
\index{composition product!tree picture of}

\begin{example}\label{ex.comp_prod_trees}
Take $p\coloneqq \yon^\2+\yon$, $q\coloneqq\yon^\2+\yon$, $p'\coloneqq\yon^\3+\yon$, and $q'\coloneqq\yon+\1$.
\[
\begin{tikzpicture}[rounded corners]
	\node (p) [draw, my-blue, "$p=$" left] {
	\begin{tikzpicture}[trees, sibling distance=5mm]
    \node["\tiny 1" below] (1) {\tiny$\blacksquare$}
      child {}
      child {};
    \node[right=.5 of 1,"\tiny 2" below] (2) {\tiny$\blacksquare$}
      child {};
  \end{tikzpicture}
  };
	\node (q) [draw, my-red, above=1 of p, "$q=$" left] {
	\begin{tikzpicture}[trees, sibling distance=5mm]
    \node["\tiny 1" below] (1) {$\bullet$}
      child {}
      child {};
    \node[right=.5 of 1,"\tiny 2" below] (2) {$\bullet$}
      child {};
  \end{tikzpicture}
  };
	\node (p') [draw, my-lavender, right=3 of p, "$p'=$" left] {
	\begin{tikzpicture}[trees, sibling distance=2.5mm]
    \node["\tiny 1" below] (1) {\tiny$\blacklozenge$}
      child {}
      child {}
      child {};
    \node[right=.5 of 1,"\tiny 2" below] (2) {\tiny$\blacklozenge$}
      child {};
  \end{tikzpicture}
  };
	\node (q') [draw, my-magenta, above=1 of p', "$q'=$" left] {
	\begin{tikzpicture}[trees, sibling distance=2.5mm]
    \node["\tiny 1" below] (1) {\tiny$\blacktriangle$}
      child {};
    \node[right=.5 of 1,"\tiny 2" below] (2) {\tiny$\blacktriangle$}
    ;
  \end{tikzpicture}
  };
\end{tikzpicture}
\]
For any pair of lenses $p\to p'$ and $q\to q'$, we have a lens $p\tri q\to p'\tri q'$. We draw $p\tri q$ and $p'\tri q'$ below.
\[
\begin{tikzpicture}[rounded corners]
	\node (p1) [draw, "$p\tri q$" above] {
	\begin{tikzpicture}[trees,
		level 1/.style={sibling distance=8mm},
	  level 2/.style={sibling distance=2.5mm}]
    \node[my-blue] (1) {\tiny$\blacksquare$}
      child[my-blue] {node[my-red] {$\bullet$}
      	child[my-red]
				child[my-red]
			}
      child[my-blue] {node[my-red] {$\bullet$}
      	child[my-red]
				child[my-red]
			};
    \node[my-blue, right=1.7 of 1] (2) {\tiny$\blacksquare$}
      child[my-blue] {node[my-red] {$\bullet$}
				child[my-red]
				child[my-red]
			}
      child[my-blue] {node[my-red] {$\bullet$}
				child[my-red]
			};
    \node[my-blue, right=1.5 of 2] (3) {\tiny$\blacksquare$}
      child[my-blue] {node[my-red] {$\bullet$}
      	child[my-red]
			}
      child[my-blue] {node[my-red] {$\bullet$}
				child[my-red]
				child[my-red]
			};
    \node[my-blue, right=1.5 of 3] (4) {\tiny$\blacksquare$}
      child[my-blue] {node[my-red] {$\bullet$}
      	child[my-red]
			}
      child[my-blue] {node[my-red] {$\bullet$}
				child[my-red]
			};
    \node[my-blue, right=1.2 of 4] (5) {\tiny$\blacksquare$}
      child[my-blue] {node[my-red] {$\bullet$}
      	child[my-red]
      	child[my-red]
			};
    \node[my-blue, right=.8 of 5] (6) {\tiny$\blacksquare$}
      child[my-blue] {node[my-red] {$\bullet$}
      	child[my-red]
			};
  \end{tikzpicture}
  };
\end{tikzpicture}
\]
\[
\begin{tikzpicture}[rounded corners]
	\node (p1) [draw, "$p'\tri q'$" above] {
	\begin{tikzpicture}[trees,
		level 1/.style={sibling distance=4mm},
	  level 2/.style={sibling distance=2.5mm}]
    \node[my-lavender] (1) {\tiny$\blacklozenge$}
      child[my-lavender] {node[my-magenta] {\tiny$\blacktriangle$}
      	child[my-magenta]
			}
      child[my-lavender] {node[my-magenta] {\tiny$\blacktriangle$}
      	child[my-magenta]
			}
      child[my-lavender] {node[my-magenta] {\tiny$\blacktriangle$}
				child[my-magenta]
			};
    \node[my-lavender, right=1.4 of 1] (2) {\tiny$\blacklozenge$}
      child[my-lavender] {node[my-magenta] {\tiny$\blacktriangle$}
      	child[my-magenta]
			}
      child[my-lavender] {node[my-magenta] {\tiny$\blacktriangle$}
      	child[my-magenta]
			}
      child[my-lavender] {node[my-magenta] {\tiny$\blacktriangle$}
			};
    \node[my-lavender, right=1.4 of 2] (3) {\tiny$\blacklozenge$}
      child[my-lavender] {node[my-magenta] {\tiny$\blacktriangle$}
      	child[my-magenta]
			}
      child[my-lavender] {node[my-magenta] {\tiny$\blacktriangle$}
			}
      child[my-lavender] {node[my-magenta] {\tiny$\blacktriangle$}
				child[my-magenta]
			};
    \node[my-lavender, right=1.4 of 3] (4) {\tiny$\blacklozenge$}
      child[my-lavender] {node[my-magenta] {\tiny$\blacktriangle$}
      	child[my-magenta]
			}
      child[my-lavender] {node[my-magenta] {\tiny$\blacktriangle$}
			}
      child[my-lavender] {node[my-magenta] {\tiny$\blacktriangle$}
			};
    \node[my-lavender, right=1.4 of 4] (5) {\tiny$\blacklozenge$}
      child[my-lavender] {node[my-magenta] {\tiny$\blacktriangle$}
			}
      child[my-lavender] {node[my-magenta] {\tiny$\blacktriangle$}
      	child[my-magenta]
			}
      child[my-lavender] {node[my-magenta] {\tiny$\blacktriangle$}
				child[my-magenta]
			};
    \node[my-lavender, below=1.2 of 1] (6) {\tiny$\blacklozenge$}
      child[my-lavender] {node[my-magenta] {\tiny$\blacktriangle$}
			}
      child[my-lavender] {node[my-magenta] {\tiny$\blacktriangle$}
      	child[my-magenta]
			}
      child[my-lavender] {node[my-magenta] {\tiny$\blacktriangle$}
			};
    \node[my-lavender, right=1.4 of 6] (7) {\tiny$\blacklozenge$}
      child[my-lavender] {node[my-magenta] {\tiny$\blacktriangle$}
			}
      child[my-lavender] {node[my-magenta] {\tiny$\blacktriangle$}
			}
      child[my-lavender] {node[my-magenta] {\tiny$\blacktriangle$}
      	child[my-magenta]
			};
    \node[my-lavender, right=1.4 of 7] (8) {\tiny$\blacklozenge$}
      child[my-lavender] {node[my-magenta] {\tiny$\blacktriangle$}
			}
      child[my-lavender] {node[my-magenta] {\tiny$\blacktriangle$}
			}
      child[my-lavender] {node[my-magenta] {\tiny$\blacktriangle$}
			};
    \node[my-lavender, below=1.2 of 4] (9) {\tiny$\blacklozenge$}
      child[my-lavender] {node[my-magenta] {\tiny$\blacktriangle$}
      	child[my-magenta]
			};
    \node[my-lavender, below=1.2 of 5] (10) {\tiny$\blacklozenge$}
      child[my-lavender] {node[my-magenta] {\tiny$\blacktriangle$}
			};
  \end{tikzpicture}
  };
\end{tikzpicture}
\]

We also pick a pair of lenses, $\varphi\colon p\to p'$ and $\psi\colon q\to q'$.
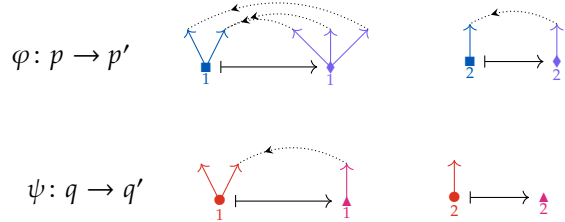
\[
\begin{tikzpicture}
	\node (p1) {\raisebox{.3cm}{$\varphi\colon p\to p'$}\qquad
	\begin{tikzpicture}[trees, sibling distance=5mm]
    \node[my-blue, "\color{my-blue} \tiny 1" below] (1) {\tiny$\blacksquare$}
      child[my-blue] {coordinate (11)}
      child[my-blue] {coordinate (12)};
    \node[right=1.5 of 1, my-lavender, "\color{my-lavender} \tiny 1" below] (2) {\tiny$\blacklozenge$}
      child[my-lavender] {coordinate (21)}
      child[my-lavender] {coordinate (22)}
      child[my-lavender] {coordinate (23)};
    \draw[|->, shorten <= 3pt, shorten >= 3pt] (1) -- (2);
    \begin{scope}[densely dotted, bend right]
      \draw[postaction={decorate}] (21) to (12);
      \draw[postaction={decorate}] (22) to (12);
      \draw[postaction={decorate}] (23) to (11);
    \end{scope}
  \end{tikzpicture}
	};
	\node (p2) [below right=-1.3 and 1 of p1] {
	\begin{tikzpicture}[trees, sibling distance=5mm]
    \node[my-blue, "\color{my-blue} \tiny 2" below] (1) {\tiny$\blacksquare$}
      child[my-blue] {coordinate (11)};
    \node[right=of 1, my-lavender, "\color{my-lavender} \tiny 2" below] (2) {\tiny$\blacklozenge$}
      child[my-lavender] {coordinate (21)};
    \draw[|->, shorten <= 3pt, shorten >= 3pt] (1) -- (2);
    \begin{scope}[densely dotted, bend right]
      \draw[postaction={decorate}] (21) to (11);
		\end{scope}
  \end{tikzpicture}
	};
	\node [below=.5 of p1] (p3) {\raisebox{.3cm}{$\psi\colon q\to q'$}\qquad
	\begin{tikzpicture}[trees, sibling distance=5mm]
    \node[my-red, "\color{my-red} \tiny 1" below] (1) {$\bullet$}
      child[my-red] {coordinate (11)}
      child[my-red] {coordinate (12)};
    \node[right=1.5 of 1, my-magenta, "\color{my-magenta} \tiny 1" below] (2) {\tiny$\blacktriangle$}
      child[my-magenta] {coordinate (21)};
    \draw[|->, shorten <= 3pt, shorten >= 3pt] (1) -- (2);
    \begin{scope}[densely dotted, bend right]
      \draw[postaction={decorate}] (21) to (12);
    \end{scope}
  \end{tikzpicture}
	};
	\node (p4) [below right=-1.05 and 1 of p3] {
	\begin{tikzpicture}[trees, sibling distance=5mm]
    \node[my-red, "\color{my-red} \tiny 2" below] (1) {$\bullet$}
      child[my-red] {coordinate (11)};
    \node[right=of 1, my-magenta, "\color{my-magenta} \tiny 2" below] (2) {\tiny$\blacktriangle$};
    \draw[|->, shorten <= 3pt, shorten >= 3pt] (1) -- (2);
  \end{tikzpicture}
	};
\end{tikzpicture}
\]
Then by \cref{exc.comp_prod_lens}, we can form the lens $\varphi\tri \psi\colon p\tri q\to p'\tri q'$ as follows.
On positions, we follow \eqref{eqn.comp_lens_pos}: for each tree $t$ in the picture of $p\tri q$, we begin by using $\varphi_\1$ to send the $p$-corolla $i$ that forms the bottom level of $t$ to a $p'$-corolla $i'$.
Then for each $p'[i']$-leaf $a'$ of $i'$, to choose which $q'$-corolla gets grafted onto $a'$, we use $\varphi^\sharp_i$ to send $a'$ back to a $p[i]$-leaf $a$.
Since $t$ has the corolla $i$ as its bottom level, $a$ is just a height-$1$ vertex of the tree $t$.
So we can take the $q$-corolla $j$ that is grafted onto $a$ in $t$, then use $\psi_\1$ to send $j$ forward to a $q'$-corolla $j'$.
This is the corolla we graft onto the $p'[i']$-leaf $a'$.
All this specifies a tree $t'$ in $p'\tri q'$ that $t$ gets sent to via $(\varphi\tri \psi)_\1$.

On directions, we follow \eqref{eqn.comp_lens_dir}: picking a direction of $t'$ consists of picking a height-$1$ vertex $a'$ and a height-$2$ leaf $b'$ emanating from $a'$.
The on-directions function $\varphi^\sharp_i$ sends $a'$ back to a height-$1$ vertex $a$ of $t$, and as we saw, the on-positions function $\psi_\1$ sends the $q$-corolla $j$ grafted onto $a$ in $t$ forward to the $q'$-corolla grafted onto $a'$.
Then $b'$ is a leaf of that $q'$-corolla, and $\psi^\sharp_j$ sends $b'$ back to a leaf $b$ emanating from $a$.
So the on-directions function $(\varphi\tri \psi)^\sharp_t$ sends the height-$2$ leaf $b'$ to the height-$2$ leaf $b$.

We draw the lens $\varphi\tri \psi\to p\tri q\to p'\tri q'$ below.
To avoid clutter, we leave out the arrows for $\psi_\1$ that show how the red corollas on the right are selected; we hope the reader can put it together for themselves.
\[
	\begin{tikzpicture}[trees]
	\begin{scope}[
		  level 1/.style={sibling distance=8mm},
      level 2/.style={sibling distance=5mm}]
    \node[my-blue] (1) {\tiny$\blacksquare$}
      child[my-blue] {node[my-red] (11') {$\bullet$}
      	child[my-red] {coordinate (11)}
				child[my-red] {coordinate (12)}
			}
      child[my-blue] {node[my-red] (12') {$\bullet$}
      	child[my-red] {coordinate (13)}
				child[my-red] {coordinate (14)}
			};
    \node[my-blue, right=5 of 1] (2) {\tiny$\blacksquare$}
      child[my-blue] {node[my-red] (21') {$\bullet$}
      	child[my-red] {coordinate (21)}
				child[my-red] {coordinate (22)}
			}
      child[my-blue] {node[my-red] (22') {$\bullet$}
      	child[my-red] {coordinate (23)}
			};
    \node[my-blue, below=1.3 of 1] (3) {\tiny$\blacksquare$}
      child[my-blue] {node[my-red] (31') {$\bullet$}
      	child[my-red] {coordinate (31)}
			}
      child[my-blue] {node[my-red] (32') {$\bullet$}
      	child[my-red] {coordinate (32)}
				child[my-red] {coordinate (33)}
			};
    \node[my-blue] at (2|-3) (4) {\tiny$\blacksquare$}
      child[my-blue] {node[my-red] (41') {$\bullet$}
      	child[my-red] {coordinate (41)}
			}
      child[my-blue] {node[my-red] (42') {$\bullet$}
      	child[my-red] {coordinate (42)}
			};
    \node[my-blue, below=1.3 of 3] (5) {\tiny$\blacksquare$}
      child[my-blue] {node[my-red] (51') {$\bullet$}
      	child[my-red] {coordinate (51)}
				child[my-red] {coordinate (52)}
			};
    \node[my-blue] at (4|-5) (6) {\tiny$\blacksquare$}
      child[my-blue] {node[my-red] (61') {$\bullet$}
      	child[my-red] {coordinate (61)}
			};
	\end{scope}
  \begin{scope}[
      level 1/.style={sibling distance=4mm},
      level 2/.style={sibling distance=5mm}]
    \node[my-lavender, right=2 of 1] (1') {\tiny$\blacklozenge$}
      child[my-lavender] {node[my-magenta] (1'1') {\tiny$\blacktriangle$}
      	child[my-magenta] {coordinate (1'1)}
			}
      child[my-lavender] {node[my-magenta] (1'2') {\tiny$\blacktriangle$}
      	child[my-magenta] {coordinate (1'2)}
			}
      child[my-lavender] {node[my-magenta] (1'3') {\tiny$\blacktriangle$}
      	child[my-magenta] {coordinate (1'3)}
			};
    \node[my-lavender, right=2 of 2] (2') {\tiny$\blacklozenge$}
      child[my-lavender] {node[my-magenta] (2'1') {\tiny$\blacktriangle$}
			}
      child[my-lavender] {node[my-magenta] (2'2') {\tiny$\blacktriangle$}
			}
      child[my-lavender] {node[my-magenta] (2'3') {\tiny$\blacktriangle$}
      	child[my-magenta] {coordinate (2'1)}
			};
    \node[my-lavender, right=2 of 3] (3') {\tiny$\blacklozenge$}
      child[my-lavender] {node[my-magenta] (3'1') {\tiny$\blacktriangle$}
      	child[my-magenta] {coordinate (3'1)}
			}
      child[my-lavender] {node[my-magenta] (3'2') {\tiny$\blacktriangle$}
      	child[my-magenta] {coordinate (3'2)}
			}
      child[my-lavender] {node[my-magenta] (3'3') {\tiny$\blacktriangle$}
			};
    \node[my-lavender, right=2 of 4] (4') {\tiny$\blacklozenge$}
      child[my-lavender] {node[my-magenta] (4'1') {\tiny$\blacktriangle$}
			}
      child[my-lavender] {node[my-magenta] (4'2') {\tiny$\blacktriangle$}
			}
      child[my-lavender] {node[my-magenta] (4'3') {\tiny$\blacktriangle$}
			};
    \node[my-lavender, right=2 of 5] (5') {\tiny$\blacklozenge$}
      child[my-lavender] {node[my-magenta] (5'1') {\tiny$\blacktriangle$}
      	child[my-magenta] {coordinate (5'1)}
			};
    \node[my-lavender, right=2 of 6] (6') {\tiny$\blacklozenge$}
      child[my-lavender] {node[my-magenta] (6'1') {\tiny$\blacktriangle$}
			};
\draw[|->, shorten <= 3pt, shorten >= 3pt] (1) -- (1');
\draw[|->, shorten <= 3pt, shorten >= 3pt] (2) -- (2');
\draw[|->, shorten <= 3pt, shorten >= 3pt] (3) -- (3');
\draw[|->, shorten <= 3pt, shorten >= 3pt] (4) -- (4');
\draw[|->, shorten <= 3pt, shorten >= 3pt] (5) -- (5');
\draw[|->, shorten <= 3pt, shorten >= 3pt] (6) -- (6');
    \begin{scope}[densely dotted, bend right=15pt]
      \draw[postaction={decorate}] (1'1') to (12');
      \draw[postaction={decorate}] (1'2') to (12');
      \draw[postaction={decorate}] (1'3') to (11');
      \draw[postaction={decorate}] (1'1) to (14);
      \draw[postaction={decorate}] (1'2) to (14);
      \draw[postaction={decorate}] (1'3) to (12);
      \draw[postaction={decorate}] (2'1') to (22');
      \draw[postaction={decorate}] (2'2') to (22');
      \draw[postaction={decorate}] (2'3') to (21');
      \draw[postaction={decorate}] (2'1) to (23);
      \draw[postaction={decorate}] (3'1') to (32');
      \draw[postaction={decorate}] (3'2') to (32');
      \draw[postaction={decorate}] (3'3') to (31');
      \draw[postaction={decorate}] (3'1) to (33);
      \draw[postaction={decorate}] (3'2) to (33);
      \draw[postaction={decorate}] (4'1') to (42');
      \draw[postaction={decorate}] (4'2') to (42');
      \draw[postaction={decorate}] (4'3') to (41');
      \draw[postaction={decorate}] (5'1') to (51');
      \draw[postaction={decorate}] (5'1) to (52);
      \draw[postaction={decorate}] (6'1') to (61');
    \end{scope}

	\end{scope}
  \end{tikzpicture}
\]
\end{example}

\begin{exercise}
With $p,q,p',q'$ and $\varphi,\psi$ as in \cref{ex.comp_prod_trees}, draw the lens $\psi\tri\varphi\colon q\tri p\to q'\tri p'$ in terms of trees as in the example.
\begin{solution}
Using the definitions, instructions, and style from \cref{ex.comp_prod_trees}, we draw $\psi\tri\varphi\colon q\tri p\to q'\tri p'$:
\[
	\begin{tikzpicture}[trees]
	\begin{scope}[
		level 1/.style={sibling distance=8mm},
	  level 2/.style={sibling distance=5mm}]
    \node[my-red] (1) {$\bullet$}
      child[my-red] {node[my-blue] (11') {\tiny$\blacksquare$}
      	child[my-blue] {coordinate (11)}
				child[my-blue] {coordinate (12)}
			}
      child[my-red] {node[my-blue] (12') {\tiny$\blacksquare$}
      	child[my-blue] {coordinate (13)}
				child[my-blue] {coordinate (14)}
			};
    \node[my-red, right=5 of 1] (2) {$\bullet$}
      child[my-red] {node[my-blue] (21') {\tiny$\blacksquare$}
      	child[my-blue] {coordinate (21)}
				child[my-blue] {coordinate (22)}
			}
      child[my-red] {node[my-blue] (22') {\tiny$\blacksquare$}
      	child[my-blue] {coordinate (23)}
			};
    \node[my-red, below=1.3 of 1] (3) {$\bullet$}
      child[my-red] {node[my-blue] (31') {\tiny$\blacksquare$}
      	child[my-blue] {coordinate (31)}
			}
      child[my-red] {node[my-blue] (32') {\tiny$\blacksquare$}
      	child[my-blue] {coordinate (32)}
				child[my-blue] {coordinate (33)}
			};
    \node[my-red] at (2|-3) (4) {$\bullet$}
      child[my-red] {node[my-blue] (41') {\tiny$\blacksquare$}
      	child[my-blue] {coordinate (41)}
			}
      child[my-red] {node[my-blue] (42') {\tiny$\blacksquare$}
      	child[my-blue] {coordinate (42)}
			};
    \node[my-red, below=1.3 of 3] (5) {$\bullet$}
      child[my-red] {node[my-blue] (51') {\tiny$\blacksquare$}
      	child[my-blue] {coordinate (51)}
				child[my-blue] {coordinate (52)}
			};
    \node[my-red] at (4|-5) (6) {$\bullet$}
      child[my-red] {node[my-blue] (61') {\tiny$\blacksquare$}
      	child[my-blue] {coordinate (61)}
			};
		\end{scope}
	\begin{scope}[
		level 1/.style={sibling distance=4mm},
	  level 2/.style={sibling distance=2.5mm}]
	    \node[my-magenta, right=2 of 1] (1') {\tiny$\blacktriangle$}
      child[my-magenta] {node[my-lavender] (1'1') {\tiny$\blacklozenge$}
      	child[my-lavender] {coordinate (1'1)}
      	child[my-lavender] {coordinate (1'2)}
      	child[my-lavender] {coordinate (1'3)}
			};
    \node[my-magenta, right=2 of 2] (2') {\tiny$\blacktriangle$}
      child[my-magenta] {node[my-lavender] (2'1') {\tiny$\blacklozenge$}
        child[my-lavender] {coordinate (2'1)}
			};
    \node[my-magenta, right=2 of 3] (3') {\tiny$\blacktriangle$}
      child[my-magenta] {node[my-lavender] (3'1') {\tiny$\blacklozenge$}
      	child[my-lavender] {coordinate (3'1)}
      	child[my-lavender] {coordinate (3'2)}
      	child[my-lavender] {coordinate (3'3)}
			};
    \node[my-magenta, right=2 of 4] (4') {\tiny$\blacktriangle$}
      child[my-magenta] {node[my-lavender] (4'1') {\tiny$\blacklozenge$}
      	child[my-lavender] {coordinate (4'1)}
			};
    \node[my-magenta, right=2 of 5] (5') {\tiny$\blacktriangle$};
    \node[my-magenta, right=2 of 6] (6') {\tiny$\blacktriangle$};
\draw[|->, shorten <= 3pt, shorten >= 3pt] (1) -- (1');
\draw[|->, shorten <= 3pt, shorten >= 3pt] (2) -- (2');
\draw[|->, shorten <= 3pt, shorten >= 3pt] (3) -- (3');
\draw[|->, shorten <= 3pt, shorten >= 3pt] (4) -- (4');
\draw[|->, shorten <= 3pt, shorten >= 3pt] (5) -- (5');
\draw[|->, shorten <= 3pt, shorten >= 3pt] (6) -- (6');
    \begin{scope}[densely dotted, bend right=15pt]
      \draw[postaction={decorate}] (1'1') to (12');
      \draw[postaction={decorate}] (1'1) to (14);
      \draw[postaction={decorate}] (1'2) to (14);
      \draw[postaction={decorate}] (1'3) to (13);
      \draw[postaction={decorate}] (2'1') to (22');
      \draw[postaction={decorate}] (2'1) to (23);
      \draw[postaction={decorate}] (3'1') to (32');
      \draw[postaction={decorate}] (3'1) to (33);
      \draw[postaction={decorate}] (3'2) to (33);
      \draw[postaction={decorate}] (3'3) to (32);
      \draw[postaction={decorate}] (4'1') to (42');
      \draw[postaction={decorate}] (4'1) to (42);
    \end{scope}

	\end{scope}
  \end{tikzpicture}
\]
\end{solution}
\end{exercise}

\begin{exercise}
Suppose $p$, $q$, and $r$ are polynomials and you're given arbitrary lenses $\varphi\colon q\to p\tri q$ and $\psi\colon q\to q\tri r$. Does the following diagram necessarily commute?\tablefootnote{When the name of an object is used in place of a morphism, we refer to the identity morphism on that object.
So for instance, $\varphi\tri r$ is the composition product of $\varphi$ with the identity lens on $r$.}
\[
\begin{tikzcd}
	q\ar[r, "g"]\ar[d, "\varphi"']&
	q\tri r\ar[d, "\varphi\:\tri\:r"]\\
	p\tri q\ar[r, "p\:\tri\:\psi"']&
	p\tri q\tri r\ar[ul, phantom, "?"]
\end{tikzcd}
\]
That is, do we have $\varphi\then (p\tri \psi)=^?\psi\then (\varphi\tri r)$?
\begin{solution}
Given arbitrary polynomials $p,q,r$ and lenses $\varphi\colon q\to p\tri q$ and $\psi\colon q\to q\tri r$, it is \emph{not} necessarily the case that $\varphi\then (p\tri \psi)=\psi\then (\varphi\tri r)$!
After all, we can let $p\coloneqq\yon$ and $q\coloneqq\2$ so that $\varphi$ is a lens $\2\to\yon\tri\2\iso\2$ (see \cref{ex.compose_yon}) and $\psi$ is a lens $\2\to\2\tri r\iso\2$ (see \cref{exc.composing_with_constants}).
Then by following the instructions for interpreting a composition product of lenses from either \cref{exc.comp_prod_lens} or \cref{ex.comp_prod_trees}, we can verify that $p\tri \psi=\yon\tri \psi$ is a lens $\2\iso\yon\tri\2\to\yon\tri\2\tri r\iso\2$ equivalent to the lens $\psi$, while $\varphi\tri r$ is a lens $\2\iso\2\tri r\to\yon\tri\2\tri r\iso\2$ equivalent to the lens $\varphi$.
If, say, we let $\varphi\colon\2\to\2$ be the function sending everything to $1\in\2$ and $\psi\colon\2\to\2$ be the function sending everything to $2\in\2$, then in this case $\varphi\then (p\tri \psi)=\varphi\then \psi\neq \psi\then \varphi=\psi\then (\varphi\tri r)$.
\end{solution}
\end{exercise}

\subsection{Dynamical systems and the composition product} \label{subsec.comon.comp.def.dyn_sys}

\index{dynamical system!composition product and}
\index{dynamical system!speeding up|see{dynamical system, composition product and}}

Back in \cref{ex.do_nothing}, we posed the question of how to model running multiple steps of dynamical system in $\poly$.
The answer lies with the composition product.

\index{interface}

Recall that a \emph{(dependent) dynamical system} is a lens $\varphi\colon S\yon^S\to p$, where $S$ is a set of \emph{states} and $p$ is a polynomial \emph{interface}.
We call $S\yon^S$ the \emph{state system} and the on-position and on-direction functions of $\varphi$ the \emph{return} and \emph{update} functions, respectively.
More generally, we saw in \cref{ex.do_nothing} that we could replace the state system with a monomial $q\coloneqq S\yon^{S'}$, where $S'$ is another set, as long as there is a function $e\colon S\to S'$ (or equivalently a section $\epsilon\colon S\yon^{S'}\to\yon$) that is bijective.

\index{state system}\index{dynamical system}

The lens models a dynamical system as follows.
Every state $s\in q(\1)=S$ returns a position $o\coloneqq\varphi_\1(s)\in p(\1)$, and every direction $i\in p[o]$ yields an updated direction $s'\coloneqq\varphi^\sharp_s(a)\in q[s]=S'$.
Then to model a second step through the system, we identify the $q[s]$-direction $s'$ with a $q$-position $e^{-1}(s')$, plug this position back into $\varphi_\1$, and repeat the process all over again.

% Equivalently, we could think of running through two steps of the system as picking an initial state $s\in q(\1)$ and using a bijection $\ol{t}\colon q[s]\to q(\1)$ to assign each direction at $s$ to a state in $q(\1)$.
% Then we apply $\varphi$ in two rounds.
% First, we use the return function to compute an output $i\coloneqq\varphi_\1(s)$ and the update function $\varphi^\sharp_s\colon p[i]\to q[s]$ to compute a direction at $s$ for each input.
% Before the next round, we use $\ol{t}$ to compute the state assigned to that direction.
% Second, we use the return function to compute an output for that state

But this is exactly what the composition product $\varphi\tri\varphi\colon q\tri q\to p\tri p$ does: by \eqref{eqn.comp_lens_pos}, its on-positions function sends the pair $(s_0,e^{-1})\in(q\tri q)(\1)$, comprised of an initial state $s_0\in q(\1)$ and the function $e^{-1}\colon q[s_0]=S'\to S=q(\1)$, to the pair
\begin{equation} \label{eqn.interface_position_pair}\index{interface}
    \left(\varphi_\1(s_0),\varphi^\sharp_{s_0}\then e^{-1}\then\varphi_\1\right)\in(p\tri p)(\1),
\end{equation}
comprised of the initial position $o_0\coloneqq\varphi_\1(s)\in p(\1)$ and a composite function
\begin{equation} \label{eqn.interface_position_fn1}
    p[o_0]\To{\varphi^\sharp_{s_0}}q[s_0]=S'\To{e^{-1}}S=q(\1)\To{\varphi_\1}p(\1),
\end{equation}
which uses the update function at $s_0$ and the return function to tell us what the next position $o_1$ will be for every possible direction $i_1$ we could select.
Then by \eqref{eqn.comp_lens_dir}, the on-directions function of $\varphi\tri\varphi$ sends each direction $(i_1,i_2)$ at the position \eqref{eqn.interface_position_pair}, comprised of an initial direction $i_1\in p[o_0]$ and (setting $o_1$ to be the function \eqref{eqn.interface_position_fn1} applied to $i_1$) a second direction $i_2\in p[o_1]$, to the pair
\[
    \left(\varphi^\sharp_{s_0}(i_1), \varphi^\sharp_{e^{-1}(\varphi^\sharp_{s_0}(i_1))}(i_2)\right),
\]
comprised of directions in $S'$ that (under $e^{-1}$) correspond to the next state $s_1$ upon selecting direction $i_1$ at state $s_0$ and the successive state $s_2$ upon selecting direction $i_2$ at $s_1$.
In summary, at certain positions, $\varphi\tri\varphi$ tells us how the dynamical system will behave when we step through it twice: starting from state $s_0$, returning position $o_0$, receiving direction $i_1$, updating its state to $s_1$, returning position $o_1$, receiving direction $i_2$, and preparing to update its state to $s_2$.
Adding another layer, $\varphi\tri\varphi\tri\varphi\colon q\tri q\tri q\to p\tri p\tri p$ will tell us how the system behaves when we step through it three times; and in general, $\varphi\tripow{n}\colon q\tripow{n}\to p\tripow{n}$ will tell us how the system behaves when we step through it $n$ times.

\index{dynamical system!stepping through}

\begin{example}[Substitution products of dynamical systems as trees] \label{ex.comp_dyn_sys_tree}
Consider the dynamical system $\varphi\colon S\yon^S\to p$ with $p\coloneqq\yon^A+\1$, corresponding to the halting deterministic state automaton \eqref{eqn.halt_dsa} from \cref{exc.halt_dsa}, depicted again here for convenience:
\[
  \begin{tikzcd}[column sep=small]
    {\color{my-lavender}\blacksquare}
      \ar[rr, bend left, my-red, dashed]
      \ar[loop left, my-blue]
      &&
    {\color{my-yellow}\blacktriangle}
      \ar[dl, bend left, my-red, dashed]\ar[ll, my-blue, bend left]
      \\&
    {\color{my-magenta}\blacklozenge}
  \end{tikzcd}
\]
Below, we draw the corolla pictures for $S\yon^S$ and for $p$.
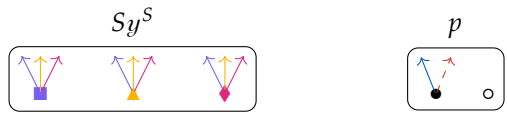
\[
\begin{tikzpicture}[rounded corners]
	\node (p1) [draw, "$S\yon^S$" above] {
	    \begin{tikzpicture}[trees, sibling distance=2.5mm]
            \node[my-lavender, below] (1) {\small$\blacksquare$}
              child[my-lavender] {}
              child[my-yellow] {}
              child[my-magenta] {};
            \node[my-yellow, right=1 of 1] (2) {\small$\blacktriangle$}
              child[my-lavender] {}
              child[my-yellow] {}
              child[my-magenta] {};
            \node[my-magenta, right=1 of 2] (3) {\small$\blacklozenge$}
              child[my-lavender] {}
              child[my-yellow] {}
              child[my-magenta] {};
        \end{tikzpicture}
    };
	\node (p2) [draw, right=2 of p1, "$p$" above] {
        \begin{tikzpicture}[trees, sibling distance=4mm]
            \node (1) {$\bullet$}
                child[my-blue] {}
                child[my-red,dashed] {};
            \node[right=0.5 of 1] (2) {$\circ$};
            ;
        \end{tikzpicture}
    };
\end{tikzpicture}
\]
In the picture for $S\yon^S$, the roots are the three states in
$\{
  {\color{my-lavender}\blacksquare},
  {\color{my-yellow}\blacktriangle},
  {\color{my-magenta}\blacklozenge}
\}$
appearing in the automaton, while the leaves of each corolla correspond to the three states as well, from left to right.
In the picture for $p$, there is one corolla whose two leaves correspond to the two arrows coming out of every state---except for the halting state, which is sent to the corolla with no leaves instead.
So the lens $\varphi\colon S\yon^S\to p$ capturing the dynamics of the automaton can be drawn as follows:
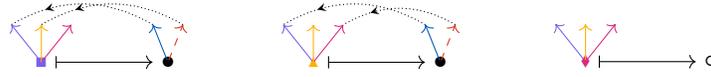
\[
\begin{tikzpicture}
	\node (p1) {
	\begin{tikzpicture}[trees, sibling distance=4mm]
        \node[my-lavender] (1) {\tiny$\blacksquare$}
            child[my-lavender]     {coordinate (11)}
            child[my-yellow]  {coordinate (12)}
            child[my-magenta]      {coordinate (13)};
        \node[right=1.5 of 1] (2) {$\bullet$}
            child[my-blue] {coordinate (21)}
            child[my-red,dashed] {coordinate (22)};
        \draw[|->, shorten <= 3pt, shorten >= 3pt] (1) -- (2);
        \begin{scope}[densely dotted, bend right]
          \draw[postaction={decorate}] (21) to (11);
          \draw[postaction={decorate}] (22) to (12);
        \end{scope}
    \end{tikzpicture}
	};
	\node [right=1 of p1] (p2) {
	\begin{tikzpicture}[trees, sibling distance=4mm]
        \node[my-yellow] (1) {\tiny$\blacktriangle$}
            child[my-lavender]     {coordinate (11)}
            child[my-yellow]  {coordinate (12)}
            child[my-magenta]      {coordinate (13)};
        \node[right=1.5 of 1] (2) {$\bullet$}
            child[my-blue] {coordinate (21)}
            child[my-red,dashed] {coordinate (22)};
        \draw[|->, shorten <= 3pt, shorten >= 3pt] (1) -- (2);
        \begin{scope}[densely dotted, bend right]
          \draw[postaction={decorate}] (21) to (11);
          \draw[postaction={decorate}] (22) to (13);
        \end{scope}
    \end{tikzpicture}
	};
	\node [right=1 of p2, yshift=-5pt] (p3) {
	\begin{tikzpicture}[trees, sibling distance=4mm]
        \node[my-magenta] (1) {\tiny$\blacklozenge$}
            child[my-lavender]     {coordinate (11)}
            child[my-yellow]  {coordinate (12)}
            child[my-magenta]      {coordinate (13)};
        \node[right=1.5 of 1] (2) {$\circ$};
        \draw[|->, shorten <= 3pt, shorten >= 3pt] (1) -- (2);
    \end{tikzpicture}
	};
\end{tikzpicture}
\]
The corolla picture tells us, for example, that from the ${\color{my-yellow}\blacktriangle}$ state, we can go in one of two directions: the solid blue direction, which leads us to the ${\color{my-lavender}\blacksquare}$ state, or the dashed red direction, which leads us to the ${\color{my-magenta}\blacklozenge}$ state.
This describes the dynamics of the automaton one step away from the ${\color{my-yellow}\blacktriangle}$ state.

But what if we want to understand the dynamics of $\varphi$ \emph{two} steps away from the ${\color{my-yellow}\blacktriangle}$ state?
Consider the following tree, corresponding to a position of the $2$-fold composite $S\yon^S\tri S\yon^S$:
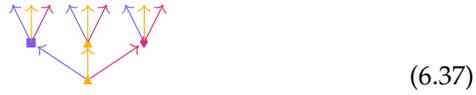
\begin{equation} \label{eqn.sys_tri2_good_tree}
\begin{tikzpicture}[trees,
  level 1/.style={sibling distance=7.5mm},
  level 2/.style={sibling distance=2.5mm}]
	\node[my-yellow] (a) {\tiny$\blacktriangle$}
		child[my-lavender] {node[my-lavender] {\tiny$\blacksquare$}
			child[my-lavender]
			child[my-yellow]
			child[my-magenta]
		}
		child[my-yellow] {node[my-yellow] {\tiny$\blacktriangle$}
			child[my-lavender]
			child[my-yellow]
			child[my-magenta]
		}
		child[my-magenta] {node[my-magenta] {\tiny$\blacklozenge$}
			child[my-lavender]
			child[my-yellow]
			child[my-magenta]
		};
\end{tikzpicture}
\end{equation}
We can follow the steps from \cref{ex.comp_prod_trees} to find that the composition product of lenses $\varphi\tri\varphi\colon S\yon^S\tri S\yon^S\to p\tri p$ acts on this tree as follows:
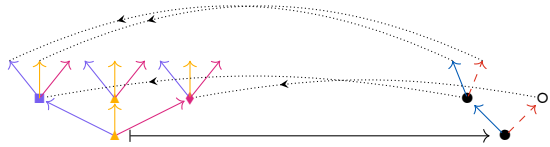
\[
\begin{tikzpicture}[trees,
  level 1/.style={sibling distance=10mm},
  level 2/.style={sibling distance=4mm}]
    \node[my-yellow] (1) {\tiny$\blacktriangle$}
		child[my-lavender] {node[my-lavender] (11) {\tiny$\blacksquare$}
			child[my-lavender] {coordinate (111)}
			child[my-yellow] {coordinate (112)}
			child[my-magenta] {coordinate (113)}
		}
		child[my-yellow] {node[my-yellow] (12) {\tiny$\blacktriangle$}
			child[my-lavender] {coordinate (121)}
			child[my-yellow] {coordinate (122)}
			child[my-magenta] {coordinate (123)}
		}
		child[my-magenta] {node[my-magenta] (13) {\tiny$\blacklozenge$}
			child[my-lavender] {coordinate (131)}
			child[my-yellow] {coordinate (132)}
			child[my-magenta] {coordinate (133)}
		};
    \node[right=5 of 1] (2) {$\bullet$}
        child[my-blue] {node[black, sibling distance=5mm] (21) {$\bullet$}
            child[my-blue] {coordinate (211)}
            child[my-red,dashed] {coordinate (212)}
        }
        child[my-red,dashed] {node[black] (22) {$\circ$}
        };

	\draw[|->, shorten <= 3pt, shorten >= 3pt] (1) -- (2);
    \begin{scope}[densely dotted, bend right=10pt]
      \draw[postaction={decorate}] (21) to (11);
      \draw[postaction={decorate}] (22) to (13);
    \end{scope}
    \begin{scope}[densely dotted, bend right=25pt]
      \draw[postaction={decorate}] (211) to (111);
      \draw[postaction={decorate}] (212) to (112);
    \end{scope}
\end{tikzpicture}
\]
Read each tree the way you would read a decision tree, and you will find that this picture tells you exactly what the dynamics of the automaton are two steps away from the ${\color{my-yellow}\blacktriangle}$ state.
Actually, it says that if we start from the ${\color{my-yellow}\blacktriangle}$ state (the root on the left, which is sent to the root on the right) and go in the dashed red direction (up from the root along the dashed red arrow to the right), the automaton will halt (as there are no more directions to follow).
But if we instead go in the solid blue direction (up from the root along the solid blue arrow to the left), we could go in the solid blue direction again (up the next solid blue arrow) to arrive at the ${\color{my-lavender}\blacksquare}$ state (as indicated by the dotted arrow above), or instead in the dashed red direction to arrive at the ${\color{my-yellow}\blacktriangle}$ state (similarly).
\end{example}
\index{decision tree}

\index{time!composition product and}

\index{dynamical system!trees}\index{interface}

This is what we meant at the start of this chapter when we said that the substitution product has to do with time.
It takes a specification $\varphi$ for how a state system and an interface can interact back-and-forth---or, indeed, any interaction pattern between wrapper interfaces---and extends it to a multistep model $\varphi\tripow{n}$ that simulates $n$ successive interaction cycles over time, accounting for all possible external directions that the interface could encounter.
Alternatively, we can think of $\varphi\tripow{n}$ as ``speeding up'' the original dynamical system $\varphi$ by a factor of $n$, as it runs $n$ steps in one---as long as whatever's connected to its new interface $p\tripow{n}$ can keep up with its pace and feed it $n$ directions of $p$ at a time!
The lens $\varphi$ tells us how the machine can run, but it is $\tri$ that makes the clock tick.
\index{clock!tick}\index{interaction!multi-step}

\subsubsection{Why this is not enough}\label{subsubsec.comon.comp.def.dyn_sys.issues}
There are several pressing issues we must address, however, before we can even begin to provide a satisfying answer to everything we asked for in \cref{ex.do_nothing}.
The first is a communication issue: as you probably sensed, our set-theoretic notation for $\varphi\tripow{n}$ is rather cumbersome, and that was just for $n=2$.
We could depict the behavior of our composition products of dynamical systems more clearly using tree pictures (see \cref{ex.comp_dyn_sys_tree}), but even that becomes infeasible in greater generality.
A concise visual representation of the back-and-forth interaction of lenses would help us reason about composition products more effectively.

The second issue is more technical: to ensure that $\varphi\tri\varphi$ behaves the way we want, when we specify a position of its domain $q\tri q$, we have to provide not only an initial state $s_0$ but also the isomorphism $e^{-1}\colon S'\to S$ to let the lens know which state in $S$ each $q[s_0]$-direction in $S'$ should lead to.
But there are many other positions $(s_0,f)$ of $q\tri q$ that we could have specified, and each $f\colon S'\to S$ associates $q[s_0]$-directions to $q$-positions in a different way.
So $\varphi\tri\varphi$ is carrying around a lot of extraneous---even misleading!---data about how our dynamical system behaves when the state system moves in the right direction but to the wrong state.
Our isomorphism $e^{-1}$ is a temporary fix, but as we pointed out in \cref{ex.do_nothing}, it relies on the set-theoretic equality of the direction-sets of $q$---there's nothing inherent to the categorical structure of $q$ in $\poly$ that encodes how directions map to states, at least not yet.
What are we missing?

\begin{example}[Composition products of dynamical systems can be misleading] \label{ex.comp_dyn_sys_mislead}
In \cref{ex.comp_dyn_sys_mislead}, instead of \eqref{eqn.sys_tri2_good_tree}, we could have picked the following tree, corresponding to a different (yet entirely valid) position of $S\yon^S\tri S\yon^S$:
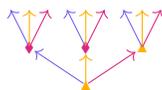
\[
\begin{tikzpicture}[trees,
  level 1/.style={sibling distance=7.5mm},
  level 2/.style={sibling distance=2.5mm}]
	\node[my-yellow] (a) {\tiny$\blacktriangle$}
		child[my-lavender] {node[my-magenta] {\tiny$\blacklozenge$}
			child[my-lavender]
			child[my-yellow]
			child[my-magenta]
		}
		child[my-yellow] {node[my-magenta] {\tiny$\blacklozenge$}
			child[my-lavender]
			child[my-yellow]
			child[my-magenta]
		}
		child[my-magenta] {node[my-yellow] {\tiny$\blacktriangle$}
			child[my-lavender]
			child[my-yellow]
			child[my-magenta]
		};
\end{tikzpicture}
\]
The composition product $\varphi\tri\varphi\colon S\yon^S\tri S\yon^S\to p\tri p$ acts on this tree like so:
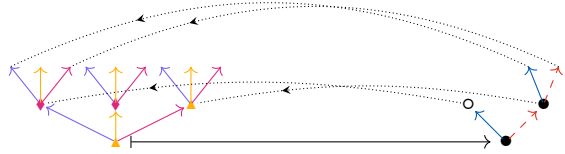
\[
\begin{tikzpicture}[trees,
  level 1/.style={sibling distance=10mm},
  level 2/.style={sibling distance=4mm}]
    \node[my-yellow] (1) {\tiny$\blacktriangle$}
    child[my-lavender] {node[my-magenta] (11) {\tiny$\blacklozenge$}
      child[my-lavender] {coordinate (111)}
      child[my-yellow] {coordinate (112)}
      child[my-magenta] {coordinate (113)}
    }
    child[my-yellow] {node[my-magenta] (12) {\tiny$\blacklozenge$}
      child[my-lavender] {coordinate (121)}
      child[my-yellow] {coordinate (122)}
      child[my-magenta] {coordinate (123)}
    }
    child[my-magenta] {node[my-yellow] (13) {\tiny$\blacktriangle$}
      child[my-lavender] {coordinate (131)}
      child[my-yellow] {coordinate (132)}
      child[my-magenta] {coordinate (133)}
    };
    \node[right=5 of 1] (2) {$\bullet$}
        child[my-blue] {node[black] (21) {$\circ$}
        }
        child[my-red,dashed] {node[black] (22) {$\bullet$}
            child[my-blue,solid] {coordinate (221)}
            child[my-red,dashed] {coordinate (222)}
        };

	\draw[|->, shorten <= 3pt, shorten >= 3pt] (1) -- (2);
    \begin{scope}[densely dotted, bend right=10pt]
      \draw[postaction={decorate}] (21) to (11);
      \draw[postaction={decorate}] (22) to (13);
    \end{scope}
    \begin{scope}[densely dotted, bend right=25pt]
      \draw[postaction={decorate}] (221) to (111);
      \draw[postaction={decorate}] (222) to (113);
    \end{scope}
\end{tikzpicture}
\]
Now the picture tells us that, starting from the ${\color{my-yellow}\blacktriangle}$ state, it is the solid blue arrow that will lead to a halting state, whereas we should somehow be able to follow the dashed red arrow twice!

Of course, this is nonsense---stemming from the fact that we have grafted the ``wrong'' corollas to each leaf when forming the position of $S\yon^S\tri S\yon^S$ on the left.
It is important to note, however, that this is nevertheless part of the data of $\varphi\tri\varphi$, and we don't yet know how to tell $\poly$ to rule it out.
\end{example}

The key to resolving both these issues lies in the next section, where we will study lenses whose codomains are composite polynomials.

%-------- Section --------%
\section{Lenses to composites}\label{sec.comon.comp.to_comp}
\index{lens!to composite}\index{polybox!maps to composition product}

Lenses to composites---that is, lenses of the form $f\colon p\to q_1\tri\cdots\tri q_n$ for some $n\in\nn$ with composites as their codomains---will be ubiquitous in the remainder of our story.
Fortunately, they have some very nice properties that make them convenient to work with, especially using polyboxes.

\subsection{Lenses to composites as polyboxes}
\index{polybox|(}
A lens $p\to q_1\tri q_2$ is an element of the set
\begin{align*}
    \poly(p, q_1\tri q_2) &\iso \poly\left(p, \sum_{j_1\in q_1(\1)}\;\prod_{b_1\in q_1[j_1]}\;\sum_{j_2\in q_2(\1)}\;\prod_{b_2\in q_2[j_2]}\yon\right) \tag*{\eqref{eqn.composite_formula}} \\
    &\iso \prod_{i\in p(\1)}\;\sum_{j_1\in q_1(\1)}\;\prod_{b_1\in q_1[j_1]}\;\sum_{j_2\in q_2(\1)}\;\prod_{b_2\in q_2[j_2]}p[i]. \tag*{\eqref{eqn.main_formula}}
\end{align*}
So we can write down the instructions for picking a lens $p\to q_1\tri q_2$ as follows.
\begin{quote}
To choose a lens $p\to q_1\tri q_2$:
\begin{enumerate}
    \item for each $p$-position $i$:
    \begin{enumerate}[label*=\arabic*.]
        \item choose a $q_1$-position $j_1$;
        \item for each $q_1[j_1]$-direction $b_1$:
        \begin{enumerate}[label*=\arabic*.]
            \item choose a $q_2$-position $j_2$;
            \item for each $q_2[j_2]$-direction $b_2$:
            \begin{enumerate}[label*=\arabic*.]
                \item choose a $p[i]$-direction $a$.
            \end{enumerate}
        \end{enumerate}
    \end{enumerate}
\end{enumerate}
\end{quote}
We could try to write out the dependent functions that these instructions correspond to.
Alternatively, we could simply draw this protocol out using polyboxes, with every ``for each'' step corresponding to a user-maintained blue box and every ``choose'' step corresponding to an automated white box:
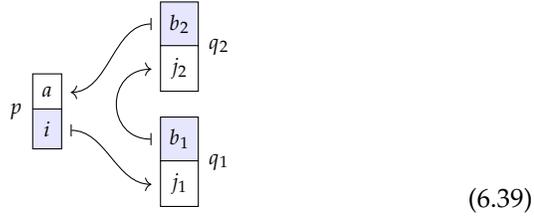
\begin{equation}\label{eqn.map_to_2ary_composite}
\begin{tikzpicture}[polybox, mapstos]
	\node[poly, dom, "$p$" left] (p) {$a$\at$i$};
	\node[poly, cod, right=1.5cm of p.south, yshift=-1ex, "$q_1$" right] (q1) {$b_1$\at$j_1$};
	\node[poly, cod, above=of q1, "$q_2$" right] (q2) {$b_2$\at$j_2$};
  	\draw (p_pos) to[first] (q1_pos);
  	\draw (q1_dir) to[climb] (q2_pos);
  	\draw (q2_dir) to[last] (p_dir);
\end{tikzpicture}
\end{equation}
Whenever we draw two pairs of polyboxes on top of each other, as we do with the polyboxes for $q_1$ and $q_2$ above on the right, we are indicating that the entire column of polyboxes depicts the composite of the polynomials depicted by each individual pair.
So the column of polyboxes on the right represents the composite $q_1\tri q_2$.
In particular, the position in the lower box of the top pair is the position associated with the direction in the upper box of the bottom pair, for the depicted position of the composite.

So a lens $p\to q_1\tri q_2$ is any protocol that will fill in the white boxes above as the user fills in the blue boxes in the direction of the arrows.
We'll see this in action in \cref{ex.map_to_comp}.

In fact, \eqref{eqn.map_to_0ary_composite}, \eqref{eqn.polybox_lens}, and \eqref{eqn.map_to_2ary_composite} are the respective polybox depictions of the $n=0, n=1,$ and $n=2$ cases of lenses $p\to q_1\tri\cdots\tri q_n$ to $n$-fold composites (we consider the monoidal unit $\yon$ of $\tri$ to be the $0$-fold composite, and a $1$-fold composite is just a polynomial on its own).
In general, for any $n\in\nn$, we can apply
\begin{align} \label{eqn.lens_to_comp}
    \poly(p, q_1\tri\cdots\tri q_n) &\iso\poly\left(p, \sum_{j_1\in q_1(\1)}\;\prod_{b_1\in q_1[j_1]}\cdots\sum_{j_n\in q_n(\1)}\;\prod_{b_n\in q_n[j_n]}\yon\right) \tag*{\eqref{eqn.composite_formula}} \\
    &\iso \prod_{i\in p(\1)}\;\sum_{j_1\in q_1(\1)}\;\prod_{b_1\in q_1[j_1]}\cdots\sum_{j_n\in q_n(\1)}\;\prod_{b_n\in q_n[j_n]}p[i], \tag*{\eqref{eqn.main_formula}}
\end{align}
so the polybox depiction of $p\to q_1\tri\cdots\tri q_n$ generalizes analogously.
For example, here are the polyboxes corresponding to a lens to a $4$-fold composite:
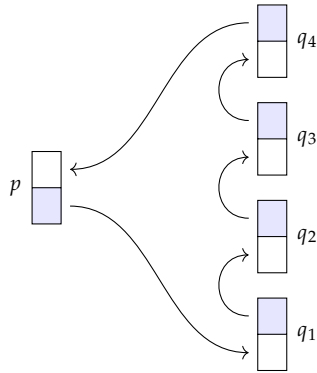
\[
\begin{tikzpicture}[polybox, tos]
	\node[poly, dom, "$p$" left] (p) {};
	\foreach \i in {1,...,4}
	{
  	\node[poly, cod, "$q_\i$" right] (q\i) at (3,1.3*\i-3.25) {};
	};
	\draw (p_pos) to[first] node[below] {} (q1_pos.west);
	\foreach \i/\j in {1/2,2/3,3/4}
	{
		\draw
			(q\i_dir.west)
			to[climb]
			node[left] {}
			(q\j_pos.west);
	};
	\draw (q4_dir) to[last] node[above left] {} (p_dir);
\end{tikzpicture}
\]

% We can use this to generalize our notation in the case $k=1$, i.e.\ for lenses $p\to q$. That is we denoted such a morphism by $\lens{f^\sharp}{f_\1}$, where $f_\1\colon p(\1)\to q(\1)$ and $f^\sharp_i\colon q[f_\1(i)]\to p[i]$. We generalize this to the $k$-ary composite case as
% \begin{equation}\label{eqn.notation_f1f2fk}
% (f_1,f_2,\ldots,f_k,f^\sharp)\colon p\too q_1\tri q_2\tri\cdots\tri q_k,
% \end{equation}
% where
% \begin{equation}\label{eqn.maps_to_comp}
% \begin{aligned}
% f_1&:p(\1)\to q_1(\1),\\
% f_2&:(i\in p(\1))\to (b_1\in q_1[f_1(i)])\to q_2(\1),\\
% f_3&:(i\in p(\1))\to (b_1\in q_1[f_1(i)])\to (b_2\in q_2[f_2(i,b_1)])\to q_3(\1),\\
% f_k&:(i\in p(\1))\to (b_1\in q_1[f_1(i)])\to  \cdots\to(b_{k-1}\in q_{k-1}[f_{k-1}(i,b_1,\ldots,b_{k-2})])\to q_k(\1),\\
% f^\sharp&:(i\in p(\1))\to (b_1\in q_1[f_1(i)])\to \cdots\to(b_{k}\in q_k[f_{k}(i,b_1,\ldots,b_{k-1})])\to p[i]
% \end{aligned}
% \end{equation}

These lenses to $n$-fold composites lend themselves to a very natural interpretation in terms of our decision-making language.
Each $p$-position is passed forward to a $q_1$-position.
For every direction at that $q_1$-position, there is then also a $q_2$-position.
Then for every direction at that $q_2$-position, there is a $q_3$-position, and so on, all the way until a direction has been selected at a $q_n$-position.
Together, all the directions chosen from $q_1,\ldots,q_n$ then inform the direction that is selected at the original $p$-position.

\index{decision!multi-step}
\index{polybox}

\slogan{A lens $p\to q_1\tri\cdots\tri q_n$ is a multi-step policy for $p$ to make decisions by asking for decisions from $q_1$, then $q_2$, etc., all the way to $q_n$, then interpreting the results.}

\begin{example}[Lenses $p\to q\tri r$]\label{ex.map_to_comp}
Consider a lens $\varphi\colon p\to q\tri r$.
Let's label the three arrows in the lens's polybox depiction:
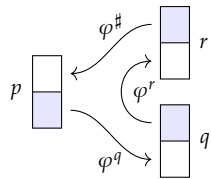
\[
\begin{tikzpicture}[polybox, tos]
	\node[poly, dom, "$p$" left] (p) {};
	\node[poly, cod, right=1.5cm of p.south, yshift=-1ex, "$q$" right] (q) {};
	\node[poly, cod, above=of q, "$r$" right] (r) {};
  	\draw (p_pos) to[first] node[below, yshift=-1mm] {$\varphi^q$} (q_pos);
  	\draw (q_dir) to[climb] node[right] {$\varphi^r$} (r_pos);
  	\draw (r_dir) to[last] node[above] {$\varphi^\sharp$} (p_dir);
\end{tikzpicture}
\]
So the on-position function of $\varphi$ can be split into two parts: a function $\varphi^q\colon p(\1)\to q(\1)$ and, for each $i\in p(\1)$, a function $\varphi^r_i\colon q[\varphi^q(i)]\to r(\1)$.
Then the on-directions function $\varphi^\sharp_i\colon (q\tri r)[\varphi_\1(i)]\to p[i]$ takes the direction of $q$ and the direction of $r$ in the two blue boxes on the right and sends them to a direction of $p$ at $i$ to fill the white box on the left.

For example, let $p\coloneqq\{A\}\yon^{\{R,S\}}+{B}\yon^{\{T\}}$, $q\coloneqq\{C\}\yon^{\{U,V,W\}}+\{D\}\yon^{\{X\}}$, and $r\coloneqq\{E\}\yon^{\{Y,Z\}}+\{F\}$.
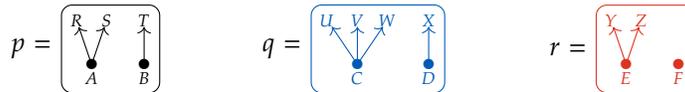
\[
\begin{tikzpicture}[rounded corners]
	\node (p) [draw, "$p=$" left] {
	\begin{tikzpicture}[trees, sibling distance=4mm]
    \node["\tiny $A$" below] (1) {$\bullet$}
      child {node[above, font=\tiny] {$R$}}
      child {node[above, font=\tiny] {$S$}};
    \node[right=.5 of 1,"\tiny $B$" below] (2) {$\bullet$}
      child {node[above, font=\tiny] {$T$}};
  \end{tikzpicture}
  };
	\node (q) [draw, my-blue, right=2 of p, "$q=$" left] {
	\begin{tikzpicture}[trees, sibling distance=4mm]
    \node["\tiny $C$" below] (1) {$\bullet$}
      child {node[above, font=\tiny] {$U$}}
      child {node[above, font=\tiny] {$V$}}
      child {node[above, font=\tiny] {$W$}};
    \node[right=.75 of 1,"\tiny $D$" below] (2) {$\bullet$}
      child {node[above, font=\tiny] {$X$}};
  \end{tikzpicture}
  };
	\node (r) [draw, my-red, right=2 of q, "$r=$" left] {
	\begin{tikzpicture}[trees, sibling distance=4mm]
    \node["\tiny $E$" below] (1) {$\bullet$}
      child {node[above, font=\tiny] {$Y$}}
      child {node[above, font=\tiny] {$Z$}};
    \node[right=.5 of 1,"\tiny $F$" below] (2) {$\bullet$};
  \end{tikzpicture}
  };
\end{tikzpicture}
\]
Here is a tree picture of a lens $\varphi\colon p\to q\tri r$:
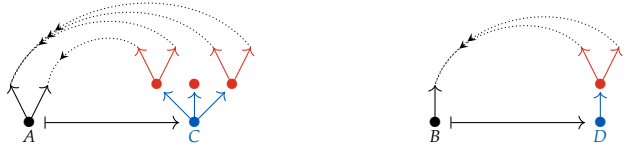
\[
	\begin{tikzpicture}[trees,
		level 1/.style={sibling distance=5mm},
	  level 2/.style={sibling distance=5mm}]
    \node["\tiny $A$" below] (1) {$\bullet$}
      child {coordinate (11')}
      child {coordinate (12')};
    \node[right=2 of 1, my-blue, "\color{my-blue} \tiny $C$" below] (1') {$\bullet$}
    	child[my-blue] {node[my-red] {$\bullet$}
				child[my-red] {coordinate (1'1)}
				child[my-red] {coordinate (1'2)}
			}
		child[my-blue] {node[my-red] {$\bullet$}
			}
		child[my-blue] {node[my-red] {$\bullet$}
				child[my-red] {coordinate (1'3)}
				child[my-red] {coordinate (1'4)}
			}
			;
    \node (2) [right=3 of 1',"\tiny $B$" below] {$\bullet$}
      child {coordinate (21')};
    \node[right=2 of 2, my-blue,"\color{my-blue} \tiny $D$" below] (2') {$\bullet$}
			child[my-blue] {node[my-red] {$\bullet$}
				child[my-red] {coordinate (2'1)}
				child[my-red] {coordinate (2'2)}
			}
			;
  \draw[|->, shorten <= 3pt, shorten >= 3pt] (1) -- (1');
  \draw[|->, shorten <= 3pt, shorten >= 3pt] (2) -- (2');
  \begin{scope}[densely dotted, bend right=60pt]
  	\draw[postaction={decorate}] (1'1) to (12');
  	\draw[postaction={decorate}] (1'2) to (11');
  	\draw[postaction={decorate}] (1'3) to (11');
  	\draw[postaction={decorate}] (1'4) to (11');
  	\draw[postaction={decorate}] (2'1) to (21');
  	\draw[postaction={decorate}] (2'2) to (21');
  \end{scope}
\end{tikzpicture}
\]
If we write $\varphi$ as the corresponding triple $(\varphi^q, \varphi^r, \varphi^\sharp)$, then we have
\begin{gather*}
\varphi^q(A)=C,\quad \varphi^q(B)=D;\\
\varphi^r_A(U)=E,\quad \varphi^r_A(V)=F,\quad \varphi^r_A(W)=E;\\
\varphi^r_B(X)=E;\\
\varphi^\sharp_A(U,Y)=S,\quad \varphi^\sharp_A(U,Z)=R,\quad \varphi^\sharp_A(W,Y)=R,\quad \varphi^\sharp_A(W,Z)=R;\\
\varphi^\sharp_B(X,Y)=T,\quad \varphi^\sharp_B(X,Z)=T.
\end{gather*}
Polyboxes display the same data in a different format:
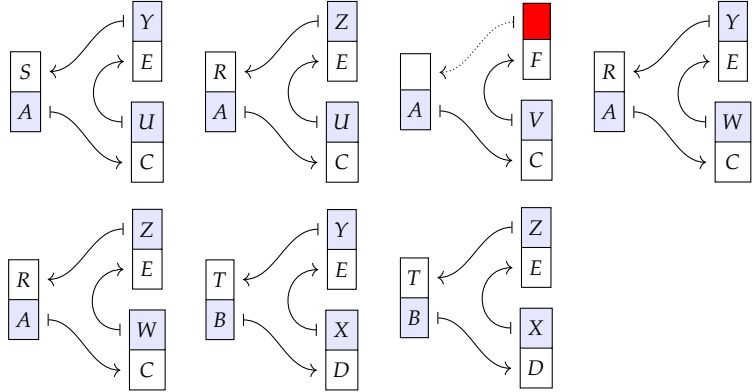
\[
\begin{tikzpicture}[polybox, mapstos, node distance=2ex and 1.4cm]
  \node (a) {
  \begin{tikzpicture}
  	\node[poly, dom] (p) {$S$\at$A$};
  	\node[poly, cod, right= of p.south, yshift=-1ex] (q) {$U$\at$C$};
  	\node[poly, cod, above=of q] (r) {$Y$\at$E$};
  	\draw (p_pos) to[first] (q_pos);
  	\draw (q_dir) to[climb] (r_pos);
  	\draw (r_dir) to[last] (p_dir);
  \end{tikzpicture}
  };
  \node[right=.3 of a] (b) {
  \begin{tikzpicture}
  	\node[poly, dom] (p) {$R$\at$A$};
  	\node[poly, cod, right= of p.south, yshift=-1ex] (q) {$U$\at$C$};
  	\node[poly, cod, above=of q] (r) {$Z$\at$E$};
  	\draw (p_pos) to[first] (q_pos);
  	\draw (q_dir) to[climb] (r_pos);
  	\draw (r_dir) to[last] (p_dir);
  \end{tikzpicture}
  };
  \node[right=.3 of b] (c) {
  \begin{tikzpicture}
  	\node[poly, dom] (p) {\at$A$};
  	\node[poly, cod, right= of p.south, yshift=-1ex] (q) {$V$\at$C$};
  	\node[poly, constant, above=of q] (r) {\at$F$};
  	\draw (p_pos) to[first] (q_pos);
  	\draw (q_dir) to[climb] (r_pos);
		\draw[densely dotted] (r_dir) to[last] (p_dir);
  \end{tikzpicture}
  };
  \node[right=.3 of c] (d) {
  \begin{tikzpicture}
  	\node[poly, dom] (p) {$R$\at$A$};
  	\node[poly, cod, right= of p.south, yshift=-1ex] (q) {$W$\at$C$};
  	\node[poly, cod, above=of q] (r) {$Y$\at$E$};
  	\draw (p_pos) to[first] (q_pos);
  	\draw (q_dir) to[climb] (r_pos);
  	\draw (r_dir) to[last] (p_dir);
  \end{tikzpicture}
	};
  \node[below=.1 of a] (e) {
  \begin{tikzpicture}
  	\node[poly, dom] (p) {$R$\at$A$};
  	\node[poly, cod, right= of p.south, yshift=-1ex] (q) {$W$\at$C$};
  	\node[poly, cod, above=of q] (r) {$Z$\at$E$};
  	\draw (p_pos) to[first] (q_pos);
  	\draw (q_dir) to[climb] (r_pos);
  	\draw (r_dir) to[last] (p_dir);
  \end{tikzpicture}
	};
  \node[below=.1 of b] (f) {
  \begin{tikzpicture}
  	\node[poly, dom] (p) {$T$\at$B$};
  	\node[poly, cod, right= of p.south, yshift=-1ex] (q) {$X$\at$D$};
  	\node[poly, cod, above=of q] (r) {$Y$\at$E$};
  	\draw (p_pos) to[first] (q_pos);
  	\draw (q_dir) to[climb] (r_pos);
  	\draw (r_dir) to[last] (p_dir);
  \end{tikzpicture}
	};
  \node[below=.1 of c] (g) {
  \begin{tikzpicture}
  	\node[poly, dom] (p) {$T$\at$B$};
  	\node[poly, cod, right= of p.south, yshift=-1ex] (q) {$X$\at$D$};
  	\node[poly, cod, above=of q] (r) {$Z$\at$E$};
  	\draw (p_pos) to[first] (q_pos);
  	\draw (q_dir) to[climb] (r_pos);
  	\draw (r_dir) to[last] (p_dir);
  \end{tikzpicture}
	};
\end{tikzpicture}
\]
As before, keep in mind that each arrow of a lens depends not only on the box it emerges from, but also on every box that came before it in our usual reading order (lower left to lower right to upper right to upper left).

\index{polybox}

The third set of polyboxes, where the left blue box has been filled with an $A$ and the lower right blue box has been filled with a $V$, is worth highlighting: as $\varphi^r_A(V)=F$, but $r[F]=\varnothing$, it is impossible to write a direction of $r$ at $F$ to go in the upper right box.
To indicate this, we color the upper right box red and leave the arrow emerging from it dashsed.
\end{example}

\index{dynamical system!composite interface}\index{interface!composite}
\begin{example}[Dynamical systems with composite interfaces] \label{ex.dyn_sys_comp_inter}
We explored dynamical systems with product interfaces in \cref{subsec.poly.dyn_sys.new.prod} and parallel product interfaces in \cref{subsec.poly.dyn_sys.new.par}.
How about dynamical systems with composite interfaces?
We now have all the tools we need to characterize them.

By the previous example, a dynamical system $\varphi\colon S\yon^S\to q\tri r$ can be drawn as
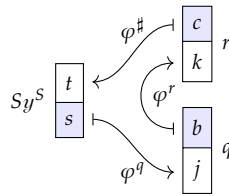
\[
\begin{tikzpicture}[polybox, mapstos]
	\node[poly, dom, "$S\yon^S$" left] (p) {$t$\at$s$};
	\node[poly, cod, right=1.5cm of p.south, yshift=-1ex, "$q$" right] (q) {$b$\at$j$};
	\node[poly, cod, above=of q, "$r$" right] (r) {$c$\at$k$};
  	\draw (p_pos) to[first] node[below, yshift=-1mm] {$\varphi^q$} (q_pos);
  	\draw (q_dir) to[climb] node[right] {$\varphi^r$} (r_pos);
  	\draw (r_dir) to[last] node[above] {$\varphi^\sharp$} (p_dir);
\end{tikzpicture}
\]
We can interpret the behavior of this system as follows.
Rather than a single interface $q\tri r$, we view $\varphi$ as having two interfaces that must be interacted with in succession, $q$ followed by $r$.\index{interface!composite}

Given the current state $s\in S$, the system feeds it into $\varphi^q$ to return a position $j$ of the first interface $q$.
Upon receiving a direction $b\in q[j]$, it then uses $\varphi^r$ to return another position $k$ (dependent on the state $s$ and the direction $b$), this time belonging to the second interface $r$.
Once a second direction $c$ is received, this time from $r[k]$, the system updates its state by feeding the current state $s$ and the pair of directions $(b,c)$ it received into $\varphi^\sharp$, yielding a new state $t\in S$.

That's a lot of words, which is why the polybox picture is so helpful: by following the arrows, we can see that a dynamical system with a composite interface actually captures a very natural type of interaction!
Mixing our metaphors a little, $\varphi$ could model a system that displays cascading menus, where selecting an option $b$ on the first menu $j$ opens up a second menu $k$.
It is only when the interacting agent selects an option $c$ from this second menu that both choices are sent back to the state system, which updates its state accordingly.

All this generalizes to $n$-fold composite interfaces exactly how you'd expect: a dynamical system with interface $q_1\tri\cdots\tri q_n$ returns a position and receives a direction through interface $q_1$, then accordingly returns a position and receives a direction through interface $q_2$, and so on, until it returns a position (according to the current state and all previous directions) and receives a direction through $q_n$, whereupon it updates its state according to the $n$ directions it received along with the current state.
\end{example}\index{interface!composite}

\subsection{The composition product of lenses as polyboxes}

A special case of a lens whose codomain is a composite is a lens that is itself the composition product of lenses.
If we draw such a lens using polyboxes by following the instructions from \eqref{eqn.comp_lens_pos} and \eqref{eqn.comp_lens_dir}, we would really just be stacking the polyboxes for the constituent lenses on top of each other.
For example, given lenses $\varphi\colon p\to q$ and $\varphi'\colon p'\to q'$, here is $\varphi\tri \varphi'$ drawn as polyboxes:
\begin{equation}
\label{eqn.comp_lens_polybox}
\begin{tikzpicture}[polybox, tos]
	\node[poly, dom, "$p$" left] (p) {};
	\node[poly, dom, above=.7 of p, "$p'$" left] (p') {};
	\node[poly, cod, right=of p, "$q$" right] (q) {};
	\node[poly, cod, above=.7 of q, "$q'$" right] (q') {};
	\draw (p_pos) -- node[below] {$\varphi_\1$} (q_pos);
	\draw (q_dir) -- node[above] {$\varphi^\sharp$} (p_dir);
	\draw (p'_pos) -- node[below] {$\varphi'_\1$} (q'_pos);
	\draw (q'_dir) -- node[above] {$(\varphi')^\sharp$} (p'_dir);
\end{tikzpicture}
\end{equation}
What differentiates this from simply writing down the polyboxes for $\varphi$ and the polyboxes for $\varphi'$ is that we are explictly associating the position that will fill the lower box of $p'$ with the direction that will fill the upper box of $p$, and likewise the position that will fill the lower box of $q'$ with the direction that will fill the upper box of $q$.
Moreover, we have the user fill out the lower set of boxes first and work their way up, so that, in particular, they can use the information they obtained from the behavior of $\varphi_\1$ and $\varphi^\sharp$ to decide what to put in the lower box of $p'$.
So this really does depict a lens $p\tri p'\to q\tri q'$.

\index{interaction}

How does \eqref{eqn.comp_lens_polybox} relate to our usual polybox depiction of a lens to a composite, as in \eqref{eqn.map_to_2ary_composite}, but with the domain also replaced with a composite?
A user who interacts with \eqref{eqn.comp_lens_polybox} can fill the lower set of polyboxes (the ones for $\varphi$) first, ignoring the upper set of polyboxes (the ones for $\varphi'$) until the entire lower half is filled.
Alternatively, after they fill in the lower box of $p$, but before they fill in anything else, they can already decide what position to put in the lower box of $p'$ for every possible direction that could end up in the upper box of $p$.
By \eqref{eqn.comp_pos}, such a choice is equivalent to picking a position of the composite $p\tri p'$.
Then by \eqref{eqn.comp_lens_pos}, following just the bottom arrow $\varphi_\1$ leads to the corresponding position of $q$ given by $\varphi\tri \varphi'$, while filling in the upper box of $q$ and following $\varphi^\sharp$, then $\varphi'_\1$ leads to the position of $q'$ that goes in the bottom box of $q'$.
Finally, once the user fills in the upper box of $q'$, following the top arrow $(\varphi')^\sharp$ completes the specification of a direction of $p\tri p'$.
In this way, \eqref{eqn.comp_lens_polybox} can be thought of as a special case of \eqref{eqn.map_to_2ary_composite}.

\begin{example}[Dynamical systems and the composition product, revisited] \label{ex.dyn_sys_comp_polyboxes}
In \cref{subsec.comon.comp.def.dyn_sys}, we explained how the $n$-fold composition product $\varphi\tripow{n}$ of a dynamical system $\varphi\colon S\yon^S\to p$ models the behavior of running through the system $n$ times, provided we choose the positions of $(S\yon^S)\tripow{n}$ appropriately.
We can visualize this behavior using polyboxes---for example, here's what the $n=3$ case looks like:
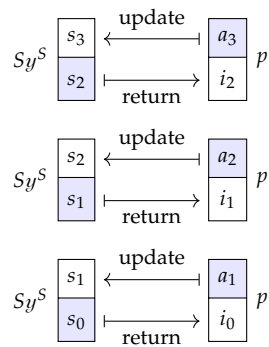
\[
\begin{tikzpicture}[polybox, mapstos]
	\node[poly, dom, "$S\yon^S$" left] (s) {$s_1$\at$s_0\vphantom{i_0}$};
	\node[poly, dom, above=.5 of s, "$S\yon^S$" left] (s') {$s_2$\at$s_1\vphantom{i_1}$};
	\node[poly, dom, above=.5 of s', "$S\yon^S$" left] (s'') {$s_3$\at$s_2\vphantom{i_2}$};

	\node[poly, cod, right=of s, "$p$" right] (p) {$a_1$\at$i_0$};
	\node[poly, cod, above=.5 of p, "$p$" right] (p') {$a_2$\at$i_1$};
	\node[poly, cod, above=.5 of p', "$p$" right] (p'') {$a_3$\at$i_2$};

	\draw (s_pos) -- node[below] {return} (p_pos);
	\draw (p_dir) -- node[above] {update} (s_dir);

	\draw (s'_pos) -- node[below] {return} (p'_pos);
	\draw (p'_dir) -- node[above] {update} (s'_dir);

	\draw (s''_pos) -- node[below] {return} (p''_pos);
	\draw (p''_dir) -- node[above] {update} (s''_dir);
\end{tikzpicture}
\]
It is now patently obvious what $\varphi\tripow3$ does from this picture, as long as we know how to read polyboxes (and we could probably make a pretty good guess even if we didn't!).
This resolves the first issue we raised in \cref{subsec.comon.comp.def.dyn_sys}, page~\pageref{subsubsec.comon.comp.def.dyn_sys.issues}: we now have a concise way of depicting the $n$-fold composite of a dynamical system.
The second issue becomes clear when we look at which boxes are blue along the left: we would really like the position $s_1$ to be entered above the direction $s_1$ automatically, the $s_2$ entered above $s_2$ automatically, etc.\ rather than having to specify the contents of those blue boxes manually.
We shouldn't even having the option to fill those blue boxes in with anything else.
We'll see how to address this issue shortly in \cref{ex.dyn_sys_trans_polyboxes}.
\end{example}

We make a big deal out of it, but \eqref{eqn.comp_lens_polybox} really is just the polyboxes of two separate lenses drawn together.
Where such polyboxes truly get interesting is when we compose them with polyboxes that look like \eqref{eqn.map_to_2ary_composite}.
That is, given a lens $g\colon r\to p\tri p'$, consider the polyboxes for $g\then(f\tri f')$:
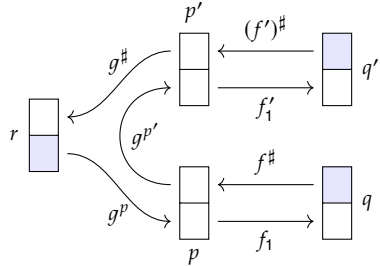
\[
\begin{tikzpicture}[polybox, tos]
	\node[poly, dom, "$r$" left] (r) {};
	\node[poly, right=1.8 of r.south, yshift=-2.5ex, "$p$" below] (p) {};
	\node[poly, above=.8 of p, "$p'$" above] (p') {};
	\node[poly, cod, right=of p, "$q$" right] (q) {};
	\node[poly, cod, above=.8 of q, "$q'$" right] (q') {};
	\draw (p_pos) -- node[below] {$f_\1$} (q_pos);
	\draw (q_dir) -- node[above] {$f^\sharp$} (p_dir);
	\draw (p'_pos) -- node[below] {$f'_\1$} (q'_pos);
	\draw (q'_dir) -- node[above] {$(f')^\sharp$} (p'_dir);
	\draw (r_pos) to[first] node[below, yshift=-1mm] {$g^{p}$} (p_pos);
	\draw (p_dir) to[climb] node[right] {$g^{p'}$} (p'_pos);
	\draw (p'_dir) to[last] node[above] {$g^\sharp$} (r_dir);
\end{tikzpicture}
\]
There is a lot going on with this lens! To fill out these polyboxes, we start from the lower box of $r$, go all the way right to the lower box of $q$, loop back left, up, and right again to the lower box of $q'$, then travel left all the way back to the upper box of $r$.

Say we knew that $g\then(f\tri f')$ were equal to some other lens $h\colon r\to q\tri q'$:
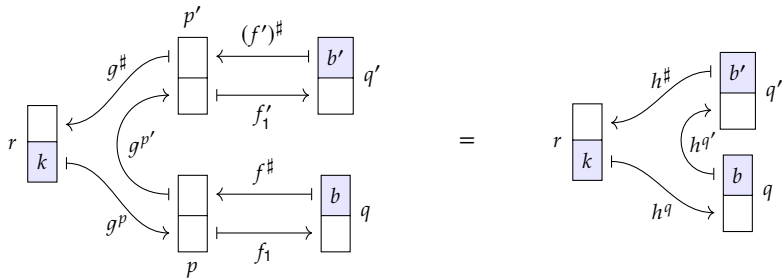
\[
\begin{tikzpicture}
	\node (1) {
  \begin{tikzpicture}[polybox, mapstos]
	\node[poly, dom, "$r$" left] (r) {\at$k$};
	\node[poly, right=1.8 of r.south, yshift=-2.5ex, "$p$" below] (p) {$\vphantom{b}$};
	\node[poly, above=.8 of p, "$p'$" above] (p') {$\vphantom{b'}$};
	\node[poly, cod, right=of p, "$q$" right] (q) {$b$};
	\node[poly, cod, above=.8 of q, "$q'$" right] (q') {$b'$};
	\draw (p_pos) -- node[below] {$f_\1$} (q_pos);
	\draw (q_dir) -- node[above] {$f^\sharp$} (p_dir);
	\draw (p'_pos) -- node[below] {$f'_\1$} (q'_pos);
	\draw (q'_dir) -- node[above] {$(f')^\sharp$} (p'_dir);
	\draw (r_pos) to[first] node[below,yshift=-1mm] {$g^p$} (p_pos);
	\draw (p_dir) to[climb] node[right] {$g^{p'}$} (p'_pos);
	\draw (p'_dir) to[last] node[above] {$g^\sharp$} (r_dir);
  \end{tikzpicture}
	};
	\node[right=1.8 of 1] (2) {
  \begin{tikzpicture}[polybox, mapstos]
  	\node[poly, dom, "$r$" left] (r) {\at$k$};
  	\node[poly, cod, right=1.8 of r.south, yshift=-1ex, "$q$" right] (q) {$b$};
  	\node[poly, cod, above=of q, "$q'$" right] (q') {$b'$};
  	\draw (r_pos) to[first] node[below,yshift=-1mm] {$h^q$} (q_pos);
  	\draw (q_dir) to[climb] node[right, xshift=-.5mm] {$h^{q'}$} (q'_pos);
  	\draw (q'_dir) to[last] node[above] {$h^\sharp$} (r_dir);
  \end{tikzpicture}
	};
	\node at ($(1.east)!.5!(2.west)$) {=};
\end{tikzpicture}
\]

We've filled in the corresponding blue boxes on either side of the equation with the same entries.
So if these two sets of polyboxes really do depict the same lens, each of the three white boxes in the domain and codomain on the left should end up with the same entry as the corresponding white box on the right (although the intermediary mechanics may differ).
Then if we follow the arrows in order on either side, matching up the white boxes in the domain and codomain along the way, we can read off three equations:
\[
    g^p\then f_\1 = h^q, \quad f^\sharp_{g^p(k)}\then g^{p'}_k\then f'_\1 = h^{q'}_k, \qand (f')^\sharp_{g^{p'}\left(f^\sharp_{g^p(k)}(b)\right)}\then g^\sharp_k = h^\sharp_k.
\]
The converse holds as well: if the three equations above all hold, then $g\then(f\tri f')=h$.
We will read equations off of polyboxes like this repeatedly in the rest of the book.%, including in the following example.

\begin{example}[Transition lenses for state systems] \label{ex.dyn_sys_trans_polyboxes}

In \cref{ex.dyn_sys_comp_polyboxes}, we saw that the $2$-fold composition product $\varphi\tripow2$ of a dynamical system $\varphi\colon S\yon^S\to p$ can be drawn as follows:
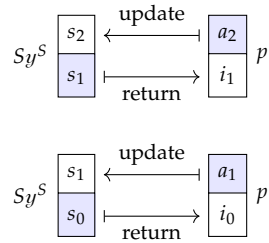
\[
\begin{tikzpicture}[polybox, mapstos]
	\node[poly, dom, "$S\yon^S$" left] (s) {$s_1$\at$s_0\vphantom{i_0}$};
	\node[poly, dom, above=.75 of s, "$S\yon^S$" left] (s') {$s_2$\at$s_1\vphantom{i_1}$};

	\node[poly, cod, right=of s, "$p$" right] (p) {$a_1$\at$i_0$};
	\node[poly, cod, above=.75 of p, "$p$" right] (p') {$a_2$\at$i_1$};

	\draw (s_pos) -- node[below] {return} (p_pos);
	\draw (p_dir) -- node[above] {update} (s_dir);

	\draw (s'_pos) -- node[below] {return} (p'_pos);
	\draw (p'_dir) -- node[above] {update} (s'_dir);
\end{tikzpicture}
\]
This \emph{almost} models the behavior of running through the system twice, except that we should really only have one blue box on the domain side---the one we fill with the initial state $s_0$.
The second blue box on the domain side, the one we fill with $s_1$, should instead be filled automatically with the same state as the direction $s_1$ in the white box below it.

In fact, it would be nice if the domain were still just $S\yon^S$.
Then we would have a lens $S\yon^S\to p\tri p$ that takes an initial state $s_0\in S$ and runs the original system $\varphi$ twice, returning two positions and receiving two directions before stopping at the new state $s_2$.
But we just learned how to take a composition product of lenses such as $\varphi\tripow2\colon S\yon^S\tri S\yon^S\to p\tri p$ and convert its domain to a new polynomial, say $S\yon^S$, with only one blue box on the domain side---just compose it with another lens $S\yon^S\to S\yon^S\tri S\yon^S$ like so:
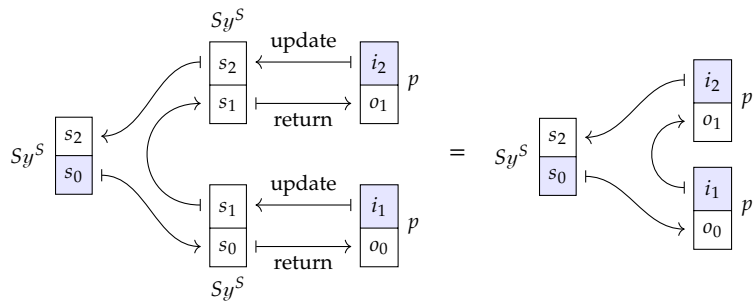
\[
\begin{tikzpicture}
	\node (1) {
  \begin{tikzpicture}[polybox, mapstos]
	\node[poly, dom, "$S\yon^S$" left] (r) {$s_2$\at$s_0$};
	\node[poly, right=1.8 of r.south, yshift=-2.5ex, "$S\yon^S$" below] (p) {$s_1\vphantom{i_1}$\at$s_0$};
	\node[poly, above=.8 of p, "$S\yon^S$" above] (p') {$s_2\vphantom{i_2}$\at$s_1$};
	\node[poly, cod, right=of p, "$p$" right] (q) {$i_1$\at$o_0\vphantom{s_0}$};
	\node[poly, cod, above=.8 of q, "$p$" right] (q') {$i_2$\at$o_1\vphantom{s_1}$};
	\draw (p_pos) -- node[below] {return} (q_pos);
	\draw (q_dir) -- node[above] {update} (p_dir);
	\draw (p'_pos) -- node[below] {return} (q'_pos);
	\draw (q'_dir) -- node[above] {update} (p'_dir);
	\draw (r_pos) to[first] node[below] {} (p_pos);
	\draw (p_dir) to[climb] node[right] {} (p'_pos);
	\draw (p'_dir) to[last] node[above] {} (r_dir);
  \end{tikzpicture}
	};
	\node[right=.5 of 1] (2) {
  \begin{tikzpicture}[polybox, mapstos]
  	\node[poly, dom, "$S\yon^S$" left] (r) {$s_2$\at$s_0$};
  	\node[poly, cod, right=1.8 of r.south, yshift=-1ex, "$p$" right] (q) {$i_1$\at$o_0$};
  	\node[poly, cod, above=of q, "$p$" right] (q') {$i_2$\at$o_1$};
  	\draw (r_pos) to[first] node[below] {} (q_pos);
  	\draw (q_dir) to[climb] node[right] {} (q'_pos);
  	\draw (q'_dir) to[last] node[above] {} (r_dir);
  \end{tikzpicture}
	};
	\node at ($(1.east)!.5!(2.west)$) {=};
\end{tikzpicture}
\]
We will denote our new lens $S\yon^S\to S\yon^S\tri S\yon^S$ on the far left by $\delta$.
Let us take a closer look at how $\delta$ behaves on its own, labeling its arrows for ease of reference:
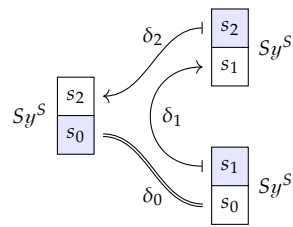
\[
\begin{tikzpicture}[polybox, mapstos]
	\node[poly, dom, "$S\yon^S$" left] (r) {$s_2$\at$s_0$};
	\node[poly, cod, right=1.8 of r.south, yshift=-2.5ex, "$S\yon^S$" right] (p) {$s_1$\at$s_0$};
	\node[poly, cod, above=.8 of p, "$S\yon^S$" right] (p') {$s_2$\at$s_1$};

	\draw[double, -] (r_pos) to[first] node[below,yshift=-1mm] {$\delta_0$} (p_pos);
	\draw (p_dir) to[climb] node[right] {$\delta_1$} (p'_pos);
	\draw (p'_dir) to[last] node[above,yshift=1mm] {$\delta_2$} (r_dir);
  \end{tikzpicture}
\]
We can see from this picture that $\delta$ arises very naturally from the structure of $S\yon^S$; indeed, every state system can be equipped with such a lens, just as every state system can be equipped with a do-nothing section from \cref{ex.do_nothing}.
The bottom arrow $\delta_0\colon S\to S$ sending $s_0\mapsto s_0$ is the identity on the position-set $S$: it sends each state to itself.
We use a double arrow to denote this equality.
While the middle arrow $\delta_1$ also looks like an identity arrow, remember that we should think of it as depending on the left blue box as well; so it is really a function $\delta_1\colon S\times S\to S$ sending $(s_0,s_1)$, a position-state $s_0$ and a direction at $s_0$ corresponding to state $s_1$, to the new position-state $s_1$.
Similarly, $\delta_2\colon S\times S\times S\to S$ sends $(s_0,s_1,s_2)$, where the last coordinate is the direction at position-state $s_1$ corresponding to state $s_2$, to the direction at position-state $s_0$ that also corresponds to state $s_2$.

Notice the crucial role that $\delta_1$ plays here: it matches up every direction at a given state to the new state that the direction in question should point to, encoding how the state system transitions from one state to the next.
We have been labeling each direction by the state that it points to, but we should really think of all the directions of $S\yon^S$ as pairs of position-states $(s,t)\in S\times S$, where $(s,t)$ is a direction at $s$ that represents the transition from state $s$ to state $t$.
The structure of the polynomial $S\yon^S$ tells us the state that each direction transitions from, but it is $\delta_1$ that tells us the new state $\delta_1(s,t)=t$ that $(s,t)$ transitions to.
For that reason, we call $\delta$ the \emph{transition lens} of the state system $S\yon^S$, and we call $\delta_1\colon S\times S\to S$ the \emph{target function}, relabeling it tgt for short, indicating where each direction leads.

This is exactly what we wanted back in \cref{ex.do_nothing}: a way to encode which directions of a state system point to which positions in the language of $\poly$. Expressions in that language are lenses such as $\delta$, and polyboxes like the ones above are the way we write them down.

We highlighted $\delta_1$ above, but the arrows $\delta_0$ and $\delta_2$ are no less important.
As an identity function, $\delta_0=\id_S$ remembers the initial state $s_0$ as we shift from its original setting $S\yon^S$, where we can only move one state away from $s_0$, to a new setting $S\yon^S\tri S\yon^S$, where we can think about all the ways to move two states away from $s_0$.
Meanwhile, $\delta_2$ shifts us from the two-state setting back to the one-state setting while ensuring a sort of transitive coherence condition: the state where we end up after moving from $s_0$ through $s_1$ to $s_2$ is the same state we would end up at if we had moved from $s_0$ directly to $s_2$.
We will call it the \emph{run function}, because it runs a sequence of two transitions together, and because it is what keeps the system running as it tells the original state system to actually move along one of its directions.

Here is another picture of the transition lens $\delta$ of $S\yon^S$ with our new arrow names:
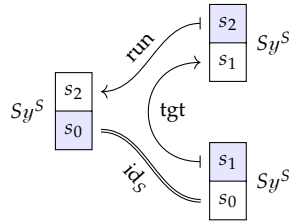
\[
\begin{tikzpicture}[polybox, mapstos]
	\node[poly, dom, "$S\yon^S$" left] (r) {$s_2$\at$s_0$};
	\node[poly, cod, right=1.8 of r.south, yshift=-2.5ex, "$S\yon^S$" right] (p) {$s_1$\at$s_0$};
	\node[poly, cod, above=.8 of p, "$S\yon^S$" right] (p') {$s_2$\at$s_1$};

	\draw[double, -] (r_pos) to[first] node[below,sloped] {$\id_S$} (p_pos);
	\draw (p_dir) to[climb] node[right] {tgt} (p'_pos);
	\draw (p'_dir) to[last] node[above,sloped] {run} (r_dir);
  \end{tikzpicture}
\]
This resolves the second issue we raised in \cref{subsec.comon.comp.def.dyn_sys}, page~\pageref{subsubsec.comon.comp.def.dyn_sys.issues} for the case of $n=2$, giving us a dynamical system $\delta\then(\varphi\tri\varphi)\colon S\yon^S\to p\tri p$ that simulates stepping through the system $\varphi$ twice.
But what about more general values of $n$?

We already have a way of talking about the $n=0$ case: that is what the do-nothing section $\epsilon\colon S\yon^S\to\yon$ models.
But there is some overlap in how $\epsilon$ matches up state-positions with directions and how $\delta$ does.
Put another way, if $\epsilon$ tells us how to do nothing, and $\delta$ tells us how to do two things, then we had better check that if one of the two things we do is nothing, then that's the same as doing just one thing.
Can we ensure that $\epsilon$ and $\delta$ agree on what they are saying about our state system?

Then for $n=3$, we would like a dynamical system $S\yon^S\to p\tri p\tri p$ that simulates stepping through $\varphi$ three times.
One way to do this would be to compose $\varphi\tripow3$ with a lens of the form $S\yon^S\to\left(S\yon^S\right)\tripow3$ that we obtain by extending $\delta$ from modeling two-step transitions to three-step transitions.
But there are two ways to derive a lens with codomain $\left(S\yon^S\right)\tripow3$ from $\delta$: we could either take $\delta\then\left(\delta\tri\id_{S\yon^S}\right)$, or we could take $\delta\then\left(\id_{S\yon^S}\tri\delta\right)$
For larger values of $n$, there are even more possibilities for what we could do.
But there should really only be one dynamical system that models stepping through $\varphi$ a fixed number of times.
How do we guarantee that all these different ways of extending $\delta$ to $n$-step transitions end up telling us the same thing?

In summary, we need some kind of compatibility condition between $\epsilon$ and $\delta$, as well as some kind of associativity condition on $\delta$ to guarantee that it can be extended coherently.
In fact, we already have all the tools we need to characterize these conditions: we'll see exactly how to state the properties we want in the next chapter.
And if this is all starting to sound suspiciously familiar, you’re not wrong—but we’ll save that surprise for the next chapter as well.
\end{example}\index{associativity}

\begin{exercise}
Let $S\coloneqq\nn$ and $p\coloneqq\rr\yon^\1$, and define $\varphi\colon S\yon^S\to p$ to be the dynamical system with return function $\varphi_\1(k)\coloneqq k$ and update function $\varphi^\sharp_k(1)\coloneqq k+1$.
\begin{enumerate}
    \item Draw the polyboxes for $\varphi$ and describe its dynamics: what does $1$ run through the system look like?
    \item Let $\delta\colon S\yon^S\to S\yon^S\tri S\yon^S$ be the transition lens of $S\yon^S$, and draw the polyboxes for the new system $\delta\then(\varphi\tri\varphi)\colon S\yon^S\to p\tri p$.
    Describe its dynamics: how does it model $2$ runs through the system? \qedhere
\end{enumerate}

\begin{solution}
Given $S\coloneqq\nn$ and $p\coloneqq\rr\yon^\1$, we have a dynamical system $\varphi\colon S\yon^S\to p$ given by $\varphi_\1(k)\coloneqq k$ and $\varphi^\sharp_k(1)\coloneqq k+1$.
\begin{enumerate}
    \item Here is the polybox picture for $\varphi$ (recall that we shade the upper box of a linear polynomial gray, as there is only one option to place there):
    \[
    \begin{tikzpicture}[polybox, mapstos]
		\node[poly, dom, "$S\yon^S$" left] (S) {$k+1$\at$k$};
		\node[poly, right=of S, linear, "$p$" right] (p) {$\vphantom{k}$\at$k$};
		\draw (S_pos) to[first] (p_pos);
		\draw (p_dir) to[last ] (S_dir);
    \end{tikzpicture}
    \]
    A single run through the system returns the current state $k\in\nn$, then increases that state by $1$.
    \item Here is the polybox picture for $\delta\then(\varphi\tri\varphi)\colon S\yon^S\to p\tri p$:
    \[
    \begin{tikzpicture}[polybox, mapstos]
		\node[poly, dom, "$S\yon^S$" left] (S) {$k+2$\at$k$};

    	\node[poly, cod, right=1.8 of S.south, yshift=-2.5ex, "$S\yon^S$" below] (S1) {$k+1$\at$k$};
    	\node[poly, cod, above=.8 of S1, "$S\yon^S$" above] (S2) {$k+2$\at$k+1$};

		\node[poly, linear, right=of S1, "$p$" right] (p1) {$\vphantom{k}$\at$k$};
		\node[poly, linear, right=of S2, "$p$" right] (p2) {$\vphantom{k}$\at$k+1$};
		\draw (S1_pos) to[first] (p1_pos);
		\draw (p1_dir) to[last] (S1_dir);
		\draw (S2_pos) to[first] (p2_pos);
		\draw (p2_dir) to[last]  (S2_dir);

		\draw (S_pos) to[first] (S1_pos);
		\draw (S1_dir) to[climb] (S2_pos);
		\draw (S2_dir) to[last] (S_dir);
    \end{tikzpicture}
    \]
    The new system has $p\tri p\iso\rr\yon\tri\rr\yon\iso\rr^\2\yon$ as its interface.
    Indeed, we see that it returns two numbers at once: the current state $k$ (what the first run through $\varphi$ would return) as well as the increased state $k+1$ (what the second run through $\varphi$ would return).
    We update the current state from $k$ to $k+1$ in one run, and from $k+1$ to $k+2$ in the next---thus increasing the current state by $2$ overall.
\end{enumerate}
\end{solution}
\end{exercise}\index{interface}

\begin{exercise}
As a lens whose domain is a state system, the transition lens $\delta\colon S\yon^S\to S\yon^S\tri S\yon^S$ of a state system $S\yon^S$ can be interpreted as a standalone dynamical system. Describe the dynamics of this system.
\begin{solution}
We give two reasonable (and of course equivalent) ways to interpret the transition lens $\delta\colon S\yon^S\to S\yon^S\tri S\yon^S$.

One way is to first evaluate its interface as $S\yon^S\tri S\yon^S\iso S\!\left(S\yon^S\right)^S\iso\left(S\times S^S\right)\yon^{S\times S}$.
Then we see that if the current state is $s_0\in S$, the system returns a position consisting of that current state $s_0$ along with a function $S\to S$, namely the identity function $s_1\mapsto s_1$.
The system then takes in a pair $(s_1,s_2)\in S\times S$, discarding $s_1$ and setting its new state to be $s_2$.

Alternatively, we can draw from \cref{ex.dyn_sys_comp_inter} to interpret $\delta$ as a dynamical system that behaves as follows.
Each run through the system is a $2$-step process: first, the current state $s_0\in S$ is returned, and a new state $s_1\in S$ is received.
Then this new state $s_1$ is immediately returned, and an ever newer state $s_2\in S$ is received.
Then the current state is updated to the newer state $s_2$.
\end{solution}
\end{exercise}\index{interface}

\index{polybox|)}

%-------- Section --------%
\section[Categorical properties of the composition product]{Categorical properties of the composition product%
  \sectionmark{Categorical properties of $\tri$}}
\sectionmark{Categorical properties of $\tri$}

\label{sec.comon.comp.prop}

We conclude this chapter by discussing several interesting properties of the composition product, many of which will come in handy in the following chapters.
We'll focus on how $\tri$ interacts with other constructions on $\poly$ that we introduced in previous chapters.

\subsection{Interaction with products and coproducts} \label{subsec.comon.comp.prop.prod}\index{coproduct!composition and}

It turns out that the composition product behaves well with products and coproducts on the left.

\index{distributive law!coproducts over composition}\index{distributive law!products over composition}

\begin{proposition}[Left distributivity of $\tri$ over $+$ and $\times$]\label{prop.left_dist_prod}
Given a polynomial $r$, the functor $(-\tri r)\colon\poly\to\poly$ that sends each $p\in\poly$ to $p\tri r$ commutes with coproducts and products (up to natural isomorphism).\index{isomorphism!natural}
That is, for any $p,q\in\poly$, we have the following natural isomorphisms:
\begin{equation}\label{eqn.comp_left_pres_plus}
    (p+q)\tri r\iso (p\tri r)+(q\tri r)
\end{equation}
and
\begin{equation}\label{eqn.comp_left_pres_times}
    pq\tri r\iso (p\tri r)(q\tri r).
\end{equation}
More generally, given a set $A$ and polynomials $(q_a)_{a\in A}$, we have the following natural isomorphisms:
\begin{equation}\label{eqn.comp_left_pres_sum}
    \left(\sum_{a\in A}q_a\right)\tri r\iso \sum_{a\in A} (q_a\tri r)
\end{equation}
and
\begin{equation}\label{eqn.comp_left_pres_prod}
    \left(\prod_{a\in A}q_a\right)\tri r\iso \prod_{a\in A} (q_a\tri r)
\end{equation}
\end{proposition}
\begin{proof}
Formally, this comes down to the fact that (co)products of functors $\smset\to\smset$ are computed pointwise (\cref{prop.presheaf_lim_ptwise}) and that (co)products in $\poly$ coincide with (co)products in $\smset^\smset$ (\cref{prop.poly_coprods,prop.poly_prods}).
One could instead give an explicit proof using \eqref{eqn.composite_formula}; this is done in \cref{exc.left_dist}.
In fact, we will see yet another proof of \eqref{eqn.comp_left_pres_prod} (and thus \eqref{eqn.comp_left_pres_times}) in \cref{exc.comp_left_pres_lim} \cref{exc.comp_left_pres_lim.prod}.
\end{proof}

\begin{exercise}\label{exc.left_dist}
Prove \cref{prop.left_dist_prod} using the explicit formula for $\tri$ given in \eqref{eqn.composite_formula} by manipulating sums and products.
\begin{solution}
To prove \cref{prop.left_dist_prod}, it suffices to verify \eqref{eqn.comp_left_pres_sum} and \eqref{eqn.comp_left_pres_prod}, as \eqref{eqn.comp_left_pres_plus} and \eqref{eqn.comp_left_pres_times} follow directly when $A\coloneqq\2$.

Given polynomials $(q_a)_{a\in A}$, recall that the position-set of the sum $\sum_{a\in A}q_a$ is $\sum_{a\in A}q_a(\1)$, while the direction-set at each position $(a,j)$ with $a\in A$ and $j\in q_a(\1)$ is $q_a[j]$.
So by \eqref{eqn.composite_formula}, we have that
\begin{align*}
    \left(\sum_{a\in A}q_a\right)\tri r &\iso \sum_{\substack{a\in A, \\ j\in q_a(\1)}}\;\prod_{b \in q_a[j]}\;\sum_{k \in r(\1)}\;\prod_{c \in r[k]}\yon \\
    &\iso \sum_{a\in A}\;\sum_{j\in q_a(\1)}\;\prod_{b \in q_a[j]}\;\sum_{k \in r(\1)}\;\prod_{c \in r[k]}\yon \\
    &\iso\sum_{a\in A}(q_a\tri r).
\end{align*}
We can also recall that the position-set of the product $\prod_{a\in A}q_a$ is $\prod_{a\in A}q_a(\1)$, while the direction-set at each position $\ol{j}\colon(a\in A)\to q_a(\1)$ is $\sum_{a\in A}q_a[\ol{j}(a)]$.
So by \eqref{eqn.composite_formula}, we have that
\begin{align*}
    \left(\prod_{a\in A}q_a\right)\tri r &\iso \sum_{\ol{j}\in\prod_{a\in A}q_a(\1)} \; \prod_{\substack{a\in A, \\ b\in q_a[\ol{j}(a)]}} \; \sum_{k \in r(\1)} \; \prod_{c \in r[j]} \yon \\
    &\iso \prod_{a\in A} \; \sum_{j\in q_a(\1)} \; \prod_{b\in q_a[j]} \; \sum_{k \in r(\1)} \; \prod_{c \in r[k]} \yon \tag*{\eqref{eqn.cat_completely_distributive}} \\
    &\iso\prod_{a\in A}(q_a\tri r).
\end{align*}
\end{solution}
\end{exercise}
\index{distributive law!products over composition}

\begin{example}[Picturing the left distributivity of $\tri$ over $\times$]\label{ex.picturing_dist}
We want an intuitive understanding of the left distributivity given by \eqref{eqn.comp_left_pres_times}.
Let $p\coloneqq\yon$, $q\coloneqq\yon+\1$, and $r\coloneqq\yon^\2+\1$, as shown here:
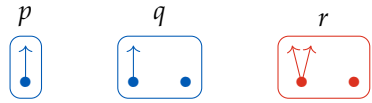
\[
\begin{tikzpicture}[rounded corners]
	\node (p) [draw, my-blue, "$p$" above] {
	\begin{tikzpicture}[trees, sibling distance=2.5mm]
    \node (p1) {$\bullet$}
      child {};
  \end{tikzpicture}
  };
	\node (q) [draw, my-blue, right=1 of p, "$q$" above] {
	\begin{tikzpicture}[trees, sibling distance=2.5mm]
    \node (q1) {$\bullet$}
      child {};
    \node[right=.5 of q1] (q2) {$\bullet$};
  \end{tikzpicture}
  };
	\node (r) [draw, my-red, right=1 of q, "$r$" above] {
	\begin{tikzpicture}[trees, sibling distance=2.5mm]
    \node (r1) {$\bullet$}
      child {}
      child {};
    \node[right=.5 of r1] (r2) {$\bullet$};
  \end{tikzpicture}
  };
\end{tikzpicture}
\]
Then $pq\cong\yon^\2+\yon$ can be drawn as follows, with each corolla comprised of a $p$-corolla and a $q$-corolla with their roots glued together:
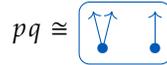
\[
\begin{tikzpicture}[rounded corners]
	\node (p) [draw, my-blue, "$pq\cong$" left] {
	\begin{tikzpicture}[trees, sibling distance=2.5mm]
        \node (q1) {$\bullet$}
          child {}
          child {};
        \node[right=.5 of q1] (q2) {$\bullet$}
          child {};
    \end{tikzpicture}
	};
\end{tikzpicture}
\]
We can therefore draw $pq\tri r$ by grafting $r$-corollas to leaves of $pq$ in every way, as follows:
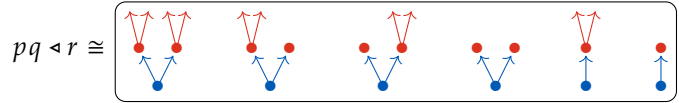
\[
\begin{tikzpicture}[rounded corners]
	\node (p) [draw, "$pq\tri r\cong$" left] {
	\begin{tikzpicture}[trees,
		level 1/.style={sibling distance=5mm},
	  level 2/.style={sibling distance=2.5mm}]
    \node[my-blue] (1) {$\bullet$}
      child[my-blue] {node[my-red] {$\bullet$}
      	child[my-red]
				child[my-red]
			}
      child[my-blue] {node[my-red] {$\bullet$}
      	child[my-red]
				child[my-red]
			};
    \node[my-blue, right=1.3 of 1] (2) {$\bullet$}
      child[my-blue] {node[my-red] {$\bullet$}
				child[my-red]
				child[my-red]
			}
      child[my-blue] {node[my-red] {$\bullet$}
			};
    \node[my-blue, right=1.3 of 2] (3) {$\bullet$}
      child[my-blue] {node[my-red] {$\bullet$}
			}
      child[my-blue] {node[my-red] {$\bullet$}
				child[my-red]
				child[my-red]
			};
    \node[my-blue, right=1.3 of 3] (4) {$\bullet$}
      child[my-blue] {node[my-red] {$\bullet$}
			}
      child[my-blue] {node[my-red] {$\bullet$}
			};
    \node[my-blue, right=1 of 4] (5) {$\bullet$}
      child[my-blue] {node[my-red] {$\bullet$}
      	child[my-red]
      	child[my-red]
			};
    \node[my-blue, right=.8 of 5] (6) {$\bullet$}
      child[my-blue] {node[my-red] {$\bullet$}
      };
  \end{tikzpicture}
	};
\end{tikzpicture}
\]
So each tree in $pq\tri r$ is obtained by grafting together the roots of a $p$-corolla and a $q$-corolla, then attaching $r$-corollas to each leaf.

Alternatively, we can compute $p\tri r$ and $q\tri r$ seperately, grafting $r$-corollas to leaves of $p$ in every way, then to leaves of $q$ in every way:
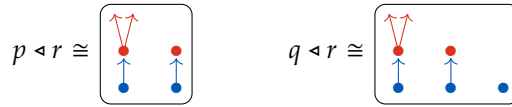
\[
\begin{tikzpicture}[rounded corners]
	\node (p) [draw, "$p\tri r\cong$" left] {
	\begin{tikzpicture}[trees,
		level 1/.style={sibling distance=5mm},
	  level 2/.style={sibling distance=2.5mm}]
    \node[my-blue] (p1) {$\bullet$}
      child[my-blue] {node[my-red] {$\bullet$}
      	child[my-red]
				child[my-red]
			};
    \node[my-blue, right=.5 of p1] (p2) {$\bullet$}
      child[my-blue] {node[my-red] {$\bullet$}
			};
  \end{tikzpicture}
  };
	\node (q) [draw, "$q\tri r\cong$" left] at (4,0) {
	\begin{tikzpicture}[trees,
		level 1/.style={sibling distance=5mm},
	  level 2/.style={sibling distance=2.5mm}]
    \node[my-blue] (q1) {$\bullet$}
      child[my-blue] {node[my-red] {$\bullet$}
      	child[my-red]
				child[my-red]
			};
    \node[my-blue, right=.5 of q1] (q2) {$\bullet$}
      child[my-blue] {node[my-red] {$\bullet$}
			};
    \node[my-blue, right=.5 of q2] (q3) {$\bullet$};
  \end{tikzpicture}
  };
\end{tikzpicture}
\]
Their product is then obtained by taking each tree from $p\tri r$ and pairing it with each tree from $q\tri r$ by gluing their roots together:
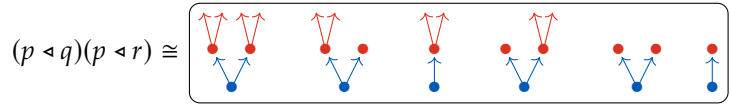
\[
\begin{tikzpicture}[rounded corners]
	\node (p) [draw, "$(p\tri q)(p\tri r)\cong$" left] {
	\begin{tikzpicture}[trees,
		level 1/.style={sibling distance=5mm},
	  level 2/.style={sibling distance=2.5mm}]
    \node[my-blue] (1) {$\bullet$}
      child[my-blue] {node[my-red] {$\bullet$}
      	child[my-red]
				child[my-red]
			}
      child[my-blue] {node[my-red] {$\bullet$}
      	child[my-red]
				child[my-red]
			};
    \node[my-blue, right=1.3 of 1] (2) {$\bullet$}
      child[my-blue] {node[my-red] {$\bullet$}
				child[my-red]
				child[my-red]
			}
      child[my-blue] {node[my-red] {$\bullet$}
			};
    \node[my-blue, right=1 of 2] (3) {$\bullet$}
      child[my-blue] {node[my-red] {$\bullet$}
      	child[my-red]
      	child[my-red]
			};
    \node[my-blue, right=1 of 3] (4) {$\bullet$}
      child[my-blue] {node[my-red] {$\bullet$}
			}
      child[my-blue] {node[my-red] {$\bullet$}
				child[my-red]
				child[my-red]
			};
    \node[my-blue, right=1.3 of 4] (5) {$\bullet$}
      child[my-blue] {node[my-red] {$\bullet$}
			}
      child[my-blue] {node[my-red] {$\bullet$}
			};
    \node[my-blue, right=.8 of 5] (6) {$\bullet$}
      child[my-blue] {node[my-red] {$\bullet$}
      };
	\end{tikzpicture}
	};
\end{tikzpicture}
\]
So each tree in $(p\tri r)(q\tri r)$ is obtained by grafting $r$-corollas to each leaf of a $p$-corolla and a $q$-corolla before gluing their roots together.

But it doesn't matter if we graft $r$-corollas onto leaves first, or if we glue the roots of $p$- and $q$-corollas together first--the processes are equivalent.
Hence the isomorphism $pq\tri r\iso(p\tri r)(q\tri r)$ holds.
\end{example}

\index{distributive law!coproducts over composition}

\begin{exercise}\label{exc.picturing_dist}
Follow \cref{ex.picturing_dist} with coproducts $(+)$ in place of products $(\times)$: use pictures to give an intuitive understanding of the left distributivity given by \eqref{eqn.comp_left_pres_plus}.
\begin{solution}
We want an intuitive understanding of the left distributivity of $\tri$ over $+$.
Let $p\coloneqq\yon^\2$, $q\coloneqq\yon+\1$, and $r\coloneqq\yon^\2+\1$, as shown here:
\[
\begin{tikzpicture}[rounded corners]
	\node (p) [draw, "$p$" above] {
	\begin{tikzpicture}[trees, sibling distance=2.5mm]
    \node (p1) {$\bullet$}
      child {}
      child {};
  \end{tikzpicture}
  };
	\node (q) [draw, my-blue, right=1 of p, "$q$" above] {
	\begin{tikzpicture}[trees, sibling distance=2.5mm]
    \node (q1) {$\bullet$}
      child {};
    \node[right=.5 of q1] (q2) {$\bullet$};
  \end{tikzpicture}
  };
	\node (r) [draw, my-red, right=1 of q, "$r$" above] {
	\begin{tikzpicture}[trees, sibling distance=2.5mm]
    \node (r1) {$\bullet$}
      child {}
      child {};
    \node[right=.5 of r1] (r2) {$\bullet$};
  \end{tikzpicture}
  };
\end{tikzpicture}
\]
Then $p+q\iso\yon^\2+\yon+\1$ can be drawn as follows, consisting of every $p$-corolla and every $q$-corolla:
\[
\begin{tikzpicture}[rounded corners]
	\node (p) [draw, "$p+q\iso$" left] {
	\begin{tikzpicture}[trees, sibling distance=2.5mm]
        \node[black] (q1) {$\bullet$}
          child[black] {}
          child[black] {};
        \node[my-blue, right=.5 of q1] (q2) {$\bullet$}
          child[my-blue] {};
        \node[my-blue, right=.5 of q2] (q3) {$\bullet$};
    \end{tikzpicture}
	};
\end{tikzpicture}
\]
We can therefore draw $(p+q)\tri r$ by grafting $r$-corollas to leaves of $p+q$ in every way, as follows:
\[
\begin{tikzpicture}[rounded corners]
	\node (p) [draw, "$(p+q)\tri r\cong$" left] {
	\begin{tikzpicture}[trees,
		level 1/.style={sibling distance=5mm},
	  level 2/.style={sibling distance=2.5mm}]
    \node[black] (1) {$\bullet$}
      child[black] {node[my-red] {$\bullet$}
      	child[my-red]
				child[my-red]
			}
      child[black] {node[my-red] {$\bullet$}
      	child[my-red]
				child[my-red]
			};
    \node[black, right=1.3 of 1] (2) {$\bullet$}
      child[black] {node[my-red] {$\bullet$}
				child[my-red]
				child[my-red]
			}
      child[black] {node[my-red] {$\bullet$}
			};
    \node[black, right=1.3 of 2] (3) {$\bullet$}
      child[black] {node[my-red] {$\bullet$}
			}
      child[black] {node[my-red] {$\bullet$}
				child[my-red]
				child[my-red]
			};
    \node[black, right=1.3 of 3] (4) {$\bullet$}
      child[black] {node[my-red] {$\bullet$}
			}
      child[black] {node[my-red] {$\bullet$}
			};
    \node[my-blue, right=1 of 4] (5) {$\bullet$}
      child[my-blue] {node[my-red] {$\bullet$}
      	child[my-red]
      	child[my-red]
			};
    \node[my-blue, right=.8 of 5] (6) {$\bullet$}
      child[my-blue] {node[my-red] {$\bullet$}
      };

    \node[my-blue, right=.8 of 6] (7) {$\bullet$};
  \end{tikzpicture}
	};
\end{tikzpicture}
\]
So each tree in $(p+q)\tri r$ is obtained by taking either a $p$-corolla or a $q$-corolla, then attaching an $r$-corolla to each leaf.

Alternatively, we can compute $p\tri r$ and $q\tri r$ separately, grafting $r$-corollas to leaves of $p$ in every way, then to leaves of $q$ in every way:
\[
\begin{tikzpicture}[rounded corners]
	\node (p) [draw, "$p\tri r\cong$" left] at (-2,0) {
	\begin{tikzpicture}[trees,
		level 1/.style={sibling distance=5mm},
	  level 2/.style={sibling distance=2.5mm}]
    \node[black] (1) {$\bullet$}
      child[black] {node[my-red] {$\bullet$}
      	child[my-red]
				child[my-red]
			}
      child[black] {node[my-red] {$\bullet$}
      	child[my-red]
				child[my-red]
			};
    \node[black, right=1.3 of 1] (2) {$\bullet$}
      child[black] {node[my-red] {$\bullet$}
				child[my-red]
				child[my-red]
			}
      child[black] {node[my-red] {$\bullet$}
			};
    \node[black, right=1.3 of 2] (3) {$\bullet$}
      child[black] {node[my-red] {$\bullet$}
			}
      child[black] {node[my-red] {$\bullet$}
				child[my-red]
				child[my-red]
			};
    \node[black, right=1.3 of 3] (4) {$\bullet$}
      child[black] {node[my-red] {$\bullet$}
			}
      child[black] {node[my-red] {$\bullet$}
			};
  \end{tikzpicture}
  };
	\node (q) [draw, "$q\tri r\cong$" left] at (4,0) {
	\begin{tikzpicture}[trees,
		level 1/.style={sibling distance=5mm},
	  level 2/.style={sibling distance=2.5mm}]
    \node[my-blue] (q1) {$\bullet$}
      child[my-blue] {node[my-red] {$\bullet$}
      	child[my-red]
				child[my-red]
			};
    \node[my-blue, right=.5 of q1] (q2) {$\bullet$}
      child[my-blue] {node[my-red] {$\bullet$}
			};
    \node[my-blue, right=.5 of q2] (q3) {$\bullet$};
  \end{tikzpicture}
  };
\end{tikzpicture}
\]
Their coproduct then consists of all the trees from $p\tri r$ and all the trees from $q\tri r$:
\[
\begin{tikzpicture}[rounded corners]
	\node (p) [draw, "$p\tri r+q\tri r\cong$" left] {
	\begin{tikzpicture}[trees,
		level 1/.style={sibling distance=5mm},
	  level 2/.style={sibling distance=2.5mm}]
    \node[black] (1) {$\bullet$}
      child[black] {node[my-red] {$\bullet$}
      	child[my-red]
				child[my-red]
			}
      child[black] {node[my-red] {$\bullet$}
      	child[my-red]
				child[my-red]
			};
    \node[black, right=1.3 of 1] (2) {$\bullet$}
      child[black] {node[my-red] {$\bullet$}
				child[my-red]
				child[my-red]
			}
      child[black] {node[my-red] {$\bullet$}
			};
    \node[black, right=1.3 of 2] (3) {$\bullet$}
      child[black] {node[my-red] {$\bullet$}
			}
      child[black] {node[my-red] {$\bullet$}
				child[my-red]
				child[my-red]
			};
    \node[black, right=1.3 of 3] (4) {$\bullet$}
      child[black] {node[my-red] {$\bullet$}
			}
      child[black] {node[my-red] {$\bullet$}
			};
    \node[my-blue, right=1 of 4] (5) {$\bullet$}
      child[my-blue] {node[my-red] {$\bullet$}
      	child[my-red]
      	child[my-red]
			};
    \node[my-blue, right=.8 of 5] (6) {$\bullet$}
      child[my-blue] {node[my-red] {$\bullet$}
      };

    \node[my-blue, right=.8 of 6] (7) {$\bullet$};
  \end{tikzpicture}
	};
\end{tikzpicture}
\]
So each tree in $p\tri r+q\tri r$ is either a $p$-corolla with an $r$-corolla attached to each leaf, or a $q$-corolla with an $r$-corolla attached to each leaf.

But it doesn't matter whether we graft $r$-corollas onto leaves first, or if we pool together corollas from $p$ and $q$ first--the processes are equivalent.
Hence the isomorphism $(p+q)\tri r\iso p\tri r+q\tri r$ holds.
\end{solution}
\end{exercise}

\begin{exercise}
Show that for any set $A$ and polynomials $p,q$, we have an isomorphism $A(p\tri q)\iso (Ap)\tri q$.
\begin{solution}
Given a set $A$ and polynomials $p,q$, the left distributivity of $\tri$ over products from \eqref{eqn.comp_left_pres_prod} implies that $(Ap)\tri q\iso(A\tri q)(p\tri q)$, while \cref{exc.composing_with_constants} \cref{exc.composing_with_constants.appl} implies that $A\tri q\iso A$.
So $(Ap)\tri q\iso A(p\tri q)$.
\end{solution}
\end{exercise}

In \cref{subsec.comon.comp.prop.lim_left}, we will see how to generalize the left distributivity of $\tri$ over products to arbitrary limits.
But first, we observe that right distributivity does not hold.

\begin{exercise} \label{exc.right_not_dist_prod}
Show that the distributivities of \cref{prop.left_dist_prod} do not hold on the other side:
\begin{enumerate}
	\item \label{exc.right_not_dist_prod.prod} Find polynomials $p,q,r$ such that $p\tri (qr)\not\iso(p\tri q)(p\tri r)$.
	\item Find polynomials $p,q,r$ such that $p\tri (q+r)\not\iso(p\tri q)+(p\tri r)$.
\qedhere
\end{enumerate}
\begin{solution}
\begin{enumerate}
    \item Let $p \coloneqq \yon + \1, q \coloneqq \1,$ and $r \coloneqq \0$.
    Then $p \tri (qr) \iso (\yon + \1) \tri \0 \iso \1$, while $(p \tri q)(p \tri r) \iso ((\yon + \1) \tri \1)((\yon + \1) \tri \0) \iso \2 \times \1 \iso \2$.
    \item Again let $p \coloneqq \yon + \1, q \coloneqq \1,$ and $r \coloneqq \0$.
    Then $p \tri (q+r) \iso (\yon + \1) \tri \1 \iso \2$, while $(p \tri q)+(p \tri r) \iso ((\yon + \1) \tri \1)+((\yon + \1) \tri \0) \iso \2 + \1 \iso \3$.
\end{enumerate}
\end{solution}
\end{exercise}

Nevertheless, there is something to be said about the relationship between $p\tri q, p\tri r,$ and $p\tri(qr)$.
We'll see this in action after we discuss how $\tri$ preserves limits on the left.

\subsection{Interaction with limits on the left} \label{subsec.comon.comp.prop.lim_left}

\index{composition product!limits and}

We saw in \cref{thm.poly_limits} that $\poly$ has all limits, and we saw in \cref{exc.refl_limits} that these limits coincide with limits in $\smset^\smset$.
Hence the argument in the proof of \cref{prop.left_dist_prod} by appealing to \cref{prop.presheaf_lim_ptwise} can be generalized to arbitrary limits.
It follows that $\tri$ preserves all limits on the left.
But we will present a proof of this fact from an alternative perspective: by appealing to the left coclosure of $\tri$.

\index{composition product!coclosure|(}
\index{left Kan extension|see{composition product, coclosure}}

\begin{proposition}[Meyers] \label{prop.comp_left_coclosed}
The composition product is left coclosed.
That is, there exists a left coclosure operation, which we denote $\lchom{-}{-}\colon\poly\op\times\poly\to\poly$,
such that there is a natural isomorphism\index{isomorphism!natural}
\begin{equation} \label{eqn.lchom_adj_iso}
    \poly(p, r\tri q)\iso\poly\left(\lchom{q}{p}, r\right).
\end{equation}
In particular, the left coclosure operation sends $q,p\in\poly$ to
\begin{equation} \label{eqn.lchom_def}
    \lchom{q}{p}\coloneqq\sum_{i\in p(\1)}\yon^{q(p[i])}.
\end{equation}
\end{proposition}
\begin{proof}
We present an argument using polyboxes; we leave it to the reader to write this proof in more standard mathematical notation in \cref{exc.comp_left_coclosed_calc}.

\index{polybox}

As in \cref{ex.map_to_comp}, a lens $\varphi\colon p\to r\tri q$ can be written as follows:
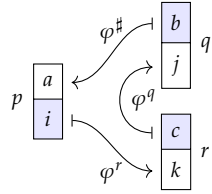
\[
\begin{tikzpicture}[polybox, mapstos]
	\node[poly, dom, "$p$" left] (p) {$a$\at$i$};
	\node[poly, cod, right=1.5cm of p.south, yshift=-1ex, "$r$" right] (r) {$c$\at$k$};
	\node[poly, cod, above=of r, "$q$" right] (q) {$b$\at$j$};
  	\draw (p_pos) to[first] node[below,yshift=-1mm] {$\varphi^r$} (r_pos);
  	\draw (r_dir) to[climb] node[right] {$\varphi^q$} (q_pos);
  	\draw (q_dir) to[last] node[above] {$\varphi^\sharp$} (p_dir);
\end{tikzpicture}
\]
But this is equivalent to the following gadget (to visualize this equivalence, imagine leaving the positions box for $p$, the arrow $\varphi^r$, and the polyboxes for $r$ untouched, while dragging the polyboxes for $q$ leftward to the directions box for $p$, merging all the data from $q$ and the arrows $\varphi^q$ and $\varphi^\sharp$ into a single on-directions arrow and directions box):
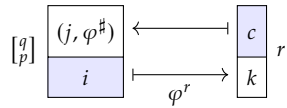
\[
\begin{tikzpicture}[polybox, mapstos]
    \node[poly, dom, "$\lchom{q}{p}$" left] (l) {$(j,\varphi^\sharp)$\at$i\vphantom{k}$};
    \node[poly, cod, "$r$" right, right=of l] (r) {$c\vphantom{(j,\varphi^\sharp)}$\at$k$};
    \draw (l_pos) -- node[below] {$\varphi^r$} (r_pos);
    \draw (r_dir) -- node[above] {} (l_dir);
\end{tikzpicture}
\]
Here the on-directions function encodes the behaviors of both $\varphi^q$ and $\varphi^\sharp$ by sending each $r[k]$-direction $c$ to both a $q$-position $j$, as $\varphi^q$ does, and a function sending $q[j]$-directions to $p[i]$-directions, as $\varphi^\sharp$ does.
So the polyboxes on the left represent a polynomial whose positions are the same as those of $p$, but whose directions at $i\in p(\1)$ are pairs $(j,\varphi^\sharp)$ consisting of a $q$-position $j$ and a function $\varphi^\sharp\colon q[j]\to p[i]$.
Such pairs are precisely the elements of $q(p[i])$, so the polynomial represented by the polyboxes on the left is indeed the one defined as $\lchom{q}{p}$ in \eqref{eqn.lchom_def}.
It follows that there is a natural isomorphism between lenses $p\to r\tri q$ and lenses $\lchom{q}{p}\to r$.
\end{proof}

\begin{exercise} \label{exc.comp_left_coclosed_calc}
Translate the polyboxes proof of \cref{prop.comp_left_coclosed} into standard mathematical notation, i.e.\ the $\sum$ and $\prod$ notation we have been using up till now.
\begin{solution}
We prove \cref{prop.comp_left_coclosed} by observing that
\begin{align*}
    \poly(p,r\tri q)&\iso\prod_{i\in p(\1)}r(q(p[i])) \tag*{\eqref{eqn.main_formula}} \\
    &\iso\poly\left(\sum_{i\in p(\1)}\yon^{q(p[i])}, r\right) \tag*{\eqref{eqn.main_formula}} \\
    &\iso\poly\left(\lchom{q}{p}, r\right) \tag*{\eqref{eqn.lchom_def}}.
\end{align*}
\end{solution}
\end{exercise}

\begin{remark}
The proof you came up with in \cref{exc.comp_left_coclosed_calc} may be more obviously rigorous and concise than the one we presented in the main text.
But polyboxes help us see right on paper exactly what is going on in this adjunction: how data on the codomain-side of a lens $p\to r\tri q$ can be simply repackaged and transferred to the domain-side of a new lens $\lchom{q}{p}\to r$.
\end{remark}

\begin{exercise} \label{exc.lchom_func}
In stating \cref{prop.comp_left_coclosed}, we implicitly assumed that $\lchom{-}{-}$ is a functor $\poly\op\times\poly\to\poly$.
Here we show that this is indeed the case.
\begin{enumerate}
    \item Given a polynomial $q$ and a lens $\varphi\colon p\to p'$, to what lens $\lchom{q}{p}\to\lchom{q}{p'}$ should the covariant functor $\lchom{q}{-}$ send $\varphi$?
    Prove that your construction is functorial.

    \item Given a polynomial $p$ and a lens $\psi\colon q'\to q$, to what lens $\lchom{q}{p}\to\lchom{q'}{p}$ should the contravariant functor $\lchom{-}{p}$ send $\psi$?
    Prove that your construction is functorial.
    \qedhere
\end{enumerate}
\begin{solution}
\begin{enumerate}
    \item Given a polynomial $q$ and a lens $\varphi\colon p\to p'$, the functor $\lchom{q}{-}$ should send $\varphi$ to a lens $\lchom{q}{p}\to\lchom{q}{p'}$; by \eqref{eqn.lchom_def}, this is a lens
    \[
        \lchom{q}{\varphi}\colon\sum_{i\in p(\1)}\yon^{q(p[i])}\to\sum_{i'\in p'(\1)}\yon^{q(p'[i'])}.
    \]
    We give $\lchom{q}{\varphi}$ the same on-positions function $p(\1)\to p'(\1)$ that $\varphi$ has.
    Then viewing $q$ as a functor, we define the on-directions function $\lchom{q}{\varphi}^\sharp_i\colon q(p'[\varphi_\1(i)])\to q(p[i])$ for each $i\in p(\1)$ to be the function obtained by applying $q$ to the corresponding on-directions function $\varphi^\sharp_i\colon p'[\varphi_\1(i)]\to p[i]$.
    Functoriality follows trivially on positions and by the functoriality of $q$ itself on directions.

    \item Given a polynomial $p$ and a lens $\psi\colon q'\to q$, the functor $\lchom{-}{p}$ should send $\psi$ to a lens $\lchom{q}{p}\to\lchom{q'}{p}$; by \eqref{eqn.lchom_def}, this is a lens
    \[
        \lchom{\psi}{p}\colon\sum_{i\in p(\1)}\yon^{q(p[i])}\to\sum_{i\in p(\1)}\yon^{q'(p[i])}.
    \]
    We let the on-positions function of $\lchom{\psi}{p}$ be the identity on $p(\1)$.
    Then viewing $\psi$ as a natural transformation, we define the on-directions function $\lchom{\psi}{p}^\sharp_i\colon q'(p[i])\to q(p[i])$ for each $i\in p(\1)$ to be the $p[i]$-component of $\psi$.
    Functoriality follows trivially on positions and because natural transformations compose componentwise on directions.
\end{enumerate}
\end{solution}
\end{exercise}

\begin{exercise}
In personal communication, Todd Trimble noted (the in-retrospect-obvious fact) that the left coclosure can be thought of as a left Kan extension
\[
\begin{tikzcd}
	\smset\ar[r, "p"]\ar[d, "q"']\ar[dr, phantom, very near start, "\Downarrow"]&
	\smset\\
	\smset\ar[ur, "\lchom{q}{p}"']&~
\end{tikzcd}
\]
Verify this.
\end{exercise}
\begin{solution}
In order for $\lchom{q}{p}$ to be a left Kan extension of $p$ along $q$, we need to first provide a natural transformation $p\to \lchom{q}{p}\tri q$, and second show that it is universal. But this is exactly the content of \cref{prop.comp_left_coclosed}: the adjunction \eqref{eqn.lchom_adj_iso} provides a unit $p\to \lchom{q}{p}\tri q$ that is universal in the appropriate way.
\end{solution}

\begin{exercise}
Let $A$ and $B$ be sets, and let $p$ and $q$ be polynomials.
\begin{enumerate}
    \item Prove that the following natural isomorphism holds:\index{isomorphism!natural}
    \begin{equation}\label{eqn.monomials_and_comp}
    	\poly(A\yon^B,p)\iso\smset(A,p(B)).
    \end{equation}

    \item Prove that the following natural isomorphism holds:
    \begin{equation}\label{eqn.flip_reps_lins}
        \poly\left(A\yon\tri p\tri \yon^B, q\right)\iso\poly\left(p,\yon^A\tri q\tri By\right).
    \end{equation}
    (Hint: Break the isomorphism down into two parts.
    You may find \eqref{eqn.two_var_adj} helpful.)
    \qedhere
\end{enumerate}
\begin{solution}
We are given $A,B\in\smset$ and $p,q\in\poly$
\begin{enumerate}
    \item By \eqref{eqn.lchom_def},
    \[
        \lchom{B}{A}=\sum_{i\in A}\yon^{B}\iso A\yon^B,
    \]
    so by \eqref{eqn.lchom_adj_iso},
    \[
        \poly(A\yon^B,p)\iso\poly(A,p\tri B).
    \]
    But $A$ and $p\tri B\iso p(B)$ are both constants, so a lens $A\to p\tri B$ is just a function $A\to p(B)$ on positions.
    Hence \eqref{eqn.monomials_and_comp} follows.

    \item We prove \eqref{eqn.flip_reps_lins} in two parts: that
    \begin{equation} \label{eqn.flipping1}
        \poly\left(A\yon\tri p, q\right)\iso\poly\left(p,\yon^A\tri q\right)
    \end{equation}
    and that
    \begin{equation} \label{eqn.flipping2}
        \poly\left(p \tri \yon^B, q\right)\cong\poly\left(p, q\tri B\yon\right).
    \end{equation}
    We have that $A\yon \tri p \iso Ap$ and $\yon^A \tri q \iso q^A$, so \eqref{eqn.flipping1} follows from \eqref{eqn.two_var_adj}.
    Meanwhile, for \eqref{eqn.flipping2}, we have by \eqref{eqn.lchom_def} that
    \[
        \lchom{B\yon}{p}=\sum_{i\in p(\1)}\yon^{Bp[i]}\iso p\tri\yon^B,
    \]
    so \eqref{eqn.flip_reps_lins} follows from \eqref{eqn.lchom_adj_iso}.
    Then combining \eqref{eqn.flipping1} and \eqref{eqn.flipping2} yields
    \[
        \poly\left(A\yon\tri p\tri\yon^B,q\right)\iso\poly\left(p\tri\yon^B,\yon^A\tri q\right)\iso\poly\left(p,\yon^A\tri q\tri B\yon\right).
    \]
\end{enumerate}
\end{solution}
\end{exercise}

% \[
% \begin{tikzpicture}
% 	\node (p1) {
%   \begin{tikzpicture}[polybox,tos]
%   	\node[poly, "$p$" left] (p) {};
%   	\node[poly, linear, below=of p, "$A\yon$" left] (Ay) {};
%   	\node[poly, pure, above=of p, "$\yon^B$" left] (yB) {};
%   	\node[poly, right=2 of p, "$q$" right] (q) {};
%   	\node at ($(p.east)!.5!(q.west)$) {$\leftrightarrows$};
%   \end{tikzpicture}
%   };
%  \node[right=4 of p1] (p2) {
%  \begin{tikzpicture}[polybox,tos]
%   	\node[poly, "$p$" left] (p) {};
%   	\node[poly, right=2 of p, "$q$" right] (q) {};
%   	\node[poly, linear, above=of q, "$B\yon$" right] (By) {};
%   	\node[poly, pure, below=of q, "$\yon^A$" right] (yA) {};
% 		\draw (p_pos) to[first] (yA_pos);
% 		\draw (yA_dir) to[climb] (q_pos);
% 		\draw (q_dir) to[climb] (By_pos);
% 		\draw (By_dir) to[last] (p_dir);
%  \end{tikzpicture}
%  };
%  \node[align=center] at ($(p1)!.5!(p2)$) {is the\\same as};
% \end{tikzpicture}
% \]
% Do you see how polyboxes with a black (one-element) part can flip upside-down to go to the other side?

\begin{example}[Dynamical systems as coalgebras] \label{ex.coalgebras}\index{coalgebra!dynamical systems as}
Taking $A=B=S\in\smset$ in \eqref{eqn.monomials_and_comp}, we find that there is a natural isomorphism between dynamical systems $S\yon^S\to p$ and functions $S\to p(S)$.\index{isomorphism!natural}
Such a function is known as a \emph{coalgebra for the functor $p$} or a \emph{$p$-coalgebra}.%
\tablefootnote{There are two versions of coalgebras we are interested in (and more that we are not) with distinct definitions: a \emph{coalgebra for a functor}, which is the version used here, and a \emph{coalgebra for a comonad}, which is a coalgebra for a functor with extra conditions that we will introduce later in \cref{sec.comon.sharp.cof.from_state}.
The version we are using will usually be clear from context---here, for example, we do not expect $p$ to be a comonad---but we will try to be explicit with our terminology whenever the interpretation may be ambiguous.}

Coalgebras as models of dynamical systems have been studied extensively in the context of computer science, most notably by Jacobs in \cite{jacobs2017introduction}.
Indeed, much of what we developed in \label{sec.poly.dyn_sys.moore,sec.poly.dyn_sys.depend_sys} stems from the theory of coalgebras.
The coalgebraic perspective has the benefit of staying in the familiar category of sets; moreover, it can be generalized to functors $\smset\to\smset$ that are not polynomial, although many of the interesting examples are.

On the other hand, we have already seen that viewing dynamical systems as lenses $S\yon^S\to p$ rather than as functions $S\to p(S)$ has the benefit of isolating the internal state system to the domain and the external interface to the codomain, aiding both intuition and functionality.
Plus, our adjunction lets us switch to the coalgebraic perspective whenever we see fit: $\poly$ lets us talk about both.
\end{example}\index{interface}

\index{composition product!coclosure|)}
\index{composition product!limits and}

\begin{proposition}[Left preservation of limits] \label{prop.left_pres_lim}
The operation $\tri$ preserves limits on the left (up to natural isomorphism).
That is, if $\cat{J}$ is a category, $p_-\colon\cat{J}\to\poly$ is a functor, and $q\in\poly$ is a polynomial, then there is a natural isomorphism
\begin{equation} \label{eqn.comp_left_pres_lim}
    \left(\lim_{j\in\cat{J}}p_j\right)\tri q\iso\lim_{j\in\cat{J}}(p_j\tri q).
\end{equation}
\end{proposition}
\begin{proof}
By \cref{prop.comp_left_coclosed}, the functor $(-\tri q)\colon\poly\to\poly$ is the right adjoint of the functor $\lchom{q}{-}\colon\poly\to\poly$, and right adjoints preserve limits.
\end{proof}

\begin{exercise} \label{exc.comp_left_pres_lim}
\begin{enumerate}
    \item Complete \cref{exc.composing_with_constants} \cref{exc.composing_with_constants.appl} using \eqref{eqn.comp_left_pres_lim} and \eqref{eqn.comp_left_pres_sum}.
    \item \label{exc.comp_left_pres_lim.prod} Deduce \eqref{eqn.comp_left_pres_prod} using \eqref{eqn.comp_left_pres_lim}.\qedhere
\end{enumerate}
\begin{solution}
\begin{enumerate}
    \item We wish to solve \cref{exc.composing_with_constants} \cref{exc.composing_with_constants.appl} using \eqref{eqn.comp_left_pres_lim} and \eqref{eqn.comp_left_pres_sum}.
    If we set $\cat{J}$ in \eqref{eqn.comp_left_pres_lim} to be the empty category, then the limit of the functor from $\cat{J}$ is just the terminal object.
    It follows that $\1\tri p\iso\1$.
    In other words, since $\tri$ preserves limits on the left, and since terminal objects are limits, $\tri$ preserves terminal objects on the left.

    Then a set $X$ can be written as a sum $\sum_{x\in X}\1$, so by \eqref{eqn.comp_left_pres_sum},
    \[
        X\tri p\iso\left(\sum_{x\in X}\1\right)\tri p\iso\sum_{x\in X}(\1\tri p)\iso\sum_{x\in X}\1\iso X.
    \]
\index{category!discrete}
    \item If we set $\cat{J}$ in \eqref{eqn.comp_left_pres_lim} to be the discrete category on the set $A$, then the limit of a functor from $\cat{J}$ is just an $A$-fold product, so \eqref{eqn.comp_left_pres_prod} follows.
    In other words, since $\tri$ preserves limits on the left, and since products are limits, $\tri$ preserves products on the left.
\end{enumerate}
\end{solution}
\end{exercise}

\subsection{Interaction with limits on the right} \label{subsec.comon.comp.prop.lim_right}

\index{composition product!connected limits and}

So $\tri$ preserves limits on the left.
How about limits on the right?
We saw in \cref{exc.right_not_dist_prod} that $\tri$ does not even preserve products on the right, so it certainly does not preserve all limits.
But it turns out that there is a special class of limits that $\tri$ does preserve on the right.

\index{limit!connected}\index{connected limit|see{limit, connected}}

\begin{definition}[Connected limit]\index{indexing category!for (co)limit}
A \emph{connected limit} is one whose indexing category $\cat{J}$ is nonempty and connected. That is, $\cat{J}$ has at least one object, and any two objects are connected by a finite zigzag of arrows.
\end{definition}

\index{category!connected}\index{equalizer}\index{pullback}\index{directed limit}\index{pullback!as connected limit}\index{equalizer!as connected limit}

\begin{example}
The following categories are connected:
\[
\fbox{$\bullet$}
\qquad
\fbox{$\bullet\tto\bullet$}
\qquad
\fbox{$\bullet\to\bullet\from\bullet$}
\qquad
\fbox{$\bullet\from\bullet\from\bullet\from\cdots$}
\]
In particular, equalizers, pullbacks, and directed limits are examples of connected limits.

The following categories are \emph{not} connected:
\[
\fbox{$\phantom{\bullet}$}
\qquad
\fbox{$\bullet\quad\bullet$}
\qquad
\fbox{$\bullet\quad\bullet\to \bullet$}
\]
In particular, terminal objects and products are \emph{not} examples of connected limits.
\end{example}

Connected limits are intimately related to slice categories, which we defined back in \cref{def.slice}.
For example, products in a slice category $\cat{C}/c$ are just pullbacks in $\cat{C}$, allowing us to view a non-connected limit as a connected one.
By relating $\poly$ to its slice categories via an adjunction, we'll be able to show that $\tri$ preserves connected limits.
(An alternative proof of this fact can be found in \cite[Proposition 1.16]{kock2012polynomial}.)
%The claim for the right side comes down to the fact that polynomials are sums of representables; representable functors commute with all limits and sums commute with connected limits in $\smset$.

\index{category!slice}

Recall that objects in a slice category $\cat{C}/c$ are just morphisms with codomain $c$.
For ease of notation, we'll often suppress the actual morphism and just write down the name of its domain when there is a canonical choice for the morphism, or when it is clear from context.
So for example, on the left hand side of \eqref{eqn.rchom_adj_iso} below, $p$ represents the lens $f\colon p\to q\tri 1$ and $q\tri r$ represents the lens $q\:\tri\:!\colon q\tri r\to q\tri\1$, both objects in the slice category $\poly/q\tri\1$.

\begin{proposition} \label{prop.comp_right_slice_coclosed}
Given polynomials $p,q,r\in\poly$ and a function $f\colon p\to q\tri\1$, there is a natural isomorphism\index{isomorphism!natural}
\begin{equation} \label{eqn.rchom_adj_iso}
    \poly/(q\tri\1)\left(p, q\tri r\right)\iso\poly\left(p\rchom{f}q, r\right),
\end{equation}
where
\begin{equation} \label{eqn.rchom_def}
    p\rchom{f}q\coloneqq\sum_{i\in p(\1)}q[f(i)]\yon^{p[i]}.
\end{equation}
\end{proposition}
\begin{proof}
Again, we present an argument using polyboxes; we leave it to the reader to write this proof in more standard mathematical notation in \cref{exc.comp_right_slice_coclosed_calc}.

\index{polybox}

Note that to consider $q\tri r$ as an object in $\poly/(q\tri\1)$, we are implicitly using the lens $q\:\tri\:!\colon q\tri r\to q\tri 1$. By definition, morphisms from $f$ to $q\:\tri\:!$ in $\poly/(q\tri\1)$ are lenses $\varphi\colon p\to q\tri r$ for which $\varphi\then(q\:\tri\:!)=f$.
We can write this equation using polyboxes:
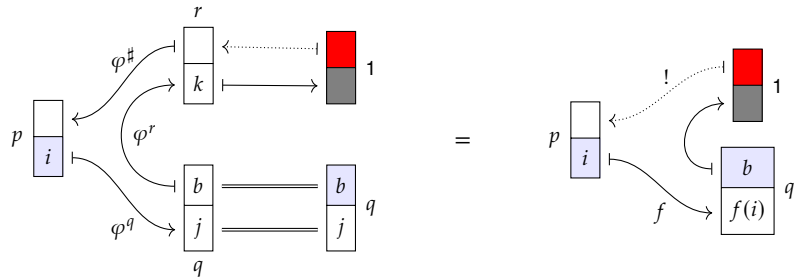
\[
\begin{tikzpicture}
	\node (a) {
  \begin{tikzpicture}[polybox, mapstos]
	\node[poly, dom, "$p$" left] (p) {\at$i$};
	\node[poly, right=1.8 of p.south, yshift=-2.5ex, "$q$" below] (q) {$b$\at$j$};
	\node[poly, above=.8 of q, "$r$" above] (r) {\at$k$};
	\node[poly, cod, right=of q, "$q$" right] (q') {$b$\at$j$};
	\node[poly, terminal, above=.8 of q', "$\1$" right] (1) {};
	\draw[double, -] (q_pos) -- node[below] {} (q'_pos);
	\draw[double, -] (q'_dir) -- node[above] {} (q_dir);
	\draw (r_pos) -- node[below] {} (1_pos);
	\draw[densely dotted] (1_dir) -- node[above] {} (r_dir);
	\draw (p_pos) to[first] node[below,yshift=-2mm] {$\varphi^q$} (q_pos);
	\draw (q_dir) to[climb] node[right] {$\varphi^r$} (r_pos);
	\draw (r_dir) to[last] node[above] {$\varphi^\sharp$} (p_dir);
  \end{tikzpicture}
	};
	\node[right=1.8 of a] (b) {
  \begin{tikzpicture}[polybox, mapstos]
  	\node[poly, dom, "$p$" left] (p) {\at$i$};
  	\node[poly, cod, right=1.8 of p.south, yshift=-1ex, "$q$" right] (q) {$b$\at$f(i)$};
  	\node[poly, terminal, above=of q, "$\1$" right] (1) {};
  	\draw (p_pos) to[first] node[below,yshift=-1mm] {$f$} (q_pos);
  	\draw (q_dir) to[climb] node[right] {} (1_pos);
  	\draw[densely dotted] (1_dir) to[last] node[above] {$!$} (p_dir);
  \end{tikzpicture}
	};
	\node at ($(a.east)!.5!(b.west)$) {=};
\end{tikzpicture}
\]
We can read off the picture that a lens $\varphi\colon p\to q\tri r$ is a morphism from $f$ to $q\:\tri\:!$ in $\poly/q\tri\1$ if and only if $\varphi^q=f_\1$.
% , so specifying such a morphism amounts to specifying a position $k$ of $r$ and an on-directions function $r[k]\to p[i]$ for each $i\in p(\1)$ and $b\in q[f_\1(i)]$.
So morphisms from $f$ to $q\:\tri\:!$ are equivalent to gadgets
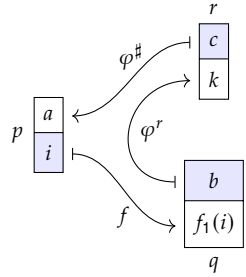
\[
\begin{tikzpicture}[polybox, mapstos]
    \node[poly, dom, "$p$" left] (p) {$a$\at$i$};
    \node[poly, cod, right=1.8 of p.south, yshift=-2.5ex, "$q$" below] (q) {$b$\at$f_\1(i)$};
    \node[poly, cod, above=.8 of q, "$r$" above] (r) {$c$\at$k$};
    \draw (p_pos) to[first] node[below,yshift=-1mm] {$f$} (q_pos);
    \draw (q_dir) to[climb] node[right] {$\varphi^r$} (r_pos);
    \draw (r_dir) to[last] node[above] {$\varphi^\sharp$} (p_dir);
\end{tikzpicture}
\]
with $f_\1$ fixed.
But this, in turn, is equivalent to the following gadget (to visualize this equivalence, imagine leaving the polyboxes for $r$, the arrow $\varphi^\sharp$, and the directions box for $p$ untouched, while dragging the polyboxes for $q$ leftward to the positions box for $p$, merging the data from $q$ and the predetermined arrow $f$ into a single positions box and adapting the arrow $\varphi^r$ into an on-positions arrow):
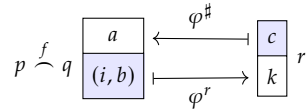
\[
\begin{tikzpicture}[polybox, mapstos]
    \node[poly, dom, "$p\rchom{f}q$" left] (h) {$a$\at$(i,b)$};
    \node[poly, cod, "$r$" right, right=of h] (r) {$c$\at$k$};
    \draw (h_pos) -- node[below] {$\varphi^r$} (r_pos);
    \draw (r_dir) -- node[above] {$\varphi^\sharp$} (h_dir);
\end{tikzpicture}
\]
Here the user can provide both the $p$-position $i$ and the $q[f(i)]$-direction $b$ right from the start, as they know what to expect from $f$ ahead of time.
Then the on-positions function encodes the behavior of $\varphi^r$.
So the polyboxes on the left represent a polynomial whose positions are pairs $(i,b)$ with $i\in p(\1)$ and $b\in q[f(i)]$, and whose directions at $(i,b)$ are the directions of the original polynomial $p$ at $i$.
This is precisely the polynomial we defined in \eqref{eqn.rchom_def}, so the isomorphism holds.
\end{proof}

\begin{remark}
As a lens $p\to q\tri\1$ can be identified with its on-positions function $p(\1)\to q(\1)$, we will use the notation $p\rchom{f}q$ interchangeably for lenses $f\colon p\to q\tri\1$ and functions $f\colon p(\1)\to q(\1)$.
\end{remark}

\begin{exercise} \label{exc.comp_right_slice_coclosed_calc}
\begin{enumerate}
    \item Translate the polyboxes proof of \cref{prop.comp_right_slice_coclosed} into standard mathematical notation.
    \item Prove that the following natural isomorphism holds:
    \begin{equation}
        \poly(p,q\tri r)\iso\sum_{f\colon p(\1)\to q(\1)}\poly\left(p\rchom{f}q,r\right).
    \end{equation}
    Thus the functor $(q\tri-)\colon\poly\to\poly$ is said to have a \emph{left multiadjoint}.
\end{enumerate}
\begin{solution}
\begin{enumerate}
	\item Note that $\poly(p,q\tri \1)\cong\smset(p(1),q(1))$. A lens $p\to q\tri r$ is an element of
    \[
    \prod_{i\in p(\1)}\sum_{j\in q(\1)}\prod_{b\in q[j]}\sum_{k\in r(\1)}\prod_{c\in r[k]}p[i]\cong\sum_{f\colon p(\1)\to q(\1)}\prod{i\in p(1)}\prod_{b\in q[j]}\sum_{k\in r(\1)}\prod_{c\in r[k]}p[i]
    \]
    Fixing the function $f\colon p(\1)\to q(\1)$ as implicit in $p$, we get
\begin{align*}
	\poly/(q\tri\1)\left(p, q\tri r\right)&\cong
	\prod{i\in p(1)}\prod_{b\in q[j]}\sum_{k\in r(\1)}\prod_{c\in r[k]}p[i]\\&\cong
	\poly\left(p\rchom{f}q, r\right),
\end{align*}
	\item This follows from the first isomorphism above.
 \end{enumerate}
\end{solution}
\end{exercise}

\begin{exercise} \label{exc.rchom_func}
In stating \cref{prop.comp_right_slice_coclosed}, we implicitly assumed that $p\rchom{f}q\in\poly$ is functorial in each variable: covariantly on the left and contravariantly on the right.
Here we show that this is indeed the case.
\begin{enumerate}
    \item Given lenses $f\colon p\to q\tri\1$ and $g\colon p'\to p$, construct the lens $p'\rchom{g\then f}q\to p\rchom{f}q$ to which the covariant functor $-\rchom{f}q$ should send $g$.
    Prove that your construction is functorial.

    \item Given lenses $f\colon p\to q\tri\1$ and $h\colon q\to q'$, to construct the lens $p\rchom{f\then(h\tri\1)}q'\to p\rchom{f}q$ to which the contravariant functor $p\rchom{f}-$ should send $h$.
    Prove that your construction is functorial.
    \qedhere
\end{enumerate}
\begin{solution}
\begin{enumerate}
    \item Given lenses $f\colon p\to q\tri\1$ and $g\colon p'\to p$, the functor $-\rchom{f}q$ should send $g$ to a lens $p'\rchom{g\then f}q\to p\rchom{f}q$; by \eqref{eqn.rchom_def}, this is a lens
    \[
        g\rchom{f}q\colon\sum_{i'\in p'(\1)}q[f_\1(g_\1(i'))]\yon^{p'[i']}\to\sum_{i\in p(\1)}q[f_\1(i)]\yon^{p[i]}.
    \]
    We give $g\rchom{f}q$ an on-positions function that is the on-positions function $p'(\1)\to p(\1)$ of $g$ on the first coordinate $i'\in p'(\1)$ and the identity on $q[f_\1(g_\1(i'))]$ on the second.
    Then we let the on-directions function at every position with first coordinate $i'\in p'(\1)$ be the on-directions function $g^\sharp_{i'}\colon p[g_\1(i')]\to p'[i']$.
    Functoriality follows trivially on both positions and directions.

    \item Given a polynomial $p$ and lenses $f\colon p\to q\tri\1$ and $h\colon q\to q'$, the functor $p\rchom{f}-$ should send $h$ to a lens $p\rchom{f\then(h\tri\1)}q'\to p\rchom{f}q$; by \eqref{eqn.rchom_def}, this is a lens
    \[
        p\rchom{f}h\colon\sum_{i\in p(\1)}q'[h_\1(f_\1(i))]\yon^{p[i]}\to\sum_{i\in p(\1)}q[f_\1(i)]\yon^{p[i]}.
    \]
    We let the on-positions function of $p\rchom{f}h$ be the identity on the first coordinate $i\in p(\1)$ and the on-directions function $h^\sharp_{f_\1(i)}\colon q'[h_\1(f_\1(i))]\to q[f_\1(i)]$ on the second.
    Then we let the on-directions function at every position with first coordinate $i\in p(\1)$ be the identity on $p[i]$.
    Functoriality follows trivially on both positions and directions.
\end{enumerate}
\end{solution}
\end{exercise}

\index{composition product!connected limits and}

\begin{theorem}[Preservation of connected limits]\label{thm.connected_limits}
The operation $\tri$ preserves connected limits on both sides.
That is, if $\cat{J}$ is a connected category, $p\colon \cat{J}\to\poly$ is a functor, and $q\in\poly$ is a polynomial, then there are natural isomorphisms\index{isomorphism!natural}
\[
	\left(\lim_{j\in \cat{J}} p_j\right)\tri q\iso \lim_{j\in \cat{J}}(p_j\tri q)
	\qqand
	q\tri\left(\lim_{j\in \cat{J}} p_j\right)\iso \lim_{j\in \cat{J}}(q\tri p_j)
\]
\end{theorem}
\begin{proof}\index{pullback!wide}
The claim for the left side is just a special case of \cref{prop.left_pres_lim}; it remains to prove the claim on the right.

By \cref{thm.poly_limits}, $\poly$ is complete, so by \cite[Theorem~4.3]{nlab:connected-limit}, it suffices to show that the functor $(q\tri-)\colon\poly\to\poly$ preserves wide pullbacks on the right.
By \cref{prop.comp_right_slice_coclosed}, the functor $(q\tri-)\colon\poly\to\poly/q\tri\1$ is a right adjoint, so it preserves limits, including wide pullbacks.
Thus $(q\tri-)$ sends a wide pullback over $r\in\poly$ to a wide pullback over the canonical lens $q\tri r\to q\tri\1$ in $\poly/q\tri\1$, corresponding to a limit in $\poly$ of a diagram consisting of arrows to $q\tri r$ and arrows to $q\tri\1$ factoring through $q\tri r$.
So up to isomorphism, this limit is just a wide pullback in $\poly$ over $q\tri r$, namely $(q\tri-)\colon\poly\to\poly$ applied to the original wide pullback.
So $\tri$ preserves wide pullbacks on the right.
\end{proof}

\begin{exercise}\label{ex.connected_limits_and_tri}
Use \cref{thm.connected_limits} in the following.
\begin{enumerate}
	\item Let $p$ be a polynomial, thought of as a functor $p\colon\smset\to\smset$. Show that $p$ preserves connected limits (of sets).
	\item Show that for any polynomials $p,q,r$ we have an isomorphism:
	\begin{equation} \label{eqn.right_prod_pullback}
	p\tri(qr)\iso (p\tri q)\times_{(p\tri\1)}(p\tri r).
	\end{equation}
	\item Take the polynomials $p,q,r$ from the counterexample you found in \cref{exc.right_not_dist_prod} \cref{exc.right_not_dist_prod.prod} and check that \eqref{eqn.right_prod_pullback} holds.
\qedhere
\end{enumerate}
\begin{solution}
\begin{enumerate}
    \item Given a polynomial functor $p \colon \smset \to \smset$, we wish to show that $p$ preserves connected limits of sets; that is, for a connected category $\cat{J}$ and a functor $X \colon \cat{J} \to \smset$, we have
    \[
        p\left(\lim_{j \in \cat{J}} X_j\right) \iso \lim_{j \in \cat{J}} p(X_j).
    \]
    But we can identify $\smset$ with the full subcategory of constant functors in $\poly$ and instead view $X$ as a functor into $\poly$.
    Then by \cref{exc.composing_with_constants} \cref{exc.composing_with_constants.appl}, the left hand side of the isomorphism we seek is isomorphic to $p\tri\left(\lim_{j \in \cat{J}} X_j\right)$, while the right hand side is isomorphic to $\lim_{j \in \cat{J}} \left(p \tri X_j\right)$.
    These are isomorphic by \cref{thm.connected_limits}.
    \item Given $p,q,r \in \poly$, we wish to show that \eqref{eqn.right_prod_pullback} holds.
    As $\1$ is terminal in $\poly$, the product $qr$ can also be written as the pullback $q \times_\1 r$.
    While products are not connected limits, pullbacks are, so by \cref{thm.connected_limits}, they are preserved by precomposition with $p$.
    Hence the desired isomorphism follows.
    \item We'll show that \eqref{eqn.right_prod_pullback} holds for $p\coloneqq\yon+\1, q\coloneqq\1,$ and $r\coloneqq\0$.
    We have $p\tri q = p\tri\1 \iso (\yon+\1)\tri\1 \iso \2$ and $p\tri r \iso (\yon+\1)\tri\0 \iso \1$, so $(p\tri q)\times_{(p\tri\1)}(p\tri r) \iso \2 \times_\2 \1$.
    We saw in \cref{ex.pullbacks_in_poly} that the position-set of a pullback in $\poly$ is just the pullback of the position-sets in $\smset$, while the direction-sets are given by a pushout of direction-sets in $\smset$.
    As our polynomials have empty direction-sets, their pullback must have an empty direction-set as well, so this pullback is just a pullback of sets: $(p\tri q)\times_{(p\tri\1)}(p\tri r) \iso \2 \times_\2 \1\iso\1$.
    And indeed we have $p\tri(qr) \iso (\yon+\1)\tri\0 \iso \1$ as well.
\end{enumerate}
\end{solution}
\end{exercise}

While we are here, it will be helpful to record the following.
\begin{proposition} %Make this an exercise, mimicking earlier stuff
For any polynomial $q\in\poly$, tensoring with $q$ (on either side) preserves connected limits. That is, if $\cat{J}$ is connected and $p\colon \cat{J}\to\poly$ is a functor, then there is a natural isomorphism:\index{isomorphism!natural}
\[
	\left(\lim_{j\in \cat{J}} p_j\right)\otimes q\cong
	\lim_{j\in \cat{J}} (p_j\otimes q).
\]
\end{proposition}

\subsection{Interaction with parallel products} \label{subsec.comon.comp.prop.par}

Before we get into how $\otimes$ interacts with $\tri$, here is a warm-up exercise.

\begin{exercise}
Let $A$ and $B$ be arbitrary sets, and let $p$ be an arbitrary polynomial.
Which of the following isomorphisms always hold?

If the isomorphism does not always hold, is there still a canonical lens in one direction or the other?
\begin{enumerate}
	\item $(A\yon)\otimes(B\yon) \iso^? (A\yon)\tri (B\yon)$.
	\item $\yon^A\otimes\yon^B\iso^?\yon^A\tri\yon^B$.
	\item $A\otimes B\iso^? A\tri B$.
	\item $B\yon\otimes p\iso^? B\yon\tri p$.
	\item $\yon^A\otimes p\iso^? \yon^A\tri p$.
	\item $p\otimes B\yon\iso^? p\tri B\yon$.
	\item $p\otimes \yon^A\iso^? p\tri\yon^A$.
	\qedhere
\end{enumerate}
What do all of the lenses you found in this exercise have in common (whether or not they were isomorphisms)?
\begin{solution}
Here $A$ and $B$ are sets and $p$ is a polynomial.
\begin{enumerate}
	\item The isomorphism always holds: we have that $(A\yon)\otimes(B\yon) \iso AB\yon \iso (A\yon)\tri (B\yon)$.
	\item The isomorphism always holds: we have that $\yon^A\otimes\yon^B \iso \yon^{AB} \iso \yon^A\tri\yon^B$.
	\item The isomorphism does not always hold: while $A\otimes B \iso AB$, we have that $A\tri B \iso A$.
	There is, however, always a canonical projection $AB\to B$; but there is not always a canonical lens $B\to AB$ (for example, take $A=\0\neq B$).
	\item The isomorphism always holds: we have that $B\yon \otimes p \iso \sum_{i \in p(\1)} B\yon \otimes \yon^{p[i]} \iso \sum_{i \in p(\1)} B\yon^{p[i]} \iso Bp \iso B\yon\tri p$.
	\item The isomorphism does not always hold: if, say, $p = B$, then $\yon^A \otimes B \iso B$, while $\yon^A \tri B \iso B^A$.
	If $A=B=\0$, then $B^A\iso\1$, so there is not always a canonical lens from right to left, either.
	There is, however, always a canonical lens from left to right: $\yon^A\otimes p\iso\sum_{i \in p(\1)}\yon^{Ap[i]}$ while $\yon^A\tri p\iso\sum_{\ol{i}\colon A\to p(\1)}\yon^{\sum_{a\in A}p[\ol{i}(a)]}$.
	So there is a lens from left to right whose on-positions function sends $i\in p(\1)$ to the constant function $A\to p(\1)$ that always evaluates to $i$; and whose on-directions function at $i$ is the identity on $Ap[i]$.
	\item The isomorphism does not always hold: if, say, $p = \yon^A$, then $\yon^A \otimes B\yon \iso B\yon^A$, while $\yon^A \otimes B\yon \iso (B\yon)^A \iso B^A\yon^A$.
	If $A=B=\0$, then $B\yon^A\iso\0$ while $B^A\yon^A\iso\0^\0\yon^\0\iso\1$, so there is not always a canonical lens from right to left, either.
	There is, however, always a canonical lens from left to right: $p\otimes B\yon\iso Bp$ while $p\tri B\yon\iso\sum_{i\in p(\1)}\sum_{\ol{b}\colon p[i]\to B}\yon^{\sum_{a\in p[i]}\1}\iso\sum_{i\in p(\1)}B^{p[i]}\yon^{p[i]}$.
	So there is a lens from left to right whose on-positions function sends $(b,i)\in Bp(\1)$ to $(i,c_b)\in\sum_{i\in p(\1)}B^{p[i]}$, where $c_b\colon p[i]\to B$ is the constant function that always evaluates to $b$; and whose on-directions function at $(b,i)$ is the identity on $p[i]$.
	\item The isomorphism always holds: we have that $p \otimes \yon^A \iso \sum_{i \in p(\1)} \yon^{Ap[i]} \iso p \tri \yon^A$.
\end{enumerate}
Every on-directions function of every lens we found in this exercise are isomorphisms, so every lens we found in this exercise is cartesian.
\end{solution}
\end{exercise}

\index{polybox}\index{indexed family}

\begin{example}[Lenses from $\otimes$ to $\tri$]
For any $p$ and $q$, there is an interesting cartesian lens $o_{p,q}\colon p\otimes q\to p\tri q$ that, stated informally, ``orders'' the operation, taking the symmetric monoidal product $\otimes$ and reinterprets it as a special case of the asymmetric monoidal product $\tri$.
Defining this lens in the usual way is rather tedious and unilluminating, but written in polyboxes, the lens looks like this (recall that positions of $p\otimes q$ are just pairs of positions of $p$ and $q$, while directions at each such pair are pairs of directions of $p$ and $q$, one at each position in the pair; we drop the parentheses around the ordered pair for readability):
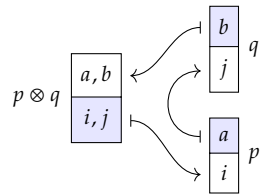
\[
\begin{tikzpicture}[polybox, mapstos]
	\node[poly, dom, "$p\otimes q$" left] (p) {$a,b$\at$i,j$};
	\node[poly, cod, "$p$" right, right= 1.5cm of p.south, yshift=-1ex] (q) {$a$\at$i$};
	\node[poly, cod, "$q$" right, above=of q] (r) {$b$\at$j$};
  	\draw (p_pos) to[first] (q_pos);
  	\draw (q_dir) to[climb] (r_pos);
  	\draw (r_dir) to[last] (p_dir);
\end{tikzpicture}
\]
Usually, the positions box of $q$ is allowed to depend on the directions box of $p$ in the polyboxes for $p\tri q$ on its own.
But in the polyboxes above, $j$ is not allowed to depend on $a$ in $p\otimes q$ on the left, so the arrow from the positions box of $q$ to the directions box of $p$ on the right doesn't actually take $a$ into account at all.
So the lens $o_{p,q}$ is in some sense the inclusion of the order-independent positions of $p\tri q$; when drawn as trees, the positions in its image are the ones whose upper-level corollas are all the same.
And of course we can flip the order using the symmetry $q\otimes p\iso p\otimes q$.
This is, we just as well have a lens $p\otimes q\to q\tri p$.

Both $\otimes$ and $\tri$ have the same monoidal unit, the identity functor $\yon$, whose identity is the unique lens $\yon\to\yon$.
In fact the lenses $o_{p,q}$ constitute a lax monoidal functor $(\poly,\yon,\otimes)\to(\poly,\yon,\tri)$.
In particular, $o_{p,q}$ commutes with associators and unitors.

This can be used in the following way. Lenses $p\to q\tri r$ into composites are fairly easy to understand (through polyboxes, for example), whereas lenses $q\tri r\to p$ are not so easy to think about. However, given such a lens, one may always compose it with $o_{q,r}$ to obtain a lens $q\otimes r\to p$.
This is quite a bit simpler to think about: they are our familiar interaction patterns from \cref{sec.poly.dyn_sys.interact}.
\end{example}\index{interaction pattern}

\index{duoidality}

It turns out that the monoidal structures $\otimes$ and $\tri$ together satisfy an interesting property known as \emph{duoidality}.
We won't give the entire definition of what it means for two monoidal structures to be duoidal here---there are a few commutative diagrams to verify for technical reasons---but the key condition is that there is a natural lens
\begin{equation} \label{eqn.duoidal}
    (p\tri p')\otimes(q\tri q')\to(p\otimes q)\tri(p'\otimes q').
\end{equation}

\begin{proposition} \label{prop.duoidal}
The monoidal structures $(\yon,\otimes)$ and $(\yon,\tri)$ together comprise a duoidal structure on $\poly$.
\end{proposition}
\begin{proof}[Idea of proof]
The key is to give the natural lens from \eqref{eqn.duoidal}, as follows.
A position of $p\tri p'$ is a pair $(i,\ol{i'})$ with $i\in p(\1)$ and $\ol{i'}\colon p[i]\to p'(\1)$; similarly, a position $q\tri q'$ is a pair $(j,\ol{j'})$ with $j\in q(\1)$ and $\ol{j'}\colon q[j]\to q'(\1)$.
So we can define the first two parts of the lens using polyboxes, like so (again we drop parentheses around some ordered pairs for readability):
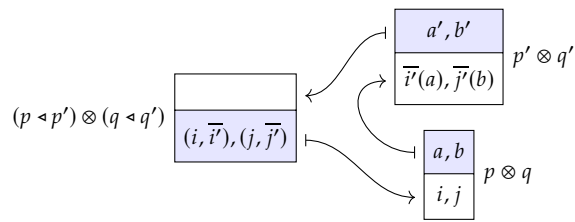
\[
\begin{tikzpicture}[polybox, mapstos]
	\node[poly, dom, "$(p\tri p')\otimes(q\tri q')$" left] (p) {\at$(i,\ol{i'}),(j,\ol{j'})$};
	\node[poly, cod, "$p\otimes q$" right, right=2.5 of p.south, yshift=-1ex] (q) {$a,b$\at$i,j$};
	\node[poly, cod, "$p'\otimes q'$" right, above=of q] (r) {$a',b'$\at$\ol{i'}(a),\ol{j'}(b)$};
  	\draw (p_pos) to[first] (q_pos);
  	\draw (q_dir) to[climb] (r_pos);
  	\draw (r_dir) to[last] (p_dir);
\end{tikzpicture}
\]
Here $(a,b)$ is a direction of $p\otimes q$ at $(i,j)$, with $a\in p[i]$ and $b\in q[j]$.

Then to fill in the remaining empty box, we need a direction of $p\tri p'$ at $(i,\ol{i'})$, which can be given by the $p[i]$-direction $a$ followed by a $p'[\ol{i'}(a)]$-direction, namely $a'$.
We also need a direction of $q\tri q'$, which can be obtained analogously:
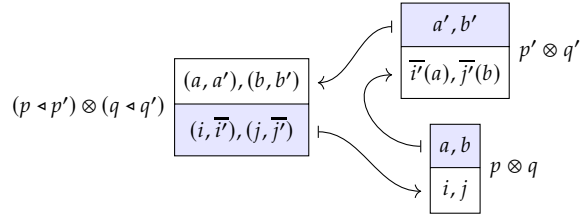
\[
\begin{tikzpicture}[polybox, mapstos]
	\node[poly, dom, "$(p\tri p')\otimes(q\tri q')$" left] (p) {$(a,a'),(b,b')$\at$(i,\ol{i'}),(j,\ol{j'})$};
	\node[poly, cod, "$p\otimes q$" right, right=2.5 of p.south, yshift=-1ex] (q) {$a,b$\at$i,j$};
	\node[poly, cod, "$p'\otimes q'$" right, above=of q] (r) {$a',b'$\at$\ol{i'}(a),\ol{j'}(b)$};
  	\draw (p_pos) to[first] (q_pos);
  	\draw (q_dir) to[climb] (r_pos);
  	\draw (r_dir) to[last] (p_dir);
\end{tikzpicture}
\]
\end{proof}

\subsection{Interaction with vertical and cartesian lenses} \label{subsec.comon.comp.prop.cart}

\index{lens!vertical}\index{lens!cartesian}

We conclude this section by examining how $\tri$ interacts with vertical and cartesian lenses, as defined in \cref{def.vert_cart}.

\index{composition product!cartesian lenses and}
\begin{proposition}[Preservation of cartesian lenses]\label{prop.comp_pres_cart}
If $\varphi\colon p\to p'$ and $\psi\colon q\to q'$ are cartesian lenses, then so is $\varphi\tri \psi\colon p\tri q \to p'\tri q'$.
\end{proposition}
\begin{proof}
We use the third characterization of cartesian lenses given in \cref{prop.cart_as_nt}, as lenses whose naturality squares are pullbacks.
For any sets $A,B$ and function $h\colon A\to B$, consider the diagram
\[
\begin{tikzcd}[column sep=small]
    p\tri q\tri A\ar[r]\ar[d] & p'\tri q\tri A\ar[r]\ar[d] & p'\tri q'\tri A\ar[d] \\
    p\tri q\tri B\ar[r] & p'\tri q\tri B\ar[r] & p'\tri q'\tri B.
\end{tikzcd}
\]
The square on the left is a pullback because $\varphi\colon p\to p'$ is cartesian.
Meanwhile, the square on the right is a pullback because $\psi\colon q\to q'$ is cartesian and $\tri$ preserves pullbacks by \cref{thm.connected_limits}.
Hence the outer rectangle is a pullback as well, implying that $\varphi\tri \psi\colon p\tri q \to p'\tri q'$ is cartesian.
\end{proof}

\index{composition product!vertical lenses and}

\begin{exercise}
Let $\varphi\colon p\to p'$ and $\psi\colon q\to q'$ be lenses.
\begin{enumerate}
	\item Show that if $\varphi$ is an isomorphism and $\psi$ is vertical, then $\varphi\tri \psi$ is vertical.
	\item Find a vertical lens $\varphi$ and a polynomial $q$ for which the lens $\varphi\tri q\colon p\tri q\to p'\tri q$ is not vertical.
\qedhere
\end{enumerate}
\begin{solution}
Here $\varphi\colon p\to p'$ and $\psi\colon q\to q'$ are lenses.
\begin{enumerate}
    \item If $\varphi$ is an isomorphism and $\psi$ is vertical, then $\psi\tri\1\colon q\tri\1\to q'\tri\1$ is an isomorphism, so $\varphi\tri \psi\tri\1\colon p\tri q\tri\1\to p'\tri q'\tri\1$ is an isomorphism as well.
    Thus $\varphi\tri \psi$ is vertical.
    \item If $\varphi$ is the unique lens $\yon\to\1$ and $q=\0$, then $\varphi$ is vertical, but since $\yon\tri\0\iso\0$ and $\1\tri\0\iso\1$, the lens $\varphi\tri\0\colon\0\to\1$ is not.
\end{enumerate}
\end{solution}
\end{exercise}

%-------- Section --------%
\section[Summary and further reading]{Summary and further reading%
  \sectionmark{Summary \& further reading}}
\sectionmark{Summary \& further reading}

In this chapter we introduced the composition (sometimes called ``substitution'') product $\tri$. Given polynomials $p,q$ thought of as functors, their composite is again polynomial and is given by $p\tri q$. We explained how it looks in terms of algebra, e.g.\ how to compute $\yon^\2\tri(\yon+\1)$; in terms of sets, e.g.\ as a sum-product-sum-product $\sum\prod\sum\prod$, which can be reduced to a single $\sum\prod$; in terms of trees, by stacking corollas on top of corollas; and in terms of polyboxes. We paid particular attention to lenses into a composite $p\to q\tri r$.

This allowed us to explain how to think of dynamical systems $S\yon^S\to p\tri q$ with composite interfaces as multi-step machines: each state produces a $p$-position $i$, and then for every $p[i]$-direction produces a $q$-position $j$, and finally for every $q[j]$-direction returns an updated state in $S$.

Finally, we discussed some facts of the composition product. For example, we showed that $-\tri q$ has a left adjoint $\lchom{q}{-}$ and that $q\tri -$ has a left multi-adjoint $-\frown q$. We also explained that $p\tri q$ preserves all with limits in the variable $p$ and all connected limits in the variable $q$. We also explained the duoidal interaction between $\otimes$ and $\tri$, i.e.\ the natural lens $\tri\otimes\tri\to\otimes\tri\otimes$, and how $\tri$ interacts with cartesian lenses.

Polynomial substitution is one of the best known aspects of polynomial functors. Again, see \cite{kock2012polynomial} for more on this. We learned of the left coclosure (see \cref{prop.comp_left_coclosed}) from Josh Meyers, though it may have already been known in the containers community.
%-------- Section --------%
\section{Exercise solutions}
\Closesolutionfile{solutions}
{\footnotesize
\input{solution-file6}}

\Opensolutionfile{solutions}[solution-file7]

\index{composition product|)}

%------------ Chapter ------------%
\chapter{Polynomial comonoids and retrofunctors}
\chaptermark{Comonoids in $\poly$ \& retrofunctors}
\label{ch.comon.sharp}
\index{polynomial comonoid|(}
\index{polynomial comonad|see{polynomial comonoid}}

\slogan{Imagine a realm where there are various positions you can be in. From every position, there are a number of moves you can make, possibly infinitely many. But whatever move you make, you'll end up in a new position. Well, technically it counts as a move to simply stay where you are, so you might end up in the same position. But wherever you move to, you can move again, and any number of moves from the original position counts as a single move. What sort of realm is this?}

The most surprising aspects of $\poly$ really begin with its comonoids.
In 2018, researchers Daniel Ahman and Tarmo Uustalu presented a characterization of comonoids in $(\poly,\yon,\tri)$ as a surprisingly familiar construct.
For us, this story will emerge naturally as we continue to expand our understanding of the humble state system of a dependent dynamical system.

%-------- Section --------%
\section{State systems, categorically}\label{sec.comon.sharp.state}

Since defining dependent dynamical systems in \cref{def.gen_moore}, we have evolved our understanding of their state systems over the course of the last few chapters.
Let's take this moment to review what we know about these state systems so far.

\index{state system}

Our original definition of a state system was as a monomial $S\yon^S$ for some set $S$.
But in \cref{ex.do_nothing}, we noted that this formulation requires us to discuss the positions and directions of a state system at the level of sets rather than in the language of $\poly$.
Instead, let's take an arbitrary polynomial $\car{s}\in\poly$ and attempt to characterize what it means for $\car{s}$ to be a state system using only the categorical machinery of $\poly$.
We will continue to refer to the positions of $\car{s}$ as \emph{states}, but we will shift from thinking of the directions of $\car{s}$ as states to thinking of them as transitions from one state to another.

\subsection{The do-nothing section}\label{subsec.comon.sharp.state.nothing}

In \cref{ex.do_nothing}, we saw that every state system $\car{s}$ is equipped with a \emph{do-nothing section}: a lens $\epsilon\colon\car{s}\to\yon$ that picks out a direction at each state that we would like to interpret as ``doing nothing'' and remaining at that state.

We drew $\epsilon$ in polyboxes in \cref{ex.do_nothing_polybox}, but that was when we let ourselves assume that the position-set of a state system was equal to each of its direction-sets.
Now all we know is that for each state $s\in\car{s}(\1)$, the do-nothing section chooses an $\car{s}[s]$-direction to signify staying at the same state; it doesn't make sense to say that this direction is literally equal to $s$.
So we need a different name for the $\car{s}[s]$-direction that $\epsilon$ identifies: call it $\id_s$, because it behaves like a sort of identity operation on the state $s$.

So the revised polyboxes for the do-nothing section $\epsilon\colon\car{s}\to\yon$ are as follows:
\begin{equation} \label{eqn.do_nothing_polybox}
\begin{tikzpicture}[polybox, mapstos]
    \node[poly, dom, "$\car{s}$" left] (S) {$\id_s$\at$s$};

    \draw (S_pos) to[climb'] node[right] {$\epsilon$} (S_dir);
\end{tikzpicture}
\end{equation}

\index{section!do nothing}

\begin{exercise}
Say I have a polynomial $\car{s}\in\poly$, and I tell you that there is a lens $\epsilon\colon\car{s}\to\yon$.
What can you say about the polynomial $\car{s}$?
\begin{solution}
Given a polynomial $\car{s}\in\poly$ equipped with a lens $\epsilon\colon\car{s}\to\yon$, we know that $\epsilon$ picks out a direction at every position of $\car{s}$.
So all we can say about $\car{s}$ is that there is at least one direction at each of its positions.
Equivalently, we could say that $\car{s}$ can be written as the product of $\yon$ and some other polynomial.
\end{solution}
\end{exercise}

\begin{example}[The do-nothing section in tree pictures] \label{ex.nothing_trees}
We have seen the do-nothing section drawn in polyboxes, but let us see what it looks like in our tree pictures.
We take $\car{s}\coloneqq
\{
{\color{my-lavender}\blacksquare},
{\color{my-yellow}\blacktriangle},
{\color{my-magenta}\blacklozenge}
\}
\yon^{
\{
{\color{my-lavender}\blacksquare},\,
{\color{my-yellow}\blacktriangle},\,
{\color{my-magenta}\blacklozenge}
\}
}$, drawn as follows, with the directions associated with $\color{my-lavender}\blacksquare$, $\color{my-yellow}\blacktriangle$, and $\color{my-magenta}\blacklozenge$ from left to right:
\[
\begin{tikzpicture}[rounded corners]
\node (p1) [draw, "$\car{s}\coloneqq$" left] {
    \begin{tikzpicture}[trees, sibling distance=4mm]
      \node[my-lavender, below] (1) {\small$\blacksquare$}
      child[my-lavender] {}
      child[my-yellow] {}
      child[my-magenta] {};
      \node[my-yellow, right=1.8 of 1] (2) {\small$\blacktriangle$}
      child[my-lavender] {}
      child[my-yellow] {}
      child[my-magenta] {};
      \node[my-magenta, right=1.8 of 2] (3) {\small$\blacklozenge$}
      child[my-lavender] {}
      child[my-yellow] {}
      child[my-magenta] {};
    \end{tikzpicture}
};
\end{tikzpicture}
\]
Then the do-nothing section $\epsilon\colon\car{s}\to\yon$ can be drawn as follows:
\[
\begin{tikzpicture}[trees, bend right]
    \foreach \i/\c/\s in {1/my-lavender/\blacksquare, 2/my-yellow/\blacktriangle, 3/my-magenta/\blacklozenge}
    {
        \node[\c] (\i) at (3*\i, 0) {\small$\s$}
            child [my-lavender] {coordinate (\i1)}
            child [my-yellow] {coordinate (\i2)}
            child [my-magenta] {coordinate (\i3)}
            ;
        \node[right=of \i] (y\i) {$\bullet$}
            child{coordinate (y\i')}
            ;
        \draw[|->, shorten <= 3pt, shorten >= 3pt] (\i) -- (y\i);
        \draw[densely dotted, postaction={decorate}] (y\i') to (\i\i);
    };
\end{tikzpicture}
\]
It picks out one direction at each position, namely the one of the same color.
\end{example}

There is not much else we can say about the do-nothing section on its own, so let us revisit the other lens that every state system is equipped with before considering the relationship between the two.

\index{lens!transition}

\subsection{The transition lens}\label{subsec.comon.sharp.state.trans}

\index{polybox}

We saw in \cref{ex.dyn_sys_trans_polyboxes} that $\car{s}$ also comes equipped with a \emph{transition lens}: a lens $\delta\colon\car{s}\to\car{s}\tri\car{s}$, which we can draw as
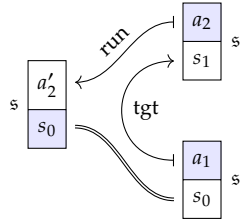
\[
\begin{tikzpicture}[polybox, mapstos]
	\node[poly, dom, "$\car{s}$" left] (r) {$a_2'$\at$s_0$};
	\node[poly, cod, right=1.8 of r.south, yshift=-2.5ex, "$\car{s}$" right] (p) {$a_1$\at$s_0$};
	\node[poly, cod, above=.8 of p, "$\car{s}$" right] (p') {$a_2$\at$s_1$};

	\draw[double, -] (r_pos) to[first] node[below] {} (p_pos);
	\draw (p_dir) to[climb] node[right] {tgt} (p'_pos);
	\draw (p'_dir) to[last] node[above,sloped] {run} (r_dir);
  \end{tikzpicture}
\]
The arrow labeled tgt is the \emph{target function}: given a state $s_0$ and a direction $a_1$ at that state, $\text{tgt}(s_0,a_1)$ tells us the new state $s_1$ that following $a_1$ from $s_0$ will lead to.
We know that when $s_0$ is fixed, the target function on the second component $a_1$ should be an isomorphism $\car{s}[s_0]\to\car{s}(\1)$; that is, there is exactly one direction at $s_0$ that leads to each state of $\car{s}$.
But this property is a little tricky to state in the language of $\poly$; in fact, we won't attempt to do so just yet.
Instead, we will use it to make a notational choice: given $s,t\in\car{s}$, we will let $s\to t$ denote the unique direction at $s$ that leads to $t$, so that $\text{tgt}(s,s\to t)=t$.
So we can redraw our transition lens as
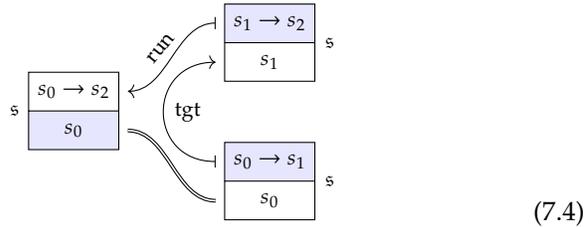
\begin{equation} \label{eqn.trans_lens_polybox}
\begin{tikzpicture}[polybox, mapstos]
	\node[poly, dom, "$\car{s}$" left] (r) {$s_0\to s_2$\at$s_0$};
	\node[poly, cod, right=2 of r.south, yshift=-2.5ex, "$\car{s}$" right] (p) {$s_0\to s_1$\at$s_0$};
	\node[poly, cod, above=.8 of p, "$\car{s}$" right] (p') {$s_1\to s_2$\at$s_1$};

	\draw[double, -] (r_pos) to[first] node[below] {} (p_pos);
	\draw (p_dir) to[climb] node[right] {tgt} (p'_pos);
	\draw (p'_dir) to[last] node[above,sloped] {run} (r_dir);
  \end{tikzpicture}
\end{equation}
In addition to the fact that $\text{tgt}(s_0,s_0\to s_1)=s_1$ as intended, this picture tells us two more properties of $\delta$.

The first is that the bottom arrow is the identity on $\car{s}(\1)$.
This is something we would like to be able to express categorically in the language of $\poly$.
We'll see that this property falls out naturally when we express how the transition lens plays nicely with the do-nothing section in \cref{subsec.comon.sharp.state.cohere}.

The second is that the run arrow, which runs the transition $s_0\to s_1$ and the transition $s_1\to s_2$ together into a transition starting at $s_0$, should have the same target as the second transition it follows: in this case, $s_2$.
Equationally, writing the left and right hand sides only in terms of the contents of the blue boxes, we have that
\begin{equation} \label{eqn.state_run_tgt}
    \text{tgt}(s_0,\text{run}(s_0,s_0\to s_1,s_1\to s_2))=s_2=\text{tgt}(\text{tgt}(s_0,s_0\to s_1),s_1\to s_2).
\end{equation}
We will see that this property arises naturally when we we generalize the transition lens to more than two steps in \cref{subsec.comon.sharp.state.coassoc}.

\begin{example}[The transition lens in tree pictures] \label{ex.trans_trees}
Continuing from \cref{ex.nothing_trees}, we draw the transition lens $\delta\colon\car{s}\to\car{s}\tri\car{s}$ of $\car{s}\coloneqq\3\yon^\3\iso\{
{\color{my-lavender}\blacksquare},
{\color{my-yellow}\blacktriangle},
{\color{my-magenta}\blacklozenge}
\}
\yon^{
  \{
  {\color{my-lavender}\blacksquare},\,
  {\color{my-yellow}\blacktriangle},\,
  {\color{my-magenta}\blacklozenge}
  \}
}$ (where directions are labeled with their targets) in tree pictures as well, recalling that the trees of $\car{s}\tri\car{s}$ are obtained by taking an $\car{s}$-corolla and grafting more $\car{s}$-corollas to each of its leaves:
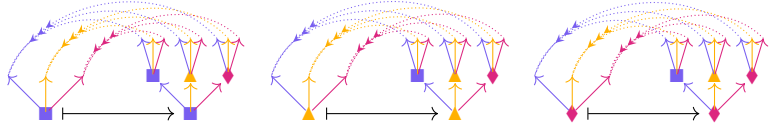
\[
\begin{tikzpicture}[trees,
  level 1/.style={sibling distance=5mm},
  level 2/.style={sibling distance=1.5mm},
	bend right=60]
	\foreach \i/\c/\s in {1/my-lavender/\blacksquare, 2/my-yellow/\blacktriangle, 3/my-magenta/\blacklozenge}
	{
  	\node[\c] (\i) at (3.5*\i, 0) {\small$\s$}
    	child [my-lavender] {coordinate (\i1)}
      child [my-yellow] {coordinate (\i2)}
      child [my-magenta] {coordinate (\i3)}
     	;
  	\node[right=1.7 of \i, \c] (SS\i) {\small$\s$}
  		child [my-lavender] {node (S\i1) {\small$\blacksquare$}
				child [my-lavender] {coordinate (\i11)}
				child [my-yellow] {coordinate (\i12)}
				child [my-magenta] {coordinate (\i13)}
				}
  		child [my-yellow] {node (S\i2) {\small$\blacktriangle$}
				child [my-lavender] {coordinate (\i21)}
				child [my-yellow] {coordinate (\i22)}
				child [my-magenta] {coordinate (\i23)}
				}
  		child [my-magenta] {node (S\i3) {\small$\blacklozenge$}
				child [my-lavender] {coordinate (\i31)}
				child [my-yellow] {coordinate (\i32)}
				child [my-magenta] {coordinate (\i33)}
				}
  		;
	\draw[|->, shorten <= 3pt, shorten >= 3pt] (\i) -- (SS\i);
	\foreach \j in {1,2,3}
	{
		\foreach \k\d in {1/my-lavender, 2/my-yellow, 3/my-magenta}
		{
			\draw[densely dotted, postaction={decorate}, \d] (\i\j\k) to (\i\k);
		};
	};
	};
\end{tikzpicture}
\]
On positions, the target function of $\delta$ tells us which root of $\car{s}$ to graft onto each leaf of $\car{s}$.
Then on directions, the run function of $\delta$ tells us how to collapse the height-$2$ leaves of the trees we obtain in $\car{s}\tri\car{s}$ down to the original height-$1$ leaves of the corollas of $\car{s}$.

We can draw what the target function is doing more compactly by taking the corollas of $\car{s}$ and ``bending the arrows'' so that they point to their targets, like so:
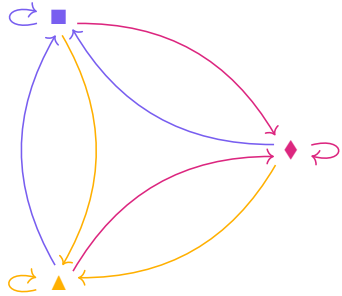
\[
\begin{tikzcd}
  {\color{my-lavender}\blacksquare} \ar[ddrrr, bend left, shift right, my-magenta] \ar[dddd, bend left, shift right, my-yellow] \ar[loop left, my-lavender] \\
  \\
  &&& {\color{my-magenta}\blacklozenge} \ar[uulll, bend left, shift right, my-lavender] \ar[ddlll, bend left, shift right, my-yellow] \ar[loop right, my-magenta] \\
  \\
  {\color{my-yellow}\blacktriangle} \ar[uuuu, bend left, shift right, my-lavender] \ar[uurrr, bend left, shift right, my-magenta] \ar[loop left, my-yellow]
\end{tikzcd}
\]
So the target function of $\delta$ turns our corolla picture of $\car{s}$ into a complete graph on its roots.
Then the run function takes any two arrows that form a path in the graph and collapses them down to a single arrow that starts (and, according to \eqref{eqn.state_run_tgt}, ends) at the same vertex as the two-arrow path.
\end{example}\index{graph!complete}

\subsection[Do-nothing section and transition lens coherence]{The do-nothing section coheres with the transition lens}\label{subsec.comon.sharp.state.cohere}

For each state $s\in\car{s}(\1)$, the do-nothing section $\epsilon\colon\car{s}\to\yon$ picks out the $\car{s}[s]$-direction $\id_s$ that ``does nothing'' and keeps the system in the same state $s$.
But it is the transition lens $\delta\colon\car{s}\to\car{s}\tri\car{s}$ that actually sets our state system in motion, specifying the target of each direction and how two directions run together.
Either of these directions could be our do-nothing direction $\id_s$, so let's try to figure out what should happen when we set each one in turn to $\id_s$.

\index{polybox}

We can draw in polyboxes what happens when we set $\id_s$, as specified by $\epsilon$, to be the first direction that $\delta$ runs together like this:
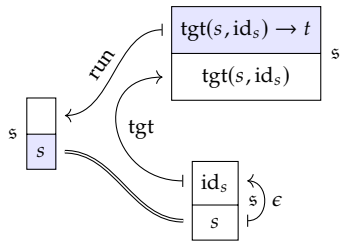
\[
\begin{tikzpicture}[polybox, mapstos]
	\node[poly, dom, "$\car{s}$" left] (r) {\at$s$};
	\node[poly, right=2 of r.south, yshift=-2.5ex, "$\car{s}$" right] (p) {$\id_s$\at$s$};
	\node[poly, cod, above=.8 of p, xshift=2.5ex, "$\car{s}$" right] (p') {tgt$(s,\id_s)\to t$\at tgt$(s,\id_s)$};

	\draw[double, -] (r_pos) to[first] node[below] {} (p_pos);
	\draw (p_pos) to[climb'] node[right] {$\epsilon$} (p_dir);
	\draw (p_dir) to[climb] node[right] {tgt} (p'_pos);
	\draw (p'_dir) to[last] node[above,sloped] {run} (r_dir);
\end{tikzpicture}
\]
Reading this picture from left to right, we see that it depicts the polyboxes of the composite lens $\delta\then(\epsilon\tri\car{s})\colon\car{s}\to\yon\tri\car{s}\iso\car{s}$ (recall that we sometimes denote the identity lens on $\car{s}$ also by $\car{s}$).
To make this interpretation more transparent, we could be a little more verbose with our polybox picture if we wanted to (omitting the contents of the boxes for clarity):
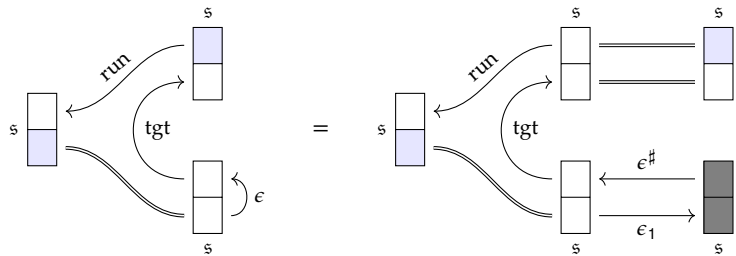
\[
\begin{tikzpicture}
    \node (1) {
        \begin{tikzpicture}[polybox, tos]
        	\node[poly, dom, "$\car{s}$" left] (r) {};
        	\node[poly, right=2 of r.south, yshift=-2.5ex, "$\car{s}$" below] (p) {};
        	\node[poly, cod, above=.8 of p, "$\car{s}$" above] (p') {};

        	\draw[double, -] (r_pos) to[first]  (p_pos.west);
        	\draw (p_pos) to[climb'] node[right] {$\epsilon$} (p_dir);
        	\draw (p_dir) to[climb] node[right] {tgt} (p'_pos);
        	\draw (p'_dir) to[last] node[above,sloped] {run} (r_dir);
        \end{tikzpicture}
	};
	\node[right=1 of 1] (2) {
	    \begin{tikzpicture}[polybox, tos]
            \node[poly, dom, "$\car{s}$" left] (yX) {};
        	\node[poly, right=2 of yX.south, yshift=-2.5ex, "$\car{s}$" below] (p) {};
            \node[poly, above=.8 of p, "$\car{s}$" above] (p') {};
            \node[poly, cod, identity, right=of p, "$\car{s}$" below] (q) {};
            \node[poly, cod, above=.8 of q, "$\car{s}$" above] (q') {};
            \draw (p_pos) -- node[below] {$\epsilon_\1$} (q_pos);
            \draw (q_dir) -- node[above] {$\epsilon^\sharp$} (p_dir);
            \draw[double, -] (p'_pos) -- (q'_pos);
            \draw[double, -] (q'_dir) -- (p'_dir);
            \draw[double, -] (yX_pos) to[first] node[below] {} (p_pos);
            \draw (p_dir) to[climb] node[right] {tgt} (p'_pos);
            \draw (p'_dir) to[last] node[above,sloped] {run} (yX_dir);
        \end{tikzpicture}
	};
	\node at ($(1.east)!.5!(2.west)$) {=};
\end{tikzpicture}
\]
Now what should $\text{tgt}(\id_s)$ be, and what should go in the direction box on the left?

If following the direction $\id_s$ from the state $s$ is really the same as doing nothing, then its target state should be the same state $s$ that it emerged from.
Moreover, running together $\id_s$ with any other direction $s\to t$ from $s$ should be no different from the direction $s\to t$ on its own.
So
\[
    \text{tgt}(s,\id_s)=s \qqand \text{run}(s,\id_s,s\to t)=s\to t.
\]
In fact, $\id_s$ should really just be the direction $s\to s$.
Pictorially, we have the equation
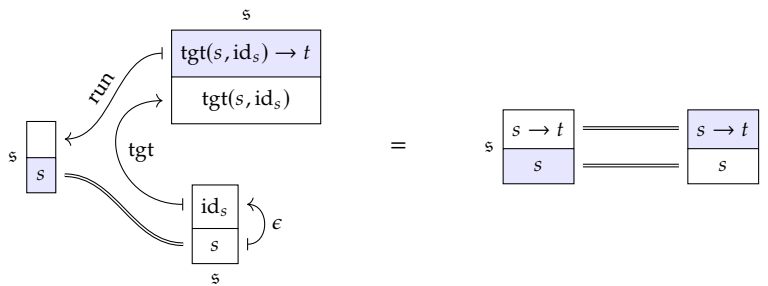
\[
\begin{tikzpicture}
	\node (1) {
        \begin{tikzpicture}[polybox, mapstos]
        	\node[poly, dom, "$\car{s}$" left] (r) {\at$s$};
        	\node[poly, right=2 of r.south, yshift=-2.5ex, "$\car{s}$" below] (p) {$\id_s$\at$s$};
        	\node[poly, cod, above=.8 of p, xshift=2.5ex, "$\car{s}$" above] (p') {tgt$(s,\id_s)\to t$\at tgt$(s,\id_s)$};

        	\draw[double, -] (r_pos) to[first] node[below] {} (p_pos);
        	\draw (p_pos) to[climb'] node[right] {$\epsilon$} (p_dir);
        	\draw (p_dir) to[climb] node[right] {tgt} (p'_pos);
        	\draw (p'_dir) to[last] node[above,sloped] {run} (r_dir);
        \end{tikzpicture}
	};
	\node[right=1.8 of 1] (2) {
        \begin{tikzpicture}[polybox, mapstos]
          	\node[poly, dom, "$\car{s}$" left] (c) {$s\to t$\at$s$};
          	\node[poly, cod, right=of c, "$\car{s}$" right] (c') {$s\to t$\at$s$};
          	\draw[double, -] (c_pos) -- (c'_pos);
          	\draw[double, -] (c'_dir) -- (c_dir);
	    \end{tikzpicture}
	};
	\node at ($(1.east)!.5!(2.west)$) {=};
\end{tikzpicture}
\]
Or, if you prefer, we might say that $\delta\then(\epsilon\tri\car{s})=\id_{\car{s}}$, or that the following diagram commutes:
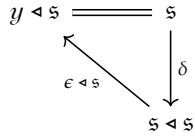
\[
\begin{tikzcd}[row sep=large]
    \yon\tri\car{s} & \car{s}\ar[d, "\delta"]\ar[l, equal] \\
    & \car{s}\tri\car{s}\ar[ul, "\epsilon\:\tri\:\car{s}"]
\end{tikzcd}
\]
This commutative diagram captures one way in which $\epsilon$ and $\delta$ always relate---and it's written entirely in the language of $\poly$, without having to talk about individual sets!

What about setting the second direction that $\delta$ runs together to what is specified by $\epsilon$, rather than the first?
To answer this, we should look at the composite lens $\delta\then(\car{s}\tri\epsilon)\colon\car{s}\to\car{s}\tri\yon\iso\car{s}$ instead.
But the do-nothing direction should still do nothing, so here's what the polybox picture should look like:
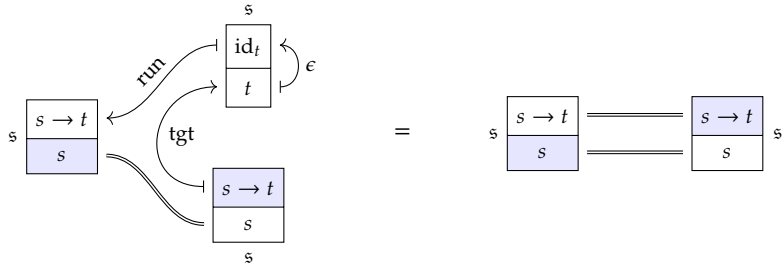
\[
\begin{tikzpicture}
	\node (1) {
        \begin{tikzpicture}[polybox, mapstos]
        	\node[poly, dom, "$\car{s}$" left] (r) {$s\to t$\at$s$};
        	\node[poly, cod, right=2 of r.south, yshift=-2.5ex, "$\car{s}$" below] (p) {$s\to t$\at$s$};
        	\node[poly, above=.8 of p, "$\car{s}$" above] (p') {$\id_t$\at $t$};

        	\draw[double, -] (r_pos) to[first] node[below] {} (p_pos);
        	\draw (p_dir) to[climb] node[right] {tgt} (p'_pos);	\draw (p'_pos) to[climb'] node[right] {$\epsilon$} (p'_dir);
        	\draw (p'_dir) to[last] node[above,sloped] {run} (r_dir);
        \end{tikzpicture}
	};
	\node[right=1.8 of 1] (2) {
        \begin{tikzpicture}[polybox, mapstos]
          	\node[poly, dom, "$\car{s}$" left] (c) {$s\to t$\at$s$};
          	\node[poly, cod, right=of c, "$\car{s}$" right] (c') {$s\to t$\at$s$};
          	\draw[double, -] (c_pos) -- (c'_pos);
          	\draw[double, -] (c'_dir) -- (c_dir);
	    \end{tikzpicture}
	};
	\node at ($(1.east)!.5!(2.west)$) {=};
\end{tikzpicture}
\]
The lens depicted on the right hand side of the equation is again the identity lens on $\car{s}$.

If we match up the two white boxes on the right hand side of the equation with the corresponding white boxes on the left, we can actually read two equations off of this polybox picture.
Matching up positions in the codomain tells us that the bottom arrow of $\delta$ on the left must send $s$ to itself: it is the identity function on $\car{s}(\1)$.
Indeed, this is exactly what we wanted to say about that arrow in \cref{subsec.comon.sharp.state.trans}.

Meanwhile, matching up directions in the domain tells us that
\[
    \text{run}(s,s\to t,\id_t)=s\to t,
\]
as we would expect: $\id_t$ is just be the direction $t\to t$.

More concisely, we can express both these facts in $\poly$ via the equation $\delta\then(\car{s}\tri\epsilon)=\id_{\car{s}}$.
The corresponding commutative diagram is as follows:
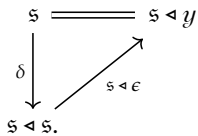
\[
\begin{tikzcd}[row sep=large]
    \car{s}\ar[d, "\delta"']\ar[r, equal] & \car{s}\tri\yon \\
    \car{s}\tri\car{s}.\ar[ur, "\car{s}\:\tri\:\epsilon"']
\end{tikzcd}
\]
We can combine this with our previous commutative diagram to say that the relationship between the do-nothing section $\epsilon\colon\car{s}\to\yon$ and the transition lens $\delta\colon\car{s}\to\car{s}\tri\car{s}$ of a state system $\car{s}$ is captured in $\poly$ by the following commutative diagram:
\begin{equation}\label{eqn.erasure_law_state}
\begin{tikzcd}[row sep=large]
	\yon\tri\car{s}&\car{s}\ar[d, "\delta" description]\ar[r, equal]\ar[l, equal]&\car{s}\tri\yon\\&
	\car{s}\tri\car{s}.\ar[ul, "\epsilon\:\tri\:\car{s}"]\ar[ur, "\car{s}\:\tri\:\epsilon"']
\end{tikzcd}
\end{equation}

\subsection{The transition lens is coassociative}\label{subsec.comon.sharp.state.coassoc}
\index{coassociativity}\index{lens!transition}\index{state system}

Toward the end of \cref{ex.dyn_sys_trans_polyboxes}, we noted that while the transition lens $\delta\colon\car{s}\to\car{s}\tri\car{s}$ gives us a canonical way to model two steps of a dynamical system with state system $\car{s}$, we have a choice of how to model three steps through the same system: we could obtain a lens $\car{s}\to\car{s}\tripow3$ that runs three directions together by taking either one of the composite lenses $\delta\then(\delta\tri\car{s})$ or $\delta\then(\car{s}\tri\delta)$.
That presents a problem for us: which one should we choose?

Happily, it turns out this choice is a false one.
If we write out the two composite lenses in polyboxes, with $\delta\then(\delta\tri\car{s})$ on the left and $\delta\then(\car{s}\tri\delta)$ on the right, we find that they are equal:
\begin{equation}\label{eqn.trans_lens_coassoc_polybox}
\scalebox{.875}{
\begin{tikzpicture}
    \node (p1) {
        \begin{tikzpicture}[polybox, mapstos, font=\tiny]
            \node[poly, dom, "$\car{s}$" left] (m') {$s_0\to s_3$\at$s_0$};
            \node[poly, right= of m'.south, yshift=-1ex, "$\car{s}$" below] (mm') {$s_0\to s_2$\at$s_0$};
            \node[poly, above=of mm', "$\car{s}$" above] (C') {$s_2\to s_3$\at$s_2$};
            \node[poly, cod, right= of mm'.south, yshift=-1ex, "$\car{s}$" right] (D') {$s_0\to s_1$\at$s_0$};
            \node[poly, cod, above=of D', "$\car{s}$" right] (mmm') {$s_1\to s_2$\at$s_1$};
            \node[poly, cod, above=of mmm', "$\car{s}$" right] (CC') {$s_2\to s_3$\at$s_2$};
            \draw[double, -] (m'_pos) to[first] (mm'_pos);
            \draw (mm'_dir) to[climb] node[right] {tgt} (C'_pos);
            \draw (C'_dir) to[last] node[above, sloped] {run} (m'_dir);
            \draw[double, -] (mm'_pos) to[first] (D'_pos);
            \draw (D'_dir) to[climb] node[right] {tgt} (mmm'_pos);
            \draw (mmm'_dir) to[last] node[above, sloped] {run} (mm'_dir);
            \draw[double, -] (C'_pos) to[first] (CC'_pos);
            \draw[double, -] (CC'_dir) to[last] (C'_dir);
        \end{tikzpicture}
	};
	\node (p2) [right=.3 of p1] {
	    \begin{tikzpicture}[polybox, mapstos, font=\tiny]
            \node[poly, dom, "$\car{s}$" left] (m) {$s_0\to s_3$\at$s_0$};
            \node[poly, right= of m.south, yshift=-1ex, "$\car{s}$" below] (D) {$s_0\to s_1$\at$s_0$};
            \node[poly, above=of D, "$\car{s}$" above] (mm) {$s_1\to s_3$\at$s_1$};
            \node[poly, cod, right= of D.south, yshift=-1ex, "$\car{s}$" right] (DD) {$s_0\to s_1$\at$s_0$};
            \node[poly, cod, above=of DD, "$\car{s}$" right] (mmm) {$s_1\to s_2$\at$s_1$};
            \node[poly, cod, above=of mmm, "$\car{s}$" right] (C) {$s_2\to s_3$\at$s_2$};
            \draw[double, -] (m_pos) to[first] (D_pos);
            \draw (D_dir) to[climb] node[right] {tgt} (mm_pos);
            \draw (mm_dir) to[last] node[above, sloped] {run} (m_dir);
            \draw[double, -] (D_pos) to[first] (DD_pos);
            \draw[double, -] (DD_dir) to[last] (D_dir);
            \draw[double, -] (mm_pos) to[first] (mmm_pos);
            \draw (mmm_dir) to[climb] node[right] {tgt} (C_pos);
            \draw (C_dir) to[last] node[above, sloped] {run} (mm_dir);
        \end{tikzpicture}
    };
	\node at ($(p1.south)!.5!(p2.north)$) {$=$};
\end{tikzpicture}
}
\end{equation}\index{polybox}
Remember: the way to read these polyboxes is to start at the lower blue square on the left and follow the path counter clockwise around the diagram; and if you reach a box with no arrows leading out of it, go up to the blue box above it and continue to follow the arrows from there.

There's a lot going on here, so let's break it down---we'll focus on the run functions first.
On the left hand side, we run together $s_0\to s_1$ and $s_1\to s_2$ to obtain $s_0\to s_2$, before running that together with $s_2\to s_3$ to obtain $s_0\to s_3$, as we see in the upper left box.
Meanwhile, on the right, we run together $s_1\to s_2$ and $s_2\to s_3$ to obtain $s_1\to s_3$, before running $s_0\to s_1$ together with our newly obtained $s_1\to s_3$ to again obtain $s_0\to s_3$ in the upper left box.
We could write this all out equationally, but all this is saying is that ``running together'' the directions of a state system is an associative operation.
When running together three directions, it doesn't matter whether we run the first two together or the last two together to start.
Not only is this guaranteed by the way in which we constructed $\delta$, it also makes intuitive sense.

\begin{exercise}
Using only the contents of the blue boxes and the target and run functions, write down the equation that we can read off of \eqref{eqn.trans_lens_coassoc_polybox} expressing the associativity of the ``running together'' operation.\index{associativity}
\begin{solution}
Following the arrows on either side of \eqref{eqn.trans_lens_coassoc_polybox} all the way to the domain's direction box, we obtain an expression for each box's contents that we can then set equal to each other.
The easiest way to actually write down these expressions is probably to start at the end with $s_0\to s_3$ and follow the arrows backward, unpacking each term until only the contents of the blue boxes remain (namely $s_0, s_0\to s_1, s_1\to s_2,$ and $s_2\to s_3$).
Here's what we get when we follow this process for the left hand side of \eqref{eqn.trans_lens_coassoc_polybox} (remember where to look for the three inputs to the run function):
\begin{align*}
    s_0\to s_3&=\text{run}(s_0,s_0\to s_2,s_2\to s_3)\\
    &=\text{run}(s_0,\text{run}(s_0,s_0\to s_1,s_1\to s_2),s_2\to s_3)\\
    &=\text{run}(s_0,\text{run}(s_0,s_0\to s_1,s_1\to s_2),s_2\to s_3);
\end{align*}
and here's what we get for the right (also remember where to look for the two inputs to the target function):
\begin{align*}
    s_0\to s_3&=\text{run}(s_0,s_0\to s_1,s_1\to s_3)\\
    &=\text{run}(s_0,s_0\to s_1,\text{run}(s_1,s_1\to s_2,s_2\to s_3))\\
    &=\text{run}(s_0,s_0\to s_1,\text{run}(\text{tgt}(s_0,s_0\to s_1),s_1\to s_2,s_2\to s_3)).
\end{align*}
Setting these equal yields our desired associativity equation:
\begin{align*}
    \text{run}(s_0,&\text{run}(s_0,s_0\to s_1,s_1\to s_2),s_2\to s_3)\\
    &=\text{run}(s_0,s_0\to s_1,\text{run}(\text{tgt}(s_0,s_0\to s_1),s_1\to s_2,s_2\to s_3)).
\end{align*}
\end{solution}
\end{exercise}
\index{associativity}

This associative property is what we get by matching up the white direction boxes on each domain side, but there are three more white position boxes on each codomain side that we can match up as well.
The fact that the lower two of these pairs coincide is a consequence of the fact that the bottom arrow of $\delta$ is the identity, which we already knew from \cref{subsec.comon.sharp.state.cohere}; so we don't learn anything new there.
On the other hand, the fact that both the upper position boxes in the codomain contain $s_2$ implies that
\begin{align*}
    \text{tgt}(s_0,\text{run}(s_0,s_0\to s_1,s_1\to s_2))&=\text{tgt}(s_0,s_0\to s_2)\\
    &=s_2\\
    &=\text{tgt}(s_1,s_1\to s_2)\\
    &=\text{tgt}(\text{tgt}(s_0,s_0\to s_1),s_1\to s_2),
\end{align*}
which is exactly what we wanted in \eqref{eqn.state_run_tgt}.
In English, this says that when we run together $s_0\to s_1$ and $s_1\to s_2$, the new direction's target is the same as the direction of $s_1\to s_2$, the latter of the two directions that we ran together.
Again, this coincides with our intuition: if we follow two directions in order, we should end up at wherever the latter direction leads us.

Hence both the associativity of running directions together and the relationship between the target and run functions from \eqref{eqn.state_run_tgt} are captured by the equality of lenses $\delta\then(\delta\tri\car{s})=\delta\then(\car{s}\tri\delta)$.
Equivalently, the following diagram in $\poly$ commutes:\index{associativity}
\begin{equation}\label{eqn.coassoc_law_states}
\begin{tikzcd}[row sep=large]
	\car{s}\ar[r, "\delta"]\ar[d, "\delta"']&
	\car{s}\tri\car{s}\ar[d, "\car{s}\:\tri\:\delta"]\\
	\car{s}\tri\car{s}\ar[r, "\delta\:\tri\:\car{s}"']&
	\car{s}\tri\car{s}\tri\car{s}.
\end{tikzcd}
\end{equation}
Another way to say this is that $\delta$ is \emph{coassociative}: while $\delta$ is only a lens $\car{s}\to\car{s}\tripow2$ as defined, the commutativity of \eqref{eqn.coassoc_law_states} tells us that the two ways of getting a lens $\car{s}\to\car{s}\tripow3$ out of $\delta$ are actually the same.
(This is dual to an \emph{associative} operation, which is a binary operation that gives rise to two identical ternary operations.)

So $\delta$ induces a canonical lens $\car{s}\to\car{s}\tripow3$, which we will call $\delta^{(3)}$, as it has $3$ copies of $\car{s}$ in its codomain.
Armed with this new lens, we can model three steps through a system $\varphi\colon\car{s}\to p$ with interface $p\in\poly$ as the composite lens
\[
    \car{s}\To{\delta^{(3)}}\car{s}\tripow3\To{\varphi\tripow3}p\tripow3.
\]

\index{coassociativity}\index{interface}

In fact, coassociativity guarantees that $\delta$ induces a canonical lens $\delta^{(n)}\colon\car{s}\to\car{s}\tripow{n}$ for every integer $n\geq2$, starting with $\delta^{(2)}\coloneqq\delta$.%
\footnote{Perhaps this notation seems a little unnatural, but it helps to think of the original $\delta\colon\car{s}\to\car{s}\tri\car{s}$ as the $n=2$ case of a generalized transition lens modeling $n$ steps through the state system.}
For concreteness, we could then define $\delta^{(n)}$ for $n>2$ inductively by $\delta^{(n)}\coloneqq\delta\then(\delta^{(n-1)}\tri\car{s})$, or just as well by $\delta^{(n)}\coloneqq\delta\then(\car{s}\tri\delta^{(n-1)})$ or even $\delta^{(n)}\coloneqq\delta\then(\delta^{(\ell)}\tri\delta^{(m)})$ for some pair of integers $\ell,m>1$ satisfying $\ell+m=n$.
Regardless, the coassociativity of $\delta$ means that it doesn't matter how we build a lens $\car{s}\to\car{s}\tripow{(n+1)}$ out of $\delta,\then,\tri,$ and identity lenses: we'll always end up with the same lens.
We will state this in more generality in \cref{prop.n_duplication}, but here's some practice with the $n=4$ case for a taste of what's to come.\index{coassociativity}

\begin{exercise}
\begin{enumerate}
    \item Say we know nothing about $\car{s}$ or $\delta$ apart from the fact that $\car{s}\in\poly$ and that $\delta$ is a lens $\car{s}\to\car{s}\tri\car{s}$.
    List all the ways to obtain a lens $\car{s}\to\car{s}\tripow4$ using only copies of $\delta,\id_\car{s}$, $\tri$, and $\then$.
    (You may write $\car{s}$ for $\id_\car{s}$.)

    \item Now assume that \eqref{eqn.coassoc_law_states} commutes.
    Show that all the lenses on your list are equal.
    (Hint: Use the fact that $(f\then g)\tri(h\then k)=(f\tri h)\then(g\tri k)$ for lenses $f,g,h,k$).
    \qedhere
\end{enumerate}
\begin{solution}
\begin{enumerate}
    \item Given $\car{s}\in\poly$ and a lens $\delta\colon\car{s}\to\car{s}\tri\car{s}$, we want all the ways to obtain a lens $\car{s}\to\car{s}\tripow4$ using $\delta,\id_\car{s}$ (i.e.\ $\car{s}$), $\tri,$ and $\then$.
    Starting with $\car{s}$, the only way to get to $\car{s}\tripow2$ is with a single $\delta\colon\car{s}\to\car{s}\tripow2$.
    From there, we can get to $\car{s}\tripow4$ directly by composing with $\delta\tri\delta$ to obtain $\delta\then(\delta\tri\delta)$.
    Alternatively, we can preserve either the first or the second $\car{s}$ using the identity, then get to $\car{s}\tripow3$ from the other $\car{s}$ in one of two ways: either $\delta\then(\delta\tri\car{s})$ or $\delta\then(\car{s}\tri\delta)$.
    This gives us $4$ more ways to write a lens $\car{s}\to\car{s}\tripow4$: either $\delta\then(\car{s}\tri(\delta\then(\delta\tri\car{s})))$ or $\delta\then(\car{s}\tri(\delta\then(\car{s}\tri\delta)))$ if we chose to preserve the first $\car{s}$, and either $\delta\then((\delta\then(\delta\tri\car{s}))\tri\car{s})$ or $\delta\then((\delta\then(\car{s}\tri\delta))\tri\car{s})$ if we chose to preserve the second.
    Here's the full list, sorted roughly by how far to the left we try to apply each $\delta$:
    \begin{enumerate}[label=(\arabic*)]
        \item $\delta\then((\delta\then(\delta\tri\car{s}))\tri\car{s})$
        \item $\delta\then((\delta\then(\car{s}\tri\delta))\tri\car{s})$
        \item $\delta\then(\delta\tri\delta)$
        \item $\delta\then(\car{s}\tri(\delta\then(\delta\tri\car{s})))$
        \item $\delta\then(\car{s}\tri(\delta\then(\car{s}\tri\delta)))$
    \end{enumerate}
    This coincides with the $5$ different ways to parenthesize a $4$-term expression.

    \item We wish to show that if \eqref{eqn.coassoc_law_states} commutes, then all the lenses on our list are equal.
    The commutativity of \eqref{eqn.coassoc_law_states} implies that $\delta\then(\delta\tri\car{s})=\delta\then(\car{s}\tri\delta)$; so (1) and (2) from our list are equal, as are (4) and (5).
    Meanwhile, since $\car{s}=\car{s}\then\car{s}$, we can rewrite (1) as
    \[
        \delta\then((\delta\then(\delta\tri\car{s}))\tri(\car{s}\then\car{s}))=\delta\then(\delta\tri\car{s})\then(\delta\tri\car{s}\tri\car{s}),
    \]
    where the associativity of $\tri$ and $\then$ allows us to drop some parentheses.\index{associativity}
    Then the commutativity of \eqref{eqn.coassoc_law_states} allows us to further rewrite this as
    \begin{align*}
        \delta\then(\car{s}\tri\delta)\then(\delta\tri\car{s}\tri\car{s})&=\delta\then((\car{s}\then\delta)\tri(\delta\then(\car{s}\tri\car{s})))\\
        &=\delta\then(\delta\tri\delta),
    \end{align*}
    so (1) and (3) are equal.
    Similarly, we can rewrite (5) as
    \begin{align*}
        \delta\then((\car{s}\then\car{s})\tri(\delta\then(\car{s}\tri\delta)))&=\delta\then(\car{s}\tri\delta)\then(\car{s}\tri\car{s}\tri\delta)\\
        &=\delta\then(\delta\tri\car{s})\then(\car{s}\tri\car{s}\tri\delta)\\
        &=\delta\then((\delta\then(\car{s}\tri\car{s}))\tri(\car{s}\then\delta))\\
        &=\delta\then(\delta\tri\delta)
    \end{align*}
    so (3) and (5) are equal.
    Hence all the lenses on our list are equal.
\end{enumerate}
\end{solution}
\end{exercise}

\subsection{Running dynamical systems}\label{subsec.comon.sharp.state.run}

\index{dynamical system!running}\index{interface}

Finally, we are ready to fulfill our promise from way back in \cref{ex.do_nothing} by using the language of $\poly$ to describe stepping through a dynamical system $n$ times for arbitrary $n\in\nn$.
Given a dynamical system $\varphi\colon\car{s}\to p$ with interface $p\in\poly$, we can construct a new dynamical system that we call $\text{Run}_n(\varphi)$, with the same state system $\car{s}$ but a new interface $p\tripow{n}$, by defining $\text{Run}_n(\varphi)\coloneqq\delta^{(n)}\then\varphi\tripow{n}$.
Visually, we define $\text{Run}_n(\varphi)$ so that the following diagram commutes:
\begin{equation*}%\label{eqn.run}
\begin{tikzcd}
	\car{s}\ar[r, "\delta^{(n)}"]\ar[rr, bend right, "{\text{Run}_n(\varphi)}"']&
	\car{s}\tripow{n}\ar[r, "\varphi\tripow{n}"]&
	p\tripow n
\end{tikzcd}
\end{equation*}
One way to think of this is that $\text{Run}_n(\varphi)$ is a sped-up version of $\varphi$: one step through $\text{Run}_n(\varphi)$ is equivalent to $n$ steps through $\varphi$.\index{interaction}
But this is just because a single interaction with the interface $p\tripow{n}$ models a sequence of $n$ interactions with the interface $p$, as detailed in \cref{subsec.comon.comp.def.dyn_sys,ex.dyn_sys_comp_polyboxes}.
So $\text{Run}_n(\varphi)$ repackages $n$ cycles through $\varphi$ into a single step.
Crucially, $\delta^{(n)}$ is what tells us how to sequence all $n$ of these steps together on the state system side.
We illustrated how $\delta$ does this for the $n=2$ case in \cref{ex.dyn_sys_trans_polyboxes}, and here's a polybox picture for the $n=3$ case:
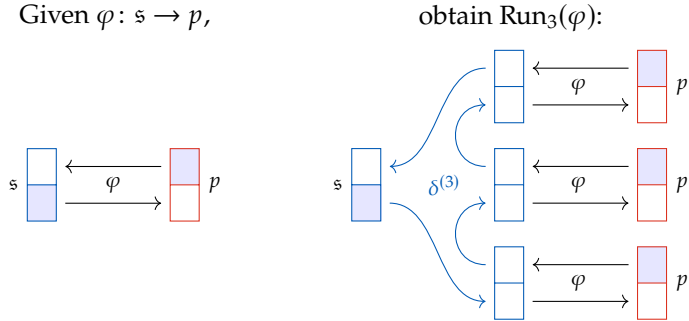
\[
\begin{tikzpicture}
	\node (given) {
	\begin{tikzpicture}[polybox, tos]
		\node[poly, dom, my-blue, "$\car{s}$" left] (S) {};
		\node[poly, cod, my-red, right=of S, "$p$" right] (p) {};
		\draw (S_pos) to[first] (p_pos);
		\draw (p_dir) to[last]  (S_dir);
		\node at ($(S.east)!.5!(p.west)$) {$\varphi$};
	\end{tikzpicture}
	};
	\node[right=of given] (obtain) {
	\begin{tikzpicture}[polybox, tos]
		\node[poly, dom, my-blue, "$\car{s}$" left] (S) {};
		\node[poly, my-blue, right=of S] (S2) {};
		\node[poly, my-blue, below=of S2] (S1) {};
		\node[poly, my-blue, above=of S2] (S3) {};
		\node[poly, my-red, cod, right=of S1, "$p$" right] (p1) {};
		\node[poly, my-red, cod, right=of S2, "$p$" right] (p2) {};
		\node[poly, my-red, cod, right=of S3, "$p$" right] (p3) {};
		\draw (S1_pos) to[first] (p1_pos);
		\draw (p1_dir) to[last] (S1_dir);
		\draw (S2_pos) to[first] (p2_pos);
		\draw (p2_dir) to[last]  (S2_dir);
		\draw (S3_pos) to[first] (p3_pos);
		\draw (p3_dir) to[last]  (S3_dir);
		\draw[my-blue] (S_pos) to[first] (S1_pos);
		\draw[my-blue] (S1_dir) to[climb] (S2_pos);
		\draw[my-blue] (S2_dir) to[climb] (S3_pos);
		\draw[my-blue] (S3_dir) to[last] (S_dir);
		\node[my-blue] at ($(S.east)!.33!(S2.west)$) {$\delta^{(3)}$};
		\node at ($(S1.east)!.33!(p1.west)$) {$\varphi$};
		\node at ($(S2.east)!.33!(p2.west)$) {$\varphi$};
		\node at ($(S3.east)!.33!(p3.west)$) {$\varphi$};
  \end{tikzpicture}
	};
	\node[above] at (obtain.north) (obtain_lab) {obtain $\text{Run}_3(\varphi)$:};
	\node at (given|-obtain_lab) {Given $\varphi\colon \car{s}\to p$,};
\end{tikzpicture}
\]
Notice that we have only defined $\delta^{(n)}$, and thus $\text{Run}_n(\varphi)$, for integers $n\geq2$.
But $n=0$ runs through $n$ is doing nothing, modeled by the do-nothing section $\epsilon\colon\car{s}\to\yon$, while $n=1$ run through $\varphi$ is modeled by $\varphi\colon\car{s}\to\car{s}$ itself.
So we want $\text{Run}_0(\varphi)=\epsilon$ and $\text{Run}_1(\varphi)=\varphi$; we can achieve this by setting $\delta^{(0)}\coloneqq\epsilon$ and $\delta^{(1)}\coloneqq\id_\car{s}$.
Here we should think of the do-nothing section $\delta^{(0)}$ as the transition lens modeling $0$ steps through our state system, and the identity $\delta^{(1)}$ as the transition lens modeling a single step.

\begin{exercise}
Verify that when $\delta^{(0)}=\epsilon$ and $\delta^{(1)}=\id_\car{s}$, if $\text{Run}_n(\varphi)$ is defined as $\delta^{(n)}\then\varphi\tripow{n}$ for all $n\in\nn$, then $\text{Run}_0(\varphi)=\epsilon$ and $\text{Run}_1(\varphi)=\id_\car{s}$.
\begin{solution}
We have $\text{Run}_n(\varphi)=\delta^{(n)}\then\varphi\tripow{n}$ for all $n\in\nn$, as well as $\delta^{(0)}=\epsilon$ and $\delta^{(1)}=\id_\car{s}$.
Then $\text{Run}_0(\varphi)=\epsilon\then\varphi\tripow0=\epsilon\then\id_\yon=\epsilon$ and $\text{Run}_1(\varphi)=\id_\car{s}\then\varphi\tripow1=\varphi\tripow1=\varphi$.
\end{solution}
\end{exercise}

\begin{example}[Returning every other position] \label{ex.long_div_skip}
In \cref{exc.long_div}, we built a dynamical system $\varphi\colon S\yon^S\to\nn\yon$ that returns natural numbers---specifically digits, alternating between $0$'s and the base-$10$ digits of $1/7$ after the decimal point like so:
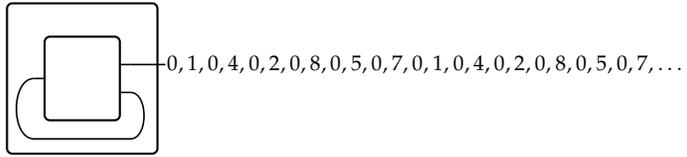
\[
\begin{tikzpicture}[oriented WD]
	\node[bb={1}{2}] (inner) {};
	\node[bb={0}{0}, inner xsep=1cm, inner ysep=1cm] (outer) {};
	\coordinate (outer_out1) at (outer.east|-inner_out1);
	\draw[shorten >=-3pt] (inner_out1) -- (outer_out1);
	\draw
		let \p1=(inner.south east), \p2=(inner.south west), \n1=\bbportlen, \n2=\bby in
		(inner_out2) to[in=0] (\x1+\n1,\y1-\n2) -- (\x2-\n1,\y1-\n2) to[out=180] (inner_in1);
		\node[right, font=\footnotesize] at (outer_out1) {$0,1,0,4,0,2,0,8,0,5,0,7,0,1,0,4,0,2,0,8,0,5,0,7,\ldots$};
\end{tikzpicture}
\]
Say we only wanted the system to return the digits of $1/7$ after the decimal point; we'd like to do away with all these $0$'s.
In other words, we want a new system $S\yon^S\to\nn\yon$ that acts like $\varphi$, except that it only returns every other position that $\varphi$ returns.

We could build such a system from scratch---or we can simply start from $\varphi$ and apply $\text{Run}_2$, yielding a system $\text{Run}_2(\varphi)\colon S\yon^S\to\nn\yon\tri\nn\yon\iso\nn^\2\yon$ that returns the positions of $\varphi$ two at a time:
\begin{align*}
    &(0,1),(0,4),(0,2),(0,8),(0,5),(0,7),\\
    &(0,1),(0,4),(0,2),(0,8),(0,5),(0,7),\ldots
\end{align*}
Then we just need to compose $\text{Run}_2(\varphi)$ with a lens $\pi_2\colon\nn^\2\yon\to\nn\yon$ equal to the second coordinate projection on positions (and the identity on directions) to extract the positions we want.
The new system $S\yon^S\to\nn\yon$ that skips over every position of $\varphi$ is therefore the following composite:
\[
\begin{tikzcd}
	S\yon^S\ar[r, "\delta^{(2)}"]\ar[rr, bend right, "{\text{Run}_2(\varphi)}"']&
    S\yon^S\tri S\yon^S\ar[r, "\varphi\tripow2"]&
	\nn\yon\tri\nn\yon\iso\nn^\2\yon\ar[r,"\pi_2"]&
	\nn\yon.
\end{tikzcd}
\]
We can apply this technique in general to skip (or otherwise act on) the positions of a dynamical system at regular intervals.
\end{example}

One drawback of the $\text{Run}_n(-)$ operation is that we need to keep track of a separate morphism $S\yon^S\to p\tripow{n}$ for every $n\in\nn$, as well as various ways to relate these morphisms for different values of $n$.
Is there a way to package all this information into a single morphism that can model arbitrarily long runs through the system?
We will answer this question in \cref{ch.comon.cofree}; but for now, let us investigate what's really going on with our state systems algebraically.

\subsection{State systems as comonoids}
\index{state system!comonoid as}

It turns out that objects equipped with morphisms like those in \cref{subsec.comon.sharp.state.nothing,subsec.comon.sharp.state.trans} that satisfy the commutative diagrams from \cref{subsec.comon.sharp.state.cohere,subsec.comon.sharp.state.coassoc} are well-known to category theorists.

\index{comonoid!definition}

\begin{definition}[Comonoid]\label{def.comonoid}
In a monoidal category $(\Cat{C},\yon,\tri)$, a \emph{comonoid} $\com{C}\coloneqq(\car{c},\epsilon,\delta)$ consists of
\begin{itemize}
    \item an object $\car{c}\in\Cat{C}$, called the \emph{carrier};
    \item a morphism $\epsilon\colon\car{c}\to\yon$ in $\Cat{C}$, called the \emph{eraser} (or the \emph{counit}); and
    \item a morphism $\delta\colon\car{c}\to\car{c}\tri\car{c}$ in $\Cat{C}$, called the \emph{duplicator} (or the \emph{comultiplication});
\end{itemize}
such that the following diagrams, collectively known as the \emph{comonoid laws}, commute:
\begin{equation}\label{eqn.erasure_law}
\begin{tikzcd}[background color=definitioncolor, row sep=large]
	\yon\tri \car{c}&\car{c}\ar[d, "\delta" description]\ar[r, equal]\ar[l, equal]&\car{c}\tri\yon\\&
	\car{c}\tri\car{c},\ar[ul, "\epsilon\:\tri\:\car{c}"]\ar[ur, "\car{c}\:\tri\:\epsilon"']
\end{tikzcd}
\end{equation}
where the left triangle is known as the \emph{left erasure} (or \emph{counit}) \emph{law} and the right triangle is known as the \emph{right erasure} (or \emph{counit}) \emph{law}; and
\begin{equation}\label{eqn.coassoc_law}
\begin{tikzcd}[row sep=large]
	\car{c}\ar[r, "\delta"]\ar[d, "\delta"']&
	\car{c}\tri\car{c}\ar[d, "\car{c}\:\tri\:\delta"]\\
	\car{c}\tri \car{c}\ar[r, "\delta\:\tri\:\car{c}"']&
	\car{c}\tri\car{c}\tri\car{c},
\end{tikzcd}
\end{equation}
known as the \emph{coassociative law}.\index{coassociativity}

We may also say that the eraser and duplicator morphisms comprise a \emph{comonoid structure} on the carrier, or we may identify a comonoid with its carrier if the eraser and duplicator can be inferred from context.

We refer to a comonoid $\com{C}$ in $(\poly,\yon,\tri)$ as a \emph{polynomial comonoid}.
\end{definition}
\index{eraser|see{polynomial comonoid, counit of}}
\index{duplicator|see{polynomial comonoid, comultiplication of}}
\index{polynomial comonoid!counit of}
\index{polynomial comonoid!comultiplication of}

\begin{remark}
The concept of a \emph{comonoid} in a monoidal category is dual to that of a \emph{monoid}, which may be more familiar.
Monoids come with \emph{unit} and \emph{multiplication} morphisms that point the other way, so named because they generalize the unit and multiplication operations of a monoid in $\smset$.
(We'll talk more about monoids in $\smset$ in \cref{ex.monoids}.)
Prepending `co-' to each term yields the corresponding terms for comonoids.

\index{monoid!unit of}\index{monoid!multiplication of}

The alternative names \emph{eraser} for the \emph{counit} and \emph{duplicator} for the \emph{comultiplication} are less standard, but we will favor them to avoid confusion between the counit of a \emph{comonoid} and the counit of an \emph{adjunction}---and so that their names match up with the Greek letters $\epsilon$ and $\delta$ that we will so often use to label them.
The word ``duplicator'' comes from the fact that $\delta\colon\car{c}\to\car{c}\tri\car{c}$ effectively turns one $\car{c}$ into two, while the ``eraser'' $\epsilon\colon\car{c}\to\yon$ erases the $\car{c}$ altogether, leaving only the monoidal unit $\yon$.
Still, it can be helpful to think of comonoids as having a \emph{coassociative} comultiplication along with a counit satisfying \emph{left and right counit laws}.
\end{remark}\index{coassociativity}

\index{comonad|see{polynomial comonoid}}

\begin{remark}
Comonoids in a functor category with respect to the composition product are generally known as \emph{comonads}.
So it would be a little more precise and familiar to refer to our polynomial comonoids as \emph{polynomial comonads}.
But since we think of our polynomials more often in terms of positions and directions than as functors, we’ll favor the term comonoid over comonad.
\end{remark}

\index{state system!as polynomial comonoid}
\begin{example}[State systems are polynomial comonoids]
Nearly all our work on state systems up until now can be summarized thusly:
\slogan{
    every state system is a polynomial comonoid,\\
    whose eraser is the do-nothing section\\
    and whose duplicator is the transition lens.
}
The comonoid structure on a state system $\car{s}$ is what allows us to write canonical lenses $\car{s}\to\car{s}\tripow{n}$ for any $n\in\nn$.
We can then model $n$ steps through a dynamical system $\varphi\colon\car{s}\to p$ with interface $p\in\poly$ by composing this canonical lens with $\varphi\tripow{n}$ to obtain a ``sped-up'' dynamical system $\text{Run}_n(\varphi)$.
This new system has the same state system $\car{s}$, but its interface is now $p\tripow{n}$.

The canonicity of $\car{s}\to\car{s}\tripow{n}$ is due to the following standard result about comonoids, which can be proved inductively.
\end{example}

\index{duplication!multifold}

\begin{proposition}[Defining $\delta^{(n)}$] \label{prop.n_duplication}
Given a comonoid $(\car{c},\epsilon,\delta)$, let $\delta^{(n)}\colon\car{c}\to\car{c}\tripow{n}$ be given as follows. Let $\delta^{(0)}\coloneqq\epsilon$ and inductively define $\delta^{(n+1)}\coloneqq\delta\then\left(\delta^{(n)}\tri\car{c}\right)$ for all $n\in\nn$.
Then we have the following:
\begin{enumerate}[label=(\alph*)]
    \item $\delta^{(n)}$ is a morphism $\car{c}\to\car{c}\tripow{n}$ for all $n\in\nn$;
    \item $\delta^{(1)}=\car{c}=\id_\car{c}$;
    \item $\delta^{(2)}=\delta$; and
    \item $\delta^{(n)}=\delta\then\left(\delta^{(k)}\tri\delta^{(n-k)}\right)$ for all $k,n\in\nn$ with $k\leq n$, so our choice of morphism $\car{c}\to\car{c}\tripow{(n+1)}$ is canonical. % More generally: $\delta^{(n)}=\delta^{(m)}\then\left(\delta^{(k_1)}\tri\cdots\tri\delta^{(k_m)}\right)$ for $m,n\in\nn$ and each $k_i\in\nn$ such that $k_1+\cdots+k_m=n$.
\end{enumerate}
\end{proposition}
 \begin{proof}
 We leave parts (a), (b), and (c) for \cref{exc.n_duplication}. Part (d) amounts to coassociativity.
\end{proof}\index{coassociativity}

We'll continue to use the notation introduced here throughout for general comonoids.

\begin{exercise} \label{exc.n_duplication}
Prove the first three parts of \cref{prop.n_duplication}. %Finish the proof of \cref{prop.n_duplication} as follows.
\begin{enumerate}
    \item Prove part (a).
    \item Prove part (b).
    \item Prove part (c).\qedhere
\end{enumerate}
\begin{solution}
We complete the proof of \cref{prop.n_duplication}, where we are given a comonoid $(\car{c},\epsilon,\delta)$ along with $\delta^{(0)}\coloneqq\epsilon$ and $\delta^{(n+1)}\coloneqq\delta\then\left(\delta^{(n)}\tri\car{c}\right)$ for all $n\in\nn$.
\begin{enumerate}
    \item We will show that $\delta^{(n)}$ is a morphism $\car{c}\to\car{c}\tripow{n}$ for every $n\in\nn$ by induction on $n$.
    We know $\delta^{(0)}=\epsilon$ is a morphism $\car{c}\to\yon=\car{c}\tripow0$, and for each $n\in\nn$, if $\delta^{(n)}$ is a morphism $\car{c}\to\car{c}\tripow{n}$, then the composite $\delta^{(n+1)}=\delta\then\left(\delta^{(n)}\tri\car{c}\right)$ is a morphism
    \[
        \car{c}\To\delta\car{c}\tri\car{c}\To{\delta^{(n)}\:\tri\:\car{c}}\car{c}\tripow{n}\tri\car{c}\iso\car{c}\tripow{(n+1)}
    \]
    Hence the result follows by induction.

    \item We have $\delta^{(1)}=\delta\then(\delta^{(0)}\tri\car{c})=\delta\then(\epsilon\tri\car{c})=\id_\car{c}$ by the left erasure law from \eqref{eqn.erasure_law}.

    \item By the previous part, we have $\delta^{(2)}=\delta\then(\delta^{(1)}\tri\car{c})=\delta\then(\car{c}\tri\car{c})=\delta$.
\end{enumerate}
\end{solution}
\end{exercise}

\index{polynomial comonoid}

\begin{example}[Not all polynomial comonoids are state systems] \label{ex.not_all_com_state}
At this point, a natural question to ask is whether everything we know about a state system $\car{s}$ is captured by the fact that state systems are polynomial comonoids.
In other words, are state systems the only polynomial comonoids there are?

The answer turns out to be no.
After all, there is one fact about state systems from \cref{subsec.comon.sharp.state.trans} that we did not encode in $\poly$: for a fixed state $s\in\car{s}(\1)$, the target function $\car{s}[s]\to\car{s}(\1)$ sending directions at $s$ to their target states is a bijection.

Nothing in our comonoid laws guarantees this bijectivity.
An arbitrary polynomial comonoid might send different directions at $s$ to the same target---given a second state $t$, there may be multiple ways to get from $s$ to $t$.
It might even send \emph{no} directions at $s$ to a target $t$, making it impossible to get from $s$ to $t$.
(We'll give an explicit example of a comonoid that is not a state system in \cref{ex.walking_arrow_com}.)
State systems as we have defined them are just the polynomial comonoids that do not allow either of these variations, for which the bijective property holds.

We consider this a feature, not a bug.
After all, it is an abstraction to say that there is exactly one way to get from any one state in a system to another.
It is perfectly plausible that the inner workings of a state system do not permit traveling between some states and differentiate ways of traveling between others.
We won't formally introduce this idea into our theory of dependent dynamical systems,
but we will often think of polynomial comonoids as a sort of generalized state system throughout the rest of the book.
\end{example}

\begin{example}[A comonoid that is not a state system]\label{ex.walking_arrow_com}
The polynomial $\yon^\2+\yon$ is not a state system: one of its direction-sets has one fewer element than its position-set.
But it can still be given a comonoid structure.
We describe that structure here, but we will go a little quickly, because we'll soon discover a much more familiar way to think about comonoids.

Define $\car{a}\coloneqq\{s\}\yon^{\{\id_s,\,a\}}+\{t\}\yon^{\{\id_t\}}\iso\yon^\2+\yon$.
Here is its tree picture:
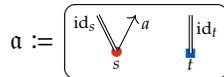
\[
\begin{tikzpicture}[rounded corners]
	\node (p1) [draw, "$\car{a}\coloneqq$" left] {
	\begin{tikzpicture}[trees, sibling distance=5mm]
    \node["\tiny $s$" below, my-red] (1) {$\bullet$}
      child  {coordinate (is) \idchild}
      child {coordinate (a)};
    \node[right=.8 of 1,"\tiny $t$" below, my-blue] (2) {\tiny$\blacksquare$}
      child  {coordinate (it) \idchild};
    \node[below left=0 of is, font=\tiny] {$\id_s$};
    \node[below right=.1 of it, font=\tiny] {$\id_t$};
    \node[below right=.1 of a, font=\tiny] {$a$};
  \end{tikzpicture}
  };
\end{tikzpicture}
\]
Notice that we have drawn one direction out of each position---$\id_s$ and $\id_t$---with a double bar.
We let these be the directions that the eraser $\epsilon\colon\car{a}\to\yon$ picks out.
The double bar is meant to evoke an equals sign from the root position to the eventual target position, which is appropriate, as these two positions should be equal for every direction that the eraser selects.
We can draw the selections that $\epsilon$ makes like so:
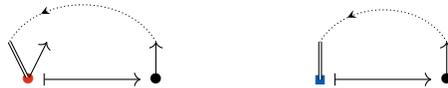
\[
\begin{tikzpicture}[trees, bend right=60]
  \node[my-red] (1) {$\bullet$}
  	child  {coordinate (11) \idchild}
    child {coordinate (12)};
  \node[right=1.5 of 1] (1y) {$\bullet$}
  	child {coordinate (1y1)};
  \node[right=2 of 1y, my-blue] (2) {\tiny$\blacksquare$}
  	child  {coordinate (21) \idchild};
  \node[right=1.5 of 2] (2y) {$\bullet$}
  	child {coordinate (2y1)};
	\draw[|->, shorten <= 3pt, shorten >= 3pt] (1) -- (1y);
	\draw[|->, shorten <= 3pt, shorten >= 3pt] (2) -- (2y);
	\draw[densely dotted, postaction={decorate}] (1y1) to (11);
	\draw[densely dotted, postaction={decorate}] (2y1) to (21);
\end{tikzpicture}
\]

Now we need a duplicator $\delta\colon\car{a}\to\car{a}\tri\car{a}$.
Before we define it, let us draw out $\car{a}\tri\car{a}$ to see what it looks like.
Remember that we need to graft corollas of $\car{a}$ onto leaves of $\car{a}$ in every possible way:
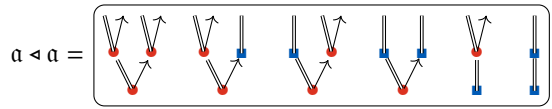
\[
\begin{tikzpicture}[rounded corners]
	\node (p1) [draw, "$\car{a}\tri\car{a}=$" left] {
	\begin{tikzpicture}[trees,
	  level 1/.style={sibling distance=5mm},
  	level 2/.style={sibling distance=2.5mm}]
    \node[my-red] (1) {$\bullet$}
      child  {
        node [my-red] {$\bullet$}
 		    child  {\idchild}
      	child {}
			\idchild
			}
      child  {
        node [my-red] {$\bullet$}
 		    child  {\idchild}
      	child {}
			};
    \node[right=1 of 1, my-red] (2) {$\bullet$}
      child  {
        node [my-red]{$\bullet$}
 		    child  {\idchild}
      	child {}
			\idchild
			}
      child {node [my-blue] {\tiny$\blacksquare$}
      	child  {\idchild}
			};
    \node[right=1 of 2, my-red] (3) {$\bullet$}
      child {node [my-blue] {\tiny$\blacksquare$}
      	child  {\idchild}
				\idchild
			}
      child  {
        node [my-red] {$\bullet$}
 		    child {\idchild}
      	child {}
			};
    \node[right=1 of 3, my-red] (4) {$\bullet$}
      child {node [my-blue] {\tiny$\blacksquare$}
      	child  {\idchild}
			\idchild
			}
      child {node [my-blue] {\tiny$\blacksquare$}
      	child  {\idchild}
			};
    \node[right=.8 of 4, my-blue] (5) {\tiny$\blacksquare$}
      child  {
        node [my-red] {$\bullet$}
 		    child  {\idchild}
      	child {}
			\idchild
			};
    \node[right=.6 of 5, my-blue] (6) {\tiny$\blacksquare$}
      child {node [my-blue] {\tiny$\blacksquare$}
      	child  {\idchild}
			\idchild
			};
  \end{tikzpicture}
  };
\end{tikzpicture}
\]
Each of these trees gives a way to match directions out of one position to positions they could lead to.
On positions, $\delta$ will decide which matchings to pick by sending the red $s$ to one of the four positions on the left and the blue $t$ to one of the two positions on the right.
We want the double-barred directions that the eraser picked out to have the same position on either end (in fact, the erasure laws guarantee this).
So the only choice to be made is whether we want the other direction $a$ at $s$ to point to $s$ or to $t$.
Let us pick $t$ for the time being, so that on positions, $\delta$ looks like this:
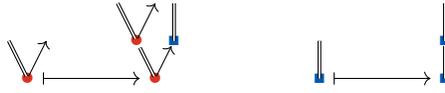
\[
\begin{tikzpicture}[trees, sibling distance=5mm,	bend right=60]
	\node (1A) [my-red] {$\bullet$}
  	child  {coordinate (1A1) \idchild}
    child {coordinate (1A2)};
  \node (2A) [right=1.5 of 1A, my-red] {$\bullet$}
      child  {
        node [my-red] {$\bullet$}
 		    child  {coordinate (2A1) \idchild}
      	child {coordinate (2A2)}
			\idchild
			}
      child {node [my-blue] {\tiny$\blacksquare$}
      	child  {coordinate (2A3) \idchild}
			};
	\draw[|->, shorten <= 3pt, shorten >= 3pt] (1A) -- (2A);
% 	\draw[densely dotted, postaction={decorate}] (2A1) to (1A1);
% 	\draw[densely dotted, postaction={decorate}] (2A2) to (1A2);
% 	\draw[densely dotted, postaction={decorate}] (2A3) to (1A2);
%
  \node[right=2 of 2A, my-blue] (1B) {\tiny$\blacksquare$}
  	child  {coordinate (1B1) \idchild};
  \node[right=1.5 of 1B, my-blue] (2B) {\tiny$\blacksquare$}
  	child {node [my-blue] {\tiny$\blacksquare$}
    child  {coordinate (2B1) \idchild}
		\idchild
	};
	\draw[|->, shorten <= 3pt, shorten >= 3pt] (1B) -- (2B);
% 	\draw[densely dotted, postaction={decorate}] (2B1) to (1B1);
\end{tikzpicture}
\]
As in \cref{ex.trans_trees}, we can interpret this as telling us how to ``bend'' the arrows of $\car{a}$ so that they point to other positions:
%\begin{equation} \label{eqn.walking_arrow_bent_cor}
%\begin{tikzcd}
%  {\huge\color{my-red}\bullet} \ar[rr] \ar[equals, loop left] && {\large\color{my-blue}\blacksquare} \ar[equals, loop right]
%\end{tikzcd}
%\end{equation}

\begin{equation} \label{eqn.walking_arrow_bent_cor}
\begin{tikzpicture}[baseline=(b)]
    \node[circle,minimum size=2cm] (b) {};

    \node[minimum size=0.45cm,draw,circle,my-red,fill=my-red] (2-1) at (b.{360/2*1}){};
    \draw[double, -] (2-1) to [in=360/2*1-30,out=360/2*1+30,loop] ();

    \node[minimum size=0.4cm,draw,my-blue,fill=my-blue] (2-2) at (b.{360/2*2}){};
    \draw[double, -] (2-2) to [in=360/2*2-30,out=360/2*2+30,loop] ();
%
%    \foreach\x/\c in {1/red, 2/blue} {
%        \node[minimum size=0.1cm,draw,circle,\c,fill=\c] (2-\x) at (b.{360/2*\x}){};
%        \draw[double, -] (2-\x) to [in=360/2*\x-30,out=360/2*\x+30,loop] ();
%        \relax
%    }
    \draw (2-1) edge[->] (2-2);
\end{tikzpicture}
\end{equation}

Meanwhile, on directions, $\delta$ should tell us how to run two directions together into one.
Fortunately, there is not much for us to do here---we know that if one of the two directions $\delta$ runs together is one of the double-barred directions that the eraser picked out, then $\delta$ should ignore that ``do-nothing'' direction and yield the other direction (again, the erasure laws ensure this).
Here is what that looks like:
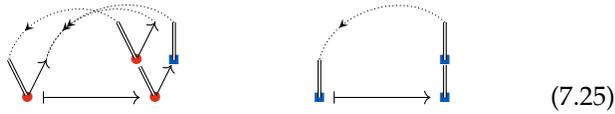
\begin{equation}\label{eqn.my_comonoid_delta}
\begin{tikzpicture}[trees, sibling distance=5mm,	bend right=60]
	\node (1A) [my-red] {$\bullet$}
  	child  {coordinate (1A1) \idchild}
    child {coordinate (1A2)};
  \node (2A) [right=1.5 of 1A, my-red] {$\bullet$}
      child  {
        node [my-red] {$\bullet$}
 		    child  {coordinate (2A1) \idchild}
      	child {coordinate (2A2)}
			\idchild
			}
      child {node [my-blue] {\tiny$\blacksquare$}
      	child  {coordinate (2A3) \idchild}
			};
	\draw[|->, shorten <= 3pt, shorten >= 3pt] (1A) -- (2A);
	\draw[densely dotted, postaction={decorate}] (2A1) to (1A1);
	\draw[densely dotted, postaction={decorate}] (2A2) to (1A2);
	\draw[densely dotted, postaction={decorate}] (2A3) to (1A2);
  \node[right=2 of 2A, my-blue] (1B) {\tiny$\blacksquare$}
  	child  {coordinate (1B1) \idchild};
  \node[right=1.5 of 1B, my-blue] (2B) {\tiny$\blacksquare$}
  	child {node [my-blue] {\tiny$\blacksquare$}
    child  {coordinate (2B1) \idchild}
		\idchild
	};
	\draw[|->, shorten <= 3pt, shorten >= 3pt] (1B) -- (2B);
	\draw[densely dotted, postaction={decorate}] (2B1) to (1B1);
\end{tikzpicture}
\end{equation}
That is all we need to specify the triple $(\car{a},\epsilon,\delta)$.

Here are $\epsilon\colon\car{a}\to\yon$ and $\delta\colon\car{a}\to\car{a}\tri\car{a}$ again, in terms of polyboxes.\index{polybox}
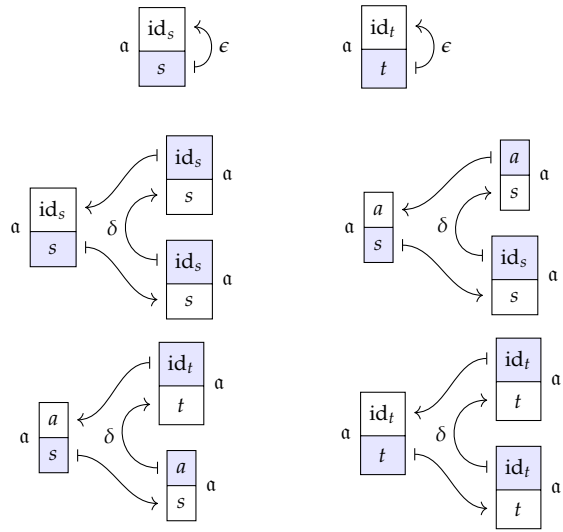
\[
\begin{tikzpicture}
	\node (p1) {
	    \begin{tikzpicture}[polybox, mapstos]
            \node[poly, dom, "$\car{a}$" left] (S) {$\id_s$\at$s$};

            \draw (S_pos) to[climb'] node[right] {$\epsilon$} (S_dir);
        \end{tikzpicture}
	};
    \node[right=1 of p1] (p2) {
        \begin{tikzpicture}[polybox, mapstos]
            \node[poly, dom, "$\car{a}$" left] (S) {$\id_t$\at$t$};

            \draw (S_pos) to[climb'] node[right] {$\epsilon$} (S_dir);
        \end{tikzpicture}
	};
\end{tikzpicture}
\]
\[
\begin{tikzpicture}
	\node (p1) {
	  \begin{tikzpicture}[polybox, mapstos]
  	\node[poly, dom, "$\car{a}$" left] (p) {$\id_s$\at$s$};
  	\node[poly, cod, right=of p.south, yshift=-1ex, "$\car{a}$" right] (q) {$\id_s$\at$s$};
  	\node[poly, cod, above=of q, "$\car{a}$" right] (r) {$\id_s$\at$s$};
  	\draw (p_pos) to[first] (q_pos);
  	\draw (q_dir) to[climb, "$\delta$" left] (r_pos);
  	\draw (r_dir) to[last] (p_dir);
  \end{tikzpicture}
	};
	\node[right=1 of p1] (p2) {
	  \begin{tikzpicture}[polybox, mapstos]
  	\node[poly, dom, "$\car{a}$" left] (p) {$a$\at$s$};
  	\node[poly, cod, right= of p.south, yshift=-1ex, "$\car{a}$" right] (q) {$\id_s$\at$s$};
  	\node[poly, cod, above=of q, "$\car{a}$" right] (r) {$a$\at$s$};
  	\draw (p_pos) to[first] (q_pos);
  	\draw (q_dir) to[climb, "$\delta$" left] (r_pos);
  	\draw (r_dir) to[last] (p_dir);
  \end{tikzpicture}
	};
	\node[below=.05 of p1] (p3) {
	  \begin{tikzpicture}[polybox, mapstos]
  	\node[poly, dom, "$\car{a}$" left] (p) {$a$\at$s$};
  	\node[poly, cod, right= of p.south, yshift=-1ex, "$\car{a}$" right] (q) {$a$\at$s$};
  	\node[poly, cod, above=of q, "$\car{a}$" right] (r) {$\id_t$\at$t$};
  	\draw (p_pos) to[first] (q_pos);
  	\draw (q_dir) to[climb, "$\delta$" left] (r_pos);
  	\draw (r_dir) to[last] (p_dir);
  \end{tikzpicture}
	};
	\node[below=.05 of p2] (p4) {
	  \begin{tikzpicture}[polybox, mapstos]
  	\node[poly, dom, "$\car{a}$" left] (p) {$\id_t$\at$t$};
  	\node[poly, cod, right= of p.south, yshift=-1ex, "$\car{a}$" right] (q) {$\id_t$\at$t$};
  	\node[poly, cod, above=of q, "$\car{a}$" right] (r) {$\id_t$\at$t$};
  	\draw (p_pos) to[first] (q_pos);
  	\draw (q_dir) to[climb, "$\delta$" left] (r_pos);
  	\draw (r_dir) to[last] (p_dir);
  \end{tikzpicture}
	};
\end{tikzpicture}
\]
Of course, we have yet to check that $(\car{a},\epsilon,\delta)$ really is a comonoid, i.e.\ that the diagrams in \eqref{eqn.erasure_law} and \eqref{eqn.coassoc_law} commute.
We leave that for \cref{exc.walking_arrow_com}.
\end{example}

\begin{exercise} \label{exc.walking_arrow_com}
Verify that $(\car{a},\epsilon,\delta)$ as defined in \cref{ex.walking_arrow_com} obeys the erasure laws in \eqref{eqn.erasure_law} and the coassociative law in \eqref{eqn.coassoc_law}.
\begin{solution}
We will eventually see that comonoids in $\poly$ are categories; this gives us intuition for what is going on here. In fact the category corresponding to the comonoid in this exercise is the walking arrow, depicted in \eqref{eqn.walking_arrow_bent_cor}. The erasure (counit) laws are verified by looking at the picture of $\delta$ in \eqref{eqn.my_comonoid_delta} and noting that it sends the double (do-nothing) lines (the identities) back to the double lines.
\end{solution}
\end{exercise}

\index{linear polynomial!comonoid structure on}

\begin{exercise}\label{exc.linear_poly_comon}
Show that if $B$ is a set, then there exists a unique comonoid structure on the linear polynomial $B\yon$.
\begin{solution}
Given a set $B$, we wish to give a comonoid structure on $B\yon$ and show that it is unique.
There is only one way to define an eraser lens $\epsilon\colon B\yon\to\yon$: the on-position function is the unique map $!\colon B\to\1$, while every on-directions function $\1\to\1$ must be the identity.
Meanwhile $B\yon\tri B\yon\iso B^\2\yon$, so to specify a duplicator lens $\delta\colon B\yon\to B^\2\yon$, it suffices to specify an on-positions function $\delta_\1\colon B\to B^\2$, and every on-directions function will again be the identity.

All the comonoid laws should then hold trivially on directions, so it suffices to consider each law on positions.
The erasure laws imply that the composite functions
\[
    B\To{\delta_\1}B\times B\To{B\:\times\:!}B\times\1\iso B,
\]
where the second map is just the canonical left projection, and
\[
    B\To{\delta_\1}B\times B\To{!\:\times\:B}\1\times B\iso B,
\]
where the second nap is just the canonical right projection, are both the identity on $B$.
The only function $\delta_\1\colon B\to B\times B$ that satisfies this condition is the diagonal $a\mapsto(a,a)$.
It is easy to verify that the diagonal is coassociative, so this does define a unique comonoid structure on $B\yon$.\index{coassociativity}
(It turns out that this is equivalent to the well-known result that there is a unique comonoid structure on every set in $(\smset,\1,\times)$.)
\end{solution}
\end{exercise}

Once you know that state systems are comonoids in $\poly$, but not the only ones, the natural question to ask is ``what are all the other comonoids in $\poly$?''
Or perhaps, as we led you through this case study of $\car{s}$, you have already suspected the truth: a polynomial comonoid---what with its directions leading from one position to another, directions that can be run together associatively among which there are directions at every position that do nothing---is just another name for a category.

%-------- Section --------%
\section{Polynomial comonoids are categories}
\index{polynomial comonoid!as category|(}
\index{Ahman-Uustalu theorem}

What Ahman and Uustalu showed was that polynomial comonoids can be identified with categories.
Every category in the usual sense is a comonoid in $\poly$, and every comonoid in $\poly$ is a category.
We find their revelation to be truly shocking, and it suggests some very different ways to think about categories.
But let's go over their result first.

\begin{theorem}[Ahman-Uustalu]\label{thm.ahman_uustalu}
There is a one-to-one isomorphism-preserving correspondence between polynomial comonoids and (small) categories.
\end{theorem}

Our goal is to spell out this correspondence so that we can justly proclaim:

\slogan{Comonoids in $\poly$ are precisely categories!}

\subsection[Translating between comonoids and categories]{Translating between polynomial comonoids and categories}

First, we describe how to translate between the carrier $\car{c}$ of a comonoid $\com{C}\coloneqq(\car{c},\delta,\epsilon)$ and the objects and morphisms of the corresponding category $\cat{C}$.
The idea is pretty simple, and you may have already guessed it: positions are objects and directions are morphisms.

\subsubsection{Positions as objects, directions as morphisms}
\index{polynomial comonoid!positions and directions of}

More precisely, the positions of $\car{c}$ are the objects of $\cat{C}$:
\begin{equation} \label{eqn.pos_obj}
    \car{c}(\1)=\Ob\cat{C}.
\end{equation}
Then for each such position or object $i$, the $\car{c}[i]$-directions are the morphisms of $\cat{C}$ with domain $i$:
\begin{equation} \label{eqn.dir_mor}
    \car{c}[i]=\sum_{j\in\Ob\cat{C}}\cat{C}(i,j).
\end{equation}
The right hand side above is a little clumsier than the left; this is because while we are used to thinking of hom-sets of categories such as $\cat{C}(i,j)$, consisting of all morphisms in $\cat{C}$ with a fixed domain and codomain, we aren't used to thinking about the collection of all morphisms in $\cat{C}$ with a fixed domain and an arbitrary codomain quite as often.%
\footnote{Except, perhaps, in the context of coslice categories.}
On the other hand, the carrier \emph{only} encodes which morphisms have each object as its domain, i.e.\ which directions are at each position.
Codomains will be encoded in the data of the comonoid elsewhere.

This is the key difference in perspective between the polynomial comonoid perspective of categories, in contrast to our usual hom-set perspective: the polynomial perspective is in a sense domain-centric, as highlighted by the following definition.

\index{category!polynomial carrier of}

\begin{definition}[Polynomial carrier]
Let $\cat{C}$ be a category.
For every object $i$ in $\cat{C}$, denote the morphisms in $\cat{C}$ with domain $i$ by $\cat{C}[i]$, so that\tablefootnote{We may also write $f\colon i\to\_$ to denote an arbitrary morphism $f\in\cat{C}[i]$, i.e.\ a morphism $f$ in $\cat{C}$ with domain $i$ and an unspecified codomain.}
\[
    \cat{C}[i]\coloneqq\sum_{j\in\Ob\cat{C}}\cat{C}(i,j).
\]
Then the \emph{polynomial carrier}, or simply \emph{carrier}, of $\cat{C}$ is the polynomial
\[
    \sum_{i\in\Ob\cat{C}}\yon^{\cat{C}[i]}.
\]
\end{definition}

So everything we have said so far about the correspondence from \cref{thm.ahman_uustalu} can be summarized by saying that it preserves carriers: the carrier of the category $\cat{C}$ is the carrier $\car{c}$ of the comonoid $\com{C}$, so that $\Ob\cat{C}=\car{c}(\1)$ and $\cat{C}[i]=\car{c}[i]$.

\begin{remark}
If we take the perspective that categories are equal if and only if their objects and morphisms are equal and obey the same laws, and similarly that polynomials are equal if and only if their position-sets and direction-sets are equal sets, then \eqref{eqn.pos_obj} and \eqref{eqn.dir_mor} really can be just strict equalities.
This is why we are comfortable naming a ``one-to-one correspondence'' in \cref{thm.ahman_uustalu} rather than just, say, some form of equivalence.
Since the positions and directions of our polynomials always form \emph{sets}, however, the categories we obtain under this correspondence are also necessarily \emph{small}: their objects form a set, as do all of their morphisms.
But we won't worry too much about size issues beyond this.
\end{remark}

\begin{exercise} \label{exc.ema_polys}
What is the carrier of each of the following categories (up to isomorphism)?
\begin{enumerate}
	\item The category
	\begin{center}
	    \boxCD{exercisecolor}{$A\Too{f}B$}
	\end{center}
	where we have drawn every morphism except for the identity morphisms.
	\item The category
	\begin{center}
	    \boxCD{exercisecolor}{$B\To{g}A\From{h}C$}
	\end{center}
	where we have drawn every morphism except for the identity morphisms.
	\item The empty category.
	\item \label{exc.ema_polys.nat_monoid} A category with exactly $1$ object and a morphism $i$, for which every morphism can be written uniquely as the $n$-fold composite of $i$ for some $n\in\nn$.
	\item \label{exc.ema_polys.nat_poset} The category
	\begin{center}
	    \boxCD{exercisecolor}{$0\to1\to2\to3\to\cdots$}
	\end{center}
	where there is a unique morphism $m\to n$ if $m\leq n$ (and no other morphisms).
	\item The category
	\begin{center}
	    \boxCD{exercisecolor}{$0\from1\from2\from3\from\cdots$}
	\end{center}
	where there is a unique morphism $m\from n$ if $m\leq n$ (and no other morphisms).\qedhere
\end{enumerate}
\begin{solution}
\begin{enumerate}
	\item The category \boxCD{white}{$A\Too{f}B$} has $2$ morphisms out of $A$, namely $\id_A$ and $f$; and $1$ morphism out of $B$, namely $\id_B$.
	So its carrier is $\{A\}\yon^{\{\id_A,f\}}+\{B\}\yon^{\{\id_B\}}\iso\yon^\2+\yon$.
	\item The category \boxCD{white}{$B\To{g}A\From{h}C$} has $1$ morphism out of $A$, namely $\id_A$; $2$ morphisms out of $B$, namely $\id_B$ and $g$; and $2$ morphisms out of $C$, namely $\id_C$ and $h$.
	So its carrier is $\{A\}\yon^{\{\id_A\}}+\{B\}\yon^{\{\id_B,g\}}+\{C\}\yon^{\{\id_C,h\}}\iso\2\yon^\2+\yon$.
	\item The empty category has no objects or morphisms, so its carrier is just $\0$.
	\item The category in question has $1$ object, and its set of morphisms is in bijection with $\nn$, so its carrier is isomorphic to $\yon^\nn$.
	(This category is the monoid $(\nn,0,+)$ viewed as a $1$-object category; see \cref{ex.monoids} for the general case.)
	\item The category in question has $\nn$ as its set of objects, and for each $m\in\nn$, the morphisms out of $m$ are determined by their codomains: there is exactly $1$ morphism $m\to n$ for every $n\in\nn$ satisfying $m\leq n$, and no other morphisms out of $m$.
	So the carrier of the category is isomorphic to
	\[
	    \sum_{m\in\nn}\yon^{\{n\in\nn\mid m\leq n\}}\iso\nn\yon^\nn,
	\]
	as $\{n\in\nn\mid m\leq n\}\iso\nn$ under the bijection $n\mapsto n-m$.
	(This category is the poset $(\nn,\leq)$.)
	\item The category in question has $\nn$ as its set of objects, and for each $n\in\nn$, the morphisms out of $n$ are again determined by their codomains: there is exactly $1$ morphism $n\to m$ for every $m\in\nn$ satisfying $m\leq n$, and no other morphisms out of $n$.
	So the carrier of the category is isomorphic to
	\[
	    \sum_{n\in\nn}\yon^{\{m\in\nn\mid m\leq n\}}\iso\sum_{n\in\nn}\yon^{\ord{n+1}}\iso\yon^\1+\yon^\2+\yon^\3+\cdots.
	\]
    (This category is the poset $(\nn,\geq)$.)
\end{enumerate}
\end{solution}
\end{exercise}\index{monoid!of natural numbers}

\index{polynomial comonoid!as category|(}

But a category $\cat{C}$ is more than its carrier polynomial, just as a comonoid $\com{C}$ is more than its carrier $\car{c}$.
In particular, we have said nothing about the codomains of morphisms in $\cat{C}$, nor anything about identity morphisms, composition, or how the laws of a category are satisfied.
Similarly, we have said nothing about the eraser $\epsilon$ or the duplicator $\delta$ of $\com{C}$, nor anything about how the comonoid laws are satisfied.
It turns out that all of these constituents and laws correspond to one another, as summarized by the following table.
Here each item in the comonoid column---either a polynomial, a lens, or a lens equation---spans two rows, with the top row corresponding to positions and the bottom row corresponding to directions.
\begin{center}
{\scriptsize\begin{tabular}{cc|cc}
  \textbf{Comonoid} & $\com{C}=(\car{c},\epsilon,\delta)$ & \textbf{Category} & $\cat{C}$ \\
  \hline
  \multirow{2}{*}{carrier} & $i\in\car{c}(\1)$ & objects & $i\in\Ob\cat{C}$ \\
  & $f\in\cat{c}[i]$ & morphisms & $f\colon i\to\_$ \\
  \hline
  \multirow{2}{*}{eraser} & $\epsilon_\1\colon\car{c}(\1)\to\1$ & --- & --- \\
  & $\epsilon_i^\sharp\colon\1\to\car{c}[i]$ & identities & $\id_i\colon i\to\_$ \\
  \hline
  \multirow{2}{*}{duplicator} & $\delta_\1\colon\car{c}(\1)\to(\car{c}\tri\car{c})(\1)$ & codomains* & $\cod\colon\cat{C}[i]\to\Ob\cat{C}$ \\
  &
  $\delta_i^\sharp\colon(\car{c}\tri\car{c})[\delta_\1(i)]\to\car{c}[i]$ & composition* & $\then$ \\ %$\displaystyle\comp\colon\sum_{f\in\cat{C}[i]}\cat{C}[\cod(f)]\to\cat{C}[i]$
  \hline
  \multicolumn{2}{c|}{\multirow{2}{*}{right erasure law}} & \multicolumn2{c}{$\ast$} \\
  \multicolumn{2}{c|}{} & \multicolumn{2}{c}{right identity law} \\
  \hline
  \multicolumn{2}{c|}{\multirow{2}{*}{left erasure law}} & codomains of identities & $\cod\id_i=i$ \\
  \multicolumn{2}{c|}{} & \multicolumn{2}{c}{left identity law} \\
  \hline
  \multicolumn{2}{c|}{\multirow{2}{*}{coassociative law}} & codomains of composites & $\cod(f\then g)=\cod g$ \\
  \multicolumn{2}{c|}{} & \multicolumn{2}{c}{associative law of composition}
\end{tabular}
}\end{center}
\index{coassociativity}
 Note that the on-positions function of $\epsilon$, being a function into the terminal set, encodes no actual data.
The asterisk $\ast$ indicates that the right erasure law on positions works together with the duplicator to ensure that codomains and composites are properly specified.

We have already covered the correspondence between the first two rows, so let us consider each of the following rows in turn.
In some sense, we have already seen each piece of this correspondence in action for state systems in \cref{sec.comon.sharp.state}, so we'll go through it a little faster this time for the general case.

\subsubsection{The eraser assigns identities}

We know that the eraser $\epsilon\colon\car{c}\to\yon$ can be identified with a dependent function $(i\in\car{c}(\1))\to\com{c}[i]$, sending each position $i\in\com{c}(\1)$ to a $\com{c}[i]$-direction.
In terms of our category $\cat{C}$, the eraser sends each object $i\in\Ob\cat{C}$ to a morphism $i\to\_$.
But this is exactly what we need to specify identity morphisms---a morphism out of each object.
So the eraser of $\car{c}$ specifies the identity morphisms of the corresponding category $\cat{C}$.
We can interpret the polybox picture for $\epsilon$ like so:
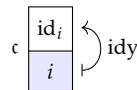
\[
 \begin{tikzpicture}[polybox, mapstos]
  	\node[poly, dom, "$\car{c}$" left] (c) {$\id_i$\at$i$};
  	\draw (c_pos) to[climb'] node[right] {$\idy$} (c_dir);
	\end{tikzpicture}
\]
Here we have given the label $\idy$ to the arrow sending objects to their identity morphisms.

Keep in mind that from the domain-centric polynomial perspective, we have not yet specified that the codomain of an identity morphism is equal to its domain; that comes later.

\subsubsection{The right erasure law on positions: some bookkeeping}

Keeping our label $\idy$ for the arrow in $\epsilon$, the right erasure law $\delta\then(\car{c}\tri\epsilon)=\id_\car{c}$ from \eqref{eqn.erasure_law} can be drawn in polyboxes like so:
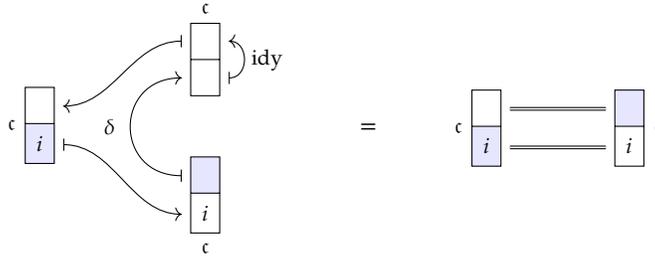
\[
\begin{tikzpicture}
	\node (1) {
        \begin{tikzpicture}[polybox, mapstos]
        	\node[poly, dom, "$\car{c}$" left] (r) {\at$i$};
        	\node[poly, cod, right=2 of r.south, yshift=-2.5ex, "$\car{c}$" below] (p) {\at$i$};
        	\node[poly, above=.8 of p, "$\car{c}$" above] (p') {};

        	\draw (r_pos) to[first] node[below] {} (p_pos);
        	\draw (p_dir) to[climb] node[left=.1] {$\delta$} (p'_pos);
        	\draw (p'_pos) to[climb'] node[right] {$\idy$} (p'_dir);
        	\draw (p'_dir) to[last] node[above=.4] {} (r_dir);
        \end{tikzpicture}
	};
	\node[right=1.8 of 1] (2) {
        \begin{tikzpicture}[polybox, mapstos]
          	\node[poly, dom, "$\car{c}$" left] (c) {\at$i$};
          	\node[poly, cod, right=of c, "$\car{c}$" right] (c') {\at$i$};
          	\draw[double,-] (c_pos) -- (c'_pos);
          	\draw[double,-] (c'_dir) -- (c_dir);
	    \end{tikzpicture}
	};
	\node at ($(1.east)!.5!(2.west)$) (3) {=};
\end{tikzpicture}
\]
We have only filled in a few of the boxes, but that is enough to interpret what the right erasure law tells us on positions: that the bottom arrow of the duplicator must be the identity function on $\car{c}(\1)$.
Equipped with this knowledge, we can focus our attention on the other two arrows of $\delta$.

\subsubsection{The duplicator assigns codomains and composites}

In fact, in the polybox picture for $\delta\colon\car{c}\to\car{c}\tri\car{c}$, the middle arrow specifies codomains, and the top arrow specifies composition.
We therefore label these arrows as follows:\footnote{Compare these labels to the names ``target'' and ``run'' that we gave to the arrows of a state system's transition lens.}
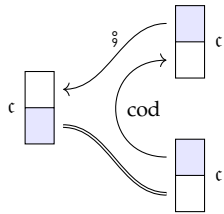
\[
\begin{tikzpicture}[polybox, tos]
    \node[poly, dom, "$\car{c}$" left] (yX) {};
    \node[poly, cod, right=1.8 of yX.south, yshift=-2.5ex, "$\car{c}$" right] (p) {};
    \node[poly, cod, above=.8 of p, "$\car{c}$" right] (p') {};

    \draw[double, -] (yX_pos) to[first] node[below] {} (p_pos);
    \draw (p_dir) to[climb] node[right] {$\cod$} (p'_pos);
    \draw (p'_dir) to[last] node[above] {$\then$} (yX_dir);
\end{tikzpicture}
\]
To check that this makes sense, we fill in the boxes:
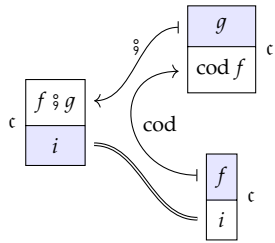
\[
\begin{tikzpicture}[polybox, mapstos]
    \node[poly, dom, "$\car{c}$" left] (yX) {$f\then g$\at$i$};
    \node[poly, cod, right=2 of yX.south, yshift=-2.5ex, "$\car{c}$" right] (p) {$f$\at$i$};
    \node[poly, cod, above=.8 of p, "$\car{c}$" right] (p') {$g$\at$\cod f$};

    \draw[double, -] (yX_pos) to[first] node[below] {} (p_pos);
    \draw (p_dir) to[climb] node[right] {$\cod$} (p'_pos);
    \draw (p'_dir) to[last] node[above] {$\then$} (yX_dir);
\end{tikzpicture}
\]
Remember: each position box contains an object of $\cat{C}$, while each direction box contains a morphism of $\cat{C}$ emanating from the object below.
So $\delta$ takes an object $i\in\Ob\cat{C}$ and a morphism $f\colon i\to\_$ in $\cat{C}$ and assigns another object $\cod f\in\Ob\cat{C}$ to be the codomain of $f$.
It then takes another morphism $g\colon\cod f\to\_$ in $\cat{C}$ and assigns a morphism $f\then g\colon i\to\_$ to be the composite of $f$ and $g$.
In this way, every morphism gets a codomain, and every pair of morphisms that can be composed (i.e.\ the codomain of one matches the domain of the other) is assigned a composite.
As with the identity morphism, we don't know what the codomain of this composite morphism is yet; but we do know that the domain of $f\then g$ matches the domain of $f$, as it should.

\subsubsection{The left erasure law on positions: codomains of identities}

As with the right erasure law, we can partially fill in the polyboxes for the left erasure law $\delta\then(\epsilon\tri\car{c})=\id_\car{c}$ from \eqref{eqn.erasure_law} to read what it says on positions:
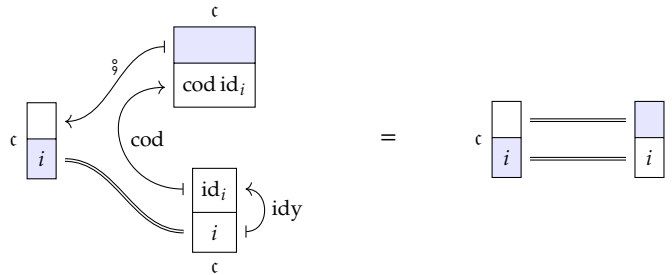
\[
\begin{tikzpicture}
	\node (1) {
        \begin{tikzpicture}[polybox, mapstos]
        	\node[poly, dom, "$\car{c}$" left] (r) {\at$i$};
        	\node[poly, right=2 of r.south, yshift=-2.5ex, "$\car{c}$" below] (p) {$\id_i$\at$i$};
        	\node[poly, cod, above=.8 of p, "$\car{c}$" above] (p') {\at$\cod\id_i$};

        	\draw[double, -] (r_pos) to[first] node[below] {} (p_pos);
        	\draw (p_pos) to[climb'] node[right] {$\idy$} (p_dir);
        	\draw (p_dir) to[climb] node[right] {$\cod$} (p'_pos);
        	\draw (p'_dir) to[last] node[above] {$\then$} (r_dir);
        \end{tikzpicture}
	};
	\node[right=1.8 of 1] (2) {
        \begin{tikzpicture}[polybox, mapstos]
          	\node[poly, dom, "$\car{c}$" left] (c) {\at$i$};
          	\node[poly, cod, right=of c, "$\car{c}$" right] (c') {\at$i$};
          	\draw[double,-] (c_pos) -- (c'_pos);
          	\draw[double,-] (c'_dir) -- (c_dir);
	    \end{tikzpicture}
	};
	\node at ($(1.east)!.5!(2.west)$) (3) {=};
\end{tikzpicture}
\]
So the left erasure law on positions guarantees that $\cod\id_i=i$ for all $i\in\Ob\cat{C}$.
It makes sense that we would find this here: the eraser assigns identities, while the duplicator assigns codomains, so a statement about codomains of identities is a coherence condition between the eraser and the duplicator.

\subsubsection{The erasure laws on directions are the identity laws}\index{polybox}

Let us finish filling in the polyboxes for the left and right erasure laws to see what they have to say on directions.
In the picture below, the left equality depicts the left erasure law (to conserve space, we'll substitute $i$ for $\cod\id_i$ on the left, which we now know we can do), while the right equality depicts the right erasure law:
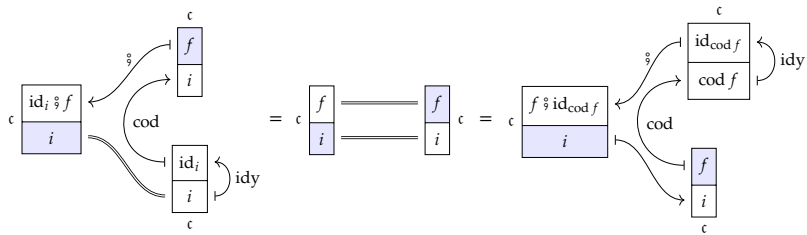
\[
\scalebox{.8}{
\begin{tikzpicture}
	\node (1) {
        \begin{tikzpicture}[polybox, mapstos]
        	\node[poly, dom, "$\car{c}$" left] (r) {$\id_i\then f$\at$i$};
        	\node[poly, right=2 of r.south, yshift=-2.5ex, "$\car{c}$" below] (p) {$\id_i$\at$i$};
        	\node[poly, cod, above=.8 of p, "$\car{c}$" above] (p') {$f$\at$i$};

        	\draw[double, -] (r_pos) to[first] node[below] {} (p_pos);
        	\draw (p_pos) to[climb'] node[right] {$\idy$} (p_dir);
        	\draw (p_dir) to[climb] node[right] {$\cod$} (p'_pos);
        	\draw (p'_dir) to[last] node[above] {$\then$} (r_dir);
        \end{tikzpicture}
	};
	\node[right=.1 of 1] (2) {
        \begin{tikzpicture}[polybox, mapstos]
          	\node[poly, dom, "$\car{c}$" left] (c) {$f$\at$i$};
          	\node[poly, cod, right=of c, "$\car{c}$" right] (c') {$f$\at$i$};
          	\draw[double,-] (c_pos) -- (c'_pos);
          	\draw[double,-] (c'_dir) -- (c_dir);
	    \end{tikzpicture}
	};
	\node[right=.25 of 2] (3) {
        \begin{tikzpicture}[polybox, mapstos]
        	\node[poly, dom, "$\car{c}$" left] (r) {$f\then\id_{\cod f}$\at$i$};
        	\node[poly, cod, right=2.1 of r.south, yshift=-2.5ex, "$\car{c}$" below] (p) {$f$\at$i$};
        	\node[poly, above=.8 of p, xshift=1.5ex, "$\car{c}$" above] (p') {$\id_{\cod f}$\at$\cod f$};

        	\draw (r_pos) to[first] node[below] {} (p_pos);
        	\draw (p_dir) to[climb] node[right] {$\cod$} (p'_pos);
        	\draw (p'_pos) to[climb'] node[right] {$\idy$} (p'_dir);
        	\draw (p'_dir) to[last] node[above] {$\then$} (r_dir);
        \end{tikzpicture}
	};
    \node at ($(1.east)!.2!(2.west)$) {=};
    \node at ($(2.east)!.5!(3.west)$) {=};
\end{tikzpicture}
}
\]
We find that on directions, the erasure laws state that for every object $i\in\Ob\cat{C}$ and morphism $f\colon i\to\_$ in $\cat{C}$,
\[
    \id_i\then f=f=f\then\id_{\cod f}.
\]
But these are precisely the identity laws of the category $\cat{C}$.

\subsubsection{The coassociative law on positions: codomains of composites}\index{coassociativity}

It remains to consider the comonoid's coassociative law \eqref{eqn.coassoc_law}, $\delta\then(\delta\tri\car{s})=\delta\then(\car{s}\tri\delta)$.
To read what it says on positions, we draw the polyboxes and fill them in, stopping just short of the uppermost direction box of the codomain:\index{polybox}
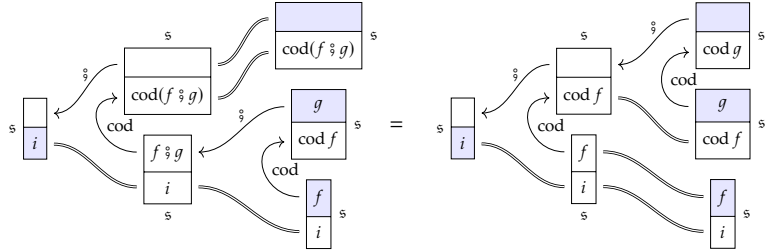
\[
\scalebox{.8}{
\begin{tikzpicture}
    \node (p1) {
        \begin{tikzpicture}[polybox, tos]
            \node[poly, dom, "$\car{s}$" left] (m') {\at$i$};
            \node[poly, right=1.8 of m'.south, yshift=-1ex, "$\car{s}$" below] (mm') {$f\then g$\at$i$};
            \node[poly, above=of mm', "$\car{s}$" above] (C') {\at$\cod(f\then g)$};
            \node[poly, cod, right=2.3 of mm'.south, yshift=-1ex, "$\car{s}$" right] (D') {$f$\at$i$};
            \node[poly, cod, above=of D', "$\car{s}$" right] (mmm') {$g$\at$\cod f$};
            \node[poly, cod, above=of mmm', "$\car{s}$" right] (CC') {\at$\cod(f\then g)$};
            \draw[double, -] (m'_pos) to[first] (mm'_pos);
            \draw (mm'_dir) to[climb] node[right] {$\cod$} (C'_pos);
            \draw (C'_dir) to[last] node[above] {$\then$} (m'_dir);
            \draw[double, -] (mm'_pos) to[first] (D'_pos);
            \draw (D'_dir) to[climb] node[right] {$\cod$} (mmm'_pos);
            \draw (mmm'_dir) to[last] node[above] {$\then$} (mm'_dir);
            \draw[double, -] (C'_pos) to[first] (CC'_pos);
            \draw[double, -] (CC'_dir) to[last] (C'_dir);
        \end{tikzpicture}
	};
	\node (p2) [right=.5 of p1] {
	    \begin{tikzpicture}[polybox, tos]
            \node[poly, dom, "$\car{s}$" left] (m) {\at$i$};
            \node[poly, right=1.8 of m.south, yshift=-1ex, "$\car{s}$" below] (D) {$f$\at$i$};
            \node[poly, above=of D, "$\car{s}$" above] (mm) {\at$\cod f$};
            \node[poly, cod, right=2.1 of D.south, yshift=-1ex, "$\car{s}$" right] (DD) {$f$\at$i$};
            \node[poly, cod, above=of DD, "$\car{s}$" right] (mmm) {$g$\at$\cod f$};
            \node[poly, cod, above=of mmm, "$\car{s}$" right] (C) {\at$\cod g$};
            \draw[double, -] (m_pos) to[first] (D_pos);
            \draw (D_dir) to[climb] node[right] {$\cod$} (mm_pos);
            \draw (mm_dir) to[last] node[above] {$\then$} (m_dir);
            \draw[double, -] (D_pos) to[first] (DD_pos);
            \draw[double, -] (DD_dir) to[last] (D_dir);
            \draw[double, -] (mm_pos) to[first] (mmm_pos);
            \draw (mmm_dir) to[climb] node[right] {$\cod$} (C_pos);
            \draw (C_dir) to[last] node[above] {$\then$} (mm_dir);
        \end{tikzpicture}
    };
	\node at ($(p1.south)!.5!(p2.north)$) {$=$};
\end{tikzpicture}
}\]
So on positions, the coassociative law states that given an object $i\in\Ob\cat C$ and morphisms $f\colon i\to\_$ and $g\colon\cod f\to\_$ in $\cat C$,
\[
    \cod(f\then g)=\cod g.
\]
Hence composites are assigned the proper codomains.

\subsubsection{The coassociative law on directions is the associative law of composition}

Finally, let us fill in the remaining polyboxes for the coassociative law (we'll substitute $\cod g$ for $\cod(f\then g)$ on the left, which we now know we can do):\index{polybox}
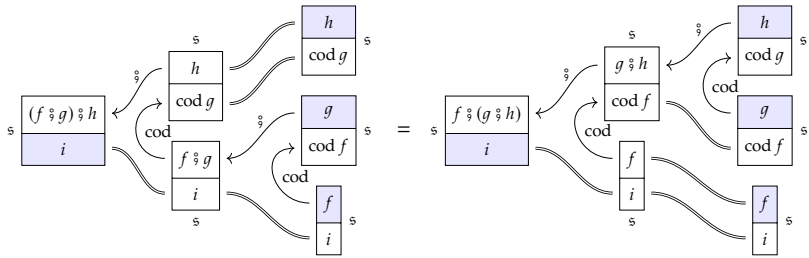
\[
\scalebox{.8}{
\begin{tikzpicture}
    \node (p1) {
        \begin{tikzpicture}[polybox, tos]
            \node[poly, dom, "$\car{s}$" left] (m') {$(f\then g)\then h$\at$i$};
            \node[poly, right=1.8 of m'.south, yshift=-1ex, "$\car{s}$" below] (mm') {$f\then g$\at$i$};
            \node[poly, above=of mm', "$\car{s}$" above] (C') {$h$\at$\cod g$};
            \node[poly, cod, right=2 of mm'.south, yshift=-1ex, "$\car{s}$" right] (D') {$f$\at$i$};
            \node[poly, cod, above=of D', "$\car{s}$" right] (mmm') {$g$\at$\cod f$};
            \node[poly, cod, above=of mmm', "$\car{s}$" right] (CC') {$h$\at$\cod g$};
            \draw[double, -] (m'_pos) to[first] (mm'_pos);
            \draw (mm'_dir) to[climb] node[right] {$\cod$} (C'_pos);
            \draw (C'_dir) to[last] node[above] {$\then$} (m'_dir);
            \draw[double, -] (mm'_pos) to[first] (D'_pos);
            \draw (D'_dir) to[climb] node[right] {$\cod$} (mmm'_pos);
            \draw (mmm'_dir) to[last] node[above] {$\then$} (mm'_dir);
            \draw[double, -] (C'_pos) to[first] (CC'_pos);
            \draw[double, -] (CC'_dir) to[last] (C'_dir);
        \end{tikzpicture}
	};
	\node (p2) [right=.5 of p1] {
	    \begin{tikzpicture}[polybox, tos]
            \node[poly, dom, "$\car{s}$" left] (m) {$f\then(g\then h)$\at$i$};
            \node[poly, right=2.2 of m.south, yshift=-1ex, "$\car{s}$" below] (D) {$f$\at$i$};
            \node[poly, above=of D, "$\car{s}$" above] (mm) {$g\then h$\at$\cod f$};
            \node[poly, cod, right=2 of D.south, yshift=-1ex, "$\car{s}$" right] (DD) {$f$\at$i$};
            \node[poly, cod, above=of DD, "$\car{s}$" right] (mmm) {$g$\at$\cod f$};
            \node[poly, cod, above=of mmm, "$\car{s}$" right] (C) {$h$\at$\cod g$};
            \draw[double, -] (m_pos) to[first] (D_pos);
            \draw (D_dir) to[climb] node[right] {$\cod$} (mm_pos);
            \draw (mm_dir) to[last] node[above] {$\then$} (m_dir);
            \draw[double, -] (D_pos) to[first] (DD_pos);
            \draw[double, -] (DD_dir) to[last] (D_dir);
            \draw[double, -] (mm_pos) to[first] (mmm_pos);
            \draw (mmm_dir) to[climb] node[right] {$\cod$} (C_pos);
            \draw (C_dir) to[last] node[above] {$\then$} (mm_dir);
        \end{tikzpicture}
    };
	\node at ($(p1.south)!.5!(p2.north)$) {$=$};
\end{tikzpicture}
}\]
Thus, on directions, the coassociative law states that given an object $i\in\Ob\cat C$ and morphisms $f\colon i\to\_, g\colon\cod f\to\_,$ and $h\colon\cod g\to\_$ in $\cat C$,
\[
    (f\then g)\then h = f\then(g\then h).
\]
But this is precisely the associative law of composition in a category.

We've seen that the data and equations of polynomial comonoids correspond exactly to the data and equations of categories.
This proves \cref{thm.ahman_uustalu}.

\index{polynomial comonoid!as category|)}

\subsubsection{Generalized duplicators as unbiased composition}

Before we move onto examples, one more note about the theory: notice that both sides of our coassociative law are given by $\delta^{(3)}\colon\car{c}\to\car{c}\tripow3$, as defined in \cref{prop.n_duplication}.\index{coassociativity}
On directions, $\delta^{(3)}$ tells us how to compose three morphisms $i\To{f}\_\To{g}\_\To{h}\_$ in $\cat{C}$ all at once to obtain $i\To{f\:\then\:g\:\then\:h}\_$, and (co)associativity ensures this is well-defined.

In general, $\delta^{(n)}\colon\car{c}\to\car{c}\tripow{n}$ on directions tells us how to compose $n$ morphisms in $\cat{C}$ for each $n\in\nn$.
After all, we have already seen that $\delta^{(2)}=\delta$ performs binary composition, that $\delta^{(1)}=\id_\car{c}$ performs ``unary'' composition (the ``unary composite'' of a single morphism $f$ is just $f$ itself), and that $\delta^{(0)}=\epsilon$ performs ``nullary'' composition (the ``nullary composite'' at any object is just its identity).
The directions of $\car{c}\tripow{n}$ at positions in the image of $\delta^{(n)}$ are exactly the sequences of composable morphisms of length $n$, and $\delta^{(n)}$ sends each sequence to the single direction that is its composite.

%---- Subsection ----%
\subsection{Examples of categories as comonoids}

\index{polynomial comonoid!examples}
Now that we know that polynomial comonoids are just categories, let's review some simple examples of categories and see how they may be interpreted as comonoids.
As we go through these examples, pay attention to how the polynomial perspective causes us to view these familiar categories somewhat differently than usual.

\index{polynomial comonoid!preorders as|(}
\index{category|seealso{polynomial comonoid}}

\subsubsection{Preorders}

A \emph{preorder} (or \emph{thin category}) is a category in which every morphism $f\colon c\to d$ is the \emph{only} morphism $c\to d$.%
\footnote{Sometimes these are also called \emph{posets}, short for \emph{partially ordered sets}, but strictly speaking the only isomorphisms in a poset are its identities, while a preorder allows objects to be isomorphic without being equal.}
Composition in preorders is easy to describe, because the composite of $c\to d$ and $d\to e$ is always just the unique arrow $c\to e$.
As such, preorders are some of the simplest examples of categories to consider---we already saw several in \cref{exc.ema_polys}---so let us interpret these as comonoids first.

\begin{example}\label{ex.walking_arrow_cat}
Let us revisit \cref{ex.walking_arrow_com}, where we first wrote down a comonoid that was not a state system.
We defined $\car{a}\coloneqq\{s\}\yon^{\{\id_s,a\}}+\{t\}\yon^{\{\id_t\}}\iso\yon^\2+\yon$ and gave it a comonoid structure, with eraser $\epsilon\colon\car{a}\to\yon$ specifying directions $\id_s$ and $\id_t$ and duplicator $\delta\colon\car{a}\to\car{a}\tri\car{a}$ pointing the direction $a$ at $t$.

Looking at the picture we drew of the comonoid in \eqref{eqn.walking_arrow_bent_cor}, it should come as no surprise that the corresponding category $\cat{A}$ is the \emph{walking arrow category}, which is a preorder with two objects and one morphism between them:
\begin{center}
    $\cat{A}\coloneqq\:$\boxCD{examplecolor}{
    $s\Too{a}t$}
\end{center}
Here we omit the identity morphisms from our picture, but we know that they exist.

The category $\cat{A}$ has two objects, the $\car{a}$-positions $s$ and $t$.
It has two morphisms with domain $s$, the $\car{a}[s]$-directions $\id_s$ and $a$; and one morphism with domain $t$, the $\car{a}[t]$-direction $\id_t$.
The morphisms $\id_s$ and $\id_t$ picked out by the erasure are the identity morphisms, and the duplicator assigns them codomains that are equal to their domains.
The duplicator also assigns $a$ the codomain $t$; and as $\cat{A}$ is then a preorder, composites are determined automatically.
\end{example}

\begin{exercise}
Let $(\car{c},\epsilon,\delta)$ be the comonoid corresponding to the preorder depicted as follows (identity morphisms omitted):
\begin{center}
    \boxCD{exercisecolor}{$B\From{f}A\To{g}C$}
\end{center}
\begin{enumerate}
    \item What is the carrier $\car{c}$?
    \item Characterize the eraser $\epsilon$.
    \item Characterize the duplicator $\delta$.\qedhere
\end{enumerate}
\begin{solution}
We are given a comonoid $(\car{c},\epsilon,\delta)$ corresponding to the preorder \boxCD{white}{$B\From{f}A\To{g}C$}.
\begin{enumerate}
    \item There are three morphisms with domain $A$, namely $\id_A,f,$ and $g$; the only other morphisms are the identity morphisms on $B$ and $C$.
    So the carrier is $\car{c}=\{A\}\yon^{\{\id_A,f,g\}}+\{B\}\yon^{\{\id_B\}}+\{C\}\yon^{\{\id_C\}}$.
    \item It suffices to specify the eraser $\epsilon\colon\car{c}\to\yon$ on directions: as always, $\epsilon_i^\sharp\colon\1\to\car{c}[i]$ picks out $\id_i$ for each $i\in\car{c}(\1)=\{A,B,C\}$.
    \item The duplicator $\delta\colon\car{c}\to\car{c}\tri\car{c}$ tells us the codomain of each morphism, as well as how every pair of composable morphisms compose (which in the case of a preorder can be deduced automatically).
    So we can completely characterize the behavior of $\delta$ using polyboxes as follows:
\[
\begin{tikzpicture}[polybox, mapstos, node distance=2ex and 1.4cm, font=\tiny]
  \node (a) {
  \begin{tikzpicture}
  	\node[poly, dom] (p) {$\id_A$\at$A$};
  	\node[poly, cod, right=of p.south, yshift=-1ex] (q) {$\id_A$\at$A$};
  	\node[poly, cod, above=of q] (r) {$\id_A$\at$A$};
  	\draw[double, -] (p_pos) to[first] (q_pos);
  	\draw (q_dir) to[climb] (r_pos);
  	\draw (r_dir) to[last] (p_dir);
  \end{tikzpicture}
  };
  \node[right=.2 of a] (b) {
  \begin{tikzpicture}
  	\node[poly, dom] (p) {$f$\at$A$};
  	\node[poly, cod, right=of p.south, yshift=-1ex] (q) {$\id_A$\at$A$};
  	\node[poly, cod, above=of q] (r) {$f$\at$A$};
  	\draw[double, -] (p_pos) to[first] (q_pos);
  	\draw (q_dir) to[climb] (r_pos);
  	\draw (r_dir) to[last] (p_dir);
  \end{tikzpicture}
  };
  \node[right=.2 of b] (c) {
  \begin{tikzpicture}
  	\node[poly, dom] (p) {$g$\at$A$};
  	\node[poly, cod, right=of p.south, yshift=-1ex] (q) {$\id_A$\at$A$};
  	\node[poly, cod, above=of q] (r) {$g$\at$A$};
  	\draw[double, -] (p_pos) to[first] (q_pos);
  	\draw (q_dir) to[climb] (r_pos);
  	\draw (r_dir) to[last] (p_dir);
  \end{tikzpicture}
  };
  \node[right=.2 of c] (d) {
  \begin{tikzpicture}
  	\node[poly, dom] (p) {$f$\at$A$};
  	\node[poly, cod, right=of p.south, yshift=-1ex] (q) {$f$\at$A$};
  	\node[poly, cod, above=of q] (r) {$\id_B$\at$B$};
  	\draw[double, -] (p_pos) to[first] (q_pos);
  	\draw (q_dir) to[climb] (r_pos);
  	\draw (r_dir) to[last] (p_dir);
  \end{tikzpicture}
  };
  \node[below=.1 of a] (e) {
  \begin{tikzpicture}
  	\node[poly, dom] (p) {$g$\at$A$};
  	\node[poly, cod, right=of p.south, yshift=-1ex] (q) {$g$\at$A$};
  	\node[poly, cod, above=of q] (r) {$\id_C$\at$C$};
  	\draw[double, -] (p_pos) to[first] (q_pos);
  	\draw (q_dir) to[climb] (r_pos);
  	\draw (r_dir) to[last] (p_dir);
  \end{tikzpicture}
  };
  \node[below=.1 of b] (f) {
  \begin{tikzpicture}
  	\node[poly, dom] (p) {$\id_B$\at$B$};
  	\node[poly, cod, right=of p.south, yshift=-1ex] (q) {$\id_B$\at$B$};
  	\node[poly, cod, above=of q] (r) {$\id_B$\at$B$};
  	\draw[double, -] (p_pos) to[first] (q_pos);
  	\draw (q_dir) to[climb] (r_pos);
  	\draw (r_dir) to[last] (p_dir);
  \end{tikzpicture}
	};
  \node[below=.1 of c] (g) {
  \begin{tikzpicture}
  	\node[poly, dom] (p) {$\id_C$\at$C$};
  	\node[poly, cod, right=of p.south, yshift=-1ex] (q) {$\id_C$\at$C$};
  	\node[poly, cod, above=of q] (r) {$\id_C$\at$C$};
  	\draw[double, -] (p_pos) to[first] (q_pos);
  	\draw (q_dir) to[climb] (r_pos);
  	\draw (r_dir) to[last] (p_dir);
  \end{tikzpicture}
	};
\end{tikzpicture}
\]
\end{enumerate}
\end{solution}
\end{exercise}

\begin{exercise} \label{exc.linear_poly_cat}
We showed in \cref{exc.linear_poly_comon} that for any set $B$, the linear polynomial $B\yon$ has a unique comonoid structure.
To what category does this comonoid correspond?
\begin{solution}
The linear polynomial $B\yon$ corresponds to a category whose objects form the set $B$ and whose only morphisms are identities: in other words, it is the discrete category on $B$.\index{category!discrete}
\end{solution}
\end{exercise}

\begin{exercise}\label{ex.star_shaped}
\begin{enumerate}
	\item Find a comonoid structure for the polynomial $p\coloneqq{\yon}^{\ord{n}+\1}+\ord{n}\yon$ whose corresponding category is a preorder.
	(It is enough to fully describe the category that it corresponds to.)
	\item Would you call your category ``star-shaped''?
\qedhere
\end{enumerate}
\begin{solution}
\begin{enumerate}
    \item The polynomial $p\coloneqq{\yon}^{\ord{n}+\1}+\ord{n}\yon$ has $n+1$ positions: $1$ with $n+1$ directions and the rest with $1$ direction each.
    So any category carried by $p$ has $n+1$ objects, of which only $1$ has any nonidentity morphisms coming out of it: in fact, it has $n$ nonidentity morphisms coming out of it.
    But if the category is to be a preorder, each of these $n$ nonidentity morphisms must have a distinct codomain.
    As there are exactly $n$ other objects, this completely characterizes the category.
    Equivalently, it is the discrete category on $\ord{n}$ adjoined with a unique initial object, so that the only nonidentity morphisms are the morphisms out of that initial object to each of the other objects exactly once.
    \item This category can be thought of as ``star-shaped'' if we picture the initial object in the center with morphisms leading out to the other $n$ objects like spokes.
\end{enumerate}\index{category!discrete}
\end{solution}
\end{exercise}
\index{polynomial comonoid!preorders as|)}

\index{state system!corresponding category of}

\begin{example}[State systems as categories] \label{ex.state_cat}
We know that every state system $\car{s}\iso S\yon^S$ with its do-nothing section $\epsilon\colon\car{s}\to\yon$ and its transition lens $\delta\colon\car{s}\to\car{s}\tri\car{s}$ is a comonoid, so what category $\cat{S}$ does $(\car{s},\epsilon,\delta)$ correspond to?

Recall from \cref{ex.not_all_com_state} that state systems are exactly those comonoids whose codomain (i.e.\ ``target'') functions $\cod\colon\car{s}[s]\to\car{s}(\1)$ for $s\in\car{s}(\1)$ are bijections.
That is, from every object $s\in\Ob\cat S=\car{s}(\1)$, there is exactly $1$ morphism to every object $t\in\Ob\cat S$.
So not only is $\cat S$ a preorder, it is the \emph{codiscrete preorder} on $\car{s}(\1)$, where there is always a morphism between every pair of objects.

Let's redraw the polyboxes for the do-nothing section of $\car{s}$ from \eqref{eqn.do_nothing_polybox} and the transition lens of $\car{s}$ from \eqref{eqn.trans_lens_polybox}, this time with our new arrow labels, as a sanity check:
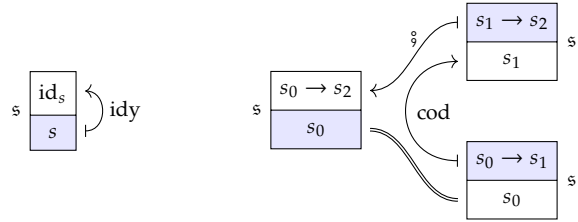
\[
\begin{tikzpicture}
\node (erase) {

\begin{tikzpicture}[polybox, mapstos]
    \node[poly, dom, "$\car{s}$" left] (S) {$\id_s$\at$s$};

    \draw (S_pos) to[climb'] node[right] {$\idy$} (S_dir);
\end{tikzpicture}

};
\node[right=of erase] (dupl) {

\begin{tikzpicture}[polybox, mapstos]
	\node[poly, dom, "$\car{s}$" left] (r) {$s_0\to s_2$\at$s_0$};
	\node[poly, cod, right=2 of r.south, yshift=-2.5ex, "$\car{s}$" right] (p) {$s_0\to s_1$\at$s_0$};
	\node[poly, cod, above=.8 of p, "$\car{s}$" right] (p') {$s_1\to s_2$\at$s_1$};

	\draw[double, -] (r_pos) to[first] node[below] {} (p_pos);
	\draw (p_dir) to[climb] node[right] {$\cod$} (p'_pos);
	\draw (p'_dir) to[last] node[above] {$\then$} (r_dir);
  \end{tikzpicture}

};
\end{tikzpicture}
\]
Indeed, we had already been writing the directions of $\car{s}$ as arrows $s\to t$, knowing that each was uniquely specified by its source $s$ and its target $t$ in $\car{s}(\1)$.
And in \cref{subsec.comon.sharp.state.cohere}, we had already noted that $\id_s$ was just the arrow $s\to s$.
So state systems have been categories with exactly one morphism between every pair of objects all along.

Other names for this category include the \emph{indiscrete preorder} and the \emph{codiscrete} or \emph{indiscrete category}.\index{category!codiscrete}
These names highlight the fact that every object of this category is isomorphic to every other object: in fact, every arrow $s\to t$ is an isomorphism with inverse $t\to s$, for these compose as $\id_s\colon s\to s$ in one direction and $\id_t\colon t\to t$ in the other.
Thus this category is also a \emph{groupoid}, and it may be called the \emph{codiscrete}, \emph{indiscrete}, or \emph{contractible groupoid}\dots but we will call it the \emph{state category on} $S$, where $S$ is the set of positions of $\car{s}$ or objects of $\cat S$.

Here are the state categories on $\3$ and on $\1\5$, with all maps (even identities) drawn:
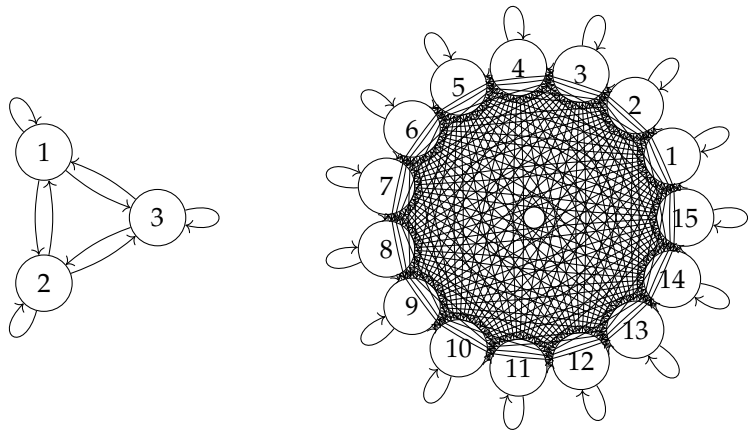
\[
\begin{tikzpicture}
\def\n{3}% how many nodes
\def\size{2cm}
\node[circle,minimum size=\size] (b) {};
\foreach\x in{1,...,\n}{
  \node[minimum size=0.75cm,draw,circle] (n-\x) at (b.{360/\n*\x}){\x};
}
\foreach\x in{1,...,\n}{
  \foreach\y in{1,...,\n}{
    \ifnum\x=\y\draw[->] (n-\x) to [in=360/\n*\x-15,out=360/\n*\x+15,loop] ();\relax\else
      \draw (n-\x) edge[->,bend right=10] (n-\y);
    \fi
  }
}
\def\n{15}% how many nodes
\def\size{4cm}
\node[circle,minimum size=\size, right=3 of b] (b) {};
\foreach\x in{1,...,\n}{
  \node[minimum size=0.75cm,draw,circle] (n-\x) at (b.{360/\n*\x}){\x};
}
\foreach\x in{1,...,\n}{
  \foreach\y in{1,...,\n}{
    \ifnum\x=\y\draw[->] (n-\x) to [in=360/\n*\x-15,out=360/\n*\x+15,loop] ();\relax\else
      \draw (n-\x) edge[->,bend right=3] (n-\y);
    \fi
  }
}
% Source: Asterix: https://tex.stackexchange.com/questions/390647/understanding-complete-graph-example-in-tikz
\end{tikzpicture}
\]
The picture on the left should look familiar: it is what we drew in \cref{ex.trans_trees} when took the corolla picture for $\3\yon^\3$ and bent the arrows to point at their targets according to its transition lens.
Notice that the graphs we obtain in this way are always complete.
\index{graph!complete}
\end{example}

\begin{exercise} \label{exc.not_state_cat_but_same_carrier}
Let $S$ be a set. Is there any comonoid structure on $S\yon^S$ other than that of the state category?
\begin{solution}
In the case of $S\coloneqq\0$, the only comonoid structure on $S\yon^S\iso\0$ is given by the empty category, the only category with no objects; and in the case of $S\coloneqq\1$, the only comonoid structure on $S\yon^S\iso\yon$ is given by the category with $1$ object and no nonidentity morphisms, again the only such category.
So in those cases, the comonoid structure on $S\yon^S$ is unique.

Now assume $|S|\geq2$.
The state category is always connected, but we can always find a comonoid structure on $S\yon^S$ given by a category that is not connected---and thus not isomorphic to the state category---as follows.
Consider a category whose object set is $S$ that has no morphisms between distinct objects, so that it is certainly not connected.
Then to specify the category, it suffices to specify a monoid associated with each object that will give the morphisms whose domain and codomain are equal to that object.
But there is always a monoid structure on a given nonempty set $S$.
If $S$ is finite, we can take the cyclic group $\zz/|S|\zz$ of order $|S|$, so that the resulting category has carrier $S\yon^{\zz/|S|\zz}\iso S\yon^S$.
On the other hand, if $S$ is infinite, we can take the free monoid on $S$, which has cardinality $\sum_{n\in\nn}|S|^\ord{n}=|\nn||S|=|S|$.
So the resulting category will again have carrier $S\yon^S$.
\end{solution}
\end{exercise}\index{monoid!structure on any nonempty set}

Not only does \cref{ex.state_cat} finally explain what our state systems really are (they're just special categories!), it illustrates two important features of our story.
One is that on positions, the duplicator $\delta\colon\car{c}\to\car{c}\tri\car{c}$ of a comonoid takes the corolla picture of $\car{c}$ and ``bends the arrows'' so that they point to other roots, yielding the underlying graph of a category.
Then $\delta$ on directions collapses two-arrow paths in the graph down to individual arrows, while the eraser $\epsilon\colon\car{c}\to\yon$ identifies empty paths with identity arrows.\index{graph!underlying a category}

Another important point is that we can view any category as a \emph{generalized state system}: its objects as \emph{states}, and its morphisms as \emph{transitions} between states.
The polynomial comonoid perspective is particularly suited for thinking about categories in this way: each object is a position that we could be in, and each morphism out of that object is a direction that we might take.
What is special about a comonoid is that each direction will always have another position at the end of it, making it reasonable to think of these directions as transitions between different states; and any sequence of transitions that we can follow is itself a transition.

Comparing these ideas, we see that they say the same thing: the first from the perspective of trees and graphs, the second from the perspective of dynamics.
We might say that
\slogan{
    a comonoid structure on a corolla forest turns\\
    roots into vertices and\\
    leaves into composable arrows between vertices;
}
or that
\slogan{
    a comonoid structure on an polynomial turns\\
    positions into states and\\
    directions into composable transitions between states.
}

\subsubsection{Monoids and monoid actions}
\index{monoid|(}
\index{polynomial comonoid!representable}

Here we use \emph{monoid} to refer to a monoid in the monoidal category $(\smset,\1,\times)$.
We denote such a monoid by $(M,e,*)$, where $M$ is the underlying set, $e\in M$ is the unit, and $*\colon M\times M\to M$ is the binary operation.

\begin{example}[Monoids as representable comonoids]\label{ex.monoids}
Recall that every monoid $(M,e,*)$ can be identified with a $1$-object category $\cat{M}$ with a single hom-set $M$, a single identity morphism $e$, and composition given by $*$.
Now we know that a $1$-object category $\cat{M}$ is also a polynomial comonoid $(\car{m},\epsilon,\delta)$ whose carrier has $1$ position, with all of the morphisms of $\cat{M}$ becoming its directions.
So the carrier of $\cat{M}$ is the representable polynomial $\yon^M$.

Then the eraser $\epsilon\colon\yon^M\to\yon$ picks out the identity morphism $e\in M$ on directions, while the duplicator $\delta\colon\yon^M\to\yon^M\tri\yon^M\iso\yon^{M\times M}$ can be identified with the binary operation $*\colon M\times M\to M$.
(We don't have to worry about codomains, since there's only one possible codomain to choose from.)
In this way, every monoid $(M,e,*)$ in $\smset$ gives rise to a representable comonoid $(\yon^M,\epsilon,\delta)$ in $\poly$.
We can just as easily invert this construction, obtaining a monoid for every representable comonoid by taking the underlying set to be the carrier's set of directions, the unit to be the direction picked out by the erasure, and the binary operation to be the duplicator's on-directions function.

% Let $(M,e,*)$ be a monoid.
% Recall that a

% Then we can construct a comonoid structure on the representable $\yon^M$. A morphism $\yon^M\to\yon$ can be identified with an element of $M$; under that identification we take $\epsilon\coloneqq e$. Similarly, $\yon^M\tri\yon^M\cong\yon^{M^\2}$ and a morphism $\yon^M\to\yon^{M^\2}$ can be identified with a function $M^\2\to M$; under that identification we take $\delta\coloneqq *$.
\end{example}

\begin{exercise}\index{associativity}
Verify \cref{ex.monoids} by showing that $(M,e,*)$ satisfies the unitality and associativity requirements of a monoid in $(\smset,\1,\times)$ if and only if $(\yon^M,\epsilon,\delta)$ satisfies the erasure and coassociativity requirements of a comonoid in $(\poly,\yon,\tri)$.
\begin{solution} %
The fact that monoids $(M,e,*)$ in $(\smset,\1,\times)$ are just comonoids $(\yon^M,\epsilon,\delta)$ in $(\poly,\yon,\tri)$, following the construction of \cref{ex.monoids}, is a direct consequence of the fully faithful Yoneda embedding $\smset\op\to\poly$ sending $A\mapsto\yon^A$ that maps $\1\mapsto\yon, A\times B\mapsto\yon^{A\times B}\iso\yon^A\tri\yon^B$ naturally, and $M\mapsto\yon^M$.\index{Yoneda lemma}

We can also state the laws and the correspondences between them explicitly, keeping in mind that $e$ and $*$ are just the on-directions functinos of $\epsilon$ and $\delta$.
The monoid's unitality condition states that $e$ is a $2$-sided unit for $*$, or that
\[
\begin{tikzcd}[row sep=large]
	\1\times M & M\ar[from=d, "*" description]\ar[from=r, equal]\ar[from=l, equal]&M\times\1\\&
	M\times M,\ar[from=ul, "e\:\times\:M"']\ar[from=ur, "M\:\times\:e"]
\end{tikzcd}
\]
commutes---equivalent to the comonoid's erasure laws, which state that
\[
\begin{tikzcd}[row sep=large]
	\yon\tri \yon^{M}&\yon^{M}\ar[d, "\delta" description]\ar[r, equal]\ar[l, equal]&\yon^{M}\tri\yon\\&
	\yon^{M}\tri\yon^{M},\ar[ul, "\epsilon\:\tri\:\yon^{M}"]\ar[ur, "\yon^{M}\:\tri\:\epsilon"']
\end{tikzcd}
\]
commutes (trivial on positions, equivalent to the monoid's unitality condition on directions).

Similarly, the monoid's associativity condition states that $*$ is associative, or that
\[
\begin{tikzcd}[row sep=large]
	M \ar[from=r, "*"']\ar[from=d, "*"]&
	M \times M \ar[from=d, "M\:\times\:*"']\\
	M \times M \ar[from=r, "*\:\times\:M"]&
	M \times M \times M,
\end{tikzcd}
\]\index{coassociativity}
commutes---equivalent to the comonoid's coassociative law, which state that
\[
\begin{tikzcd}[row sep=large]
	\yon^{M}\ar[r, "\delta"]\ar[d, "\delta"']&
	\yon^{M}\tri\yon^{M}\ar[d, "\yon^{M}\:\tri\:\delta"]\\
	\yon^{M}\tri \yon^{M}\ar[r, "\delta\:\tri\:\yon^{M}"']&
	\yon^{M}\tri\yon^{M}\tri\yon^{M},
\end{tikzcd}
\]
commutes (trivial on positions, equivalent to the monoid's associativity condition on directions).
\end{solution}
\end{exercise}

\index{list!cyclic}

\begin{example}[Cyclic lists]
For any $n\in\nn$, consider $\zz/n\zz$, the cyclic group of order $n$, viewed as a monoid or, equivalently, a $1$-object category.
Its carrier is $\yon^{\zz/n\zz}$.

As a polynomial functor, $\yon^{\zz/n\zz}$ sends each set $X$ to the set of length-$n$ tuples in $X$.
But the comonoid structure lets us think of these tuples as \emph{cyclic lists}: once we reach the last element, we can loop back around to the first element.
Indeed, as a natural transformation, $\epsilon\colon\yon^{\zz/n\zz}\to\yon$ picks out the ``current'' element via its $X$-component $\epsilon\tri X\colon\yon^{\zz/n\zz}\tri X\to\yon\tri X$, which is just a function $\epsilon_X\colon X^{\zz/n\zz}\to X$; and $\delta$ lets us move around the list.

We will see later that comonoids are closed under coproducts, so $\sum_{n\in\nn}\yon^{\zz/n\zz}$ is also a comonoid.
\end{example}

\index{monoid!action of}

\begin{example}[Monoid actions]\label{ex.monoid_action}
Suppose that $(M,e,*)$ is a monoid, $S$ is a set, and $\alpha\colon S\times M\to S$ is a \emph{(right) monoid action}.
That is, for all $s\in S$ we have $\alpha(s,e)=s$ and $\alpha(s,m*n)=\alpha(\alpha(s,m),n)$ for $m,n\in M$; equivalently, the diagrams
\[
\begin{tikzcd}
	S\times\1\ar[rd,equals]\ar[r,"S\:\times\:e"] & S\times M\ar[d,"\alpha"]\\
	& S
\end{tikzcd}
\qqand
\begin{tikzcd}[row sep=large]
    S\times M\times M\ar[r,"S\:\times\:*"]\ar[d,"\alpha\:\times\:M"] & S\times M\ar[d,"\alpha"]\\
    S\times M\ar[r,"\alpha"] & S
\end{tikzcd}
\]
commute.

Then there is an associated category $\cat{M\!A}$ with objects in $S$ and morphisms $s\To{m}\alpha(s,m)$ for each $s\in S$ and $m\in M$.
This, in turn, corresponds to a comonoid $(S\yon^M,\epsilon,\delta)$, as we will see in the next exercise.
\end{example}

\begin{exercise}
With notation as in \cref{ex.monoid_action}, characterize the comonoid structure on $S\yon^M$.
\begin{enumerate}
    \item How can we define the erasure $\epsilon$?
    \item How can we define the duplicator $\delta$?
    \item Verify that the erasure laws hold.
    \item Verify that the coassociative law holds.
    \item Describe the corresponding category $\cat{M\!A}$.
    In particular, what are the morphisms between any fixed pair of objects, what are the identity morphisms, and how do morphisms compose?
	\item $M$ always acts on itself by multiplication. Is the associated comonoid structure on $M\yon^M$ the same or different from the one coming from \cref{ex.state_cat}?
\qedhere
\end{enumerate}
\begin{solution}
Here $(M,e,*)$ is a monoid, $S$ is a set, $\alpha\colon S\times M\to S$ is a monoid action, and $\cat{M\!A}$ is the associated category, whose corresponding comonoid is $(S\yon^M,\epsilon,\delta)$.
We also know that for each $s\in S$ and $m\in M$, there is a morphism $s\To{m} \alpha(s,m)$ in $\cat{M\!A}$.
\begin{enumerate}
    \item The erasure $\epsilon\colon S\yon^M\to\yon$ picks out an element of $m\in M$ for every element $s\in S$ that will play the role of an identity, which in particular should also have $s$ as its codomain.
    Since we want the codomain of the morphism $m$ out of $s$ to be $\alpha(s,m)$, we can take $m=e$ to guarantee that its codomain will be $\alpha(s,e)=s$.
    So we let $\epsilon$ be the lens whose on-directions function at each $s\in S$ is $\epsilon^\sharp_s\colon\1\to M$ always maps to $e$.
    \item The duplicator $\delta\colon S\yon^M\to S\yon^M\tri S\yon^M$ is determined by what we want the codomain of each morphism to be and how we want the morphisms to compose.
    We already know that we want the morphism $m\in M$ out of each $s\in S$ to have the codomain $\alpha(s,m)$.
    If we then have another morphism $n\in M$ out of $\alpha(s,m)$, its codomain will be $\alpha(\alpha(s,m),n)=\alpha(s,m*n)$, the same as the codomain of the morphism $m*n$ out of $s$.
    So it makes sense for the composite $s\To{m}\alpha(s,m)\To{n} \alpha(\alpha(s,m),n)$ to be the morphism $s\To{m*n}\alpha(s,m*n)$.
    Thus, we can define $\delta$ in polyboxes as
    \[
    \begin{tikzpicture}[polybox, mapstos, font=\tiny]
        \node[poly, dom] (p) {$m*n$\at$s$};
        \node[poly, cod, right=1.5 of p.south, yshift=-1ex] (q) {$m$\at$s$};
        \node[poly, cod, above=of q, xshift=3] (r) {$n$\at$\alpha(s,m)$};
        \draw[double, -] (p_pos) to[first] (q_pos);
        \draw (q_dir) to[climb] node[right] {$\cod$} (r_pos);
        \draw (r_dir) to[last] node[above] {$\then$} (p_dir);
    \end{tikzpicture}
    \]
    \item We constructed $\delta$ above so that its bottom arrow is an identity function, so verifying the erasure laws amounts to checking that the direction $e\in M$ that $\epsilon$ picks out at each position $s\in S$ really does function as an identity morphism $s\To{e}\alpha(s,e)$ under the codomain and composition operations specified by $\delta$.
    We have already ensured that the codomain of $e$ at $s$ is $\alpha(s,e)=s$; meanwhile, given $m\in M$ we have that the composite of $s\To{m}\alpha(s,m)\To{e}\alpha(s,m)$ is $m*e=m$ and that the composite of $s\To{e}s\To{m}\alpha(s,m)$ is $e*m=m$ by the monoid's own unit laws.
    So the erasure laws hold.
    \item Verifying the coassociativity of $\delta$ amounts to checking that composition plays nicely with codomains and is associative.
    We already checked the former when defining $\delta$, and the latter follows from the monoid's own associativity laws: given $m,n,p\in M$, we have $(m*n)*p=m*(n*p)$.
    \item The associated category $\cat{M\!A}$ is a category whose objects are the elements of the set $S$ being acted on, and whose morphisms $s\to t$ for each $s,t\in S$ are the elements of the monoid $m\in M$ that send $s$ to $t$, i.e.\ $\alpha(s,m)=t$.
    The identity morphism on each object is just the unit $e\in M$, while morphisms compose via the multiplication $*$.

    \item It is the same iff $M$ is a group, i.e.\ if every $m\in M$ has an inverse. Indeed, the comonoid structure on $M\yon^M$ from \cref{ex.state_cat} corresponds to a category in which every map is an isomorphism, so for each $m\in M$, the left-action $\alpha(m,-)\colon M\to M$ would need to be a bijection, and this is the case iff $M$ is a group.
\end{enumerate}
\end{solution}
\end{exercise}

\index{monoid|)}

\index{streams!category of}
\begin{example}[The category of $B$-streams]\label{ex.streams_category}
Fix a set $B$.
The set $B^\nn$ consists of countable sequences of elements in $B$, which we will call \emph{$B$-streams}.
We can write an $B$-stream $\ol{b}\in B^\nn$ as
\[
    \ol{b}\coloneqq(b_0\to b_1\to b_2\to b_3\to\cdots),
\]
with $b_n\in B$ for each $n\in\nn$.

Then there is a monoid action $\tau\colon B^\nn\times\nn\to B^\nn$ for which
\[
    \tau(\ol{b},n)\coloneqq(b_n\to b_{n+1}\to b_{n+2}\to b_{n+3}\to\cdots).
\]
Roughly speaking, $\nn$ acts on $B$-streams by shifting them forward by a natural number of steps.
We can check that this is a monoid action by observing that $\tau(\ol{b},0)=\ol{b}$ and that $\tau(\ol{b},m+n)=\tau(\tau(\ol{b},m),n)$.

So by \cref{ex.monoid_action}, the corresponding comonoid is carried by $B^\nn\yon^\nn$.
Each $B$-stream $\ol{b}$ is a position, and each $n\in\nn$ is a direction at $\ol{b}$ that can be visualized as the sequence of $n$ arrows starting from $b_0$ and ending at $b_n$.
Then at the end of the direction $n$ is a new $B$-stream: the rest of $\ol{b}$ starting at $b_n$.
Indeed, this $B$-stream is exactly $\tau(\ol{b},n)$, the codomain assigned to the direction $n$ at $\ol{b}$.

Alternatively, if we shift from the domain-centric perspective to the usual hom-set perspective, this comonoid corresponds to a category whose objects are $B$-streams and whose morphisms $\ol{b}\to\ol{b'}$ consist of every way in which $\ol{b'}$ can be viewed as a continguous substream of $\ol{b}$: that is, there is a morphism $n\colon\ol{b}\to\ol{b'}$ for each $n\in\nn$ satisfying
\[
    (b_n\to b_{n+1}\to\cdots)=(b'_0\to b'_1\to\cdots).
\]
The identity on $\ol{b}$ is given by $0\colon\ol{b}\to\ol{b}$; and the composite of two morphisms is the sum of the corresponding natural numbers, as a substream of a substream of $\ol{b}$ is just a substream of $\ol{b}$ shifted by the appropriate amount.

We will see this category again in \cref{ex.streams_cofree}.
\end{example}

\begin{exercise}
Let $\rr/\zz\iso[0,1)$ be the quotient of $\rr$ by the $\zz$-action sending $(r,n)\mapsto r+n$.
More concretely, it is the set of real numbers between 0 and 1, including 0 but not 1.
\begin{enumerate}
	\item Find a comonoid structure on $(\rr/\zz)\yon^\rr$.
	\item Is the corresponding category a groupoid?
\qedhere
\end{enumerate}
\begin{solution}
\begin{enumerate}
    \item Since $\rr$ acts on $\rr/\zz$ by addition modulo 1, (e.g.\ $\alpha(.7,5.4)=.1$), we obtain a comonoid structure on $(\rr/\zz)\yon^\rr$ by \cref{ex.monoid_action}. For example, the erasure $(\rr/\zz)\yon^\rr\to\yon$ sends everything to $0\in\rr$, because $0$ is the identity in $\rr$.
    \item Yes it is a groupoid because $\rr$ is a group: every element is invertible.
\end{enumerate}
\end{solution}
\end{exercise}

\subsubsection{The degree of an object}

We could continue to list examples of polynomial comonoids, but of course any list of small categories is already a list of such comonoids.
So instead, we conclude this section with some terminology that the polynomial perspective on a category affords.

\index{category!degree of an object}

\begin{definition}[Degree, linear]
Let $\cat{C}$ be a category and $c\in\Ob\cat{C}$ an object. The \emph{degree of $c$}, denoted $\deg(c)$, is the set of arrows in $\cat{C}$ that emanate from $c$.

If $\deg(c)\iso\1$, we say that $c$ is \emph{linear}.
If $\deg(c)\iso\ord{n}$ for $n\in\nn$, we say $c$ has \emph{degree $n$}.
\end{definition}

\begin{exercise}
\begin{enumerate}
	\item If every object in $\cat{C}$ is linear, what can we say about $\cat{C}$?
	\item Is it possible for an object in $\cat{C}$ to have degree $0$?
	\item Find a category that has an object of degree $\nn$.
	\item Up to isomorphism, how many categories are there that have just one linear and one quadratic (degree 2) object?
	\item Is the above the same as asking how many comonoid structures on $\yon^\2+\yon$ there are?\qedhere
\end{enumerate}
\begin{solution}
\begin{enumerate}
    \item If every object in $\cat{C}$ is linear, then the only morphisms in $\cat{C}$ are the identity morphisms, so $\cat{C}$ must be a discrete category.\index{category!discrete}
    \item It is not possible for an object in $\cat{C}$ to have degree $0$, as every object must have at least an identity morphism emanating from it.
    \item Some possible examples of categories with objects of degree $\nn$ are the monoid $(\nn, 0, +)$ (see \cref{exc.ema_polys} \cref{exc.ema_polys.nat_monoid}), the poset $(\nn, \leq)$ (see \cref{exc.ema_polys} \cref{exc.ema_polys.nat_poset}), and the state category on $\nn$ (see \cref{ex.state_cat}).
    \item Up to isomorphism, there are $3$ categories with just one linear and one quadratic object.
    They can be distinguished by the behavior of the single nonidentity morphism.
    Either its domain and its codomain are distinct, in which case we have the walking arrow category; or its domain and its codomain are the same, in which case it can be composed with itself to obtain either itself or the identity.
    So there are $3$ possible categories in total.
    \item Yes: since (isomorphic) categories correspond to (isomorphic) comonoids, there are as many categories with one linear and one quadratic object up to isomorphism as there are comonoid structures on $\yon^\2+\yon$.
\end{enumerate}
\end{solution}
\end{exercise}

\index{polynomial comonoid!as category|)}

%-------- Section --------%
\section[Morphisms of polynomial comonoids are retrofunctors]{Morphisms of polynomial comonoids are retrofunctors%
  \sectionmark{Comonoid morphisms in $\poly$ are retrofunctors}}
\sectionmark{Comonoid morphisms in $\poly$ are retrofunctors}

\label{sec.comon.sharp.cof}

\index{retrofunctor}\index{polynomial comonoid!morphism of|see{retrofunctor}}

Now that we have characterized the comonoids of $\poly$, let us consider the morphisms between them.
These turn out to correspond to a rather odd kind of map between categories known as a \emph{retrofunctor}.
\footnote{Many authors have referred to these as \emph{cofunctors}, including ourselves in other work and in early versions of this book. However, the prefix \emph{co} in category theory is very special---having to do with taking opposites---and we will see in \cref{rem.retrofunctor_pare} that comonoid homomorphisms are not just opposite-functors.. Thus to keep the prefix \emph{co} more pristine, and in solidarity with other researchers, we have decided to use the term \emph{retrofunctor}, which is an appropriate usage of the term defined by Bob Par\'e \cite{pare2023retrocells}.}

%---- Subsection ----%
\subsection[Comonoid morphisms and retrofunctors]{Introducing comonoid morphisms and retrofunctors}

First, let us define morphisms of comonoids in the most general setting.
If you've seen the definition of a monoid homomorphism (or even a group homomorphism), then this definition may look familiar.

\index{comonoid!morphism of}

\begin{definition}[Comonoid morphism]\label{def.morphism_comonoids}
Given a monoidal category $(\cat{C},\yon,\tri)$ with comonoids $\com{C}\coloneqq(\car{c},\epsilon,\delta)$ and $\com{C}'\coloneqq(\car{c'},\epsilon',\delta')$, a \emph{comonoid morphism} (or \emph{morphism of comonoids}) $\com{C}\to\com{C}'$ is a morphism $F\colon \car{c}\to\car{c'}$ in $\cat{C}$ for which the following diagrams commute:
\begin{equation}\label{eqn.pres_era}
\begin{tikzcd}
    \car{c}\ar[r, "F"]\ar[d, "\epsilon"']&
    \car{c'}\ar[d, "\epsilon'"]\\
    \yon\ar[r, equal]&
    \yon,
\end{tikzcd}
\end{equation}
called the \emph{eraser preservation law} (we say $F$ \emph{preserves erasure}); and
\begin{equation}\label{eqn.pres_dup}
\begin{tikzcd}
    \car{c}\ar[r, "F"]\ar[d, "\delta"']&
    \car{c'}\ar[d, "\delta'"]\\
    \car{c}\tri\car{c}\ar[r, "F\:\tri\:F"']&
    \car{c'}\tri\car{c'}
\end{tikzcd}
\end{equation}
called the \emph{duplicator preservation law} (we say $F$ \emph{preserves duplication}).
We may also say that $F$ \emph{preserves the comonoid structure}.

When the monoidal structure on $\cat{C}$ can be inferred, we let $\comon(\cat{C})$ denote the subcategory of $\cat{C}$ whose objects are comonoids in $\cat{C}$ and whose morphisms are comonoid morphisms.
\end{definition}

So when our monoidal category of interest is $(\poly,\yon,\tri)$, a morphism between polynomial comonoids is just a special kind of lens between their carriers that preserves erasure and duplication.

\begin{exercise}
There is something to be proved in the definition above: that comonoids and comonoid morphisms really do form a category.
Using the notation from \cref{def.morphism_comonoids}, verify the following:
\begin{enumerate}
    \item The identity morphism on a comonoid is a comonoid morphism.
    \item The composite of two comonoid morphisms is a comonoid morphism.
\end{enumerate}
This will show that $\comon(\cat{C})$ is indeed a subcategory of $\cat{C}$.
\begin{solution}
As in \cref{def.morphism_comonoids}, we have a monoidal category $(\cat{C},\yon,\tri)$ with comonoids $\com{C}\coloneqq(\car{c},\epsilon,\delta)$ and $\com{C}'\coloneqq(\car{c'},\epsilon',\delta')$ and a comonoid morphism $F\colon\com{C}\to\com{C'}$ (really a morphism $F\colon\car{c}\to\car{c'}$ in $\cat{C}$).
Let's throw in another comonoid $\com{C}''\coloneqq(\car{c''},\epsilon'',\delta'')$ and another comonoid morphism $G\colon\com{C'}\to\com{C''}$ (really a morphism $G\colon\car{c'}\to\car{c''}$ in $\cat{C}$).
\begin{enumerate}
    \item To show that the identity morphism $\id_\car{c}\colon\car{c}\to\car{c}$ is a comonoid morphism $\com{C}\to\com{C}$, we must check that it preserves erasure by showing that \eqref{eqn.pres_era} commutes, then check that it preserves duplication by showing that \eqref{eqn.pres_era} commutes:
    \[
    \begin{tikzcd}
        \car{c}\ar[r, "\id_\car{c}", equal]\ar[d, "\epsilon"']&
        \car{c}\ar[d, "\epsilon"]\\
        \yon\ar[r, equal]&
        \yon
    \end{tikzcd}
    \hspace{.5in}
    \begin{tikzcd}[column sep=large]
        \car{c}\ar[r, "\id_\car{c}", equal]\ar[d, "\delta"']&
        \car{c}\ar[d, "\delta"]\\
        \car{c}\tri\car{c}\ar[r, "\id_\car{c}\:\tri\:\id_\car{c}", equal]&
        \car{c}\tri\car{c}.
    \end{tikzcd}
    \]
    But they do commute, since $\id_\car{c}$ is the identity on $\car{c}$ and $\id_\car{c}\tri\id_\car{c}$ is the identity on $\car{c}\tri\car{c}$.

    \item To show that the composite $F\then G\colon\car{c}\to\car{c''}$ of the two comonoid morphisms $F$ and $G$ is itself a comonoid morphism $\com{C}\to\com{C''}$, we check that it preserves erasure by showing that \eqref{eqn.pres_era} commutes, then check that it preserves duplication by showing that \eqref{eqn.pres_era} commutes:
    \[
    \begin{tikzcd}
        \car{c}\ar[r, "F\:\then\:G"]\ar[d, "\epsilon"']&
        \car{c''}\ar[d, "\epsilon''"]\\
        \yon\ar[r, equal]&
        \yon
    \end{tikzcd}
    \hspace{.5in}
    \begin{tikzcd}[column sep=huge]
        \car{c}\ar[r, "F\:\then\:G"]\ar[d, "\delta"']&
        \car{c''}\ar[d, "\delta''"]\\
        \car{c}\tri\car{c}\ar[r, "(F\:\then\:G)\:\tri\:(F\:\then\:G)"]&
        \car{c''}\tri\car{c''}.
    \end{tikzcd}
    \]
    But we can rewrite these squares like so, using the fact that $(F\then G)\tri(F\then G)=(F\tri F)\then(G\tri G)$ on the right:
    \[
    \begin{tikzcd}
        \car{c}\ar[r, "F"]\ar[d, "\epsilon"']&
        \car{c}\ar[r, "G"]\ar[d, "\epsilon'"']&
        \car{c''}\ar[d, "\epsilon''"]\\
        \yon\ar[r, equal]&
        \yon\ar[r, equal]&
        \yon
    \end{tikzcd}
    \hspace{.5in}
    \begin{tikzcd}
        \car{c}\ar[r, "F"]\ar[d, "\delta"']&
        \car{c'}\ar[r, "G"]\ar[d, "\delta"']&
        \car{c''}\ar[d, "\delta'"]\\
        \car{c}\tri\car{c}\ar[r, "F\:\tri\:F"]&
        \car{c'}\tri\car{c'}\ar[r, "G\:\tri\:G"]&
        \car{c''}\tri\car{c''}.
    \end{tikzcd}
    \]
    Then the left square in each diagram commutes because $F$ is a comonoid morphism, while the right square in each diagram commutes because $G$ is a comonoid morphism.
    So both diagrams commute.
\end{enumerate}
\end{solution}
\end{exercise}

Notice that we were very careful in how we stated \cref{thm.ahman_uustalu}: while we asserted the existence of an isomorphism-preserving one-to-one correspondence between the objects of $\comon(\poly)$ and $\smcat$, we never claimed that these two categories are isomorphic or even equivalent.
The strange truth of the matter is that they are not: polynomial comonoid morphisms correspond not to functors, but to different maps of categories called \emph{retrofunctors}.

How exactly do these maps behave?
If we specify \cref{def.morphism_comonoids} to the case of $(\poly,\yon,\tri)$, we can write the eraser preservation law \eqref{eqn.pres_era} in polyboxes as
\begin{equation}\label{eqn.pres_era_draw}
\begin{tikzpicture}
	\node (id1) {
	\begin{tikzpicture}[polybox, mapstos]
		\node[poly, dom, my-blue, "$\car{c}\vphantom{\car{c'}}$" above] (p) {\at$c$};
		\draw[my-blue,double,-] (p_pos) to[climb'] node[right] {$\idy$} (p_dir);
	\end{tikzpicture}
	};
	\node[right=of id1] (id2) {
	\begin{tikzpicture}[polybox, mapstos]
		\node[poly, dom, my-blue, "$\car{c}$" above] (p) {\at$c$};
		\node[poly, my-red, right=1 of p, "$\car{c'}$" above] (q) {};
		\draw (p_pos) to[first] (q_pos);
		\draw (q_dir) to[last] (p_dir);
		\draw[my-red] (q_pos) to[climb'] node[right] {$\idy$} (q_dir);
		\node at ($(p.east)!.3!(q.west)$) {$F$};
	\end{tikzpicture}
	};
	\node at ($(id1.east)!.3!(id2.west)-(0,6pt)$) {$=$};
\end{tikzpicture}
\end{equation}
and the duplicator preservation law \eqref{eqn.pres_dup} in polyboxes as
\begin{equation}\label{eqn.pres_dup_draw}
\begin{tikzpicture}
	\node (sp1) {
	\begin{tikzpicture}[polybox, mapstos]
		\node[poly, dom, my-blue, "$\car{c}$" above] (c) {\at$c$};
		\node[poly, my-blue, right=2 of c.south, yshift=-1ex, "$\car{c}$" below] (c1) {$\vphantom{g}$};
		\node[poly, my-blue, above=of c1, "$\car{c}$" above] (c2) {$\vphantom{h}$};
		\node[poly, cod, my-red, right=1 of c1, "$\car{c'}$" below] (c'1) {$g$};
		\node[poly, cod, my-red, right=1 of c2, "$\car{c'}$" above] (c'2) {$h$};
		\node at ($(c1.east)!.5!(c'1.west)$) {$F$};
		\node at ($(c2.east)!.5!(c'2.west)$) {$F$};
		\draw[my-blue,double,-] (c_pos) to[first] (c1_pos);
		\draw[my-blue] (c1_dir) to[climb] node[right] {$\cod$} (c2_pos);
		\draw[my-blue] (c2_dir) to[last] node[above] {$\then$} (c_dir);
		\draw (c1_pos) to[first] (c'1_pos);
		\draw (c'1_dir) to[last] (c1_dir);
		\draw (c2_pos) to[first] (c'2_pos);
		\draw (c'2_dir) to[last] (c2_dir);
    \end{tikzpicture}
	};
	\node[right=of sp1] (sp2) {
	\begin{tikzpicture}[polybox, mapstos]
		\node[poly, dom, my-blue, "$\car{c}$" above] (c) {\at$c$};
		\node[poly, my-red, right=1 of c, "$\car{c'}$" above] (c') {};
		\node[poly, cod, my-red, right=2 of c'.south, yshift=-1ex, "$\car{c'}$" below] (c'1) {$g$};
		\node[poly, cod, my-red, above=of c'1, "$\car{c'}$" above] (c'2) {$h$};
		\node at ($(c.east)!.3!(c'.west)$) {$F$};
		\draw (c_pos) to[first] (c'_pos);
		\draw (c'_dir) to[last] (c_dir);
		\draw[my-red,double,-] (c'_pos) to[first] (c'1_pos);
		\draw[my-red] (c'1_dir) to[climb] node[right] {$\cod$} (c'2_pos);
		\draw[my-red] (c'2_dir) to[last] node[above] {$\then$} (c'_dir);
	\end{tikzpicture}
	};
	\node at ($(sp1.east)!.5!(sp2.west)-(0,4pt)$) {$=$};
\end{tikzpicture}
\end{equation}
If we read off the equations from these polyboxes, interpreting polynomial comonoids as categories, we derive the following definition of a retrofunctor.
(Here \eqref{eqn.pres_id} is equivalent to \eqref{eqn.pres_era_draw}, while \eqref{eqn.pres_cod} and \eqref{eqn.pres_comp} are together equivalent to \eqref{eqn.pres_dup_draw}.)

\begin{definition}[Retrofunctor]\label{def.retrofunctor}
Let $\cat{C}$ and $\cat{C'}$ be (small) categories.
A \emph{retrofunctor} $F\colon\cat{C}\cof\cat{C'}$ consists of
\begin{itemize}
    \item a function $F\colon\Ob\cat{C}\to\Ob\cat{C'}$ \emph{forward on objects}\tablefootnote{In keeping with standard functor notation, we omit the usual subscript $\1$ that we include for on-positions (in this case, on-objects) functions. We often omit parentheses when applying these functions as well.} and
    \item a function $F^\sharp_c\colon\cat{C'}[Fc]\to\cat{C}[c]$ \emph{backward on morphisms} for each $c\in\Ob\cat{C}$,
\end{itemize}
satisfying the following conditions, collectively known as the \emph{retrofunctor laws}:
\begin{enumerate}[itemsep=0pt, label=\roman*.]
	\item $F$ \emph{preserves identities}:
	\begin{equation} \label{eqn.pres_id}
	    F^\sharp_c\,\id_{Fc}=\id_c
	\end{equation}
	for each $c\in\Ob\cat{C}$;
	\item $F$ \emph{preserves codomains}:
	\begin{equation} \label{eqn.pres_cod}
	    F\cod F^\sharp_c g=\cod g
	\end{equation}
	for each $c\in\Ob\cat{C}$ and $g\in\cat{C'}[Fc]$;
	\item $F$ \emph{preserves composites}:%
	\tablefootnote{In particular, the codomains of either side of \eqref{eqn.pres_comp} are equal.
	This isn't actually guaranteed by the other laws, so it is worth noting on its own; see for example the proof of \cref{prop.retrofunctors_isos}.}
	\begin{equation} \label{eqn.pres_comp}
	    F^\sharp_c g\then F^\sharp_{\cod F^\sharp_c g} h=F^\sharp_c\left(g\then h\right)
	\end{equation}
	for each $c\in\Ob\cat{C}, g\in\cat{C'}[Fc],$ and $h\in\cat{C'}[\cod g]$.
\end{enumerate}
We let $\catsharp\iso\comon(\poly)$ denote the category of (small) categories and retrofunctors.
\end{definition}

\begin{remark}\label{rem.retrofunctor_pare}
For experts, we explain the term \emph{retrofunctor} from \cref{def.retrofunctor}. Let $\cat{C},\cat{C'}$ be categories, and consider them as monads in $\mathbb{S}\Cat{pan}$. A functor between them consists of a function $F_\mathrm{Ob}\colon\Ob(\cat{C})\to\Ob(\cat{C'})$ and a 2-cell, as shown left
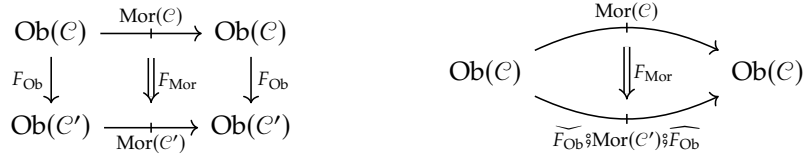
\[
\begin{tikzcd}[column sep=large]
	\Ob(\cat{C})\ar[r, tick, "\mathrm{Mor}(\cat{C})", ""' name=MM]\ar[d, "F_\mathrm{Ob}"']&
	\Ob(\cat{C})\ar[d, "F_\mathrm{Ob}"]\\
	\Ob(\cat{C'})\ar[r, tick, "\Mor(\cat{C'})"', "" name=MM']&
	\Ob(\cat{C'})
	\ar[from=MM, to=MM', Rightarrow, shorten=4pt, "F_{\mathrm{Mor}}"]
\end{tikzcd}
\hspace{.7in}
\begin{tikzcd}[column sep=70pt]
	\Ob(\cat{C})
		\ar[r, tick, bend left=25pt, "\mathrm{Mor}(\cat{C})", ""' name=MM]
		\ar[r, tick, bend right=25pt, "\wc{F_\mathrm{Ob}}\then\mathrm{Mor}(\cat{C'})\then \wh{F_\mathrm{Ob}}"', "" name=MM']&
	\Ob(\cat{C})
  	\ar[from=MM, to=MM', Rightarrow, shorten=2pt, "F_{\mathrm{Mor}}"]
\end{tikzcd}
\]
satisfying the properties of a monad homomorphism.

In \cite[Definition 6.1]{pare2023retrocells}, Par\'e defines \emph{retromorphism of monads} in double categories like $\mathbb{S}\Cat{pan}$. We will not discuss the more general definition here, but in a framed bicategory (equipment), where we have companions $\wc{f}$ and conjoints $\wh{f}$ of tight maps $f$, it is easy to check that monads lift along tight morphisms in the sense that the horizontal cell $\left(\wc{F_\mathrm{Ob}}\then\mathrm{Mor}(\cat{C'})\then \wh{F_\mathrm{Ob}}\right)\colon \Ob(\cat{C})\tickar \Ob(\cat{C})$ is a monad if $\mathrm{Mor}(\cat{C'})$ is. Hence, a functor is equivalently given by a function $F_\mathrm{Ob}\colon\Ob(\cat{C})\to\Ob(\cat{C'})$ such that the 2-cell shown above right is a monad homomorphism.

A retromorphism of monads---in our case, a \emph{retrofunctor}---is simply a monad map going the other way:%
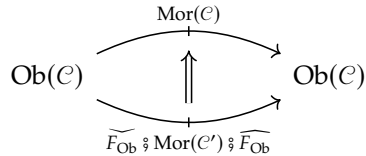
\footnote{
As a hint to the connection in terms of the (yet-undefined) double category $\mathbb{C}\Cat{at}^\sharp$, note that when the monads in question are left adjoints, their right adjoints will automatically be comonads, and a morphism between these comonads will be a retromorphism between the monads.
}
\[
\begin{tikzcd}[column sep=70pt]
	\Ob(\cat{C})
		\ar[r, tick, bend left=25pt, "\mathrm{Mor}(\cat{C})", ""' name=MM]
		\ar[r, tick, bend right=25pt, "\wc{F_\mathrm{Ob}}\:\then\:\mathrm{Mor}(\cat{C'})\:\then\:\wh{F_\mathrm{Ob}}"', "" name=MM']&
	\Ob(\cat{C})
  	\ar[from=MM', to=MM, Rightarrow, shorten=2pt]
\end{tikzcd}
\]
\end{remark}

\index{$\catsharp$}\index{comonoids|seealso{$\catsharp$}}

Henceforth we will identify the category $\catsharp$ with the isomorphic category $\comon(\poly)$, eliding the difference between comonoids in $\poly$ and categories.

Since each retrofunctor includes a lens between its carrier polynomials, retrofunctors compose the way lenses do.

\begin{exercise}
Let $\cat{C},\cat{D},\cat{E}$ be categories, and let $F\colon\cat{C}\cof\cat{D}$ and $G\colon\cat{D}\cof\cat{E}$ be retrofunctors between them.
\begin{enumerate}
    \item Characterize the behavior of the identity retrofunctor $\id_\cat{D}$ on $\cat{D}$.
    Where does it send each object?
    Where does it send each morphism?
    \item Characterize the behavior of the composite retrofunctor $F\then G$.
    Where does it send each object and morphism?\qedhere
\end{enumerate}
\begin{solution}
Here $F\colon\cat{C}\cof\cat{D}$ and $G\colon\cat{D}\cof\cat{E}$ are retrofunctors in $\catsharp$.
\begin{enumerate}
    \item The identity retrofunctor $\id_\cat{D}$ on $\cat{D}$ should correspond to the identity lens on the carrier of $\cat{D}$, which is the identity on both positions (objects) and directions (morphisms).
    So $\id_\cat{D}$ sends each object $d\in\Ob\cat{D}$ to itself, while $\left(\id_\cat{D}\right)^\sharp_d\colon\cat{D}[d]\to\cat{D}[d]$ sends each morphism out of $d$ to itself as well.
    \item The composite retrofunctor $F\then G\colon\cat{C}\cof\cat{D}$ should correspond to the composite of $F$ as a lens between the carriers of $\cat{C}$ and $\cat{D}$ with $G$ as a lens between the carriers of $\cat{D}$ and $\cat{E}$.
    So on objects, $F\then G$ sends each $c\in\Ob\cat{C}$ to $G(Fc)\in\Ob\cat{E}$.
    Then given $c\in\Ob\cat{C}$, the on-morphisms function $(F\then G)^\sharp_c\colon\cat{E}[G(Fc)]\to\cat{C}[c]$ is the composite of on-directions functions
    \[
        \cat{E}[G(Fc)]\To{G^\sharp_{Fc}}\cat{D}[Fc]\To{F^\sharp_c}\cat{C}[c],
    \]
    sending each morphism $h$ with domain $G(Fc)$ to $F^\sharp_c\left(G^\sharp_{Fc}h\right)$.
\end{enumerate}
\end{solution}
\end{exercise}

\index{retrofunctor!similarity with functor}

On the surface, functors and retrofunctors have much in common: both send objects to objects and morphisms to morphisms in a way that preserves domains and codomains as well as identities and composites.
The main difference is that functors send morphisms \emph{forward}, while retrofunctors send morphisms \emph{backward}.
As we work with retrofunctors, it will be helpful to remember the following:
\slogan{
    A retrofunctor $F$ goes forward on objects and backward on morphisms.\\
    Codomains are objects, so $F$ preserves them going forward.\\
    Identities and composites are morphisms, so $F$ preserves them going backward.
}

Before we explore just how different functors and retrofunctors can be, let us note a few more similarities that these two kinds of maps between categories share.
For example, retrofunctors, like functors, preserve isomorphisms.

\begin{proposition}\label{prop.retrofunctors_isos}
Let $F\colon\cat{C}\cof\cat{D}$ be a retrofunctor, $c\in\cat{C}$ be an object, and $g\colon Fc\to\_$ be an isomorphism in $\cat{D}$.
Then $F^\sharp_c g$ is also an isomorphism in $\cat{C}$.
\end{proposition}
\begin{proof}
Let $c'\coloneqq\cod F^\sharp_c g$, so that $Fc'=\cod g$ by \eqref{eqn.pres_cod}, and let $g^{-1}\colon Fc'\to Fc$ be the inverse of $g$.
Then
\begin{align*}
	\id_c&=
	F^\sharp_c\,\id_{Fc} \tag*{\eqref{eqn.pres_id}}\\&=
	F^\sharp_c\left(g\then g^{-1}\right)\\&=
	F^\sharp_c g\then F^\sharp_{c'}\left(g^{-1}\right), \tag*{\eqref{eqn.pres_comp}}
\end{align*}
so in particular $c=\cod\id_c=\cod F^\sharp_{c'}\left(g^{-1}\right)$, and
\begin{align*}
	\id_{c'}&=
	F^\sharp_{c'}\,\id_{Fc'} \tag*{\eqref{eqn.pres_id}}\\&=
	F^\sharp_{c'}\left(g^{-1}\then g\right)\\&=
	F^\sharp_{c'}\left(g^{-1}\right)\then F^\sharp_{\cod F^\sharp_{c'}\left(g^{-1}\right)}g \tag*{\eqref{eqn.pres_comp}}\\&=
	F^\sharp_{c'}\left(g^{-1}\right)\then F^\sharp_c g.
\end{align*}
Hence $F^\sharp_c g$ and $F^\sharp_{c'}\left(g^{-1}\right)$ are inverses, and the result follows.
\end{proof}

Moreover, isomorphisms in $\smcat$ correspond to isomorphisms in $\catsharp$.
\index{$\catsharp$!isomorphisms in}

\index{retrofunctor!isomorphism}

\begin{exercise}
We've justified the ``isomorphism-preserving'' part of \cref{thm.ahman_uustalu} implicitly, but let's make it explicit.

Recall that two categories $\cat{C}$ and $\cat{D}$ are isomorphic in $\smcat$ if there exist functors $F\colon\cat{C}\to\cat{D}$ and $G\colon\cat{D}\to\cat{C}$ that are mutually inverse, i.e.\ $F\then G$ and $G\then F$ are identity functors on $\cat{C}$ and $\cat{D}$.
Similarly, $\cat{C}$ and $\cat{D}$ are isomorphic in $\catsharp$ if there exist retrofunctors $H\colon\cat{C}\cof\cat{D}$ and $K\colon\cat{D}\cof\cat{C}$ that are mutually inverse, i.e.\ $H\then K$ and $K\then H$ are identity retrofunctors on $\cat{C}$ and $\cat{D}$.
Show that $\cat{C}$ and $\cat{D}$ are isomorphic in $\smcat$ if and only if they are isomorphic in $\catsharp$.
\begin{solution}
We want to show that categories $\cat{C}$ and $\cat{D}$ are isomorphic in $\smcat$ if and only if they are isomorphic in $\catsharp$.
First, assume that $\cat{C}$ and $\cat{D}$ are isomorphic in $\smcat$, so that there exist mutually inverse functors $F\colon\cat{C}\to\cat{D}$ and $G\colon\cat{D}\to\cat{C}$.
Then we can define a retrofunctor $H\colon\cat{C}\cof\cat{D}$ such that for each $c\in\cat{C}$ we have $Hc\coloneqq Fc\in\cat{D}$, and for each $g\in\cat{D}[Fc]$ we have $H^\sharp_c g\coloneqq Gg\in\cat{C}[GFc]=\cat{C}[c]$.
We can verify that $H$ really is a retrofunctor: it preserves identities and composition because $G$ does, and it preserves codomains because $H\cod H^\sharp_c g=F\cod Gg=FG\cod g=\cod g$.
Analogously, we can define a retrofunctor $K\colon\cat{D}\cof\cat{C}$ with $Kd\coloneqq Gd$ and $K^\sharp_d f\coloneqq Ff$ for each $d\in\cat{D}$ and $f\in\cat{C}[Gd]$.
Then $H\then K$ is equal to $F\then G$ both on objects and on morphisms, so it is the identity retrofunctor on $\cat{C}$; analogously, $K\then H$ is the identity retrofunctor on $\cat{D}$.
Thus $\cat{C}$ and $\cat{D}$ are isomorphic in $\catsharp$.

Conversely, assume $\cat{C}$ and $\cat{D}$ are isomorphic in $\catsharp$, so that there exist mutually inverse retrofunctors $H\colon\cat{C}\cof\cat{D}$ and $K\colon\cat{D}\cof\cat{C}$.
Given objects $c,c'$ and a morphism $f\colon c\to c'$ in $\cat{C}$, we have that $KHc=c$, so $K^\sharp_{Hc}$ is a function $\cat{C}[c]\to\cat{D}[Hc]$.
In particular, $K^\sharp_{Hc}f$ is a morphism in $\cat{D}$ whose domain is $Hc$ and whose codomain satisfies $K\cod K^\sharp_{Hc}f=\cod f=c'$, and thus $\cod K^\sharp_{Hc}f=Hc'$.
Hence we can define a functor $F\colon\cat{C}\to\cat{D}$ such that for each $c\in\cat{C}$ we have $Fc\coloneqq Hc\in\cat{D}$, and for each morphism $f\colon c\to c'$ in $\cat{C}$ we have $Ff\coloneqq K^\sharp_{Hc}f\colon Hc\to Hc'$.
Functoriality follows from the fact that $K$ preserves identities and composition.
Analogously, we can define a functor $G\colon\cat{D}\to\cat{C}$ with $Gd\coloneqq Kd$ for each $d\in\cat{D}$ and $Gg\coloneqq H^\sharp_{Kd}g$ for each $g\colon d\to d'$ in $\cat{D}$.
Then $F\then G$ is equal to $H\then K$ on objects; on morphisms, it sends each $f\colon c\to c'$ in $\cat{C}$ to $H^\sharp_{KHc}\left(K^\sharp_{Hc}f\right)=H^\sharp_c\left(K^\sharp_{Hc}f\right)=(K\then H)^\sharp_cf=f$ itself.
So $F\then G$ is the identity retrofunctor on $\cat{C}$.
Analogously, $G\then F$ is the identity functor on $\cat{D}$.
Thus $\cat{C}$ and $\cat{D}$ are isomorphic in $\smcat$.
\end{solution}
\end{exercise}
But while isomorphisms in $\catsharp$ are the same as those in $\smcat$, the non-isomorphisms can be very different.
We'll see this in the examples to come.

%---- Subsection ----%
\subsection{Examples of retrofunctors}
\index{retrofunctor!examples of|(}

From a realm where functors reign supreme, the back-and-forth behavior of retrofunctors can seem foreign and  counterintuitive.
Whereas a functor $\cat{C}\to\cat{D}$ can be thought of as a \emph{diagram}---a picture in the shape of $\cat{C}$, drawn with the objects and arrows of $\cat{D}$---retrofunctors are much more like the dynamical systems of \cref{ch.poly.dyn_sys}.\footnote{In fact, we will see in \cref{ch.comon.cofree} that retrofunctors generalize our dynamical systems.}

That is, a retrofunctor $F\colon\cat{C}\cof\cat{D}$ is a way of interacting with the states (objects) and transitions (morphisms) within $\cat{C}$ by way of $\cat{D}$.
Imagine the retrofunctor as a box, with $\cat{C}$ on the inside and $\cat{D}$ on the outside.
Some $c\in\cat{C}$ may be the current state inside the box, but all anyone outside the box can see is the object $Fc\in\cat{D}$ that the box chooses to display in lieu of $c$.
Still, any transition $g$ out of $Fc$ can be selected from the outside; the box guarantees that whatever $c$ is on the inside, there is a corresponding transition $F^\sharp_c g$ out of that $c$.
As $g$ is followed from $Fc$ to $\cod g$ on the outside, $F^\sharp_c g$ is followed from $c$ to $\cod F^\sharp_c g$ on the inside.
But codomain preservation guarantees that the new state $\cod g$ on the outside is equal to what the box would want to display in lieu of the new state $\cod F^\sharp_c g$ on the inside, as $\cod g=F\cod F^\sharp_c g$,
Then the process repeats in a manner compatible with identities and composition.

Here we give a variety of examples of retrofunctors to get a better handle on them.
Often we will denote a category by its carrier when its comonoid structure can be inferred from context, and $\cat{C}$ will be a category throughout with carrier $\car{c}$.

\subsubsection{Retrofunctors to preorders}
\index{retrofunctor!to preorders}

Given a retrofunctor from $\cat{C}$ to a preorder $\cat{P}$, we can think of $\cat{P}$ as providing a simplified model or abstraction of the states and transitions possible in $\cat{C}$, picking canonical transitions in $\cat{C}$ along the way to exhibit the model.
While the transitions in a general category may be more complex, all that a preorder tells you is whether you can get from one state to another or not.
Let's see some examples.

\index{retrofunctor!to discrete categories}

\begin{example}[Retrofunctors to discrete categories] \label{ex.cof_to_discrete}
The discrete category on a set $S$ is the category with objects in $S$ and only identity morphisms; its carrier is $S\yon$.\index{category!discrete}
So a retrofunctor $F\colon\cat C\cof S\yon$ is completely determined by its behavior on objects: to preserve identities, it can only send the morphisms in $S\yon$ back to the identity morphisms in $\cat C$.
We can identify $F$ with a function $\Ob\cat C\to S$, assigning each state in $\cat C$ a label in $S$ without revealing anything about the transitions between them.
\end{example}

\begin{exercise}\label{exc.yon_comonoid}
\begin{enumerate}
	\item Show that $\yon$ has a unique comonoid structure.
	\item Show that $\yon$ with its comonoid structure is terminal in $\catsharp$.
	\item Explain why $\yon$ is terminal using the language of states and transitions.
\qedhere
\end{enumerate}
\begin{solution}
\begin{enumerate}
    \item We actually already showed that $\yon$ has a unique comonoid structure, corresponding to the category with $1$ object and no nonidentity morphisms (which we will also denote by $\yon$), in \cref{exc.not_state_cat_but_same_carrier}, for the case of $S\coloneqq\1$.
    \item For any category $\cat{C}$, there is a unique retrofunctor $\cat{C}\cof\yon$: it sends every object in $\cat{C}$ to the only object in $\yon$, and it sends the only morphism in $\yon$, an identity morphism, to each identity morphism in $\cat{C}$.
    \item By \cref{ex.cof_to_discrete}, a retrofunctor from a category $\cat C$ to the discrete category on $S\in\smset$ is a way of assigning each state in $\cat C$ a label in $S$.
    In this case, $\yon$ is the discrete category on $\1$, so there is only $1$ label to choose from; hence there is always just $1$ way to assign the labels.
\end{enumerate}\index{category!discrete}
\end{solution}
\end{exercise}

\begin{example}[Retrofunctors to the walking arrow] \label{ex.cof_to_walking_arrow}
Consider a retrofunctor $F\colon\cat{C}\cof\cat{A}$, where $\cat{A}$ is the walking arrow category
\begin{center}
    $\cat{A}\coloneqq\:$\boxCD{examplecolor}{
    $s\Too{a}t$}
\end{center}
from \cref{ex.walking_arrow_cat}.
On objects, $F$ is a function $\Ob\cat{C}\to\{s,t\}$, so each object of $\cat{C}$ lies in either $\cat{C}_s\coloneqq F^{-1}s$ or $\cat{C}_t\coloneqq F^{-1}t$ (but not both).
Then on morphisms, preservation of identities determines where $F^\sharp$ sends $\id_s$ and $\id_t$, while preservation of codomains ensures that for each $c\in\cat{C}_s$, the morphism $F^\sharp_c a\colon a\to\_$ that $F^\sharp$ sends $a$ back to must satisfy
\[
    F\cod F^\sharp_c a=\cod a=t
\]
and thus $\cod F^\sharp_c a\in\cat{C}_t$.
In particular, for every object $c\in\cat{C}$ that $F$ sends to $s$, there must be at least one morphism from $c$ to an object that $F$ sends to $t$, so that one of those morphisms can be $F^\sharp_c a$.
As there are no nontrivial composites in $\cat{A}$, the retrofunctor $F$ automatically preserves composites.

In summary, a retrofunctor $F\colon\cat{C}\cof\cat{A}$ divides the objects of $\cat{C}$ between $\cat{C}_s$ and $\cat{C}_t$ and fixes a morphism from each object in $\cat{C}_s$ to some object in $\cat{C}_t$.
We can think of $F$ as separating the states of $\cat{C}$ into source states and target states, modeled by the $s$ state and the $t$ state in $\cat{A}$, respectively; then every source state is assigned a target state and a way of getting to that target state via a transition in $\cat{C}$.
\end{example}

Given a retrofunctor $F\colon\cat{C}\cof\cat{D}$ and an object $d\in\cat{D}$, we will continue to use the notation $\cat{C}_d\coloneqq F^{-1}d$ to denote the set of objects in $\cat{C}$ that $F$ sends to $d$.

\begin{exercise}
Let $F\colon\cat{C}\cof\cat{A}$ be a retrofunctor from $\cat{C}$ to the walking arrow category $\cat{A}$, as in \cref{ex.cof_to_walking_arrow}.
If $Fc=s$ for all $c\in\cat{C}$, what can we say about $\cat{C}$?
\begin{solution}
Given a retrofunctor $F\colon\cat{C}\cof\cat{A}$, where $\cat{A}$ is the walking arrow category as in \cref{ex.cof_to_walking_arrow}, assume $Fc=s$ for all $c\in\cat{C}$.
By \cref{ex.cof_to_walking_arrow}, if $Fc=s$ for $c\in\cat{C}$, then there must be a morphism from $c$ to an object in $\cat{C}$ that $F$ sends to $t$ for $a\colon s\to t$ in $\cat{A}$ to be sent back to via $F^\sharp$.
But there are no objects in $\cat{C}$ that $F$ sends to $t$.
So the only way such a retrofunctor could be defined is if there are no objects in $\cat{C}$ that it sends to $s$, either: we conclude that $\cat{C}$ is the empty category.
\end{solution}
\end{exercise}

\begin{exercise}
\begin{enumerate}
	\item Recall the star-shaped category ${\yon}^{\ord{n}+\1}+\ord{n}\yon$ from \cref{ex.star_shaped}. Describe retrofunctors to it.
	\item Describe retrofunctors to the preorder $(\nn,\leq)$, viewed as a category: its objects are natural numbers, and there is a morphism $m\to n$ if and only if $m\leq n$.
	\item Describe retrofunctors to the preorder $(\nn,\geq)$: its objects are natural numbers, and there is a morphism $n\to m$ if and only if $n\geq m$.
\qedhere
\end{enumerate}
\begin{solution}
The star-shaped category has $n+1$ objects, one of which is ``central'' in the sense that it maps uniquely to every object; the other objects have only identity maps.
\begin{enumerate}
    \item A retrofunctor $\cat{C}\cof\yon^{\ord{n}+\1}+\ord{n}\yon$ comprises an assignment of a label to each object in $\cat{C}$: either it assigns ``center'' or it assigns an element of $\{1,2,\ldots,n\}$. If it assigns ``center'', the object is equipped with $n$-many morphisms, with the $i$th one having as its codomain an object labeled $i$.
    \item A retrofunctor $\cat{C}\to(\nn,\leq)$ comprises an assignment of a natural number label to each object in $\cat{C}$, as well as a choice of morphism in $\cat{C}$ from each object labeled $n$ to some object labeled $n+1$.
    \item A retrofunctor $\cat{C}\to(\nn,\leq)$ comprises an assignment of a natural number label to each object in $\cat{C}$, as well as a choice of morphism in $\cat{C}$ from each object labeled $n+1$ to some object labeled $n$.
\end{enumerate}
\end{solution}
\end{exercise}

\begin{example}[Retrofunctors to the walking commutative square] \label{ex.cof_to_comm_sq}
Consider a retrofunctor $F\colon\cat{C}\cof\cat{C\!S}$, where $\cat{C\!S}$ is the \emph{walking commutative square category}
\begin{center}
    $\cat{C\!S}\coloneqq\:$\boxCD{examplecolor}{\[
    \begin{tikzcd}[ampersand replacement=\&]
        w\ar[d,"f"']\ar[r,"h"] \& y\ar[d,"k"] \\
        x\ar[r,"g"'] \& z
    \end{tikzcd}
    \]
    $f\then g=h\then k$}
\end{center}
On objects, $F$ is a function $\Ob\cat{C}\to\{w,x,y,z\}$, so each object of $\cat{C}$ lies in exactly one of $\cat{C}_w,\cat{C}_x,\cat{C}_y,$ and $\cat{C}_z$.
Then on morphisms, out of every object $X\in\cat{C}_x$ there is a morphism $F^\sharp_Xg\colon X\to\_$ to an object in $\cat{C}_z$, and out of every object $Y\in\cat{C}_y$ there is a morphism $F^\sharp_Yk\colon Y\to\_$ also to an object in $\cat{C}_z$.
Finally, out of every object $W\in\cat{C}_w$ there is a morphism $F^\sharp_Wf\colon W\to X_W$ to an object $X_W\in\cat{C}_x$ and a morphism $F^\sharp_Wh\colon W\to Y_W$ to an object in $Y_W\in\cat{C}_y$.
As $F$ preserves composites, these must all then satisfy
\[
    F^\sharp_Wf\then F^\sharp_{X_W}g=F^\sharp(f\then g)=F^\sharp(h\then k)=F^\sharp_Wh\then F^\sharp_{Y_W}k;
\]
in particular, $F^\sharp_{X_W}g$ and $F^\sharp_{Y_W}k$ must share a common codomain $Z_W\in\cat{C}_z$, yielding the following commutative square in $\cat{C}$:
\[
\begin{tikzcd}[ampersand replacement=\&]
    W\ar[d,"F^\sharp_Wf"']\ar[r,"F^\sharp_Wh"] \& Y_W\ar[d,"F^\sharp_{Y_W}k"] \\
    X_W\ar[r,"F^\sharp_{X_W}g"'] \& Z_W.
\end{tikzcd}
\]
\end{example}

\begin{exercise}
Let $\cat{A}$ denote the walking arrow category, as in \cref{ex.cof_to_walking_arrow}, and let $\cat{C\!S}$ denote the walking commutative square category, as in \cref{ex.cof_to_comm_sq}.
\begin{enumerate}
    \item List the retrofunctors $\cat{C\!S}\cof\cat{A}$.
    \item List the retrofunctors $\cat{A}\cof\cat{C\!S}$.\qedhere
\end{enumerate}
\begin{solution}
\begin{enumerate}
    \item Every object $c\in \cat{C\!S}$ in the commutative square must be labeled $s$ or $t$.
    If $c$ is labeled $t$, there are no further restrictions, since the only arrow emanating from $t$ is the identity.
    If $c$ is labeled $s$, then we must choose an outgoing arrow from $c$ to an object labeled $t$. So the possible retrofunctors are
\[
    \begin{tikzcd}[sep=10pt]
        t&t\\
	t&t
    \end{tikzcd}
	\quad\;
    \begin{tikzcd}[sep=10pt]
        s\ar[r]&t\\
	t&t
    \end{tikzcd}
	\quad\;
    \begin{tikzcd}[sep=10pt]
        s\ar[d]&t\\
	t&t
    \end{tikzcd}
	\quad\;
    \begin{tikzcd}[sep=10pt]
        s\ar[dr]&t\\
	t&t
    \end{tikzcd}
	\quad\;
    \begin{tikzcd}[sep=10pt]
        t&s\ar[d]\\
	t&t
    \end{tikzcd}
	\quad\;
    \begin{tikzcd}[sep=10pt]
        t&t\\
	s\ar[r]&t
    \end{tikzcd}
\]
\[
    \begin{tikzcd}[sep=10pt]
        s\ar[d]&s\ar[d]\\
	t&t
    \end{tikzcd}
	\quad\;
    \begin{tikzcd}[sep=10pt]
        s\ar[dr]&s\ar[d]\\
	t&t
    \end{tikzcd}
	\quad\;
    \begin{tikzcd}[sep=10pt]
        s\ar[r]&t\\
        s\ar[r]&t
    \end{tikzcd}
	\quad\;
    \begin{tikzcd}[sep=10pt]
        s\ar[dr]&t\\
	s\ar[r]&t
    \end{tikzcd}
        \quad\;
    \begin{tikzcd}[sep=10pt]
        t&s\ar[d]\\
        s\ar[r]&t
    \end{tikzcd}
        \quad\;
    \begin{tikzcd}[sep=10pt]
        s\ar[dr]&s\ar[d]\\
        s\ar[r]&t
    \end{tikzcd}
\]
    \item Every object $a\in\cat{A}$ in the walking arrow must be labeled $w,x,y,$ or $z$.
    If $a$ is labeled $z$, there are no further restrictions.
    If it is labeled $x$ or $y$, then $a$ must have a map to an element labeled $z$; in particular this implies that $a$ must be the source object $s\in\cat{A}$.
    Finally, $a$ cannot be labeled $w$ because then $a$ would need too many outgoing arrows.
    So the possible retrofunctors are
\[
    \begin{tikzcd}[sep=15pt]
        z&z
    \end{tikzcd}
	\qquad\;
    \begin{tikzcd}[sep=15pt]
        x\ar[r]&z
    \end{tikzcd}
	\qquad\;
    \begin{tikzcd}[sep=15pt]
        y\ar[r]&z
    \end{tikzcd}
	\qquad\;
\]
\end{enumerate}
\end{solution}
\end{exercise}

\begin{exercise}
\begin{enumerate}
	\item What does a retrofunctor from $\yon$ to a poset represent?
	\item Consider the chain poset $[n]\cong\sum_{i=0}^n\yon^{i+1}$.
	How many retrofunctors are there from $[m]\to[n]$ for all $m\leq n$?
\qedhere
\end{enumerate}
\begin{solution}
\begin{enumerate}
    \item Let $\cat{P}$ be a poset. A retrofunctor $\yon\cof\cat{P}$ represents an maximal element $p\in\cat{P}$. Indeed, if there were some $p\leq p'$ with $p'\neq p$, then the retrofunctor would have to send the map $p\to p'$ to some morphism in $\yon$ whose codomain is sent to $p'$; this is impossible. But if $p$ is maximal, then there is no obstruction to sending the unique object of $\yon$ to $p$.
   \item Let $m=\fbox{$\bullet^0\to\cdots\to\bullet^{m-2}\to\bullet^{m-1}$}$. By the same reasoning as above, a map $[m]\to[n]$ must send the final object of $[m]$ to the final object of $[n]$. It can send the penultimate object $m-2$ to either the penultimate object $n-2$ or to the final object $n-1$ in $[n]$. Repeating in this way, we see that there $2^m$ many retrofunctors. For example, the retrofunctors $[2]\to[2]$ are those labeled by $(2,2,2)$, $(1,2,2)$, $(1,1,2)$, and $(0,1,2)$.
 \end{enumerate}
\end{solution}
\end{exercise}

\index{retrofunctor!to monoids}

\subsubsection{Retrofunctors to monoids}

When a monoid $(M,e,*)$, viewed as a $1$-object category $\yon^M$, is the codomain of a retrofunctor $\cat{C}\cof\yon^M$, it plays the role of a joystick: an ``input device'' that ``reports its\dots direction to the device it is controlling.''\index{control}%
\footnote{Description from Wikipedia.} % https://en.wikipedia.org/wiki/Joystick
Like a joystick, $\yon^M$ stays in one ``place''---a single state---but has a number of directions it can take that are reported back to $\cat{C}$, controlling the way it moves through its transitions.
As we string together a sequence of directions in $M$, we chart a course through the transitions of $\cat C$.
We make this analogy concrete in the following examples.

\begin{example}[Arrow fields]\label{ex.arrow_field}\index{arrow field}
Consider the monoid $(\nn,0,+)$ viewed as a category $\yon^\nn$.
Retrofunctors $\cat{C}\cof\yon^\nn$ have been called \emph{admissible sections} \cite{aguiar1997internal}.
We prefer to call them \emph{arrow fields} (on $\cat{C}$), for they turn out to resemble vector fields---but with arrows in $\cat{C}$ instead of vectors.%
\tablefootnote{After all, a vector field is a section of a vector bundle.\index{bundle!vector}
Our arrow fields will be sections of $\cat{C}$'s carrier, viewed as a bundle of directions over positions.}
We'll have more to say about these in \cref{thm.catsharp_to_mon}, but our goal here is simply to unpack the definition.

To specify a retrofunctor $A\colon\cat{C}\cof\yon^\nn$, we first say what it does on objects, but this is already decided: there is only one object in $\yon^\nn$, so every object of $\cat{C}$ is sent to it.
This also means that codomains are automatically preserved.
So as will be the case for all retrofunctors to monoids, $A$ is characterized by its behavior on morphisms: for each object $c\in\cat{C}$, the retrofunctor assigns each $n\in\nn$ a morphism $A^\sharp_c n$ of $\cat{C}$ emanating from $c$.
That's a lot of data, but we still have two retrofunctor laws to pare it down:
\[
    A^\sharp_c0=\id_c
        \qqand
    A^\sharp_c(m+n)=A^\sharp_c m\then A^\sharp_{\cod A^\sharp_c m}n.
\]
Then for each $c\in\cat{C}$, since every $n\in\nn$ is a sum of 1's, the morphism $A^\sharp_i n$ can be decomposed into $n$ copies of $A^\sharp_{c_j}1$ for objects $c_0\coloneqq c,c_1,\ldots,c_n\in\cat{C}$, as follows:
\begin{equation} \label{eqn.arrow_field_composite}
    c=c_0\To{A^\sharp_{c_0}1}c_1\To{A^\sharp_{c_1}1}\cdots\To{A^\sharp_{c_{n-1}}1}c_n.
\end{equation}
Here each $c_{j+1}\coloneqq\cod A^\sharp_{c_j}1$.

Thus an arrow field of $\cat{C}$ is given by independently choosing a morphism emanting from each $c\in\cat{C}$ to be $A^\sharp_c1$: an arrow (morphism) out of each object, like how a vector field has a vector out of each point.
Indeed, any such choice uniquely determines the retrofunctor $A\colon\cat{C}\cof\yon^\nn$: as every object is assigned an arrow coming out of it, we can follow the arrow out of $c$ to an object $c_1$, then following the arrow out of $c_1$ to an object $c_2$, and so on until we have followed the $n$ arrows in \eqref{eqn.arrow_field_composite}, which then compose to yield $A^\sharp_cn$.
With our joystick analogy, a single flick sends us from a state $c\in\cat{C}$ along its assigned arrow $A^\sharp_c1$, while $n$ flicks send us through $n$ arrows along the arrow field.

Here is an example of an arrow field on the product of preorders $(\nn,\leq)\times(\nn,\leq)$:
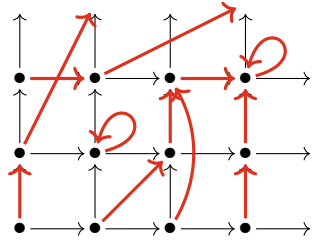
\[
\begin{tikzpicture}[shorten <=4pt, shorten >=4pt, tips=proper]
	\foreach \i in {0,1,2,3}
	{
		\foreach \j in {0,1,2}
		{
			\node (\i\j) at (\i,\j) {$\bullet$};
			\draw[->] (\i,\j) -- (\i+1,\j);
			\draw[->] (\i,\j) -- (\i,\j+1);
		};
	};
	\begin{scope}[my-red, very thick]
		\draw[->] (0,0) -- (0,1);
		\draw[->] (0,1) -- (1,3);
		\draw[->] (0,2) -- (1,2);
		\draw[->] (1,0) -- (2,1);
		\draw[->] (1,1) edge[in=75, out=15, distance=10mm] (1,1);
		\draw[->] (1,2) -- (3,3);
		\draw[->] (2,0) to[bend right] (2,2);
		\draw[->] (2,1) -- (2,2);
		\draw[->] (2,2) -- (3,2);
		\draw[->] (3,0) -- (3,1);
		\draw[->] (3,1) -- (3,2);
		\draw[->] (3,2) edge[in=75, out=15, distance=10mm] (3,2);
	\end{scope}
\end{tikzpicture}
\]
Every object has been assigned an emanating morphism drawn in red, but there need not be any rhyme or reason to our choice.
\end{example}

\index{retrofunctor!arrow field}

\begin{exercise}
How many arrow fields on the category \fbox{$\bullet\to\bullet$} are there?
\begin{solution}
We seek the number of arrow fields on the category \fbox{$\bullet\to\bullet$}.
There are $2$ choices of morphisms emanating from the object on the left, and $1$ choice of morphism emanating from the object on the right, for a total of $2 \cdot 1 = 2$ arrow fields.
\end{solution}
\end{exercise}

% Show that arrow fields can be composed
% maybe just bring all the arrow field stuff up here!

We will see later in \cref{prop.traj_mon_poly} that the arrow fields on a category form a monoid, and that this operation $\catsharp\to\Cat{Mon}\op$ is functorial and in fact an adjoint.

\begin{exercise}
Consider the monoid of integers $(\zz,0,+)$ as a $1$-object category $\yon^\zz$, and let $\yon^\nn$ be the monoid of natural numbers $(\nn,0,+)$ viewed as a $1$-object category as above.
\begin{enumerate}
	\item Describe the data of a retrofunctor $\cat{C}\cof\yon^\zz$.
	\item What would you say is the canonical retrofunctor $\yon^\zz\cof\yon^\nn$?
\qedhere
\end{enumerate}
\begin{solution}
\begin{enumerate}
    \item A retrofunctor $\cat{C}\cof\yon^\zz$ assigns to each object $c\in\cat{C}$ and integer $n\in\zz$ an emanating arrow $c.n\in\cat{C}$ to some other object, with the property that $c.0=c$ and $c.n.n'=c.(n+n')$. But this is overkill. Indeed, it is enough to assign the $c.1$ arrow and to check that for every object $c'$ there exists a unique object $c$ with $\cod(c.1)=c'$.
    \item We seek a canonical retrofunctor $\yon^\zz\cof\yon^\nn$.
    The canonical inclusion $i\colon\nn\inj\zz$ gives rise to a lens $\iota$ from $\yon^\zz$ to $\yon^\nn$, whose sole on-directions function $\iota^\sharp\colon\nn\inj\zz$ coincides with $i$.
    We verify that $\iota$ is a retrofunctor: it preserves identities, as $\iota^\sharp0=0$; it automatically preserves codomains; and it preserves composites, given by addition in either monoid, as $\iota^\sharp(m+n)=m+n=\iota^\sharp(m)+\iota^\sharp(n)$.
\end{enumerate}
\end{solution}
\end{exercise}

\begin{exercise}
\begin{enumerate}
	\item Suppose that $M,N$ are monoids (each is a category with one object).
	Are retrofunctors between them related to monoid homomorphisms? If so, how?
	\item Suppose $\cat{C}$ and $\cat{D}$ are categories and $F\colon\cat{C}\cof\cat{D}$ is a retrofunctor.
	Does there necessarily exist a retrofunctor $\cat{C}\op\cof\cat{D}\op$ that acts the same as $F$ on objects?
\qedhere
\end{enumerate}
\begin{solution}
\begin{enumerate}
    \item Retrofunctors $\yon^M\to\yon^N$ are the same as monoid homomorphisms $N\to M$. See \cref{prop.monoids_ff}.
    \item No! This is the weirdest thing about retrofunctors. For example, there is a unique retrofunctor $\yon\to\yon^\2+\yon$ from the walking object to the walking arrow, so if we reverse the arrows, there is no longer a retrofunctor that acts the same on objects.
\end{enumerate}
\end{solution}
\end{exercise}

\begin{exercise}[Monoid actions]\label{exc.monoid_action}
Recall from \cref{ex.monoid_action} that every monoid action $\alpha\colon S\times M\to S$, where $S$ is a set and $(M,e,*)$ is a monoid, gives rise to a category carried by $S\yon^M$.
Show that the projection $S\yon^M\to\yon^M$ is a retrofunctor.
\begin{solution}
We are given a monoid $(M,e,*)$, a set $S$, and an $M$-action $\alpha\colon M\times S\to S$. The category $S\yon^M$ has $e$ as the identity on each $s\in S$; the codomain of the map labeled $m$ emanating from $s$ is $\alpha(s,m)$, and composition is given by $*$. The projection $S\yon^M\to\yon^M$ sends the identity back to the identity, trivially preserves codomains, and also preserves composition, so it is a retrofunctor.
\end{solution}
\end{exercise}

\begin{example}\label{ex.BGEG}
Let $(G,e,*)$ be a group and $(\yon^G,\epsilon,\delta)$ the corresponding comonoid. There is a retrofunctor $G\yon^G\cof\yon^G$ given by
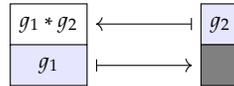
\[
\begin{tikzpicture}[polybox, mapstos]
	\node[poly, dom] (p) {$g_1*g_2$\at$g_1$};
	\node[poly, pure cod, right=of p] (q) {$g_2$\at\vphantom{$g_1$}};
	\draw (p_pos) to[first] (q_pos);
	\draw (q_dir) to[last] (p_dir);
\end{tikzpicture}
\]
To see this is a retrofunctor, we check that identities, codomains, and compositions are preserved. For any $g_1$, the identity $e$ is passed back to $g_1*e=g_1$, and this is the identity on $g_1$ in $G\yon^G$. Codomains are preserved because there is only one object in $\yon^G$. Composites are preserved because for any $g_2,g_3$, we have $g_1*(g_2*g_3)=(g_1*g_2)*g_3$.
\end{example}

\begin{exercise}\label{exc.BGEG}
Does the idea of \cref{ex.BGEG} work when $G$ is merely a monoid, or does something go subtly wrong somehow?
\begin{solution}
This works!
\end{solution}
\end{exercise}

\begin{proposition}\label{prop.monoids_ff}
There is a fully faithful functor $\Cat{Mon}\op\to\catsharp$, whose image consists of all categories whose carriers are representable.
\end{proposition}
\begin{proof}
Given a monoid $(M,e,*)$, we think of it as a category with one object; its carrier $\yon^M$ is representable. A retrofunctor between such categories carries no data in its on-objects part, and codomains are automatically preserved. Retrofunctors $\yon^M\to\yon^N$ simply carry elements of $N$ to elements of $M$, preserving identity and composition, exactly the description of monoid homomorphisms.
\end{proof}

\begin{proposition}\index{adjunction!between $\smset$ and $\catsharp$}
There is an adjunction
\[
\catsharp(\cat{C},A\yon)\cong\smset(\Ob\cat{C},A)
\]
for $\cat{C}\in\catsharp$ and $A\in\smset$.
\end{proposition}
\begin{proof}
In the solution to \cref{exc.linear_poly_cat}, we saw that a category is discrete iff its carrier is a linear polynomial: this occurs when the only arrow emanating from each object is its identity. Thus $A\yon$ corresponds to a discrete category. A retrofunctor from any category to a discrete category needs to say what happens on objects, but the rest of the data is determined because identities need to be sent back to identities. This is the content of the proposition.
\end{proof}\index{category!discrete}

\begin{exercise}[Continuous arrow fields]
Suppose we say that a \emph{continuous arrow field} on $\cat{C}$ is a retrofunctor $\cat{C}\cof\yon^\rr$, viewing $\yon^\rr$ as the monoid of real numbers with addition.

Describe continuous arrow fields in $\cat{C}$ using elementary terms, i.e.\ to someone who doesn't know what a retrofunctor is and isn't yet ready to learn.
\begin{solution}
A continuous arrow field on $\cat{C}$ assigns to each object $c\in\cat{C}$ and each real number $r\in\rr$ a morphism $c.r$ emanating from $c$. These have the property that $c.0$ is the identity on $c$ and that $(c.r).r'=c.(r+r')$. In other words, you can evolve $c$ forward or backward in time by any $r\in\rr$, and this works as expected.
\end{solution}
\end{exercise}

\index{Systems of ODEs}

\begin{example}[Systems of ODEs]
A system of ordinary differential equations (ODEs) in $n$ variables, e.g.
\begin{align*}
    \dot{x}_1 &= f_1(x_1, \ldots, x_n) \\
    \dot{x}_2 &= f_2(x_1, \ldots, x_n) \\
    & \; \; \; \vdots \\
    \dot{x}_n &= f_n(x_1, \ldots, x_n),
\end{align*}
can be understood as a vector field on $\rr^n$.
We are often interested in integrating this vector field to get flow lines, or integral curves.
In other words, for each $x=(x_1, \ldots, x_n)\in\rr^n$, viewed as a point, and each $t\in\rr$, viewed as a quantity of time, we can begin at $x$ and move along the vector field for time $t$, arriving at a new point $x^{+t}$. These satisfy the equations
\begin{equation} \label{eqn.retrofunctor_ode}
    x^{+0} = x \qqand x^{+t_1+t_2} = (x^{+t_1})^{+t_2}.
\end{equation}
Let's call such things \emph{differentiable dynamical systems} with time domain $(T, 0, +)$; above, we used $T\coloneqq\rr$, but any monoid will do.

\index{dynamical system!as retrofunctor}

Differential dynamical systems as defined above are retrofunctors $F \colon \rr^n \yon^{\rr^n} \cof \yon^T$.
In order to say this, we first need to say how both $\cat{C} := \rr^n \yon^{\rr^n}$ and $\yon^T$ are being considered as categories.
The category $\cat{C}$ has objects $\rr^n$, and for each object $x \in \rr^n$ and outgoing arrow $v \in \rr^n$, the codomain of $v$ is $x + v$; in other words, $v$ is a vector emanating from $x$.
The identity is $v = 0$, and composition is given by addition.
The category $\yon^T$ is the monoid $T$ considered as a category with one object, $\bullet$.

The retrofunctor assigns to every object $x \in \rr^n$ the unique object $F(x) = \bullet$, and to each element $t \in T$ the morphism $F^\sharp(x, t) = x^{+t} - x \in \rr^n$, which can be interpreted as a vector emanating from $x$.
Its codomain is $\cod F^\sharp(x, t) = x^{+t}$, and we will see that \eqref{eqn.retrofunctor_ode} ensures the retrofunctoriality properties.

The codomain law ii is vacuously true, since $\yon^T$ only has one object.
Law i follows because $F^\sharp(x, 0) = x^{+0} - x = 0$, and law iii follows as
\begin{align*}
    F^\sharp(x^{+t_1}, t_2) + F^\sharp(x, t_1) &= \\
    (x^{+t_1})^{+t_2} - x^{+t_1} + x^{+t_1} - x &= \\
    x^{+t_1 + t_2} - x &= \\
    F^\sharp(x, t_1+t_2).
\end{align*}
\end{example}

\subsubsection{Retrofunctors from state categories}
\index{state category|see{state system}}\index{state system!retrofunctor from|(}

By now we should be very familiar with lenses from state categories, which are our original dynamical systems.
A retrofunctor from a state category, then, is just a dynamical system that satisfies the retrofunctor laws.
It turns out that retrofunctors from state categories are particularly noteworthy: just as a polynomial comonoid $\cat{C}$ can be identified with a category, a retrofunctor out of $\cat{C}$ can be identified with a number of equivalent categorical constructions on $\cat{C}$, perhaps the most familiar being a functor $\cat{C}\to\smset$.
But these equivalences deserve their own subsection to examine in full; we'll defer them to \cref{sec.comon.sharp.cof.from_state}.
For now, let's look at some examples of retrofunctors out of state categories.

\begin{example}[Retrofunctors $S\yon^S\cof\cat{C}$ are $\cat{C}$-coalgebras]\index{coalgebra!of polynomial comonad}\index{polynomial comonoid!coalgebra of}

Recall from \cref{ex.coalgebras} that for a set $S$, lenses $S\yon^S\to p$ correspond to functions $S\to p(S)$ known as coalgebras for the functor $p$.
As a retrofunctor $S\yon^S\cof\cat{C}$ is just a special kind of lens from $S\yon^S$ to $\car{c}$, the carrier of $\cat{C}$, it should correspond to a special kind of coalgebra $S\to\car{c}(S)$ for the functor $\car{c}$.

% that plays nicely with the comonoid structure of $\cat{C}$, the function $S\to\car{c}(S)$ that plays nicely with the comonoid structure carried by $\car{c}$.

Taking $A=B=S\in\smset$ in \eqref{eqn.monomials_and_comp}, we find that there is a natural isomorphism between dynamical systems $S\yon^S\to p$ and functions $S\to p(S)$, also known as a \emph{coalgebra for the functor} $p$ or a $p$\emph{-coalgebra}.\tablefootnote{There are two versions of coalgebras we are interested in (and more that we are not) with distinct definitions: a \emph{coalgebra for a functor}, which is the version used here, and a \emph{coalgebra for a comonad}, which is a coalgebra for a functor with extra conditions that we will introduce later.}
\end{example}

\index{lens!very well-behaved|seealso{retrofunctor, between state categories}}\index{retrofunctor!between state categories}\index{functional programming!lenses in}

\begin{example}[Retrofunctors between state categories are very well-behaved lenses] \label{ex.very_well_behaved_lenses}
Our familiar state category on $S$ from \cref{ex.state_cat} is the category with objects in $S$ and exactly $1$ morphism between every pair of objects; when we label each morphism with its codomain, its carrier is $S\yon^S$, the identity of $s\in S$ is $s$, and (disregarding domains) $s\then s'=s'$ for composable $s,s'\in S$.

Then a retrofunctor $S\yon^S\cof T\yon^T$ between two state categories corresponds to what is known to functional programmers as a \emph{very well-behaved lens}.
We actually defined this way back in \cref{ex.lens_get_put}, where we called the on-objects (on-positions) function of such a retrofunctor $\lensget\colon S\to T$, and the on-morphisms (on-directions) function $\lensput\colon S\times T\to S$.%
\tablefootnote{More precisely, we are treating the on-morphism functions $T\to S$ for each $s\in S$ of a retrofunctor $S\yon^S\cof T\yon^T$ as a single function $S\times T\to S$.}
Then the retrofunctor laws are as follows:
\begin{enumerate}
    \item Preservation of identities \eqref{eqn.pres_id} becomes
    \[
        \lensput(s,\lensget(s))=s,
    \]
    for all $s\in S$, known as the \emph{get-put law} (named in diagrammatic order: we apply $\lensget$ before we apply $\lensput$).

    \item Preservation of codomains \eqref{eqn.pres_cod} becomes
    \[
        \lensget(\lensput(s,t))=t,
    \]
    for all $s\in S$ and $t\in T$, known as the \emph{put-get law}.

    \item Preservation of composition becomes
    \[
        \lensput(\lensput(s,t),t')=\lensput(s,t')
    \]
    for all $s\in S$ and $t,t'\in T$, known as the \emph{put-put law}.
\end{enumerate}

In fact, it turns out that these laws can be satisfied if and only if $\lensget$ is a product projection!
For example, if the cardinalities $|S|$ and $|T|$ of $S$ and $T$ are finite and $|S|$ is not divisible by $|T|$, then there are no retrofunctors $S\yon^S\cof T\yon^T$.
A stringent condition, no?
We'll explore it in  \cref{exc.how_many_vwbls} below.

Let's explore why retrofunctors between state categories are just product projections.
A product projection $A\times B\to A$ always has a second factor $B$; if every retrofunctor between state categories is a product projection, what is the second factor?
It turns out to be
\[
    U\coloneqq\{u\colon T\to S\mid \forall t,t'\in T,\;\lensget(u(t))=t\text{ and }\lensput(u(t),t')=u(t')\}.
\]
In other words, we will show that if $(\lensget,\lensput)$ defines a retrofunctor $S\yon^S\cof T\yon^T$, then there is a bijection $S\iso T\times U$ making $\lensget\colon S\to T$ a product projection.
We then prove the converse in \cref{exc.well_behaved_boring}.

Assume $(\lensget,\lensput)\colon S\yon^S\cof T\yon^T$ is a retrofunctor, so that it satisfies the enumerated laws.
First, we define a function $\alpha\colon S\to T\times U$ as follows.
Given $s\in S$, the function $u_s\colon T\to S$ defined by
\[
    u_s(t)=\lensput(s,t)
\]
lies in $U$: we check that it satisfies
\[
    \lensget(u_s(t))=\lensget(\lensput(s,t))=t
\]
by the put-get law for $t\in T$ and
\[
    \lensput(u_s(t),t')=\lensput(\lensput(s,t),t')=\lensput(s,t')=u_s(t')
\]
by the put-put law for $t,t'\in T$.
We can therefore define a function $\alpha\colon S\to T\times U$ by
\[
    \alpha(s)=\left(\lensget(s),u_s\right).
\]
In the other direction, we have a function $\beta\colon T\times U\to S$ given by
\[
    \beta(t,u)=u(t).
\]
The two functions $\alpha$ and $\beta$ are mutually inverse: $\alpha\then\beta\colon S\to S$ is the identity because
\[
    \beta(\alpha(s))=u_s(\lensget(s))=\lensput(s,\lensget(s))=s
\]
by the get-put law, while $\beta\then\alpha\colon T\times U\to T\times U$ is the identity because
\[
    \alpha(\beta(t,u))=\left(\lensget(u(t)),u_{u(t)}\right)=(t,t'\mapsto\lensput(u(t),t'))=(t,u),
\]
as $\lensget(u(t))=t$ and $\lensput(u(t),t')=u(t')$ by construction for $u\in U$.
Thus $S\iso T\times U$, and the product projection $S\To\alpha T\times U\to T$ sends $s\mapsto\lensget(s)$, as desired.

We have therefore shown that for every retrofunctor $S\yon^S\cof T\yon^T$, there exists a set $U$ for which $S\iso T\times U$ and the on-positions function $\lensget$ is the product projection $S\iso T\times U\to T$.
Notice that the on-directions function $\lensput$ can be uniquely recovered from the bijection $S\iso T\times U$ we constructed: it is determined by the functions $u_s\colon T\to S$ for $s\in S$, which in turn is determined by the second projection $T\times U\to U$.

More precisely, composing $\alpha\colon S\to T\times U$ with the projection to $U$ yields a map $S\to U$ sending $s\mapsto u_s$; then $\lensput(s,t)$ is given by $u_s(t)$.
Of course, a priori $U$ could just be a set---we may not know how to interpret its elements as a functions $T\to S$.
This is where $\beta\colon T\times U\to S$ comes in: we know $\beta(t,u_s)=u_s(t)$.
So $\lensput(s,t)$ must be $\beta(t,u_s)$.
\end{example}

\begin{exercise}\label{exc.well_behaved_boring}
Let $S,T,U$ be sets for which we have a bijection $S\iso T\times U$.
Show that there exists a unique retrofunctor $S\yon^S\cof T\yon^T$ whose on-positions function $S\iso T\times U\to T$ is given by the product projection.
\begin{solution}
Given a bijection of sets $S\iso T\times U$, we seek a unique retrofunctor $S\yon^S\cof T\yon^T$ whose on-positions function $\lensget$ is given by the product projection $S\iso T\times U\to T$.
Such a retrofunctor should also have an on-directions function $\lensput\colon S\times T\to S$ such that the three laws from \cref{ex.very_well_behaved_lenses} are satisfied.
With $S\iso T\times U$, such a function is uniquely determined by its components
\[
    S\times T\To\lensput S\To\lensget T \qqand S\times T\To\lensput S\To\pi U,
\]
where $\pi\colon S\iso T\times U\to U$ is the other product projection.
The left component is determined by the put-get law, which specifies the behavior of the composite $\lensput\then\lensget$: it should send $(s,t)\mapsto t$.
Identifying $S$ with $T\times U$, the get-put law for $s=(t,u)\in T\times U$ (so $\lensget(t,u)=t$) reads as
\[
    \lensput((t,u),t)=(t,u),
\]
while the put-put law reads as
\[
    \lensput(\lensput((t,u),t'),t'')=\lensput((t,u),t'')
\]
for $t',t''\in T$.
Applying $\lensget$ to both sides, we observe that the first coordinate (in $T$) of either side of each equation are automatically equal when the put-get law holds.
So we are really only concerned with the second coordinates of either side (in $U$); applying $\pi$ to both sides yields
\[
    \pi(\lensput((t,u),t))=u \qqand \pi(\lensput(\lensput((t,u),t'),t''))=\pi(\lensput((t,u),t''))
\]
\end{solution}
\end{exercise}

\begin{exercise}\label{exc.how_many_vwbls}
\begin{enumerate}
	\item Suppose $|S|=3$. How many retrofunctors are there $S\yon^S\to S\yon^S$?
	\item Suppose $|S|=4$ and $|T|=2$. How many retrofunctors are there $S\yon^S\cof T\yon^T$?
\qedhere
\end{enumerate}
\begin{solution}
\begin{enumerate}
    \item There are 6 product projections $\3\to \3$, namely the three automorphisms, so the answer is 6.
    \item There are 6 product projections $\4\to\2$, so the answer is 6.
\end{enumerate}
\end{solution}
\end{exercise}

\begin{example}
We have a state category $\car{c}(\1)\yon^{\car{c}(\1)}$ on the set of objects of $\cat{C}$.
Define a lens $\car{c}(\1)\yon^{\car{c}(\1)}\to\car{c}$ by
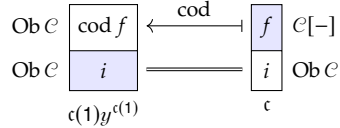
\[
\begin{tikzpicture}[polybox, mapstos]
	\node[poly, dom, "$\car{c}(\1)\yon^{\car{c}(\1)}$" below] (s) {$\cod f$\at$i$};
	    \node[left=0pt of s_pos] {$\Ob\cat{C}$};
        \node[left=0pt of s_dir] {$\Ob\cat{C}$};

	\node[poly, cod, right=of s, "$\car{c}$" below] (c) {$f\vphantom{d}$\at$i$};
	    \node[right=0pt of c_pos] {$\Ob\cat{C}$};
        \node[right=0pt of c_dir] {$\cat{C}[-]$};

	\draw[double, -] (s_pos) -- (c_pos);
	\draw (c_dir) -- node[above] {$\cod$} (s_dir);
\end{tikzpicture}
\]
sending each object $i\in\cat{C}$ to itself on positions and, at $i$, sending each morphism $f\colon i\to\_$ to its codomain $\cod f$ on directions.

This lens is a retrofunctor $\car{c}(\1)\yon^{\car{c}(\1)}\cof\cat{C}$ because it sends identities back to identities, codomains forward to codomains, and preserves composition (trivially, since each morphism in $c(\1)\yon^{c(\1)}$ is determined by its domain and codomain.
\end{example}

\begin{exercise}
Fix an object $i\in\car{c}(\1)=\Ob\cat{C}$.
Then we have a state category $\car{c}[i]\yon^{\car{c}[i]}$ on the set $\car{c}[i]=\cat{C}[i]$ of morphisms out of $i$ in $\cat{C}$.
Define a lens $\car{c}[i]\yon^{\car{c}[i]}\to\car{c}$ by
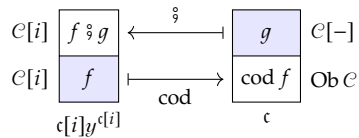
\[
\begin{tikzpicture}[polybox, mapstos]
	\node[poly, dom, "{$\car{c}[i]\yon^{\car{c}[i]}$}" below] (s) {$f\then g$\at$\vphantom{\cod}f$};
        \node[left=0pt of s_pos] {$\cat{C}[i]$};
        \node[left=0pt of s_dir] {$\cat{C}[i]$};

	\node[poly, cod, right=of s, "$\car{c}$" below] (c) {$\vphantom{f\then}g$\at$\cod f$};
	    \node[right=0pt of c_pos] {$\Ob\cat{C}$};
        \node[right=0pt of c_dir] {$\cat{C}[-]$};

	\draw (s_pos) -- node[below] {$\cod$} (c_pos);
	\draw (c_dir) -- node[above] {$\then$} (s_dir);
\end{tikzpicture}
\]
sending each morphism $f\colon i\to\_$ to its codomain on positions and, at $f$, sending each morphism $g\colon\cod f\to\_$ to the composite $f\then g\colon i\to\_$ on directions.
Is this lens a retrofunctor $\car{c}[i]\yon^{\car{c}[i]}\cof\cat{C}$?
\begin{solution}
Yes, it is a retrofunctor. One sees easily that it sends identities back to identities and composites back to composites. The codomain of $f\then g$ is simply $f\then g$ as an element of $\car{c}[i]$, and it is sent forward to $\cod(f\then g)$, so the map preserves codomains.
\end{solution}
\end{exercise}

We'll revisit retrofunctors from state categories in \cref{sec.comon.sharp.cof.from_state}.

\index{state system!retrofunctor from|)}

\subsubsection{Other retrofunctors}

\begin{example}[Objects aren't representable in $\catsharp$]\label{ex.rep_objects}
In the world of categories and the usual functors between them, the terminal category $\cat{T}\coloneqq\fbox{$\bullet$}$ with one object and one morphism \emph{represents objects}, in the sense that functors $\cat{T}\to\cat{C}$ naturally correspond to objects in $\cat{C}$.

Unfortunately, the same cannot be said for retrofunctors: we'll see in \cref{exc.rep_objects} that there does not exist a fixed category $\cat{U}$ for which retrofunctors $\cat{U}\cof\cat{C}$ are in bijection with objects in $\cat{C}$ for every category $\cat{C}$.

Retrofunctors $\cat{T}\cof\cat{C}$ are somewhat strange beasts: because they must preserve codomains, they can be identified with objects $c\in\cat{C}$ for which the codomain of every emanating morphism $c\to c'$ is $c'=c$ itself.
\end{example}

\begin{exercise}\label{exc.rep_objects}
We saw in \cref{exc.linear_poly_comon} that $\2\yon$ has a unique comonoid structure.
\begin{enumerate}
	\item Show that for any category $\cat{U}$, retrofunctors $\cat{U}\cof\2\yon$ are in bijection with the set $\2^{\Ob\cat{U}}$.
	\item Use the case of $\cat{C}\coloneqq\2\yon$ to show that if retrofunctors $\cat{U}\cof\cat{C}$ are always in bijection with objects in $\cat{C}$, then $\cat{U}$ must have exactly one object.
	\item Now use a different category $\cat{D}$ to show that if retrofunctors $\cat{U}\cof\cat{D}$ are in bijection with objects in $\cat{D}$, then $\cat{U}$ must have more than one object.
	Conclude that objects are not reprsentable in $\catsharp$ the way they are in $\smcat$.
	\item Is there a fixed category $\cat{V}$ for which retrofunctors $\cat{E}\cof\cat{V}$ are in bijection with objects in $\cat{E}$ for every category $\cat{E}$?
	If there is, find it; if there isn't, prove there isn't.
\qedhere
\end{enumerate}
\begin{solution}
\begin{enumerate}
    \item Every object in $\cat{U}$ must be labeled with either 1 or 2, but there are no other requirements.
    \item If retrofunctors $\cat{U}\to\cat{C}$ are always in bijection with objects in $\cat{C}$, then with $\cat{C}=\2\yon$, we have a bijection $\2\cong\2^{\Ob(\cat{U})}$ by part 1, so $\cat{U}$ has one object.
    \item Take $\cat{D}$ to be the walking arrow. Then if $\cat{U}$ has only one object, it cannot be sent to the source object of $\cat{D}$ because the emanating morphism would have nowhere to go. Hence the unique object of $\cat{U}$ must be sent to the target object of $\cat{D}$. But there is only one such retrofunctor, whereas there are two objects in $\cat{D}$. We conclude that objects are not reprsentable in $\catsharp$ the way they are in $\smcat$.
    \item No, there is no such $\cat{V}$. There is exactly one retrofunctor $0\cof\cat{V}$, but there are no objects of $0$.
\end{enumerate}
\end{solution}
\end{exercise}

\begin{example}\label{ex.cof_to_rr} %this is a retrofunctor a monoid action... section for those??
Consider the category $\rr\yon^\rr$, where the codomain of $r$ emanating from $x$ is $x+r$, identities are $0$, and composition is given by addition. What are retrofunctors into $\rr\yon^\rr$?

Let $\cat{C}$ be a category and $|\cdot|\colon\cat{C}\cof\rr\yon^\rr$ a retrofunctor. It assigns to every object $c$ both a real number $|c|\in\rr$ and a choice of emanating morphism $|c|^\sharp(r)\colon c\to c_r$ such that $|c|+r=|c_r|$. This assignment satisfies some laws. Namely we have $c_0=c$ and, given reals $r,s\in\rr$, we have $(c_r)_s=c_{r+s}$.
\end{example}

%  Interpret above. Possible names: $(\rr,0,+)$-action on the objects of $\cat{C}$, filtration, valuation

\begin{exercise}
How many retrofunctors
\[
    \fbox{$s\To{a}t$}\cof\fbox{$u\Tto{b}{c}v$}
\]
are there from the walking arrow category $\cat{A}$, drawn above on the left, to the walking parallel-arrows category $\cat{P\!A}$, drawn above on the right?
\begin{solution}
There are two: one of which sends $s\mapsto u$ and $t\mapsto v$ and the other of which sends both $s,t\mapsto v$.
\end{solution}
\end{exercise}

\begin{exercise}
\begin{enumerate}
	\item For any category $\cat{C}$ with carrier $\car{c}$, find a category with carrier $\car{c}\yon$.
	\item Show that your construction is functorial; i.e.\ assign each retrofunctor $\cat{C}\cof\cat{D}$ a retrofunctor $\car{c}\yon\cof\car{d}\yon$ in a way that preserves identities and composites.
	\item Is your functor a monad on $\catsharp$, a comonad on $\catsharp$, both, or neither?
\qedhere
\end{enumerate}
\begin{solution}
\begin{enumerate}
    \item Take the category $\cat{C}'$ that has the same objects as $\cat{C}$ and almost the same arrows, except that it has one more arrow $i_c$ from each object $c\in\cat{C}'$ to itself. This new arrow is the identity. The composites of all the old arrows in $\cat{C}'$ are exactly as they are in $\cat{C}$. The old identity $\id_c^\tn{old}$ is no longer an identity because $\id_c^\tn{old}\then i=\id_c^\tn{old}\neq i$.
    \item Given a retrofunctor $\varphi\colon\car{c}\cof\car{d}$, we get a retrofunctor $\car{c}\yon\cof\car{d}\yon$ that acts the same on objects and all the old arrows and sends the new identities back to the new identities. It preserves identities and composites going backward and codomains going forward, so it's a retrofunctor.
    \item It is a monad: there is an obvious unit map $\car{c}\to\car{c}\yon$ and a multiplication map $\car{c}\yon\yon\to\car{c}\yon$ that sends the new identity back to the newest identity.
\end{enumerate}
\end{solution}
\end{exercise}

\begin{exercise}
Suppose $\car{c},\car{d},\car{e}$ are polynomials, each with a comonoid structure, and that $f\colon\car{c}\to\car{d}$ and $g\colon\car{d}\to\car{e}$ are lenses.
\begin{enumerate}
	\item If $f$ and $f\then g$ are each retrofunctors, is $g$ automatically a retrofunctor?
	If so, sketch a proof; if not, sketch a counterexample.
	\item If $g$ and $f\then g$ are each retrofunctors, is $f$ automatically a retrofunctor?
	If so, sketch a proof; if not, sketch a counterexample.
\qedhere
\end{enumerate}
\begin{solution}
\begin{enumerate}
    \item Counterexample: take $\car{c}=0$. Then the unique lens $f\colon\car{c}\to\car{d}$ is a retrofunctor and so is $f\then g$ for any lens $g\colon\car{d}\to\car{e}$, but some lenses are not retrofunctors.
    \item Counterexample: take $\car{e}=\yon$. Then there is a unique retrofunctor $\epsilon\colon\car{c}\to\yon$. As long as $f$ is copointed, meaning it sends identities backward to identities, then $f\then g$ will be a retrofunctor. But not every copointed lens is a retrofunctor.
\end{enumerate}
\end{solution}
\end{exercise}

In the next chapter, we will delve deeper into the categorical structure and properties of $\catsharp$.
We'll encounter many more categories and retrofunctors along the way.
But first, we'll conclude this chapter with several alternative characterizations of retrofunctors out of state categories.

\index{retrofunctor!examples of|)}

\subsection[Retrofunctors from state categories]{Equivalent characterizations of retrofunctors from state categories} \label{sec.comon.sharp.cof.from_state}

Fix a category $\cat{C}$ throughout with polynomial carrier $\car{c}$.
How can view the data of a retrofunctor from a state category $S\yon^S\cof\cat{C}$?
This is actually a very natural categorical concept---we'll see some equivalent ways to express this data below, then state and prove even more equivalences next chapter when we have the machinery to do so.

\subsubsection{As coalgebras}
\index{coalgebra!as retrofunctor from state category}\index{functor!coalgebra for}

Recall from \cref{ex.coalgebras} that for a set $S$, lenses $S\yon^S\to p$ correspond to functions $S\to p(S)$ known as coalgebras for the functor $p$.
As a retrofunctor $S\yon^S\cof\cat{C}$ is just a special kind of lens from $S\yon^S$ to $\car{c}$, it should correspond to a special kind of coalgebra $S\to\car{c}(S)$ for the functor $\car{c}$.
Indeed, whenever $\car{c}$ carries a comonoid $\cat{C}$ with respect to the composition product (i.e.\ a comonad), there is a special notion of a $\cat{C}$-coalgebra (i.e.\ a coalgebra for the comonad $\cat{C}$), as follows.

\begin{definition}[Coalgebra for a polynomial comonoid]\label{def.coalgebra}
Let $\cat{C}=(\car{c},\epsilon,\delta)$ be a polynomial comonoid.
A \emph{$\cat{C}$-coalgebra} $(S,\alpha)$ is
\begin{itemize}
    \item a set $S$, called the \emph{carrier}, equipped with
    \item a function $\alpha\colon S\to\car{c}\tri S$,
\end{itemize}
such that the following diagrams, collectively known as the \emph{coalgebra laws}, commute:
\begin{equation} \label{eqn.coalg_laws}
\begin{tikzcd}
    S \ar[r, "\alpha"] \ar[dr, equal] &
    \car{c}\tri S \ar[d, "\epsilon\:\tri\:S"] \\
    & S
\end{tikzcd}
\hspace{.6in}
\begin{tikzcd}
    S \ar[d, "\alpha"] \ar[r, "\alpha"] &
    \car{c}\tri S \ar[d, "\delta\:\tri\:S"] \\
    \car{c}\tri S \ar[r, "\car{c}\:\tri\:\alpha"'] &
	\car{c}\tri\car{c}\tri S.
\end{tikzcd}
\end{equation}

A \emph{morphism} of $\cat{C}$-coalgebras $(S,\alpha)\to(T,\beta)$ is a function $h\colon S\to T$ such that the following diagram commutes:\index{coalgebra!morphism of}
\[
\begin{tikzcd}
    S \ar[d, "h"'] \ar[r, "\alpha"] &
    \car{c}\tri S \ar[d, "\car{c}\:\tri\:h"] \\
    T \ar[r, "\beta"'] &
    \car{c}\tri T
\end{tikzcd}
\]
\end{definition}

\begin{proposition}
Retrofunctors $S\yon^S\cof\cat{C}$ can be identified (up to isomorphism) with $\cat{C}$-coalgebras carried by $S$.\index{coalgebra!carrier of}
\end{proposition}
\begin{proof}
Let $\car{c}$ be the carrier of $\cat{C}$.
In \cref{ex.coalgebras}, we showed that \eqref{eqn.monomials_and_comp} gives a natural correspondence between lenses $\Phi\colon S\yon^S\to\car{c}$ and functions $\varphi\colon S\to\car{c}\tri S$.
We can unravel this correspondence via the proof of \cref{prop.comp_left_coclosed} as follows.
A lens $\Phi\colon S\yon^S\to\car{c}$ can be drawn like so (we will adopt our former convention of identifying each morphism $s\to t$ from $S\yon^S$ with its codomain $t$):
\[
\begin{tikzpicture}[polybox, mapstos]
    \node[poly, dom, "$S\yon^S$" left, my-blue] (l) {$t\vphantom{f}$\at$s\vphantom{i}$};
    \node[poly, cod, "$\car{c}$" right, right=of l, my-red] (r) {$f$\at$i$};
    \draw (l_pos) -- node[below] {$\Phi_\1$} (r_pos);
    \draw (r_dir) -- node[above] {$\Phi^\sharp$} (l_dir);
\end{tikzpicture}
\]
Meanwhile, the corresponding function $\varphi\colon S\to\car{c}\tri S$, equivalently a lens between constants, can be drawn thusly (recall that we color a box red when it is impossible to fill, i.e.\ when it can only be filled by an element of the empty set):
\[
\begin{tikzpicture}[polybox, mapstos]
	\node[poly, constant dom, "$S$" left, my-blue] (p) {\at$s$};
	\node[poly, cod, right=1.5cm of p.south, yshift=-1.25ex, "$\car{c}$" right, my-red] (r) {$f$\at$i$};
	\node[poly, above=of r, "$S$" right, constant, my-blue] (q) {\at$t$};
  	\draw (p_pos) to[first] node[below] {$\Phi_\1$} (r_pos);
  	\draw (r_dir) to[climb] node[right] {$\Phi^\sharp$} (q_pos);
\end{tikzpicture}
\]

Then it suffices to show that $\Phi$ satisfies the retrofunctor laws if and only if $\varphi$ satisfies the coalgebra laws.\index{coalgebra}
We can verify this using polyboxes.
From \eqref{eqn.pres_era_draw}, the eraser preservation law for $\Phi$ would state the following (remember that the arrow in the eraser for the state category $S\yon^S$ is just an equality):\index{polybox}\index{state system}
\[
\begin{tikzpicture}
	\node (id1) {
	\begin{tikzpicture}[polybox, mapstos]
		\node[poly, dom, my-blue, "$S\yon^S$" left] (p) {$s$\at$s$};
		\draw[my-blue,double,-] (p_pos) to[climb'] node[right] {$\idy$} (p_dir);
	\end{tikzpicture}
	};
	\node[right=of id1] (id2) {
	\begin{tikzpicture}[polybox, mapstos]
		\node[poly, dom, my-blue, "$S\yon^S$" left] (p) {\at$s$};
		\node[poly, my-red, right=1 of p, "$\car{c}$" below] (q) {};
		\draw (p_pos) to[first] node[below] {$\Phi_\1$} (q_pos);
		\draw (q_dir) to[last] node[above] {$\Phi^\sharp$} (p_dir);
		\draw[my-red] (q_pos) to[climb'] node[right] {$\idy$} (q_dir);
	\end{tikzpicture}
	};
	\node at ($(id1.east)!.3!(id2.west)$) {$=$}; % -(0,6pt)
\end{tikzpicture}
\]
Meanwhile, the commutative triangle on the left of \eqref{eqn.coalg_laws} can be written as follows:
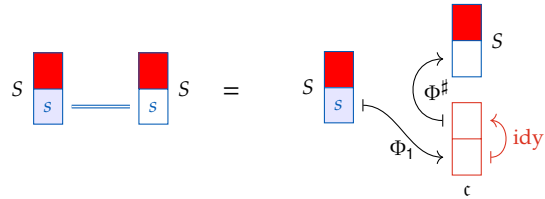
\[
\begin{tikzpicture}
	\node (id1) {
	\begin{tikzpicture}[polybox, mapstos]
    	\node[poly, constant dom, "$S$" left, my-blue] (p) {\at$s$};
    	\node[poly, constant, "$S$" right, my-blue, right=1 of p] (q) {\at$s$};
		\draw[double,my-blue,-] (p_pos) to[first] node[below] {} (q_pos);
	\end{tikzpicture}
	};
	\node[right=of id1, yshift=-1.25ex] (id2) {
	\begin{tikzpicture}[polybox, mapstos]
    	\node[poly, constant dom, "$S$" left, my-blue] (p) {\at$s$};
    	\node[poly, right=1.5cm of p.south, yshift=-1.25ex, "$\car{c}$" below, my-red] (r) {};
    	\node[poly, above=of r, "$S$" right, constant, my-blue] (q) {};

      	\draw (p_pos) to[first] node[below] {$\Phi_\1$} (r_pos);
      	\draw (r_dir) to[climb] node[right] {$\Phi^\sharp$} (q_pos);

		\draw[my-red] (r_pos) to[climb'] node[right] {$\idy$} (r_dir);
    \end{tikzpicture}
	};
	\node at ($(id1.east)!.3!(id2.west)$) {$=$};
\end{tikzpicture}
\]
But these polybox equations are entirely equivalent.

Then from \eqref{eqn.pres_dup_draw}, the duplicator preservation law for $\Phi$ would state the following (remember that the three arrows in the duplicator for the state category $S\yon^S$ are all equalities)
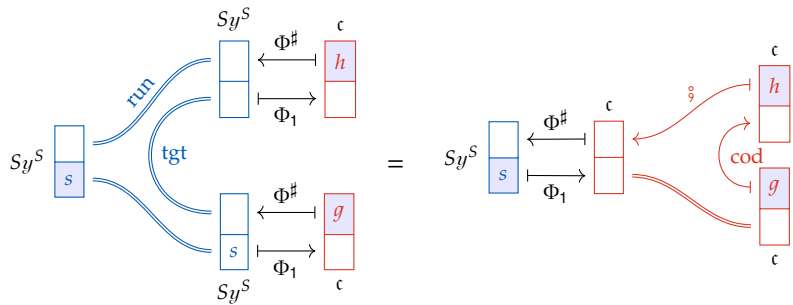
\[
\begin{tikzpicture}
	\node (sp1) {
	\begin{tikzpicture}[polybox, mapstos]
		\node[poly, dom, my-blue, "$S\yon^S$" left] (c) {\at$s$};
		\node[poly, my-blue, right=2 of c.south, yshift=-2.75ex, "$S\yon^S$" below] (c1) {$\vphantom{g}$\at$s$};
		\node[poly, my-blue, above=1 of c1, "$S\yon^S$" above] (c2) {$\vphantom{h}$};
		\node[poly, cod, my-red, right=1 of c1, "$\car{c}$" below] (c'1) {$g$\at$\vphantom{s}$};
		\node[poly, cod, my-red, right=1 of c2, "$\car{c}$" above] (c'2) {$h$};
		\draw[my-blue,double,-] (c_pos) to[first] (c1_pos);
		\draw[my-blue,double,-] (c1_dir) to[climb] node[right] {tgt} (c2_pos);
		\draw[my-blue,double,-] (c2_dir) to[last] node[above,sloped] {run} (c_dir);
		\draw (c1_pos) to[first] node[below] {$\Phi_\1$} (c'1_pos);
		\draw (c'1_dir) to[last] node[above] {$\Phi^\sharp$} (c1_dir);
		\draw (c2_pos) to[first] node[below] {$\Phi_\1$} (c'2_pos);
		\draw (c'2_dir) to[last] node[above] {$\Phi^\sharp$} (c2_dir);
    \end{tikzpicture}
	};
	\node[right=.8 of sp1] (sp2) {
	\begin{tikzpicture}[polybox, mapstos]
		\node[poly, dom, my-blue, "$S\yon^S$" left] (c) {\at$s$};
		\node[poly, my-red, right=1 of c, "$\car{c}$" above] (c') {};
		\node[poly, cod, my-red, right=2 of c'.south, yshift=-1ex, "$\car{c}$" below] (c'1) {$g$};
		\node[poly, cod, my-red, above=of c'1, "$\car{c}$" above] (c'2) {$h$};
		\draw (c_pos) to[first] node[below] {$\Phi_\1$} (c'_pos);
		\draw (c'_dir) to[last] node[above] {$\Phi^\sharp$} (c_dir);
		\draw[my-red,double,-] (c'_pos) to[first] (c'1_pos);
		\draw[my-red] (c'1_dir) to[climb] node[right] {$\cod$} (c'2_pos);
		\draw[my-red] (c'2_dir) to[last] node[above] {$\then$} (c'_dir);
	\end{tikzpicture}
	};
	\node at ($(sp1.east)!.5!(sp2.west)-(0,4pt)$) {$=$};
\end{tikzpicture}
\]
Meanwhile, the commutative triangle on the right of \eqref{eqn.coalg_laws} can be written as follows:
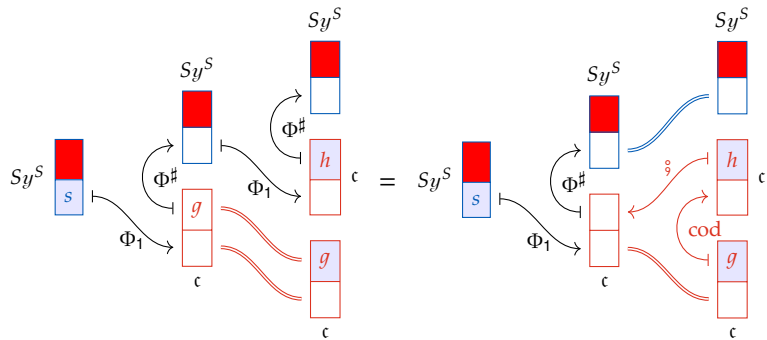
\[
\begin{tikzpicture}
    \node (p1) {
	    \begin{tikzpicture}[polybox, mapstos]
            \node[poly, constant dom, "$S\yon^S$" left, my-blue] (m) {$\vphantom{h}$\at$s$};

            \node[poly, right= of m.south, yshift=-1ex, "$\car{c}$" below, my-red] (D) {$g$\at};
            \node[poly, constant, above=of D, "$S\yon^S$" above, my-blue] (mm) {};

            \node[poly, cod, right= of D.south, yshift=-1ex, "$\car{c}$" below, my-red] (DD) {$g$};
            \node[poly, cod, above=of DD, "$\car{c}$" right, my-red] (mmm) {$h$\at$\vphantom{s}$};
            \node[poly, constant, above=of mmm, "$S\yon^S$" above, my-blue] (C) {};
            \draw (m_pos) to[first] node[below] {$\Phi_\1$} (D_pos);
            \draw (D_dir) to[climb] node[right] {$\Phi^\sharp$} (mm_pos);

            \draw[my-red,double,-] (D_pos) to[first] (DD_pos);
            \draw[my-red,double,-] (DD_dir) to[last] (D_dir);

            \draw (mm_pos) to[first] node[below] {$\Phi_\1$} (mmm_pos);
            \draw (mmm_dir) to[climb] node[right] {$\Phi^\sharp$} (C_pos);
        \end{tikzpicture}
	};
	\node (p2) [right=.3 of p1] {
        \begin{tikzpicture}[polybox, mapstos]
            \node[poly, constant dom, "$S\yon^S$" left, my-blue] (m') {$\vphantom{h}$\at$s$};

            \node[poly, right= of m'.south, yshift=-1ex, "$\car{c}$" below, my-red] (mm') {};
            \node[poly, constant, above=of mm', "$S\yon^S$" above, my-blue] (C') {};

            \node[poly, cod, right= of mm'.south, yshift=-1ex, "$\car{c}$" below, my-red] (D') {$g$};
            \node[poly, cod, above=of D', "$\car{c}$" right, my-red] (mmm') {$h$\at$\vphantom{s}$};
            \node[poly, constant, above=of mmm', "$S\yon^S$" above, my-blue] (CC') {};
            \draw (m'_pos) to[first] node[below] {$\Phi_\1$} (mm'_pos);
            \draw (mm'_dir) to[climb] node[right] {$\Phi^\sharp$} (C'_pos);

            \draw[my-red,double,-] (mm'_pos) to[first] (D'_pos);
            \draw[my-red] (D'_dir) to[climb] node[right] {cod} (mmm'_pos);
            \draw[my-red] (mmm'_dir) to[last] node[above] {$\then$} (mm'_dir);

            \draw[my-blue,double, -] (C'_pos) to[first] (CC'_pos);
        \end{tikzpicture}
    };
	\node at ($(p1.south)!.5!(p2.north)-(0,4pt)$) {$=$};
\end{tikzpicture}
\]
But these polybox equations are equivalent as well.\index{polybox}
Hence the retrofunctor laws for $\Phi$ are equivalent to the coalgebra laws for $\varphi$.
\end{proof}

So a retrofunctor from a state category to $\cat{C}$ bears the same data as a $\cat{C}$-coalgebra.
The equivalences don't stop there, however.

\subsubsection{As discrete opfibrations}

\index{discrete opfibration!as retrofunctor}\index{functor!discrete opfibration|see{discrete opfibration}}
The concept of a $\cat{C}$-coalgebra is in turn equivalent to a better known categorical construction on $\cat{C}$, which we introduce here.

\begin{definition}[Discrete opfibration]\label{def.dopf}
Let $\cat{C}$ be a category. A pair $(\cat{S},\pi)$, where $\cat{S}$ is a category and $\pi\colon\cat{S}\to\cat{C}$ is a functor, is called a \emph{discrete opfibration over $\cat{C}$} if it satisfies the following condition:
\begin{itemize}
	\item for every object $s\in\cat{S}$, object $c'\in\cat{C}$, and morphism $f\colon \pi(s)\to c'$ there exists a unique object $s'\in\cat{S}$ and morphism $\ol{f}\colon s\to s'$ such that $\pi(\ol{f})=f$. Note that in this case $\pi(s')=c$.
\end{itemize}
\[
\begin{tikzcd}
  s\ar[r, dashed, "\ol{f}"]\ar[d, |->, "\pi"']&
  s'\ar[d, |->, "\pi"]\\
  \pi(s)\ar[r, "f"']&
  \_\vphantom{()}
\end{tikzcd}
\]
A \emph{morphism} of discrete opfibrations $(\cat{S},\pi)\to(\cat{S}',\pi')$ over $\cat{C}$ is a functor $F\colon\cat{S}\to\cat{S}'$ making the following triangle commute:
We refer to $\ol{f}$ as the \emph{lift} of $f$ to $s$.

A \emph{morphism} $(\cat{S},\pi)\to(\cat{S}',\pi')$ between discrete opfibrations over $\cat{C}$ is a functor $F\colon\cat{S}\to\cat{S}'$ making the following triangle commute:
\begin{equation}\label{eqn.dopf_triangle}
\begin{tikzcd}[column sep=small]
	\cat{S}\ar[dr, "\pi"']\ar[rr, "F"]&&
	\cat{S}'\ar[dl, "\pi'"]\\&
	\cat{C}
\end{tikzcd}
\end{equation}
We denote the category of discrete opfibrations over $\cat{C}$ by $\Cat{dopf}(\cat{C})$.
\end{definition}\index{discrete opfibration!morphism of}\index{discrete opfibration!category of}

\begin{exercise}
Show that if $F\colon \cat{S}\to\cat{S}'$ is a functor making the triangle \eqref{eqn.dopf_triangle} commute, where both $\pi$ and $\pi'$ are discrete opfibrations, then $F$ is also a discrete opfibration.
\begin{solution}
In general, a functor $F\colon\cat{C}\to\cat{D}$ is a discrete opfibration iff, for every object $c\in\cat{C}$, the induced map $F[c]\colon\cat{C}[c]\to\cat{D}[Fc]$ is a bijection.

In our case, if $\pi$ and $\pi'$ are discrete opfibrations then for any $s\in\cat{S}$ we have that both the second map and the composite in $\cat{S}[s]\To{F[s]}\cat{S}'[Fs]\To{\pi'[Fs]}\cat{C}[\pi s]$ are bijections, so the first map is too.
\end{solution}
\end{exercise}

\begin{exercise}\label{exc.dopf_cof}
Suppose $\pi\colon\cat{S}\to\cat{C}$ is a discrete opfibration and $i\in\cat{S}
$ is an object. With notation as in \cref{def.dopf}, show the following:
\begin{enumerate}
	\item Show that the lift $\ol{\id_{\pi(i)}}=\id_i$ of the identity on $\pi(i)$ is the identity on $i$.
	\item Show that for $f\colon\pi(i)\to c$ and $g\colon c\to c'$, we have $\ol{f}\then\ol{g}=\ol{f\then g}$.
	\item Show how $\pi$ could instead be interpreted as a retrofunctor.% (We will characterize exactly which retrofunctors can arise in this way in \cref{**}.)
\qedhere
\end{enumerate}
\begin{solution}
\begin{enumerate}
    \item The identity map $\id_i$ satisfies $\pi(\id_i)=\id_{\pi(i)}$, so $\ol{\id_{\pi(i)}}=\id_i$ by the uniqueness-of-lift condition in \cref{def.dopf}.
    \item The morphism $\ol{f}\then\ol{g}$ satisfies $\pi(\ol{f}\then\ol{g})=\pi(\ol{f})\then\pi(\ol{g})=f\then g$, so again this follows by uniqueness of lift.
    \item To see that $\pi$ is a retrofunctor, we need to understand and verify conditions about its action forward on objects and backward on morphisms. Its action forward on objects is that of $\pi$ as a functor. Given an object $s\in\cat{S}$, the action of $\pi$ backward on morphisms is the lift operation $\pi^\sharp_s(f)\coloneqq\ol{f}$. This preserves identity by part 1, composition by part 2, and codomains because $\pi(\cod(\pi^\sharp_s(f)))=\cod(\pi(\pi^\sharp_s(f)))=\cod(f)$.
\end{enumerate}
\end{solution}
\end{exercise}

As it turns out, not only do $\cat{C}$-coalgebras carry the same data as discrete opfibrations over $\cat{C}$, they in fact comprise isomorphic categories.

\begin{proposition}\index{coalgebra}
The category of $\cat{C}$-coalgebras is isomorphic to the category $\Cat{dopf}(\cat{C})$ of discrete opfibrations over $\cat{C}$.
\end{proposition}\index{discrete opfibration!as retrofunctor}
%\begin{proof}
%**
%\end{proof}
\index{copresheaf!as retrofunctor}\index{functor!copresheaf|see{copresheaf}}

\subsubsection{As copresheaves}

It is well-known in the category theory literature that the category of discrete opfibrations over $\cat{C}$ is equivalent to yet another familiar category: the category $\smset^{\cat{C}}$, whose objects are functors $\cat{C}\to\smset$.
Such a functor is known as a \emph{copresheaf} on $\cat{C}$ for short. These are relevant, e.g.\ in the theory of categorical databases \cite{spivak2012functorial}.
Here we review what is needed to understand this equivalence.
We begin by giving a standard construction on any copresheaf.\index{database}

\begin{definition}[Category of elements]\label{def.cat_elements}\index{elements!category of elements}\index{category!of elements}
Given a copresheaf on $\cat{C}$, i.e.\ a functor $I\colon\cat{C}\to\smset$, its \emph{category of elements} $\elts^\cat{C}I$ is defined to have objects
\[
    \Ob\elts^\cat{C}I\coloneqq\{(c,x)\mid c\in\cat{C}, x\in Ic\}
\]
and a morphism $f\colon(c,x)\to(c',x')$ for every morphism $f\colon c\to c'$ from $\cat{C}$ satisfying
\[
    (If)(x)=x'.
\]
Identities and composites in $\elts^\cat{C}I$ are inherited from $\cat{C}$; they obey the usual category laws by the functoriality of $I$.

The category is so named because its objects are the elements of the sets that the objects of $\cat{C}$ are sent to by $I$.
Each morphism $f\colon c\to\_$ in $\cat{C}$ then becomes as many morphisms in $\elts^\cat{C}I$ as there are elements of $Ic$, tracking where $If$ sends each such element.
\end{definition}

The next exercise shows how this construction turns every copresheaf into a discrete opfibration.\index{copresheaf!and discrete opfibration}

\begin{exercise} \label{exc.copre_to_dopf}\index{functor!set-valued}\index{category!of elements}
Let $I\colon\cat{C}\to\smset$ be a functor, and let $\int^{\cat{C}}I$ be as in \cref{def.cat_elements}.
\begin{enumerate}
    \item Show that there is a functor $\pi\colon\elts^\cat{C}I\to\cat{C}$ sending objects $(c,x)\mapsto c$ and morphisms $f\colon(c,x)\to(c',x')$ to $f\colon c\to c'$.
    \item Show that $\pi$ is in fact a discrete opfibration. \qedhere
\end{enumerate}
\begin{solution}
\begin{enumerate}
    \item The functor action on objects and morphisms is defined in the exercise statement, so it suffices to check that it preserves identities and composites. But this too is obvious: $\pi(\id)=\id$ and $\pi(f\then g)$ is $f\then g$.
    \item Given an object $(c,x)\in\int^\cat{C}I$ and $f\colon c\to c'$, we need to check that there is a unique object $s\in\int^{\cat{C}}I$ and morphism $\ol{f}\colon(c,x)\to s$ with $\pi(\ol{f})=f$. Using $x'\coloneqq I(f)(x)$ and $s\coloneqq(c',x')$ and $\ol{f}\coloneqq f$, we do indeed get $\pi(\ol{f})=f$ and we find that it is the only possible choice.
\end{enumerate}
\end{solution}
\end{exercise}

In fact, the assignment of a discrete opfibration to every copresheaf given above is functorial, as the next exercise shows.

\begin{exercise}\label{exc.elts_functor}\index{natural transformation!as copresheaf morphism}
Suppose that $I,J\colon\cat{C}\to\smset$ are functors and $\alpha\colon I\to J$ is a natural transformation.
\begin{enumerate}
	\item Show that $\alpha$ induces a functor $(\elts^\cat{C}I)\to(\elts^\cat{C}J)$.
	\item Show that it is a morphism of discrete opfibrations in the sense of \cref{def.dopf}.
	\item Have you now verified that there is a functor
	\[
	\elts^\cat{C}\colon\smset^{\cat{C}}\to\Cat{dopf}(\cat{C})
	\]
	or is there something left to do?
	\end{enumerate}
\begin{solution}
\begin{enumerate}
    \item Given a natural transformation $\alpha\colon I\to J$, we need a functor $\int^\cat{C}\alpha\colon\int^\cat{C}I\to\int^\cat{C}J$. On objects have it send $(c,x)\mapsto(c,\alpha_cx)$, where $x\in I(c)$ so $\alpha_cx\in J(c)$; on morphisms have it send $f\mapsto f$, a mapping which clearly preserves identities and composition.
    \item To see that $\int^\cat{C}\alpha$ is a morphism between discrete opfibrations, one only needs to check that it commutes with the projections to $\cat{C}$, but this is obvious: the mapping $(c,x)\mapsto(c,\alpha_cx)$ and the mapping $f\mapsto f$ preserve the objects and morphisms of $\cat{C}$.
    \item It is clear that $\int^\cat{C}$ preserves identities and composition in $\smset^\cat{C}$, so we have verified that this is functor.
\end{enumerate}
\end{solution}
\end{exercise}\index{category!of elements}

\begin{exercise}\label{exc.elts_free_grph}
Let $G$ be a graph, and let $\cat{G}$ be the free category on it. Show that for any functor $S\colon\cat{G}\to\smset$, the category $\elts^\cat{G} S$ of elements is again free on a graph.
\begin{solution}
Given a graph $G\coloneqq(\fun{src},\fun{tgt}\colon A\tto V)$ and functor $S\colon\cat{G}\to\smset$, define a new graph $H$ as the following diagram of sets:
\[
\begin{tikzcd}[column sep=50pt]
	\sum_{a\in A} S(\fun{src}(a))\ar[r, shift left, "{(\fun{src},\id)}"]\ar[r, shift right, "{(\fun{\tgt},S(a))}"']&
	\sum_{v\in V} S(v)
\end{tikzcd}
\]\index{graph}
In other words, it is a graph for which a vertex is a pair $(v,s)$ with $v\in V$ a $G$-vertex and $s\in S(v)$, and for which an arrow is a pair $(a, s)$ with $a\in A$ a $G$-arrow, say $a\colon v\to v'$, and $s\in S(v)$ is an element over the source of that arrow. With this notation, the arrow $(a,s)$ has as source vertex $(v,s)$ and has as target vertex $(v',S(a)(s))$ where $S(a)\colon S(v)\to S(v')$ is the function given by the functor $S\colon\cat{G}\to\smset$.

It remains to see that the free category on $H$ is isomorphic to $\int^{\cat{G}}S$. We first note that it has the same set of objects, namely $\sum_{v\in V}S(v)$. A morphism in $\int^\cat{G}S$ can be identified with an object $(v,s)$ together with a morphism $f\colon v\to v'$ in $\cat{G}$, but this is just a length-$n$ sequence $(a_1,\ldots,a_n)$ of arrows, with $v=\fun{src}(a_1)$, $\fun{tgt}(a_{i})=\fun{src}(a_{i+1})$ for $1\leq i<n$, and $\fun{tgt}(a_n)=v'$. This is the same data as a morphism in the free category on $H$.
\end{solution}
\end{exercise}\index{functor!set-valued}

\begin{proposition}\label{prop.tfae_dopf}\index{discrete opfibration!as copresheaf}\index{category!of elements}
The category $\smset^{\cat{C}}$ of copresheaves on $\cat{C}$ is equivalent to the category of discrete opfibrations over $\cat{C}$.
\end{proposition}
\begin{proof}
By \cref{exc.elts_functor} we have a functor $\elts^\cat{C}\colon\smset^{\cat{C}}\to\Cat{dopf}(\cat{C})$. There is a functor going back: given a discrete opfibration $\pi\colon\cat{S}\to\cat{C}$, we define a functor $\partial\pi\colon\cat{C}\to\smset$ on objects by sending each $c\in\cat{C}$ to the set of objects in $\cat{S}$ that $\pi$ maps to $c$; that is,
\[
    (\partial\pi)(c)\coloneqq\{s\in\cat{S}\mid\pi(s)=c\}.
\]
Then on morphisms, for each $f\colon c\to\_$ in $\cat{C}$ and $s\in(\partial\pi)(c)$ we have $\pi(s)=c$, so by \cref{def.dopf} there exists a unique morphism $\ol{f}\colon s\to\_$ for which $\pi(\ol{f})=f$.
As $\pi(\cod\ol{f})=\cod f$, we have $\cod\ol{f}\in(\partial\pi)(\cod f)$, so we can define
\[
    (\partial\pi)(f)(s)\coloneqq\cod\ol{f}.
\]
On objects, the roundtrip $\smset^{\cat{C}}\to\smset^{\cat{C}}$ sends $I\colon\cat{C}\to\smset$ to the functor
\begin{align*}
	c&
	\mapsto\{s\in\elts^\cat{C}I\mid \pi(s)\\&
	=\{(c,x)\mid x\in I(c)\}&=I(c).
\end{align*}
The roundtrip $\Cat{dopf}(\cat{C})\to\Cat{dopf}(\cat{C})$ sends $\pi\colon\cat{S}\to\cat{C}$ to the discrete opfibration whose object set is $\{(c,s)\in\Ob(\cat{C})\times\Ob(\cat{S})\mid\pi(s)=c\}$ and this set is clearly in bijection with $\Ob(\cat{S})$. Proceeding similarly, one defines an isomorphism of categories $\cat{S}\cong\elts^\cat{C}\partial\pi$.
\end{proof}

\begin{proposition}\label{prop.ds_dopf}
Up to isomorphism, discrete opfibrations into $\cat{C}$ can be identified with dynamical systems on $\cat{C}$.
\end{proposition}
In case it isn't clear, this association is only functorial on the groupoid of objects and isomorphisms.
\begin{proof}
Given a discrete opfibration $\pi\colon\cat{S}\to\cat{C}$, take $S\coloneqq\Ob(\cat{S})$ and define $(\varphi_1,\varphi^\sharp)\colon S\yon^S\to\car{c}$ by $\varphi_1=\pi$ and with $\varphi^\sharp$ given by the lifting: $\varphi(g)\coloneqq \hat{g}$ as in \cref{def.dopf}. One checks using \cref{exc.dopf_cof} that this defines a retrofunctor.

Conversely, given a retrofunctor $(\varphi_1,\varphi^\sharp)\colon S\yon^S\to\car{c}$, the function $\varphi_1$ induces a lens $S\yon\to\car{c}$, and we can factor it as a vertical followed by a cartesian $S\yon\to\car{s}\To{\psi}\car{c}$. We can give $\car{s}$ the structure of a category such that $\psi$ is a retrofunctor; see \cref{exc.ds_dopf}.
\end{proof}

\begin{exercise}\label{exc.ds_dopf}
With notation as in \cref{prop.ds_dopf}, complete the proof as follows.
\begin{enumerate}
	\item Check that $(\varphi,\varphi^\sharp)$ defined in the first paragraph is indeed a retrofunctor.
	\item Find a comonoid structure on $\car{s}$ such that $\psi$ is a retrofunctor, as stated in the second paragraph.
	\item Show that the two directions are inverse, up to isomorphism.
\qedhere
\end{enumerate}
\end{exercise}

Given a dynamical system $S\yon^S\to p$, we extend it to a retrofunctor $\varphi\colon S\yon^S\cof\cofree{p}$. By \cref{prop.tfae_dopf,prop.ds_dopf}, we can consider it as a discrete opfibration over $\cofree{p}$. By \cref{exc.elts_free_grph} the category $\elts \varphi$ is again free on a graph. It is this graph that we usually draw when depicting the dynamical system, e.g.\ in \eqref{eqn.dyn_sys_misc573}.

\index{graph!category free on}\index{category!free}

%\begin{exercise}
%Give an example of a dynamical system on $p\coloneqq\yon^\2+\yon$ for which the corresponding copresheaf has in its image a set with at least two elements.
%\begin{solution}
%**
%\end{solution}
%\end{exercise}

%\begin{exercise}
%Given a retrofunctor $F\colon S\yon^S\cof\yon$, what does the corresponding copresheaf look like?
%\begin{solution}
%**
%\end{solution}
%\end{exercise}

To summarize, we have four equivalent notions:\index{coalgebra!as equivalent to other notions}
\begin{enumerate}[label=(\arabic*)]
    \item retrofunctors $F\colon S\yon^S\cof\cat{C}$;
    \item $\cat{C}$-coalgebras $(S,\alpha)$, with $\alpha\colon S\to\car{c}\tri S$;
    \item discrete opfibrations $\pi\colon \cat{S}\to\cat{C}$, with $\Ob\cat{S}=S$;
    \item copresheaves $I\colon\cat{C}\to\smset$, with $\Ob\elts^{\cat{C}}I=S$.
\end{enumerate}
Moreover, (2), (3), and (4) form equivalent categories.
Translating between these notions yields different perspectives on familiar categorical concepts.
In the next chapter, we will discover even more characterizations of the same data within $\poly$ (see \cref{prop.cart_dopf}).

\index{state system!retrofunctor from|(}

%-------- Section --------%
\section[Summary and further reading]{Summary and further reading%
  \sectionmark{Summary \& further reading}}
\sectionmark{Summary \& further reading}

In this chapter we began by showing that for every set $S$, the thing that makes a dynamical system like $S\yon^S\to p$ actually run is the fact that $S\yon^S$ has the structure of a comonoid.  We then explained Ahman-and-Uustalu's result that comonoids, i.e.\ polynomials $p$, equipped with a pair of lenses $\epsilon\colon p\to\yon$ and $\delta\colon p\to p\tri p$, are exactly categories \cite{ahman2016directed}. We explained how $\epsilon$ picks out an identity for each object and how $\delta$ picks out a codomain for each morphism and a composite for each composable pair of morphisms. In particular we showed that the category corresponding to $S\yon^S$ is the \emph{contractible groupoid on $S$}, i.e.\ the category with $S$-many objects and a unique morphism between any two.

We then discussed how comonoid morphisms $\car{c}\to\car{d}$ are not functors but \emph{retrofunctors}: they map forwards on objects and backwards on morphisms. Retrofunctors were first defined by Marcelo Aguiar \cite{aguiar1997internal}, though his definition was opposite to ours; he referred to these as \emph{cofunctors}. Retrofunctors are the morphisms of a category we notate as $\catsharp$. We showed that retrofunctors between state categories $S\yon^S\cof T\yon^T$ are what are known in the functional programming community as \emph{very well-behaved lenses}. For more on lenses, see \cite{nlab2022lens}.

\index{retrofunctor!Aguiar's definition}

%-------- Section --------%
\section{Exercise solutions}
\Closesolutionfile{solutions}
{\footnotesize
\input{solution-file7}}

\Opensolutionfile{solutions}[solution-file8]

%------------ Chapter ------------%
\chapter{Categorical properties of polynomial comonoids}
\chaptermark{Properties of comonoids in $\poly$}
\label{ch.comon.cofree}

While we defined the category $\catsharp$ of categories and retrofunctors in the last chapter, we have only begun to scratch the surface of the properties it satisfies.
As the category of comonoids in $\poly$, equipped with a canonical forgetful functor $U\colon\catsharp\to\poly$ sending each polynomial to its carrier, $\catsharp$ inherits many of the categorical properties satisfied by $\poly$.
Underpinning this inheritance is the fact that $U$ has a right adjoint, a cofree functor $\poly\to\catsharp$.
We will introduce this adjunction in the first section of this chapter.
Then we will discuss how many of the categorical properties of $\poly$, including much of what we covered in \cref{ch.poly.bonus}, play nicely with restriction to $\catsharp$.
Finally, we will touch on other constructions we can make over the comonoids in $\poly$, such as their comodules and coalgebras.\index{coalgebra}

\index{$\catsharp$!categorical properties of|(}\index{functor!forgetful $\catsharp\to\poly$}

%-------- Section --------%
\section{Cofree comonoids} \label{sec.comon.cofree.cons}

\index{comonoid!cofree|see{cofree comonoid}}
\index{cofree comonoid|(}\index{interface}

Consider a dynamical system $\varphi\colon\car{s}\to p$ with state system $\car{s}$ and interface $p$.
In \cref{subsec.comon.sharp.state.run}, we posed the question of whether there was a single morphism that could capture all the information encoded by the family of lenses $\text{Run}_n(\varphi)\colon\car{s}\to p\tripow{n}$ for all $n\in\nn$, defined as the composite
\[
    \car{s}\To{\delta^{(n)}}\car{s}\tripow{n}\To{\varphi\tripow{n}}p\tripow{n}
\]
that models $n$ runs through the system $\varphi$.

It turns out that there is: the key is that there is a natural way to interpret every polynomial $p$ as a category $\cofree{p}$, so that retrofunctors into $\cofree{p}$ are exactly lenses into $p$.
In other words, $\cofree{p}$ will turn out to be the \emph{cofree comonoid} (or \emph{cofree category}) \emph{on $p$}.
Cofree comonoids in $\poly$ are beautiful objects, both in their visualizable structure as a category and in the metaphors we can make about them. They allow us to replace the interface of a dynamical system with a category and get access to a rich theory that exists there.

We'll go through the construction of the cofree comonoid and its implications in this section, featuring a purely formal proof of the fact the forgetful functor $U\colon\catsharp\to\poly$ has a right adjoint $\cofree{-}\colon\poly\to\catsharp$, where for each $p\in\poly$, the carrier $\car{t}_p\coloneqq U\cofree{p}$ of the category $\cofree{p}$ is given by the limit of the following diagram:\index{functor!forgetful $\catsharp\to\poly$}
\begin{equation} \label{eqn.cofree_diagram}\index{limit!cofree comonoid as}
\begin{tikzcd}
	\yon \ar[d] &
	p \ar[d] &
	p\tripow2 \ar[d] &
	[10pt] p\tripow3 \ar[d] &
	\cdots\\
	\1 &
	p\tri\1 \ar[l, "!"'] &
	p\tripow2\tri\1 \ar[l, "p\:\tri\:!"'] &
	p\tripow3\tri\1 \ar[l, "p\tripow2\:\tri\:!"'] &
	\cdots.\ar[l]
\end{tikzcd}
\end{equation}
Thus, we will show in \cref{thm.cofree} that $U$ and $\cofree{-}$ form a forgetful-cofree adjunction, making $\cofree{p}$ the cofree comonoid on $p$.
But first, let us concretely characterize the canonical comonoid structure on the limit of \eqref{eqn.cofree_diagram}, before showing that it is indeed cofree.

%---- Subsection ----%
\subsection{The carrier of the cofree comonoid} \label{subsec.comon.cofree.cons.car}

\index{cofree comonoid!carrier of|(}

Let $\car{t}_p$ be the limit of the diagram \eqref{eqn.cofree_diagram} in $\poly$; it will turn out to be the carrier of the cofree comonoid on $p$ (where $p$ will be an arbitrary polynomial throughout).
We could compute this limit directly, but we will be able to describe it more concretely in terms of what we call \emph{trees on $p$} or \emph{$p$-trees}: trees comprised of $p$-corollas.
In doing so, we will formalize the tree pictures we have been using to describe polynomials all along.

\subsubsection{Trees on polynomials}

\index{tree!on a polynomial}

\begin{definition}[Tree on a polynomial] \label{def.poly_tree}
Let $p\in\poly$ be a polynomial.
A \emph{tree on $p$}, or a \emph{$p$-tree}, is a rooted tree whose every vertex $v$ is assigned a $p$-position $i$ and a bijection from the children of $v$ to $p[i]$.
We denote the set of $p$-trees by $\tr_p$.
\end{definition}
We can think of a $p$-tree as being ``built'' out of $p$-corollas according to these instructions:
\begin{quote}
To choose a $p$-tree in $\tr_p$:
\begin{enumerate}
    \item choose a $p$-corolla:
    \begin{itemize}
        \item its root $i_0\in p(\1)$ will be the tree's root, and
        \item its leaves in $p[i_0]$ will be the edges out of the root;
    \end{itemize}
    \item for each $p[i_0]$-leaf $a_1$:
    \begin{enumerate}[label*=\arabic*.]
        \item choose a $p$-corolla:
        \begin{itemize}
            \item its root $i_1\in p(\1)$ will be the vertex adjoined to $a_1$, and
            \item its leaves in $p[i_1]$ will be the edges out of that vertex;
        \end{itemize}
        \item for each $p[i_1]$-leaf $a_2$:
        \begin{enumerate}[label*=\arabic*.]
            \item choose a $p$-corolla:
            \begin{itemize}
                \item its root $i_2\in p(\1)$ will be the vertex adjoined to $a_2$, and
                \item its leaves in $p[i_2]$ will be the edges out of that vertex;
            \end{itemize}
            \item for each $p[i_1]$-leaf $a_2$:

            $\cdots$
        \end{enumerate}
    \end{enumerate}
\end{enumerate}
\end{quote}
Of course, there may eventually be multiple copies of any one $p$-root as a vertex or $p$-leaf as an edge in our $p$-tree, and these vertices and edges are not literally the same.
So we should really think of the positions and directions of $p$ involved each step as \emph{labels} for the vertices and edges of a $p$-tree---although crucially, each $p[i]$-direction must be used as a label exactly once among the edges emanating from a given vertex labeled with $i$.

Although these instructions continue forever, we could abbreviate them by writing them recursively:
\begin{quote}
To choose a $p$-tree in $\tr_p$:
\begin{enumerate}
    \item choose a $p$-corolla:
    \begin{itemize}
        \item its root $i_0\in p(\1)$ will be (the label of) the tree's root, and
        \item its leaves in $p[i_0]$ will be (the labels of) the edges out of the root;
    \end{itemize}
    \item for each $p[i_0]$-leaf $a_1$:
    \begin{enumerate}[label*=\arabic*.]
        \item choose a $p$-tree in $\tr_p$:
        \begin{itemize}
            \item it will be the subtree whose root is adjoined to $a_1$.
        \end{itemize}
    \end{enumerate}
\end{enumerate}
\end{quote}

We would like to draw some examples of $p$-trees, but note that a $p$-tree can have infinite height---in fact, it always will unless every one of its branches terminates at a $p$-corolla with no leaves, i.e.\ a position with an empty direction set.
This means that plenty of $p$-trees cannot be drawn, even when all the position- and direction-sets of $p$ are finite; but we can instead consider finite-height portions of them that we will call \emph{pretrees}.

\subsubsection{Pretrees on polynomials}
\index{tree!pretree}

Before we define pretrees, let's give some examples.

\begin{example}[A few example $p$-pretrees] %TODO
Let $p\coloneqq\{\bul[my-red],\bul[my-blue]\}\yon^\2+\{\bul[black]\}\yon+\{\bul[my-yellow]\}$.
Here are four partially constructed $p$-trees:
\begin{equation}\label{eqn.some_trees_misc58}
\begin{tikzpicture}[trees,
  level 1/.style={sibling distance=10mm},
  level 2/.style={sibling distance=5mm},
  level 3/.style={sibling distance=2.5mm}]
	\node[my-red] (a) {$\bullet$}
		child {node[my-red] {$\bullet$}
			child {node[my-yellow] {$\bullet$}}
			child {node[my-yellow] {$\bullet$}}
		}
		child {node[black] {$\bullet$}
			child {node[black] {$\bullet$}
				child
			}
		}
		;
	\node[black, right=2 of a] (b) {$\bullet$}
		child {node[my-blue] {$\bullet$}
			child {node[my-blue] {$\bullet$}
				child
				child
			}
			child {node[my-red] {$\bullet$}
				child
				child
			}
		};
	\node[my-yellow, right=2 of b] (c) {$\bullet$};
	\node[my-red, right=2 of c] {$\bullet$}
		child {node[my-blue] {$\bullet$}
			child {node[my-red] {$\bullet$}
				child
				child
			}
			child {node[my-red] {$\bullet$}
				child
				child
			}
		}
		child {node[my-red] {$\bullet$}
			child {node[my-red] {$\bullet$}
				child
				child
			}
			child {node[my-red] {$\bullet$}
				child
				child
			}
		};
\end{tikzpicture}
\end{equation}
Here only the third one---the single yellow dot---would count as an element of $\tr_p$.
After all, in \cref{def.poly_tree}, when we speak of a tree on $p$, we mean a tree for which every vertex is a position in $p$ with all of its emanating directions filled by another position in $p$.
Since three of the four trees shown in \eqref{eqn.some_trees_misc58} have leaves emanating from the top that have not been filled by any $p$-corollas, these trees are not elements of $\tr_p$.

However, each of these trees can be extended to an actual element of $\tr_p$ by continually filling in each open leaf with another $p$-corolla.
These might continue forever---or, if you're lazy, you could just cap them all off with the direction-less yellow dot.
\end{example}

\begin{exercise}
Let $q\coloneqq\yon^\2+\3\yon^\1$.
Are there any finite $q$-trees?
If not, could there be any vertices of a given $q$-tree with finitely many descendents?
\begin{solution}
With $q\coloneqq\yon^\2+\3\yon^\1$, every vertex of a $q$-tree starting from the root has either $1$ or $2$ outgoing edges, from which it follows that the tree is infinite and that every vertex has infinitely many descendents.
\end{solution}
\end{exercise}

The trees in \eqref{eqn.some_trees_misc58} can all be obtained by following just the first $3$ levels of instructions for building a $p$-tree (in fact, exactly as many instructions as we initially wrote out).
On the other hand, we know from \cref{subsec.comon.comp.def.corolla} that such a tree represents a position of $p\tripow3$, whose directions are its height-$3$ leaves---and this is true for any $n\in\nn$ in place of $3$.
So these trees are still important to our theory; but since they are not always complete $p$-trees, we will call them something else.

\index{tree!stages}

\begin{definition} \label{def.tripow_pretree}
Given $p\in\poly$, a \emph{stage-$n$ pretree} (or \emph{$p\tripow{n}$-pretree}) is defined to be an element of $p\tripow{n}(\1)$. For each $i\in p\tripow{n}(\1)$ the \emph{height-$n$ leaves} of $i$ is defined to be the set $p\tripow{n}[i]$.
\end{definition}

\begin{remark} \label{rmk.stage}
In \cref{def.tripow_pretree}, we use \emph{stage} instead of \emph{height} to allow for the fact that a $p\tripow{n}$-pretree may not reach its maximum height $n$ if all of its branches terminate early.
For example, the yellow dot in \eqref{eqn.some_trees_misc58} is a $p\tripow1$-pretree, but it is also a $p\tripow2$-pretree, a $p\tripow{50}$-pretree, and indeed a $p$-tree.

Note that for any polynomial $p\in\poly$ there is exactly one stage-$0$ pretree on $p$ for any $p\in\poly$, because $p\tripow0(\1)=\yon(\1)=\1$.
\end{remark}

\begin{example}[Trimming pretrees] \label{ex.trim}
Since $p\tripow1(\1)\iso p(\1)$ and $p\tripow0(\1)\iso\yon(\1)\iso\1$, the unique function $!\colon p(\1)\to\1$ can be thought of as a function from $p\tripow1$-positions to $p\tripow0$-positions, or equivalently a function from stage-$1$ pretrees (i.e.\ corollas) to stage-$0$ pretrees on $p$.
We can interpret this function as taking a corolla and ``stripping away'' its leaves along with the position-label on its root, leaving only a single unlabeled root: a stage-$0$ pretree.

This deceptively simple function has a surprising amount of utility when combined with other maps.
For any $n\in\nn$, we can take the composition product in $\poly$ of the identity on $p\tripow{n}$ and $!$, interpreted as a lens between constant polynomials, to obtain a lens $p\tripow{n}\tri \:!\colon p\tripow{n}\tri p(\1)\to p\tripow{n}\tri\1$, or equivalently a function $p\tripow{n}(!)\colon p\tripow{n+1}(\1)\to p\tripow{n}(\1)$ from stage-$(n+1)$ pretrees to stage-$n$ pretrees on $p$.\index{constant polynomial}

We can deduce the behavior of this function on stage-$(n+1)$ pretrees from what we know about how the composition product interacts with pretrees on $p$.
The identity lens on $p\tripow{n}$ keeps the lower $n$ levels of each $p\tripow{(n+1)}$-pretree intact, while $!$ will ``strip away'' the $p\tripow{(n+1)}$-pretree's height-$(n+1)$ leaves, along with all the position-labels on its height-$n$ vertices.
Thus $p\tripow{n}(!)$ is the function sending every stage-$(n+1)$ pretree on $p$ to its stage-$n$ pretree, effectively \emph{trimming it down} a level.

We can go even further: composing several such functions yields a composite function
\begin{equation} \label{eqn.trim}
    p\tripow{n}(\1) \From{p\tripow{n}(!)}
	p\tripow{(n+1)}(\1) \from\cdots\from
	p\tripow{(n+k-1)}(\1) \From{p\tripow{(n+k-1)}(!)}
	p\tripow{(n+k)}(\1)
\end{equation}
sending every stage-$(n+k)$ pretree on $p$ to its stage-$n$ pretree by trimming off its top $k$ levels.
We will see these functions again shortly.
\end{example}

\index{tree!as limit of pretrees}
\subsubsection{Trees as a limit of pretrees}

Before we go any further in the theory of $p$-trees, let us look at some more examples.

\begin{example}[A few more actual $p$-trees] \label{ex.imagining_trees}
Keeping $p\coloneqq\{\bul[my-red],\bul[my-blue]\}\yon^\2+\{\bul[black]\}\yon+\{\bul[my-yellow]\}$, here are some elements of $\tr_p$ that we could imagine (or even draw, at least in part):
\begin{itemize}
	\item The binary tree that's ``all red all the time.''
	\item The binary tree where odd layers are red and even layers are blue.
	\item The binary tree whose root is red, but after which every left child is red and every right child is blue.%
	\tablefootnote{To formalize the notions of ``left'' and ``right,'' we could think of the direction-sets of the red and blue dots as $\2\iso\{\const{left},\const{right}\}$, so that out of every vertex there is an edge labeled $\const{left}$ and an edge labeled $\const{right}$.}
	\item The tree where all the nodes are red, except for the rightmost branch, which (apart from the red root) is always green.
	\item Any finite tree, where every branch terminates in a yellow dot.
	\item A completely random tree: for the root, randomly choose either red, blue, green, or yellow, and at every leaf, loop back to the beginning, i.e.\ randomly choose either red, blue, green, or yellow, etc.
\end{itemize}
In fact, there are uncountably many trees in $\tr_p$ (even just $\tr_{2\yon}$ has cardinality $\2^\nn$), but only countably many can be uniquely characterized in a finite language like English (and of course only finitely many can be uniquely characterized in the time we have!).
Thus most elements of $\tr_p$ cannot even be described.
\end{example}

\begin{exercise}\label{exc.p_tree_sets}
% \begin{enumerate}
% 	\item Interpret each of the five tree examples imagined in \cref{ex.imagining_trees} by drawing three or four layers (your choice) of it.
% \end{enumerate}
For each of the following polynomials $p$, characterize the set of trees $\tr_{p}$.
\begin{enumerate}[resume]
	\item $p\coloneqq\1$.
	\item $p\coloneqq\2$.
	\item \label{exc.p_tree_sets.unary} $p\coloneqq\yon$.
	\item \label{exc.p_tree_sets.binary} $p\coloneqq\yon^\2$.
	\item $p\coloneqq{\2\yon}$.
	\item $p\coloneqq{\yon+\1}$.
	\item \label{exc.p_tree_sets.monomial} $p\coloneqq{B\yon^A}$ for some sets $A,B\in\smset$.
\qedhere
\end{enumerate}
\begin{solution}
\begin{enumerate}
    \item To select a $\1$-tree, we must select a position from $\1$\dots and then we are done, because there are no directions. So there is a unique $\1$-tree given by a root and no edges, implying that $\tr_\1\iso\1$.
    \item To select a $\2$-tree, we must select a position from $\2$.
    Then we are done, because there are no directions.
    So every $\2$-tree is a root with no edges, and the root is labeled with one of the elements of $\2$.
    Hence $\tr_\2\iso\2$.
    \item To select a $\yon$-tree, we must select a position from $\1$; then for the unique direction at that position, we must select a position from $\1$; and so forth.
    Since we only ever have $1$ position to choose from, there is only $1$ such $\yon$-tree we can build: an infinite ray extending from the root, where every vertex has exactly $1$ child.
    Hence $\tr_\yon\iso\1$.
    \item To select a $\yon^\2$-tree, we must select a position from $\1$; then for each direction at that position, we must select a position from $\1$; and so forth.
    We only ever have $1$ position to choose from, so there is only $1$ such $\yon^\2$-tree: an infinite binary tree, where every vertex has exactly $2$ children.
    Hence $\tr_{\yon^\2}\iso\1$.
    \item To select a $\2\yon$-tree, we must select a position from $\2$; then for the unique direction at that position, we must select a position from $\2$; and so forth, one position out of $\2$ for each level of the tree.
    So every $\2\yon$-tree is an infinite ray, where every vertex is labeled with an element of $\2$ and has exactly $1$ child.
    There is then a unique vertex of height $n$ for each $n\in\nn$, so there is a bijection between $2\yon$-trees $T$ and functions $f\colon\nn\to\2$, where $f(n)\in\2$ is the label of the height-$n$ vertex of $T$.
    Hence $\tr_{\2\yon}\iso\2^\nn$.
    \item To select a $(\yon+\1)$-tree, we choose either a position with $1$ direction or a position with no directions: if the position has no directions, we stop, and if the position has $1$ direction, we repeat our choice for the next level.
    So a choice of $(\yon+\1)$-tree is equivalent to a choice of when, if ever, to stop.
    That is, for each $n\in\nn$, there is a unique $(\yon+\1)$-tree of height-$n$ consisting of $n+1$ vertices along a single path, so that every vertex aside from the height-$n$ leaf has exactly $1$ child.
    Then there is exactly one more $(\yon+\1)$-tree: an infinite ray, obtained by always picking the position with $1$ directions and never the one with no directions.
    Hence we could write $\tr_{\yon+\1}\iso\nn\cup\{\infty\}$ (although this set is in bijection with $\nn$).
    \item To select a $B\yon^A$-tree for fixed $A,B\in\smset$, we must select a position from $B$; then for each direction in $A$ at that position, we must select a position from $B$; and so forth.
    This yields a tree in which every vertex has $A$-many children and a label from $B$.
    Such a tree has $1$ root and $|A|$ times as many vertices in one level than the level below it, so its height-$n$ vertices for each $n\in\nn$ are in bijection with the set $A^\ord{n}$: each $n$-tuple in $A^\ord{n}$ gives a length-$n$ rooted path of directions in $A$ to follow up the tree, uniquely specifying a height-$n$ vertex.
    So the set of all vertices is given by $\sum_{n\in\nn}A^\ord{n}=\lst(A)$.
    Then specifying a $B\yon^A$-tree amounts to assigning a label from $B$ to every vertex via a function $\lst(A)\to B$.
    Hence $\tr_{B\yon^A}\iso B^{\lst(A)}$.
\end{enumerate}
\end{solution}
\end{exercise}

\begin{exercise} \label{exc.n_ary_trees}
\begin{enumerate}
    \item Say we were interested in \emph{$\ord{n}$-ary trees}: infinite (unless $n=0$) rooted trees in which every vertex has $n$ children.
    Is there a polynomial $p$ for which $\tr_p$ is the set of $\ord{n}$-ary trees?
    \item \label{exc.n_ary_trees.l_labeled} Now say we wanted to assign each vertex of an $\ord{n}$-ary tree a label from a set $L$.
    Is there a polynomial $q$ for which $\tr_q$ is the set of $L$-labeled $\ord{n}$-ary trees?\qedhere
\end{enumerate}
\begin{solution}
\begin{enumerate}
    \item Yes: by \cref{exc.p_tree_sets} \cref{exc.p_tree_sets.monomial} in the case of $B=1$, or by analogy with \cref{exc.p_tree_sets} \cref{exc.p_tree_sets.unary,exc.p_tree_sets.binary}, we have that $\tr_{\yon^\ord{n}}$ is the set of $\ord{n}$-ary trees.
    \item Yes: by \cref{exc.p_tree_sets} \cref{exc.p_tree_sets.monomial}, we have that $\tr_{L\yon^\ord{n}}$ is the set of $L$-labeled $\ord{n}$-ary trees.
\end{enumerate}
\end{solution}
\end{exercise}

From here, a natural question to ask is how the set $\tr_p$ of $p$-trees is related to the set of $p\tripow{n}$-pretrees $p\tripow{n}(\1)$ for each $n\in\nn$.
A look back at \cref{def.tripow_pretree} gives a clue: every $p$-tree has a stage-$n$ pretree obtained by removing all vertices of height greater than $n$, yielding a function $\tr_p\to p\tripow{n}(\1)$ that we will denote by $\pi^{(n)}$.

Moreover, since the stage-$n$ pretree of a given $p$-tree agrees with the stage-$n$ pretree of any stage-$(n+k)$ pretree of the $p$-tree, the functions $\pi^{(n)}$ should commute with our tree-trimming functions $p\tripow{(n+k)}(\1)\to p\tripow{n}(\1)$ from \eqref{eqn.trim} in \cref{ex.trim}.
In particular, the following diagram commutes for all $n\in\nn$:
\[
\begin{tikzcd}[column sep=large]
    \tr_p \ar[d,"\pi^{(n)}"'] \ar[dr,"\pi^{(n+1)}"]%,outer sep=-2pt]
    \\
    p\tripow{n}(\1) & p\tripow{(n+1)}(\1). \ar[l,"p\tripow{n}(!)"]
\end{tikzcd}
\]
All this says is that if we trim a $p$-tree down until only $n+1$ levels are left via $\pi^{(n+1)}$, then trimmed off one more level via $p\tripow{n}(!)$, it would be the same as if we had trimmed it down to $n$ levels from the start via $\pi^{(n)}$.
This is summarized by the following larger commutative diagram, which contains every function of the form \cref{ex.trim}:
\begin{equation} \label{eqn.cofree_pos_with_lim}
\begin{tikzcd}%[sep=large]
    \tr_p \ar[dd,"\pi^{(0)}"'] \ar[ddr,"\pi^{(1)}"',outer sep=-2pt] \ar[ddrr,"\pi^{(2)}"',outer sep=-1.5pt] \ar[ddrrrr,"\pi^{(3)}"',outer sep=-1.5pt]\\
    &&&\cdots\\
    \1 &
	p(\1) \ar[l, "!"] &
	p\tripow2(\1) \ar[l, "p(!)"] &&
% 	[10pt]
	p\tripow3(\1) \ar[ll, "p\tripow2(!)"] &
	\cdots.\ar[l]
\end{tikzcd}
\end{equation}
% \[
% \begin{tikzcd}[sep=large]
%     \tr_p \ar[dd,"\pi^{(0)}"'] \ar[ddr,"\pi^{(1)}"',outer sep=-2pt] \ar[ddrr,"\pi^{(2)}"',outer sep=-1.5pt] \ar[ddrrrr,"\pi^{(3)}"',outer sep=-1.5pt]\\
%     &&&\ldots\\
%     \1 &
% 	p(\1) \ar[l, "!"] &
% 	p\tri p(\1) \ar[l, "p(!)"] &&
% % 	[10pt]
% 	p\tri p\tri p(\1) \ar[ll, "p\:\tri\:p(!)"] &
% 	\cdots.\ar[l]
% \end{tikzcd}
% \]
So $\tr_p$ with the functions $\pi^{(n)}$ forms a cone over the the bottom row---in fact, it is the universal cone.
Intuitively, this is because a $p$-tree carries exactly the information that a compatible sequence of $p\tripow{n}$-pretrees does: no more, no less.
But we can prove it formally as well.

\begin{exercise}\label{exc.tree_as_limit}
Prove that $\tr_p$ with the functions $\pi^{(n)}\colon\tr_p\to p\tripow{n}(\1)$ is the limit of the diagram
\begin{equation} \label{eqn.cofree_pos_diagram}
\begin{tikzcd}
    \1 &
	p(\1) \ar[l, "!"'] &
	p\tripow2(\1) \ar[l, "p(!)"'] &
 	[10pt]
	p\tripow3(\1) \ar[l, "p\tripow2(!)"'] &
	\cdots.\ar[l]
\end{tikzcd}
\end{equation}
You may use the fact that \eqref{eqn.cofree_pos_with_lim} commutes.
\begin{solution}
We defined $\tr_p$ in \cref{def.poly_tree} as a rooted tree $T$ where every vertex $v$ is labeled with an element of $p(1)$, together with a choice of bijection between the set of $v$'s children and the set $p[i]$. We need to translate between that description and the description as the set-theoretic limit
\[
1\from p(1)\from p(p(1))\from\cdots
\]
An element of this set consists of an element of $1$, which is determined, an element $i\in p(1)$, an element $(i_1,i_2)\in\sum_{i\in p(1)}\prod_{d\in p[i]}p(1)$ that ``agrees'' with $i$, an element of $\sum_{i_1\in p(1)}\prod_{d_1\in p[i_1]}\sum_{i_2\in p(1)}\prod_{d_2\in p[i_2]}p(1)$ that agrees, etc. But all this agreement is easy to simplify: we need elements
\begin{align*}
	 i_1&\in p(1),\\
	 i_2&\in\prod_{d_1\in p[i_1]}p(1),\\
	 i_3&\in \prod_{d_1\in p[i_1]}\prod_{d_2\in p[i_2 d_1]}p(1),\\
	 i_4&\in\prod_{d_1\in p[i_1]}\prod_{d_2\in p[i_2 d_1]}\prod_{d_3\in p[i_3 d_1 d_2]}p(1)\\
	 &\vdots
\end{align*}

To go from the former description to the latter, take such a rooted tree $T$ and let $i_1$ be the label on the root vertex. Since by assumption there is a bijection between the set children of that vertex and the set $p[i_1]$, take $i_2\colon p[i_1]\to p(1)$ to assign to each $d_1\in p[i_1]$ the label on the associated child vertex. Repeat this indefinitely.

To go from the latter description to the former, simply invert the description in the previous paragraph.
\end{solution}
\end{exercise}

\index{coalgebra!terminal}\index{functor!coalgebra for}

If you know about coalgebras for functors, as mentioned in \cref{ex.coalgebras}, then you might know the limit $\car{t}_p(\1)$ of \eqref{eqn.cofree_pos_diagram} by a different name: it is the \emph{terminal coalgebra for the functor} $p$, or the \emph{terminal $p$-coalgebra}, because it is terminal in the category of $p$-coalgebras, as the following exercise shows.

\index{coalgebra!morphism of coalgebras}
\begin{exercise}
A \emph{$p$-coalgebra morphism} between $p$-coalgebras $\varphi\colon S\to p(S)$ and $\psi\colon T\to p(T)$ (as in \cref{ex.coalgebras}) is a function $f\colon S\to T$ such that the square
\[
\begin{tikzcd}
    S \ar[r,"\varphi"] \ar[d,"f"'] & p(S) \ar[d,"p(f)"] \\
    T \ar[r,"\psi"] & p(T)
\end{tikzcd}
\]
commutes.
\begin{enumerate}
	\item Choose what you think is a good function $\tr_p\to p(\tr_p)$.
	\item Show that $\tr_p$ equipped with your function $\tr_p\to p(\tr_p)$ is the terminal object in the category of $p$-coalgebras and the morphisms between them.\index{category!of coalgebras}\index{coalgebra!category of}
	\item Show that the function $\tr_p\to p(\tr_p)$ you chose is a bijection.
\end{enumerate}
\begin{solution}
\begin{enumerate}
	\item To give a function $t\colon\tr_p\to p(\tr_p)$ we need to choose, for each tree $T\in \tr_p$, an element $i\in p(1)$ and a function $T.-\colon p[i]\to\tr_p$. Choose $i$ to be the label on the root of $T$, and choose $T.d$ to be the tree obtained by following the branch labeled $d$ for any $d\in p[i]$.
	\item Given an arbitrary coalgebra $(S,f)$, where $S\in\smset$ and $f\colon S\to p\tri S$, we need to show that there is a unique function $u\colon S\to\tr_p$ that commutes with $f$ and $t$, i.e.\ with $t\circ u=p(f)$.

	For each $n$, we have a function $f_n\colon S\to p\tripow{n}\tri S$ by induction: for $n=0$ use the identity function $S\to S$ and given $S\to p\tripow{n}\tri S$ we compose with $p\tripow{n}\tri f$ to get $S\to p\tripow{n}\tri p\tri S$. Let $f_n'\colon S\to p\tripow{n}(1)$ be given by $f_n'\coloneqq f_n\then(p\tripow{n}\tri!$. We know by \cref{exc.tree_as_limit} that $\tr_p$ is the limit of $p\tripow{n}(1)$, so we get a unique map $(f_n')_{n\in\nn}\colon S\to\tr_p$. It commutes with $f$ and $t$ by construction, and the fact that it is constructed using the universal property of limits implies that it is appropriately unique.
		\item The function $\tr_p\to p\tri\tr_p$ from the first part is invertible. Indeed, given $i\in p(1)$ and $T'\colon p[i]\to\tr_p$, we construct an element of $\tr_p$ by taking its root to be labeled with $i$, its set of children to be $p[i]$, and for each $d\in p[i]$ taking the remaining tree to be $T'(d)$.
\end{enumerate}
\end{solution}
\end{exercise}

\subsubsection{Positions of the cofree comonoid}

\index{cofree comonoid!positions of}

Now recall from \cref{ex.compute_limits} that in $\poly$, the positions of a limit are the limit of the positions.
Moreover, in \eqref{eqn.cofree_diagram}, every vertical lens $p\tripow{n}\to p\tripow{n}\tri\1$ is, in fact, a \emph{vertical} lens in the sense of \cref{def.vert_cart}: an isomorphism on positions.
So \eqref{eqn.cofree_diagram} on positions collapses down to its bottom row, viewed as a diagram in $\smset$.
Yet this is precisely the diagram from \eqref{eqn.cofree_pos_diagram}, whose limit is $\tr_p$.
So $\tr_p$ is the position-set of the limit $\car{t}_p$ of \eqref{eqn.cofree_diagram}: we write $\car{t}_p(\1)\iso\tr_p$.

So, as we said above, $\car{t}_p(\1)$ is the set carrying the terminal coalgebra of $p$, but  we prefer to think of it as the set of $p$-trees, for it gives us a concrete way to realize the limit and its projections, as well as a natural interpretation of the directions of $\car{t}_p$ at each position.\index{coalgebra!terminal}

\subsubsection{The directions of the cofree comonoid}
\index{cofree comonoid!directions of}

Given that $\car{t}_p$ is the limit of \eqref{eqn.cofree_diagram}, \cref{ex.compute_limits} also tells us how to compute its directions: the directions of a limit are the colimit of the directions.\index{colimit}
But every polynomial in the bottom row of \eqref{eqn.cofree_diagram} has an empty direction-set, and there are no arrows between polynomials in the top row.
So the directions of $\car{t}_p$ are given by a coproduct of the directions of each $p\tripow{n}$.

More precisely, given a $p$-tree $T\in\tr_p$ whose stage-$n$ pretree for $n\in\nn$ is $\pi^{(n)}T$, the direction-set $\car{t}_p[t]$ is given by the following coproduct:
\[
    \car{t}_p[T]\coloneqq\sum_{n\in\nn}p\tripow{n}[\pi^{(n)}T].
\]
But by \cref{def.tripow_pretree}, each $p\tripow{n}[\pi^{(n)}T]$ is the set of height-$n$ leaves of the $p\tripow{n}$-pretree $\pi^{(n)}T$, which in turn is the stage-$n$ pretree of $T$.
So its height-$n$ leaves coincide with the height-$n$ vertices of $T$.
Therefore, we can identify $p\tripow{n}[\pi^{(n)}T]$ with the set of height-$n$ vertices of $T$; we denote this set by $\vtx_n(T)$.

Since the coproduct above ranges over all $n\in\nn$, it follows that $\car{t}_p[T]$ is the set of \emph{all} vertices in $T$; we denote this set by $\vtx(T)$.
The on-directions function at $T$ of each projection $\car{t}_p\to p\tripow{n}$ from the limit must then be the canonical inclusion
\[
    p\tripow{n}[\pi^{(n)}T]\iso \vtx_n(T)\inj \vtx(T)\iso\car{t}_p[T],
\]
sending height-$n$ leaves of $\pi^{(n)}T$ to height-$n$ vertices of $T$.

Here is alternative way to think about the directions of $\car{t}_p$ and each $p\tripow{n}$ that will be helpful.
A defining feature of a rooted tree is that its vertices are in bijection with its finite rooted paths: each vertex gives rise to a unique path to that vertex from the root, and every finite rooted path arises this way.
So the directions of $p\tripow{n}$ at a given $p\tripow{n}$-pretree correspond in turn to the rooted paths of that pretree leading to its height-$n$ leaves; indeed, each such direction is comprised of a sequence of $n$ directions of $p$, which together specify a length-$n$ rooted path up the pretree.
Then for $T\in\tr_p$, the direction-set $\car{t}_p[T]$ consists of every finite rooted path of $T$.

Since only finite rooted paths correspond to vertices, all our paths will be assumed to be finite from here on out. This is our preferred way to think about directions in $\car{t}_p$.
When we wish to refer to what one might call an infinite (rooted) ``path,'' we will instead call it a (rooted) \emph{ray}.

\begin{exercise}\label{exc.p_tree_polys}
For each of the following polynomials $p$, characterize the polynomial $\car{t}_p$.
You may choose to think of the directions of $\car{t}_p$ either as vertices or as rooted paths.
(Note that these are the same polynomials from \cref{exc.p_tree_sets}.)
\begin{enumerate}
	\item \label{exc.p_tree_polys.1} $p\coloneqq\1$.
	\item $p\coloneqq\2$.
	\item \label{exc.p_tree_polys.unary} $p\coloneqq\yon$.
	\item \label{exc.p_tree_polys.binary} $p\coloneqq\yon^\2$.
	\item \label{exc.p_tree_polys.2_label_ray} $p\coloneqq{\2\yon}$.
	\item \label{exc.p_tree_polys.maybe_unary} $p\coloneqq{\yon+\1}$.
	\item \label{exc.p_tree_polys.monomial} $p\coloneqq{B\yon^A}$ for some sets $A,B\in\smset$.
\qedhere
\end{enumerate}
\begin{solution}
For the given values of $p$, we have already characterized the position-set $\car{t}_p(\1)\iso\tr_p$ of $\car{t}$ in \cref{exc.p_tree_sets}, so it remains to characterize the directions at each position, i.e.\ the vertices of each $p$-tree.
\begin{enumerate}
    \item We saw that $\tr_\1\iso\1$ and that the unique $\1$-tree has a single vertex, its root.
    Equivalently, it has only $1$ rooted path: the empty path from the root to itself.
    So $\car{t}_\1\iso\1\yon^\1\iso\yon$.
    \item We saw that $\tr_\2\iso\2$ and that each $\2$-tree has a single vertex, its root.
    Equivalently, its only rooted path is the empty path.
    So $\car{t}_\2\iso\2\yon^\1\iso\2\yon$.
    \item We saw that $\tr_\yon\iso\1$ and that the unique $\yon$-tree is a single ray extending from the root, where every vertex has exactly $1$ child.
    So there is exactly $1$ vertex of height-$n$ for each $n\in\nn$, yielding a bijection between the vertices of the $\yon$-tree to $\nn$.
    Equivalently, the ray has exactly $1$ rooted path of length $n$ for each $n\in\nn$.
    So $\car{t}_\yon\iso\yon^\nn$.
    \item We saw that $\tr_{\yon^\2}\iso\1$ and that the unique $\yon^\2$-tree is an infinite binary tree, where every vertex has exactly $2$ children.
    So the vertices of height-$n$ are in bijection with $\2^\ord{n}$, yielding a bijection between the vertices of the $\yon^\2$-tree to $\sum_{n\in\nn}\2^\ord{n}\iso\lst(\2)$.
    Equivalently, the rooted paths of an infinite binary tree are just finite binary sequences, which comprise the set $\lst(\2)$.
    Hence $\car{t}_{\yon^\2}\iso\yon^{\lst(\2)}$.
    \item We saw that $\tr_{\2\yon}\iso\2^\nn$, and that every $2\yon$-tree is an infinite ray, whose vertices (or rooted paths) are in bijection with $\nn$.
    Hence $\car{t}_{\2\yon}\iso\2^\nn\yon^\nn$.
    \item We saw that $\tr_{\yon+\1}\iso\nn\cup\{\infty\}$.
    Here the $(\yon+\1)$-tree corresponding to $n\in\nn$ consists of $n+1$ vertices along a single path, so that every vertex aside from the height-$n$ leaf has exactly $1$ child; and the $(\yon+\1)$-tree corresponding to $\infty$ is an infinite ray.
    Thus the direction-set of $\car{t}_{\yon+\1}$ at $n\in\nn$ is $\ord{n}+\1$, while its direction-set at $\infty$ is $\nn$.
    Hence $\car{t}_{\yon+\1}\iso\{\infty\}\yon^\nn+\sum_{n\in\nn}\yon^{\ord{n}+\1}$.
    \item We saw that $\tr_{B\yon^A}\iso B^{\lst(A)}$, where the height-$n$ vertices of a given $B\yon^A$-tree are in bijection with the set $A^\ord{n}$, so the set of all vertices of a $B\yon^A$-tree is given by $\sum_{n\in\nn}A^\ord{n}=\lst(A)$.
    Equivalently, the rooted paths of the $B\yon^A$-tree are just finite sequences in $A$, which comprise the set $\lst(A)$.
    Hence $\car{t}_{B\yon^A}\iso B^{\lst(A)}\yon^{\lst(A)}$.
\end{enumerate}
\end{solution}
\end{exercise}

We summarize the results of this section in the following proposition, thus concretely characterizing the carrier $\car{t}_p$ of the cofree comonoid on $p$ in terms of $p$-trees.
For reasons that will become clear shortly, we will denote each projection from the limit of \eqref{eqn.cofree_diagram} by $\epsilon^{(n)}_p\colon\car{t}_p\to p\tripow{n}$. We denote $\epsilon$ simply by $\epsilon$.

\begin{proposition} \label{prop.cofree_carrier_lim}
For $p\in\poly$, let
\[
    \car{t}_p\coloneqq\sum_{T\in\tr_p}\yon^{\vtx(T)}
\]
be the polynomial whose positions are $p$-trees and whose directions at each $p$-tree are the rooted paths.
Then $\car{t}_p$ is the limit of the diagram \eqref{eqn.cofree_diagram}, with projections $\epsilon^{(n)}_p\colon\car{t}_p\to p\tripow{n}$ for every $n\in\nn$ making the following diagram commute:
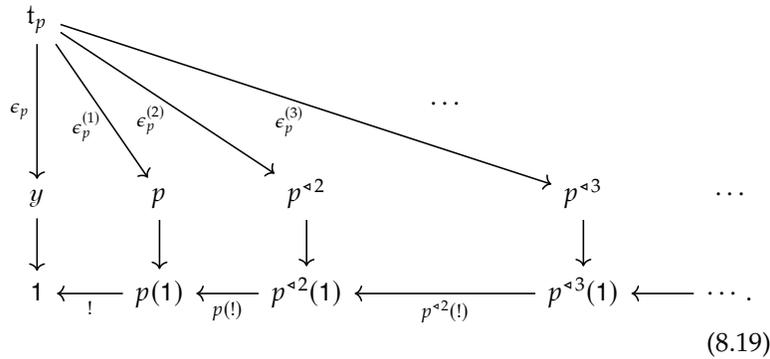
\begin{equation} \label{eqn.cofree_with_lim}
\begin{tikzcd}%[sep=large]
    \car{t}_p \ar[dd,"\epsilon_p"'] \ar[ddr,"\epsilon^{(1)}_p"',outer sep=-2pt] \ar[ddrr,"\epsilon^{(2)}_p"',outer sep=-1.5pt] \ar[ddrrrr,"\epsilon^{(3)}_p"',outer sep=-1.5pt]\\
    &&&\cdots\\
    \yon \ar[d] &
	p \ar[d] &
	p\tripow2 \ar[d] &&
% 	[10pt]
	p\tripow3 \ar[d] &
	\cdots \\
    \1 &
	p(\1) \ar[l, "!"] &
	p\tripow2(\1) \ar[l, "p(!)"] &&
% 	[10pt]
	p\tripow3(\1) \ar[ll, "p\tripow2(!)"] &
	\cdots.\ar[l]
\end{tikzcd}
\end{equation}
The lens $\epsilon^{(n)}_p\colon\car{t}_p\to p\tripow{n}$ sends each $p$-tree $T\in\tr_p$ to its stage-$n$ pretree $\pi^{(n)}T$ on positions and each height-$n$ leaf of $\pi^{(n)}T$ to the corresponding height-$n$ rooted path of $T$ on directions.
\end{proposition}

\begin{example}[Drawing $\epsilon^{(n)}_p$ in polyboxes] \label{ex.eps_n_polybox}
Here is $\epsilon^{(n)}_p\colon\car{t}_p\to p\tripow{n}$ drawn in polyboxes, where we continue to denote the stage-$n$ pretree of a $p$-tree by $\pi^{(n)}T$:
\begin{equation} \label{eqn.eps_n_polybox}
\begin{tikzpicture}
  \node (f) {
    \begin{tikzpicture}[polybox, mapstos]
  	  \node[poly, dom, "$\car{t}_p$" below] (p) {$v$\at$T$};
  	  \node[left=0pt of p_pos] {$\tr_p$};
  	  \node[left=0pt of p_dir] {$\vtx(-)$};

  	  \node[poly, cod, right=of p, "$p\tripow{n}$" below] (q) {$v$\at$\pi^{(n)}T$};
  	  \node[right=0pt of q_pos] {$p\tripow{n}(\1)$};
	  \node[right=0pt of q_dir] {$\vtx_n(-)$};

  	  \draw (p_pos) -- node[below] {} (q_pos);
  	  \draw (q_dir) -- node[above] {} (p_dir);
    \end{tikzpicture}
  };
\end{tikzpicture}
\end{equation}
On the right hand side, $v$ is a height-$n$ leaf of $\pi^{(n)}T$; on the left, $v$ is identified with the corresponding height-$n$ vertex of $T$.

This isn't the only way we can write this lens in polyboxes, however; polyboxes have special notation for lenses to composites, allowing us to write, say, the $n=4$ case like so:
\begin{equation} \label{eqn.eps_1_1_1_1_polybox}
\begin{tikzpicture}[polybox, mapstos]
	\node[poly, dom, "$\car{t}_p$" left] (p) {$a_1\leadsto a_2\leadsto a_3\leadsto a_4$\at$T$};
	\foreach \i/\j in {1/0,2/1,3/2,4/3}
	{
  	\node[poly, cod, "$p$" right] (q\i) at (4,1.3*\i-3.25) {$a_\i$\at$i_\j$};
	};
	\draw (p_pos) to[first] node[below] {} (q1_pos.west);
	\foreach \i/\j in {1/2,2/3,3/4}
	{
		\draw
			(q\i_dir.west)
			to[climb]
			node[left] {}
			(q\j_pos.west);
	};
	\draw (q4_dir) to[last] node[above left] {} (p_dir);
\end{tikzpicture}
\end{equation}
Here we are unpacking the construction of the first $4$ levels of $T$, according to our nested instructions for building $p$-trees.
The $p$-position $i_0$ is the label on the root of $T$, while the $p[i_0]$-direction $a_1$ specifies one of the edges coming out of it---leading to a height-$1$ vertex of $T$ labeled $i_1$, and so on.

The contents of the position-boxes and the arrows going up on the codomain side carry all the data of the bottom $4$ levels of $T$: namely the label $i_0$ of the root, the label $i_1$ of the vertex at the end of every direction $a_1$ out of the root, and so on until $i_3$.
All this specifies a unique $p\tripow4$-position, a $p\tripow4$-pretree, which is the same as the $p\tripow4$-pretree $T$.
Indeed, we can think of $\pi^{(4)}T$ as a shorthand for the gadget comprised of the $4$ polyboxes on the right hand side of \eqref{eqn.eps_1_1_1_1_polybox} when their blue boxes are yet to be filled.
So the position and direction depicted on the codomain side of \eqref{eqn.eps_1_1_1_1_polybox} is equivalent to the position and direction depicted on the codomain side of \eqref{eqn.eps_n_polybox}.
This generalizes to all values of $n$. The direction $a_4$ emanating from $i_3$, together with all the data below it, corresponds to the rooted path we've denoted $a_1\leadsto a_2\leadsto a_3\leadsto a_4$ in $T$.

There are even more ways to express $\epsilon^{(4)}_p\colon\car{t}_p\to p\tripow4$ in polyboxes, however.
After all, $p\tripow4\iso p\tripow2\tri p\tripow2$.
So we ought to be able to draw $\epsilon^{(4)}_p$ as follows:
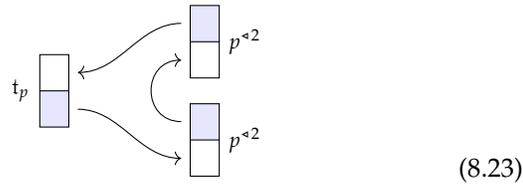
\begin{equation} \label{eqn.eps_2_2_polybox_empty}
\begin{tikzpicture}[polybox, tos]
    \node[poly, dom, "$\car{t}_p$" left] (r) {};
    \node[poly, cod, right=1.8 of r.south, yshift=-1ex, "$p\tripow2$" right] (q) {};
    \node[poly, cod, above=of q, "$p\tripow2$" right] (q') {};

    \draw (r_pos) to[first] node[below] {} (q_pos);
    \draw (q_dir) to[climb] node[right] {} (q'_pos);
    \draw (q'_dir) to[last] node[above] {} (r_dir);
\end{tikzpicture}
\end{equation}
Let's think about how we should fill these boxes.
We can still put a $p$-tree $T$ in the lower left position box.
From the commutativity of \eqref{eqn.cofree_with_lim}, we know the on-positions function $\pi^{(4)}$ of $\epsilon^{(4)}_p$ factors through the on-positions function $\pi^{(2)}$ of $\epsilon^{(2)}_p$, which tells us that $p\tripow2$-pretree that should go in the lower right position box should be $\pi^{(2)}T$: the $p\tripow2$-pretree of $T$.

Another way to think about this is that the polyboxes for the lower $p\tripow2$ on the right side of \eqref{eqn.eps_2_2_polybox_empty} are equivalent to the polyboxes for the $2$ lower copies of $p$ on the right side of \eqref{eqn.eps_1_1_1_1_polybox}---just like how the polyboxes for $p\tripow{n}$ in \eqref{eqn.eps_n_polybox} are equivalent to the polyboxes for all $n=4$ copies of $p$ in \eqref{eqn.eps_1_1_1_1_polybox}.
There, the position could be represented with a single $p\tripow{n}$-pretree $\pi^{(n)}T$, and the direction $(a_1,\ldots,a_n)$ is one of its length-$n$ rooted paths.
So here, too, we can package the $2$ lower copies of $p$ into a single pair of polyboxes for $p\tripow2$, if we let $u$ be the height-$2$ vertex of $T$ at the end of the rooted path $(a_1,a_2)$:
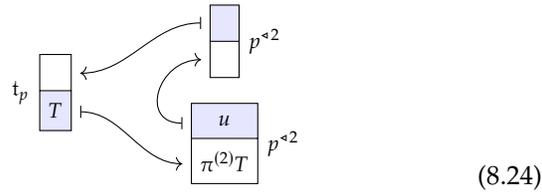
\begin{equation} \label{eqn.eps_2_2_polybox_partly}
\begin{tikzpicture}[polybox, mapstos]
    \node[poly, dom, "$\car{t}_p$" left] (r) {\at$T$};
    \node[poly, cod, right=1.8 of r.south, yshift=-1ex, "$p\tripow2$" right] (q) {$u$\at$\pi^{(2)}T$};
    \node[poly, cod, above=of q, "$p\tripow2$" right] (q') {};

    \draw (r_pos) to[first] node[below] {} (q_pos);
    \draw (q_dir) to[climb] node[right] {} (q'_pos);
    \draw (q'_dir) to[last] node[above] {} (r_dir);
\end{tikzpicture}
\end{equation}
But then the polyboxes for the $2$ upper copies of $p$ from \eqref{eqn.eps_1_1_1_1_polybox} should also collapse down to a single pair of polyboxes for $p\tripow2$, with a $p\tripow2$-pretree as its position and a height-$2$ leaf of that pretree as its direction.
Indeed, once we have followed the directions $(a_1,a_2)$ up to the height-$2$ vertex $u$, the subtree of $T$ rooted at $u$ is itself a $p$-tree: it has a label $i_2$ on its root, out of which we can follow the direction $a_3$ to reach a height-$1$ vertex labeled $i_3$, then follow the direction $a_4$ to reach a height-$2$ vertex that we call $w$.
Of course, the vertex labeled $i_3$ is actually a height-$3$ vertex of the whole tree $T$; likewise $w$ corresponds to a height-$4$ vertex of $T$.
However, from the perspective of the upper copy of $p\tripow2$ in \eqref{eqn.eps_2_2_polybox_partly}, we are starting over from $u$ and moving up the subtree of $T$ rooted at $u$---or, more precisely, the $p\tripow2$-pretree of the subtree rooted at $u$.
So the polyboxes end up looking like this:
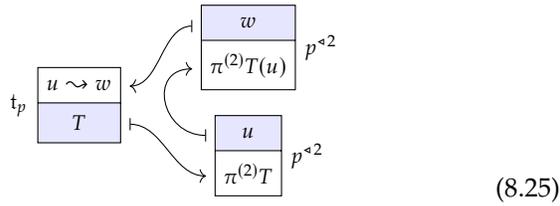
\begin{equation} \label{eqn.eps_2_2_polybox}
\begin{tikzpicture}[polybox, mapstos]
    \node[poly, dom, "$\car{t}_p$" left] (r) {$u\leadsto w$\at$T$};
    \node[poly, cod, right=1.8 of r.south, yshift=-1ex, "$p\tripow2$" right] (q) {$u$\at$\pi^{(2)}T$};
    \node[poly, cod, above=of q, "$p\tripow2$" right] (q') {$w$\at$\pi^{(2)}T(u)$};

    \draw (r_pos) to[first] node[below] {} (q_pos);
    \draw (q_dir) to[climb] node[right] {} (q'_pos);
    \draw (q'_dir) to[last] node[above] {} (r_dir);
\end{tikzpicture}
\end{equation}
Here $T(u)$ denotes the $p$-tree equal to the subtree of $T$ rooted at its vertex $u$, and $w$ is one of $T(u)$'s height-$2$ vertices.
When viewed as the subtree of $T$ rooted at its the height-$2$ vertex $u$, the $p$-tree $T(u)$ has its height-$2$ vertex $w$ identified with a height-$4$ vertex of $T$, a descendent of $u$ that we denote by $u\leadsto w$.
Alternatively, $u$ corresponds to the rooted path $(a_1,a_2)$ of $T$, while $w$ corresponds to the rooted path $(a_3,a_4)$ of $T(u)$, so $u\leadsto w$ corresponds to the concatenated rooted path $(a_1,a_2,a_3,a_4)$ of $T$.
As concatenation is associative, $\leadsto$ is associative as well.

Both the notation $T(u)$, for what we will call the $p$-\emph{subtree} of the $p$-tree $T$ rooted at $u\in\vtx(T)$, and the notation $u\leadsto w$, for the vertex in $\vtx(T)$ that $w\in\vtx(T(u))$ coincides with when $T(u)$ is identified with the subtree of $p$ rooted at $u$, will turn out to be temporary---we will soon justify why we can express these concepts in much more familiar terms.%
\tablefootnote{As a reminder, due to the definition of a $p$-tree $T$, a $p$-subtree of $T$ is still itself a $p$-tree, whereas the $p\tripow{n}$-pretree of $T$ is not a $p$-tree unless the height of $T$ is strictly less than $n$.}
As a sneak preview, they will be crucial to defining the categorical structure of the cofree comonoid on $p$.

Here are three more ways to depict $\epsilon^{(4)}_p$, viewing its codomain in turn as $p\tripow3\tri p\tripow1$, as $p\tripow1\tri p\tripow2\tri p\tripow1$, or as $p\tripow1\tri p\tripow1\tri p\tripow1\tri p\tripow1$:
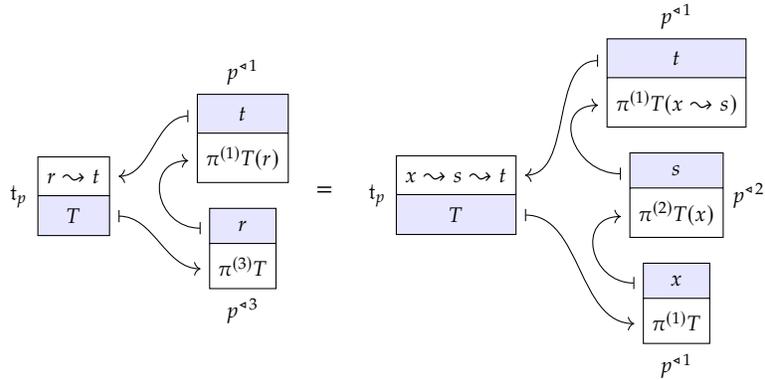
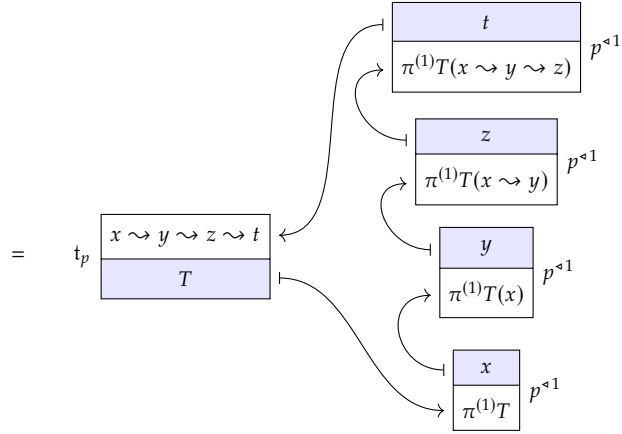
\[
\begin{tikzpicture}
\node (1) {
\begin{tikzpicture}[polybox, mapstos]
    \node[poly, dom, "$\car{t}_p$" left] (r) {$r\leadsto t$\at$T$};
    \node[poly, cod, right=1.8 of r.south, yshift=-1ex, "$p\tripow3$" below] (q) {$r$\at$\pi^{(3)}T$};
    \node[poly, cod, above=of q, "$p\tripow1$" above] (q') {$t$\at$\pi^{(1)}T(r)$};

    \draw (r_pos) to[first] node[below] {} (q_pos);
    \draw (q_dir) to[climb] node[right] {} (q'_pos);
    \draw (q'_dir) to[last] node[above] {} (r_dir);
\end{tikzpicture}
};
\node[right=.7 of 1] (2) {
\begin{tikzpicture}[polybox, mapstos]
	\node[poly, dom, "$\car{t}_p$" left] (S) {$x\leadsto s\leadsto t$\at$T$};
	\node[poly, cod, right=of S, "$p\tripow2$" right] (S2) {$s$\at$\pi^{(2)}T(x)$};
	\node[poly, cod, below=of S2, "$p\tripow1$" below] (S1) {$x$\at$\pi^{(1)}T$};
	\node[poly, cod, above=of S2, "$p\tripow1$" above] (S3) {$t$\at$\pi^{(1)}T(x\leadsto s)$};
	\draw (S_pos) to[first] (S1_pos);
	\draw (S1_dir) to[climb] (S2_pos);
	\draw (S2_dir) to[climb] (S3_pos);
	\draw (S3_dir) to[last] (S_dir);
\end{tikzpicture}
};
\node at ($(1.east)!.5!(2.west)$) {=};
\end{tikzpicture}
\]
\[
\begin{tikzpicture}[polybox, mapstos]
	\node[poly, dom, "$=\qquad\car{t}_p$" left] (p) {$x\leadsto y\leadsto z\leadsto t$\at$T$};
  	\node[poly, cod, "$p\tripow1$" right] (q1) at (4,1.5*1-3.25) {$x$\at$\pi^{(1)}T$};
  	\node[poly, cod, "$p\tripow1$" right] (q2) at (4,1.5*2-3.25) {$y$\at$\pi^{(1)}T(x)$};
  	\node[poly, cod, "$p\tripow1$" right] (q3) at (4,1.5*3-3.25) {$z$\at$\pi^{(1)}T(x\leadsto y)$};
  	\node[poly, cod, "$p\tripow1$" right] (q4) at (4,1.5*4-3.25) {$t$\at$\pi^{(1)}T(x\leadsto y\leadsto z)$};

	\draw (p_pos) to[first] node[below] {} (q1_pos.west);
	\foreach \i/\j in {1/2,2/3,3/4}
	{
		\draw
			(q\i_dir.west)
			to[climb]
			node[left] {}
			(q\j_pos.west);
	};
	\draw (q4_dir) to[last] node[above left] {} (p_dir);
\end{tikzpicture}
\]
Here $r=x\leadsto s$ and $s=y\leadsto z$.
Notice that the last depiction is just another way to write \eqref{eqn.eps_1_1_1_1_polybox}, with each $p\tripow1$-pretree in place of the $p$-position that labels its root and height-$1$ vertices in place of their corresponding directions, which are just length-$1$ rooted paths.

Finally, all this could be generalized to other values of $n$.
Here are different ways to draw polyboxes for $\epsilon^{(n)}_p\colon\car{t}_p\to p\tripow{n}$, viewing its codomain as $p\tripow\ell\tri p\tripow{m}$ with $n=\ell+m$ or $p\tripow{i}\tri p\tripow{j}\tri p\tripow{k}$ with $n=i+j+k$:
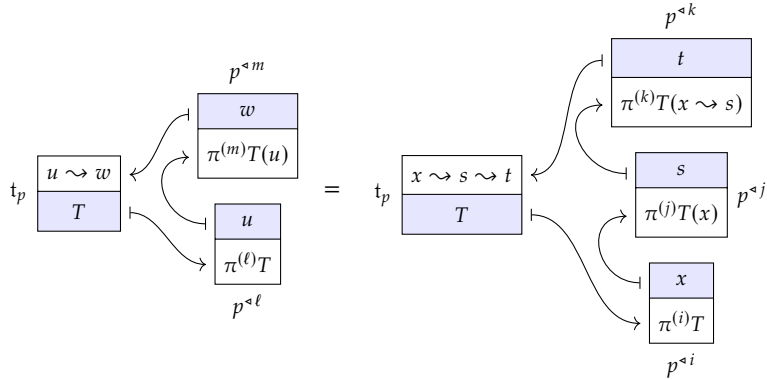
\[
\begin{tikzpicture}
\node (1) {
\begin{tikzpicture}[polybox, mapstos]
    \node[poly, dom, "$\car{t}_p$" left] (r) {$u\leadsto w$\at$T$};
    \node[poly, cod, right=1.8 of r.south, yshift=-1ex, "$p\tripow\ell$" below] (q) {$u$\at$\pi^{(\ell)}T$};
    \node[poly, cod, above=of q, "$p\tripow{m}$" above] (q') {$w$\at$\pi^{(m)}T(u)$};

    \draw (r_pos) to[first] node[below] {} (q_pos);
    \draw (q_dir) to[climb] node[right] {} (q'_pos);
    \draw (q'_dir) to[last] node[above] {} (r_dir);
\end{tikzpicture}
};
\node[right=.65 of 1] (2) {
\begin{tikzpicture}[polybox, mapstos]
	\node[poly, dom, "$\car{t}_p$" left] (S) {$x\leadsto s\leadsto t$\at$T$};
	\node[poly, cod, right=of S, "$p\tripow{j}$" right] (S2) {$s$\at$\pi^{(j)}T(x)$};
	\node[poly, cod, below=of S2, "$p\tripow{i}$" below] (S1) {$x$\at$\pi^{(i)}T$};
	\node[poly, cod, above=of S2, "$p\tripow{k}$" above] (S3) {$t$\at$\pi^{(k)}T(x\leadsto s)$};
	\draw (S_pos) to[first] (S1_pos);
	\draw (S1_dir) to[climb] (S2_pos);
	\draw (S2_dir) to[climb] (S3_pos);
	\draw (S3_dir) to[last] (S_dir);
\end{tikzpicture}
};
\node at ($(1.east)!.5!(2.west)$) {=};
\end{tikzpicture}
\]
\end{example}

\index{polybox}
\index{cofree comonoid!carrier of|)}

%---- Subsection ----%
\subsection{Cofree comonoids as categories} \label{subsec.comon.cofree.cons.cat}
\index{category!cofree comonoids as}\index{cofree comonoid!as category}

We have now characterized the carrier $\car{t}_p$ of the comonoid $\cofree{p}$ that will turn out to be cofree on $p$; but we have yet to describe the comonoid structure of $\cofree{p}$ that $\car{t}_p$ carries.
This structure will allow us to interpret $\cofree{p}$ as a category, whose objects will be $p$-trees and whose morphisms will be the vertices of each $p$-tree.
We could go ahead and describe this category right now, but as category theorists, let us show that the eraser and duplicator of $\cofree{p}$, as well as the comonoid laws that they satisfy, arise naturally from our construction of $\car{t}_p$ as the limit of \eqref{eqn.cofree_diagram}.

\subsubsection{The eraser of the cofree comonoid}
\index{cofree comonoid!eraser of}

The eraser for a comonoid $\cofree{p}$ carried by $\car{t}_p$ should be a lens of the form $\car{t}_p\to\yon$.
Conveniently, since $\yon$ appears in \eqref{eqn.cofree_diagram}, its limit $\car{t}_p$ is already equipped with a canonical lens $\epsilon_p\colon\car{t}_p\to\yon$, as seen in \eqref{eqn.cofree_with_lim}.
This lens will turn out to be the eraser of the cofree comonoid on $p$.

The eraser picks out a direction at each position of the carrier to be the identity morphism on that object.
As each $\car{t}_p$-position is a tree $T\in\tr_p$, and each $\car{t}_p[T]$-direction is a vertex of $T$, the eraser $\epsilon_p$ should pick out a single vertex of every $p$-tree.
Even without looking at \eqref{eqn.cofree_with_lim}, you may already have your suspicions as to which vertex this will turn out to be, as there is, after all, only one sensible way to choose a canonical vertex that every $p$-tree is guaranteed to have: choose its root.
Indeed, \cref{prop.cofree_carrier_lim} tells us that $\epsilon_p\colon\car{t}_p\to\yon$ sends each $p$-tree $T\in\tr_p$ to its stage-$0$ pretree on positions, stripping away everything except for the single unlabeled root of $T$; then on-directions, $\epsilon_p$ picks out the unique height-$0$ vertex of $T$ in $\vtx_0(T)\ss\vtx(T)$, which is just its root.

So in the category $\cofree{p}$, where morphisms out of objects are vertices of trees, the identity morphism on a tree is its root---in fact, we will henceforth denote the root of a $p$-tree by $\id_T$ (so that $\vtx_0(T)=\{\id_T\}$).
Equivalently, we can identify $\id_T$ with the unique length-$0$ rooted path of $T$: the \emph{empty rooted path}, which starts and ends at the root.
The rest of the categorical structure of $\cofree{p}$ will be determined by its duplicator.

\subsubsection{The duplicator of the cofree comonoid}
\index{cofree comonoid!duplicator of}

The duplicator for a comonoid $\cofree{p}$ carried by $\car{t}_p$ should be a lens $\delta_p\colon\car{t}_p\to\car{t}_p\tri\car{t}_p$.
But before we can specify such a lens, we need to figure out what kind of polynomial $\car{t}_p\tri\car{t}_p$ is.

Here is where our work from \cref{subsec.comon.comp.prop.lim_left,subsec.comon.comp.prop.lim_right} comes in handy: since the diagram in \eqref{eqn.cofree_diagram} is connected, its limit $\car{t}_p$ is a connected limit, and we showed in \cref{thm.connected_limits} that $\tri$ preserves connected limits.
Therefore $\car{t}_p\tri\car{t}_p$ is itself the limit of some diagram, and a morphism to $\car{t}_p\tri\car{t}_p$ is just a cone over that diagram.

Which diagram is it?
First, we can use the fact that $\tri$ preserves connected limits on the left to expand $\car{t}_p\tri\car{t}_p$ as the limit of the following diagram, where we have applied the functor $-\tri\car{t}_p$ to the diagram from \eqref{eqn.cofree_diagram}:
\[
\begin{tikzcd}
	\yon\tri\car{t}_p \ar[d] &
	p\tri\car{t}_p \ar[d] &
	p\tripow2\tri\car{t}_p \ar[d] &
	p\tripow3\tri\car{t}_p \ar[d] &
	\cdots\\
	\1\tri\car{t}_p &
	p\tri\1\tri\car{t}_p \ar[l] &
	p\tripow2\tri\1\tri\car{t}_p \ar[l] &
	p\tripow3\tri\1\tri\car{t}_p \ar[l] &
	\cdots.\ar[l]
\end{tikzcd}
\]
But we have $p\tripow{n}\tri\1\tri\car{t}_p\iso p\tripow{n}\tri\1$, so we can simplify this diagram like so:
\begin{equation} \label{eqn.cofree_tri2_lim_expand_left}
\begin{tikzcd}
	\yon\tri\car{t}_p \ar[d] &
	p\tri\car{t}_p \ar[d] &
	p\tripow2\tri\car{t}_p \ar[d] &
	[10pt] p\tripow3\tri\car{t}_p \ar[d] &
	\cdots\\
	\1 &
	p\tri\1 \ar[l] &
	p\tripow2\tri\1 \ar[l] &
	p\tripow3\tri\1 \ar[l] &
	\cdots.\ar[l]
\end{tikzcd}
\end{equation}
(This works because according to \cref{prop.left_pres_lim}, $\tri$ actually preserves \emph{all} limits on the left---including the terminal object $\1$.)
So $\car{t}_p\tri\car{t}_p$ is the limit of \eqref{eqn.cofree_tri2_lim_expand_left}.\index{limit!preserved by $\tri$}\index{limit!connected}

Then we use the fact that $\tri$ preserves connected limits on the right to expand each $p\tripow\ell\tri\car{t}_p$ in \eqref{eqn.cofree_tri2_lim_expand_left} as the limit of the following:
\begin{equation} \label{eqn.cofree_tri2_lim_expand_right}
\begin{tikzcd}[column sep=small]
	p\tripow\ell\tri\yon \ar[d] &
	p\tripow\ell\tri p \ar[d] &
	p\tripow\ell\tri p\tripow2 \ar[d] & p\tripow\ell\tri p\tripow3 \ar[d] &
	\cdots\\
	p\tripow\ell\tri\1 &
	p\tripow\ell\tri p\tri\1 \ar[l] &
	p\tripow\ell\tri p\tripow2\tri\1 \ar[l] &
	p\tripow\ell\tri p\tripow3\tri\1 \ar[l] &
	\cdots.\ar[l]
\end{tikzcd}
\end{equation}
So the limit of \eqref{eqn.cofree_tri2_lim_expand_left} is the limit of a larger diagram, where we have ``plugged in'' a copy of \eqref{eqn.cofree_tri2_lim_expand_right} in place of each $p\tripow\ell\tri\car{t}_p$ that appears.

\index{functor!limit-preserving}

It's worth being a little careful, though, when we draw this diagram: each $p\tripow\ell\tri\car{t}_p$ in \eqref{eqn.cofree_tri2_lim_expand_left} appears with a lens to $p\tripow\ell\tri\1$, and we need to work out what arrow should go in its place.
The lens in question is given by applying $p\tripow\ell\tri-$ to the unique lens $\car{t}_p\to\1$.
But $\1$ actually shows up in the lower left corner of \eqref{eqn.cofree_diagram}, so the unique lens $\car{t}_p\to\1$ is just the projection from the limit of \eqref{eqn.cofree_diagram} to that $\1$.
Then once we apply the connected limit-preserving functor $p\tripow\ell\tri-$ to \eqref{eqn.cofree_diagram}, yielding \eqref{eqn.cofree_tri2_lim_expand_right}, we obtain the lens $p\tripow\ell\tri\car{t}_p\to p\tripow\ell\tri\1$ we desire as the projection from the limit of \eqref{eqn.cofree_tri2_lim_expand_right} to the $p\tripow\ell\tri\1$ in the lower left corner.
So if we want to replace each $p\tripow\ell\tri\car{t}_p$ in \eqref{eqn.cofree_tri2_lim_expand_left} with a copy of \eqref{eqn.cofree_tri2_lim_expand_right}, without changing the limit of the whole diagram, we should replace the arrow $p\tripow\ell\tri\car{t}_p\to p\tripow\ell\tri\1$ with an identity arrow from the $p\tripow\ell\tri\1$ in the lower left corner of \eqref{eqn.cofree_tri2_lim_expand_right} to the $p\tripow\ell\tri\1$ in \eqref{eqn.cofree_tri2_lim_expand_left}.
Of course, we can then collapse these identity arrows down without changing the limit, so the diagram we are left with should look like this:
\begin{equation} \label{eqn.cofree_tri2_diagram}
\scalebox{.9}{
\begin{tikzcd}[row sep=small, column sep=tiny, ampersand replacement=\&]
    p\tripow3 \ar[rd] \&\& p\tripow3\tri p \ar[rd] \&\& p\tripow3\tri p\tripow2 \ar[rd] \&\& \cdots \\
    \& p\tripow3\tri\1 \ar[dd] \&\& p\tripow3\tri p\tri1 \ar[dd] \&\& p\tripow3\tri p\tripow2\tri\1 \ar[dd] \&\& \cdots \\
    p\tripow2 \ar[rd] \&\& p\tripow2\tri p \ar[rd] \&\& p\tripow2\tri p\tripow2 \ar[rd] \&\& \cdots \\
    \& p\tripow2\tri\1 \ar[dd] \&\& p\tripow2\tri p\tri1 \ar[dd] \&\& p\tripow2\tri p\tripow2\tri\1 \ar[dd] \&\& \cdots \\
    p \ar[rd] \&\& p\tri p \ar[rd] \&\& p\tri p\tripow2 \ar[rd] \&\& \cdots \\
    \& p\tri\1 \ar[dd] \&\& p\tri p\tri1 \ar[dd] \&\& p\tri p\tripow2\tri\1 \ar[dd] \&\& \cdots \\
    \yon \ar[rd] \&\& p \ar[rd] \&\& p\tripow2 \ar[rd] \&\& \cdots \\
    \& \1 \&\& p\tri\1 \ar[ll] \&\& p\tripow2\tri\1 \ar[ll] \&\& \cdots.\ar[ll]
\end{tikzcd}
}
\end{equation}
Here the bottom row of \eqref{eqn.cofree_tri2_lim_expand_left} is still the bottom row of \eqref{eqn.cofree_tri2_diagram}, but in the place of each $p\tripow\ell\tri\car{t}_p$ in the top row, we have grafted in a copy of \eqref{eqn.cofree_tri2_lim_expand_right}.
Yet the limit of the diagram is preserved: the limit of \eqref{eqn.cofree_tri2_diagram} is still $\car{t}_p\tri\car{t}_p$.
In summary, for purely formal reasons, \cref{prop.left_pres_lim,thm.connected_limits} ensure that since the limit of \eqref{eqn.cofree_diagram} is $\car{t}_p$, the limit of \eqref{eqn.cofree_tri2_diagram} is $\car{t}_p\tri\car{t}_p$.

While \eqref{eqn.cofree_tri2_diagram} has become rather unwieldy to draw, it is easy to characterize: it is a diagram with a copy each of
\[
    p\tripow\ell\tri p\tripow{m} \qqand p\tripow\ell\tri p\tripow{m}\tri\1
\]
for every $\ell,m\in\nn$, with arrows
\begin{equation} \label{eqn.cofree_tri2_diagram_arr_diag}
    p\tripow\ell\tri p\tripow{m}\to p\tripow\ell\tri p\tripow{m}\tri\1
\end{equation}
between them; along with arrows
\begin{equation} \label{eqn.cofree_tri2_diagram_arr_vert}
    p\tripow\ell\tri p\tripow{(m+1)}\tri\1\to p\tripow\ell\tri p\tripow{m}\tri\1,
\end{equation}
drawn vertically in \eqref{eqn.cofree_tri2_diagram}; and
\begin{equation} \label{eqn.cofree_tri2_diagram_arr_horiz}
    p\tripow{\ell+1}\tri p\tripow0\tri\1\iso p\tripow{\ell+1}\tri\1\to p\tripow\ell\tri\1\iso p\tripow\ell\tri p\tripow0\tri\1,
\end{equation}
drawn horizontally in \eqref{eqn.cofree_tri2_diagram} along the bottom row.
Of course, we could write each $p\tripow\ell\tri p\tripow{m}$ as $p\tripow{(\ell+m)}$, but the notation $p\tripow\ell\tri p\tripow{m}$ helps us distinguish it as the object appearing $\ell$ rows above the bottom and $m$ columns to the right in \eqref{eqn.cofree_tri2_diagram}, in contrast with any other $p\tripow{\ell'}\tri p\tripow{m'}$ with $(\ell,m)\neq(\ell',m')$, even if $\ell+m=\ell'+m'$ guarantees that there are isomorphisms between these objects.

Nevertheless, it is these isomorphisms that induce a canonical lens $\car{t}_p\to\car{t}_p\tri\car{t}_p$, by giving a map from the diagram \eqref{eqn.cofree_diagram} to the diagram \eqref{eqn.cofree_tri2_diagram} as follows.
For all $\ell,m,n\in\nn$ with $n=\ell+m$, we have a canonical isomorphism $p\tripow{n}\To\iso p\tripow\ell\tri p\tripow{m}$ sending the $p\tripow{n}$ that appears in \eqref{eqn.cofree_diagram} to the $p\tripow\ell\tri p\tripow{m}$ that appears in \eqref{eqn.cofree_tri2_diagram}, and similarly sending $p\tripow{n}\tri\1$ in \eqref{eqn.cofree_diagram} to $p\tripow\ell\tri p\tripow{m}\tri\1$ in \eqref{eqn.cofree_tri2_diagram}.
Then all the arrows that appear in \eqref{eqn.cofree_tri2_diagram}---i.e.\ all the arrows in \eqref{eqn.cofree_tri2_diagram_arr_diag}, \eqref{eqn.cofree_tri2_diagram_arr_vert}, and \eqref{eqn.cofree_tri2_diagram_arr_horiz}---can be identified with arrows that appear in \eqref{eqn.cofree_diagram}, so everything commutes.
We therefore induce a lens from the limit $\car{t}_p$ of \eqref{eqn.cofree_diagram} to the limit $\car{t}_p\tri\car{t}_p$ of \eqref{eqn.cofree_tri2_diagram}, which we call $\delta_p\colon\car{t}_p\to\car{t}_p\tri\car{t}_p$.
This turns out to be the duplicator of $\cofree{p}$.

How does $\delta_p$ behave in terms of $p$-trees?
First, we characterize $\car{t}_p\tri\car{t}_p$ and its projections to each $p\tripow\ell\tri p\tripow{m}$ in \eqref{eqn.cofree_tri2_diagram}.
Concretely, we know that a $\car{t}_p\tri\car{t}_p$-position is just a $p$-tree $T\in\tr_p$ along with a $p$-tree $U(v)\in\tr_p$ associated with every vertex $v\in\vtx(T)$.
Then a direction at that position is a choice of vertex $v\in\vtx(T)$ and another vertex $w\in\vtx(U(v))$ of the $p$-tree $U(v)$ associated with $v$.
Each projection from $\car{t}_p\tri\car{t}_p$ to $p\tripow\ell\tri p\tripow{m}$ in \eqref{eqn.cofree_tri2_diagram} is then the composition product of the projections $\epsilon^{(\ell)}_p\colon\car{t}_p\to p\tripow\ell$ and $\epsilon^{(m)}_p\colon\car{t}_p\to p\tripow{m}$, which by \cref{ex.eps_n_polybox} we can draw in polyboxes like so:
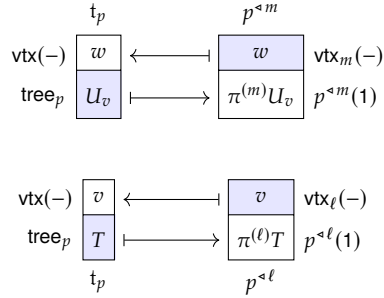
\[
\begin{tikzpicture}[polybox, mapstos]
	\node[poly, dom, "$\car{t}_p$" below] (p) {$v$\at$\vphantom{\pi^{(\ell)}}T$};
        \node[left=0pt of p_pos] {$\tr_p$};
        \node[left=0pt of p_dir] {$\vtx(-)$};
	\node[poly, dom, above=.8 of p, "$\car{t}_p$" above] (p') {$w$\at$\vphantom{\pi^{(m)}}U_v$};
        \node[left=0pt of p'_pos] {$\tr_p$};
        \node[left=0pt of p'_dir] {$\vtx(-)$};

	\node[poly, cod, right=of p, "$p\tripow\ell$" below] (q) {$v$\at$\pi^{(\ell)}T$};
	    \node[right=0pt of q_pos] {$p\tripow\ell(\1)$};
        \node[right=0pt of q_dir] {$\vtx_\ell(-)$};
	\node[poly, cod, above=.8 of q, "$p\tripow{m}$" above] (q') {$w$\at$\pi^{(m)}U_v$};
	    \node[right=0pt of q'_pos] {$p\tripow{m}(\1)$};
        \node[right=0pt of q'_dir] {$\vtx_m(-)$};

	\draw (p_pos) -- node[below] {} (q_pos);
	\draw (q_dir) -- node[above] {} (p_dir);
	\draw (p'_pos) -- node[below] {} (q'_pos);
	\draw (q'_dir) -- node[above] {} (p'_dir);
\end{tikzpicture}
\]
On positions, $\epsilon^{(\ell)}_p\tri\epsilon^{(m)}_p$ sends the $p$-tree $T$ to its stage-$\ell$ pretree $\pi^{(\ell)}T$, and it sends each $p$-tree $U_v$ associated to a height-$\ell$ vertex $v\in\vtx_\ell(T)$ to its stage-$m$ pretree $\pi^{(m)}U_v$, to be $p\tripow{m}$-pretree associated with the height-$\ell$ leaf $v$ of $\pi^{(m)}T$.
Equivalently, this specifies a $p\tripow{(\ell+m)}$-pretree on the right: its bottom $\ell$ levels coincide with the bottom $\ell$ levels of $T$, and its top $m$ levels coincide with the bottom $m$ levels of the $U_v$'s for $v\in\vtx_\ell(T)$.
Then on directions, $\epsilon^{(\ell)}_p\tri\epsilon^{(m)}_p$ is the canonical inclusion of vertices $\vtx_\ell(T)\inj\vtx(T)$ sending $v\mapsto v$, followed by the canonical inclusion of vertices $\vtx_m(U_v)\inj\vtx(U_v)$.
These lenses comprise the universal cone over \eqref{eqn.cofree_tri2_diagram}.

Meanwhile, the other cone we formed over \eqref{eqn.cofree_tri2_diagram} is comprised of lenses of the form $\epsilon^{(\ell+m)}_p\colon\car{t}_p\to p\tripow{(\ell+m)}\iso p\tripow\ell\tri p\tripow{m}$, each of which should factor through $\car{t}_p\tri\car{t}_p$, as depicted in the following commutative diagram:
\begin{equation} \label{eqn.delta_p_tri_epsilons}
\begin{tikzcd}[row sep=large]
    \car{t}_p \ar[r,dashed,"\delta_p"] \ar[dr,"\epsilon^{(\ell+m)}_p"'] & \car{t}_p\tri\car{t}_p \ar[d,"\epsilon^{(\ell)}_p\:\tri\:\epsilon^{(m)}_p"] \\
    & p\tripow{\ell}\tri p\tripow{m}.
\end{tikzcd}
\end{equation}
Indeed, by the universal property of $\car{t}_p\tri\car{t}_p$, our $\delta_p$ is the unique lens for which the above diagram commutes for all $\ell,m\in\nn$.
We will use the equation given by the commutativity of \eqref{eqn.delta_p_tri_epsilons} repeatedly in what follows, whenever we work with $\delta_p$.

Expressing \eqref{eqn.delta_p_tri_epsilons} as an equation of polyboxes, using our usual labels for the arrows of the duplicator $\delta_p$ on the left and our depiction of the projection $\epsilon^{(\ell+m)}_p\colon\car{t}_p\to p\tripow\ell\tri p\tripow{m}$ from \cref{ex.eps_n_polybox} on the right, we have
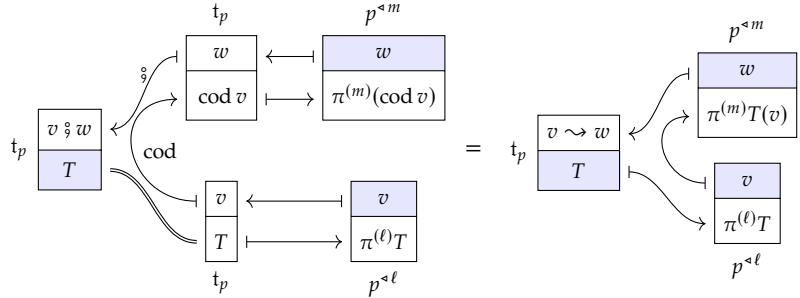
\[
\begin{tikzpicture}
	\node (1) {
  \begin{tikzpicture}[polybox, mapstos]
	\node[poly, dom, "$\car{t}_p$" left] (r) {$v\then w$\at$T$};
	\node[poly, right=1.8 of r.south, yshift=-2.5ex, "$\car{t}_p$" below] (p) {$v$\at$\vphantom{\pi^{(\ell)}}T$};
	\node[poly, above=.8 of p, "$\car{t}_p$" above] (p') {$w$\at$\vphantom{\pi^{(m)}()}\cod v$};
	\node[poly, cod, right=of p, "$p\tripow\ell$" below] (q) {$v$\at$\pi^{(\ell)}T$};
	\node[poly, cod, above=.8 of q, "$p\tripow{m}$" above] (q') {$w$\at$\pi^{(m)}(\cod v)$};

	\draw (p_pos) -- node[below] {} (q_pos);
	\draw (q_dir) -- node[above] {} (p_dir);
	\draw (p'_pos) -- node[below] {} (q'_pos);
	\draw (q'_dir) -- node[above] {} (p'_dir);
	\draw[double, -] (r_pos) to[first] node[below] {} (p_pos);
	\draw (p_dir) to[climb] node[right] {$\cod$} (p'_pos);
	\draw (p'_dir) to[last] node[above] {$\then$} (r_dir);
  \end{tikzpicture}
	};
	\node[right=.5 of 1] (2) {
  \begin{tikzpicture}[polybox, mapstos]
  	\node[poly, dom, "$\car{t}_p$" left] (r) {$v\leadsto w$\at$T$};
  	\node[poly, cod, right=1.8 of r.south, yshift=-1ex, "$p\tripow\ell$" below] (q) {$v$\at$\pi^{(\ell)}T$};
  	\node[poly, cod, above=of q, "$p\tripow{m}$" above] (q') {$w$\at$\pi^{(m)}T(v)$};

  	\draw (r_pos) to[first] node[below] {} (q_pos);
  	\draw (q_dir) to[climb] node[right] {} (q'_pos);
  	\draw (q'_dir) to[last] node[above] {} (r_dir);
  \end{tikzpicture}
	};
	\node at ($(1.east)!.5!(2.west)$) {=};
\end{tikzpicture}
\]
Recall from \cref{ex.eps_n_polybox} that $T(v)$ denotes the $p$-tree equal to the subtree of $T$ rooted at $v\in\vtx(T)$, while $v\leadsto w\in\vtx(T)$ for $w\in\vtx(T(v))$ is the descendent of $v$ that coincides with $w$ when $T(v)$ is identified with the subtree of $T$ rooted at $v$.
Equivalently, we can identify $v\in\vtx(T)$ with the rooted path in $T$ that ends at the vertex $v$ and $w\in\vtx(T(v))$ with the rooted path of $T(v)$ that ends at the vertex $w$, so $v\leadsto w$ becomes the rooted path in $T$ obtained by concatenating $v$ and $w$.
Then for this equality to hold over all $\ell,m\in\nn$, we want $\cod v\coloneqq T(v)$ and $v\then w\coloneqq v\leadsto w$; in fact, $\cod v$ and $v\then w$ will henceforth be our preferred notation for $T(v)$ and $v\leadsto w$.

\subsubsection{Verifying the comonoid laws}

Putting together our constructions of the carrier, the eraser, and the duplicator of the $\cofree{p}$ yields the following result.
\begin{proposition} \label{prop.cofree_as_cat}
As defined above, $\left(\car{t}_p,\epsilon_p,\delta_p\right)$ is a polynomial comonoid corresponding to a category $\cofree{p}$ characterized as follows.
\begin{itemize}
    \item An object in $\cofree{p}$ is a \textit{$p$-tree} in $T\in\tr_p$.
    \item A morphism emanating from $T$ is a \textit{rooted path} in $T$; its codomain is the $p$-subtree rooted at the end of the path.
    \item The identity morphism on $T$ is its empty rooted path.
    \item Composition is given by concatenating rooted paths: given a rooted path $v$ in $T$ and a rooted path $w$ in $\cod v$, the $p$-subtree rooted at the end of $v$, we identify $w$ with the corresponding path in $T$ that starts where $v$ ends, then concatenate the two paths in $T$ to obtain the composite morphism $v\then w$, another rooted path in $T$.
\end{itemize}
\end{proposition}
\begin{proof}
We have already shown how to construct the given carrier, eraser, and duplicator purely diagrammatically, and we have given them concrete interpretations in terms of $p$-trees, their vertices, and their $p$-subtrees.
So it remains to verify that the category laws hold for $\cofree{p}$, or equivalently that the comonoid laws hold for $\left(\car{t}_p,\epsilon_p,\delta_p\right)$.
Again, our argument will be purely formal, although it is not too hard to see that our concrete characterization above in terms of $p$-trees satisfies the laws for a category.

First, we verify the left erasure law: that
\[
\begin{tikzcd}
    \yon\tri\car{t}_p & \car{t}_p \ar[l,equals] \ar[d, "\delta_p"] \\
    & \car{t}_p\tri\car{t}_p \ar[ul, "\epsilon_p\:\tri\:\car{t}_p"]
\end{tikzcd}
\]
commutes.
We know the lens $\epsilon_p\tri\car{t}_p\colon\car{t}_p\tri\car{t}_p\to\yon\tri\car{t}_p$ is characterized by its components $\epsilon_p\tri\epsilon^{(n)}_p\colon\car{t}_p\tri\car{t}_p\to p\tripow0\tri p\tripow{n}$, a projection out of the limit $\car{t}_p\tri\car{t}_p$, for all $n\in\nn$.
Then by our construction of $\delta_p$, each composite $\delta_p\then\left(\epsilon_p\tri\epsilon^{(n)}_p\right)\colon\car{t}_p\to p\tripow0\tri p\tripow{n}$ is the component of $\delta_p$ equal to $\epsilon^{(n)}_p=\epsilon^{(0+n)}_p\colon\car{t}_p\to p\tripow0\tri p\tripow{n}\iso p\tripow{n}$ (this is just the commutativity of \eqref{eqn.delta_p_tri_epsilons} in the case of $(\ell,m)=(0,n)$).
Together, these characterize the composite lens $\delta_p\then\left(\epsilon_p\tri\car{t}_p\right)\colon\car{t}_p\to\yon\tri\car{t}_p\iso\car{t}_p$ as the lens whose components are the projections $\epsilon^{(n)}_p\colon\car{t}_p\to p\tripow{n}$ from the limit.
It follows from the universal property of $\car{t}_p$ that $\delta_p\then\left(\epsilon_p\tri\car{t}_p\right)$ can only be the identity lens on $\car{t}_p$.
Hence the left erasure law holds.

The right erasure law, that
\[
\begin{tikzcd}
    \car{t}_p \ar[r,equals] \ar[d, "\delta_p"'] & \car{t}_p\tri\yon \\
    \car{t}_p\tri\car{t}_p \ar[ur, "\car{t}_p\:\tri\:\epsilon_p"']
\end{tikzcd}
\]
commutes, follows similarly: the composite $\delta_p\then\left(\car{t}_p\tri\epsilon_p\right)$ is characterized by its components over $n\in\nn$ of the form $\delta_p\then\left(\epsilon^{(n)}_p\tri\epsilon_p\right)\colon\car{t}_p\to p\tripow{n}\tri p\tripow0$, which we know by \eqref{eqn.delta_p_tri_epsilons} is equal to $\epsilon^{(n+0)}_p=\epsilon^{(n)}_p$, the projection out of the limit $\car{t}_p$.
It follows from the universal property of this limit that $\delta_p\then\left(\car{t}_p\tri\epsilon_p\right)$ must be the identity lens on $\car{t}_p$.

Finally, we check the coassociative law: that
\[
\begin{tikzcd}[column sep=large]
	\car{t}_p\ar[r, "\delta_p"]\ar[d, "\delta_p"']&
	\car{t}_p\tri\car{t}_p\ar[d, "\car{t}_p\:\tri\:\delta_p"]\\
	\car{t}_p\tri\car{t}_p\ar[r, "\delta_p\:\tri\:\car{t}_p"']&
	\car{t}_p\tri\car{t}_p\tri\car{t}_p,
\end{tikzcd}
\]
commutes.
Because $\tri$ preserves connected limits, we can write $\car{t}_p\tri\car{t}_p\tri\car{t}_p$ as a limit the way we did with $\car{t}_p\tri\car{t}_p$: it is the limit of diagram consisting of arrows
\[
    p\tripow{i}\tri p\tripow{j}\tri p\tripow{k}\to p\tripow{i}\tri p\tripow{j}\tri p\tripow{k}\tri\1
\]
for each $i,j,k\in\nn$, with additional arrows between the position-sets.
So a lens to $\car{t}_p\tri\car{t}_p\tri\car{t}_p$, such as those in the square above, is uniquely determined by its components to $p\tripow{i}\tri p\tripow{j}\tri p\tripow{k}$ for all $i,j,k\in\nn$, obtained by composing it with the projections $\epsilon^{(i)}_p\tri\epsilon^{(j)}_p\tri\epsilon^{(k)}_p$ out of the limit.
Then by \eqref{eqn.delta_p_tri_epsilons},
\begin{align*}
    \delta_p\then\left(\car{t}_p\tri\delta_p\right)\then\left(\epsilon^{(i)}_p\tri\epsilon^{(j)}_p\tri\epsilon^{(k)}_p\right) &=
    \delta_p\then\left(\epsilon^{(i)}_p\tri\left(\delta_p\then\left(\epsilon^{(j)}_p\tri\epsilon^{(k)}_p\right)\right)\right) \\ &=
    \delta_p\then\left(\epsilon^{(i)}_p\tri\epsilon^{(j+k)}_p\right) \tag*{\eqref{eqn.delta_p_tri_epsilons}} \\ &=
    \epsilon^{(i+j+k)}_p \tag*{\eqref{eqn.delta_p_tri_epsilons}} \\ &=
    \delta_p\then\left(\epsilon^{(i+j)}_p\tri\epsilon^{(k)}_p\right) \tag*{\eqref{eqn.delta_p_tri_epsilons}} \\ &=
    \delta_p\then\left(\left(\delta_p\then\left(\epsilon^{(i)}_p\tri\epsilon^{(j)}_p\right)\right)\tri\epsilon^{(k)}_p\right) \tag*{\eqref{eqn.delta_p_tri_epsilons}} \\ &=
    \delta_p\then\left(\delta_p\tri\car{t}_p\right)\then\left(\epsilon^{(i)}_p\tri\epsilon^{(j)}_p\tri\epsilon^{(k)}_p\right).
\end{align*}
for all $i,j,k\in\nn$.
So by the universal property of $\car{t}_p\tri\car{t}_p\tri\car{t}_p$, coassociativity holds.
\end{proof}\index{coassociativity}

We call the category $\cofree{p}$ corresponding to the comonoid $\left(\car{t}_p,\epsilon_p,\delta_p\right)$ the \emph{category of $p$-trees}, the \emph{category of trees on $p$}, or the \emph{$p$-tree category}.

\index{category!of $p$-trees|see{cofree comonoid}}\index{indexed set}

\begin{example}[The category of $p$-trees: states and transitions for (co)free]
Let's step back and think about how $p$-tree categories relate to the original polynomial $p$ from the perspective of the states and transitions of dynamical systems.
Before we build any trees or categories out of it, a polynomial $p$ is just a family of sets of directions indexed over a set of positions.
The directions emerge from positions, but don't point to anywhere in particular.
If we want to interpret the positions of $p$ as states and the directions of $p$ as composable transitions, first we would need to point the directions to other positions by assigning them codomains.
There are many ways to do this, but we would like to do so ``freely,'' without having to make any choices along the way that would bias us one way or another, so as to give a canonical way to interpret $p$ as a category.

So to avoid making choices, rather than fixing a codomain for each direction, we allow every possible direction to point to every possible position once.
But a single direction of any one position cannot point to multiple positions, which is when we need to make several copies of each position: at least one copy for every possible combination of codomains that can be assigned to its directions.
This is how we get from $p$ to $p\tripow2$: the $p\tripow2$-pretrees represent all the ways we could assign codomains to the directions in $p$.
Essentially, we have freely refined our positions into more specific states to account for everywhere their directions could lead.
Each $p\tripow2$-pretree still remembers which $p\tripow1$-pretree it grew from, giving us our canonical trimming operation $p\tripow2(\1)\to p\tripow1(\1)$.

Yet even that is not enough: sure, we've assigned codomains to the directions of $p$ in every possible way, but now that we're building a category, we want our directions to be composable transitions.
So each pair of composable directions of $p$—each length-$2$ rooted path of the $p\tripow2$-pretrees—is now a possible transition as well, another morphism in our category.
To keep everything canonical, we still want to avoid making any actual choices; we can't simply say that two directions of $p$ compose to a third direction of $p$ that we already have.
Each pair of composable directions must be an entirely new morphism---and every new morphism needs a codomain.

You can tell where this is going.
To avoid actually choosing codomains for the new morphisms, we need to refine our $p\tripow2$-pretrees by making copies of them to account for every possibility, building $p\tripow3$-pretrees as a result.
Then their length-$3$ rooted paths are new morphisms, too, and they need new codomains, and so forth, ad infinitum.
Indeed, this process cannot terminate in a finite number of steps, but that's okay---we can take the limit, yielding the ultimate free refinement of positions into states and free composites of directions as transitions.

This is where \eqref{eqn.cofree_diagram} comes from; it captures the infinite process of turning sequences of directions into morphisms by giving them codomains in every possible way, then dealing with the longer sequences of directions that emerge as a result.
Note that we also need to account for the fact that every object needs an identity; to avoid making a choice, we don't set it to be any direction or composite of multiple directions of $p$, reserving it instead for the empty sequence of directions that every pretree has.
Then the limit $\car{t}_p$ of \eqref{eqn.cofree_diagram} is all we need: the $p$-trees are states representing every possible way that sequences of directions emerging from a $p$-position can lead to other $p$-positions, and the rooted paths of each $p$-tree are transitions accounting for every finite sequence of composable directions from the corresponding state.

Since composites were freely generated, directions and thus entire rooted paths compose by concatenation---making the empty rooted path the correct identity.
And the codomain is whatever the direction at the end of a rooted path has been freely assigned to point to---not just the $p$-position there but the whole $p$-subtree, an entire state representing all the sequences of directions one could take and the positions to which they lead, starting from the end of the path just followed.
This gives the comonoid structure on $\car{t}_p$, completing the cofree construction of the category of $p$-trees $\cofree{p}$.
\end{example}

Soon we will prove that $\cofree{p}$ really is the cofree category on $p$, but first we give some examples to make all this concrete.

\subsubsection{Examples of $p$-tree categories}
\index{cofree comonoid!examples of|(}

In what follows, we will freely switch between the morphisms-as-vertices and the morphisms-as-rooted-paths perspectives of $p$-tree categories given in \cref{prop.cofree_as_cat} whenever convenient.

\begin{example}[The category of $\1$-trees]
Let's start by taking $p\coloneqq\1$.
In \cref{exc.p_tree_polys} \cref{exc.p_tree_polys.1}, we showed that there is a unique $\1$-tree: a tree with only $1$ vertex, its root.
So $\car{t}_\1\iso\yon$ is the carrier of the category of $\1$-trees $\cofree\1$.

Up to isomorphism, there is only one category carried by $\yon$: the terminal category with $1$ object and no morphisms aside from the $1$ identity.
When we think of this as the category of $\1$-trees, we can characterize it as follows.
\begin{itemize}
    \item The $1$ object is the tree with $1$ vertex: call the tree $\bullet$, because that's what it looks like.
    \item The $1$ morphism $\bullet\to\_$ is the $1$ vertex of $\bullet$, its root.
    After all, there are no directions of $\1$ to freely compose, so the only morphism generated is the identity.
    Since the $\1$-subtree of $\bullet$ rooted as its root is just the entire $\1$-tree $\bullet$, the codomain of this morphism is still $\bullet$, which makes sense---it's the only object we have.
    \item The identity morphism $\id_\bullet\colon\bullet\to\bullet$ is the root of $\bullet$.
    This is the same morphism we just mentioned.
    Equivalently, it is the empty rooted path.
    \item The composite morphism $\id_\bullet\then\id_\bullet=\id_\bullet\then\id_{\cod(\id_\bullet)}\colon\bullet\to\_$ is obtained by concatenating the empty rooted path $\id_\bullet$ of $\bullet$ with the empty rooted path $\id_{\cod(\id_\bullet)}$ of $\cod(\id_\bullet)$, which is again just the empty rooted path of $\bullet$.
    Hence $\id_\bullet\then\id_{\cod\id_\bullet}=\id_\bullet$, as expected---it's the only morphism we have.
\end{itemize}
This is certainly too trivial an example to say much about, but it demonstrates how we can interpret the category carried by $\yon$ as the category of $p$-trees for $p\coloneqq\1$.
Moreover, it helps us see concretely why taking the identity to be the root gives it the right codomain and compositional behavior.

\end{example}

\begin{exercise}[The category of trees on a constant]
Let $B$ be a set, viewed as a constant polynomial.\index{constant polynomial}
\begin{enumerate}
	\item What is the polynomial $\car{t}_B$?
	\item Characterize the $B$-tree category $\cofree{B}$.
\qedhere
\end{enumerate}
\begin{solution}
\begin{enumerate}
    \item We compute $\car{t}_B$ as follows.
    To select a $B$-tree, we must select a position from $B$.
    Then we are done, because there are no directions.
    So every $B$-tree consists of $1$ vertex labeled with an element of $B$.
    Hence the set of $B$-trees $\car{t}_B(\1)=\tr_B$ is in bijection with $B$ itself, yielding $\car{t}_B\iso B\yon$.
    (We could have also applied \cref{exc.p_tree_polys} \cref{exc.p_tree_polys.monomial} in the case of $A=\0$, yielding $\car{t}_{B\yon^\0}\iso B^{\lst(\0)}\yon^{\lst(\0)}\iso B\yon$, as $\lst(\0)\iso\1+\0^\1+\0^\2+\cdots\iso\1$---the empty sequence is the only sequence with elements in $\0$.)
    \item We characterize the $B$-tree category $\cofree{B}$ as follows.
    We saw in \cref{exc.linear_poly_comon} and \cref{exc.linear_poly_cat} that (up to isomorphism) there is a unique comonoid structure on $B\yon$, corresponding to the discrete category on $B$.
    So $\cofree{B}$ must be isomorphic to the discrete category on $B$.
    In the language of $B$-trees, the objects of $\cofree{B}$ are trees with exactly $1$ vertex, labeled with an element of $B$.
    Each tree has an identity morphism, its root; but there are no vertices aside from the root, so there are no morphisms aside from the identities.
    In particular, the only $B$-subtree of any one $B$-tree is the entire $B$-tree itself, so there are no morphisms from one $B$-tree to a different $B$-tree.
    So the category of $B$-trees can indeed be identified with the discrete category on $B$.
\end{enumerate}
\end{solution}
\end{exercise}

\begin{example}[The category of $\yon$-trees] \label{ex.yon_tree_nn}
Now consider $p\coloneqq\yon$.
In \cref{exc.p_tree_polys} \cref{exc.p_tree_polys.unary}, we showed that there is a unique $\yon$-tree: a single ray extending from the root in which every vertex has exactly $1$ child, so that there is exactly $1$ height-$n$ vertex---and $1$ length-$n$ rooted path---for every $n\in\nn$.
Then $\car{t}_\yon\iso\yon^\nn$ is the carrier of the category of $\yon$-trees $\cofree\yon$.
In fact, we can identify the set of rooted paths of this $\yon$-tree with $\nn$, so that $n\in\nn$ is the $\yon$-tree's unique length-$n$ rooted path.

We know from \cref{ex.monoids} that comonoids with representable carriers correspond to $1$-object categories, which can be identified with monoids.
In particular, a comonoid structure on $\yon^\nn$ corresponds to a monoid structure on $\nn$.
There is more than one monoid structure on $\nn$, though, so which one corresponds to the category $\cofree\yon$?

We can characterize $\cofree\yon$ in terms of $\yon$-trees as follows.
\begin{itemize}
    \item The $1$ object is the unique $\yon$-tree, a single ray: call it $\uparrow$.
    Here is a picture of this ray, to help you visualize its vertices, rooted paths, and subtrees:
    \begin{gather*}
    \vdots\\
    \begin{tikzpicture}[trees]
    	\node {$\bullet$}
    		child {node {$\bullet$}
    		    child {node {$\bullet$}
    		        child {node {$\bullet$}
    		            child {node {$\vphantom\bullet$}}
    		        }
    	        }
    	    };
    \end{tikzpicture}
    \end{gather*}
    \item The morphisms are the rooted paths of $\uparrow$; they comprise the set $\nn$, where each $n\in\nn$ is the unique rooted path of $\uparrow$ of length $n$.
    Each rooted path represents the free composite of $n$ copies of the sole direction of $\yon$.
    Since $\uparrow$ is an infinite ray, the $\yon$-subtree of $\uparrow$ rooted at any of its vertices is still just a copy of $\uparrow$.
    So the codomain of each morphism is still $\uparrow$, which makes sense---it's the only object we have.
    \item The identity morphism $\id_\uparrow\colon{\uparrow}\to{\uparrow}$ is the empty rooted path, which has length $0$; so $\id_\uparrow=0$.
    In the corresponding monoid structure on $\nn$, the element $0\in\nn$ must then be the unit.
    \item The composite morphism $m\then n\colon{\uparrow}\to{\uparrow}$ for $m,n\in\nn$ is obtained by concatenating $m$, i.e.\ the length-$m$ rooted path of $\uparrow$, with $n$, i.e.\ the length-$n$ rooted path of $\uparrow$, translating the latter path so that it begins where the former path ends.
    The result is then the rooted path of $\uparrow$ of length $m+n$; so $m\then n=m+n$.
    In the corresponding monoid structure on $\nn$, the binary operation must then be given by addition.
\end{itemize}
Hence $\cofree{\yon}$ is the monoid $(\nn,0,+)$ viewed as a $1$-object category.
\end{example}

\begin{example}[$B\yon$-trees are $B$-streams]\label{ex.streams_cofree}\index{streams!cofree comonoid and}
Let $B$ be a set, and consider $p\coloneqq B\yon$.
Generalizing \cref{exc.p_tree_polys} \cref{exc.p_tree_polys.2_label_ray} for $B$ instead of $\2$, or applying \cref{exc.p_tree_polys} \cref{exc.p_tree_polys.monomial} for the case of $A\coloneqq\1$, we can deduce that a $B\yon$-tree consists of a single ray for which every vertex is given a label from $B$.
So the vertices of a $B\yon$-tree are in bijection with $\nn$, and $B\yon$-trees are in bijection with functions $\nn\to B$ assigning each vertex a label.
Hence $\car{t}_{B\yon}\iso B^\nn\yon^\nn$ is the carrier of the category of $B\yon$-trees $\cofree{B\yon}$.
As in \cref{ex.yon_tree_nn}, we identify the set of rooted paths of a given $B\yon$-tree with $\nn$, so that $n\in\nn$ is the $B\yon$-tree's unique length-$n$ rooted path.

In fact, we have already seen the category $\cofree{B\yon}$ once before: it is the category of $B$-streams from \cref{ex.streams_category}.\index{streams}

\begin{itemize}
    \item Recall that a $B$-stream is an element of $B^\nn$ interpreted as a countable sequence of elements $b_n\in B$ for $n\in\nn$, written like so:
    \[
        \ol{b}\coloneqq(b_0\to b_1\to b_2\to b_3\to\cdots).
    \]
    This is just a $B\yon$-tree lying on its side!
    Each arrow is a copy of the sole direction at each position of $B\yon$, pointing to one of the positions in $B$ for that direction to lead to next.
    \item Given a $B$-stream $\ol{b}\in B^\nn$ above, recall that a morphism out of $\ol{b}$ is a natural number $n\in\nn$, and its codomain is the substream of $B$ starting at $b_n$:
    \[
        \cod(\ol{b}\To{n}\_)=(b_n\to b_{n+1}\to b_{n+2}\to b_{n+3}\to\cdots).
    \]
    This coincides with how we view morphisms as rooted paths and codomains as subtrees rooted at the end of those paths in the category of $B\yon$-trees.
    \item Whether we view $\ol{b}$ as a $B$-stream or a $B\yon$-tree, its identity morphism $\id_{\ol{b}}\colon\ol{b}\to\ol{b}$ corresponds to $0\in\nn$; from the latter perspective, it is the length-$0$ path from the root to itself.
    \item Whether we view $\ol{b}$ as a $B$-stream or a $B\yon$-tree, composition is given by addition; from the latter perspective, concatenating a length-$m$ rooted path with a length-$n$ path yields a length-$(m+n)$ rooted path.
\end{itemize}
\end{example}

\begin{example}[The category of $\nn$-labeled binary trees] \label{ex.nn_labeled_binary_trees}
Consider $p\coloneqq\nn\yon^\2$.
By \cref{exc.p_tree_polys} \cref{exc.p_tree_polys.monomial}, or \cref{exc.n_ary_trees} \cref{exc.n_ary_trees.l_labeled} in the case of $L\coloneqq\nn$ and $\ord{n}\coloneqq\2$, an $\nn\yon^\2$-tree is an infinite binary tree with vertices labeled by elements of $\nn$.
Here's how such an $\nn\yon^\2$-tree might start:
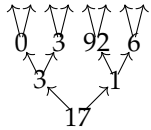
\[
\begin{tikzpicture}[trees,
  level 1/.style={sibling distance=10mm},
  level 2/.style={sibling distance=5mm},
  level 3/.style={sibling distance=2.5mm}]
	\node {$17$}
		child {node {$3$}
			child {node {$0$}
				child
				child
			}
			child {node {$3$}
				child
				child
			}
		}
		child {node {$1$}
			child {node {$92$}
				child
				child
			}
			child {node {$6$}
				child
				child
			}
		};
\end{tikzpicture}
\]
We can characterize the category $\cofree{\nn\yon^\2}$ as follows.
\begin{itemize}
    \item The objects are $\nn$-labeled binary trees from the set $\nn^{\lst(\2)}$.
    \item A morphism out of an $\nn$-labeled binary tree is a binary sequence: a finite list of directions in $\2$, whose elements could be interpreted as ``left'' and ``right,'' thus uniquely specifying a rooted path in a binary tree.
    They comprise the set $\lst(\2)$.
    The codomain of each rooted path is the $\nn$-labeled binary subtree rooted at the end of the path.
    \item The identity morphism on a given $\nn$-labeled binary tree is its empty rooted path.
    \item The composite of two binary sequences is obtained by concatenation.
\end{itemize}
\end{example}

\begin{exercise}[The category of $B$-labeled $A$-ary trees] \label{exc.b-lab_a-ary_trees}
Characterize the $B\yon^A$-tree category $\cofree{B\yon^A}$.
\begin{solution}
We characterize $\cofree{B\yon^A}$, following \cref{ex.nn_labeled_binary_trees}.
\begin{itemize}
    \item The objects are $B\yon^A$-trees.
    By \cref{exc.p_tree_polys} \cref{exc.p_tree_polys.monomial}, or \cref{exc.n_ary_trees} \cref{exc.n_ary_trees.l_labeled} with $L\coloneqq B$ and $\ord{n}$ replaced with an arbitrary set $A$, a $B\yon^A$-tree is a (infinite, assuming $|A|,|B|>0$) tree whose every vertex is labeled by an element of $B$ and whose children are in bijection with $A$.
    We can call these $B$-labeled $A$-ary trees; they are in bijection with the set $B^{\lst(A)}$.
    \item A morphism out of a $B$-labeled $A$-ary tree is a finite list of directions in $A$, or a rooted path in an $A$-ary tree.
    They comprise the set $\lst(A)$.
    The codomain of each rooted path is the $B$-labeled $A$-ary subtree rooted at the end of the path.
    \item The identity morphism on a given $B$-labeled $A$-ary tree is its empty rooted path.
    \item The composite of two lists in $A$ is obtained by concatenation.
\end{itemize}
\end{solution}
\end{exercise}

\begin{exercise}
Characterize the $\left(\yon+\1\right)$-tree category $\cofree{\yon+\1}$.
\begin{solution}
We characterize $\cofree{\yon+\1}$ as follows.
\begin{itemize}
    \item The objects are $(\yon+\1)$-trees.
    By \cref{exc.p_tree_polys} \cref{exc.p_tree_polys.maybe_unary}, a $(\yon+\1)$-tree is either a single length-$n$ path for some $n\in\nn$, which we denote by $[n]$, or a ray, which we denote by $[\infty]$.
    So the set of objects is $\tr_{\yon+\1}=\{[n]\mid n\in\nn\cup\{\infty\}\}\iso\nn\cup\{\infty\}$.
    \item A finite $(\yon+\1)$-tree $[n]$ for $n\in\nn$ has $n+1$ rooted paths: $1$ each of length $0$ through $n$, inclusive.
    So a morphism out of $[n]$ can be identified with an element of the set $\{0,\ldots,n\}$.

    Meanwhile, the infinite $(\yon+\1)$-tree $[\infty]$, a ray, has exactly $1$ rooted path of every length $\ell\in\nn$.
    So a morphism out of $[\infty]$ can be identified with an element of the set $\nn$.
    Every $(\yon+\1)$-subtree of the ray is still a ray, so the codomain of each of these morphisms is still $[\infty]$.
    \item The identity morphism on a given $B$-labeled $A$-ary tree is its empty rooted path.
    \item The composite of two lists in $A$ is obtained by concatenation.
\end{itemize}
\end{solution}
\end{exercise}

\begin{exercise}
Let $p\coloneqq \{\const{a},\const{b},\const{c},\ldots,\const{z},\,\text{\textvisiblespace}\}\yon+\{\bullet\}$.
\begin{enumerate}
	\item Describe the objects of the cofree category $\cofree{p}$, and draw one.
	\item For a given such object, describe the set of emanating morphisms.
	\item Describe how to take the codomain of a morphism.
\qedhere
\end{enumerate}
\begin{solution}
\begin{enumerate}
    \item An object in $\cofree{p}$ is a stream of letters $\const{a},\const{b},\ldots,\const{z}$ and spaces $\text{\textvisiblespace}$, that may go on forever or may end with a large period, $\bullet$. So an example is the infinite stream $\const{aaaaa}\cdots$. Another example is $\const{hello}\text{\textvisiblespace}\const{world}\bullet$.
    \item The object $\const{aaaaa}\cdots$ has $\nn$ as its set of emanating morphisms. The object $\const{hello}\text{\textvisiblespace}\const{world}\bullet$ has $\{0,\ldots,11\}$ as its set of emanating morphisms.
    \item The codomain of any morphism out of $\const{aaaaa}\cdots$ is again $\const{aaaaa}\cdots$. The codomain of any morphism $0\leq i\leq 11$ out of $\const{hello}\text{\textvisiblespace}\const{world}\bullet$ is the string one obtains by removing the first $i$ letters of $\const{hello}\text{\textvisiblespace}\const{world}\bullet$.
\end{enumerate}
\end{solution}
\end{exercise}

\begin{exercise} %%TODO
Let $p\coloneqq\{\bul[my-red],\bul[my-blue]\}\yon^\2+\{\bul[black]\}\yon+\{\bul[my-yellow]\}$ as in \cref{ex.imagining_trees}.
\begin{enumerate}
	\item Choose an object $t\in \tr_p$, i.e.\ a tree in $p$, and draw a finite approximation of it (say four layers).
	\item What is the identity morphism at $t$?
	\item Choose a nonidentity morphism $f$ emanating from $t$ and draw it.
	\item What is the codomain of $f$? Draw a finite approximation of it.
	\item Choose a morphism emanating from the codomain of $f$ and draw it.
	\item What is the composite of your two morphisms? Draw it on $t$.
\qedhere
\end{enumerate}
\begin{solution}
\begin{enumerate}
    \item
\[
\begin{tikzpicture}[trees,
  level 1/.style={sibling distance=10mm},
  level 2/.style={sibling distance=5mm},
  level 3/.style={sibling distance=2.5mm}]
	\node {$\bul[my-red]$}
		child {node {$\bul[my-blue]$}
			child {node {$\bul[black]$}
				child {node {$\bul[my-red]$}
					child
					child
				}
			}
			child {node {$\bul[my-yellow]$}
			}
		}
		child {node {$\bul[black]$}
			child {node {$\bul[my-blue]$}
				child {node {$\bul[my-yellow]$}
				}
				child {node {$\bul[my-blue]$}
					child
					child
				}
			}
		};
\end{tikzpicture}
\]
    \item The identity morphism at $t$ is the trivial path starting at the root node.
    \item We indicate a morphism $f$ using thick arrows:
    \[
\begin{tikzpicture}[trees,
  level 1/.style={sibling distance=10mm},
  level 2/.style={sibling distance=5mm},
  level 3/.style={sibling distance=2.5mm}]
	\node {$\bul[my-red]$}
		child[thin] {node {$\bul[my-blue]$}
			child {node {$\bul[black]$}
				child {node {$\bul[my-red]$}
					child
					child
				}
			}
			child {node {$\bul[my-yellow]$}
			}
		}
		child[very thick] {node {$\bul[black]$}
			child {node {$\bul[my-blue]$}
				child[thin] {node {$\bul[my-yellow]$}
				}
				child[thin] {node {$\bul[my-blue]$}
					child
					child
				}
			}
		};
\end{tikzpicture}
\]
    \item The codomain of $f$ is the tree rooted at the target of the path, namely
   \[
\begin{tikzpicture}[trees,sibling distance=2.5mm]
	\node {$\bul[my-blue]$}
				child[thin] {node {$\bul[my-yellow]$}
				}
				child[thin] {node {$\bul[my-blue]$}
					child
					child
				};
\end{tikzpicture}
\]
    \item Here's a morphism emanating from the codomain of $f$:
   \[
\begin{tikzpicture}[trees,sibling distance=2.5mm]
	\node {$\bul[my-blue]$}
				child[very thick] {node {$\bul[my-yellow]$}
				}
				child[thin] {node {$\bul[my-blue]$}
					child
					child
				};
\end{tikzpicture}
\]
    \item The composite morphism on the original tree is:
    \[
\begin{tikzpicture}[trees,
  level 1/.style={sibling distance=10mm},
  level 2/.style={sibling distance=5mm},
  level 3/.style={sibling distance=2.5mm}]
	\node {$\bul[my-red]$}
		child[thin] {node {$\bul[my-blue]$}
			child {node {$\bul[black]$}
				child {node {$\bul[my-red]$}
					child
					child
				}
			}
			child {node {$\bul[my-yellow]$}
			}
		}
		child[very thick] {node {$\bul[black]$}
			child {node {$\bul[my-blue]$}
				child[very thick] {node {$\bul[my-yellow]$}
				}
				child[thin] {node {$\bul[my-blue]$}
					child
					child
				}
			}
		};
\end{tikzpicture}
\]
\end{enumerate}
\end{solution}
\end{exercise}

\begin{exercise}
Let $p$ be a polynomial, let $Q\coloneqq\{q\in\qq\mid q\geq 0\}$ and consider the monoid $\yon^Q$ of nonnegative rational numbers under addition. Is it true that any retrofunctor $\varphi\colon\cofree{p}\cof\yon^\qq$ is constant, i.e.\ that it factors as
\[
\cofree{p}\cof\yon\cof\yon^\qq?
\]
\begin{solution}
Take $p\coloneqq\2\yon$, and consider the object $x\in\cofree{p}(\1)$ given by the stream
\[
x\coloneqq(2\,12\,112\,1112\,11112\,111112\ldots)
\]
(with spaces only for readability); note that every morphism emanating from $x$ has a different codomain. We need to give $\varphi^\sharp_i(q)$ for every $i\in\cofree{p}(\1)$ and $q\geq 0$. Define
\[
	\varphi^\sharp_i(q)\coloneqq
	\begin{cases}
		i&\tn{ if }i\neq x \tn{ or } q=0\\
		x'&\tn{ if }i=x\tn{ and } q>0
	\end{cases}
\]
where $x'\coloneqq(12\,112\,1112\,11112\,111112\ldots)$. There are three retrofunctor conditions to check, namely identity, codomains, and composition. The codomain condition is vacuous since $\yon^Q$ has one object, and the identity condition holds by construction, because we always have $\varphi^\sharp_i(0)=i$. Now take $q_1,q_2\in Q$; we need to check that
\[\varphi^\sharp_{\cod\varphi^\sharp_i(q_1)}(q_2)=^?\varphi^\sharp_i(q_1+q_2)\]
holds. If $i\neq x$ or $q_1=q_2=0$, then it holds because both sides equal $i$. If $i=x$ and either $q_1>0$ or $q_2>0$, it is easy to check that both sides equal $x'$, so again it holds.
\end{solution}
\end{exercise}

\index{cofree comonoid!examples of|)}

%---- Subsection ----%
\subsection{Exhibiting the forgetful-cofree adjunction}

\index{forgetful-cofree adjunction}\index{cofree comonoid!adjunction with $\poly$}
\index{adjunction!between $\poly$ and $\catsharp$}

We are now ready to give a diagrammatic proof of the main result of this section: as promised, the category $\cofree{p}$ we constructed is the cofree comonoid on $p$.

\begin{theorem}[Cofree comonoid] \label{thm.cofree}\index{cofree comonoid}\index{functor!forgetful $\catsharp\to\poly$}
The forgetful functor $U\colon \catsharp\to\poly$ has a right adjoint $\cofree{-}\colon\poly\to\catsharp$, giving rise to an adjunction
\[
    \adj{\catsharp}{U}{\cofree{-}}{\poly},
\]
such that for each $p\in\poly$, the carrier $\car{t}_p\coloneqq U\cofree{p}$ of the category $\cofree{p}$ is given by the limit of the diagram \eqref{eqn.cofree_diagram}, repeated here:
\[
\begin{tikzcd}
	\yon \ar[d] &
	p \ar[d] &
	p\tripow2 \ar[d] &
	[10pt] p\tripow3 \ar[d] &
	\cdots\\
	\1 &
	p\tri\1 \ar[l, "!"'] &
	p\tripow2\tri\1 \ar[l, "p\:\tri\:!"'] &
	p\tripow3\tri\1 \ar[l, "p\tripow2\:\tri\:!"'] &
	\cdots.\ar[l]
\end{tikzcd}
\]
That is, for any category $\cat{C}\in\catsharp$ with carrier $\car{c}\coloneqq U\cat{C}$, there is a natural isomorphism\index{isomorphism!natural}
\[
    \poly(\car{c},p)\iso\catsharp(\cat{C},\cofree{p}).
\]
\end{theorem}\index{functor!forgetful}

\begin{proof}
To show that $\cofree{-}$ is the right adjoint of $U$, it is enough to show that for every lens $\varphi\colon\car{c}\to p$, there exists a unique retrofunctor $F\colon\cat{C}\cof\cofree{p}$ for which
\begin{equation} \label{eqn.cofree_universal_lens_counit}
\begin{tikzcd}
    \car{c} \ar[d,dashed,"F"'] \ar[dr,"\varphi"] \\
    \car{t}_p \ar[r,"\epsilon^{(1)}_p"'] & p
\end{tikzcd}
\end{equation}
commutes.
Here the projection $\epsilon^{(1)}_p\colon U\cofree{p}\iso\car{t}_p\to p$ serves as the counit of the adjunction, and we identify the retrofunctor $F$ with its underlying lens $UF\colon\car{c}\to\car{t}_p$.

First, we construct $F$ from $\varphi$ as follows.
If we let $\epsilon$ and $\delta$ be the eraser and duplicator of $\cat{C}$, the diagram
\begin{equation} \label{eqn.phi_seed}
\begin{tikzcd}[column sep=large]
    \car{c} \ar[d, "\epsilon"'] \ar[r, "\delta"] \ar[dr, "\varphi"'] &
    \car{c}\tri\car{c} \ar[d, "\epsilon\:\tri\:\varphi"] \\
    \yon \ar[d] & p \ar[d] \\
    \1 & p\tri\1 \ar[l]
\end{tikzcd}
\end{equation}
commutes: the pentagon in the lower left commutes trivially, while the triangle in the upper right commutes by the left erasure law of $\cat{C}$, as
\begin{align*}
    \delta\then\left(\epsilon\tri\varphi\right)&=
    \delta\then\left(\epsilon\tri\car{c}\right)\then\varphi \\ &=
    \id_\car{c}\then\varphi \tag{Left erasure law} \\ &=
    \varphi.
\end{align*}
Then by induction, the larger diagram
\begin{equation} \label{eqn.phi_expand}
\begin{tikzcd}[column sep=large]
    \car{c} \ar[d,"\epsilon"'] \ar[r,"\delta"] \ar[dr,"\varphi"'] & \car{c}\tri\car{c} \ar[d,"\epsilon\:\tri\:\varphi"] \ar[r,"\delta\:\tri\:\car{c}"] & \car{c}\tripow3 \ar[d, "\epsilon\:\tri\:\varphi\tripow2"] \ar[r,"\delta\:\tri\:\car{c}\tripow2"] & \car{c}\tripow4 \ar[d,"\epsilon\:\tri\:\varphi\tripow{3}"] \ar[r] & \cdots \\
    \yon \ar[d] & p \ar[d] & p\tripow2 \ar[d] & p\tripow3 \ar[d] & \cdots \\
    \1 & p\tri\1 \ar[l] & p\tripow2\tri\1 \ar[l] & p\tripow3\tri\1 \ar[l] & \cdots \ar[l]
\end{tikzcd}
\end{equation}
commutes as well: its leftmost rectangle is \eqref{eqn.phi_seed}, and taking the composition product of each rectangle in \eqref{eqn.phi_expand} with the commutative rectangle
\[
\begin{tikzcd}[column sep=large]
    \car{c} \ar[d, "\varphi"'] \ar[r, equals] & \car{c} \ar[d, "\varphi"] \\
    p \ar[d, equals] & p \ar[d, equals] \\
    p & p \ar[l, equals]
\end{tikzcd}
\]
yields the rectangle to its right.
As $\car{t}_p$ is the limit of the bottom two rows of \eqref{eqn.phi_expand}, it follows that there is an induced lens $F\colon\car{c}\to\car{t}_p$ that, when composed with each projection $\epsilon^{(n)}_p\colon\car{t}_p\to p\tripow{n}$, yields the lens depicted in \eqref{eqn.phi_expand} from $\car{c}$ to $p\tripow{n}$.
This lens is the composite of the lens $\car{c}\to\car{c}\tripow{(n+1)}$ in the top row, which by \cref{prop.n_duplication} is the canonical lens $\delta^{(n+1)}$ associated with the comonoid $\cat{C}$, composed with the lens $\epsilon\tri\varphi\tripow{n}\colon\car{c}\tripow{(n+1)}\to p\tripow{n}$.
%** check that it satisfies retrofunctor laws

Next, we prove uniqueness: that a retrofunctor $\cat{C}\cof\cofree{p}$ with underlying lens $f\colon\car{c}\to\car{t}_p$ is completely determined by the value of $f\then\epsilon^{(1)}_p$.
It suffices to show that we can recover the $n^\text{th}$ component of $g$ from its first component.
**
% We do so by demonstrating that the diagram
% \begin{diag} \label{diag.g_components}
% \begin{tikzcd}[sep=large]
%     C \ar[d, "\d^{n-1}"'] \ar[r, "g"] & C_p \ar[d, "\d^{n-1}_p"'] \ar[dr, "\e^n_p"] \\
%     C\circpow{n} \ar[r, "g\circpow{n}"'] & C_p\circpow{n} \ar[r, "\left(\e^1_p\right)\circpow{n}"'] & p\circpow{n}
% \end{tikzcd}
% \end{diag}
% commutes, which would show that $g \sm \e^n_p = \d^{n-1} \sm \left( g \sm \e^1_p \right)\circpow{n}$.

% The square on the left side of \cref{diag.g_components} commutes because $g$ respects the comonoid structure, and we can show that the triangle on the right side of \cref{diag.g_components} commutes via induction on $n$.
% When $n = 0$, the triangle commutes because $\delta^{-1}_p$ and $\e^0_p$ are both just the counit of the comonoid $C_p$; and if the triangle commutes for some fixed $n$, then by \cref{eq.d_pow_n,eq.d_m_n},
% \[
% \begin{split}
%     \d^n_p \sm \left(\e^1_p\right)\circpow{(n+1)} &= \d_p \sm \left(C_p \circ \d^{n-1}_p\right) \sm \left(\e^1_p\right)\circpow{(n+1)} \\
%     &= \d_p \sm \left(\e^1_p \circ \left(\d^{n-1}_p \sm \left(\e^1_p\right)\circpow{n}\right)\right) \\
%     &= \d_p \sm \left(\e^1_p \circ \e^n_p \right) \\
%     &= \e^{n+1}_p.
% \end{split}
% \]
% Hence the triangle commutes by induction, so \cref{diag.g_components} commutes, as desired.
\end{proof}

% ** include examples

%---- Subsection ----%
\subsection{The many (inter)faces of the cofree comonoid} \label{subsec.comon.cofree.cons.faces}
\index{cofree comonoid!interfaces for}\index{interface}

The forgetful-cofree adjunction of \cref{thm.cofree} tells us that given a category $\cat{C}\in\catsharp$ with carrier $\car{c}\coloneqq U\cat{C}$ and a polynomial $p\in\poly$, there is a natural isomorphism\index{isomorphism!natural}
\[
    \poly(\car{c},p)\iso\catsharp(\cat{C},\cofree{p}).
\]
So every lens $\varphi\colon\car{c}\to p$ has a corresponding retrofunctor $F\colon\cat{C}\cof\cofree{p}$ that we call its \emph{mate}.

We can view $\varphi$ as a dynamical system with an interface $p$ and a generalized state system $\car{c}$, carrying an arbitrary category $\cat{C}$ of states and transitions.
Then the lens $\text{Run}_n(\varphi)\colon\car{c}\to p\tripow{n}$ for $n\in\nn$, defined as the composite
\[
    \car{c}\To{\delta^{(n)}}\car{c}\tripow{n}\To{\varphi\tripow{n}}p\tripow{n},
\]
models $n$ runs through the system $\varphi$.

Meanwhile, as the limit of the diagram \eqref{eqn.cofree_diagram} with the row of polynomials of the form $p\tripow{n}$ across the top, the carrier $\car{t}_p$ of the cofree comonoid on $p$ comes equipped with a lens $\epsilon^{(n)}_p\colon\car{t}_p\to p\tripow{n}$ for each $n\in\nn$.
Since the retrofunctor $F\colon\cat{C}\cof\cofree{p}$ has an underlying lens between carriers $f\coloneqq UF\colon\car{c}\to\car{t}_p$, we can obtain another lens $\car{c}\to\tripow{n}$ as the composite
\[
    \car{c}\To{f}\car{t}_p\To{\epsilon^{(n)}_p}p\tripow{n}.
\]

It follows from our forgetful-cofree adjunction that these two composites are equal.

\begin{proposition}
With the definitions above, the following diagram commutes:
\[
\begin{tikzcd}
    \car{c} \ar[r,"f"] \ar[d, "\delta^{(n)}"'] & \car{t}_p \ar[d, "\epsilon^{(n)}_p"] \\
    \car{c}\tripow{n} \ar[r, "\varphi\tripow{n}"'] & p\tripow{n}.
\end{tikzcd}
\]
\end{proposition}
%\begin{proof}
%**
%\end{proof}

Now we have a better sense of what we mean when we say that $F\colon\cat{C}\cof\cofree{p}$ captures all the information that the lenses $\text{Run}_n(\varphi)\colon\car{c}\to p\tripow{n}$ encode: the category $\cofree{p}$ is carried by a polynomial $\car{t}_p$ equipped with lenses $\epsilon^{(n)}_p\colon\car{t}_p\to p\tripow{n}$, each of which exposes a part of the category as the $n$-fold interface $p\tripow{n}$.
All together, $\car{t}_p$ acts as a giant interface that captures the $n$-fold behavior of $p\tripow{n}$ for every $n\in\nn$.
But to see all this explicitly, let's consider some examples.

\begin{example}\label{ex.halt_dsa_accept}\index{deterministic state automaton!halting}
Let $S\coloneqq\{\bul[black],\bul[my-yellow],\bul[my-red]\}$ and $p\coloneqq\yon^\2+\1$, and consider the dynamical system $\varphi\colon S\yon^S\to p$ modeling the halting deterministic state automaton from \cref{exc.halt_dsa}, depicted here again for your convenience:
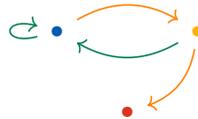
\[ %TODO
\begin{tikzcd}[column sep=small]
	\bul[my-blue]\ar[rr, bend left, orange]\ar[loop left, dgreen]&&
	\bul[my-yellow]\ar[dl, bend left, orange]\ar[ll, dgreen, bend left]\\&
	\bul[my-red]
\end{tikzcd}
\]
Under the forgetful-cofree adjunction, the lens $\varphi$ coincides with a retrofunctor $F\colon S\yon^S\cof\cofree{p}$ from the state category on $S$ to the category of $p$-trees.
The retrofunctor sends each state in $S=\{\bul[dgreen],\bul[my-yellow],\bul[my-red]\}$ to a $p$-tree; these $p$-trees are drawn below (the first two are infinite).
Then the vertices of each $p$-tree are sent back to a morphism in the state category; the color of each vertex indicates the codomain of the morphism to which that vertex is sent.
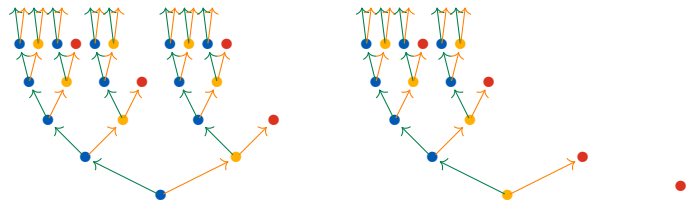
\[ %%COLOR
\begin{tikzpicture}[trees, scale=1,
  level 1/.style={sibling distance=20mm},
  level 2/.style={sibling distance=10mm},
  level 3/.style={sibling distance=5mm},
  level 4/.style={sibling distance=2.5mm},
  level 5/.style={sibling distance=1.25mm}]
  \node[my-blue] (a) {$\bullet$}
    child[dgreen] {node[my-blue] {$\bullet$}
    	child[dgreen] {node[my-blue] {$\bullet$}
    		child[dgreen] {node[my-blue] {$\bullet$}
  				child[dgreen] {node[my-blue] {$\bullet$}
    				child[dgreen] {}
    				child[orange] {}
    			}
  				child[orange] {node[my-yellow] {$\bullet$}
    				child[dgreen] {}
    				child[orange] {}
    			}
  			}
    		child[orange] {node[my-yellow] {$\bullet$}
				child[dgreen] {node[my-blue] {$\bullet$}
      			    child[dgreen] {}
      			    child[orange] {}
     			}
    			child[orange]  {node[my-red] {$\bullet$}}
  			}
    	}
    	child[orange] {node[my-yellow] {$\bullet$}
    		child[dgreen] {node[my-blue] {$\bullet$}
  				child[dgreen] {node[my-blue] {$\bullet$}
    				child[dgreen] {}
    				child[orange] {}
    			}
  				child[orange] {node[my-yellow] {$\bullet$}
    				child[dgreen] {}
    				child[orange] {}
    			}
  			}
    		child[orange]  {node[my-red] {$\bullet$}}
    	}
    }
    child[orange] {node[my-yellow] {$\bullet$}
    	child[dgreen] {node[my-blue] {$\bullet$}
    		child[dgreen] {node[my-blue] {$\bullet$}
  				child[dgreen] {node[my-blue] {$\bullet$}
    				child[dgreen] {}
    				child[orange] {}
    			}
  				child[orange] {node[my-yellow] {$\bullet$}
    				child[dgreen] {}
    				child[orange] {}
    			}
  			}
    		child[orange] {node[my-yellow] {$\bullet$}
				child[dgreen] {node[my-blue] {$\bullet$}
          			child[dgreen] {}
          			child[orange] {}
     			}
    			child[orange] {node[my-red] {$\bullet$}}
  			}
  		}
  		child[orange] {node[my-red] {$\bullet$}
  		}
  	}
  ;
\end{tikzpicture}
\hspace{10mm}
\begin{tikzpicture}[trees, scale=1,
  level 1/.style={sibling distance=20mm},
  level 2/.style={sibling distance=10mm},
  level 3/.style={sibling distance=5mm},
  level 4/.style={sibling distance=2.5mm},
  level 5/.style={sibling distance=1.25mm}]
  \node[my-yellow] (a) {$\bullet$}
    child[dgreen] {node[my-blue] {$\bullet$}
    	child[dgreen] {node[my-blue] {$\bullet$}
    		child[dgreen] {node[my-blue] {$\bullet$}
  				child[dgreen] {node[my-blue] {$\bullet$}
    				child[dgreen] {}
    				child[orange] {}
    			}
  				child[orange] {node[my-yellow] {$\bullet$}
    				child[dgreen] {}
    				child[orange] {}
    			}
  			}
    		child[orange] {node[my-yellow] {$\bullet$}
					child[dgreen] {node[my-blue] {$\bullet$}
      			child[dgreen] {}
      			child[orange] {}
     			}
    			child[orange]  {node[my-red] {$\bullet$}}
  			}
    	}
    	child[orange] {node[my-yellow] {$\bullet$}
    		child[dgreen] {node[my-blue] {$\bullet$}
  				child[dgreen] {node[my-blue] {$\bullet$}
    				child[dgreen] {}
    				child[orange] {}
    			}
  				child[orange] {node[my-yellow] {$\bullet$}
    				child[dgreen] {}
    				child[orange] {}
    			}
  			}
    		child[orange]  {node[my-red] {$\bullet$}}
    	}
    }
    child[orange] {node[my-red] {$\bullet$}
  	}
  ;
\end{tikzpicture}
\hspace{10mm}
\begin{tikzpicture}
	\node[my-red] {$\bullet$};
\end{tikzpicture}
\]
Each $p$-tree above, with its vertices colored, thus encodes all the ways to navigate the automaton.
In particular, the \emph{maximal rooted paths} of these trees (i.e.\ those that terminate at a leaf) trace out all the ways in which the automaton can halt, and therefore all the words that the automaton accepts.
Notice, too, that every $p$-subtree of any one of these $p$-trees is another one of these three $p$-trees---do you see why?

In general, given a lens $\varphi\colon S\yon^S\to\yon^A+\1$ modeling a halting deterministic state automaton and an initial state $s_0\in S$, the $(\yon^A+\1)$-tree to which the corresponding retrofunctor $F\colon S\yon^S\cof\cofree{\yon^A+\1}$ sends $s_0$ encodes the set of words accepted by the automaton with that initial state in its maximal rooted paths.
\end{example}
\index{deterministic state automaton!language of}

\begin{example}[Languages recognized by deterministic state automata] \label{ex.dsa_lang_recog}
Recall from \cref{prop.dsa} that a deterministic state automaton with a set of states $S$ and a set of symbols $A$ can be identified with a lens $\yon\to S\yon^S$, indicating the initial state $s_0\in S$, and a lens $\varphi\colon S\yon^S\to\2\yon^A$, indicating the subset of accept states $F\coloneqq\varphi_1^{-1}(2)\ss S$ and the update function $u\colon S\times A\to S$ via $u(s,a)=\varphi^\sharp_s(a)$.

Under the forgetful-cofree adjunction, the lens $\varphi$ coincides with a retrofunctor $F\colon S\yon^S\cof\cofree{\2\yon^A}$.
By \cref{exc.b-lab_a-ary_trees}, $\cofree{\2\yon^A}$ is the category of $\2$-labeled $A$-ary trees.
So $F(s_0)$ is an element of $\tr_{\2\yon^A}\iso\2^{\lst(A)}$: an $A$-ary tree where each rooted path corresponds to a list of elements of $A$ and bears one of the elements of $\2$, indicating whether the course through the automaton corresponding to that rooted path ends at an accept state or not.
Equivalently, an element of $\2^{\lst(A)}$ is a subset of $\lst(A)$; in the case of $F(s_0)$, it is the subset containing exactly the words in $\lst(A)$ accepted by the automaton.
So on objects, $F$ sends each possible start state $s_0\in S$ of the automaton to the subset of words that the automaton would then accept!
Backward on morphisms, $F^\sharp_{s_0}$ sends every possible word to the state the automaton would reach if it followed that word starting from $s_0$.
\end{example}

\index{cofree comonoid!and Moore machines}

\begin{example}[Direction sequences to position sequences] \label{ex.input_output}
We interpret our dynamical systems as converting sequences of directions to sequences of positions.
The forgetful-cofree adjunction allows us to express this conversion formally in the language of $\poly$.
For convenience, we'll focus on the example of an $(A,B)$-Moore machine $\varphi\colon S\yon^S\to B\yon^A$, although of course we can generalize this beyond monomial interfaces.

\index{interface!monomial}

The lens $\varphi$ corresponds to a retrofunctor $F\colon S\yon^S\cof\cofree{B\yon^A}$.
By \cref{exc.b-lab_a-ary_trees}, $\cofree{B\yon^A}$ is the category of $B$-labeled $A$-ary trees; in paticular, its carrier is $\car{t}_{B\yon^A}\iso B^{\lst(A)}\yon^{\lst(A)}$.

Then for every initial state $s_0\in S$, the $B$-labeled $A$-ary tree $F(s_0)$ can be interpreted as a decision tree of all the possible sequences of directions that the system may receive.
The label in $B$ corresponding to the vertex (or rooted path) specified by each sequence in $\lst(A)$ tells us the final position that the system returns when that sequence is fed in as directions.
Put another way, $F(s_0)\in B^{\lst(A)}$ can be interpreted as a function $\lst(A)\to B$.
So if the direction sequence is $(a_1,\ldots,a_n)\in A^n\ss\lst(A)$, then the corresponding position sequence $(b_0,\ldots,b_n)\in B^{n+1}$ is given (non-recursively!) by
\[
    b_i\coloneqq F(s_0)(a_1,\ldots,a_i).
\]
Finally, $F^\sharp_{s_0}(a_1,\ldots,a_i)\in S$ then corresponds to the system's state after that sequence of directions is given.
\end{example}\index{decision tree}

\begin{example}
  In \cref{ex.halt_dsa_accept} we had a dynamical system with $S\coloneqq\{\bul[dgreen],\bul[my-yellow],\bul[my-red]\}$ and $p\coloneqq\yon^\2+\1$, and $\varphi\colon S\yon^S\to p$ from \cref{exc.halt_dsa}, depicted here again for your convenience:
  
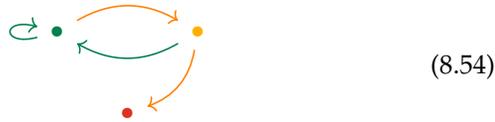
\begin{equation}\label{eqn.dyn_sys_misc573}
    \begin{tikzcd}[column sep=small]
      \bul[dgreen]\ar[rr, bend left, orange]\ar[loop left, dgreen]&&
      \bul[my-yellow]\ar[dl, bend left, orange]\ar[ll, dgreen, bend left]\\&
      \bul[my-red]
    \end{tikzcd}
  \end{equation}
  Under the forgetful-cofree adjunction (\cref{thm.cofree}), the lens $\varphi$ coincides with a retrofunctor $F\colon S\yon^S\cof\cofree{p}$ from the state category on $S$ to the category of $p$-trees.
  We can now see this as a copresheaf on the category $\cofree{p}$ itself.\index{copresheaf!on cofree category}

  The cofree category $\cofree{t}$ is actually the free category on a graph, as we saw in \cref{prop.cofree_free_on_graph}, and so the schema is easy. There is one table for each tree (object in $\cofree{p}$), e.g.\ we have a table associated to this tree:
  \[
  \treepic
  \]
  The table has two columns, say left and right, corresponding to the two arrows emanating from the root node. The left column refers back to the same table, and the right column refers to another table (the one corresponding to the yellow dot).

  Again, there are infinitely many tables in this schema. Only three of them have data in them; the rest are empty. We know in advance that this instance has three rows in total, since $|S|=3$.
\end{example}

%---- Subsection ----%
\subsection{Morphisms between cofree comonoids}
\index{cofree comonoid!morphisms between}\index{functor!cofree comonad}

Given a lens $\varphi \colon p \to q$, the cofree functor gives us a comonoid morphism $\cofree{\varphi}\colon\cofree{p}\to\cofree{q}$ as follows.

An object $t \in \tr_p$ is a tree; the tree $u \coloneqq \cofree{\varphi} (t) \in \tr_q $ is constructed recursively as follows. If the root of $t$ is $i \in p(\1)$ then the root of $u$ is $j \coloneqq \varphi_1 (i)$. To each branch $b \in q[j]$, we need to assign a new tree, and we use the one situated at $\varphi_i ^ \sharp (b)$.

\begin{exercise}
Let $p \coloneqq \{\bul[dgreen]\}\yon ^\2 + \{\bul[my-red]\}$ and $q \coloneqq \{\bul[dgreen],\bul[my-yellow]\} \yon + \{\bul[my-red],\bul[black]\}$.
\begin{enumerate}
    \item Choose a lens $\varphi\colon p\to q$, and write it out.
    \item Choose a tree $T\in\tr_p$ with at least height $3$.
    \item What is $\cofree{\varphi}(T)$?
    \qedhere
\end{enumerate}
\begin{solution}
Recall that $p \coloneqq \{\bul[dgreen]\}\yon ^\2 + \{\bul[my-red]\}$ and $q \coloneqq \{\bul[dgreen],\bul[my-yellow]\} \yon + \{\bul[my-red],\bul[black]\}$.
\begin{enumerate}
    \item Have $\varphi\colon p\to q$ send forward $\bul[dgreen]\mapsto\bul[dgreen]$, with the unique element of $q[\bul[dgreen]]$ sent back to the left branch of $p[\bul[dgreen]]$, and send forward $\bul[my-red]\mapsto\bul[my-red]$.
    \item Here is a sample $T\in\tr_p$.
  \[
\begin{tikzpicture}[trees, scale=1.5,
  level 1/.style={sibling distance=10mm},
  level 2/.style={sibling distance=5mm}]
  \node[dgreen] (a) {$\bullet$}
    child {node[dgreen] {$\bullet$}
    	child {node[my-red] {$\bullet$}
			}
			child{node[dgreen] {$\bullet$}
				child
				child
			}
		}
    child {node[my-red] {$\bullet$}
    }
  ;
\end{tikzpicture}
\]
    \item Here is $\cofree{\varphi}(T)$:
  \[
\begin{tikzpicture}[trees, scale=1.5,
  level 1/.style={sibling distance=10mm},
  level 2/.style={sibling distance=5mm}]
  \node[dgreen] (a) {$\bullet$}
    child {node[dgreen] {$\bullet$}
    	child {node[my-red] {$\bullet$}
			}
		}
  ;
\end{tikzpicture}
\]
\end{enumerate}
\end{solution}
\end{exercise}

% The following exercise is useful when considering the (topos-theoretic) logic of dynamical systems. Namely, it will allow us to specify legal subtrees of height $n$. %** what does this mean?

\begin{exercise}
Let $p$ be a polynomial.
\begin{enumerate}
    \item Show there is an induced retrofunctor $\cofree{p} \to \cofree{p\tripow{n}}$ for all $n\in\nn$.
    \item When $n\geq1$, is the induced retrofunctor is an isomorphism? \qedhere
\end{enumerate}
\begin{solution}
\begin{enumerate}
    \item By the adjunction from \cref{thm.cofree}, retrofunctors $\cofree{p} \to \cofree{p\tripow{n}}$ are in bijection with lenses $\car{t}_{p}\to p\tripow{n}$, where $\car{t}_{p}=U\cofree{p}$ is the carrier.

    From \cref{prop.n_duplication} we have a lens $\delta^{(n)}\colon\cofree{p}\to\cofree{p}\tripow{n}$, and we also have the counit of the cofree adjunction which we'll temporarily call $r\colon\cofree{p}\to p$. Then the desired lens is $(\delta^{(n)}\then r\tripow{n})\colon\car{t}_{p}\to p\tripow{n}$.
    \item When $n=1$ the induced retrofunctor is an isomorphism, but for $n\geq 2$ it is not. However, we should start by saying that the function on objects $\tr_p\to\tr_{p\tripow{n}}$ is a bijection. It sends a $p$-tree to the $p\tripow{n}$ tree that simply compresses every height-$n$ segment into a single vertex, labeled by that height-$n$ segment. To go back, just decompress the segments. But this is not a bijection on maps, because every rooted path on $\tr_{p\tripow{n}}$ would correspond to a rooted path on $\tr_p$ whose length is a multiple of $n$. When $n=2$ and $p=\yon+1$, some rooted paths in $\tr_p$ have odd lengths.
\end{enumerate}
\end{solution}
\end{exercise}

%---- Subsection ----%
\subsection{Some categorical properties of cofree comonoids}

\index{cofree comonoid!as free category on a graph}\index{graph!free category on}

\begin{proposition}\label{prop.cofree_free_on_graph}
For every polynomial $p$, the cofree category $\cofree{p}$ is free on a graph. That is, there is a graph $G_p$ whose associated free category in the usual sense (the category of vertices and paths in $G_p$) is isomorphic to $\cofree{p}$.
\end{proposition}
\begin{proof}
For vertices, we let $V_p$ denote the set of $p$-trees,
\[V_p\coloneqq\tr_p(\1).\]
For arrows we use the counit lens $\pi\colon\tr_p\to p$ from \cref{thm.cofree} to define
\[
A_p\coloneqq\sum_{t\in\tr_p(\1)}p[\pi_1(t)]
\]
In other words $A_p$ is the set $\{d\in p[\pi_1(t)]\,\mid\,t\in\tr_p\}$ of directions in $p$ that emanate from the root corolla of each $p$-tree. The source of $(t,d)$ is $t$ and the target is $\cod(\pi^\sharp_t(d))$. It is clear that every morphism in $\cofree{t}$ is the composite of a finite sequence of such morphisms, completing the proof.
\end{proof}

\begin{corollary}
Let $p$ be a polynomial and $\cofree{p}$ the cofree comonoid. Every morphism in $\cat{C}_p$ is both monic and epic.
\end{corollary}
\begin{proof}
The free category on a graph always has this property, so the result follows from \cref{prop.cofree_free_on_graph}.
\end{proof}\index{graph!free category on}

\begin{proposition}\label{prop.ynn_monoid}
The additive monoid $\yon^\nn$ of natural numbers has a $\times$-monoid structure in $\catsharp$.
\end{proposition}
\begin{proof}
The right adjoint $p\mapsto\cofree{p}$ preserves products, so $\yon^{\List(\ord{n})}\cong\cofree{\yon^{\ord{n}}}$ is the $n$-fold product of $\yon^\nn$ in $\catsharp$. We thus want to find retrofunctors $e\colon \yon\to\yon^\nn$ and $m\colon\yon^{\List(\2)}\to\yon^\nn$ that satisfy the axioms of a monoid.

The unique lens $\yon\to\yon^\nn$ is a retrofunctor (it is the mate of the identity $\yon\to\yon$). We take $m$ to be the mate of the lens $\yon^{\List(\2)}\to\yon$ given by the list $[1,2]$. One can check by hand that these definitions make $(\yon^\nn,e,m)$ a monoid in $(\catsharp,\yon,\times)$.
\end{proof}

Recall from \cref{ex.arrow_field} that an arrow field of a category $\cat{C}$ is a retrofunctor $\cat{C}\cof\yon^\nn$.\index{arrow field}

\begin{corollary}
For any category $\cat{C}$, the set $\catsharp(\cat{C},\yon^\nn)$ of arrow fields has the structure of a monoid. Moreover, this construction is functorial
\[\catsharp(-,\yon^\nn)\colon\catsharp\to\Cat{Mon}\op\]
\end{corollary}
\begin{proof}
We saw that $\yon^\nn$ is a monoid object in \cref{prop.ynn_monoid}.
\end{proof}

A retrofunctor $\cat{C}\cof\yon^\nn$ is a policy in $\cat{C}$: it assigns an outgoing morphism to each object of $\cat{C}$. Any two such trajectories can be multiplied: we simply do one and then the other; this is the monoid operation. The policy assigning the identity to each object is the unit of the monoid.

We use the notation $\cat{C}\mapsto\vec{\cat{C}}$ for the monoid of arrow fields.\index{arrow field!monoid of arrow fields}

\index{arrow fields!as adjoint}

\begin{theorem}\label{thm.catsharp_to_mon}\index{functor!arrow fields}
The arrow fields functor
\[\catsharp\to\Cat{Mon}\op\]
is right adjoint to the inclusion $\Cat{Mon}\op\to\catsharp$ from \cref{prop.monoids_ff}.
\end{theorem}
\begin{proof}
Let $\cat{C}$ be a category and $(M,e,*)$ a monoid. A retrofunctor $F\colon\cat{C}\cof\yon^M$ has no data on objects; it is just a way to assign to each $c\in \cat{C}$ and $m\in M$ a morphism $F^\sharp_c(m)\colon c\to c'$ for some $c'\coloneqq\cod(F^\sharp_c(m))$. This assignment must send identities to identities and composites to composites: given $m'\in M$ we have $F^\sharp_c(m\then m')=F^\sharp_c(m)\then F^\sharp_{c'}(m')$. This is exactly the data of a monoid morphism $M\to \vec{\cat{C}}$: it assigns to each $m\in M$ an arrow field $\cat{C}$, preserving unit and multiplication.
\end{proof}

\begin{proposition}\label{prop.traj_mon_poly}
There is a commutative square of left adjoints
\[
\begin{tikzcd}
	\Cat{Mon}\op\ar[r, "U"]\ar[d, "\yon^-"']&
	\smset\op\ar[d, "\yon^-"]\\
	\catsharp\ar[r, "U"']&
	\poly
\end{tikzcd}
\]
where the functors denoted $U$ are forgetful functors.
\end{proposition}\index{functor!forgetful}
\begin{proof}
Using the fully faithful functor $\yon^-\colon\Cat{Mon}\op\fromto\catsharp$ from \cref{prop.monoids_ff}, it is easy to check that the above diagram commutes.

The free-forgetful adjunction $\smset\fromto\Cat{Mon}$ gives an opposite adjunction $\smset\op\fromto\Cat{Mon}\op$, where $U$ is now left adjoint. We saw that $\yon^-\colon\smset\op\to\poly$ is a left adjoint in \cref{prop.yoneda_left_adjoint}, that $U\colon\catsharp\to\poly$ is a left adjoint in \cref{thm.cofree}, and that $\yon^-\colon\Cat{Mon}\to\catsharp$ is a left adjoint in \cref{thm.catsharp_to_mon}.
\end{proof}

\index{cofree comonoid|)}

\section{More categorical properties of $\catsharp$}

Many of the properties of $\poly$ we covered in \cref{ch.poly.cat,ch.poly.bonus} have analogues in $\catsharp$; we review these here.

%-------- Section --------%
\subsection{Other special comonoids and adjunctions}

We begin by highlighting a few other adjunctions involving $\catsharp$, as well as the special comonoids in $\catsharp$ they provide.

\index{polynomial comonoid!representable}

\begin{proposition}
The functor $\yon^-$ from \cref{prop.traj_mon_poly} factors through an isomorphism of categories
\[
    \catsharp_{\text{rep}}\iso\Cat{Mon}\op,
\]
where $\Cat{Mon}$ is the category of monoids and monoid homomorphisms and  $\catsharp_{\text{rep}}$ is the full subcategory of $\catsharp$ consisting of categories with representable carriers $\yon^M$ for some $M\in\smset$.
\end{proposition}
\begin{proof}
Let $\cat{C}$ be a category. It has only one object iff its carrier $\car{c}$ has only one position, i.e.\ $\car{c}\cong\yon^M$ for some $M\in\smset$, namely where $M$ is the set of morphisms in $\cat{C}$. It remains to show that retrofunctors between monoids are dual---opposite---to morphisms between monoids.

A retrofunctor $f\colon\yon^M\to\yon^N$ involves a single function $f^\sharp\colon N\to M$ that must satisfy a law coming from unitality and one coming from composition, as in \cref{def.morphism_comonoids}. The result can now be checked by hand, or seen formally as follows. Each object in the two diagrams of \eqref{def.morphism_comonoids} is representable by \cref{exc.composites_of_specials}. The Yoneda embedding $\smset\op\to\poly$ is fully faithful, so these two diagrams are equivalent to the unit and composition diagrams for monoid homomorphisms.
\end{proof}\index{Yoneda lemma}

\index{polynomial comonoid!linear}

\begin{exercise}\label{exc.lin_comon_set}
Let $\catsharp_{\text{lin}}$ be the full subcategory of $\catsharp$ consisting of categories with linear carriers $S\yon$ for some $S\in\smset$.
Show that there is an isomorphism of categories
\[
\catsharp_{\text{lin}}\iso\smset.=
\]
\begin{solution}
A category $S\yon$ with linear carrier is a discrete category on $S$, so we need to show that a retrofunctor between the discrete category on $S$ and the discrete category on $T$ is the same thing as a function $S\to T$. But this is clear: a retrofunctor $S\yon\to T\yon$ is a function $S\to T$ on objects and a function backwards on morphisms, and the latter is unique because each object in $S\yon$ has a unique outgoing morphism.
\end{solution}
\end{exercise}

\begin{proposition}[Discrete categories]
The inclusion $\catsharp_{\text{lin}}\to\catsharp$ has a left adjoint sending each $(\car{c},\epsilon,\delta)\in\catsharp$ to the unique comonoid carried by $(\car{c}\tri\1)\yon$ in $\catsharp_{\text{lin}}$.
\end{proposition}
\begin{proof}
We need to show that for any comonoid $(\car{c},\epsilon,\delta)$ and set $A$, we have a natural isomorphism\index{isomorphism!natural}
\[
  \catsharp(\car{c},A\yon).
  \cong^?
  \catsharp((\car{c}\tri\1)\yon,A\yon)
\]
But every morphism in $A\yon$ is an identity, so the result follows from the fact that every retrofunctor must pass identities back to identities.
\end{proof}

%**
%* monoids, etc.

%-------- Section --------%
\subsection{Vertical-cartesian factorization of retrofunctors}

\index{factorization system!vertical cartesian}
\index{retrofunctor!vertical}\index{retrofunctor!cartesian}

A retrofunctor is called \emph{cartesian} if the underlying lens $f\colon\car{c}\to\car{d}$ is cartesian (i.e.\ for each position $i\in\car{c}(\1)$, the map $f^\sharp_i\colon\car{d}[f_1(i)]\to\car{c}[i]$ is an isomorphism). %** refer back to earlier section where we defined this...also vertical cartesian...and make this its own section

\begin{proposition}\label{prop.factor_retrofunctor}
Every retrofunctor $f\colon\cat{C}\cof\cat{D}$ factors as a vertical morphism followed by a cartesian morphism
\[
\cat{C}\overset{\text{vert}}{\cof}\cat{C}'\overset{\text{cart}}{\cof}\com{D}.
\]
\end{proposition}
\begin{proof}
A retrofunctor $\cat{C}\cof\cat{D}$ is a lens $\car{c}\to\car{d}$ satisfying some properties, and any lens $f\colon\car{c}\to\car{d}$ can be factored as a vertical morphism followed by a cartesian morphism
\[
	\car{c}\To{g}\car{c'}\To{h}\car{d}.
\]
For simplicity, assume $g_1\colon\car{c}(\1)\to\car{c'}(\1)$ is identity (rather than merely isomorphism) on positions and similarly that for each $i\in\car{c}$ the map $h^\sharp_i\colon\car{c'}[i]\to\car{d}[h_1(i)]$ is identity (rather than merely isomorphism) on directions.

It suffices to show that the intermediate object $\car{c'}$ can be endowed with the structure of a category such that $g$ and $h$ are retrofunctors. Given an object $i\in\car{c'}(\1)$, assign its identity to be the identity on $h_1(i)=f(i)$; then both $g$ and $h$ preserve identities because $f$ does. Given an emanating morphism $m\in\car{c'}[i]=\car{d}[f(i)]$, assign its codomain to be $\cod(m)\coloneqq\cod(f^\sharp_i(m))$, and given an emanating morphism $m'\in\car{c'}[\cod(m)]$, assign the composite $m\then m'$ in $\car{c'}$ to be $m\then m'$ in $\car{d}$. In \cref{exc.factor_retrofunctor} we will check that with these definitions, $\car{c'}$ is a category and both $g$ and $h$ are retrofunctors.
\end{proof}

\begin{exercise}\label{exc.factor_retrofunctor}
We will complete the proof of \cref{prop.factor_retrofunctor}, using the same notation.
\begin{enumerate}
	\item Show that composition is associative and unital in $\car{c'}$.
	\item Show that $g$ preserves codomains.
	\item Show that $g$ preserves compositions.
	\item Show that $h$ preserves codomains.
	\item Show that $h$ preserves compositions.
\qedhere
\end{enumerate}
\begin{solution}
\begin{enumerate}
    \item Given an object $i\in\car{c'}(1)=\car{c}(1)$, the set of emanating morphisms is $\car{c'}[i]\coloneqq\car{d}[fi]$, and they compose according to the structure of $\car{d}$. Since morphisms in $\car{d}$ compose associatively and unitaly, so do morphisms in $\car{c'}$.
    \item Given $i\in\car{c}(1)$ and $m\in \car{c'}[i]=\car{d}[fi]$, we have $g(\cod(g^\sharp_i(m)))=\cod(g^\sharp_i(m))=\cod(f^\sharp_i(m))=\cod(m)$ by definition.
    \item For any $i\in\car{c'}(1)$, the lens $g^\sharp_i$ preserves compositions of morphisms in $\car{c'}[i]=\car{d}[fi]$ because it agrees with $f^\sharp_i$, which preserves compositions.
    \item Given $i\in\car{c'}(1)$ and $m\in\car{d}[hi]$, we have $h(\cod(h^\sharp_i(m)))=f(\cod(f^\sharp_i(m)))=m$.
    \item The map $h^\sharp_i$ was chosen to be the identity for each $i$, so it certainly preserves compositions.
\end{enumerate}
\end{solution}
\end{exercise}

\index{discrete opfibration!as cartesian retrofunctor}

\begin{proposition} \label{prop.cart_dopf}
The wide%
\tablefootnote{A subcategory $\Cat{D}$ of a category $\Cat{C}$ is \emph{wide} if every object of $\Cat{C}$ is in $\Cat{D}$.}
subcategory of cartesian retrofunctors in $\catsharp$ is isomorphic to the wide subcategory of discrete opfibrations in $\smcat$.
\end{proposition}
\begin{proof}
Suppose that $\cat{C}$ and $\cat{D}$ are categories. Both a functor and a retrofunctor between them involve a map on objects, say $f\colon\Ob\cat{C}\to\Ob\cat{D}$. For any object $c\in\Ob\cat{C}$, a functor gives a function, say $f_\sharp\colon\cat{C}[c]\to\cat{D}[f(c)]$ whereas a retrofunctor gives a function $f^\sharp\colon\cat{D}[f(c)]\to\cat{C}[c]$. The retrofunctor is cartesian iff $f^\sharp$ is an iso, and the functor is a discrete opfibration iff $f_\sharp$ is an iso. We thus transform our functor into a retrofunctor (or vice versa) by taking the inverse function on morphisms. It is easy to check that this inverse appropriately preserves identities, codomains, and compositions.
\end{proof}

% ** connect the above to earlier discrete opfibration stuff; also see below

% The above correspondence is well-known; it remains to address the relationship between (2) and (3). A cartesian retrofunctor $(\varphi_1,\varphi^\sharp)\colon\cat{S}\cof\cat{C}$ gives a function $\varphi\colon\Ob(\cat{S})\to\Ob(\cat{C})$ and for each $s\in\cat{S}$ an isomorphism
% \[
%   \varphi_s^\sharp\colon\cat{C}[\varphi_1(s)]\To{\cong}\cat{S}[s]
% \]
% between the set of $\cat{S}$-morphisms emanating from $s$ and the set of $\cat{C}$-morphisms emanating from $\varphi_1(s)$. This isomorphism respects identities, codomains, and composites. As such we can define a functor that acts as $\varphi_1$ on objects and $(\varphi^\sharp)\inv$ on morphisms, and it is easily checked to be a discrete opfibration.

% Finally, given a discrete opfibration $\pi\colon\cat{S}\to\cat{C}$, we define $\varphi_1\coloneqq\Ob(\pi)$ to be its on-objects part, and for any $s\in\Ob(\cat{S})$ and emanating morphism $f\in\cat{C}[\varphi_1(s)]$, we define $\varphi^\sharp_s(f)\coloneqq\ol{f}$ to be the lift guaranteed by \cref{def.dopf}. The conversions between discrete opfibrations and cartesian retrofunctors are easily seen to be functorial and the roundtrips are identities.

\index{retrofunctor!vertical}
\begin{proposition}\label{prop.com_vert_cat_boo}
The wide subcategory of vertical maps in $\catsharp$ is isomorphic to the opposite of the wide subcategory bijective-on-objects maps in $\smcat$:
\[
\catsharp_{\text{vert}}\cong(\smcat_{\text{boo}})\op.
\]
\end{proposition}
\begin{proof}
Let $\cat{C}$ and $\cat{D}$ be categories. Given a vertical retrofunctor $F\colon\cat{C}\cof\cat{D}$, we have a bijection $F_1\colon\Ob\cat{C}\to\Ob\cat{D}$; let $G_1$ be its inverse. We define a functor $G\colon\cat{D}\to\cat{C}$ on objects by $G_1$ and, for any $f\colon d\to d'$ in $\cat{D}$ we define $G(f)\coloneqq F^\sharp_{G_1(d)}$. It has the correct codomain: $\cod(G(f))=G_1(F_1(\cod(G(f)))=G_1(\cod f)$. And it sends identities and compositions to identities and compositions by the laws of retrofunctors.

The construction of a vertical retrofunctor from a bijective-on-objects functor is analogous, and it is easy to check that the two constructions are inverses.
\end{proof}

\begin{exercise}
Let $S$ be a set and consider the state category $\cat{S}\coloneqq(S\yon^S,\epsilon,\delta)$. Use \cref{prop.com_vert_cat_boo} to show that categories $\cat{C}$ equipped with a vertical retrofunctor $\cat{S}\cof\cat{C}$ can be identified with categories whose set of objects has been identified with $S$.
\begin{solution}
The category corresponding to $S\yon^S$ is the contractible groupoid $K_s$ on $S$. A vertical retrofunctor $K_S\cof\cat{C}$ includes a bijection on objects, so we can assume that $\Ob(\cat{C})=S$. But the rest of the retrofunctor assigns to each map $s_1\to s_2$ in $\cat{C}$ some choice of morphism $s_1\to s_2$ in $K_S$, and there is exactly one. Thus the retrofunctor includes no additional data.
\end{solution}
\end{exercise}

\begin{exercise}\index{functor!bijective on objects (boo)}
Consider the categories $\cat{C}\coloneqq\fbox{$\bullet\tto\bullet$}$ and $\cat{D}\coloneqq\fbox{$\bullet\to\bullet$}$. There is a unique bijective-on-objects (boo) functor $F\colon\cat{C}\to\cat{D}$ and two boo functors $G_1,G_2\colon\cat{D}\to\cat{C}$. These have corresponding retrofunctors going the other way.
\begin{enumerate}
	\item Write down the morphism $\car{d}\to\car{c}$ of carriers corresponding to $F$.
	\item Write down the morphism $\car{c}\to\car{d}$ of carriers corresponding to either $G_1$ or $G_2$.
\qedhere
\end{enumerate}
\begin{solution}
Let's take the two nonidentity arrows in $\cat{C}$ to be labeled $s,t$ and the unique nonidentity arrow in $\cat{D}$ to be labeled $f$. Then the carrier of $\cat{C}$ is $\car{c}\coloneqq\yon^{\{\id_1,s,t\}}+\yon^{\{\id_2\}}$ and that of $\cat{D}$ is $\car{d}\coloneqq\yon^{\{\id_1,f\}}+\yon^{\{\id_2\}}$.
\begin{enumerate}
    \item The boo functor $F\colon\cat{C}\to\cat{D}$ can be identified with a vertical retrofunctor $\car{d}\to\car{c}$ that sends back $\id_1\mapsto\id_1$ and $s\mapsto f$ and $t\mapsto f$ and $\id_2\mapsto\id_2$.
    \item For $G_1\colon\cat{D}\to\cat{C}$, which we take to send $f\mapsto s$, the corresponding vertical retrofunctor $\car{c}\to\car{d}$ sends back $\id_1\mapsto\id_1$ and $f\mapsto s$ and $\id_2\mapsto\id_2$.
\end{enumerate}
\end{solution}
\end{exercise}

We record the following proposition here.%; it will be useful in \cref{cor.cartesian_cof_extra_adjoint}.

\index{cofree comonoid!preserves cartesian maps}

\begin{proposition}
If $\varphi\colon p\to q$ is a cartesian lens, then $\cofree{\varphi}\colon\cofree{p}\to\cofree{q}$ is a cartesian retrofunctor: that is, for each $t\in\tr_p$, the on-morphisms function
\[
    \left(\cofree{\varphi}\right)^\sharp_t\colon\cofree{q}[\cofree{\varphi}t]\To{\iso}\cofree{p}[t]
\]
is a bijection.
\end{proposition}
\begin{proof}
Given $\varphi\colon p\to q$ cartesian, each tree $T\in\tr_p$ is sent to a tree in $\tr_q$ with the same branching profile. A morphism emanating from it is just a finite rooted path, and the set of these is completely determined by the branching profile. Thus we have the desired bijection.
\end{proof}

%-------- Section --------%
\subsection{Limits and colimits of comonoids}
\index{polynomial comonoid!colimit of comonoids}\index{colimit!in $\catsharp$}

We saw in \cref{thm.poly_limits,thm.poly_colimits} that $\poly$ has all small limits and colimits.
It turns out that $\catsharp$ has all small limits and colimits as well.
We start by discussing colimits in $\catsharp$, as they are somewhat easier to get a handle on.

\subsubsection{Colimits in $\catsharp$}
\index{polynomial comonoid!colimits of|(}

It is a consequence of the forgetful-cofree adjunction from \cref{thm.cofree} that $\catsharp$ inherits all the colimits from $\poly$.
As these results are somewhat technical, relying on a property of functors known as \emph{comonadicity}, we defer their proofs to references.

\begin{proposition}[Porst]\index{forgetful functor|seealso{functor, forgetful}}\index{functor!comonadic}
The forgetful functor $\catsharp\to\poly$ is comonadic.
\end{proposition}
\begin{proof}
The fact that a forgetful functor $\catsharp\iso\Cat{Comon}(\poly)\to\poly$ is comonadic if it has a right adjoint follows from Beck's monadicity theorem via a straightforward generalization of an argument given by Par{\'e} in \cite[pp.~138-9]{pare1969absolute}, as pointed out by Porst in \cite[Fact~3.1]{porst2019colimits}.
So the result follows from \cref{thm.cofree}.
\end{proof}

\begin{corollary} \label{cor.comon_cocomp}\index{functor!forgetful $\catsharp\to\poly$}
The category $\catsharp$ has all small colimits.
They are created by the forgetful functor $\catsharp\to\poly$.
\end{corollary}
\begin{proof}
A comonadic functor creates all colimits that exist in its codomain (see \cite{nlab:created-limit}), and by \cref{thm.poly_colimits}, the category $\poly$ has all small colimits.
\end{proof}

\index{polynomial comonoid!coproduct of comonoids}\index{coproduct!of polynomial comonoids}

\begin{example}[Coproducts in $\catsharp$] \label{ex.comon_coprod}
Probably the most familiar kind of colimit in $\catsharp$ is the coproduct, as \cref{cor.comon_cocomp} tells us that it agrees with the usual coproduct from $\smcat$.
Here's why.

For concreteness, let $I$ be a set and $(\cat{C}_i)_{i\in I}$ be categories with carriers $(\car{c}_i)_{i\in I}$.
Then the coproduct of $(\cat{C}_i)_{i\in I}$ in $\smcat$ is the category $\sum_{i\in I}\cat{C}_i$, whose objects are given by the disjoint union of the objects in each summand, so that
\[
    \Ob\left(\sum_{i\in I}\cat{C}_i\right)\iso\sum_{i\in I}\Ob\cat{C}_i=\sum_{i\in I}\car{c}_i(\1)\iso\left(\sum_{i\in I}\car{c}_i\right)(\1),
\]
and whose morphisms out of each $c\in\Ob\cat{C}_j\ss\Ob\sum_{i\in I}\cat{C}_i$ are just the morphisms out of $c$ in the summand $\cat{C}_j$, so
\[
    \left(\sum_{i\in I}\cat{C}_i\right)[c]\iso\cat{C}_j[c]=\car{c}_j[c]\iso\left(\sum_{i\in I}\car{c}_i\right)[c]
\]
Hence $\sum_{i\in I}\cat{C}_i$ is carried by the coproduct of polynomials $\sum_{i\in I}\car{c}_i$.

It remains to show that $\sum_{i\in I}\cat{C}_i$ is also the coproduct of $(\cat{C}_i)_{i\in I}$ in $\catsharp$.
We already know from \cref{cor.comon_cocomp} that it has the right carrier: the coproduct of carriers $\sum_{i\in I}\car{c}_i$. It also has the right morphisms, for any object $(i,x)$ with $x\in\car{c}_i(1)$, the set of emanating morphisms is---and should be---$\car{c}_i[x]$.

% ** give inclusions as retrofunctors

% For each $j\in J$, we define an inclusion retrofunctor $\iota_j\colon\cat{C}_j\cof\sum_{i\in I}\cat{C}_i$ as follows: it is the canonical inclusion on objects, **. %**show all fully faithful functors can be interp as retrofunctors first??

% Suppose we have a category $\cat{D}$ and retrofunctors $F_i\colon\cat{C}_i\cof\cat{D}$ for $i\in I$.
% For every object $c\in\Ob\cat{C}_j\ss\Ob\sum_{i\in I}\cat{C}_i$, we have an object $F_j c\in\Ob\cat{D}$, so we can set $Fc\coloneqq F_j c$.
% Then for every morphism $g\colon Fc\to\_$ in $\cat{D}$, we have a morphism $\left(F_j\right)^\sharp_c g\colon c\to\_$, so we can set $F^\sharp_c g\coloneqq\left(F_j\right)^\sharp_c g$.
% As each $F_j$ satisfies the retrofunctor laws, so does $F$.

% % ** check retrofunctor laws
% % ** check commutes
% % ** check uniqueness

% % It is easy to check that the retrofunctor laws hold for $f$. Uniqueness of $f$ given $f_1,f_2$ is also straightforward.
\end{example}

\begin{exercise}\label{exc.0_initial_com}\index{functor!forgetful $\catsharp\to\poly$}
\begin{enumerate}
	\item Show that $\0$ has a unique comonoid structure.
	\item Explain why $\0$ with its comonoid structure is initial in $\catsharp$ in two ways: by explicitly showing it has the required universal property, or by invoking \cref{cor.comon_cocomp}.
\qedhere
\end{enumerate}
\begin{solution}
\begin{enumerate}
    \item We actually already showed that $\0$ has a unique comonoid structure, corresponding to the empty category (which we will also denote by $\0$), in \cref{exc.not_state_cat_but_same_carrier}, for the case of $S\coloneqq\0$.
    \item We show that $\0$ has the universal property of the initial object in $\catsharp$ as follows.
    For any category $\cat{D}$, there is a unique retrofunctor $\0\cof\cat{D}$: it vacuously sends each object in $\0$ to an object in $\cat{D}$, and since there are no objects in $\0$, it does nothing to morphisms.

    Alternatively, we know by \cref{cor.comon_cocomp} that $\catsharp$ has an initial object, and that the forgetful functor $\catsharp\to\poly$ sends it to the initial object in $\poly$, which is $\0$ (by \cref{prop.poly_coprods} and the following discussion).
    So the initial object in $\catsharp$ must be the unique comonoid carried by $\0$.
\end{enumerate}
\end{solution}
\end{exercise}

\begin{exercise}
Given a comonoid $(\car{c},\epsilon,\delta)\in\catsharp$, show that there is an induced comonoid structure on the polynomial $\2\car{c}$.
\begin{solution}
The coproduct of $\car{c}$ with itself is again a comonoid (as a category it is the disjoint union of $\car{c}$ with itself) and it has the right carrier polynomial, since the carrier functor $U\colon\catsharp\to\poly$ is a left adjoint.
\end{solution}
\end{exercise}
\index{polynomial comonoid!colimits of|)}

\subsubsection{Limits in $\catsharp$}

As in the case of colimits, there is a rather technical result showing that the forgetful-cofree adjunction implies the existence of all small limits of comonoids, which we summarize here.

\begin{corollary} \label{cor.comon_compl}\index{functor!forgetful $\catsharp\to\poly$}
The category $\catsharp$ has all small limits.
\end{corollary}
\begin{proof}
By \cref{thm.cofree}, the forgetful functor $U\colon\catsharp\to\poly$ has a right adjoint, and by \cref{thm.poly_limits}, $\poly$ itself has all small limits.
Furthermore equalizers in $\poly$ are connected limits, so by \cref{thm.connected_limits}, they are preserved by $\tri$ on either side.
Then the result follows from \cite[Fact~3.4]{porst2019colimits}.
\end{proof}\index{equalizer}

%-------- Section --------%
\subsection{Parallel product comonoids}

The usual product of categories is not the categorical product in $\catsharp$. %**show this with small counterexample
It is, however, a \emph{monoidal} product on $\catsharp$, coinciding with the parallel product $\otimes$ on $\poly$.

\index{parallel product!of polynomial comonoids}\index{polynomial comonoid!parallel product of}

\begin{proposition}\label{prop.parallel_on_catsharp}\index{functor!forgetful $\catsharp\to\poly$}
The parallel product $(\yon,\otimes)$ on $\poly$ extends to a monoidal structure $(\yon,\otimes)$ on $\catsharp$, such that the forgetful functor
$U\colon\catsharp\to\poly$
is strong monoidal with respect to $\otimes$.
The parallel product of two categories is their product in $\smcat$.
\end{proposition}
\begin{proof}
Let $\cat{C},\cat{D}\in\catsharp$ be categories corresponding to comonoids $(\car{c},\epsilon_{\car c},\delta_{\car c})$ and $(\car{d},\epsilon_{\car d},\delta_{\car d})$.

The carrier of $\cat{C}\otimes\cat{D}$ is defined to be $\car{c}\otimes\car{d}$. A position in it is a pair $(c,d)$ of objects, one from $\cat{C}$ and one from $\cat{D}$; a direction there is a pair $(f,g)$ of a morphism emanating from $c$ and one emanating from $d$.

We define $\epsilon_{\cat{C}\otimes\cat{D}}\colon\car{c}\otimes\car{d}\to\yon$ as
\[
\car{c}\otimes\car{d}\To{\epsilon_{\car{C}}\otimes\epsilon_{\car{D}}}\yon\otimes\yon\cong\yon.
\]
This says that the identity at $(c,d)$ is the pair of identities.

We define $\delta_{\cat{C}\otimes\cat{D}}\colon(\car{c}\otimes\car{d})\to(\car{c}\otimes\car{d})\tri(\car{c}\otimes\car{d})$ using the duoidal property:\index{duoidality}
\[
\car{c}\otimes\car{d}\To{\delta_{\car{c}}\otimes\delta_{\car{d}}}(\car{c}\tri\car{c})\otimes(\car{d}\tri\car{d})\to(\car{c}\otimes\car{d})\tri(\car{c}\otimes\car{d}).
\]
One can check that this says that codomains and composition are defined coordinate-wise, and that $(\car{c}\otimes\car{d},\epsilon_{\car{c}\otimes\car{d}},\delta_{\car{c}\otimes\car{d}})$ forms a comonoid. One can also check that this is functorial in $\cat{C},\cat{D}\in\catsharp$. See \cref{exc.parallel_on_catsharp}.
\end{proof}

\begin{exercise}\label{exc.parallel_on_catsharp}
We complete the proof of \cref{prop.parallel_on_catsharp}.
\begin{enumerate}
	\item Show that $(\car{c}\otimes\car{d},\epsilon_{\car{c}\otimes\car{d}},\delta_{\car{c}\otimes\car{d}})$, as described in \cref{prop.parallel_on_catsharp}, forms a comonoid.
	\item Check that the construction $(\cat{C},\cat{D})\mapsto\cat{C}\otimes\cat{D}$ is functorial in $\cat{C},\cat{D}\in\catsharp$.
\qedhere
\end{enumerate}
\begin{solution}
\begin{enumerate}
    \item The highbrow way to check that $\car{c}\otimes\car{d}$ forms a comonoid is to say that
    \[\otimes\colon(\poly,\yon,\tri)\times(\poly,\yon,\tri)\too(\poly,\yon,\tri)\]
     is a colax monoidal functor, so it sends comonoids to comonoids. A lowbrow way to check it is to see that the category corresponding to $\car{c}\otimes\car{d}$ is just the usual cartesian product $\cat{C}\times\cat{D}$ in $\smcat$. This is a category, hence a comonoid in $\poly$.
    \item Because $\otimes$ is symmetric, it suffices to show that for any retrofunctor $\varphi\colon\cat{C}\cof\cat{C}'$, there is an induced retrofunctor $\varphi\otimes\cat{D}\colon\cat{C}\otimes\cat{D}\cof\cat{C}'\otimes\cat{D}$.
\end{enumerate}
\end{solution}
\end{exercise}

The cofree construction works nicely with this monoidal product.

\index{parallel product!and the cofree construction}

\begin{proposition}\label{prop.cofree_lax_monoidal}
The cofree functor $p\mapsto\cofree{p}$ is lax monoidal; in particular there is a lens $\yon\to\car{t}_\yon$, and for any $p,q\in\poly$ there is a natural lens
\[
	\car{t}_p\otimes\car{t}_q\to\car{t}_{p\otimes q}.
\]
satisfying the usual conditions.
\end{proposition}
\begin{proof}
By \cref{prop.parallel_on_catsharp}, the left adjoint $U\colon\catsharp\to\poly$ is strong monoidal. A consequence of Kelly's doctrinal adjunction theorem \cite{kelly1974doctrinal} says that the right adjoint of an oplax monoidal functor is lax monoidal.
\end{proof}

\begin{exercise}
\begin{enumerate}
	\item What polynomial is $\car{t}_\yon$?
	\item What is the lens $\yon\to\car{t}_\yon$ from \cref{prop.cofree_lax_monoidal}?
	\item Explain in words how to think about the lens $\car{t}_p\otimes\car{t}_q\to\car{t}_{p\otimes q}$ from \cref{prop.cofree_lax_monoidal}, for arbitrary $p,q\in\poly$.
\qedhere
\end{enumerate}
\begin{solution}
\begin{enumerate}
    \item The carrier of $\cofree{\yon}$ is $\car{t}_\yon\coloneqq\yon^\nn$.
    \item The lens $\yon\to\car{t}\yon$ is the unique lens.
    \item Given a $p$-tree $T$ and a $q$-tree $U$, we get a $(p\otimes q)$-tree by having its root be labeled by the pair of root labels for $T$ and $U$, and for each child there---the set of which is the set of pairs consisting of a child in $T$ and a child in $U$---recursively using the above formula to combine these two children.
\end{enumerate}
\end{solution}
\end{exercise}

\section{Comodules over polynomial comonoids}

\index{comodule|(}\index{comodule!bicomodule|see{bicomodule}}

Just as we can define a category of comonoids within any monoidal category, we can further define a notion of comodules over such comonoids.
And much like how polynomial comonoids are categories, these comodules can also be described in categorical terms we are already familiar with.
There is much more to say about comodules over polynomial comonoids than we have room for here---we will merely give a glimpse of the theory and applications.

\subsection{Left and right comodules}
\index{comodule!left}

When the monoidal category is not symmetric, left comodules and right comodules may differ significantly, so we define them separately (but notice that the diagrams for one are analogous to the diagrams for the other).

\begin{definition}[Left comodule]\label{def.left_comod}
In a monoidal category $(\Cat{C},\yon,\tri)$, let $\cat{C}=(\car{c},\epsilon,\delta)$ be a comonoid.
A \emph{left $\cat{C}$-comodule} is
\begin{itemize}
    \item an object $m\in\Cat{C}$, called the \emph{carrier}, equipped with
    \item a morphism $\car{c}\tri m\From{\lambda}m$ called the \emph{left coaction},
\end{itemize}
such that the following diagrams, collectively known as the \emph{left comodule laws}, commute:
\begin{equation} \label{eqn.left_comod_laws}
\begin{tikzcd}[ampersand replacement=\&]
	\car{c}\tri m\ar[d, "\epsilon\:\tri\:m"']\&
	m\ar[l, "\lambda"']\ar[dl, equal, bend left]\\
	m
\end{tikzcd}
\hspace{.6in}
\begin{tikzcd}[ampersand replacement=\&]
	\car{c}\tri m\ar[d, "\delta\:\tri\:m"']\&
	m\ar[l, "\lambda"']\ar[d, "\lambda"]\\
	\car{c}\tri\car{c}\tri m\&
	\car{c}\tri m\ar[l, "\car{c}\:\tri\:\lambda"]
\end{tikzcd}
\end{equation}
When referring to a left $\cat{C}$-comodule, we may omit its coaction if it can be inferred from context (or simply unspecified), identifying the comodule with its carrier.

\index{comodule!morphism of}

A \emph{morphism} of left $\cat{C}$-comodules $m$ and $n$ is a morphism $\alpha\colon m\to n$ such that the following diagram commutes:
\[
\begin{tikzcd}
	\car{c}\tri m \ar[d, "\car{c}\:\tri\:\alpha"'] &
	m \ar[l] \ar[d, "\alpha"] \\
	\car{c}\tri n &
	n \ar[l]
\end{tikzcd}
\]
Here the top and bottom morphisms are the left coactions of $m$ and $n$.
\end{definition}

\begin{exercise} \label{exc.coalg_is_const_l_comod}\index{coalgebra}
Show that the category of $\cat{C}$-coalgebras from \cref{def.coalgebra} is exactly the category of \textit{constant} left $\cat{C}$-comodules---i.e.\ the full subcategory of the category of left $\cat{C}$-comodules spanned by those left $\cat{C}$-comodules whose carriers are constant polynomials.\index{constant polynomial!as carrying bicomodule}
\begin{solution}
In \cref{def.left_comod}, if the carrier of the left $\cat{C}$-comodule is always chosen to be a constant polynomial, i.e.\ a set, then we recover \cref{def.coalgebra}.
Hence constant left $\cat{C}$-comodules are precisely left $\cat{C}$-coalgebras; their notions of morphisms coincide as well, so the categories they form are isomorphic.
\end{solution}
\end{exercise}

\index{comodule!right}

\begin{definition}[Right comodule]\label{def.right_comod}
In a monoidal category $(\Cat{C},\yon,\tri)$, let $\cat{D}=(\car{d},\epsilon,\delta)$ be a comonoid.
A \emph{right $\cat{D}$-comodule} is
\begin{itemize}
    \item an object $m\in\Cat{C}$, called the \emph{carrier}, equipped with
    \item a morphism $m\To{\rho}m\tri\car{d}$ called the \emph{right coaction},
\end{itemize}
such that the following diagrams, collectively known as the \emph{right comodule laws}, commute:
\begin{equation} \label{eqn.right_comod_laws}
\begin{tikzcd}[ampersand replacement=\&]
	m\ar[r, "\rho"]\ar[dr, equal, bend right]\&
	m\tri\car{d}\ar[d, "m\:\tri\:\epsilon"]\\\&
	m
\end{tikzcd}
\hspace{.6in}
\begin{tikzcd}[ampersand replacement=\&]
	m\ar[r, "\rho"]\ar[d, "\rho"']\&
	m\tri\car{d}\ar[d, "m\:\tri\:\delta"]\\
	m\tri\car{d}\ar[r, "\rho\:\tri\:\car{d}"']\&
	m\tri\car{d}\tri\car{d}
\end{tikzcd}
\end{equation}
When referring to a right $\cat{D}$-comodule, we may omit its coaction if it can be inferred from context (or unspecified), identifying the comodule with its carrier.

A \emph{morphism} of right $\cat{D}$-comodules $m$ and $n$ is a morphism $\alpha\colon m\to n$ such that the following diagram commutes:
\[
\begin{tikzcd}
	m \ar[r] \ar[d, "\alpha"'] &
	m\tri\car{d} \ar[d, "\alpha\:\tri\:\car{d}"] \\
	n \ar[r] &
	n\tri\car{d}
\end{tikzcd}
\]
Here the top and bottom morphisms are the right coactions of $m$ and $n$.
\end{definition}

\begin{exercise}\label{exc.left_right_comodules}
\begin{enumerate}
	\item Draw the equations of \eqref{eqn.left_comod_laws} using polyboxes.
	\item Draw the equations of \eqref{eqn.right_comod_laws} using polyboxes.
\qedhere
\end{enumerate}
\begin{solution}
\begin{enumerate}
	\item Here is the first equation in polyboxes:
\[
\begin{tikzpicture}[polybox, mapstos]
  \node[poly, dom, "$m$" left] (l) {$x$\at$a\vphantom{a}$};
  \node[poly, cod, "$m$" right, right=of l] (r) {$x$\at$a$};
  \draw[double,-] (l_pos) -- (r_pos);
  \draw[double,-] (r_dir) -- (l_dir);
	\node[poly, dom, right=2.5 of r, "$m$" left] (p) {$\lambda^\sharp_a(\id_{|a|},x)$\at$a$};
	\node[poly, right=2.2cm of p.south, yshift=-1.25ex] (r') {$\id_{|a|}$\at$|a|$};
	\node[poly, cod, above=of r', "$m$" right] (q) {$x$\at$a.\id_{|a|}$};
  	\draw (p_pos) to[first] (r'_pos);
  	\draw (r'_dir) to[climb] node[left] {$\lambda$} (q_pos);
		\draw (q_dir) to[last] (p_dir);
		\draw (r'_pos) to[climb'] node[right] {$\epsilon_\car{c}$} (r'_dir);
	\node at ($(r.east)!.5!(p.west)$) {$=$};
\end{tikzpicture}
\]
The equation says that for any $a\in m(1)$ we have $a=a.\id_{|a|}\coloneqq\lambda(a,\id_{|a|})$ and that for any $x\in m[a]$ we have $x=x'\coloneqq\lambda^\sharp_a(\id_{|a|},x)$.

Here is the second equation in polyboxes:
\[
\begin{tikzpicture}
    \node (p1) {
        \begin{tikzpicture}[polybox, mapstos]
            \node[poly, dom, "$m$" left] (m') {$\lambda^\sharp_a((f\then f'),x)$\at$a$};
            \node[poly, right=2.3 of m'.south, yshift=-1ex, "$\car{c}$" below] (mm') {$f\then f'$\at$|a|$};
            \node[poly, above=of mm', "$m$" above] (C') {$x$\at$a.(f\then f')$};
            \node[poly, cod, right=2 of mm'.south, yshift=-1ex, "$\car{c}$" right] (D') {$f$\at$|a|$};
            \node[poly, cod, above=of D', "$\car{c}$" right] (mmm') {$f'$\at$\cod f$};
            \node[poly, cod, above=of mmm', "$m$" right] (CC') {$x$\at$a.(f\then f')$};
            \draw (m'_pos) to[first] (mm'_pos);
            \draw (mm'_dir) to[climb] node[left] {$\lambda$} (C'_pos);
            \draw (C'_dir) to[last]  (m'_dir);
            \draw[double,-] (mm'_pos) to[first] (D'_pos);
            \draw (D'_dir) to[climb] node[left] {$\delta$} (mmm'_pos);
            \draw (mmm'_dir) to[last] (mm'_dir);
            \draw[double,-] (C'_pos) to[first] (CC'_pos);
            \draw[double,-] (CC'_dir) to[last] (C'_dir);
        \end{tikzpicture}
	};
	\node (p2) [below=1 of p1] {
	    \begin{tikzpicture}[polybox, tos]
            \node[poly, dom, "$m$" left] (m) {$\lambda^\sharp_a(f,\lambda^\sharp_{a.f}(f',x))$\at$a$};
            \node[poly, right=3 of m.south, yshift=-1ex, "$\car{c}$" below] (D) {$f$\at$|a|$};
            \node[poly, above=of D, "$m$" above] (mm) {$\lambda^\sharp_{a.f}(f',x)$\at$a.f$};
            \node[poly, cod, right=2 of D.south, yshift=-1ex, "$\car{c}$" right] (DD) {$f$\at$|a|$};
            \node[poly, cod, above=of DD, "$\car{c}$" right] (mmm) {$f'$\at$|a.f|$};
            \node[poly, cod, above=of mmm, "$m$" right] (C) {$x$\at$(a.f).f'$};
            \draw (m_pos) to[first] (D_pos);
            \draw (D_dir) to[climb] node[left] {$\lambda$} (mm_pos);
            \draw (mm_dir) to[last] (m_dir);
            \draw[double,-] (D_pos) to[first] (DD_pos);
            \draw[double,-] (DD_dir) to[last] (D_dir);
            \draw (mm_pos) to[first] (mmm_pos);
            \draw (mmm_dir) to[climb] node[left] {$\lambda$} (C_pos);
            \draw (C_dir) to[last] (mm_dir);
        \end{tikzpicture}
    };
	\node at ($(p1.south)!.5!(p2.north)$) {$=$};
\end{tikzpicture}
\]
These say that $\cod(f)=|a.f|$, that $a.(f\then f')=(a.f).f'$, and that $\lambda^\sharp_a((f\then f'),x)=\lambda^\sharp_a(f,\lambda^\sharp_{a.f}(f',x))$.
	\item Here is the first equation in polyboxes:
\[
\begin{tikzpicture}[polybox, mapstos]
  \node[poly, dom, "$m$" left] (l) {$x$\at$a\vphantom{a}$};
  \node[poly, cod, "$m$" right, right=of l] (r) {$x$\at$a$};
  \draw[double,-] (l_pos) -- (r_pos);
  \draw[double,-] (r_dir) -- (l_dir);
	\node[poly, dom, right=2 of r, "$m$" left] (p) {$\rho^\sharp_a(x,\id_{|a'.x|})$\at$a$};
	\node[poly, cod, right=2.5cm of p.south, yshift=-1.25ex, "$m$" right] (r') {$x$\at$a'$};
	\node[poly, above=of r'] (q) {$\id_{|a'.x|}$\at$|a'.x|$};
  	\draw (p_pos) to[first] (r'_pos);
  	\draw (r'_dir) to[climb] node[left] {$\rho$} (q_pos);
		\draw (q_dir) to[last] (p_dir);
		\draw (q_pos) to[climb'] node[right] {$\epsilon_\car{d}$} (q_dir);
	\node at ($(r.east)!.5!(p.west)$) {$=$};
\end{tikzpicture}
\]
The equation says that for any $a\in m(1)$ we have $a=a'\coloneqq\rho_1(i)$ and that for any $f\in m[a]$ we have $x=\rho^\sharp_a(x,\id_{|a'.x|})$.

Here is the second equation in polyboxes:
\[
\begin{tikzpicture}
    \node (p1) {
        \begin{tikzpicture}[polybox, mapstos]
            \node[poly, dom, "$m$" left] (m') {$\rho^\sharp_a(\rho_a^\sharp(x,g),g')$\at$a$};
            \node[poly, right=2.7 of m'.south, yshift=-1ex, "$m$" below] (mm') {$\rho^\sharp_i(x,g)$\at$a$};
            \node[poly, above=of mm', "$\car{d}$" above] (C') {$g'$\at$|a.\rho^\sharp_a(x,g)|$};
            \node[poly, cod, right=2 of mm'.south, yshift=-1ex, "$m$" right] (D') {$x$\at$a$};
            \node[poly, cod, above=of D', "$\car{d}$" right] (mmm') {$g$\at$|a.x|$};
            \node[poly, cod, above=of mmm', "$\car{d}$" right] (CC') {$g'$\at$|a.\rho^\sharp_a(x,g)|$};
            \draw[double,-] (m'_pos) to[first] (mm'_pos);
            \draw (mm'_dir) to[climb] node[left] {$\rho$}  (C'_pos);
            \draw (C'_dir) to[last]  (m'_dir);
            \draw[double,-] (mm'_pos) to[first] (D'_pos);
            \draw (D'_dir) to[climb] node[left] {$\rho$} (mmm'_pos);
            \draw (mmm'_dir) to[last] (mm'_dir);
            \draw[double,-] (C'_pos) to[first] (CC'_pos);
            \draw[double,-] (CC'_dir) to[last] (C'_dir);
        \end{tikzpicture}
	};
	\node (p2) [below=1 of p1] {
	    \begin{tikzpicture}[polybox, tos]
            \node[poly, dom, "$m$" left] (m) {$\rho^\sharp_a(x,(g\then g'))$\at$a$};
            \node[poly, right=2.5 of m.south, yshift=-1ex, "$m$" below] (D) {$x$\at$a$};
            \node[poly, above=of D, "$\car{d}$" above] (mm) {$g\then g'$\at$|a.x|$};
            \node[poly, cod, right=2 of D.south, yshift=-1ex, "$m$" right] (DD) {$x$\at$a$};
            \node[poly, cod, above=of DD, "$\car{d}$" right] (mmm) {$g$\at$|a.x|$};
            \node[poly, cod, above=of mmm, "$\car{d}$" right] (C) {$g'$\at$\cod g$};
            \draw[double,-] (m_pos) to[first] (D_pos);
            \draw (D_dir) to[climb]  node[left] {$\rho$} (mm_pos);
            \draw (mm_dir) to[last] (m_dir);
            \draw[double,-] (D_pos) to[first] (DD_pos);
            \draw[double,-] (DD_dir) to[last] (D_dir);
            \draw[double,-] (mm_pos) to[first] (mmm_pos);
            \draw (mmm_dir) to[climb]  node[left] {$\delta$} (C_pos);
            \draw (C_dir) to[last] (mm_dir);
        \end{tikzpicture}
    };
	\node at ($(p1.south)!.5!(p2.north)$) {$=$};
\end{tikzpicture}
\]
These equations say that $\cod(g)=|a.\rho^\sharp_a(x,g)|$ and that $\rho_a^\sharp(\rho_a^\sharp(x,g),g')=\rho^\sharp_a(x,(g\then g'))$.
\end{enumerate}
\end{solution}
\end{exercise}

\begin{exercise}\label{exc.unique_yon_comodule}
Recall from \cref{exc.yon_comonoid} that $\yon$ has a unique category structure.
\begin{enumerate}
	\item Show that the category of left $\yon$-comodules is isomorphic to $\poly$.
	\item Show that the category of right $\yon$-comodules is isomorphic to $\poly$. If it is similar, just say ``similar''; if not, explain.
\end{enumerate}
\begin{solution}
The comonoid structure on $\yon$ has $\epsilon\colon\yon\to\yon$ and $\delta\colon\yon\to\yon\tri\yon=\yon$ the unique lenses.
\begin{enumerate}
	\item We first show that for any polynomial $m$, there is a unique $\yon$-comodule structure on $m$. By \eqref{eqn.left_comod_laws}, a $\yon$-comodule structure is a lens $\lambda\colon m\to\yon\tri m=m$ such that $\lambda\then(!\tri m)=\id_m$, but $!\tri m=\id_m$, so $\lambda=\id_m$ is forced. A morphism $(m,\id)\to(n,\id_n)$ of $\yon$-comodules is also easily seen to be just a lens $m\to n$.
	\item Similar.
\end{enumerate}
\end{solution}
\end{exercise}

We can characterize left comodules as something far more familiar.

\begin{proposition} \label{prop.left_comod}
Let $\cat{C}$ be a category.
The category of left $\cat{C}$-comodules is equivalent to the category of functors $\cat{C}\to\poly$.
\end{proposition}
\begin{proof}[Sketch of proof]%CHECK
Let $\car{c}$ be the carrier of the category $\cat{C}$. Given a left $\cat{C}$-comodule $m\to \car{c}\tri m$, we saw in \cref{exc.left_right_comodules} that there is an induced function $|-|\colon m(1)\to\car{c}(\1)$, so to each object $i\in\Ob(\cat{C})$ we can associate a set $P^m_i\coloneqq\{a\in m(\1)\mid |a|=i\}$. We also have for each $a\in P_i\ss m(\1)$ a set $m[a]$, and so we can consider the polynomial
\[
p^m_i\coloneqq\sum_{a\in P^m_i}\yon^{m[a]}.
\]
The polynomial $p^m_i\in\poly$ is easily seen to be functorial in the comodule $m$. Moreover, given a morphism $f\colon i\to i'$ in $\cat{C}$, we have for each $a\in p^m_i(1)$ an element $a.f\in m(1)$, with $|a.f|=i'$, i.e.\ $a.f\in  P^m_{i'}=p^m_{i'}(\1)$. Similarly for each $x\in m[a.f]=p^m_{i'}[a.f]$ we have $\lambda_a^\sharp(f,x)\in m[a]=p^m_i[a]$, and thus we have a lens $p^m_i\to p^m_{i'}$. In this way we have obtained we obtain a functor $p^m_-\colon\cat{C}\to\poly$, and this functor is itself functorial in $m$.

We leave the proof that this construction has an inverse to the reader.
\end{proof}

Similarly, a right comodule is not an altogether foreign concept.

\begin{proposition} \label{prop.right_comod}
Let $\cat{D}$ be a category.
A right $\cat{D}$-comodule $m$ can be identified (up to isomorphism) with a functor $\cat{D}\to\smset^{m(\1)}$.
\end{proposition}
\begin{proof}[Sketch of proof]
Let $\car{d}$ be the carrier of $\cat{D}$. Given a right $\cat{D}$-comodule $m\to m\tri\car{d}$ we need a functor $F\colon\cat{D}\times m(1)\to\smset$. We again rely heavily on \cref{exc.left_right_comodules}.
Given $j\in\cat{D}$ and $a\in m(1)$, define
\[F(j,d)\coloneqq\{x\in m[a]\mid|a.f|=j\}.\]
Given $g\colon j\to j'$ we get $\rho_a^\sharp(x,g)\in m[a]$ with $|a.\rho_a^\sharp(x,g)|=j'$, and thus $g$ induces a function $F(g,d)\colon F(j,d)\to F(j',d)$. Again by laws from \cref{exc.left_right_comodules} these are functorial in $g$, completing the proof sketch.
\end{proof}

\begin{proposition}\label{prop.all_free_modules}
Let $\com{C}=(\car{c},\epsilon,\delta)$ be a comonoid in $\poly$. For any set $G$, the polynomial $\yon^G\tri\car{c}$ has a natural right $\com{C}$-comodule structure.
\end{proposition}
\begin{proof}
We use the lens $(\yon^G\tri\delta)\colon(\yon^G\tri\car{c})\to(\yon^G\tri\car{c}\tri\car{c})$. It satisfies the unitality and associativity laws because $\car{c}$ does.\index{associativity}
\end{proof}

We can think of elements of $G$ as ``generators''. Then if $i'\colon G\to\car{c}\tri\1$ assigns to every generator an object of a category $\cat{C}$, then we should be able to find the free $\cat{C}$-set that $i'$ generates.

\begin{proposition}
Functions $i'\colon G\to\car{c}\tri\1$ are in bijection with positions $i\in\yon^G\tri\car{c}\tri\1$. Let $m\coloneqq i^*(\yon^G\tri\car{c})$ and let $\rho_i$ be the induced right $\com{C}$-comodule structure. %from \cref{prop.right_modules_as_sums}.
Then $\rho_i$ corresponds to the free $\cat{C}$-set generated by $i'$.
\end{proposition}
\begin{proof}
The polynomial $m$ has the following form:
\[
m\cong\yon^{\sum_{g\in G}\car{c}[i'(g)]}
\]
In particular $\rho_i$ is a representable right $\com{C}$-comodule, and we can identify it with a $\cat{C}$-set by \cref{thm.tfae_c_sets}. The elements of this $\cat{C}$-set are pairs $(g, f)$, where $g\in G$ is a generator and $f\colon i'(g)\to\cod(f)$ is a morphism in $\cat{C}$ emanating from $i'(g)$. It is easy to see that the comodule structure induced by \cref{prop.all_free_modules} is indeed the free one.
\end{proof}

\begin{exercise}
Let $\cat{C}$ be a category and $i\in\cat{C}$ an object.
\begin{enumerate}
    \item Consider $i$ as a lens $\yon\to\car{c}$.
    Show that the vertical-cartesian factorization of this lens is $\yon \to \yon^{\car{c}[i]} \To{\varphi} \car{c}$.
    \item Use \cref{prop.comp_pres_cart} to show that $\yon^{\car{c}[i]} \tri \car{c} \to \car{c} \tri \car{c}$ is cartesian.
    \item Show that there is a commutative square
    \[
    \begin{tikzcd}
        \yon^{\car{c}[i]} \ar[r, "\delta^i"] \ar[d, "\varphi"'] & \yon^{\car{c}[i]} \tri \car{c} \ar[d, "\text{cart}"] \\
        \car{c} \ar[r, "\delta"'] & \car{c} \tri \car{c} \ar[ul, phantom, very near end, "\lrcorner"]
    \end{tikzcd}
    \]
    \item Show that this square is a pullback, as indicated.
    \item Show that $\delta^i$ makes $\yon^{\car{c}[i]}$ a right $\com{C}$-comodule.
    \qedhere
\end{enumerate}
\end{exercise}\index{polynomial functors!pullback of polynomials}

The lens $\delta^i$ can be seen as the restriction of $\delta \colon \car{c} \to \car{c} \tri \car{c}$ to a single starting position.

We can extend this to a functor $\cat{C} \to \bimod{\yon}{\cat{C}}$ that sends the object $i$ to $\yon^{\car{c}[i]}$. Given a morphism $f \colon i \to i'$ in $\cat{C}$, we get a function $\car{c}[i'] \to \car{c}[i]$ given by composition with $f$, and hence a lens $\yon^{\car{c}[f]} \colon \yon^{\car{c}[i]} \to \yon^{\car{c}[i']}$.

\begin{exercise}
\begin{enumerate}
    \item Show that $\yon^{\car{c}[f]}$ is a right $\com{C}$-comodules morphism.
    \item Show that the construction $\yon^{\car{c}[f]}$ is functorial in $f$. \qedhere
\end{enumerate}
\end{exercise}

\begin{proposition}\label{prop.break_up_right_mods}
Let $\car{c}$ be a comonoid. For any set $I$ and right $\car{c}$-comodules $(m_i)_{i\in I}$, the coproduct $m\coloneqq \sum_{i\in I}m_i$ has a natural right-comodule structure. Moreover, each representable summand in the carrier $m$ of a right $\car{c}$-comodule is itself a right-$\car{c}$ comodule and $m$ is their sum.
\end{proposition}
\begin{proof}[Sketch of proof]
This follows from the fact in \cref{prop.left_dist_prod} that $-\tri\car{c}$ commutes with coproduct, i.e.\ $\left(\sum_{i\in I}m_i\right)\tri\car{c}\cong\sum_{i\in I}(m_i\tri\car{c})$.
\end{proof}

\begin{proposition}
If $m\in\poly$ is equipped with both a right $\com{C}$-comodule and a right $\com{D}$-comodule structure, we can naturally equip $m$ with a $(\com{C}\times\com{D})$-comodule structure.
\end{proposition}
\begin{proof}
It suffices by \cref{prop.break_up_right_mods} to assume that $m=\yon^M$ is representable. But a right $\com{C}$-comodule with carrier $\yon^M$ can be identified with a retrofunctor $M\yon^M\to\com{C}$.

Thus if $\yon^M$ is both a right-$\com{C}$ comodule and a right-$\com{D}$ comodule, then we have comonoid morphisms $\com{C}\from M\yon^M\to\com{D}$. This induces a unique comonoid morphism $M\yon^M\to(\com{C}\times\com{D})$ to the product, and we identify it with a right-$(\com{C}\times\com{D})$ comodule on $\yon^M$.
\end{proof}

\subsection{Bicomodules}
\index{bicomodule}

We take special note of any object that is both a left comodule and a right comodule in compatible ways.

\begin{definition}[Bicomodule]\label{def.bicomodule}
In a monoidal category $(\Cat{C},\yon,\tri)$, let $\cat{C}$ and $\cat{D}$ be comonoids with carriers $\car{c}$ and $\car{d}$.
A \emph{$(\cat{C},\cat{D})$-bicomodule} is
\begin{enumerate}
	\item an object $m\in\Cat{C}$ that is both
	\item a left $\cat{C}$-comodule, with left coaction $\car{c}\tri m\From{\lambda}m$, and
	\item a right $\cat{D}$-comodule, with right coaction $m\To{\rho}m\tri\car{d}$,
\end{enumerate}
such that the following diagram, known as the \emph{coherence law}, commutes:
\begin{equation} \label{eqn.bimod_coherence}
\begin{tikzcd}[sep=small]
    & m \ar[dl, "\lambda"'] \ar[dr, "\rho"] & \\
    \car{c}\tri m \ar[dr, "\car{c}\:\tri\:\rho"'] & &
    m\tri\car{d} \ar[dl, "\lambda\:\tri\:\car{d}"] \\
    & \car{c}\tri m\tri\car{d}.
\end{tikzcd}
\end{equation}
We may denote such a bicomodule by $\cat{C}\bimodfrom[m]\cat{D}$ when its left and right coactions have yet to be specified or may be inferred from context---or even by $\car{c}\bimodfrom[m]\car{d}$ when the comonoid structures on $\car{c}$ and $\car{d}$ can be inferred from context as well.%
\tablefootnote{The notation $\car{c}\bimodfrom\car{d}$ has a mnemonic advantage, as each $\tri$ goes in the correct direction: $\car{c}\tri m\from m$ and $m\to m\tri\car{d}$.
But it also looks like an arrow going backward from $\car{d}$ to $\car{c}$, which will turn out to have a semantic advantage as well.%; see \cref{**}.
}

A \emph{morphism} of $(\cat{C},\cat{D})$-bicomodules is one that is a morphism of left $\cat{C}$-comodules and a morphism of right $\cat{D}$-comodules.
\end{definition}
\index{bicomodule!morphism of}\index{polybox}

We draw the commutativity of \eqref{eqn.bimod_coherence} using polyboxes.
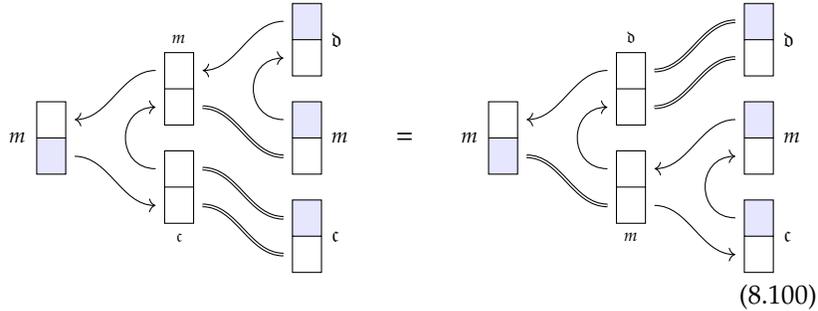
\begin{equation}\label{eqn.bimod_coherence_polybox}
\begin{tikzpicture}[baseline=(current bounding box.center)]
	\node (p1) {
  \begin{tikzpicture}[polybox, tos]
  	\node[poly, dom, "$m$" left] (m) {};
  	\node[poly, right= of m.south, yshift=-1ex, "\tiny$\car{c}$" below] (D) {};
  	\node[poly, above=of D, "\tiny$m$" above] (mm) {};
  	\node[poly, cod, right= of D.south, yshift=-1ex, "$\car{c}$" right] (DD) {};
  	\node[poly, cod, above=of DD, "$m$" right] (mmm) {};
  	\node[poly, cod, above=of mmm, "$\car{d}$" right] (C) {};
		\draw (m_pos) to[first] (D_pos);
		\draw (D_dir) to[climb] (mm_pos);
		\draw (mm_dir) to[last] (m_dir);
		\draw[double,-] (D_pos) to[first] (DD_pos);
		\draw[double,-] (DD_dir) to[last] (D_dir);
		\draw[double,-] (mm_pos) to[first] (mmm_pos);
		\draw (mmm_dir) to[climb] (C_pos);
		\draw (C_dir) to[last] (mm_dir);
	\end{tikzpicture}
	};
	\node (p2) [right=of p1] {
  \begin{tikzpicture}[polybox, tos]
  	\node[poly, dom, "$m$" left] (m') {};
  	\node[poly, right= of m'.south, yshift=-1ex, "\tiny$m$" below] (mm') {};
  	\node[poly, above=of mm', "\tiny$\car{d}$" above] (C') {};
  	\node[poly, cod, right= of mm'.south, yshift=-1ex, "$\car{c}$" right] (D') {};
  	\node[poly, cod, above=of D', "$m$" right] (mmm') {};
  	\node[poly, cod, above=of mmm, "$\car{d}$" right] (CC') {};
		\draw[double,-] (m'_pos) to[first] (mm'_pos);
		\draw (mm'_dir) to[climb] (C'_pos);
		\draw (C'_dir) to[last] (m'_dir);
		\draw (mm'_pos) to[first] (D'_pos);
		\draw (D'_dir) to[climb] (mmm'_pos);
		\draw (mmm'_dir) to[last] (mm'_dir);
		\draw[double,-] (C'_pos) to[first] (CC'_pos);
		\draw[double,-] (CC'_dir) to[last] (C'_dir);
	\end{tikzpicture}
	};
	\node at ($(p1.east)!.5!(p2.west)$) {$=$};
\end{tikzpicture}
\end{equation}
This polybox equation implies that we can unambiguously write
\[
\begin{tikzpicture}[polybox, tos]
	\node[poly, dom, "$m$" left] (m) {};
	\node[poly, cod, right=of m, "$m$" right] (mm) {};
	\node[poly, cod, above=of mm, "$\car{d}$" right] (C) {};
	\node[poly, cod, below=of mm, "$\car{c}$" right] (D) {};
	\draw (m_pos) to[out=0, in=180] (D_pos);
	\draw (D_dir) to[climb] (mm_pos);
	\draw (mm_dir) to[climb] (C_pos);
	\draw (C_dir) to[last] (m_dir);
\end{tikzpicture}
\]
for any bicomodule $\car{c}\bimodfrom[m]\car{d}$.

\begin{exercise}
Let $\cat{C}=(\car{c},\epsilon,\delta)$ be a category. Recall from \cref{exc.yon_comonoid} that $\yon$ has a unique category structure.
\begin{enumerate}
	\item Show that a left $\cat{C}$-module is the same thing as a $(\cat{C},\yon)$-bimodule.
	\item Show that a right $\cat{C}$-module is the same thing as a $(\yon,\cat{C})$-bimodule.
	\item Show that every polynomial $p\in\poly$ has a unique $(\yon,\yon)$-bimodule structure.
	\item Show there is an isomorphism of categories $\poly\cong\bimod{\yon}{\yon}$.
\qedhere
\end{enumerate}
\begin{solution}
\begin{enumerate}
    \item By \cref{exc.unique_yon_comodule}, it suffices to check that the diagram \eqref{eqn.bimod_coherence} commutes vacuously when $\car{d}=\yon$.
    \item By \cref{exc.unique_yon_comodule}, it suffices to check that the diagram \eqref{eqn.bimod_coherence} commutes vacuously when $\car{c}=\yon$.
    \item This follows from \cref{exc.unique_yon_comodule}.
    \item This follows from \cref{exc.unique_yon_comodule}.
\end{enumerate}
\end{solution}
\end{exercise}

\subsection{More equivalences}\index{retrofunctor!from state category}\index{coalgebra}\index{discrete opfibration}

We have seen that a retrofunctor from a state category to a category $\cat{C}$ carries the same data as a $\cat{C}$-coalgebra, which is in turn equivalent to a discrete opfibration over $\cat{C}$ or a copresheaf on $\cat{C}$.\index{copresheaf}
We then showed that cartesian retrofunctors to $\cat{C}$ is an equivalent notion as well. % TODO: refs?
We are now ready to show that several kinds of comodules are also equivalent to all of these concepts.

\begin{theorem}\label{thm.tfae_c_sets}
Given a polynomial comonoid $\cat{C}=(\car{c},\epsilon,\delta)$, the following comprise equivalent categories:
\begin{enumerate}
	\item functors $\cat{C}\to\smset$;
	\item discrete opfibrations over $\cat{C}$;
	\item cartesian retrofunctors to $\com{C}$;
	\item $\com{C}$-coalgebras (sets with a coaction by $\com{C}$);
	\item constant left $\com{C}$-comodules;
	\item $(\com{C},\0)$-bicomodules;
	\item linear left $\com{C}$-comodules; and
	\item representable right $\com{C}$-comodules (opposite).
\end{enumerate}
In fact, all but the first comprise isomorphic categories; and up to isomorphism, any one of these can be identified with a retrofunctor from a state category to $\cat{C}$.
\end{theorem}
\begin{proof}
\begin{description}
	\item[$1\simeq 2\cong 3$:] This was shown in \cref{prop.tfae_dopf}.
	\item[$3\cong 4$:] Given a cartesian retrofunctor $(\pi_1,\pi^\sharp)\colon\cat{S}\cof\cat{C}$, let $S\coloneqq\Ob(\cat{S})$ and define a $\car{c}$-coalgebra structure $\alpha\colon S\to\car{c}\tri S$ on an object $s\in\Ob(\cat{S})$ and an emanating morphism $f\colon\pi_1(s)\to c'$ in $\cat{C}$ by
	
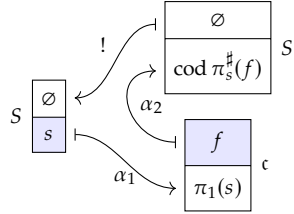
\[
  \begin{tikzpicture}[polybox, mapstos]
  	\node[poly, dom, "$S$" left] (yX) {$\varnothing$\at$s$};
  	\node[poly, cod, right=1.8 of yX.south, yshift=-1ex, "$\car{c}$" right] (p) {$f$\at$\pi_1(s)$};
  	\node[poly, above=of p, "$S$" right] (p') {$\varnothing$\at$\cod\pi^\sharp_s(f)$};
  	\draw (yX_pos) to[first] node[below] {$\alpha_1$}(p_pos);
  	\draw (p_dir) to[climb] node[right] {$\alpha_2$} (p'_pos);
  	\draw (p'_dir) to[last] node[above left] {$!$} (yX_dir);
  \end{tikzpicture}
	\]
	We check that this is indeed a coalgebra using properties of retrofunctors. For identities in $\com{C}$, we have
	\begin{align*}
  	\alpha_2(s,\id_{\pi_1(s)})&=
  	\cod\pi^\sharp_s(\id_{\pi_1(s)})\\&=
		\cod\id_{s}=s
	\end{align*}
	and for compositions in $\com{C}$ we have
	\begin{align*}
		\alpha_2(s, f\then g)&=
		\cod\left(\pi_s^\sharp(f\then g)\right)\\&=
		\cod\left(\pi_s^\sharp(f)\then\pi^\sharp_{\cod\pi^\sharp_s(f)}g\right)\\&=
		\cod\left(\pi^\sharp_{\cod\pi^\sharp_s(f)}g\right)\\&=
		\alpha_2(\alpha_2(s,f),g).
	\end{align*}
	Going backward, if we're given a coalgebra $\alpha\colon S\to\car{c}\tri S$, we obtain a function $\alpha_1\colon S\to\car{c}\tri\1$ and we define $\car{s}\coloneqq \alpha_1^*\car{c}$ and the cartesian lens $\varphi\coloneqq(\alpha_1,\id)\colon\car{s}\to\car{c}$ to be the base change from \cref{prop.basechange}. We need to show $\car{s}$ has a comonoid structure $(\car{s},\epsilon,\delta)$ and that $\varphi$ is a retrofunctor. We simply define the eraser $\epsilon\colon\car{s}\to\yon$ using $\alpha_1$ and the eraser on $\car{c}$:
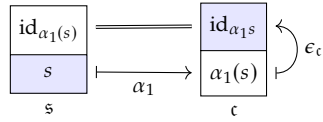
	\[
	\begin{tikzpicture}[polybox, mapstos]
  	\node[poly, dom, "$\car{s}$" below] (c) {$\id_{\alpha_1(s)}$\at$s\vphantom{s_1}$};
  	\node[poly, cod, right=of c, "$\car{c}$" below] (c') {$\id_{\alpha_1{s}}$\at$\alpha_1(s)$};
  	\draw (c_pos) -- node[below] {$\alpha_1$} (c'_pos);
  	\draw[double,-] (c'_dir) -- (c_dir);
		\draw (c'_pos) to[climb'] node[right] {$\epsilon_\car{c}$} (c'_dir);
	\end{tikzpicture}
	\]
 We give the duplicator $\delta\colon\car{s}\to\car{s}\tri\car{s}$ using $\alpha_2$ for the codomain, and using the composite $\then$ from $\car{c}$:
	\[
  \begin{tikzpicture}[polybox, mapstos]
  	\node[poly, dom, "$\car{s}$" left] (yX) {$f\then g$\at$s$};
  	\node[poly, cod, right=2.5 of yX.south, yshift=-1ex, "$\car{s}$" right] (p) {$f$\at$s$};
  	\node[poly, cod, above=of p, "$\car{s}$" right] (p') {$g$\at$\alpha_2(s,f)$};
  	\draw[double,-] (yX_pos) to[first] (p_pos);
  	\draw (p_dir) to[climb] node[right] {$\cod$} (p'_pos);
  	\draw (p'_dir) to[last] node[above left] {$\comp$} (yX_dir);
  \end{tikzpicture}
  \]
  In \cref{exc.tfae_c_sets} we check that $(\car{s},\epsilon,\delta)$ really is a comonoid, that $(\alpha_1,\id)\colon\car{s}\cof\car{c}$ is a retrofunctor, that the roundtrips between cartesian retrofunctors and coalgebras are identities, and that these assignments are functorial.
	\item[$4\cong5$:]This is straightforward and was mentioned in \cref{def.coalgebra}.
	\item[$5\cong6$:]A right $\0$-comodule is in particular a polynomial $m\in\poly$ and a lens $\rho\colon m\to m\tri\0$ such that $(m\tri\epsilon)\circ\rho=\id_m$. This implies $\rho$ is monic, which itself implies by \cref{prop.monics_in_poly} that $m$ must be constant since $m\tri\0$ is constant. This makes $m\tri\epsilon$ the identity, at which point $\rho$ must also be the identity. Conversely, for any set $M$, the corresponding constant polynomial is easily seen to make the diagrams in \eqref{eqn.right_comod_laws} commute.\index{constant polynomial}
	\item[$5\cong7$:] By the adjunction in \cref{thm.adjoint_quadruple} and the fully faithful inclusion $\smset\to\poly$ of sets as constant polynomials, \cref{prop.ff_const_set_to_poly}, we have isomorphisms
	\[\poly(S\yon,\car{c}\tri S\yon)\cong\smset(S,\car{c}\tri S\yon\tri\1)=\smset(S,\car{c}\tri S)\cong\poly(S,\car{c}\tri S).\]
	One checks easily that if $S\yon\to\car{c}\tri S\yon$ corresponds to $S\to\car{c}\tri S$ under the above isomorphism, then one is a left comodule if and only if the other is.
	\item[$7\cong8$:] By \eqref{eqn.flip_reps_lins} we have a natural isomorphism\index{isomorphism!natural}
	\[
		\poly(S\yon,\car{c}\tri S\yon)\cong\poly(\yon^S,\yon^S\tri\car{c}).
	\]
	In pictures,
	
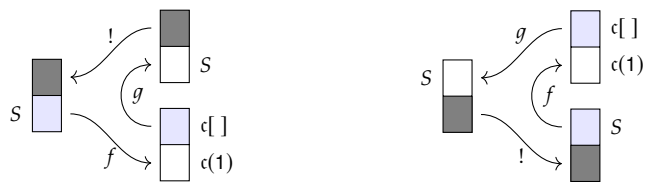
\[
	\begin{tikzpicture}
		\node (p1) {
		\begin{tikzpicture}[polybox, tos]
			\node[poly, linear dom] (s) {};
			\node[poly, cod, right=of s.south, yshift=-1ex] (c) {};
			\node[poly, linear cod, above=of c] (s') {};
    	\node[left=0pt of s_pos] {$S$};
    	\node[right=0pt of c_pos] {$\car{c}(\1)$};
    	\node[right=0pt of c_dir] {$\car{c}[\ ]$};
    	\node[right=0pt of s'_pos] {$S$};
			\draw (s_pos) to[first] node[below] {$f$} (c_pos);
			\draw (c_dir) to[climb] node[right] {$g$} (s'_pos);
			\draw (s'_dir) to[last] node[above] {$!$} (s_dir);
		\end{tikzpicture}
		};
		\node (p2) [right=2 of p1] {
		\begin{tikzpicture}[polybox, tos]
			\node[poly, pure dom] (s) {};
			\node[poly, pure cod, right=of s.south, yshift=-1ex] (s') {};
			\node[poly, cod, above=of s'] (c) {};
    	\node[left=0pt of s_dir] {$S$};
    	\node[right=0pt of c_pos] {$\car{c}(\1)$};
    	\node[right=0pt of c_dir] {$\car{c}[\ ]$};
    	\node[right=0pt of s'_dir] {$S$};
			\draw (s_pos) to[first] node[below] {$!$} (s'_pos);
			\draw (s'_dir) to[climb] node[right] {$f$} (c_pos);
			\draw (c_dir) to[last] node[above] {$g$} (s_dir);
		\end{tikzpicture}
		};
	\end{tikzpicture}
	\]
\end{description}
\noindent
The last claim was proven in \cref{prop.ds_dopf}.
\end{proof}

\begin{exercise}\label{exc.tfae_c_sets}
Complete the proof of \cref{thm.tfae_c_sets} ($3\cong 4$) by proving the following.
\begin{enumerate}
	\item Show that $(\car{s},\epsilon,\delta)$ really is a comonoid.
	\item Show that $(\alpha_1,\id)\colon\car{s}\cof\car{c}$ is a retrofunctor.
	\item Show that the roundtrips between cartesian retrofunctors and coalgebras are identities.
	\item Show that the assignment of a $\com{C}$-coalgebra to a cartesian retrofunctor over $\cat{C}$ is functorial.
	\item Show that the assignment of a cartesian retrofunctor over $\cat{C}$ to a $\com{C}$-coalgebra is functorial.
\qedhere
\end{enumerate}
\end{exercise}

Let $\cat{C}$ be a category. Under the above correspondence, the terminal functor $\cat{C}\to\smset$ corresponds to the identity discrete opfibration $\cat{C}\to\cat{C}$, the identity retrofunctor $\com{C}\to\com{C}$, a certain left $\com{C}$ comodule with carrier $\com{C}(\1)\yon$ which we call the \emph{canonical left $\com{C}$-comodule}, a certain constant left $\com{C}$ comodule with carrier $\com{C}(\1)$ which we call the \emph{canonical $(\com{C},\0)$-bicomodule}, and a certain representable right $\com{C}$-comodule with carrier $\yon^{\com{C}(\1)}$ which we call the \emph{canonical right $\cat{C}$-comodule}.

\begin{exercise}
For any object $c\in \cat{C}$, consider the representable functor $\cat{C}(c,-)\colon\cat{C}\to\smset$. What does it correspond to as a
\begin{enumerate}
	\item discrete opfibration over $\cat{C}$?
	\item cartesian retrofunctor to $\com{C}$?
	\item linear left $\com{C}$-comodule?
	\item constant left $\com{C}$-comodule?
	\item $(\com{C},\0)$-bicomodule?
	\item representable right $\com{C}$-comodule?
	\item dynamical system with comonoid interface $\com{C}$?
\qedhere
\end{enumerate}
\end{exercise}\index{interface!comonoid}

\index{copresheaf!as bicomodule}\index{category!of bicomodules}\index{category!copresheaf}

\begin{exercise}
We saw in \cref{thm.tfae_c_sets} that the category $\bimod{\cat{C}}{\0}$ of $(\cat{C},\0)$-bicomodule has a very nice structure: it's the topos of copresheaves on $\cat{C}$.
\begin{enumerate}
	\item What is a $(\0,\cat{C})$-bicomodule?
	\item What is $\bimod{\0}{\cat{C}}$?
\qedhere
\end{enumerate}
\end{exercise}

\subsection{Bicomodules are parametric right adjoints}
\index{parametric right adjoint}\index{bicomodule!as parametric right adjoint}

There is an equivalent characterization of bicomodules between a pair of polynomial comonoids due to Richard Garner.\index{Garner, Richard} Thus we attribute the foundational theory of $\catsharp$ to Ahman-Uustalu-Garner.

\begin{proposition}[Garner]\label{prop.prafunctor}
Let $\cat{C}$ and $\cat{D}$ be categories.
Then the following can be identified, up to isomorphism:
\begin{enumerate}
	\item a $(\cat{C},\cat{D})$-bicomodule.
% 	\item a profunctor + discrete opfibration
	\item a parametric right adjoint $\smset^{\cat{D}}\to\smset^{\cat{C}}$.
	\item a connected limit-preserving functor $\smset^{\cat{D}}\to\smset^{\cat{C}}$.
\end{enumerate}
\end{proposition}

Parametric right adjoints model \emph{data migrations} between categorical databases; see \cite{Spivak.Wisnesky:2015a}. \index{database!data migration}

When a polynomial
\[
m\coloneqq\sum_{i\in m(\1)}\yon^{m[i]}
\]
is given the structure of a $(\cat{D},\cat{C})$-bicomodule, the symbols in that formula are given a hidden special meaning:
\[
  m(\1)\in\smset^{\cat{D}}
  \qqand
	m[i]\in\smset^{\cat{C}}
\]
% Thus $m(1)$ is a copresheaf on $\cat{D}$; in particular, each position in $i\in m(\1)$ is **. And each $m[i]$ is a database instance on $\cat{C}$; in particular, each direction $d\in m[i]$ is a row in that instance.

Before we knew about bicomodule structures, what we called positions and directions were understood as each forming an ordinary set. In the presence of a bicomodule structure, the positions $m(\1)$ have been organized into a $\cat{D}$-set and the directions $m[i]$ have been organized into a $\cat{C}$-set for each position $i$. We are listening for $\cat{C}$-sets and positioning ourselves in a $\cat{D}$-set.

\begin{definition}[Prafunctor]
Let $\cat{C}$ and $\cat{D}$ be categories. A \emph{prafunctor} (also called a \emph{parametric right adjoint functor}) $\smset^{\cat{C}}\to\smset^{\cat{D}}$ is one satisfying any of the conditions of \cref{prop.prafunctor}.
\end{definition}

\subsection{Bicomodules in dynamics}\index{bicomodule!dynamical systems and bicomodules}\index{dynamical system!as bicomodule}

We conclude with a peek at how bicomodules can model dynamical systems themselves.

\index{cellular automata}

\begin{example}[Cellular automata]\label{ex.cellular_bimod}
In \cref{ex.graph_interaction,exc.conway} we briefly discussed cellular automata; here we will discuss another way that cellular automata show up, this time in terms of bicomodules.

Suppose that $\src,\tgt\colon A\tto V$ is a graph, and consider the polynomial
\[
  g\coloneqq\sum_{v\in V}\yon^{\src\inv(v)}
\]
so that positions are vertices and directions are emanating arrows. It carries a natural bicomodule structure
\[
V\yon\bimodfrom[g]V\yon
\]
where the right structure lens uses $\tgt$; see \cref{exc.cellular_bimod} for details. A bicomodule
\[
V\yon\bimodfrom[T]\0
\]
can be identified with a functor $T\colon V\to\smset$, i.e.\ it assigns to each vertex $v\in V$ a set. Let's call $T(v)$ the color set at $v$; for many cellular automata we will put $T(v)\cong\2$ for each $v$.

Then a cellular automata on $g$ with color sets $T$ is given by a map
\[
\begin{tikzcd}
	V\yon\ar[r, bimr-biml, "g"]\ar[rr, bend right=40pt, bimr-biml, "T"', "" name=T]&
	V\yon\ar[r, bimr-biml, "T"]\ar[to=T, Rightarrow, "\alpha"]&
	\0
\end{tikzcd}
\]
Indeed, for every vertex $v\in V$ the map $\alpha$ gives a function
\[
\prod_{\src(a)=v}T(\tgt(a))\Too{\alpha_v} T(v),
\]
which we call the \emph{update} function. In other words, given the current color at the target of each arrow emanating from $v$, the function $\alpha_v$ returns a new color at $v$.

Note that if $V\in\bimod{V\yon}{\0}$ is the terminal object, then the composite $V\yon\bimodfrom[g]V\yon\bimodfrom[V]\0$ is again $V$.
\end{example}

\begin{exercise}\label{exc.cellular_bimod}
Let $\src,\tgt\colon A\tto V$ and $g$ and $T$ be as in \cref{ex.cellular_bimod}.
\begin{enumerate}
	\item Give the structure lens $\lambda\colon g\to V\yon\tri g$
	\item Give the structure lens $\rho\colon g\to g\tri V\yon$.
	\item Give the set $T$ and the structure lens $T\to V\yon\tri T$ corresponding to the functor $V\to\smset$ that assigns $\2$ to each vertex.
\qedhere
\end{enumerate}
\begin{solution}
For each $v\in V$, let $A_v\coloneqq\fun{src}\inv(v)$, the set of arrows emanating from $v$.
\begin{enumerate}
    \item We need a lens $\sum_{v\in V}\yon^{A_v}\to V\sum_{v\in V}\yon^{A_v}$, and we use the obvious ``diagonal'' lens, which on positions sends $v\mapsto (v,v)$ and on directions is the identity.
    \item We need a lens $\sum_{v\in V}\yon^{A_v}\to\sum_{v\in V}\prod_{a\in A_v}\sum_{v\in V}\yon$, which is the same as an element of
    \[
    \prod_{v\in V}\sum_{v'\in V}\prod_{a\in A_{v'}}\sum_{v''\in V}A_v
    \]
    and we use $v\mapsto (v, a\mapsto\fun{tgt}(a), a)$. Most of this is forced on us, and the only interesting part is the function $A_v\to V$, which we can take to be anything, the choice being exactly the choice of ``target map'' that defines our graph.
    \item Let $T=\2V$ and take the map $\2V\to V\times \2V$ to be $(i,v)\mapsto(v,i,v)$ for any $i\in \2$ and $v\in V$.
\end{enumerate}
\end{solution}
\end{exercise}

\index{cellular automata!running}

\begin{example}[Running a cellular automaton]\index{graph!cellular automaton on}
Let $g$ be a graph on vertex set $V$, let $T$ assign a color set to each $v\in V$, and let $\alpha$ be the update function for a cellular automaton. As in \cref{ex.cellular_bimod}, this is all given by a diagram
\[
\begin{tikzcd}
	V\yon\ar[r, bimr-biml, "g"]\ar[rr, bend right=40pt, bimr-biml, "T"', "" name=T]&
	V\yon\ar[r, bimr-biml, "T"]\ar[to=T, Rightarrow, "\alpha"]&
	\0
\end{tikzcd}
\]
To run the cellular automaton, one simply chooses a starting color in each vertex. We call this an initialization; it is given by a map of bicomodules
\begin{equation}\label{eqn.starting_color_v}
\begin{tikzcd}
	V\yon
		\ar[r, bimr-biml, bend left, "V"]
		\ar[r, bimr-biml, bend right, "T"']\ar[r, phantom, "\hphantom{\scriptstyle\sigma}\Downarrow\scriptstyle\sigma"]&
	0\hphantom{\yon}
\end{tikzcd}
\end{equation}

Note that $g$ is a profunctor, i.e.\ for each $v\in V$ the summand $g_v=\yon^{\fun{src}\inv(v)}$ is representable, so it preserves the terminal object. In other words $g\tri_{V\yon}V\cong V$.

Now to run the cellular automaton on that initialization for $k\in\nn$ steps is given by the composite
\[
\begin{tikzcd}
	V\yon\ar[d, equal]\ar[r, bimr-biml, "g"]\ar[rrrr, bimr-biml, bend left, "V", ""' name=V']&
	\cdots\ar[r, bimr-biml, "g"]&|[alias=VV]|
	V\yon\ar[d, equal]\ar[r, bimr-biml, "g"]&
	V\yon\ar[d, equal]\ar[r, bimr-biml, "V"]&
	\0\ar[d, equal]\ar[dl, phantom, "\Downarrow\scriptstyle\sigma"]\\
	V\yon\ar[d, equal]\ar[r, bimr-biml, "g"]&
	\cdots\ar[r, bimr-biml, "g"]&
	V\yon\ar[d, equal]\ar[r, bimr-biml, "g"]&
	V\yon\ar[r, bimr-biml, "T"]&
	\0\ar[d, equal]\ar[dll, phantom, "\Downarrow\scriptstyle\alpha"]\\
	V\yon\ar[d, equal]\ar[r, bimr-biml, "g"]&
	\cdots\ar[r, bimr-biml, "g"]&|[alias=V]|
	V\yon\ar[rr, bimr-biml, "T"]&&
	\0\ar[d, equal]\\
	V\yon\ar[rrrr, bimr-biml, "T"', "" name=T]&&&&
	\0
	\ar[from=V, to=V|-T, Rightarrow, "\alpha\circ\cdots\circ\alpha"]
	\ar[from=V', to=VV, equal]
\end{tikzcd}
\]
\end{example}

\begin{exercise}
Explain why \eqref{eqn.starting_color_v} models an initialization, i.e.\ a way to choose a starting color in each vertex.
\begin{solution}
Recall that the $(V\yon,0)$-bicomodule $T$ assigns to each vertex $v\in V$ some set $T_v$ of ``colors'', whereas the $(V\yon,0)$-bicomodule $V$ assigns each vertex a singleton set. A map $V\to T$ between them chooses one color for each vertex, which we're calling the ``starting color''.
\end{solution}
\end{exercise}

\index{comodule|)}

%-------- Section --------%
\section[Summary and further reading]{Summary and further reading%
  \sectionmark{Summary \& further reading}}
\sectionmark{Summary \& further reading}

In this chapter we gave more of the theory of $\tri$-comonoids, sometimes called polynomial comonads. We began by defining an adjunction
\[
\adj{\catsharp}{U}{\cofree{-}}{\poly}
\]
where $\cofree{p}$ is the cofree comonoid $\cofree{p}$ associated to any polynomial $p$, and we gave intuition for it in terms of $p$-trees. A $p$-tree is a (possibly infinite) tree for which each vertex is labeled by a position $i\in p(\1)$ and its outgoing arrows are each labeled by an element of $p[i]$. The category corresponding to $\cofree{p}$ has $p$-trees as its objects; the morphisms emanating from such a $p$-tree are the finite rooted paths up the tree, and the codomain of such a path is the tree rooted at its endpoint.

We then briefly discussed some other properties of $\catsharp$ including a formal proof of the fact that it has all limits and colimits. Colimits in $\catsharp$ are created by, i.e.\ fit nicely with, colimits in $\poly$, but limits are quite strange; we did not really discuss them here, but for example the product in $\catsharp$ of the walking arrow \fbox{$\bullet\to\bullet$} with itself has infinitely-many objects!

We then moved on to left-, right-, and bicomodules between polynomial comonoids. In particular, we showed that left $\cat{C}$-comodules can be identified with functors $\cat{C}\to\poly$, and that $(\cat{C},\cat{D})$-bicomodules correspond to parametric right adjoint functors (prafunctors) $\smset^\cat{D}\to\smset^\cat{C}$. This idea was due to Richard Garner; it is currently unpublished, but can be found in video form here: \url{https://www.youtube.com/watch?v=tW6HYnqn6eI}. What we call prafunctors are sometimes called familial functors between (co-)presheaf categories; see \cite{weber2007theory}, \cite{garner2018shapely}, and \cite{shapiro2021familial} for more on this notion.

\index{polynomial comonoid|)}
\index{$\catsharp$!categorical properties of|)}

%-------- Section --------%
\section{Exercise solutions}
\Closesolutionfile{solutions}
{\footnotesize
\input{solution-file8}}

\setcounter{chapter}{8}%Just finished 8.

%------------ Chapter ------------%
\chapter{New horizons} \label{ch.comon.vistas}

In this brief chapter, we lay out some questions that whose answers may or may not be known, but which were not known to us at the time of writing. They vary from concrete to open-ended, they are not organized in any particular way, and are in no sense complete. Still we hope they may be useful to some readers.

\begin{enumerate}
  \item What can you say about comonoids in the category of all functors $\smset\to\smset$, e.g. ones that aren't polynomial.
  \item What can you say about the internal logic for the topos $[\cofree{p},\smset]$ of dynamical systems with interface $p$, in terms of $p$?
  \item How does the logic of the topos $\cofree{p}$ help us talk about issues that might be useful in studying dynamical systems?
  \item Morphisms $p\to q$ in $\poly$ give rise to left adjoints $\cofree{p}\to\cofree{q}$ that preserve connected limits. These are not geometric morphisms in general; in some sense they are worse and in some sense they are better. They are worse in that they do not preserve the terminal object, but they are better in that they preserve every connected limit not just finite ones. How do these left adjoints translate statements from the internal language of $p$ to that of $q$?
  \item Consider the $\times$-(co)monoids and $\otimes$-(co)monoids in three categories: $\poly$, $\smcat^\sharp$, and $\bimod{}{}$. Find examples of these (co)monoids, and perhaps characterize them or create a theory of them.
  \item The category $\poly$ has pullbacks, so one can consider the bicategory of spans in $\poly$. Is there a functor from that to $\bimod{}{}$ that sends $p\mapsto\cofree{p}$?
  \item Databases are static things, whereas dynamical systems are dynamic; yet we see them both in terms of $\poly$. How do they interact? Can a dynamical system read from or write to a database in any sense?
  \item Can we do database aggregation in a nice dynamic way?
  \item In the theory of polynomial functors, sums of representable functors $\smset\to\smset$, what happens if we replace sets with homotopy types: how much goes through? Is anything improved?
  \item Are there any functors $\smset\to\smset$ that aren't polynomial, but which admit a comonoid structure with respect to composition $(\yon,\tri)$?
  \item Characterize the monads in poly. They're generalizations of one-object operads (which are the Cartesian ones), but how can we think about them?
  \item Describe the limits that exist in $\smcat^\sharp$ combinatorially.
  \item Since the forgetful functor $U\colon\smcat^\sharp\to\poly$ is faithful, it reflects monomorphisms: if $f\colon\cat{C}\cof\cat{D}$ is a retrofunctor whose underlying map on carriers is monic, then it is monic. Are all monomorphisms in $\smcat^\sharp$ of this form?
  \item Are there polynomials $p$ such that one use something like G\"odel numbers to encode logical propositions from the topos $[\tr_p,\smset]$ into a ``language'' that $p$-dynamical systems can ``work with''?
\end{enumerate}

\appendix
\begingroup
\footnotesize

\backmatter

\printbibliography

@inproceedings{abbott2003derivatives,
  title={Derivatives of containers},
  author={Abbott, Michael and Altenkirch, Thorsten and Ghani, Neil and McBride, Conor},
  booktitle={International Conference on Typed Lambda Calculi and Applications},
  pages={16--30},
  year={2003},
  organization={Springer}
}

@article{abbott2005containers,
title = "Containers: Constructing strictly positive types",
journal = "Theoretical Computer Science",
volume = "342",
number = "1",
pages = "3 - 27",
year = "2005",
note = "Applied Semantics: Selected Topics",
author = "Michael Abbott and Thorsten Altenkirch and Neil Ghani",
}

@phdthesis{abbot2003categoriesthesis,
  author       = {Michael Gordon Abbott},
  title        = {Categories of Containers},
  school       = {University of Leicester},
  year         = 2003,
  month        = 8,
}

@incollection{abou2016reflections,
  title={Reflections on monadic lenses},
  author={Abou-Saleh, Faris and Cheney, James and Gibbons, Jeremy and McKinna, James and Stevens, Perdita},
  booktitle={A List of Successes that can Change the World},
  pages={1--31},
  year={2016},
  publisher={Springer}
}

@book{aguiar1997internal,
  title={Internal categories and quantum groups},
  author={Aguiar, Marcelo},
  year={1997},
  publisher={Cornell University}
}

@article{ahman2016directed,
Author = {Danel Ahman and Tarmo Uustalu},
Title = {Directed Containers as Categories},
Year = {2016},
Eprint = {arXiv:1604.01187},
journal = {EPTCS 207, 2016, pp. 89-98},
% Doi = {10.4204/EPTCS.207.5},
}

@article{beurier2019memoryless,
  title={Memoryless systems generate the class of all discrete systems},
  author={Beurier, Erwan and Pastor, Dominique and Spivak, David I},
  journal={International Journal of Mathematics and Mathematical Sciences},
  volume={2019},
  year={2019},
  publisher={Hindawi}
}

@inproceedings{bohannon2006relational,
  title={Relational lenses: a language for updatable views},
  author={Bohannon, Aaron and Pierce, Benjamin C and Vaughan, Jeffrey A},
  booktitle={Proceedings of the twenty-fifth ACM SIGMOD-SIGACT-SIGART symposium on Principles of database systems},
  pages={338--347},
  year={2006},
  organization={ACM}
}

@book {Borceux:1994a,
    AUTHOR = {Borceux, Francis},
     TITLE = {Handbook of categorical algebra 1},
    SERIES = {Encyclopedia of Mathematics and its Applications},
    VOLUME = {50},
      NOTE = {Basic category theory},
 PUBLISHER = {Cambridge University Press, Cambridge},
      YEAR = {1994},
}

@misc{capucci2022diegetic,
  url = {https://arxiv.org/abs/2206.12338},
  author = {Capucci, Matteo},
  title = {Diegetic representation of feedback in open games},
  publisher = {arXiv},
  year = {2022},
 }

@book{cheng_2022, place={Cambridge}, title={The Joy of Abstraction: An Exploration of Math, Category Theory, and Life}, DOI={10.1017/9781108769389}, publisher={Cambridge University Press}, author={Cheng, Eugenia}, year={2022}}

@book{conway2012regular,
  title={Regular algebra and finite machines},
  author={Conway, John Horton},
  year={2012},
  publisher={Courier Corporation}
}

@book{fong2019seven,
  title={Seven Sketches in Compositionality: An Invitation to Applied Category Theory},
  author={Fong, Brendan and Spivak, David I.},
  year={2019},
  publisher={Cambridge University Press},
}

@article{garner2018shapely,
  title={Shapely monads and analytic functors},
  author={Garner, Richard and Hirschowitz, Tom},
  journal={Journal of Logic and Computation},
  volume={28},
  number={1},
  pages={33--83},
  year={2018},
  publisher={Oxford University Press}
}

@article{gibbons2012relating,
  title={Relating algebraic and coalgebraic descriptions of lenses},
  author={Gibbons, Jeremy and Johnson, Michael},
  journal={Electronic Communications of the EASST},
  volume={49},
  year={2012}
}

@misc{hedges2018limits,
Author = {Jules Hedges},
Title = {Limits of bimorphic lenses},
Year = {2018},
Eprint = {arXiv:1808.05545},
}

@misc{hedges2016compositionality,
  title={Compositionality and string diagrams for game theory},
  author={Hedges, Jules and Shprits, Evguenia and Winschel, Viktor and Zahn, Philipp},
  eprint={arXiv:1604.06061},
  year={2016}
}

@misc{hedges2017coherence,
  title={Coherence for lenses and open games},
  author={Hedges, Jules},
  eprint={arXiv:1704.02230},
  year={2017}
}

@article{hedges2018morphisms,
  title={Morphisms of open games},
  author={Hedges, Jules},
  journal={Electronic Notes in Theoretical Computer Science},
  volume={341},
  pages={151--177},
  year={2018},
  publisher={Elsevier}
}

@book{jacobs2017introduction,
  title={Introduction to Coalgebra},
  author={Jacobs, Bart},
  volume={59},
  year={2017},
  publisher={Cambridge University Press}
}

@article{johnson2012lenses,
  title={Lenses, fibrations and universal translations},
  author={Johnson, Michael and Rosebrugh, Robert and Wood, Richard J},
  journal={Mathematical Structures in Computer Science},
  volume={22},
  number={1},
  pages={25--42},
  year={2012},
  publisher={Cambridge University Press}
}

@inproceedings{kelly1974doctrinal,
  title={Doctrinal adjunction},
  author={Kelly, G Max},
  booktitle={Category seminar},
  pages={257--280},
  year={1974},
  organization={Springer}
}

@article{kock2012polynomial,
   title={Polynomial functors and polynomial monads},
   volume={154},
   number={1},
   journal={Mathematical Proceedings of the Cambridge Philosophical Society},
   publisher={Cambridge University Press (CUP)},
   author={Gambino, Nicola and Kock, Joachim},
   year={2012},
   month={09},
   pages={153–192}
}

@unpublished{kock2016,
	 title={Notes on Polynomial Functors},
  author={Kock, Joachim},
  howpublished={\url{https://mat.uab.cat/~kock/cat/polynomial.pdf}},
  note={Accessed 2023/12/04}
 }

@book{Leinster:2014a,
  title={Basic category theory},
  author={Leinster, Tom},
  volume={143},
  year={2014},
  publisher={Cambridge University Press}
}

@book{MacLane:1998a,
  Author = {Mac Lane, Saunders},
  Title = {Categories for the working mathematician},
  Edition = {2},
  %Isbn = = {0-387-98403-8},
  Location = {New York},
  Publisher = {Springer-Verlag},
  Series = {Graduate Texts in Mathematics},
  Number = {5},
  Year = {1998}
}

@book{macLane1992sheaves,
author = {MacLane, Saunders and Moerdijk, Ieke},
%ISBN = = {0387977104},
publisher = {Springer},
title = {Sheaves in Geometry and Logic: A First Introduction to Topos Theory},
year = {1992},
}

@article{mcbride2001derivative,
  title={The derivative of a regular type is its type of one-hole contexts},
  author={McBride, Conor},
  journal={Unpublished manuscript},
  pages={74--88},
  year={2001},
  publisher={Citeseer}
}

@book{jaz,
   author = {David Jaz Myers},
   year = {2022},
   title = {Categorical Systems Theory}
}

@misc{nlab2022lens,
  author = {{nLab authors}},
  title = {lens (in computer science)},
  howpublished = {\url{http://ncatlab.org/nlab/show/lens\%20\%28in\%20computer\%20science\%29}},
  note = {\href{http://ncatlab.org/nlab/revision/lens\%20\%28in\%20computer\%20science\%29/26}{Revision 26}},
  month = aug,
  year = 2022
}

@misc{nlab:connected-limit,
    author = "nLab, Contributors To",
    title = "Connected limit --- nLab",
    year = "2019",
    url = "https://ncatlab.org/nlab/show/connected+limit",
}

@misc{nlab:created-limit,
    author = "nLab, Contributors To",
    title = "Created limit --- nLab",
    year = "2018",
    url = "https://ncatlab.org/nlab/show/created+limit",
}

@misc{oconnor2011functor,
  title={Functor is to lens as applicative is to biplate: Introducing multiplate},
  author={O'Connor, Russell},
  eprint={arXiv:1103.2841},
  year={2011}
}

@misc{pare2023retrocells,
      title={Retrocells}, 
      author={Robert Paré},
      year={2023},
      eprint={2306.06436},
      archivePrefix={arXiv},
      primaryClass={math.CT}
}

@incollection{pare1969absolute,
    author="Par{\'e}, Robert",
    editor="Hilton, Peter J.",
    title="Absolute coequalizers",
    booktitle="Category Theory, Homology Theory and their Applications I",
    year="1969",
    publisher="Springer",
    address="Berlin, Heidelberg",
    pages="132--145",
    isbn="978-3-540-36095-7"
}

@book{Pierce:1991,
  title = {Basic Category Theory for Computer Scientists},
  author = {Benjamin C. Pierce},
  year = {1991},
  publisher = {MIT Press},
}

@article{porst2019colimits,
    title={Colimits of monoids},
    author={Porst, Hans-E},
    journal={Theory and Applications of Categories},
    volume={34},
    number={17},
    pages={456--467},
    year={2019}
}

@book{Riehl:2017a,
  title={Category theory in context},
  author={Riehl, Emily},
  year={2017},
  publisher={Courier Dover Publications}
}

@misc{shapiro2021familial,
  url = {https://arxiv.org/abs/2111.14796},
  author = {Shapiro, Brandon},
  title = {Familial Monads as Higher Category Theories},
  publisher = {arXiv},
  year = {2021},
}

@article{shulman2008framed,
  Author = {Shulman, Michael},
  Date = {2008},
%  Eprint = {0706.1286},
%  Eprintclass = {math.CT},
%  Eprinttype = {arXiv},
  Journaltitle = {Theory and Applications of Categories},
  Pages = {Paper No. 18, 650-738},
  Title = {Framed bicategories and monoidal fibrations},
  Volume = {20}
}

@article {spivak2012functorial,
    AUTHOR = {Spivak, David I.},
     TITLE = {Functorial data migration},
%   JOURNAL = {Inform. and Comput.},
%  FJOURNAL = {Information and Computation},
  JOURNAL = {Information and Computation},
    VOLUME = {217},
      YEAR = {2012},
     PAGES = {31--51},
      %ISSN = {0890-5401},
   MRCLASS = {68P15 (18A40)},
  MRNUMBER = {2943977},
       % DOI = {10.1016/j.ic.2012.05.001},
%       URL = {https://doi.org/10.1016/j.ic.2012.05.001},
}

@inproceedings{Spivak.Wisnesky:2015a,
  author = {Spivak, David I. and Wisnesky, Ryan},
  title = {Relational Foundations for Functorial Data Migration},
  booktitle = {Proceedings of the 15th Symposium on Database Programming Languages},
  series = {DBPL},
  year = {2015},
  %ISBN = = {978-1-4503-3902-5},
  location = {Pittsburgh, PA},
  pages = {21-28},
  %doi = {10.1145/2815072.2815075},
  publisher = {ACM},
}

@article{spivak2017nesting,
  title={Nesting of dynamical systems and mode-dependent networks},
  author={Spivak, David I and Tan, Joshua},
  journal={Journal of Complex Networks},
  volume={5},
  number={3},
  pages={389--408},
  year={2017},
  publisher={Oxford University Press}
}

@misc{spivak2019generalized,
Author = {David I. Spivak},
Title = {Generalized Lens Categories via functors $\mathcal{C}^{\mathrm op}\to\mathsf{Cat}$},
Year = {2019},
Eprint = {arXiv:1908.02202},
}

@misc{spivak2022polynomial,
      title={Polynomial functors and Shannon entropy},
      author={David I. Spivak},
      year={2022},
      eprint={2201.12878},
      archivePrefix={arXiv},
      primaryClass={math.CT}
}

@misc{spivak2023functorial,
      title={Functorial aggregation}, 
      author={David I. Spivak},
      year={2023},
      eprint={2111.10968},
      archivePrefix={arXiv},
      primaryClass={math.CT}
}

@article{weber2007theory,
  Date = {2007},
  Journaltitle = {Theory and Applications of Categories},
  Volume = {18},
%  Number = {22},
  Pages = {Paper No. 22, 665--732},
  Author = {Weber, Mark},
  Title = {Familial 2-Functors and Parametric Right Adjoints},
}
\printindex

%\numnotes

\end{document}